\documentclass[10.6pt, reqno]{amsart}
\usepackage{tikz}
\usepackage{tabularx}
\usepackage{caption}
\usepackage{graphicx}
\usepackage{cancel}
\usepackage{subfigure}
\usepackage{slashed}
\usepackage{amsfonts}
\usepackage{amssymb}
\usepackage{accents}
\usepackage{amsmath}
\usepackage{amsxtra}
\usepackage{epstopdf}
\usepackage{latexsym}
\usepackage{mathrsfs}
\usepackage{leftidx}
\usepackage{tensor}
\usepackage{enumitem}

\newtheorem{theorem}{Theorem}[section]
\newtheorem{problem}[theorem]{Open Problem}
\newtheorem{lemma}[theorem]{Lemma}
\newtheorem{corollary}[theorem]{Corollary}
\newtheorem{proposition}[theorem]{Proposition}
\newtheorem{definition}[theorem]{Definition}

\newtheorem{assumption}[theorem]{Assumption}

\theoremstyle{remark}
\newtheorem{remark}[theorem]{Remark}

\numberwithin{equation}{section}

\newcommand{\Poincare}{Poincar\'e }

\def\sC{\mathscr{C}}

\def\sta#1#2{\stackrel{#1}{#2}}
\def\stc#1{\sta{\circ}{#1}}

\def\rp#1{^{\!(#1)}}

\def\sN{\mathscr{N}}

\def\sB{\mathscr{B}}

\def\fm#1{\mathfrak{m}[#1]}

\def\bp{\boldsymbol{\partial}}

\def\T{\mathcal{T}}

\def\sDc{\stackrel{\circ}\sD}

\def\tt{{t'}}

\def\bb{{\mathbf{b}}}

\def\Er{\mbox{Er}}

\def\sn{{\slashed{\nabla}}}

\def\sQ{\mathscr{Q}}

\def\eh{\hat{\eta}}

\def\zb{{\underline{\zeta}}}

\def\J{{\mathcal{J}}}

\def\M{{\mathcal{M}}}

\def\bT{{\textbf{T}}}

\def\bR{{\textbf{R}}}

\def\bd{{\textbf{D}}}
\def\ti{\tilde}

\def\bg{\mathbf{g}}

\def\hk{{\hat{k}}}

\def\beaa{\begin{eqnarray*}}
\def\eeaa{\end{eqnarray*}}

\def\ba{\begin{array}}
\def\ea{\end{array}}
\def\d{\delta}
\def\be#1{\begin{equation} \label{#1}}
\def \eeq{\end{equation}}

\newcommand{\nn}{\nonumber}

\def\l{\langle}
\def\r{\rangle}

\def\pih{\hat{\pi}}
\def\cir{\overset\circ}
\def\nn{\nonumber}
\def\S{{\mathcal S}}

\def\ud#1{\underline{#1}}
\def\zb{\ud{Z}}

\def\S2{{\mathbb S}^2}

\def\A{\mathcal {A}}

\def\E{{\mathcal E}}

\def\ze{{\zeta}}

\def\Lie{{\mathcal L}}

\def\tr{\mbox{tr}}

\def\D{{\mathcal D}}
\def\H{{\mathcal H}}

\def\N{{\mathcal N}}

\def\La{{\Lambda}}

\def\P{{\mathcal P}}

\def\c{\cdot}
\def\hot{\widehat{\otimes}}
\def\sig{\sigma}

\def\a{\alpha}
\def\b{\beta}

\def\ep{{\epsilon}}
\def\ve{{{\textbf{$\varepsilon$}}}}
\def\l{\langle}
\def\r{\rangle}
\def\ga{\gamma}
\def\Ga{\Gamma}

\def\O{\mathcal{O}}

\def\p{\partial}
\def\P{{\mathcal P}}

\def\nab{\nabla}
\def\hb{{\ud h}}
\def\fR{\mathfrak{R}}

\def\C{{\mathcal C}}

\def\Lb{{\underline{L}}}
\def\aaa{{\mathbf a}}

\def\div{\mbox{\,div\,}}
\def\curl{\mbox{\,curl\,}}
\def\divc{\sta{\circ}{\sl\div}}
\def\tr{\mbox{tr}}
\def\Tr{\mbox{Tr}}

\def\tir{{\tilde r}}

\def\wt{\widetilde}
\def\f14{\frac{1}{4}}
\def\f12{{\frac{1}{2}}}

\def\t1a{t^{-\frac{1}{a}}}

\def\bm{{\bf m}}
\def\sl{\slashed}
\def\sD{\slashed{\Delta}}

\def\sn{{\slashed{\nabla}}}

\def\eh{\hat{\eta}}

\def\zb{{\underline{\zeta}}}

\def\J{{\mathcal{J}}}

\def\M{{\mathcal{M}}}
\def\fY{\mathfrak{Y}}
\def\bT{{\emph{\bf{T}}}}

\def\bR{{\emph{\bf{R}}}}

\def\bd{{\emph{\bf{D}}}}
\def\ti{\tilde}

\def\hk{{\hat{k}}}

\def\beaa{\begin{eqnarray*}}
\def\eeaa{\end{eqnarray*}}

\def\ba{\begin{array}}
\def\ea{\end{array}}
\def\be#1{\begin{equation} \label{#1}}
\def \eeq{\end{equation}}
\def\nn{\nonumber}

\def\l{\langle}
\def\r{\rangle}

\def\pih{\hat{\pi}}
\def\cir{\overset\circ}
\def\nn{\nonumber}
\def\S{{\mathcal S}}

\def\S2{{\mathbb S}^2}

\def\A{\mathcal {A}}

\def\piS{{}\rp{S}\pi}
\def\piSh{{}\rp{S}\pih}
\def\piShb{{}\rp{S}\ud\pih}
\def\piT{{}\rp{\bT}\pi}
\def\pio{{}\rp{\Omega}\pi}

\def\pioh{{}\rp{\Omega}\pih}
\def\piohb{{}\rp{\Omega}\ud\pih}

\def\E{{\mathcal E}}

\def\ze{{\zeta}}

\def\Lb{\underline{L}}
\def\tr{\mbox{tr}}

\def\bA{{\emph{\textbf{A}}}}
\def\D{{\mathcal D}}
\def\H{{\mathcal H}}

\def\La{{\Lambda}}

\def\P{{\mathcal P}}

\def\c{\cdot}
\def\hot{\widehat{\otimes}}
\def\sig{\sigma}

\def\a{\alpha}
\def\b{\beta}

\def\l{\langle}
\def\r{\rangle}
\def\ga{\gamma}

\def\Ga{\Gamma}
\def\la{\lambda}

\def\p{\partial}
\def\P{{\mathcal P}}

\def\nab{\nabla}

\def\Lb{{\underline{L}}}
\def\aaa{{\mathbf a}}

\def\div{\mbox{\,div\,}}
\def\curl{\mbox{\,curl\,}}

\def\tr{\mbox{tr}}
\def\Tr{\mbox{Tr}}

\def\tir{{\tilde r}}

\def\wt{\widetilde}
\def\f14{\frac{1}{4}}
\def\f12{{\frac{1}{2}}}

\def\t1a{t^{-\frac{1}{a}}}

\def\bm{{\bf m}}
\def\sl{\slashed}
\def\sD{\slashed{\Delta}}

\def\ckk{\check}
\def\ckc{\ckk c}

\newcommand{\bea}{\begin{eqnarray}}
\newcommand{\eea}{\end{eqnarray}}

\def\nn{\nonumber}

\def\gaz{\gamma^{(0)}}

\newcommand{\chih}{\hat{\chi}}
\newcommand{\chib}{\underline{\chi}}

\newcommand{\chibh}{\underline{\hat{\chi}}\,}
\newcommand{\les}{\lesssim}
\newcommand{\ges}{\gtrsim}

\def\gac{{\stackrel{\circ}\ga}}
\def\Kc{{\stackrel{\circ}K}}
\def\thetac{{\stackrel{\circ}\theta}}

\def\fB{\mathfrak{B}}
\def\bN{{\mathbf{N}}}

\def\S{\mathcal{S}}

\def\cir#1{\stackrel{\circ}{#1}}

\def\vs{\varsigma}

\def\ss{\mathfrak{S}}

\def\bff{\mathbf{f}}
\def\bbf{\sta{\diamond}{\mathbf{f}}}

\def\sX{\mathscr{X}}
\def\sB{\mathscr{B}}
\def\sG{{\mathscr{G}}}

\def\sF{\mathscr{F}}

\def\bJ{{\mathbf{J}}}

\def\sP{{\mathscr{P}}}

\def\sQ{{\mathscr{Q}}}

\def\ud#1{\underline{#1}}

\def\fw{\mathfrak{w}}
\def\be{{(e)}}

\def\cb{\bar c}
\def\sdiv{{\sl{\div}}}
\def\scurl{\sl{\curl}}
\def\divsc{{\stackrel{\circ}{\sdiv}}}
\def\curlsc{{\stackrel{\circ}{\scurl}}}
\def\Delsc{{\stackrel{\circ}{\sD}}}
\def\snc{{\sta{\circ}\sn}}
\def\thetac{\sta{\circ}\theta}

\def\vc{{\sta{\circ}v}}

\def\Pic{\stackrel{\circ}\Pi}
\def\ckr{\ckk r}

\def\sk{\sl{\stc{k}}}

\def\al{\aleph}

\def\hN{{\hat\bN}}
\def\bAn{{\bA^\natural}}

\def\Sc{\mbox{Sc}}
\def\Ac{\mbox{Ac}}
	
\def\Osc{\mbox{\bf osc}}
\def\zsn{\sta{(0)}\sn}
\def\zga{\ga^{(0)}}
\def\tzga{\ti\ga^{(0)}}

\begin{document}
\title[]{On global dynamics of $3$-D irrotational compressible fluids}
\author{Qian Wang}
\address{
Oxford PDE center, Mathematical Institute, University of Oxford, Oxford, OX2 6GG, UK}
  \email{qian.wang@maths.ox.ac.uk}
  \date{\today}
\begin{abstract}
We consider global-in-time evolution of irrotational, isentropic, compressible Euler flow in $3$-D, for a broad class of $H^4$ classical Cauchy data without assuming symmetry, prescribed on an annulus surrounded by a constant state in the exterior. By giving a sufficient expansion condition on the initial data  and using the nonlinear structure of the compressible Euler equations, we show that the decay rate of the first order transversal derivative of the normalized density is better than that of the same derivative of a free wave, provided that the perturbation arising from the tangential derivatives can be properly controlled for all $t$ by using a bootstrap argument. Building on this critical analysis, we construct global exterior solutions in $H^4$ for the broad class of data, with a rather general subclass forming rarefaction at null infinity. Our result does not require smallness on the transversal derivatives of classical data, thus applies to data with a total energy of any size.
\end{abstract}
\maketitle
\tableofcontents
\section{Introduction}
Compression and expansion are fundamental physical phenomena in fluid mechanics, where compression refers to an increase in density due to a reduction in volume, and expansion refers to a decrease in density due to an increase in volume. These physical phenomena give rise to two types of solutions in the classical one-dimensional conservation law equations: shock waves and rarefaction waves. Shock occurs when characteristics intersect, leading to a blow-up in the distribution density of characteristics within finite time, while rarefaction corresponds to the density of characteristics vanishing. These phenomena  have been extensively studied in one dimension.

In the one-dimensional case, for a genuinely nonlinear 
 quasilinear hyperbolic system, singularity is known to form for compactly supported small data. This was proved by Lax \cite{Lax1, Lax2}, Glimm and Lax \cite{Glimm_Lax} for scalar equations and $2\times 2$ systems, and by John \cite{FJohn} for $n\times n$ system,  extended by Liu \cite{TPLiu} and others, due to the existence of compression in such data.  Sideris \cite{Sideris} gave the first blow up result for the compressible Euler equations in three
dimensions. One can refer to the works of Alinhac \cite{Alin_3}-\cite{Alin2} for results on formation of singularity for $2$-D compressible Euler equations and for quasi-linear wave equations. In particular for the latter, the results in Alinhac \cite{Alin1} and \cite{Alin2} were established due to the mechanism of the collapse of the foliation of characteristic hypersurfaces, without any symmetry assumptions on the data.  

In recent years, there have been significant developments in multi-dimensional analysis leading to the construction of small, smooth initial data that form shock within finite time for the 3-D relativistic Euler flow. This was demonstrated in Christodoulou \cite{shock_demetrios} by imposing a non-degenerate compression condition on the initial data constructed in an annulus. A similar result was proved for the irrotational isentropic compressible Euler equations in Christodoulou-Miao \cite{Miao_thesis}. These results have been extended to various scenarios:
shock formation for the geometric wave equation in $3$-D in Speck \cite{Spck_shock_1}; shock formation along the incoming acoustical cones in Miao-Yu \cite{Pin-Shuang} for some quasilinear wave equation verifying the classical null conditions; and for the $2$-D or $3$-D compressible Euler equations, including the non-irrotational case in Luk-Speck \cite{Spck-luk_2}  and the non-isentropic case in Luk-Speck \cite{Spck-luk_3}. In Buckmaster-Shkoller-Vicol \cite{Buck_2}-\cite{Buck_4},  by perturbing a Burgers shock, the authors
 construct shock with vorticity and entropy.  One can also find a different type of singularity, constructed in Merle-Rapha\"el-Rodnianski-Szeftel  \cite{Implo_1} and \cite{Implo_2} for the compressible $3$-D Navier-Stokes and Euler equations in the self-similar setting  under a
suitable regime of barotropic laws, with 
the density becoming infinity at the blow-up point.  

 For compressible Euler equations in $1$-D, one can refer to  Chen-Pan-Zhu \cite{Chen_G}  and Chen-Chen-Zhu \cite{Chen} for if-and-only-if results for singularity formation, classified by the existence of initial compression. In multi-dimensional cases, achieving a comprehensive classification result is out of reach. To complement the shock formation results in Christodoulou  \cite{shock_demetrios} and Christodoulou-Miao \cite{Miao_thesis} in this paper, we explore the global dynamics of the classical compressible Euler flow in $3$-D, initiated from compactly supported data (modulo constant states). In $3$-D, one has the favorable linear behavior of the free wave solution $\phi\approx 1/t$,  which lies at the border line between making nonlinear analysis work or fail, while such a decay property is absent in lower dimensions. 
Therefore for the global dynamics for the irrotational ideal gas flow in $3$-D, we propose the following:
\begin{problem}
 For the isentropic, irrotational compressible Euler equations in $3$-space-dimension \begin{footnote}{We focus on the perfect fluid with pressure satisfying $\gamma$-law for $\ga\ge 1$.}\end{footnote},  there is a certain class of classical Cauchy data, nontrivial in a unit ball surrounded by non-vacuum constant states, that may not satisfy any symmetry assumption, from which the classical solutions extend uniquely and globally in time, forming rarefaction with time approaching infinity.
 \end{problem}

 Unlike the $1$-D case, due to the decay property of the solution of the linearized problem, either confirming or disproving the above statement is open.  It is worth noting that in Serre \cite{Serre} and Grassin \cite{Grassin}, 
 the authors proved decent global existence results for expansion waves for the compressible Euler equations in the multi-dimensional case. In these works, global solutions were constructed as perturbations of the global solution of a force-free Burgers' equation when the density of fluid and its spatial derivatives are sufficiently small and the initial derivatives of velocity verify an expansion condition which forces the initial velocity to be infinite at the spatial infinity. Such an expansion condition in initial data obviously excludes the nontrivial data surrounded with constant states (and of non-zero density). We refer to Sideris \cite{Expansion1} and Hadžić-Jang \cite{Expansion2} for global expansion solution of Euler equations for free-boundary problems.

 Motivated by the open problem, we construct the global exterior solution of the irrotational compressible Euler equations by providing a sufficient condition on the classical initial data  prescribed in an annulus, without requiring the total initial energies to be small. Initiating from a standard radial foliation, we provide a further condition on data at a point in the initial slice, under which we prove that rarefaction occurs at the infinity of the corresponding null geodesic, as indicated by the distribution density of the characteristic surfaces vanishing to zero.
 
Rarefaction is often studied in the context of Riemann problems, particularly in one dimension. Setting it as the type of singularity in the initial data, the results in literature are usually to construct local solutions from such data. One can refer to Alinhac \cite{Alin} for the remarkable result in the general multi-dimensional case. A recent breakthrough on this topic can be found in Luo-Yu \cite{Pin_Luo1, Pin_Luo2} on the irrotational isentropic compressible Euler equation in $2$-D without losing derivatives, which also provides a deep geometric perspective. To clarify the difference, we emphasize that our goal is to form rarefaction of classical solution from classical (normal) Cauchy data, i.e. regular initial data that are not close to rarefaction, in a general multi-dimensional setting. 

The results of shock formation in both the irrotational relativistic and the irrotational compressible Euler equations in $3$-D in Christodoulou  \cite{shock_demetrios} and Christodoulou-Miao \cite{Miao_thesis} suggest that it is impossible for a global-in-time existence result for the exterior solution to hold for general classical data given on an annulus, surrounded by a non-vacuum constant state, even with smallness assumptions. When it comes to constructing a rarefaction wave, singularities are not expected to form in finite time. (We will provide a heuristic example to illustrate this shortly.) Hence to confirm the dynamic formation of rarefaction from Cauchy data, we have to construct global-in-time solutions, as well as proving the desired decay property for these solutions. The framework in \cite{shock_demetrios} and \cite{Miao_thesis} is devised for providing semi-global results in time for sufficiently small perturbation of the surrounding constant states. Achieving global-in-time solution for fluids in our setting is a significant challenge. We will give the model argument to explain the mechanism shortly. 

Starting from Klainerman \cite{K-commu} and Klainerman-Rodnianski \cite{KRduke, KR1}, there established a systematical geometric-analytical method for treating the lower regularity local well-posedness problem of quasi-linear wave equations. In Wang \cite{WangCMCSH, Wangrough, rough_fluid}, the geometric energy regimes were developed directly in the rough physical background. Thus, more geometric structures were uncovered in particular for treating quasi-linear wave and compressible Euler equations, whose metrics lack favorable Ricci curvature for geometric analysis. The lower order derivatives of the background geometric quantities exhibit better asymptotic behavior, which is advantageous for establishing energy inequalities without either losing decay or requiring higher order control on geometry. In this paper, to carry out the energy estimates, we crucially rely on the weighted energy regime in Wang \cite{Wangrough, rough_fluid}  which was established originally for solving the difficulty from the rough background geometry.  More importantly, we further extend the geometric structures obtained in \cite[Section 8]{rough_fluid} for the acoustical spacetime. For completing the global-in-time analysis, we need more detailed geometric structures of the acoustical metric to gain either better asymptotic behavior or smallness. For instance,  by virtue of the geometric structures given in Lemma \ref{dcom_s},  we gain the crucial hidden smallness in certain sense for the lowest order transversal derivative of the density (see (\ref{12.6.2.23}) in  Proposition \ref{12.21.1.21}).  Moreover, we provide the intrinsic formulation for the null forms in the Euler equations (see Proposition \ref{geonul_5.23_23}), and also identify null forms in Ricci decomposition (see Proposition \ref{ric44}), etc. These delicacies are in particular crucial for the global dynamics.

In this paper, we regard spatial velocity derivative $\p v$ as the second fundamental form $\f12\Lie_\bT \delta_e$, with $\delta_e$ the Euclidean metric in $\mathbb R^3$. From this perspective, influenced by the seminal work Christodoulou-Klainerman \cite{CK} for constructing global solutions to Einstein vacuum equations, we derive the Hodge-system and radial transport equations for various components of the spatial derivative of velocity. Using these equations, we gain the necessary transversal control or smallness for the asymptotic behaviors of various components of $\p v$ by merely assuming smallness for the initial  tangential (along characteristic cones) derivatives of the solution in $L^2$ and for the initial energy of angular derivatives of $\varrho$ and $\div v$. Together with our fully intrinsic formulation of null forms in the Euler equations (\ref{4.10.1.19}) and (\ref{4.10.2.19}) and the fact that our analysis is completed in $H^4$ for velocity and density, we reduce significantly the cumbersome analysis for geometric comparison and for controlling higher order derivatives in literature. More importantly, it does not require smallness assumptions on the total derivatives of the classical data, and the whole work is structured to cope with the lack of smallness in the total energy.

Next, we will begin by presenting the geometric set-up and the main result of this paper. We then will introduce the model problem we consider, and discuss the main steps, difficulties and ideas in our proof.
 \subsection{Basic set-up and the main result}\label{10.25.5.23}

We consider the compressible Euler equations of $3$ space dimension for a perfect fluid under a barotropic equation of state, that is, the
pressure $p$ is a function of the density  $\rho:{\mathbb R}^{1+3}\rightarrow (0,\infty)$,
\begin{equation}\label{10.12.1.19}
p=p(\rho).
\end{equation}
 We fix outside of the Euclidean sphere $\{|x|=R\}$ at $t=0$ a constant  background  density $\bar\rho>0$. Define the normalized density
\begin{equation}\label{10.12.2.19}
\varrho=\ln (\rho/\bar \rho)
\end{equation}
 and the sound speed
 \begin{equation*}
  c=\sqrt{\frac{dp}{d\rho}}.
 \end{equation*}
Clearly, due to (\ref{10.12.1.19}),   $c=c(\varrho)$. For convenience, we assume the pressure function
\begin{equation*}
p(\rho)=A\rho^\ga,
\end{equation*}
where the constant $A>0$, and $\ga\ge 1$, which ensures the constant  $\wp=c^{-1}c'+1=\frac{\ga+1}{2}\ge 1$. Keeping $\wp$ nonvanishing is important to our analysis.

Let $v$ be the velocity of the compressible fluid, $v:{\mathbb R}^{3+1}\rightarrow {\mathbb R}^3$. At $t=0$, let $\{|x|\le R\}$ be surrounded by the constant states $(\bar\rho, \bar v)$. Using the Galilean transform of coordinates $x^a\rightarrow x^a-{\bar v}^a t$, $v\rightarrow v-\bar v$
 we may assume without loss of generality that
\begin{equation}\label{cst}
\varrho=0,\quad c=c_*=(A \ga)^\f12 (\bar\rho)^{\frac{\ga-1}{2}}>0,  \quad  v=0, \mbox{ for }|x|\ge R, t=0
\end{equation} 
 and set $R=5$ for convenience. 
  
We define the acoustic metric $\bg$ as
\begin{equation}\label{metric}
\bg:=-dt\otimes dt+c^{-2} \sum_{a=1}^3(d x^a-v^adt)\otimes (dx^a-v^a dt).
\end{equation}

The inverse metric $\bg^{-1}$ can be written as
\begin{equation*}
\bg^{-1}=-\bT\otimes \bT+c^2 \Sigma_{a=1}^3 \p_a \otimes \p_a,
\end{equation*}
where $\bT$ is the future directed, time-like unit normal of $\Sigma_t$, i.e. the level set of $t$. And the component of  $\bg^{-1}$ will be denoted by
$\bg^{\a\b}$. \begin{footnote}{We adopt Einstein summation convention in this article. The range of  the indices of Greek letters such as
$\a,\b,\mu,\nu$ is $0,\cdots, 3$, and the range of the Latin letters $i,j,k,l,m,n,a,b$ is $1,2,3$, unless specified differently. We also fix the convention that
$\p_0=\p_t$.}\end{footnote}
 Relative to the Cartesian coordinates, $\bT$ is written as
\begin{equation*}
\bT=\p_t+v^a \p_a.
\end{equation*}

 The  induced metric on $\Sigma_t= {\mathbb R}^3\times \{t\}$ takes the form $g_{ij}=c^{-2}\delta_{ij}$ , where $\delta_{ij}$ is the kronecker
 delta. Define the second fundamental form
 \begin{equation*}
 k_{ij}=-\f12 \Lie_\bT g_{ij}, \qquad \Tr k =g^{ij} k_{ij},\qquad \hk_{ij}=k_{ij}-\frac{1}{3} \Tr k g_{ij}
 \end{equation*}
 where $\Lie_X$ denotes the Lie derivative by the vector field $X$.
Let $\cir{k}_{ij}=-\f12 \Lie_\bT \delta_{ij}$. Thus $\Tr \cir{k}:=\delta^{ij} \cir{k}_{ij}= -\p_i v^i$. 

Now we introduce the compressible Euler equations with (\ref{10.12.1.19}) for $\varrho$ and $v$,
\begin{equation}\label{4.23.1.19}
\left\{
\begin{array}{lll}
\bT \varrho=-\div v\\
\bT v^i=-c^2\delta^{ia} \p_a \varrho,
\end{array}
\right.
\end{equation}
where $\div v= \p_i v^i$ and  $\varrho$ is the normalized density function in (\ref{10.12.2.19}).

 Let $\tensor{\ud\ep}{_i^j^k}$, $i,j,k=1,2,3$, be the standard volume form on ${\mathbb R}^3$. We define the vorticity to be
$
\fw_i=\tensor{\ud\ep}{_i^j_k}\p_j v^k.
$  \begin{footnote}{ The indices of the tensor fields and differentiation are lifted and lowered by the Euclidean metric.}\end{footnote}

In the isentropic irrotational case, $\fw=0$ for all $t$. The compressible Euler equations (\ref{4.23.1.19}) can be reduced to
\begin{align}
& \Box_\bg v^i= \sQ^i, \label{4.10.1.19}\\
& \Box_\bg \varrho= \sQ^0,\label{4.10.2.19}
\end{align}
where $\Box_\bg$ is the Laplace-Beltrami operator of the Lorentzian metric $\bg$, and the  quadratic terms are
\begin{align*}
&\sQ^i:=-(1+c^{-1}c')\bg^{\a\b}\p_\a\varrho\p_\b v^i, \\
&\sQ^0:=-3c^{-1}c' \bg^{\a\b}\p_\a \varrho \p_\b \varrho+2\sum_{1\le a<b\le 3}\big(\p_a v^a\p_b v^b-\p_b v^a \p_a v^b\big).
\end{align*}
See the equations from the work of Luk-Speck \cite[Page 13]{Jared_Luk}. Different from the known treatment in literature, c.f. Luk-Speck \cite{Spck-luk_3}, in this paper, we will give an invariant formulation for $\sQ^i$ and $\sQ^0$, 
see Proposition \ref{geonul_5.23_23}. This observation plays a crucial role for our analysis together with the structure equations for $\p v$ components derived in Lemma \ref{dcom_s}.  

Define the optical function $u$ to be the solution of the Eikonal equation
\begin{equation}\label{optical}
\bg^{\a\b}\p_\a u \p_\b u=0.
\end{equation}
We refer to the level sets of the optical function $u$, denoted by $\H_u$, as the acoustical null cones. The global optical function $u$ is
constructed as follows.

 For each  $\omega\in {\mathbb
S}^2$,  define the null vector field $L'$  to be the generator of  the null geodesic $\Upsilon_\omega$ in the acoustical spacetime by
\begin{equation}\label{6.29.2.19}
\bd_{L'} {L'}=0, \quad \frac{d}{ds}\Upsilon_\omega(s)=L',\quad  L'(s)=1.
\end{equation}
 We denote $\bb^{-1}=-\l L',\bT\r$. $L=\bb L'=\bT+\bN$ where $\bN=\frac{\nab u}{|\nab u|}$ with $\nab=\nab_g$ the Levi-Civita connection of $g$. Thus
 \begin{equation*}
 \bT u=-\bb^{-1}=-\bN u.
 \end{equation*}


On the initial slice, let $u=r$ and set $u_*=5$ at $r=5$, and let $2\le u_0<u_*$. It is straightforward to obtain
\begin{equation*}
\bb=c^{-1}, \quad\mbox{at } t=0,\quad u_0\le u\le u_*.
\end{equation*}
By abuse of notation, we set $\Sigma_t=\{t'=t, u_0\le u\le u_*\}$, since our analysis is only on the exterior region,  and set 
\begin{equation}\label{b0}
\bb_0:=\inf_{\Sigma_0} c^{-1}.
\end{equation}

 We define $S_{t,u}=\Sigma_t\cap \H_u$.
For each $u$ the geodesic generator (\ref{6.29.2.19}) defines a smooth one-to-one mapping from $S_{0,u}$ to $S_{t,u}$. Moreover since $S_{0,u}=\{|x|=u\}$, we can assign any point $p$ in $S_{t,u}$ a pull back coordinates $\omega=(\omega_1,\omega_2)\in \mathbb S^2$, namely by following the null geodesic $\Upsilon_{\omega, u}$ originated at $q\in S_{0,u}$ with $q=(u,\omega_1, \omega_2)$.  This defines the acoustical coordinates for the point $p$, i.e. $p=(t,u, \omega_1, \omega_2)$. We note the integral curves of $L$ are curves being constant in $\omega$ and $u$, along which
\begin{equation}\label{11.27.1.23}
L=\frac{\p}{\p t}.
\end{equation}

Let $\ga$ be the induced metric on $S_{t,u}$. We denote by $\sn$ its Levi-Civita connection, and the area  $|S_{t,u}|_\ga=\int_{S_{t,u}} 1 d\mu_\ga$. Similarly, we denote  by $\snc$ and $ |S_{t,u}|_\gac$ the Levi-Civita connection and the area  with respect to $\gac=c^2\ga$.
For scalar functions $f$, we denote the spherical mean of $f$ by  $\bar f(t,u)=\frac{1}{|S_{t,u}|_{\gac}}\int_{S_{t,u}} f d\mu_\gac$, and denote 
$$\ckc=\cb(0,u), \quad \tir:=t+(\ckc)^{-1}u, \quad t\ge 0.$$

Let ${}\rp{a}\O=\tensor{\ud\ep}{^a_j_k}x^j\p_k, a=1,2,3$ and $\Pi^i_j=\delta_j^i-\bN^i\bN_j$. Define
\begin{equation}\label{5.14.2.23}
{}\rp{a}\Omega^m={}\rp{a}\O^k\Pi_k^m,\quad S=\tir L.
\end{equation}
Let $\p$ represent the spatial derivative $\p_i, i=1,2,3$, and $\bp$ stand for $\p$ and $\bT$.
 We will regard $\Phi^0=\varrho,$ and $\Phi^i=v^i$ for convenience. $\bp\Phi$ means applying $\bp$ to  the scalar component of $\Phi$. Other derivatives on $\Phi$ are understood similarly. It is important to distinguish different components of $\bp \Phi$. For this purpose, we adopt the convention for the scalar components:
 \begin{align*}
 [X \Phi]=X v^i\bN^j \bg_{ij} , X\varrho
 \end{align*}  
 where $X$ is a linear differential operator and for the repeated indices we adopt the Einstein summation convention to sum over $i, j=1,2,3$. 
 
 By the finite speed of propagation property, exterior to $\H_{u*}$, the spacetime is fully constant, with the values of $(\varrho, c, v)=(0, c_*, 0)$ as in (\ref{cst}).
 Throughout this paper, we suppose the initial data is free of vorticity in $\Sigma_0$.\begin{footnote}{In the domain of dependence external to $\{u=u_0\}$, vorticity vanishes wherever the classical solution extends to. Thus it suffices to consider $\varrho$ and $v$ as the solutions of (\ref{4.10.1.19}) and (\ref{4.10.2.19}).}\end{footnote} Moreover, for the initial data we make the following boundedness assumption.
\begin{assumption}\label{A1}
Let $2\le u_0<5$ be fixed. For $\Phi=\varrho, v$, suppose at $t=0$, there exists a positive constant $\A_0>0$
\begin{equation}\label{1.12.1.22}\tag{\bf A1}
\int_{\Sigma_0} |\bp  \bp^{\le 3}\Phi|^2 dx\le \A_0^2<\infty\mbox{ at } t=0
\end{equation}
and assume
\begin{equation}\label{rarefied}\tag{\bf A2}
q_0=\inf_{u\in [u_0, u_*], \omega\in \mathbb S^2} (-\ckc^{-1}u\Lb\varrho+\varrho)(0, u, \omega)>0, 
\end{equation}
where $\ckc$ denotes the spherical mean of $c$ on $\{|x|=u\}$.
\end{assumption}
 
Note that $\bb^{-1}$ is the density of the foliation of the null cones $\H_u$. At the point where a shock occurs, $\bb^{-1}\rightarrow \infty$. On the contrary, $\bb^{-1}\rightarrow 0$ as rarefaction occurs.  In this paper we will show for a certain class of initial data  the solution spacetimes can be uniquely extended for all $t>0$ and $u_0\le u\le u_*$. And in a subclass class of such spacetimes, $\bb^{-1}\rightarrow 0$ as $t\rightarrow \infty$ along null geodesics, which confirms that the phenomenon of rarefaction does occur at the null infinity even if $\bb(0)\approx 1$.

\begin{theorem}\label{mainthm1}(Main theorem) Let $\Phi$ be $(v^i, \varrho)$ in the equation system (\ref{4.10.1.19})-(\ref{4.10.2.19}). Let the initial data  of $\Phi$ satisfy (\ref{1.12.1.22}) and (\ref{rarefied}). Suppose 
\begin{equation}\label{9.22.1.22}\tag{\bf A3}
\|X^{1+\le 3}\Phi\|_{L^2(\Sigma_0)}+\sum_{\ell=0}^1\sum_{m=\ell}^2\|\bp\Omega^{1+\le m-\ell}\bT^\ell\varrho\|_{L^2(\Sigma_0)}\le \ve\le\ve_0
\end{equation}
 where $X$ represents all the vector fields in $\{S, \Omega\}$, and $\ve_0$ is sufficiently small, depending on $\A_0, c_*$. It yields the following rough bound for $q_0$, \begin{footnote}
{
For  $A\les B$, we mean $A\le C B$ where  $C>0$ is a universal constant. Here if a constant is called  universal, it means the constant depends only on $c_*, \A_0$. $A\ges B$ is understood in the similar way.}
\end{footnote} $$q_0\les\ve^\f12.$$ 

There hold the following results:  
 \begin{itemize}
\item[(1)]  There exists a universal constant $C_1>1$ such that as long as the data verify
\begin{equation}\label{exist}\tag{\bf A4}
q_0>C_1\ve, 
\end{equation}
the solution $(\varrho, v)$ for the irrotational compressible Euler equation system (\ref{4.10.1.19})-(\ref{4.10.2.19}) exists globally and uniquely for all $2\le u_0\le u\le u_*$ and $t>0$.

\item[(2)] Under the assumptions of (\ref{1.12.1.22})-(\ref{exist}), there exists a universal constant $C_2$, if at some $u_1, \omega_1$ with  $u_0\le u_1\le u_*$ and $\omega_1\in \mathbb S^2$,  the initial density verifies the condition
\begin{equation}\label{rarif}\tag{\bf A5}
\ckc^{-1} c^{-1}r\Lb \varrho(0,u_1, \omega_1)<-C_2 \ve^\f12,
\end{equation}
then the null lapse $\bb\rightarrow\infty$ as $t\rightarrow \infty$ along the null geodesic $\Upsilon_{\omega_1, u_1}(t)$, which indicates the rarefaction occurs at $(u_1,\omega_1)$ of the null infinity. More precisely, with constant $C''>0$  and the constant $C'\approx \A_0$, 
\begin{equation*}
\bb(t, u_1, \omega_1)\ges \ve^\f12\big(\frac{1}{C'}+ C'' \wp\log (\f12\l t\r)\big).
\end{equation*}
\end{itemize}
\end{theorem}
\begin{remark}
Technically, the rough upper bound of $q_0$  and the bound assumed in (\ref{rarif}) can both be slightly refined, with a more refined analysis using the same method we employed to prove our main result. 
\end{remark}
\begin{remark}
 For convenience, we will assume $\A_0=1$ throughout the proof, which can be achieved by rescaling.
\end{remark}
\begin{remark}
Detailed asymptotic behaviors of the solution can be found in Proposition \ref{1steng}, Proposition \ref{8.29.8.21}, Proposition \ref{9.8.6.22}, Proposition \ref{10.30.4.21}, Proposition \ref{3.14.4.24} and Proposition \ref{imp_decay} and Proposition \ref{10.13.5.23} with $\La_0, \Delta_0$ therein bounded by $C\ve$. Here $C>1$ is a universal constant.  
\end{remark}
\begin{remark}\label{5.19.1.23}
 The second assumption  in (\ref{9.22.1.22}) is only for $\varrho$. We can not assume the same estimate for $v$, since it is incompatible with the fact that $\Lb\Phi, \Lb^2\Phi$ may not be small in our data.  Indeed the angular component of $\p_r \Omega v$ contains the term $\p_r v(\p_r)\c r^{-1}\Omega$ (see (\ref{5.26.1.23})), which could be large at $t=0$. Thus we do not have small standard energy on $\Omega v$  even at the initial slice. 
  Under the assumptions (\ref{1.12.1.22})-(\ref{9.22.1.22}), we will give the hierarchy on the amplitude of various derivatives of $\Phi$ at $t=0$ in Table (\ref{5.24.1.23}), proved in Section \ref{9.25.2.22}.
\end{remark}
In what follows we will justify the assumptions made in Theorem \ref{mainthm1}. For this purpose, we first will give a simplest example.
\subsection{Motivation}
\subsubsection{A heuristic example}
Consider the initial value problem of Burgers' equation in $1$-D.
$$
\phi_t + (\phi^2/2)_x =0, \quad \phi(x, 0) = \phi_0(x) :=\left\{ \begin{array}{lll}
1,  & x<0,\\
1-x,  & 0\le x\le 1,\\
0, & x>1.
\end{array}\right.
$$

The characteristic line crossing $x$-axis at $x_0$ is given by
$$
x(t)= x_0 + t \phi_0(x_0), \qquad  x_0 \in {\mathbb R},
$$
and on this line
$$
\phi = \phi_0(x_0).
$$
Since all characteristics starting at $(x_0, 0)$ with $0\le x_0\le 1$ cross at $(1, 1)$, where they form  a shock, $\phi(x, t)$ can not be smooth
for $t\ge 1$.

 Therefore, for $t<1$,
\begin{equation}\label{4.11.1.23}
\phi(x, t) = \left\{ \begin{array}{lll}
1, & x<t,\\
(1-x)/(1-t), & t\le x\le 1,\\
0, & x>1.
\end{array}\right.
\end{equation}

In order to use the geometric formulation in (\ref{optical}) for the characteristics, heuristically, we reduce the Burgers' equation to a second order equation by a straightforward calculation:
 \begin{equation*}
 -\phi_{tt}+\phi^2 \phi_{xx}=-\phi\c \phi_x^2+\phi_t \phi_x.
 \end{equation*}
Hence, the second order differential operator $-\p_t^2+\phi^2\p_x^2$ on the left-hand side differs by quadratic lower order terms from the Laplace-Beltrami operator in  a spacetime $\mathbb R\times [0,1)$ equipped with metric
$$
\bg_{\mu\nu}=\mbox{diag}(-1,\phi^{-2}), \,\mu,\nu=0,1.
$$
We can adopt this heuristic geometric formulation wherever $\phi(t,x)\neq 0$.

In line with our definition (\ref{optical}), the optical function $u$ for Burgers' equation takes the form
 $$
 -(\p_t u)^2+\phi^2(\p_x u)^2=0.
 $$
 To construct the optical function, we usually need to prescribe initial or boundary conditions. For Burgers' equation, we can explicitly write down the optical function such that $u(x,0)=x$. Corresponding to (\ref{4.11.1.23}) we have
 \begin{itemize}
 \item In the region $\phi=1$, $u=x-t\le 0$, $\bb^{-1}=1$.
 \item In the region $\phi=0$, $u=x\ge 1$, $\bb^{-1}=0$.
 \item In the region $\phi=\frac{1-x}{1-t}$, $u=\frac{x-t}{1-t}=x_0\in(0,1)$, $\bb^{-1}=\frac{1-x}{(1-t)^2}=\frac{1-u}{1-t}$.
  As $t\rightarrow 1$ along the characteristics,  $\bb^{-1}\rightarrow\infty$.
 \end{itemize}
 Clearly, the shock forms at the point $(1,1)$ where $\bb^{-1}$ blows up. (See FIGURE 1.)
 
 To create a counter scenario to the above solution, we give the following initial data to the Burgers' equation
 $$
 \phi(x, 0) = \phi_0(x) :=\left\{ \begin{array}{lll}
0,  & x<0,\\
x,  & 0\le x\le 1,\\
1, & x>1.
\end{array}\right.
$$
Following the characteristic approach, we can derive the solution
 $$
\phi(x, t) = \left\{ \begin{array}{lll}
0, & x<0,\\
\frac{x}{1+t}, & 0<x\le t+1,\\
1, & x>t+1.
\end{array}\right.
$$
Moreover, setting $u(x, 0)=x$, we have (see FIGURE 2)
\begin{itemize}
\item In the region $\phi=0$, $u=x$, $\bb^{-1}=0$.\\
\item In the region $\phi=\frac{x}{1+t}$, $u=\frac{x}{1+t}=x_0\in (0,1)$, $\bb^{-1}=\frac{u}{1+t}$, $\bb^{-1}\rightarrow 0$ as $t\rightarrow \infty$ along the characteristics.
\item In the region $\phi=1$, $u=x-t$, $\bb^{-1}=1$.
\end{itemize}
\begin{figure}[ht]
 \begin{minipage}[t]{.5\linewidth}
\begin{tikzpicture}[scale=1.1]
\draw[->] (-0.4,0) --(4,0) node[right] {$x$};
\draw[->] (0,-0.4) --(0,3) node[above] {$t$};
\draw (0,0)--(1,1);
\draw (1,1)--(2,3) node[right] {$x = 1 + \frac{t-1}{2}$};
\draw (1,0)--(1,1);
\draw [densely dashed] (0,1)--(3.5,1) node[right] {$t = 1$};
\node at(2,0.5) {$\phi = 0$};
\node at(2.5,2) {$\phi = 0$};
\node at(0.8,2) {$\phi = 1$};
\node at(0.4,0.78) {\tiny $\phi = 1$};
\node at(2, -0.4) {$\phi(x,t) = \frac{1-x}{1-t}$};
\draw[->] (1,-0.3)--(0.8,0.2);
\draw [densely dashed] (0.45,0)--(1,1);
\node at(0.44,-0.2) {$x_0$};
\end{tikzpicture}
\caption{Shock formation}
\end{minipage}\hfill
 \begin{minipage}[t]{.5\linewidth}
\begin{tikzpicture}[scale=0.95]
\draw[->] (-0.7,0) --(5.5,0) node[right] {$x$};
\draw[->] (0,-0.4) --(0,4) node[above] {$t$};
\node at(-0.2,-0.2) {$0$};
\node at(-0.6, 1.4) {$\phi=0$};
\draw (1,0)--(5,4);
\node at(4,1.4) {$\phi = 1$};
\node at(1,-0.2) {$1$};
\draw[->] (1/4,0)--(5/4,4);
\node at(1/2,-0.2) {$x_0$};
\draw[->] (1/2,0)--(5/2,4) node[above] {\tiny $x = x_0(1+t)$};
\draw[->] (3/4,0)--(14/4,4);
\node at(-0.8,3.6) {$\phi = \frac{x}{1+t}$};
\draw[->] (-1.0,3.4)--(0.5, 2.5);
\end{tikzpicture}
\caption{Formation of rarefaction}
 \end{minipage}
\end{figure}
 

This is an elementary example of Burgers' equation in $1$-D showing that rarefaction forms along characteristics as time approaches infinity for continuous data. Our goal is to extend this example to $3$-D irrotational compressible Euler flow,  without assuming any symmetric property in data. Theorem \ref{mainthm1} shows the classical exterior solution exists uniquely and globally and rarefaction can form at the null infinity of the acoustical spacetime in $3$-D for a rather general class of classical data, with the initial null lapse satisfying $\bb^{-1}\approx 1$.
\subsubsection{Remarks on the assumptions in Theorem \ref{mainthm1}}\label{3.24.4.24}

In $3$-D, Christodoulou  \cite{shock_demetrios} and Christodoulou-Miao \cite{Miao_thesis} provide sufficient conditions for shock formation in a finite time, which requires a non-degenerate compression assumption on small and smooth data defined in a proper sense.

 Here we briefly explain the mechanism in their works, by using the quantities and notations in our work (for ease of exposition). Based on the key observation that $L(\bb\tir \Lb\Phi)$ is a quantity linked to perturbations, the quantity $\bb \tir\Lb \Phi$ is expected to be of the same amplitude as initially. Assuming smallness of the initial energies up to very high order derivatives of all directions, propagating the smallness of the initial energies up to where the first singularity is expected to occur, it gives the semi-global existence result and  more importantly  the necessary decay properties of $\bb\bp\Phi$. While $\Lb \Phi$ blows up at the singularity, $\bb\tir \Lb\Phi$ stays close to its initial value due to the aforementioned key observation. Assuming $\bb\tir\Lb \Phi$ to be non-degenerate initially, then $\bb\tir \Lb\Phi$ is proved to be a finite non-zero quantity when approaching the singularity along characteristics. This forces $\bb\rightarrow 0$, which justifies the singularity is indeed a shock. The main difficulty in \cite{shock_demetrios} and \cite{Miao_thesis} is to achieve long-time energy propagation, with the key quantities, $\Lb \Phi$ and its derivatives, which determine the dynamics, blowing up as approaching the singularity. To balance with the blowing up terms, it is important to keep the vanishing quantity $\bb$ paired to energies properly throughout the long-time analysis. 
 
 Different from the above set-up, our assumptions on initial total energies are merely up to the order of $H^4$ (for $\Phi$) and of arbitrary size. To construct the global exterior solution in $3$-D, we certainly have to work in the complement of the data set of \cite{Miao_thesis}. 
      (\ref{1.12.1.22}) is a standard assumption even for local-in-time energy propagation.  (\ref{rarefied}) and (\ref{exist}) are expansion type conditions, which exclude initial compression.  We will show that the evolution of both $\bb\tir\Lb \varrho$ and the quantity $\tir \Lb\varrho-\varrho$ can be determined by their initial data with error contributed by the tangential perturbations. Inspired by the fact that solutions of $1$-D Burgers' equations are constant along characteristics, we therefore impose the smallness assumption (\ref{9.22.1.22}) merely on $L^2$ norm of tangential (to the characteristic surfaces) derivatives of $\Phi$ and the energies of the angular derivatives of $\varrho$ and $\div v$ at $t=0$, consistent with the order of $H^4(\Sigma_0)$. Combining (\ref{rarefied}) and (\ref{9.22.1.22}), we will show the bound $q_0\les \ve^\f12$ in Proposition \ref{12.21.1.21}.  In our setting  and with the mechanism to be explained shortly, for data satisfying (\ref{1.12.1.22})-(\ref{exist}), we can show that the transversal derivative of the solution decays at a rate bounded by $\frac{1}{\l t\r\log \l t\r}$ with $\l t\r=t+2$, although it does not necessarily maintain its initial amplitude. This decay property, in principle, allows us to propagate the tangential perturbations given in (\ref{9.22.1.22}) for all $t$ by means of a bootstrap argument, which is the main task of our analysis. The  condition (\ref{rarif}) is posed for guaranteeing the formation of rarefaction, analogous to the conditions for shock formation. We remark that, from the assumption (\ref{9.22.1.22}) and due to the sign requirement in (\ref{rarefied}), we can gain some hidden smallness for $\Lb \varrho(0)$ in $L_u^2 L_\omega^\infty$ initially (weaker than the amplitude of the tangential perturbation assumed in (\ref{9.22.1.22})). No smallness can be gained for the higher order general energies in (\ref{1.12.1.22}). 
 
In what follows, we will run the model argument on the asymptotic equation for $\varrho$ \begin{footnote}{We refer to \cite{Lind_Rod} for the definition of the asymptotic equation.}\end{footnote} to explain the mechanism for proving the main result.


\subsection{The model argument}\label{5.9.1.23} 
The breakdown of $C^1$ solutions for the equations of $1$-D hyperbolic system  has been extensively studied over decades. The approach largely is based on treating the Riccati equations of functions constructed from the Riemann invariants, see Lax \cite{Lax1} for instance. To generalize this idea in $3$-D, it is  important to find the proper evolution equation to work on. 

\subsubsection{Global exterior solution} The model argument is based on analysing the evolution of the solution $y\approx (t+u)\Lb \varrho$  of the asymptotic equation of (\ref{4.10.2.19}), which will be given symbolically shortly. The asymptotic behavior of $y$ that we are about to show is better than the asymptotic behavior of $r\Lb\phi$ with $\phi$ being free wave. It is achieved due to the initial condition and the nonlinear structure of Euler equation, with the latter modeled by the asymptotic equation as a vastly simplified form.  Controlling this quantity is the key to carry out energy estimates in our work.

Suppose the classical solution of the asymptotic equation breaks down first at some $0<t_*<\infty$.  Note,  due to (\ref{1.12.1.22}) and (\ref{rarefied})
  \begin{equation}\label{9.23.5.23}\tag{\bf H1}
-\varrho< -y(t)\le \M_0 \big(1+\f12\wp\log(\frac{\l t\r}{2})\big)^{-1} 
\end{equation}
holds at least at $t=0$, where $\M_0>1$ is a fixed constant comparable to $\A_0$, depending on $c_*$, to be specified.
Let $T_*=\sup\{t': (\ref{9.23.5.23}) \mbox{ holds  on } (0, t')\}.$ Due to continuity, $T_*>0$. 
 In the sequel, we will bound $\varrho=O(\ve^\f12 \l t\r^{-1+\delta})$ on $[0, T_*)$ and we can prove $T_*=t_*$, as long as the energy bound $\ve_0$ in (\ref{9.22.1.22}) is sufficiently small. This implies that $t_*=\infty$. Next we sketch he proof in three steps.

{\bf Step 1} Suppose $T_*<t_*$, then on $[0, T_*)$, (\ref{9.23.5.23}) holds. We reduce from the wave equation (\ref{4.10.2.19}) for $\varrho$ to the asymptotic equation for $y\approx (t+u)\Lb\varrho$, which takes the form schematically \begin{footnote}{The actual factor appeared after $y^2$ is closer to  $(t+u)^{-1}$ than to $\l t\r^{-1}$. We use the latter in the model equation to keep the argument simple.}\end{footnote}:
\begin{equation}\label{9.10.1.22}
y'-\f12 \wp y^2 \l t\r^{-1}=G
\end{equation}
where $\l t\r=t+2$; and for the null lapse $\bb$, we have
\begin{equation*}
(\log \bb)'=-\f12 \wp y(t+u)^{-1}+G_1,
\end{equation*} 
 with $|G, G_1|< C_0\ve^\f12\l t\r^{-\frac{7}{4}+\delta}$, $C_0$ universal. Here and in the sequel,  $\delta>0$ is a fixed constant, any close to $0$. We can run analysis at $t=0$ to see $\bb_0>0$ merely depends on $c_*$ and $\A_0$. Using the above equation for $\log \bb$ and (\ref{9.23.5.23}),  provided that $\varrho=O(\ve^\f12 \l t\r^{-1+\delta})$, and as long as $\ve_0>0$ is sufficiently small, chosen depending only on $\A_0$ and $c_*$, we will obtain that $\bb>\frac{1}{4}\bb_0$ on $(0,T_*)$. 
  
  We remark that the presence of $\l t\r^{-1}$ in (\ref{9.10.1.22}) is due to volume element in $3$-D in the nearly-spherical model. This factor is absent in $1$-D and the plane-wave scenario. It plays a crucial role in our evolution argument. 

 {\bf Step 2}
 Next, with\begin{footnote}{The actual version is more complicated since we have to pair $y$ with an integrating factor. }\end{footnote} $\ti y=y-\varrho$, \begin{footnote}{Here we ignore the factor of $\ckc$ in the definition of $q_0$ for simplicity, in order to give the heuristic argument.}\end{footnote}  we show there is a positive-valued function $\C(t,u)$ decreasing in $t$ with $u$ fixed, and a positive function $f(t)=C \ve\l t\r^{-\frac{3}{4}+\delta}$ with $C>0$ a universal constant  such that   
 \begin{equation}\label{9.23.2.23}
 -\ti y>f(t)+\C(t,u)>0, \quad 0<t<T_*.
 \end{equation}
 Indeed, in $(0,T_*)$, again from (\ref{4.10.2.19}) we derive 
 \begin{equation}\label{9.28.1.23}
 \ti y'=\f12\wp {\ti y}^2\l t\r^{-1}+G_2
 \end{equation}
 with $G_2$ being the refined error terms. 
 Taking the tangential perturbation with the bound in (\ref{9.22.1.22}) being $\ve$, we will obtain
 \begin{equation}\label{9.23.4.23} 
 |G_2|<\f12 C {\ve}\l t\r^{-\frac{7}{4}+\delta}. 
 \end{equation}
Due to the assumptions (\ref{rarefied}) and (\ref{exist}), we can choose $\a_1(t)=2C \ve\l t\r^{-\frac{3}{4}+\delta}$ such that 
\begin{equation*}
\ti y+\a_1(0)<0.
\end{equation*} 
This is achievable as long as the constant in (\ref{exist}) is greater than $2C$ in the above.  
From (\ref{9.28.1.23}), we derive
\begin{align*}
(\ti y+\a_1)'&=\f12\wp\{(\ti y+\a_1)^2-2\a_1(\ti y+\a_1)+\a_1^2\}\l t\r^{-1}+\a_1'+G_2\\
&=\f12\wp\{(\ti y+\a_1)^2-2\a_1(\ti y+\a_1)\}\l t\r^{-1}+\a_1'+G_2+\f12 \wp \a_1^2 \l t\r^{-1}.
\end{align*}
Note that it follows due to (\ref{9.23.4.23}) that $$\a_1'+G_2+\f12 \wp\a_1^2\l t\r^{-1}=2(-\frac{3}{4}+\delta)C \ve \l t\r^{-\frac{7}{4}+\delta}+2 \wp C^2{\ve}^2 \l t\r^{-\frac{5}{2}+2\delta}+G_2<0.$$
Hence with $z=\ti y+\a_1$, we derive
\begin{equation*}
z'<\wp(\f12 z^2-\a_1 z)\l t\r^{-1}.
\end{equation*} 
Suppose $t_1=\sup\{t: t\le T_*, z<0 \mbox{ on } (0, t)\}$. It follows by continuity that $0<t_1\le T_*$. In $(0, t_1)$, 
since $z<0$, 
\begin{equation*}
-z'z^{-2}> -\f12 \wp(1-2\a_1 z^{-1}) \l t\r^{-1}. 
\end{equation*}
This implies
\begin{align*}
z^{-1}\exp\int_0^t- \wp \a_1 \l t'\r^{-1}>z^{-1}(0)-\f12\wp\int_0^t\l t'\r^{-1} \exp(\int_0^{t'} -\wp \a_1 \l t''\r^{-1}) dt'.  
\end{align*}
In $(0, t_1)$, we have 
\begin{align*}
-z \exp\int_0^t\wp \a_1 \l t'\r^{-1}&>\Big(-z^{-1}(0)+\f12\wp\int_0^t\l t'\r^{-1} \exp(\int_0^{t'} -\wp \a_1 \l t''\r^{-1}) dt'\Big)^{-1}\\
&>(-z^{-1}(0)+\f12\wp\int_0^t\l t'\r^{-1})^{-1}.
\end{align*}
Note that there is a constant bound $\ti C>0$ such that  $1<\exp\int_0^t\wp \a_1 \l t'\r^{-1}\le \ti C$. 
Denoting by $\ti C\c \C(t,u)$ the function on the right-hand side,  we conclude $-z>\C(t,u)$, i.e. 
\begin{align*}
\ti y+\a_1<-\C(t,u)\le-\C(T_*,u)<0, 
\end{align*}
where  the last inequality holds since $0<t\le t_1\le T_*<\infty$. Hence $t_1=T_*$. Thus (\ref{9.23.2.23}) is proved, which improves the left-hand inequality in (\ref{9.23.5.23}) if $t<T_*$. It remains to improve the right-hand side inequality. 
 
Moreover, in view of $L-\Lb=2\bN$ and $\bb\bN\approx \p_u$ with a shift term neglected for simplicity, we will integrate the above result on $[u, u_*]$ to derive,  with some error $G_3$ verifying $\|G_3\|_{L^2_u}\les \l t\r^{-1+\delta}\ve$, that 
\begin{equation*}
\varrho(u)\les \int_u^{u_*}\{-\l t\r^{-1}(\C(t,u')+\a_1(t))+G_3\}.
\end{equation*}
Hence, we conclude  that
\begin{equation}\label{11.11.1.23}
- y\ges (\inf_{u_0\le u\le u_*}\C(t,u)+\a_1(t))\l t\r^{-1}+\a_1(t)-O(\l t\r^{-1+\delta}\ve), \, \forall\, 0<t<T_*.
\end{equation}

{\bf Step 3.} Next we improve the upper bound in (\ref{9.23.5.23}). Let $\a=2 C_0\ve^\f12 \l t\r^{-\frac{3}{4}+\delta}$, with $C_0$ specified in {\bf Step 1}. We derive from (\ref{9.10.1.22}) that
\begin{align*}
(-y+\a)'&=-\f12 \wp y^2\l t\r^{-1}-G+\a'\\
&=-\f12\wp\left((-y+\a)^2-2(-y+\a)\a+\a^2\right)\l t\r^{-1}-G+\a'.
\end{align*}
Since $|G|<C_0\ve^\f12 \l t\r^{-\frac{7}{4}+\delta}$, we derive
\begin{equation}\label{9.10.5.22}
\a'-G-\f12 \wp\a^2\l t\r^{-1}=2(-\frac{3}{4}+\delta)C_0\ve^\f12\l t\r^{-\frac{7}{4}+\delta}-G-\f12\wp\a^2\l t\r^{-1}<0.
\end{equation}
Consequently,
\begin{equation*}
(-y+\a)'+\f12\wp\big((-y+\a)^2-2(-y+\a)\a)\l t\r^{-1}<0.
\end{equation*}
Now let $z=-y+\a$. It is crucial that  due to (\ref{11.11.1.23}), $z>0$ for $0<t<T_*$ as long as $\ve>0$ is sufficiently small. Consequently,
\begin{equation*}
(z^{-1})'=-z^{-2}z'>\f12 \wp (1-2z^{-1}\a) \l t\r^{-1},\quad 0<t<T_*.
\end{equation*}
Hence
\begin{align*}
(z^{-1}\exp\int_0^t\wp \a \l t'\r^{-1} )'&=\big((z^{-1})'+\wp\a z^{-1}\l t\r^{-1}\big) \exp\int_0^t \wp\a \l t'\r^{-1}\\
&>\f12 \wp\l t\r^{-1}\exp\int_0^t \wp \a \l t'\r^{-1}> \f12 \wp \l t\r^{-1}.
\end{align*}
Integrating the last inequality in $t$ gives
\begin{equation*}
z^{-1}>(\exp\int_0^t \wp \a \l t'\r^{-1}dt')^{-1}(z^{-1}(0)+\f12\wp\log (\f12\l t\r)).
\end{equation*}
We will check in Proposition \ref{12.21.1.21} that $y(0)\les \ve$, which gives $z(0)>0$. Let $\bff(t)=\exp\int_0^t \wp \a \l t'\r^{-1}dt'$. Then $1\le \bff<\bff(\infty)=\ti C_1$, with $\ti C_1$ a constant slightly greater than $1$. We thus have
\begin{equation*}
z^{-1}>\ti C_1^{-1}(z^{-1}(0)+\f12\wp\log (\f12\l t\r))
\end{equation*}
i.e.
\begin{equation*}
(-y+\a)^{-1}>\ti C_1^{-1}\left((-y(0)+\a(0))^{-1}+\f12\wp \log (\f12 \l t\r)\right)
\end{equation*}
which implies on $(0, T_*)$
\begin{equation}\label{9.10.6.22}
0<-y+\a<\ti C_1\left((-y(0)+\a(0))^{-1}+\f12\wp \log (\f12 \l t\r)\right)^{-1}.
\end{equation}
Since $0<\a(0)-C\ve<-y(0)+\a(0)<M_1+1$  with $M_1>0$ depending on $\A_0$, and due to the fact that $\ve>0$ is sufficiently small, we derive 
\begin{equation*}
-y<\ti C_1\left((M_1+1)^{-1}+\f12\wp \log (\f12 \l t\r)\right)^{-1}\le\ti C_1(M_1+1)\left(1+\f12 \wp \log (\f12 \l t\r)\right)^{-1}.
\end{equation*}
With $\M_0=2\ti C_1(M_1+1)$, the above estimate  improves the right-hand side inequality in (\ref{9.23.5.23}) on $(0, T_*)$. Hence (\ref{9.23.5.23}) holds beyond $T_*$ due to the principle of continuity, which implies $T_*=t_*$. This shows that the solution is bounded in $(0, t_*)$ which contradicts to the definition of $t_*$ if $t_*<\infty$. Consequently, $t_*=\infty$. 

Thus we proved the result for the model problem, analogous to  (1) in Theorem \ref{mainthm1}. It remains to show the result analogous to (2) in Theorem \ref{mainthm1}. 
 
\subsubsection{Rarefaction at null infinity} We consider in the region $\{u_0\le u\le u_*\}$, where the global solution has been constructed from the given data satisfying (\ref{1.12.1.22})-(\ref{exist}). Along the outgoing null geodesics in this region, schematically, using the asymptotic equation of (\ref{4.10.2.19}),  we can obtain for $y\approx (t+u) \Lb \varrho$ that
\begin{equation}\label{9.10.7.22}
|(\bb y)'|\le C_3\ve^\f12 \l t\r^{-\frac{7}{4}+\delta},
\end{equation}
with $C_3$ a constant depending only on $\A_0$ and $c_*$. Let $C=C_3\int_0^\infty \l t\r^{-\frac{7}{4}+\delta} dt$. Integrating in $t$, noting that $u=r$, $\bb=c^{-1}$ at $t=0$  with $C_2=2C$ in (\ref{rarif}), we can obtain
\begin{equation*}
\bb y(t, u_1, \omega_1)<-C \ve^\f12.
\end{equation*}
 $\bb$ is positive and bounded from below due to {\bf Step 1}. By the asymptotic behavior of $y$, $\bb^{-1}\rightarrow 0$ as $y\rightarrow 0$ when $t\rightarrow \infty,$ which shows the rarefaction occurs at the infinity of null geodesic $\Upsilon_{\omega_1, u_1}(t)$. In Section \ref{brate}, we will give the explicit lower bound of $\bb$ which diverges to $\infty$ as $t\rightarrow \infty$ along the geodesic $\Upsilon_{\omega_1, u_1}(t)$. 

\subsection{Main steps of the proof}
We outline main steps of the proof.
\begin{enumerate}
\item Reduction to an asymptotic equation: We reduce (\ref{4.10.2.19}) to (\ref{9.10.1.22}) and (\ref{9.28.1.23}), with the errors $G$ and $G_2$. To control $G$ and $G_2$, it requires obtaining the full sets of weighted energy estimates. 

\item Control of weighted energies: This includes the standard energy of the angular derivatives of $\varrho$ and $\div v$ obtained by applying rotation vector-fields, and the weighted energies of tangential derivatives of $\Phi$ obtained by applying the vector-fields in $\{\Omega, S\}$. 
 In Assumption \ref{5.13.11.21+} we make bootstrap assumptions for $0<t<T_*$ on decay properties of various components of $\bp\Phi$, higher order derivatives of $\Phi$ and other geometric quantities. The decay of the key quantity $\Lb \varrho$, assumed in (\ref{6.5.1.21}), is expected to improve by analyzing the asymptotic equation of $\varrho$ after the weighted energies are controlled. All remaining bootstrap assumptions are on quantities controllable by the weighted energies using various weighted Sobolev inequalities. 
 
  The boundedness of weighted energies is achieved through establishing two types of hierarchies:
\begin{enumerate}
    \item[(a)] Hierarchies in the asymptotic behaviors of energies across different types and orders.
\item[(b)] Hierarchies in the amplitude across various 
geometric quantities and energies of different types.
\end{enumerate}
To be more specific on (b),  it is important to note that the large transversal derivatives in the initial data poses the major difficulty to the nonlinear analysis. The standard energies of $\bp\Phi$ are not small initially.  We only commute tangential vector-fields $\Omega$ and $S$ with the wave operator to derive higher order weighted energies, since these energies have smallness at the initial slice.  
 In principle, for the quantities which are possibly large either pointwisely or in $L^4(S_{0,u})$ initially, such as $\Lb \varrho$, and $\bp^2\Phi,$ we rely on transport equations to represent them in terms of the initial data;  for those having  smallness property initially we obtain their smallness for all $t<T_*$ by propagating the weighted energies, using Sobolev inequality and geometric structures. Nevertheless, commuting the vector-fields $S$ and $\Omega$ with the wave operator generates various terms without smallness properties, which are not controllable by the weighted energies, and therefore are harmful to nonlinear analysis. We adopt two strategies to cope with this issue: one is to decompose such non-small quantities into the higher order part which has smallness property, and the remaining lower order part which can be represented by $\Lb \varrho$; the other is to cancel the large term by carrying out proper normalization. 

 We will explain the main technicalities in this step shortly.    
 
\item Control of the total energy: We use the full set of weighted energies, the asymptotic behavior of $\Lb\varrho$ due to bootstrap assumption, and the bound of $\|\bp^2\Phi\|_{L^4(S_{t,u})}$ obtained in (2) to control the total energy for all $t<T_*$ by their initial data, which then extends the solution beyond $T_*$ by the classical well-posedness result. 

\item Control of $\tir \Lb \varrho$  by extending the model arguments in Section \ref{5.9.1.23} : with errors $G$ and $G_2$ controlled by using the decay properties given by the weighted energies, we improve the bootstrap assumption for $\Lb \varrho$ for $t<T_*$. Finally, we use (\ref{rarif}) to show that rarefaction occurs as $t\rightarrow \infty$ along the corresponding null geodesics.
    \end{enumerate}

\subsection{Key technicalities for controlling the weighted energies}
We have briefly given a simplified argument in Section \ref{5.9.1.23}, which determines the key behaviors of the solution.  To make the strategy work, it is crucial to control the error terms and ensure that they have sufficient decay and smallness properties provided by the initial data, particularly by the smallness assumption (\ref{9.22.1.22}).  This requires us to control the weighted energies, which dominates the whole paper.

We carry out two types of analysis for controlling the weighted energies: 

(1) Control the spacetime geometry and the deformation tensors of $S$ and $\Omega$, as well as their derivatives. In particular we establish a set of estimates of geometric comparison to compare connection coefficients with derivatives of $\Phi$, which are important for bounding the weighted energy estimates; 

(2) Run energy estimates by using the geometric estimates of deformation tensors, which  relies crucially on the nonlinear structures of the Euler equations and acoustical geometry. 

These two sets of estimates are integrated together using a bootstrap argument. We face three primary challenges: (1) the lack of smallness for closing the energy argument; (2) the lack of control over transversal derivatives of geometric quantities; (3) Loss of decay in the acoustical causal geometry for closing the top order energies.  In the sequel, we overview the techniques for solving these issues.
\subsubsection{Establish the hierarchy of amplitude} 

 In a standard small-data-global-existence result, the amplitudes of asymptotic estimates are expected to be consistent with their initial value. The difference in asymptotic behaviors of various quantities lies in the decay rate. In our case,  both the amplitudes and decay rates could vary across different quantities, due to the distinctive structures associated with each.
Therefore to establish the hierarchy of amplitude,
 we  exploit geometric structures of of various quantities with the help of the Euler equations. In this regard, we emphasize two novel elements:
 \begin{enumerate}
     \item[$\bullet$] Geometric null forms in Euler equations 
     \item[$\bullet$]  The application of the geometric decomposition of the $v$ derivatives. 
 \end{enumerate}

 Note that the terms on the right-hand side of the equations (\ref{4.10.1.19}) and (\ref{4.10.2.19}) are quadratic in terms of $\bp\Phi$, with at least one factor being tangential derivatives of $\Phi$. This is regarded as the general null condition in this paper.  In \cite{shock_demetrios, Miao_thesis}, despite working on different unknowns, the right-hand side of Euler equations therein verify null conditions in the same sense.
Nevertheless, in either formulation, Euler equations do not satisfy the typical null conditions for quasi-linear wave given in Klainerman \cite{klinvar}.  
  This can be clearly seen in the asymptotic equations (\ref{9.10.1.22}) and (\ref{9.28.1.23}). The term $\f12\wp y^2\l t\r^{-1}$ breaks the standard null condition for achieving the small-data-global-existence result. This term is hidden in the quasi-linear wave operator (see similar discussions in Speck \cite{Spck_shock_1} and Luk-Speck \cite{Jared_Luk}). In our setting the transversal derivatives are not necessarily small even initially, thus the null condition is even more important, (for instance,) for controlling the error terms $G$ and $G_2$ stated in Section \ref{5.9.1.23}.  
 
In Proposition \ref{geonul_5.23_23}, we summarize the quadratic structures, appeared both on the right-hand side of the Euler equations and in geometric quantities such as the nonlinear terms in $\bR_{44}$, into an intrinsic, and much more refined formulation than a general null condition.  This more precise geometric formulation guarantees the necessary level of smallness for our analysis. Additionally, employing the intrinsic formulation alongside the set of geometric structure equations for the components of $\bp v$, derived in Lemma \ref{dcom_s}, allows us to circumvent the technical complexities involved in comparing the acoustical geometry with the Minkowski geometry.  

To establish the hierarchy of the amplitude, it is important to distinguish the types of components of $\Phi$-derivative. If taking transversal derivative over the $v$-derivatives, different components may have completely different behaviors. For example, the angular component of $\hN\Omega v$ contains a large term $\approx\hN v_\hN$ (see Lemma \ref{2.9.3.23}), while its radial component $\hN\Omega v\c\hN$ can be represented by geometric quantities having smallness property.
 Regarding $\p v$ as the components of the second fundamental form $\Lie_\bT \delta_e$, 
we decompose $\p v$ into 
\begin{equation}\label{10.4.14.23}
 \ep:=-\p_A v^i \hN^i,\, \eta_{AB}:=-\p_A v^i \hat e_B^i=\eh_{AB}+\f12\tr\eta\delta_{AB}, \, \stc\delta:=-\hN v^i \hN^i,  
\end{equation}
where  $\hN=c^{-1}\bN$ and $\{\hat e_A, \hat e_B\}_{A,B=1,2}$ forms the orthonormal basis on $(S_{t,u}, \gac)$. 

 We further derive the geometric structure $\tr\eta=[L\Phi]$ 
in (\ref{7.04.9.19}), and the Hodge system for $\eh$ and the radial transport equations for $\ep$ and $\eh$, (\ref{5.30.3.23})-(\ref{12.16.1.23}).
We observe that the derivatives of $\varrho$ have similar behavior as the radial component of the same derivatives of $v$, hence decompose 
general one-derivatives of $\Phi$  into 
 $$[\sn\Phi], [L\Phi], [\Lb \Phi], \eh.$$ 
Using (\ref{9.22.1.22}) and the boundedness in (\ref{1.12.1.22}), we obtain in Proposition \ref{12.21.1.21}:
\begin{equation*}
|[L \Phi]|(0)\les \ve^\f12,\, |[\sn\Phi], \eh|(0)\les \ve, |[\Lb\Phi](0)|\les 1.
\end{equation*}
  By analysing the initial data, we observe the hierarchy of amplitude for various quantities after being differentiated (modulo their expected decay rates), due to the geometric structure of their derivatives. In  
  Definition  \ref{5.23.2.23}, we introduce the index of amplitude via which quantities such as $\bp \Phi$ and the important connection coefficients are classified in three sets $\al_0, \al_\f12, \al_1$. For instance,  
 \begin{equation*}
 [\Lb \Phi]\in \al_0, [L\Phi]\in \al_\f12, [\sn\Phi], \eh\in \al_1.
 \end{equation*}
 We provide in Table (\ref{5.24.1.23}) their classification after being differentiated by different derivatives. The analysis in this paper is structured to show such hierarchy of the amplitude holds for all $t$ by using both energy estimates and geometric structure equations. 

The main advantage of using the decomposition of $\p v$ in (\ref{10.4.14.23}) lies in that the estimates obtained by using the set of structure equations for $\ep$ and $\eh$ complement crucially the control from the weighted energies, important for proving the strong derivative  estimates of $\al_1$ in Table (\ref{5.24.1.23}), (see for instance Section \ref{12.6.1.23}). This decomposition allows us to give geometric decomposition for the null forms in Euler equations in Proposition \ref{geonul_5.23_23}, which are some particular quadratic forms satisfying the general null condition. Decomposing the null forms without relying on the Cartesian coordinates avoids significantly the comparison analysis between the frames in the acoustical spacetime with their Minkowski counter parts.
 
As a direct application of the geometric formulation of null form (given in  Proposition \ref{geonul_5.23_23}),  it gives sufficient smallness to important geometric quantities. For instance, $\bR_{44}$ appears in (\ref{s1}), which is the most important equation for controlling the acoustical geometry. It is important to identify (in Proposition \ref{ric44}) the key nonlinear error terms  in $\bR_{44}$  taking the null form of $\N(\Phi,\bp\Phi)$, a particular quadratic combinations of the scalar components of $\bp\Phi$ together with the good term $\eh\c \eh$. This is, for instance, used to gain sufficient smallness for $\sn\tr\chi$ and its derivatives, which are crucial for controlling both geometry and energies. For the similar reason, in Section \ref{geostru2}, we provide detailed decompositions for the Riemann curvature terms appeared in Proposition \ref{6.29.1.24}. 


\subsubsection{New decompositions of commutators between tangential vector fields and the wave operator}\label{1.5.2}  For quasi-linear wave equations without verifying null conditions, one usually does not expect to have global solution in the full exterior region $\{u\ge u_0\}$ even assuming small initial total energies (see Wang \cite{Wang_Exterior}).
To construct global exterior solutions and to propagate the smallness of the weighted energies, the most obvious challenges arise from that the Euler equations do not verify the standard null condition, and that our assumption on the general initial total energy is only bounded. Typically, for obtaining global stability  result of the trivial solutions for a nonlinear wave equation, the core analysis is to recover to the solution the linear behavior of a free wave. The method to obtain the desired decay properties is either by Klainerman-Sobolev inequalities, or by representation via characteristics.
 Even for quasi-linear wave equations that verify null conditions, proving the global stability of a trivial solution requires smallness assumptions on the initial total energies (see Klainerman \cite{klinvar, SKNull} and Christodoulou \cite{DM}). Usually the major effort is on bounding weighted energies, which will give the decay property by using Klainerman-Sobolev inequalities. The smallness assumptions are crucial for bounding the weighted energies. One could refer to Yang-Yu \cite{Yu-Yang}, Klainerman-Wang-Yang \cite{KWY} and Fang-Wang-Yang \cite{FWY} for fully or partially large data results for achieving asymptotic decays on Maxwell Klein-Gordon equations, which are semi-linear and verifying null conditions.%

As seen in the model equation (\ref{9.10.1.22}), the irrotational isentropic compressible Euler equations are quasi-linear, and do not satisfy the standard null conditions. Due to the term which breaks the null condition,  by making the initial assumptions (\ref{1.12.1.22})-(\ref{exist}), we can extend the model argument to show that the key behavior of the solution, i.e. the decay of $\Lb \varrho$, is slightly stronger than $\bp\phi$ if $\phi$ is the free wave,  provided that the tangential derivatives decay much better and have sufficient smallness. To control the tangential derivatives, we employ various Klainerman-Sobolev-type inequalities, which relies on the smallness of the weighted energies for all $t$ in the exterior region $\{u\ge u_0\}$. This necessitates energy arguments using the wave equations (\ref{4.10.1.19}) and (\ref{4.10.2.19}). 
 
Due to lacking smallness in the total energy and the lowest order weighted energies, we only commute tangential derivatives with the wave operator in Euler equations and commute $\Omega$ with the wave operator of the equation for $\Box_\bg \bT \varrho$ to gain smallness for weighted energies from the initial conditions. As discussed, the right-hand side of the equations (\ref{4.10.1.19}) and (\ref{4.10.2.19}) verify the null conditions. We need to check whether the geometric quadratic null forms are preserved in the analysis of higher order energies, particularly when commuting the tangential derivatives $S=\tir L$ and $\Omega$ with the wave operators $\Box_\bg$.

Recall the commutator formula 
\begin{equation}\label{5.10.2.23}
[\Box_\bg, X]f={}\rp{X}\pi^{\a\b} \bd^2_{\a\b} f+(\bd^\a {}\rp{X}\pi_\a^\la-\f12 \bd^\la \Tr{}\rp{X}\pi)\p_\la f, \quad X=S\, \mbox {or}\, \Omega,
\end{equation}
where the deformation tensor of $X$ is defined by ${}\rp{X}\pi_{\mu\nu}:=\l \bd_\mu X, \p_\nu\r+\l \bd_\nu X, \p_\mu\r$.
The standard method of manipulating commutators easily breaks smallness, resulting in several sets of large terms that are challenging to be grouped back together to yield the desired structure for smallness. 
Note the standard scaling vector field $t\partial_t + r\partial_r$ is a conformal Killing vector in Minkowski space, which is transversal to the Minkowskian null cone. We choose $S$ to be tangential to the null cone, nevertheless, our $S$ is not conformal Killing even in Minkowski space. Therefore, we do not expect to have smallness from the traceless part of ${}\rp{S}\pi$. 

In our setting, the components of deformation tensors ${}\rp{S}\pi$ have non-small terms $\tir\tr\chi$ and $\tir k_{\bN\bN}$, with the latter vanishing in Minkowski space. To control the weighted energies of $f=S^n\Phi$ with $n=1,2,3$,  in commutators there are the completely large terms ${}\rp{S}\pi^{L\Lb} \bd^2_{L\Lb} S^{n-1}\Phi$ and ${}\rp{S}\pi^{A B} \bd^2_{A B} S^{n-1}\Phi$, for which we do not expect to have smallness in amplitude, since none of the factors are.
 
Adopting the decomposition for ${}\rp{S}\pi$ into trace and traceless parts still generates completely large terms. To solve this issue, we give a new decomposition of the commutators, in Proposition \ref{11.12.2.22} in Section \ref{error_manu}. It is based on decompositions in both ${}\rp{X}\pi$ and the paired derivatives of $\Phi$. 

 Different from  the standard trace-traceless decomposition, we decompose $${}\rp{X}\pi={}\rp{X}\pih-\f12{}\rp{X}\pi_{L\Lb}\bg$$ so that the factor ${}\rp{X}\pih_{L\Lb}=0$ is paired with $\bd^2_{L\Lb} f$, the factor paired with $\Lb f$ has sufficient smallness, with the potentially large term ${}\rp{X}\pi_{L\Lb}$ paired with $\Box_\bg f$. This allows us to take advantage of the geometric null forms inductively. It works particularly well when commuting $S$ with $\Box_\bg$, and higher-order commutations for $S$ can be done similarly. 

 Remarkably, with careful analysis, for decomposing $[\Box_\bg, \Omega]$, the above new decomposition of the commutators introduced originally to treat $[S, \Box_\bg]$ does not bring unwanted terms such as the higher order derivatives of $\sn\log \bb$ to the calculation,  see in Proposition \ref{error_terms} for
 the cancellation in achieving $\bJ[\Omega]_B$ (the term paired to $\sn f$ in the commutator in (\ref{5.10.2.23})). This helps crucially for closing the top order energies. Therefore, we adopt such decomposition for commuting both $S$ and $\Omega$ with the wave operator. 
 
 To control the weighted energies of $f=\Omega^n\Phi$ with $n=1,2,3$, we consider the terms generated due to the commutator $[\Omega, \Box_\bg]$. By using the new decomposition, for instance, we have cancelled terms involving $\bd^2_{L\Lb}\Omega^{n-1}\Phi$ in the commutator and  the decomposition of $\bd^2_{AB}\Omega^{n-1}\Phi$ allows us to pair  small factors to $\sn^2\Omega^{n-1} \Phi$ and $\Lb \Omega^{n-1}\Phi$  since $\pioh$ has small amplitude.  If $\Phi=v$, $\Lb \Omega^{n-1}\Phi$ is not a small term, and in the terms (for instance) $\sum_{l=0}^{n-1}\sn_\Omega^{l}{}\pioh_{AL}\bd^2_{\Lb A}\Omega^{n-1-l}\Phi$,  the second factors are mixed derivatives of $\Phi$,  containing non-small components. All of theses terms are not bounded directly by the weighted energies. We decompose such mixed derivatives of $\Phi$ into two parts: the good part, small, of higher order, controllable by weighted energies, and the large part representable by $\Lb \varrho$. (See Lemma \ref{9.18.5.23}.)  We then close the energy estimates by using the estimates of $\pioh$, $\bJ[\Omega]$ and the energy inequalities with the help of  the decay of $\Lb \varrho$ that is controlled separately. 

 For the weighted energies of  derivatives  $X^n\Phi$ with $X\in\{ S, \Omega\}$ generally, and the energy of $\Omega^n\bT \varrho$, we adopt the same strategy. 

\subsubsection{Energy hierarchy}
The energy hierarchy for the weighted energies (including the standard energy of $\Omega^{n-a}\bT^a\varrho, a=0,1, n=1,2,3$) is formulated for controlling $\tr\chi-\frac{2}{\tir}$ and $\ze$, (see (\ref{ricc_def}) for the definitions of $\tr\chi$ and $\ze$), which are crucial geometric quantities to close the energy arguments. Bounding up to the third order weighted energies is necessitated by obtaining the pointwise estimate of $\sD\varrho$, crucial for achieving the desired decay for $\tr\chi-\frac{2}{\tir}$.  
 
For the important top order quantities $\Omega^n\tr\chi$ and $\sn_\Omega^n\ze$, with $n\le 3$, it is expected to incur the loss of derivatives,  due to the presence of nontrivial Ricci terms (effectively, higher-order derivatives) in the transport equation (\ref{s1}) and the equation for $L\mu$, which is normalized into (\ref{8.31.4.19}). As mentioned in Section \ref{1.5.2}, due to cancellation, the highest order derivative of $\ze$ does not appear in the most crucial way in the commutator $[\Box_\bg, \Omega^3]$ (see the derivation of $\bJ[\Omega]_A$ in Proposition \ref{error_terms}). Hence, the primary challenge from acoustical geometry lies in the loss of derivatives in bounding $\tr\chi$ and getting the additional factor of $\bb$, for controlling the highest order energies.

We observe that, with $n=3$, the hardest error terms for bounding the top order weighted energies are
\begin{align*}
\E_1&=\int_0^t \int_\Sigma \Omega^n \tr\chi\c \Lb \varrho \c \bT \Omega^n\varrho d\mu_g dt'\\
 \E_2&=\int_0^t \int_{\Sigma}\Omega^n \tr\chi \Lb \Phi (L+\f12\tr\chi)\Omega^n\Phi \tir^m d\mu_g dt'
\end{align*}
where $\E_1$ arises from controlling $E[\Omega^3\varrho](t)$ \begin{footnote}{See Definition \ref{3.24.2.24} for the definitions of various types of energies.}\end{footnote}, the standard energy of $\Omega^3 \varrho$, and $\E_2$ is for bounding $W_m[\Omega^3\Phi](t)$, the weighted energy of $\Omega^3\Phi$, with the weight depending on the choice of $m$. To solve the issue of loss of derivatives, we decompose  $\tr\chi=\sX-\Xi_4$, with $\Xi_4=\Lb\varrho+[L\Phi]$ symbolically (defined in (\ref{ricc6.7.2})). The terms contributed by the part of $\Xi_4$ in $\E_1$ and $\E_2$ are sharp terms in the energy argument, provided that the commutation $[\Lb, \Omega^3]\varrho$ can be well controlled.

For $\Omega^3\sX$, thanks to the renormalized equation (\ref{12.5.1.21}), we derive  
\begin{align*}
L\Big(\tir^3\sn_\Omega^2(\bb^{-1}\sn\sX)\Big)=\tir^3 \sn_\Omega^2\Big(\bb^{-1}\big(\sn \Xi_4(\Xi_4+\f12 \sX)\big)\Big)+\mbox{ better terms}.
\end{align*}
 Integrating this transport equations along null cones leads to 
\begin{equation}\label{3.24.1.24}
\|\tir^2\sn_\Omega^2\big(\bb^{-1}\sn(\tr\chi+\Xi_4)\big)\|_{L^2_u L_\omega^2}\les\l t\r^{-1}(\ve+ \int_0^t \|\tir^2\sn_\Omega^2(\bb^{-1}\sn\Xi_4)\|_{L^2_u L_\omega^2}+\cdots).
\end{equation}
Due to $\Xi_4=\Lb\varrho+[L\Phi]$ symbolically, with  $\bb^{-1}[\Lb,\Omega^3]\varrho$ treated properly, the leading term on the right-hand side is roughly bounded by 
$$
\sup_{t'\le t}\|\bb^{-\frac{3}{2}}\bT\Omega^3\varrho(t)\|_{L^2_\Sigma} \les \sup_{t'\le t}E[\Omega^3\varrho]^\f12(t).
$$
The $\E_1$ term in principle could be treated by using Gronwall's inequality and the decay for $\Lb \varrho$ in (\ref{9.10.6.22}) (in the vastly simplified case). Nevertheless, there is a crucial additional growth in the factor $\bb$ due to the presence of $\bb^{-1}$ paired with $\sn\sX$ in the transport equation. Since $\bb$ could be of $\log \l t\r$ growth, such additional growth of $\bb$ does not allow us to treat $\E_1$ directly by Gronwall's inequality for bounding the standard energy $E[\Omega^3\varrho](t)$. On the other hand, due to the expected bound for the right-hand side of (\ref{3.24.1.24}), if $n=3$, we choose $m=1$ in $\E_2$. While for lower order weighted energies, i.e, $n=1,2$, we set $m=2$. Again there is also certain growth in $\bb$ for bounding $\E_2$ due to the same reason if $n=3$. Such growth fails the energy estimate as well. 

 In Christodoulou \cite{shock_demetrios}, Christodoulou-Miao \cite{Miao_thesis} and Miao-Yu \cite{Pin-Shuang}, as well as in various other studies on shock formation, the phenomenon of additional growth arising from the appearance of $\bb^{-1}$ in the top-order energy control is inevitable. 
 The method of treating this issue, often referred to as the descent scheme, requires controlling very high-order energies with growth to gain additional decay to bound lower order energies. 

  Since we only have $H^4$ data to work with,  we treat the  growth due to the additional $\bb$ factor differently. 
    Note that schematically $$\Omega^3 \sX\Lb\Phi\approx \tir\sn_\Omega^2(\bb^{-1}\sn\sX)\c\bb\Lb\varrho+\mbox{better terms}.
    $$ Also  taking advantage of the structure mentioned in Section \ref{3.24.4.24} that 
 \begin{equation}\label{3.25.1.24}
 L(\bb\tir\Lb \varrho)=\mbox{ small good terms} 
 \end{equation}
 we have an important observation
   \begin{equation}\label{12.23.1.24}
   L(\tir^4\sn_\Omega^2(\bb^{-1}\sn\sX)\c\bb\Lb\varrho)\approx \tir^2 \Omega^3\Xi_4 \Lb \varrho+\mbox{better terms}.   
   \end{equation}
   We choose a positive constant $
   \M$  depending on $\M_0$ in (\ref{9.23.5.23}), pair  a $(\log (\tir+3))^{-\M}$ and a weight of $(\log (\tir+3))^{-\M-1}$  to the energy densities of $E[\Omega^3\varrho](t)$ and $W_1[\Omega^3\Phi](t)$ respectively.
     Using (\ref{9.23.5.23}) and (\ref{12.23.1.24}), 
 we avoid the loss of the weight of $\bb$ by carrying out bilinear estimate for the treating the error $\E_1$ and $\E_2$ (see Lemma \ref{3.8.3.24}).  
In this way the original top order energies $E[\Omega^3\varrho](t)$ and $W_1[\Omega^3\Phi](t)$ are controlled with the additional growth of $(\log \l t\r)^\M$ or $(\log \l t\r)^{\M+1}$, which suffices to close the bootstrap argument. 

Moreover, to bound the right-hand side of (\ref{12.23.1.24}), we need to bound $\|\Omega^3\Lb\varrho \log (\tir+3)^{-\frac{\M}{2}}\|^2_{L^2_\Sigma}$ by $E_{-\M}[\Omega^3\varrho](t)$ and the $L^2_u L_\omega^2$-type norm for $\Omega^3(\bb^{\f12})$  by an integral of $E_{-\M}[\Omega^3\varrho]^\f12(t)$.  In Lemma \ref{3.24.3.24}, we establish a series of delicate estimates to minimize the loss of decay due to commutation $[\Lb, \Omega]$, which requires better control of $\sn\log \bb$ and its angular derivatives than the bootstrap assumptions for them. 
  Using (\ref{3.25.1.24}), we gain  better decay in Lemma \ref{5.13.11.21} (5) for the lower derivative bounds on $\bb\Lb \varrho$ and  $\sn\bb$  over the bootstrap assumptions. However, even a lowest order estimate such as $\|\Omega\bb\|_{L_\omega^4}\les \log \l t\r\Delta_0$ is not bounded in terms of $t$, and $\Omega\log\bb$ has even stronger $t$-growth in $L^\infty_\omega$. Hence, with each commutation of $[\Lb, \Omega]$, a certain loss is inevitable. To counteract this growth stemming from commutation, it is important to enhance the lower-order energy bounds before bounding the top order energies. In addition to employing this strategy, we have a crucial observation that (\ref{9.10.6.22}) (due to (\ref{rarefied})) provides the correct signs for the first terms (i.e., the highest-order term) on the right-hand sides in the commutation formulas in Lemma \ref{3.25.2.24} and in Lemma \ref{6.29.2.24}. This enables us to control the commutator without incurring significant loss of decay rate. 
These two strategies collectively allow us to circumvent loss of decay when bounding the right-hand side of (\ref{12.23.1.24}).

In summary, due to the loss of derivatives for controlling $\Omega^n\tr\chi$, we can only improve the first order energy estimates over the bootstrap assumption directly. For the second order, we are not able to directly improve the growth rate for the weighted energies $W_2[\Omega^2\Phi](t)$. We improve its amplitude allowing an additional $\log \l t\r$ growth compared with the bootstrap assumption. Using this bound, we improve the remaining second order energies, and the majority of the decay properties assumed in Assumption \ref{5.13.11.21+} up to the second order. For the top order energies, we reduce the weight in $\tir$ in the weighted energies of $X^3\Phi$, with $X\in\{\Omega, S\}$, due to the loss of derivatives for controlling $\Omega^3\tr\chi$. Using the obtained lower order energy bounds,  we can reduce the $\l t\r$-growth in error estimates for controlling the top order energies. Once the top order energies are obtained, we use it to improve the second order weighted energies of $\Omega^2\Phi$, and several related decay estimates. Thus we complete the weighted energy estimates (including the standard energy estimates of $\Omega^3\varrho$ and $\Omega^2 \bT\varrho$).

Finally, to bound the total energy for all $t$, we rely on the $L^4_{S_{t,u}}$ bound of $\bp^2\Phi$ (see (6) in Lemma \ref{5.13.11.21}), which, together with other decay properties of $\bp\Phi$,  works as the criterion to extend the solution by continuity. Note the key part of the quantity $\bp^2\Phi$ is $\Lb^2\varrho$. To bound this quantity by Sobolev inequalities, we need the bound $\|(\bN, \tir\sn)\Lb^2\varrho\|_{L^2_\Sigma}$ or $\|(S,\tir\sn)\Lb^2\varrho\|_{L^2(\H_u)}$, which relies on the uncontrolled part of the total energy.  
  However the total energy is not obtained a priori, we instead directly derive the $L^4_{S_{t,u}}$ estimate of $\Lb^2\varrho$ from initial data  by integrating along null cone $\H_u$  using the transport equation given in Lemma \ref{12.21.3.23}. 
  
\subsection{Organization}  
\begin{itemize}
\item In Section \ref{geosetup}, we give the geometric basics: in Lemma \ref{dcom_s} we derive a set of crucial structures on $k=-\f12\Lie_\bT g$ and the structure equations for components of $\p v$; we give the notion of the geometric null condition for the Euler equations in Proposition \ref{geonul_5.23_23}, the crucial decomposition of the Ricci component $\bR_{44}$ in Proposition \ref{ric44} and introduce the index of amplitude in Definition \ref{5.23.2.23} and give in Table (\ref{5.24.1.23}) the hierarchy on it. 
\item In  Section \ref{causal_geo}, we give the null structure equations, as well as derive the normalized ones. We also derive the transport equation for $\Lb \sX$ in Lemma  \ref{10.8.1.23}, which is crucial for giving the control on $\|\sn^3 \ze\|_{L_u^2 L_\omega^2}$ in Proposition \ref{gmtrc_high_od}; and derive the transport equations for $\bb \Lb (\tir \Lb\varrho)$ in Lemma \ref{12.21.3.23}, which is crucial for controlling $\Lb^2\varrho$ in Proposition \ref{6.27.1.24}. 
 \item In Section \ref{rotation}, we recall the definition from \cite{shock_demetrios}  of the rotation vector-fields and give basic properties. Since our acoustical metric is conformal to the metric in \cite{shock_demetrios}, we enclose proofs of the properties in Section \ref{append} for completeness. \item In Section \ref{BA_decay}, we set up the bootstrap assumptions (Assumption \ref{5.13.11.21+}) and derive an important set of basic decay properties in Lemma \ref{5.13.11.21}. For setting up Assumption \ref{5.13.11.21+}, we  give in Proposition \ref{12.21.1.21} the estimates for various quantities at $t=0$  with the detailed proof given in Section \ref{9.25.2.22}. In Section \ref{framecmp} we give basic comparison results of frames and vector fields.  In Section \ref{trans_control}, we give the control on various transversal derivatives of solutions and of connection coefficients. We also give a set of important estimates on $\mho=\Lb \log \tir-\f12 \tr\chib$ and its derivatives, which gives basic bounds on $\Lb\tir$  and is used in controlling the commutators between vectorfields and the wave operator in later sections.
       \item In Section \ref{mul_1}, using the basic decay estimates obtained in Section \ref{BA_decay}, we derive the fundamental energy inequalities in Proposition \ref{10.10.3.22} for the standard energy and Proposition \ref{MA2} for the weighted energy. 
       \item We give a set of elliptic estimates in Section \ref{Elliptic} to treat the Hodge systems on spheres $S_{t,u}$. 
        We then give the lower order estimates of connection coefficients in Section \ref{low_ricci}: in Proposition \ref{11.4.1.22} for comparing derivatives of $\tr\chi, \chih, \zeta$ with derivatives of $\Phi$, and in Proposition \ref{8.12.1.23} and Proposition \ref{10.4.4.23} for derivative estimates of $\eh$. They are obtained by using the transport  null structure equations and the Hodge systems on spheres $S_{t,u}$ given in Section \ref{geosetup} and Section \ref{causal_geo}. 
           Finally in Proposition \ref{7.15.5.22} we summarize the estimates we obtain by using Assumption \ref{5.13.11.21+}, the aforementioned estimates, and Lemma \ref{5.13.11.21}.  
       \item In Section \ref{error_manu}, we give a careful manipulation on the commutator $[\Box_\bg, X]$ for $X=\Omega, S$. In Section \ref{Jformulas}, we derive the formulas of  ${}\rp{X}\pih$ and $\bJ[X]$.  Under the bootstrap assumptions (\ref{3.12.1.21})-(\ref{6.5.1.21}), we derive a set of estimates on ${}\rp{X}\pih$ and their derivatives in Proposition \ref{5.24.1.21} with the help of results in Section \ref{BA_decay}-\ref{low_ricci}. Moreover, in Section \ref{comp_sec}, as an important step for bounding the weighted energy, we give the comparison estimates between the derivatives of $\Phi$ appearing in the commutator $X^a[\Box_\bg, X]X^{n-a-1}\Phi$, with $0\le a\le n-1, n=1,2$ with the weighted energies. 
       \item In Section \ref{low_energy}, we bound the weighted energies for $X\Phi$ with $X=S, \Omega$ and provide decay estimates as its consequence. The main result is presented in Proposition \ref{1steng}.
       \item In Section \ref{2nd} and Section \ref{top_eng}, we obtain the second and the third order estimates of the weighted energies. The estimates of the second order energies and the related decay estimates are given in Proposition \ref{8.29.8.21} and Proposition \ref{9.8.6.22}. The weighted energy estimates of $\Omega^2\Phi$ we obtain in this section  have not improved the corresponding bootstrap assumptions in terms of decay rate. We gain smallness in amplitude for $WFIL_2[\Omega^2\Phi]$, while slightly losing the decay rate compared with the bootstrap assumptions. With the help of such gain, all other desired second order energies are improved over their bootstrap assumptions. We provide the top order energy estimates, the improved second order estimates and the top order decay estimates in  Proposition \ref{10.30.4.21}-Proposition \ref{imp_decay}. This allows us to improve Assumptions \ref{5.13.11.21+} except for (\ref{6.5.1.21}). 
             
     \item In Section \ref{10.24.2.23}, we use the decay properties  on the fully tangential derivatives and mixed derivatives derived so far, and the key control of $\Lb^2\varrho$ in Proposition \ref{6.27.1.24} to obtain the complete control on total energies for all $t<T_*$ in Proposition \ref{10.13.5.23}, thus extend the solution beyond $T_*$ by the standard local well-posedness result. 
     \item In Section \ref{10.24.1.23}, we adapt the model argument to prove (\ref{6.5.1.21}) in Assumption \ref{5.13.11.21+}, and confirm the formation of rarefaction under the additional assumption (\ref{rarif}).  Hence the proof of our main result, Theorem \ref{mainthm1}, is complete.
     \item In Section \ref{9.25.2.22}, we provide the proof of Proposition \ref{12.21.1.21}. We rely on the formulas and basic comparisons established in other sections. Nevertheless, the analysis involved in this section is independent of time evolution.  
         \item In Section \ref{geostru2}, we give the structures of Christoffel symbols and crucial components of Riemann curvature.
         \item In Section \ref{geocal}, we give basic calculations for treating $\pioh_b$ and its derivatives,  with the help of the result in Lemma \ref{5.13.11.21} under the assumptions of (\ref{3.12.1.21})-(\ref{1.25.1.22}). 
         \item In Section \ref{append}, we provide the proofs of Proposition \ref{2.19.4.22} and Proposition \ref{3.22.6.21} for completeness.                  
\end{itemize}

\section{Geometric set-up and the geometric null conditions}\label{geosetup}
We start with recalling the basic geometric set-up defined by using the null tetrad $\{L, \Lb, e_A, A=1,2\}$,\begin{footnote}{Alternatively, $e_4=L, e_3=\Lb$. }\end{footnote} which appeared in \cite{CK, KRduke, KR1, Wangrough, rough_fluid}.
Let $\{e_A, e_B\}$ with $A,B=1,2$  be the orthonormal basis of the tangent bundle on $S_{t,u}$.
 The null second fundamental forms $\chi$ and $\chib$,
the torsion $\zeta$, and the Ricci coefficient $\zb$ of the foliation $S_{t,u}$ are defined by
\begin{equation}\label{ricc_def}
\begin{split}
\chi_{AB}=\bg(\bd_A e_4, e_B), &\qquad \chib_{AB}=\bg (\bd_A e_3, e_B),\\
\zeta_A=\f12 \bg(\bd_3 e_4, e_A), &\qquad \zb_A=\f12 \bg (\bd_4 e_3, e_A).
\end{split}
\end{equation}
 We denote by $\tr\chi$ and $\chih$ the trace and traceless part of $\chi$ taken by the metric $\ga$, and apply the same convention to $\chib$.
The second fundamental form of $S_{t,u}$ in $\Sigma_t$  is given by
\begin{equation}\label{7.15.7.19}
\theta(X, Y)=\l \sn_X \bN, Y\r
\end{equation}
for any vector fields $X, Y$ tangent to $S_{t,u}$. The trace of $\theta$ is defined by $\tr\theta=\ga^{AB} \theta_{AB}$, and the traceless part of $\theta$ is denoted by $\hat \theta$.
Let
 $\hN=c^{-1}\bN$ be the unit normal to $S_{t,u}\subset\Sigma_t$ with respect to the Euclidean metric $\delta_e$. Note that $\gac$ is induced metric of $\delta_e^{ij}$ on $S_{t,u}.$ Using $\hN$, with respect to $\gac$, we can define on the tangent space of  $S_{t,u}$ the second fundamental form $\thetac$,  $\tr\thetac$ and $\hat\thetac$ in the same way as their counter parts $\theta, \tr\theta$ and $\hat \theta$ relative to the metric $\ga$.
  
\begin{proposition}\label{6.7con}
\begin{align}
\bd_A e_4=\chi_{AB} e_B-k_{A\bN} e_4 &\qquad  \bd_A e_3=\chib_{AB}e_B+k_{A\bN} e_3 \label{7.21.3.19}\\
\bd_4 e_4=-k_{\bN\bN} e_4   &\qquad  \bd_4 e_3=2\zb_A e_A+  k_{\bN\bN} e_3\label{6.29.5.19} \\
\bd_3 e_4=2\zeta_A e_A+k_{\bN\bN} e_4 & \qquad \bd_4 e_A=\sn_L e_A+\zb_A e_4 \label{6.29.6.19}\\
\bd_B e_A=\sn_B e_A+\f12 \chi_{AB} e_3+\f12 \chib_{AB} e_4& \qquad \bd_3 e_3 = (-2\zeta_A+2k_{\bN A}) e_A -k_{\bN\bN} e_3\label{6.29.7.19}\\
\chi_{AB}=\theta_{AB}-k_{AB}, \quad\zb^A=-k^A_\bN &\qquad  \zeta^A=\sn^A\log \bb+k^A_\bN\label{3.19.1}\\
 \bd_\Lb e_B=\sn_\Lb e_B -\zeta_B \Lb+(\zeta_B-k_{\bN B})L.&\label{5.01.3.21}
\end{align}
\end{proposition}
As a direct consequence of (\ref{7.21.3.19})-(\ref{6.29.7.19}), the following decomposition holds under the null tetrad.
\begin{corollary}\label{decom_wave}
Let $h=\f12 \tr\chi$ and $\hb=\f12 \tr\chib$. For a scalar function $f$, there hold
\begin{align}
\Box_\bg f&=\sD f-\Lb L f-(\hb-k_{\bN\bN}) L f-h\Lb f+2\zeta^A \sn_A f\label{6.30.1.19}\\
\Box_\bg f&=\sD f-L \Lb f-(h-k_{\bN\bN})\Lb f-\hb L f+2\zb^A \sn_A f\label{6.30.2.19}\\
L(v_t^\f12\bb \Lb f)&=\bb v_t^\f12(\sD f-\hb L f+2\zb^A \sn_A f-\Box_\bg f).\label{1.14.4.22}
\end{align}
\end{corollary}
Recall from (\ref{metric}) that the components of $\bg$ and $\bg^{-1}$ are
\begin{align*}
&\bg_{00}=-1+c^{-2}|v|^2, \quad \bg_{0i}=-c^{-2} v_i, \quad \bg_{ij}=c^{-2}\delta_{ij}\\
&\bg^{00}=-1, \quad \bg^{0i}=-v^i \quad \bg^{ij}=c^2 \delta^{ij}-v^i v^j.
\end{align*}

We can also write the second fundamental form $k(X,Y):=-\l \bd_X \bT, Y\r$ for $\Sigma$-tangent \begin{footnote}
{We call a vector field $F$  $\Sigma$-tangent if $\l F, \bT\r=0$. We can directly check  that $F^0=0$ iff $\l F, \bT\r=0$.
In general any tensor field $F$ is called $\Sigma$-tangent, if the contractions of $F$ by $\bT$ using the acoustical metric $\bg$ all vanish.}
\end{footnote} vector field $X, Y$. For geometric structure of various components of $k$ relative to the orthonormal frame $\{\bN, e_A, A=1,2\}$ in $\T \Sigma$, we recall the following geometric identities from \cite[Prop 2.1]{rough_fluid}.
\begin{lemma}[Crucial decomposition for second fundamental forms]\label{dcom_s}
Let $\wp=1+(\log c)'(\varrho)$ and  $\Xi_\mu=\Ga^\eta_{\a\b}(\bg) \bg^{\a\b}\bg_{\eta \mu}$, where $\Ga(\bg)$ is the Christoffel symbol of $\bg$.
Let $\Xi$ be the 1-form
\begin{equation}\label{ricc6.7.2}
\Xi_\ga=(\Ga_{\a\b}^\eta-{\hat \Ga}_{\a\b}^\eta)\bg^{\a\b}\bg_{\ga\eta},
\end{equation}
with $\hat\Ga $ being the Christoffel symbol of a smooth reference metric  $\hat\bg$, which is chosen to be the Minkowski metric $\bm$.
 We denote by $\{\hat e_A\}$ the orthonormal frames relative to the restricted metric $\gac$ on $\T S_{t,u}$.\begin{footnote}{We may use $\hat A, \hat B=1,2$ as short-hand notations for $\hat e_A, \hat e_B$.}\end{footnote}
Relative to this frame, let $\eta_{AB}=\stc{k}_{AB}$ and $\ep_A=\stc{k}_{A\hat\bN}$ and 
 we decompose $\eta$ as 
\begin{equation*}
\tr \eta=\delta^{CD}\eta_{CD}, \quad \eta_{AB}=\eh_{AB}+\f12 \gac_{AB}\tr \eta,
\end{equation*}
where $A,B,C,D=1,2$.

There hold
\begin{align}
k_{ij}&=c^{-2}(\bT(\log c)\delta_{ij}-\p_i v_j) \label{k1}\\
\Tr k&=3 \bT \log c-\div v\label{1.26.7.22}\\
\Xi_\mu\bT^\mu&=\Tr k; \, \Xi_j=\p_j(\log c-\varrho)\label{6.14.1.19}\\
\displaybreak[0]
2\bT\varrho&=\wp^{-1}\{\Xi_L-L (\log c-\varrho)\} \label{7.04.7.19} \\
\delta:&=k_{\bN\bN}=\f12 (\Xi_L-L(\log c+\varrho)-2L v_\bN)\label{7.04.8.19}\\
&=\wp \bT \varrho-L \varrho-L v_\bN.\label{3.22.1.21}\\
\displaybreak[0]
-k_{A\bN}&=e_A(v^i) \bg_{ij}\bN^j=\bN v^i \bg_{ij} e_A^j\label{7.04.11.19}\\
\bN v^i \bN^j \bg_{ij}&= -\bT \varrho-e_A(v^i) e_A^j g_{ij};\quad L v^i \bN^j \bg_{ij}+L\varrho=\tr\eta \label{7.04.9.19}\\
\displaybreak[0]
\p_i \stc{k}_{ij}&=- \p_j \div v\label{5.24.2.19}\\
c\tr\thetac&=\tr\theta+2\bN \log c=\tr\chi+2L \log c+\tr\eta\label{5.17.1.21}\\
\tr\chi&=c\tr\sta{\circ}\theta-2L \log c+\f12 (\p_A v_B+\p_B v_A) \ga^{AB}\label{1.6.1.21}
\end{align}
\begin{align}
&\snc_{\hat \bN}\ep_{\hat B}+\thetac_{\hat B\hat A}\ep_{\hat A}=\snc_{\hat B}(\stc{k}_{\hat\bN\hat \bN})+\snc_{\hat A}\log(\bb c)k_{\hN\hN}-\snc_{\hat A}\log (\bb c) \eta_{\hat B\hat A}\label{5.30.3.23}\\
&\snc_{\hat A} (\eh_{\hat A\hat B}-\f12 \tr\eta\gac_{\hat A\hat B})=-\tr\thetac \ep_{\hat B}\label{6.1.1.23}\\
&\snc_{\hat\bN}\eh+\f12\tr\thetac \eh=\snc\hot \ep-\hat\thetac_{\hat A\hat C}\hot\eh_{\hat B\hat C}+\snc\log (\bb c)\hot \ep+(\stc k_{\hN\hN}-\f12\tr\eta)\hat\thetac_{\hat A\hat B}\label{12.16.1.23}
\end{align}
where $\hat\thetac_{\hat A\hat C}\hot\eh_{\hat B\hat C}=\hat \theta_{\hat A\hat C}\hat \eta_{\hat C \hat B}-\f12\gac_{\hat A\hat B}\hat\thetac\c \eh$.  
\begin{align}\label{2.18.4.24}
\nab_g^i k_{ij}&=-2\p^i \log c (\bT\log c \delta_{ij}-\f12 (\p_i v_j +\p_j v_i))+\p_j \bT \log c\nn\\
&-\f12 (\Delta_e v_j+\p_j \div v)-c^2\p_l \log c k_{lj}+\p_j \log c \Tr k
\end{align}
\begin{align}
&\overline{\bb^\a(\hN v_{\hN}+\tr\thetac v_\hN+\snc_{\hat A}\log (c\bb) v_{\hat A})}=-\overline{\bb^\a(\bT \varrho+\snc_{\hat A} (v^{\hat A}))}\label{9.29.5.23}\\
&-\overline{\bb^\a\tr\thetac v_{\hN}}=\overline{\bb^\a(\tr\eta-\a\snc_{\hat A}\log \bb v^{\hat A})}, \,\a\ge 0,\label{9.29.6.23}
\end{align}
where for $S_{t,u}$ tangent tensor-field $F$,  $\snc_{\hN} F=\hN^i (\nab_e)_i F$ with $\nab_e=\p$ representing the connection of $\delta_e$ on $\mathbb R^3$ and any type of derivative for $v^i$ is understood as the derivative of the  $v^i$-component. 
\end{lemma}
\begin{remark}
Note $\tr\eta=[L\Phi]$ due to (\ref{7.04.9.19}). A rough calculation implies $\sn_\Lb \eta$ is not necessarily small initially.  We gain various estimates with smallness for $\eh$, $\sn_\Lb \eh$ and their derivatives, by using the Hodge system (\ref{6.1.1.23}) and (\ref{12.16.1.23})  in Section \ref{12.6.1.23}.
\end{remark}
\begin{proof}
(\ref{k1})-(\ref{7.04.8.19}) can be found or derived by using \cite[Proposition 2.1, Proposition 7.7]{rough_fluid}.
(\ref{3.22.1.21}) follows by using (\ref{7.04.7.19}) and (\ref{7.04.8.19}).
 (\ref{7.04.11.19}) is from the definition of $k_{ij}$ and $\p_i v_j=\p_j v_i$. The two identities in (\ref{7.04.9.19}) are obtained by using (\ref{4.23.1.19}). In particular, for the second one, we calculate
\begin{align*}
-\gac^{\hat A\hat B}\p_{\hat A} v_{\hat B}&=\hN v^i \hN^i+\bT\varrho=c^{-2}\bN v^i \bN^i+\bT\varrho\\
&=c^{-2}L v^i \bN^i-c^{-2}\bT v^i \bN^i+\bT \varrho=c^{-2}L v^i \bN^i+L\varrho.
\end{align*}

(\ref{5.24.2.19}) is obtained by $\p_i v_j =\p_j v_i$.  (\ref{5.17.1.21}) can be derived by using (\ref{k1}). In (\ref{1.6.1.21}), we write (\ref{5.17.1.21}) into a slightly different form. 

 We decompose the left-hand side of (\ref{5.24.2.19}) using the radial frame $\hat\bN, \hat{e}_A$. Next we show
\begin{align}
\snc_{\hat A}\eta_{\hat A\hat B}&+\frac{3}{2}\tr\thetac \ep_{\hat B}+\hat{\thetac}_{\hat A\hat B}\c \ep_{\hat A}+\sn_{\hat\bN}\ep_{\hat B}+\snc_{\hat A} \log(\bb c) \eta_{\hat A\hat B}\nn\\
&-\stc{k}_{\hat\bN\hat\bN} \snc_{\hat B}\log (\bb c)=-\p_{\hat B} \div v\label{5.30.2.23}
\end{align}
Note
$\l \p_{\hat \bN}\hat\bN, \hat e_A\r_e=-\hat e_A\log (\bb c).$
Recall from \cite[Corollary 3.2.3.1]{CK} that
\begin{equation}\label{5.30.1.23_com}
[\snc_{\hat\bN}, \snc]f=-\thetac\c\snc f+\snc\log (\bb c) \hat\bN f.
\end{equation}
To see (\ref{5.30.2.23}), with $\Pic^{ij}=\delta^{ij}-\hN^i \hN^j$, we write the left-hand side of (\ref{5.24.2.19}) by the radial decomposition
\begin{align*}
-\p_i \stc{k}_{ij}{\hat e}^j_B=&
\snc_{\hat A}(\snc_{\hat A} v^j) \Pic_{jl}{\hat e}_B^l+\tr\thetac \snc_{\hat \bN} v^j \Pic_{jl}{\hat e}_B^l+\hat \bN(\hat \bN v^j) \Pic_{jl} {\hat e}_B^l \\
&-\snc_{\hat \bN}{\hat \bN}^m \p_m v^j \Pic_{jl} {\hat e}_B^l\\
\displaybreak[0]
&=\snc_{\hat A}(\snc_{\hat A} v^i \Pic_i^j+\snc_{\hat A} v^i \hat \bN^i \hat\bN^j)\Pic_{j\hat B}-\tr\thetac \ep_{\hat B}+\hat \bN(\hat \bN v^j) \Pic_{jl} {\hat e}_B^l-\snc_{\hat \bN}{\hat \bN}^m \p_m v^j \Pic_{j\hat B}\\
&=-(\sl{\div}\eta_{\hat B}+\ep_{\hat A} \thetac_{\hat A\hat B}+\snc_{\hat A}\log(\bb c) \eta_{\hat A\hat B})-\tr\thetac \ep_{\hat B}+\hat \bN(\hat \bN v^j) \Pic_{jl} {\hat e}_B^l.
\end{align*}
For the last term, we further derive
\begin{align*}
\hN(\hN v^j) \Pic_{j\hat B}&=\hN(\hN v^j \Pic_{jl})\hat e_B^l -\hN v^j \hN(\Pic_{jl}) \hat e_B^l\\
&=-(\snc_\hN\ep_{\hat B}+\stc{k}_{\hN\hN}\snc_\hN \hN_{\hat B})=-\snc_\hN \ep_{\hat B}+\stc{k}_{\hN\hN}\snc_{\hat B}\log (\bb c).
\end{align*}
Combining the above two calculations implies (\ref{5.30.2.23}).

Using (\ref{5.30.1.23_com}) we derive
\begin{align*}
\snc_\hN \ep_{\hat B}& =-\snc_\hN(\snc_{\hat B} v_i \hN^i)=\snc_{\hat B}(\stc{k}_{\hN\hN})-[\snc_\hN, \snc_{\hat B}]v_i \hN^i-\snc_{\hat A}\log(\bb c)\eta_{\hat B\hat A}\\
&=\snc_{\hat B}(\stc{k}_{\hN\hN})-\thetac_{\hat B\hat A} \ep_{\hat A}+\snc_{\hat A}\log(\bb c)\stc{k}_{\hN\hN}-\snc_{\hat A}\log (\bb c)\eta_{\hat B\hat A}.
\end{align*}
This gives (\ref{5.30.3.23}).

Combining (\ref{5.30.2.23}) with (\ref{5.30.3.23}) we find
\begin{equation*}
\snc_{\hat A} \eh_{\hat A\hat B}=-\tr\thetac \ep_{\hat B}-\snc_{\hat B}(\stc{k}_{\hN\hN}+\div v+\f12\tr\eta).
\end{equation*}
 Note that the terms in the bracket gives $-\f12\tr\eta$. Hence we proved (\ref{6.1.1.23}).
Next, we prove (\ref{12.16.1.23}). It is straightforward to derive
\begin{align*}
-\snc_{\hN} \eta_{\hat A\hat B}&=\snc_{\hN}\p_i v_j \Pic^i_{i'}\Pic^j_{j'}\hat e_A^{i'}\hat e_B^{j'}+\p_i v_j \snc_\hN(\Pic^j_{j'}\Pic_{i'}^i)\hat e_A^{i'}\hat e_B^{j'}\\
&=\p_i(\hN^m \p_m v_j \Pic^j_{j'})\Pic_{i'}^i \hat e_A^{i'}\hat e_B^{j'}-\hN^m \p_m v_j \p_i \Pic^{j}_{j'} \Pic^i_{i'}\hat e_A^{i'}\hat e_B^{j'}-\p_i \hN^m \p_m v_j \Pic^{j}_{j'}\Pic^i_{i'}\hat e_A^{i'} \hat e_B^{j'}\\
&-\ep_{\hat A}\snc_{\hat B}\log (\bb c) 
\end{align*}
and $\p_i \Pic_{j'}^{j}\hat e_A^{i'}\hat e_B^{j'}\Pic^i_{i'}=-\hN^j\thetac_{\hat A\hat B}$. Thus 
\begin{equation*}
\snc_\hN\eta_{\hat A\hat B}=\snc_{\hat A}\ep_{\hat B}-\thetac_{\hat A\hat C}\eta_{\hat C\hat B}+\thetac_{\hat A\hat B}\stackrel{\circ}k_{\hN\hN}+\ep_{\hat A}\snc_{\hat B}\log (\bb c).
\end{equation*}
 Taking the traceless part yields (\ref{12.16.1.23}).
 
  Using (\ref{k1}), it is direct to calculate
\begin{align*}
\nab_g^i k_{ij}&= c^2 \delta^{mi}\p_m\{ c^{-2}\big(\bT(\log c)\delta_{ij}-\f12 (\p_i v_j+\p_j v_i)\big)\}-c^2\p_l \log c k_{lj}+\p_j \log c \Tr k\\
&=-2\p^i \log c (\bT\log c \delta_{ij}-\f12 (\p_i v_j +\p_j v_i))+\p_j \bT \log c-\f12 (\Delta_e v_j+\p_j \div v)\\
&-c^2\p_l \log c k_{lj}+\p_j \log c \Tr k
\end{align*}
as stated in (\ref{2.18.4.24}).

 (\ref{9.29.5.23}) is obtained by radially decomposing (\ref{4.23.1.19})  and integrating on spheres. Similarly by radially decomposing $\tr\eta$ and integrating on spheres, we obtain (\ref{9.29.6.23}).  
\end{proof}
Next we derive a set of important formulas on the derivatives of $\bg$. 
 \begin{lemma}\label{dg}
There hold $\Ga_{X\bN}^i=0$ with $X=L, \Lb$, and the following formulas for $Y$ a linear differential operator,  
\begin{equation}\label{10.6.1.22}
\begin{split}
Y\bg_{\mu\nu} L^\mu L^\nu&=-2c^{-2}\bN^i Y v_i-2 Y \log c, \quad Y \bg_{\mu\nu} \Lb^\mu \Lb^\nu=2c^{-2} Y v^i \bN^i-2Y\log c\\
Y \bg_{\mu\nu}L^\mu \Lb^\nu&=-2 Y\log c\\
Y \bg_{\mu\nu}\Lb^\mu e_A^\nu&, Y \bg_{\mu\nu}L^\mu e_A^\nu=-c^{-2} Y v^i e_A^i,\quad Y \bg_{\mu\nu} e_B^\mu e_A^\nu=-2\ga_{AB} Y\log c.
\end{split}
\end{equation}
In short, expressing $Y\bg$ relative to null tetrad takes the form of $Y (\Phi)\c Z$, where $\c Z$ denotes multiplication by $\Pi$, $c^{m}$, and $\bN^i$.
\end{lemma}
Indeed,
\begin{align*}
L^\mu L^\nu Y \bg_{\mu\nu}&=Y \bg_{00}+2L^0 L^i Y \bg_{0i}+L^i L^j  Y \bg_{ij}\\
&=Y(-1+c^{-2} |v|^2)+2(\bN^i+v^i)Y(-c^{-2}v_i)+(\bN^i+v^i)(\bN^j+v^j)Y(c^{-2}\delta_{ij})\\
&=Y(c^{-2} |v|^2)+2(\bN^i+v^i)Y(-c^{-2}v_i)+|v|^2 Y(c^{-2})+2 \bN^i v^j\delta_{ij}Y(c^{-2})\\
&+\bN^i \bN^j \delta_{ij} Y(c^{-2})\\
&=-2\bN^i Y(c^{-2}v_i)+2\bN^i v_j \delta_{ij} Y(c^{-2})+c^2 Y (c^{-2})\\
&=-2c^{-2}\bN^i Y v_i-2 Y \log c;\\
Y \bg_{\mu\nu}\Lb^\mu \Lb^\nu&= Y \bg_{00}\Lb^0\Lb^0+2 Y\bg_{0i}\Lb^0\Lb^i+Y\bg_{ij}\Lb^i \Lb^j\\
&=Y(-1+c^{-2}|v|^2)-2 Y(c^{-2}v^i)(v^i-\bN^i)+Y(c^{-2})\delta_{ij}(v^i-\bN^i)(v^j-\bN^j)\\
&=c^{-2} Y(|v|^2)-2c^{-2}Y v^i(v^i-\bN^i)+Y(c^{-2}) c^2\\
&=2c^{-2} Y v^i \bN^i-2Y\log c;\\
 \displaybreak[0]
Y \bg_{\mu\nu}L^\mu \Lb^\nu&=Y \bg_{00}L^0 \Lb^0+Y\bg_{ij}L^i \Lb^j+ Y\bg_{i0}L^i \Lb^0+Y\bg_{i0} L^0 \Lb^i\\
&=Y(-1+c^{-2}|v|^2)+Y(c^{-2}) \delta_{ij}(\bN^i+v^i)(v^j-\bN^j)+2Y(-c^{-2}v^i) v^i\\
&=-2Y\log c;\\
Y \bg_{\mu\nu}L^\mu e_A^\nu&=Y\bg_{0i} L^0 e_A^i+Y \bg_{ij} e_A^i L^j\\
 &=Y(-c^{-2}v_i)e_A^i+Y(c^{-2})\delta_{ij} (v^j+\bN^j) e_A^i=-c^{-2}Y(v_i)e_A^i;\\
Y \bg_{\mu\nu}\Lb^\mu e_A^\nu&=Y \bg_{0i}\Lb^0 e_A^i+Y \bg_{ij}e_A^i \Lb^j\\
&=Y(-c^{-2}v_i)e_A^i+Y(c^{-2})\delta_{ij} (v^j-\bN^j) e_A^i=-c^{-2}Y(v_i)e_A^i;\\
Y\bg_{\mu\nu}e_B^\mu e_A^\nu&=Y\bg_{ij}e_B^i e_A^j=Y(c^{-2})e_B^i e_A^i=-2\ga_{AB} Y\log c.
\end{align*}
Thus the proof of (\ref{10.6.1.22}) is complete.
\subsection{Symbolic notations}
For convenience, we fix the following conventions:
\begin{itemize}
\item $[X F]$:  $X F^i \bN^j \bg_{ij}$ for  3-vectorfields $F$ and $X(F)$ for scalar functions $F$
\item $X F^\dagger$:  $[X F]$, $XF^\|:=X F\c \Pi$ for 3-vectorfields $F$
\end{itemize}
where $X$ is one of the differential operators $L, \Lb, S, \Omega$,  and the covariant derivative $\sn$ on $(S_{t,u}, \ga)$.

We classify the important geometric quantities according to their asymptotic behaviors:
\begin{equation}\label{symlst}\tag{SYB}
\begin{split}
&\fB: k_{\bN\bN},\, [\bp \Phi]; \quad \bA_b:\tr\chi-\frac{2}{\tir},\,\tr\thetac-\frac{2}{\ckr}\\
&\bA_{g,1}: [\sn\Phi],\, \eh, \Lb v^i\Pi_i^j,\,L v^i\Pi_i^j,\, \ud\bA: \ze, \bA_{g,1}\\
&\bA_{g,2}: \chih,\, \hat \theta_{AB};\quad \bA^\natural: [L\Phi],\bA_{g,1}\\
&\bA_g=\bA_{g,1}, \bA_{g,2};\quad \bA:\bA_g, \bA_b, \bA^\natural
\end{split}
\end{equation}
where $\ckr=c\tir$. Since a factor of $c^m$, with $m\in\mathbb Z$, does not change the asymptotic behavior, in future application of the above symbolic convention, we may ignore such a factor.
  
Next we give two important facts about the pair of quantities $[X v]$ and $X v^\|$.
\begin{proposition}\label{6.7.1.23}
\begin{enumerate}
\item
We have by using (\ref{4.23.1.19}) and (\ref{7.04.11.19}) that
\begin{align*}
&-\p_A v_B=\eta_{AB}=\eh_{AB}+\f12\delta_{AB} c^2\tr\eta \\
&Xv^i\Pi_{ij}e_A^j=k_{A\bN}\pm \sn \varrho=[\sn\Phi], \,X=L, \Lb\\
&\tr \eta=[L\Phi].
\end{align*}
Hence for $X$ being any operator among $\{L, \Lb, \sn\}$, symbolically we decompose $X v^\dagger$ into
\begin{equation*}
X v^\dagger: [\sn \Phi],\, [\Lb\Phi], \,\eta
\end{equation*}
where we dropped the $c^m$ factor, which do not affect the analysis in application. 
 
\item Let $X$ be one of the differential operators $L, \Lb, S, \Omega$ or the covariant differentiation $\sn$ on $(S_{t,u}, \ga)$.
By using Proposition \ref{6.7con}, it follows directly that for a vector-field $Y$ in the spacetime, and the 3-vector field $F^i$  there hold
\begin{equation}\label{5.23.1.23}
\begin{split}
&Y[XF]=[YX F]+(XF^\|)_A Y\bN^A+Y\log c[XF]\\
&\sn_Y (XF^\|)_A=e^j_A\sn_Y(XF^i\c \Pi_{ij})=(Y X F)^\|_A+[XF]Y \bN_A+Y\log c(XF^\|)_A\\
&YXF^i=(YX F)^A e_A^i+YXF\c\bN  \bN^i
\end{split}
\end{equation}
where the contraction in the above is taken with respect to $\bg_{ij}$ and  we neglected the actual value of the constant coefficients of the last terms in the first two lines.
\item 
Define  
\begin{align*}
	\Sc(X^{n}v)&= X^{n-1}[X_1 v], X^{n-2}[X_2 X_1 v], \cdots, [X^n v]\\
	\Ac(X^n v)&=\sn_X^{n-1}(X_1 v^\|), \sn_X^{n-2}( (X_2 X_1 v)^\|), \cdots, (X^n v)^\|.
\end{align*}

If $Y\in\{\Omega, S\}$, in view of Proposition \ref{6.7con} and calculating by using $\Ga_{\bN Y}^i=0$ with if $Y$ is null (from Lemma \ref{dg}), symbolically,
\begin{equation}\label{11.30.2.23}
Y \bN^A=(1-\vs(Y))c\thetac^A(Y)+\vs(Y)\tir[\sn\Phi],\quad \Lb  \bN^A=\zb^A+\ze^A.
\end{equation}
\end{enumerate}
\end{proposition}

\begin{definition}
\begin{enumerate}
\item  Let $f$ be a nowhere trivial, smooth scalar function. We define the signature of the differential operators $f L, f\Lb, f\sn, \Omega$ by
 $$\vs(f L)=1, \vs( f\Lb)=-1, \vs(f\sn)=\vs(\Omega)=0,$$
 and  for $Y_1,\cdots, Y_n$ among the above differential operators
 $$\vs(Y_1\cdots Y_n)=\vs(Y_1)+\cdots+\vs(Y_n).$$
 Moreover, we set
\begin{equation}\label{5.15.2.23}
 \vs^-(X^n)=\min\{\vs(X_1),\cdots, \vs(X_n)\}, \quad \vs^+(X^n)=\max\{\vs(X_1),\cdots, \vs(X_n)\}.
\end{equation}
 \end{enumerate}
\end{definition}
\subsection{Geometric null conditions}
\subsubsection{Null conditions in the irrotational Euler equation}
For $\Sigma$-tangent vector-fields or scalars $W$ and $V$, we denote scalar bilinear forms
\begin{align*}
\sG(W,V)&=\bg^{\a\b}\p_\a W \p_\b V \\
\sB(W,V)&=\sum_{1\le c<d\le 3}(\p_c W^c \p_d V^d-\p_c W^d \p_d V^c).
\end{align*}
Recast (\ref{4.10.1.19}) and (\ref{4.10.2.19}) by using the above conventions, symbolically
\begin{equation}\label{8.25.2.22}
\Box_\bg v^i=\sG(\varrho, v^i); \quad \Box_\bg \varrho=\sG(\varrho, \varrho)+\sB(v, v).
\end{equation}
It is straightforward to derive 
\begin{align*}
\sB(W,V)&=\sl{\p}W \sl{\p}V+\sl{\p}V\p W+\p V\c \sl{\p}W\\
\sG(W,V)&=-\f12(\Lb W LV+\Lb V LW)+\sn W\sn V.
\end{align*}
where  $\p_i=\sl{\p}_i+\l\p_i, \hat \bN\r_\delta \hat \bN$ and $\hat \bN=c^{-1}\bN$. With $\bar\bp \Phi=\sn \Phi, L \Phi$ and $Z$ either in null tetrad, or being $c^m$ with $m\in \mathbb Z$,
 we can symbolically write
\begin{equation*}
\Box_\bg \Phi=(\bar\bp \Phi, \sl{\p}\Phi)\c \bp \Phi\c Z.
\end{equation*}
To avoid treating $\sl{\p}$ in this paper, we will present a geometric formulation of the above two bilinear forms.

Now we derive the geometric decomposition for $\sB(v,v)$, which is important for our analysis.
\begin{lemma}
The following quadratic forms for $\bp v$ take the following form
\begin{align}
\sG(v, v), [\sG(v, \varrho)], \sB(v,v)=\sum_{\vs(Y_1Y_2)=0, 2}[Y_1 \Phi][Y_2\Phi]+|\eh|^2 \label{null}
\end{align}
\begin{footnote}{We use $[F]$ to represent either the radial component $F_\bN$ if $F$ is a $3$-vector, or $F$ itself if $F$ is a scalar. }\end{footnote}where $\sG(v, v)=\sum_{i=1}^3\sG(v^i, v^i)$,  $Y_1, Y_2$ are in null tetrad, and  the factors of the type $c^m$ on the right-hand side are neglected. 
\end{lemma}
Indeed
\begin{align*}
\sB(v,v)&=-\f12(\stc{k}\c\stc{k}-(\Tr\stc{k})^2)\\
&=-\f12(\sk_{AB}\sk_{AB}+2\stc{k}_{A\hat\bN} \stc{k}_{A\hat \bN}+\stc{k}_{\hat\bN\hat \bN}^2-(\tr \sk+\stc{k}_{\bN\bN})^2)\\
 &=-\f12(|\eh|^2+2|\ep|^2-\f12(\tr\eta)^2-2\tr\eta \stc{k}_{\bN\bN}).
\end{align*}
Then using (\ref{7.04.9.19}), we obtain the formula of $\sB(v,v)$ in (\ref{null}). It is direct to see that \begin{footnote}{We hide the summation over $i=1,2,3$ by using Einstein convention.}\end{footnote}
\begin{align*}
\sG(v, v)&=-\Lb v^i L v^i+\sn v^i \sn v^i=\sum_{\vs(Y_1Y_2)=0, 2}[Y_1 \Phi][Y_2\Phi]+|\eh|^2.
\end{align*}
The identity for $[\sG(v, \varrho)]$ follows similarly. Thus (\ref{null}) is proved.

\subsubsection{Decomposition of $\bR_{44}$}
We prove an crucial decomposition for $\bR_{44}$, which is the key structure for controlling $\tr\chi-\frac{2}{\tir}$ and its derivatives. 
\begin{proposition}\label{ric44}
There holds the following important decomposition for $\bR_{44}$,
\begin{equation}\label{3.7.6.21}
\bR_{44}=L(\Xi_4)-k_{\bN\bN}\Xi_4+ \N(\Phi, \bp \Phi)
\end{equation}
where $\N(\Phi, \bp \Phi)$ symbolically represents a set of scalars, with $n,m\in \mathbb Z$, taking the following form 
\begin{equation}\label{10.8.1.22}
\N(\Phi, \bp\Phi)=\sum_{\vs(Y_1 Y_2)=0,2} c^n [Y_1\Phi][Y_2\Phi]+c^m |\eh|^2.
\end{equation}
\end{proposition}
 \begin{proof}
We first recall that relative to any coordinate, the Ricci curvature of $\bg$ can be decomposed as
\begin{equation}\label{ricc6.7.1}
\bR_{\a\b}=-\f12 \Box_\bg (\bg_{\a\b})+\f12 (\bd_\a \Xi_\b+\bd_\b \Xi_\a)+S_{\a\b}
\end{equation}
where $\Xi$ is the 1-form defined in (\ref{ricc6.7.2}).
The term
$S_{\a\b}$ is quadratic in $\bp \bg$, (see \cite[Page 25]{Andersson_Moncrief})
\begin{equation}\label{s44}
\begin{split}
S_{\mu\nu}&=\f12 \bg^{\b\a}\bg^{\sigma \ga}(\p_\nu \bg_{\a \ga}\p_\sigma \bg_{\mu \b}+\p_\mu \bg_{\b \sigma} \p_\ga \bg_{\nu\a}
-\f12 \p_\nu \bg_{\a\ga}\p_\mu \bg_{\b\sigma}\\
&+\p_\b \bg_{\mu\sigma}\p_\a\bg_{\nu \ga}-\p_\b \bg_{\mu\sigma} \p_\ga\bg_{\nu\a})-\f12 \Xi^\a\p_\a \bg_{\mu\nu},
\end{split}
\end{equation}
where  the reference connection is already taken with respect to the Minkowski metric.

Recall from \cite[Proposition 7.5]{rough_fluid} that
\begin{align}
\bR_{44}&=\delta^{ij} c^{-2} \bN^j \sQ^i-\f12 c^2 \Box_\bg (c^{-2})+\bd_L \Xi_L +S_{44}\nn\\
&-(c^{-2} \bd^\a v^i \bd_\a v^j \delta_{ij}-2\delta_{ij}\bN^j \bd^\a(c^{-2}) \bd_\a(v^i))\label{6.23.2.19}.
\end{align}
In view of (\ref{4.10.1.19}), (\ref{4.10.2.19}) and using (\ref{6.29.5.19}) we derive
\begin{equation}\label{3.7.7.21}
\bR_{44}=L(\Xi_4)+k_{\bN\bN} \Xi_4+S_{44}+c^{-2}(\sG(v^i, v^j)\delta_{ij}+\sG(\varrho, v^i) \bN^i)+\sG(\varrho, \varrho)+\sB(v,v).
\end{equation}

 Next we claim
\begin{equation}\label{3.7.8.21}
S_{44}=-4(k_{\bN\bN})^2+\N(\Phi,\bp \Phi).
\end{equation}
Using (\ref{7.04.8.19}) which gives $\Xi_4-2k_{\bN\bN}=[L\Phi]$, also using (\ref{3.22.1.21}), we obtain (\ref{3.7.6.21}) by combining the above two identities.

Now we prove (\ref{3.7.8.21}). Using (\ref{s44}), we first prove
\begin{equation}\label{10.6.2.22}
S_{44}=\N(\Phi, \bp \Phi)-\f12 (\bg^{\b\a}\bg^{\sigma \ga}\p_\b \bg_{\mu\sigma} \p_\ga\bg_{\nu\a}+\Xi^\a\p_\a \bg_{\mu\nu})L^\mu L^\nu.
\end{equation}
Indeed, in view of the last line in Lemma \ref{dg}, writing terms in the first line of (\ref{s44}) relative to null tetrad gives a scalar $L \Phi \c Y \Phi\c Z$, with $Y$ in null tetrad. Due to being scalars, here are the only possibilities:
\begin{equation*}
  L\log c \c \tr \eta, \quad L \Phi^A Y\Phi_A, \quad [L\Phi][Y\Phi], \quad Y=L, \Lb.
\end{equation*}
 Hence, using $\tr\eta=[L\Phi]$ from (\ref{7.04.9.19}), they can only take the form of $\sum_{\vs(Y_1 Y_2)=0,2} [Y_1 \Phi][Y_2\Phi]$, in view of Proposition \ref{6.7.1.23} (1).
 
 It remains to consider the first term of the second line in (\ref{s44}), which is a scalar, and a full contraction $[\sG(\Phi,\Phi)]$.  By applying the formula for such terms in (\ref{null}), and noting that $\sC(\varrho, \varrho)$ is obviously of the desired form, we conclude (\ref{10.6.2.22}).
 
Next we compute the second term on the right-hand side of (\ref{10.6.2.22}). Using (\ref{10.6.1.22}), it is direct to see that $L^\mu L^\nu \bg^{\la s} \bg^{\rho s'}\p_\rho \bg_{s\mu}\p_\la \bg_{s' \nu}$ takes the form
\begin{itemize} 
\item $\sum_{\vs(Y_1 Y_2)=0, 2}c^m [Y_1\Phi][Y_2\Phi]$, if the indices $\rho$ or $\la=L$;
\item $c^m \eta\c \eta+c^n[\sn\Phi]^2,$  if both $\rho, \la$ are contracted by $e_A, e_B$;
\item $c^m[\sn\Phi][\sn\Phi]$ if one of $\rho, \la=\Lb,$ while the other is contracted by $e_A$,    
\end{itemize}
where $m,n\in\mathbb Z$ all are some integers, whose particular values are not important for our analysis.    
 Note that applying the first line in (\ref{10.6.1.22}) to $Y=\Lb$, also using (\ref{k1}), we have
\begin{equation}\label{3.8.1.21}
L^\mu L^\nu \Lb \bg_{\mu\nu}=-4k_{\bN\bN}+2(c^{-2}\bN^j L v^i +L \log c).
\end{equation}
Hence, we obtain
\begin{align*}
\begin{split}
-\frac{1}{2}L^\mu L^\nu \bg^{\la s} \bg^{\rho s'}\p_\rho \bg_{s\mu}\p_\la \bg_{s' \nu}
&=-\f12 L^\mu L^\nu  \bg^{34}\bg^{34} \Lb \bg_{\a\mu} \Lb \bg_{\b\nu}L^\a L^\b+\N(\Phi, \bp \Phi)\\
&=-\frac{1}{8} (L^\nu L^s \Lb \bg_{s\nu})^2+\N(\Phi, \bp \Phi)\\
&=-2k_{\bN\bN}^2+\N(\Phi, \bp \Phi).
 \end{split}
\end{align*}

Hence the second term of (\ref{10.6.2.22}) contains a bad term $-2(k_{\bN\bN})^2$ apart from good quadratic terms $\N(\Phi, \bp\Phi)$.

For the third term on the right-hand side of (\ref{10.6.2.22}), we derive by using (\ref{7.04.8.19}), (\ref{3.22.1.21}) and (\ref{3.8.1.21})
\begin{equation*}
-\f12\Xi^\a \p_\a \bg_{\mu\nu}L^\mu L^\nu=\f12 k_{\bN\bN}\Lb \bg_{LL}-\f12(\Xi^L L\bg_{\mu\nu}+\Xi^A \sn_A \bg_{\mu\nu})L^\mu L^\nu+[L\Phi]([\Lb \Phi]+[L\Phi]).
\end{equation*}
By using (\ref{10.6.1.22}) for $Y=\Lb, L, e_A$ respectively, also using (\ref{6.14.1.19}) we obtain
\begin{equation*}
-\f12\Xi^\a \p_\a \bg_{\mu\nu}L^\mu L^\nu=\f12 k_{\bN\bN}\Lb \bg_{LL}+\N(\Phi, \bp\Phi)=-2k_{\bN\bN}^2+\N(\Phi, \bp\Phi).
\end{equation*}
Using (\ref{3.8.1.21}) again, hence we conclude (\ref{3.7.8.21}).
\end{proof}

In view of (\ref{null}), we summarize the important null forms appeared in the geometric calculations and Euler equations.
\begin{proposition}\label{geonul_5.23_23}
For  $\N(\Phi, \bp\Phi)$ defined in (\ref{10.8.1.22}) and $\widetilde{\N}(\Phi, \bp\Phi):=\sQ^0, \sQ^i$ which appeared in (\ref{8.25.2.22}), there hold  
\begin{align}\label{5.22.2.23}
\sB(v,v), \sG(v,v)&=\sum_{\vs(Y_1Y_2)=0, 2}[Y_1 \Phi][Y_2\Phi]+|\eh|^2\nn\\
\sG(v, \varrho)&=-\f12(L\varrho\Lb v+\Lb \varrho Lv)+\sn \varrho \sn v\nn\\
\displaybreak[0]
\sG(v, \varrho)(\hat \bN), \sG(\varrho, \varrho)&=\sum_{\vs(Y_1 Y_2)= 0}[Y_1\Phi][Y_2 \Phi]\\
\N(\Phi, \bp\Phi)&=\sum_{\vs(Y_1 Y_2)=0,2} [Y_1\Phi][Y_2\Phi]+|\eh|^2\nn\\
\widetilde\N(\Phi, \bp\Phi)&=\N(\Phi, \bp\Phi)+\sG(v, \varrho)\nn
\end{align}
where $Y_1, Y_2$ are in null tetrad, and we dropped the factors of $c^n$ in front of the terms $[Y_1\Phi][Y_2\Phi]$ or $|\eh|^2$ in the above.  
\end{proposition}
 In application, we can drop the terms with $\vs(Y_1Y_2)=2$  since they are much better terms than those with $\vs(Y_1 Y_2)=0$.

Due to our small-large regime, we give the following signature rule to classify terms according to their amplitudes.
\begin{definition}\label{5.23.2.23}
We assign an index of amplitude $\al(F)$ for a tensor-field or scalar $F$, according to the maximal exponent of $\ve$ used in this paper in the following  bound at the initial slice,
\begin{equation*}
\|F\|_{L_\omega^4(S_{0,u})}\les \ve^{\al(F)}.
\end{equation*}
We use $\al$-value to classify quantities into three sets $\al_0, \al_\f12, \al_1$.
\end{definition} 
 Based on the above convention, we have
$$ [\Lb\Phi], \fB\in \al_0, \quad  [L\Phi],\tr\eta, \bA_b\in \al_\f12, \quad [\sn\Phi],\zb, \eh, \bA_g\in  \al_1,$$
which can be directly justified by analyzing initial data (see Proposition \ref{12.21.1.21}). 
\begin{itemize}
\item There hold the rules of products and sums:
\begin{equation*}
 \al(F_1\c F_2\cdots F_n)=\sum_{i=1}^n\al(F_i), \al(F_1+F_2+\cdots+F_n)=\min_{i=1}^n(\al(F_i)).
\end{equation*}
\item There holds the following differentiation rule: each column in the following table represents the class which the quantity belongs to after being differentiated by the corresponding operator at the beginning of the row.
\begin{equation} \label{5.24.1.23}\tag{\bf sig}
\begin{array}{|l|l|l|l|}
\hline
  & \al_0 & \al_\f12 & \al_1 \\
 \hline
\Lb & \al_0 & \al_0 & \al_1 \\
\hline
\sn & \al_1 & \al_1 & \al_1\\
\hline
L & \al_0 &\al_\f12 & \al_1  \\
\hline
\end{array}
\end{equation}
\end{itemize}

We have a few remarks about the index of amplitude.
\begin{remark}\label{5.24.3.23}
The table applies to all the elements in each set of $\al_0, \al_\f12, \al_1$, not depending on the particular form of the quantities. We use it inductively for higher order derivatives.  
\end{remark}
\begin{remark}\label{5.24.4.23}
The calculation rule and differentiation rule are summarized based on the bounds achieved from the initial data. 
We will refer to the table as a guide to establish bootstrap assumptions shortly. In this paper, we will show that the amplitude index of the quantities in the bootstrap assumption remains the same over time.
\end{remark}
\section{Preliminaries of causal geometry}\label{causal_geo}
\subsubsection{Commutation formulas}
We recall the following commutation relations used in \cite[Section 5]{Wangrough} (see also in \cite{KRduke}).
\begin{proposition}
(1) There hold for the scalar functions $f$ that
\begin{align}
[L,\bT]f&=\f12[L, \Lb]f=(\zb^A-\zeta^A)\sn_A f-k_{\bN\bN} \bN f\label{3.19.2}\\
[L, \bb \Lb]f&=2\bb(\zb_A-\zeta_A) \sn_A f-\bb \delta L f.\label{8.31.5.19}
\end{align}

Moreover for $S_{t,u}$ tangent tensor field $U_A$,
\begin{equation}\label{7.04.5.21}
\f12[\sn_L, \sn_\Lb]U_A=(\zb-\zeta)\sn U_A-k_{\bN\bN} \sn_\bN U_A+R_{ABL\Lb} U_B+\zeta\c \zb\c U
\end{equation}
where the last two terms are symbolic terms.
%

(2) There holds for $S_{t,u}$-tangent $m$-covariant tensor fields $U_A$ that
\begin{equation}\label{cmu2}
\begin{split}
&\sn_L\sn_B U_A-\sn_B \sn_L U_A  \\
&=-\chi_{BC}\c \sn_C U_A+\sum_{i}(\chi_{A_i B} \zb_C-\chi_{BC} \zb_{A_i}
+\bR_{{A_i}C4 B})U_{A_1\cdots\ckk C\cdots A_m}
\end{split}
\end{equation}
and for any scalar function $f$ there holds
\begin{equation}\label{cmu_2}
[L,\sn_A] f=- \chi_{AB} \sn_B f.
\end{equation}
Thus, for an $S$-tangent one tensor $U$
\begin{equation}\label{6.22.16.19}
[\sn_L, \sl{\div}]U=-\chi\c \sn U+\chi\c\zb\c U+\bR_{AC4B} \delta^{AB} U_C
\end{equation}
\begin{equation}\label{2.18.1.22}
[\sn_L, \sl{\curl}]U=-\chih\c \sn U-\f12 \tr\chi\c \sl{\curl} U+\chi\c \zb\c U+\bR_{AC4B}\ep^{AB}U_C.
\end{equation}
Consequently, for any scalar function $f$ there holds
\begin{align}
L\sD f+  \tr\chi \sD f&= \sD L f-2\chih\c \sn^2 f-\sn_A\chi_{AC}\sn_C f\nn
\\&+(\tr\chi \zb_C-\chi_{AC}\zb_{A}-\delta^{AB}\bR_{CA4B}) \sn_C f.\label{tran1}
\end{align}
(3)
 \begin{equation}\label{7.03.1.19}
\begin{split}
[\sn_\Lb, \sn_B]U_A&=-\chib_{BC} \sn_C U_A+(\zeta_B-k_{B\bN})(\sn_\Lb-\sn_L) U_A\\
&+\Big(\chi\c (-\zeta+k_{A\bN})+\chib\c \zeta+\bR_{AC\Lb B}\Big)U_C.
\end{split}
\end{equation}
\begin{equation}\label{4.6.1.24}
\begin{split}
[\bb\sn_\Lb, \sn_B]U_A&=-\bb\chib_{BC} \sn_C U_A-\bb(\zeta_B-k_{B\bN})\sn_L U_A\\
&+\bb\Big(\chi\c (-\zeta+k_{A\bN})+\chib\c \zeta+\bR_{AC\Lb B}\Big)U_C.
\end{split}
\end{equation}
(4) If $X$ is $S_{t,u}$ tangent,
\begin{align}\label{8.10.1.21}
[\sn_X, \sn_A]F_B&=(X^BF_A-\ga_{BA} X^D F_D)K-\sn_A X^D \sn_D F_B.
\end{align}
\end{proposition}
Indeed, by using (\ref{lb}) and (\ref{3.19.2}), we derive
\begin{equation*}
[L, \bb \Lb]f=L\bb \Lb f+\bb[L,\Lb]f=-\bb k_{\bN\bN} \Lb f+2\bb(\zb-\zeta)\sn f-2\bb k_{\bN\bN} \bN f
\end{equation*}
(\ref{8.31.5.19}) follows by using $2\bN=L-\Lb$. (\ref{4.6.1.24}) is derived by using (\ref{7.03.1.19}) and $\ze=\sn\log \bb+k_{A\bN}$. 

(\ref{8.10.1.21}) follows from the calculation below
\begin{align*}
[\sn_X, \sn_A]F_B&=X^C(\sn_C \sn_A F_B-\sn_A \sn_C F_B)-\sn_A X^C \sn_C F_B\\
&=X^D R_{BCD A}F_C-\sn_A X^D \sn_D F_B\\
&=X^D K(\ga_{BD}\ga_{CA}-\ga_{BA}\ga_{CD})F_C-\sn_A X^D \sn_D F_B.
\end{align*}
All other commutation formulas have been given in the references. 

\subsubsection{Conformal change of metric}
In
$\Sigma_t$, we will frequently use the standard Euclidean metric in calculation for convenience. Recall that $\stc\ga$ is the induced metric of ${\delta_e}_{ij}$ on $S_{t,u}$, and $\snc$ is its Levi-Civita connection on $S_{t,u}$.  Since $g_{ij}=c^{-2}\delta_{ij}$, the induced metric of $g_{ij}$ on $S_{t,u}$ verifies $\ga_{ij}=c^{-2}\gac_{ij}$ in the tangent space of $S_{t,u}$. We have the following formulas for $S_{t,u}$-tangent vector-field $F$,  $S_{t,u}$- tangent  1-form $G$ and scalar functions $f$, 
\begin{equation}\label{1.27.1.22}
\sl{\div}F=\divsc F+\snc\log c \c F, \quad \scurl G=c^2 \curlsc G, \quad \sD f=c^2 \Delsc f,
\end{equation}
where $\divsc$, $\curlsc$, and $\Delsc$ are the divergence, curl and Laplace-Beltrami operators associated to $(\snc, \gac)$.

\subsubsection{Null Structure Equations}

 We will rely heavily on the following structure equations for analysing the connection coefficients
on null hypersurfaces $\H_u$ (see \cite[Section 5]{Wangrough} and \cite[Section 7]{rough_fluid}):
\begin{proposition}[Transport equations and Hodge systems for connection coefficients]\label{6.29.1.24}
\begin{align}
&L \bb=-\bb { k}_{\bN\bN}, \label{lb}\\
&L\tr\chi+\f12 (\tr\chi)^2=-|\chih|^2-{ k}_{\bN\bN} \tr\chi-\bR_{44}, \label{s1} \displaybreak[0]\\
&\sn_L \chih_{AB}+ \tr\chi \chih_{AB}=-{k}_{\bN\bN} \chih_{AB}-(\bR_{4A4B}-\f12 \bR_{44} \delta_{AB}), \label{s2} \displaybreak[0]\\
&L \tr\chib+\f12 \tr\chi \tr\chib=2\sl{\div} \zb+k_{\bN\bN} \tr\chib-\chih\c \chibh+2|\zb|^2+\delta^{AB}\bR_{A34B}, \label{mub} \displaybreak[0]\\
&\sn_L \zeta+\f12\tr\chi \zeta=-(k_{B\bN}+\zeta_B) \chih_{AB}-\f12 \tr\chi k_{A\bN}-\f12 \bR_{A4 43}, \label{tran2} \displaybreak[0]\\
&(\sl{\div} \chih)_A+\chih_{AB}\c k_{B\bN}=\f12(\sn \tr\chi+k_{A\bN} \tr\chi)+\bR_{B4BA}, \label{dchi} \displaybreak[0]\\
&\sl{\div} \zeta=\f12(\mu-k_{\bN\bN} \tr\chi-2|\zeta|^2-|\chih|^2-2k_{AB}\chih_{AB})-\f12\delta^{AB}\bR_{A43B}, \label{dze} \displaybreak[0]\\
&\sl{\curl} \zeta=-\f12 \chih\wedge \chibh+\f12 \ep^{AB}\bR_{B43A}, \label{dcurl} \displaybreak[0]\\
&\sl{\curl} \zb=-2\zb_A \sn_B\log c\ep^{AB}+\ep^{BC}\eh_{AC}\hat\theta_B^A,\label{7.04.18.19}\displaybreak[0]\\
&\sn_\Lb \chih_{AB}+\f12 \tr\chib \chih_{AB}=-\f12 \tr\chi\chibh_{AB}+2\sn_A\zeta_B-\sl{\div} \ze \delta_{AB}+k_{\bN\bN} \chih_{AB} \label{3chi}
\displaybreak[0]\\
&\quad\quad\quad\quad\quad \quad\quad+(2\zeta_A\ze_B-|\ze|^2\delta_{AB})+\bR_{A43B}-\f12 \delta^{CD}\bR_{C43D}\delta_{AB},\nn
\end{align}
where the mass aspect function $\mu:=\Lb \tr\chi+\f12 \tr\chi \tr\chib$; and for an $S_{t,u}$-tangent tensor field $F$,  $\sn_L F:=L^\mu\bd_\mu
F$, with $\bd$ the covariant derivative of $(\M,\bg)$.

\end{proposition}
Indeed, it suffices to check  (\ref{7.04.18.19}). Other identities have been checked in \cite{KRduke, Wangrough}, etc.  
Note that by using (\ref{7.04.11.19}) and the formula for $\zb$ in Proposition \ref{6.7con}
\begin{equation*}
\zb=\sn v^i \bN^j g_{ij}.
\end{equation*}
 It is straightforward to derive
\begin{align*}
\sl{\curl} \zb=\ep^{AB}\sn_B (\sn_A  v^i \bg_{ij}\bN^j)&=\ep^{AB} \p_A v_C \hat\theta_B^C-2\zb_A \sn_B\log c\ep^{AB}.
\end{align*}
This gives (\ref{7.04.18.19}).
	

Next we derive a set of structure equations by using our important observations on geometric structures.
\begin{lemma}[Normalized structural equations]
\begin{enumerate}
\item Let $\sX=\tr\chi+\Xi_4$ and $h=\f12\tr\chi$.
\begin{align}
&L \tr\chi+\f12 (\tr\chi)^2=-|\chih|^2-k_{\bN\bN}(\tr\chi-\Xi_4)-L(\Xi_4)+\N(\Phi, \bp \Phi)\label{8.31.6.19}\\
&L\sX+\f12 \sX^2=\f12\Xi_4^2+\sX(\Xi_4-k_{\bN\bN})-|\chih|^2+\N(\Phi, \bp\Phi)\label{3.8.3.21}
\end{align}
\item
\begin{equation}\label{6.3.1.23}
\begin{split}
 L(\tr\chi-\frac{2}{\tir})&+\frac{2}{\tir}(\tr\chi-\frac{2}{\tir})+\f12 (\tr\chi-\frac{2}{\tir})^2=L \tr\chi+\f12 (\tr\chi)^2\\
&=-|\chih|^2-\widetilde{L \Xi_4}-\tr\chi(k_{\bN\bN}-\f12 \Xi_4)+\N(\Phi, \bp \Phi),
\end{split}
\end{equation}
where
\begin{align}
\widetilde{L \Xi_4}:&=L (\Xi_4)+(h-k_{\bN\bN})\Xi_4\nn\\
&=\wp(\sD \varrho-\hb L \varrho+2\zb^A \sn_A \varrho-\Box_\bg \varrho)+2(\wp-1)(L+(h-k_{\bN\bN}))L \varrho.\label{3.20.1.22}
\end{align}

\begin{align}
\sn_L \sn \tr\chi+\frac{3}{2}\tr\chi\sn\tr\chi&=-\chih\c \sn \tr\chi-2\chih\c \sn \chih-\sn\widetilde{L(\Xi_4)}-\sn(\tr\chi[L\Phi])\nn\\
&+\sn\N(\Phi, \bp \Phi);\label{9.15.3.22}
\end{align}
With $\wt{\tr\chi}=\tr\chi+2(\wp-1) L\varrho$, 
\begin{align}\label{1.26.1.23}
\begin{split}
\sn_L \sn \wt {\tr\chi}+\frac{3}{2}\tr\chi\sn \wt{\tr\chi}&=(-2\sn \chih-\sn\wt{\tr\chi})\c \chih-\wp\sn(\sD \varrho-\hb L \varrho+2\zb^A \sn_A \varrho-\Box_\bg \varrho)\\
&+  [L \Phi] (\sn \wt{\tr\chi}+\sn k_{\bN\bN}+\sn L \varrho)+ \tr\chi \sn [L \Phi]+\sn \N(\Phi, \bp\Phi)
\end{split}
\end{align}
\item With
\begin{equation}\label{3.31.5.22}
\sF=\bb^{-1}\sn \sX,
\end{equation}
there holds
\begin{equation}\label{12.5.1.21}
\begin{split}
L\sF+\frac{3}{2}\tr\chi\sF&=-\chih\c \sF+\bb^{-1}\Big(\sn \Xi_4(\Xi_4+\f12 \sX)+\sX\sn(\f12 \Xi_4-k_{\bN\bN})-2\sn \chih\c \chih\\
&+\sn \N(\Phi, \bp \Phi)\Big).
\end{split}
\end{equation}
With the right-hand side of (\ref{12.5.1.21}) denoted by $G_2$, there holds the following equation for the derivatives of $\sF$, 
\begin{align}
(\sn_L+2\tr\chi)\sn_A\sF_B=-\chih\c\sn \sF+\chi\c \zb\c \sF+\bR_{BCLA}\sF_C+\sn G_2.\label{12.5.2.21}
\end{align}
\item There hold the following equations for $\Lb \tir$ and $\mho=\Lb \log \tir-\f12 \tr\chib.$
\begin{align}
&L(\bb \Lb \tir-\bb)=0\label{8.11.8.22}\\
&L(\bb \tir \mho-\bb)=-\f12 \bb\tir(\ud\mu-k_{\bN\bN}\tr\chib)+\frac{\bb\tir}{4}\tr\chib(\tr\chi-\frac{2}{\tir})\label{8.11.7.22}\\
&(\sn_\Lb +\f12\tr\chib)\big(\tr\chi-\frac{2}{\tir}\big)=\mu-k_{\bN\bN}\tr\chi-\varpi+k_{\bN\bN}\tr\chi+\varpi+2\tir^{-1}\mho\label{8.11.9.22}\\
&\Lb(\tir \tr\chi)=\tir\tr\chi \mho+\tir\mu\label{8.13.5.22}\\
&\sn_A\Lb \tir=(\zeta-k_{A\bN})\c(1-\Lb\tir)\label{4.9.1.24}\\
\displaybreak[0]
&\sn \mho=(\zeta-k_{A\bN})(\tir^{-1}-\f12\tr\chib-\mho)-\f12\sn \tr\chib\label{8.13.3.22}\\
&\sn_L \sn \log\bb+\chi\c \sn \log \bb=-\sn k_{\bN\bN}.\label{1.27.6.24}\\
&\begin{array}{lll}
\ud\mu-k_{\bN\bN}\tr\chib-\varpi=\sn \bA_{g,1}+\bA_g\c\bA_g, \\
\mu-\tr\chi k_{\bN\bN}-\varpi=\sn\ud \bA+\ud \bA\c \ud \bA+\bA_g\c\bA_g
\end{array}\label{8.13.4.22}
\end{align}
where $\varpi=(\bA+\frac{1}{\tir}+\fB)\fB$. 
\end{enumerate}
\end{lemma}
\begin{proof}
We first combine (\ref{s1}) and (\ref{3.7.6.21}) to derive (\ref{8.31.6.19}). (\ref{6.3.1.23}) follows from (\ref{8.31.6.19}) directly. 

Applying (\ref{7.04.8.19}) to $f=\varrho$ gives 
\begin{equation}\label{6.22.2.21}
L \Lb \varrho+h \Lb \varrho-k_{\bN \bN}\Lb \varrho=\sD \varrho+\hb L \varrho+2\zb^A \sn_A \varrho-\Box_\bg \varrho.
\end{equation}
From (\ref{7.04.7.19}), we have
\begin{equation*}
\Xi_L=\wp \Lb \varrho+2(\wp-1) L \varrho.
\end{equation*}
By direct substitution and using (\ref{6.30.2.19}), we derive
\begin{align*}
L(\Xi_4)+(h-k_{\bN\bN})\Xi_4&=(L +(h-k_{\bN\bN}))(\wp\Lb \varrho+2(\wp-1) L \varrho)\\
&=\wp(L+(h-k_{\bN\bN}))\Lb \varrho+2(\wp-1)(L+(h-k_{\bN\bN}))L \varrho\nn\\
&=\wp(\sD \varrho-\hb L \varrho+2\zb^A \sn_A \varrho-\Box_\bg \varrho)+2(\wp-1)(L+(h-k_{\bN\bN}))L \varrho,\nn
\end{align*}
as stated in (\ref{3.20.1.22}).

(\ref{9.15.3.22}) and (\ref{1.26.1.23}) can be obtained by differentiating (\ref{6.3.1.23}) with the help of (\ref{cmu_2}) and (\ref{7.04.8.19}).

Using (\ref{8.31.6.19}) and (\ref{7.04.8.19}),  we infer
\begin{align*}
L\sX+\f12 \sX^2&=\Xi_4^2+\tr\chi(\Xi_4-k_{\bN\bN})-|\chih|^2+\N(\Phi, \bp\Phi)\\
&=\f12\Xi_4^2+\f12 \Xi_4^2-\Xi_4(\Xi_4-k_{\bN\bN})+\sX(\Xi_4-k_{\bN\bN})-|\chih|^2+\N(\Phi, \bp\Phi)\\
&=\f12\Xi_4^2+\sX(\Xi_4-k_{\bN\bN})-|\chih|^2+\N(\Phi, \bp\Phi).
\end{align*}
This gives (\ref{3.8.3.21}).

Differentiating (\ref{3.8.3.21}) by using (\ref{cmu_2}) implies
\begin{align*}
\sn_L \sn \sX+(\frac{3}{2}\tr\chi+k_{\bN\bN})\sn \sX&=\sn\Xi_4\c \Xi_4+\sX\sn(\Xi_4-k_{\bN\bN})-2\sn \chih\c \chih\\
&+\sn \N(\Phi, \bp \Phi)-\chih\c\sn \sX \\
&=\sn \Xi_4(\Xi_4+\f12 \sX)+\sX\sn(\f12 \Xi_4-k_{\bN\bN})-2\sn \chih\c \chih\\
&+\sn \N(\Phi, \p \Phi)-\chih\c \sn \sX.
\end{align*}
In view of (\ref{lb}), we obtain
\begin{align*} 
L(\bb^{-1}\sn \sX)+\frac{3}{2}\tr\chi\bb^{-1}\sn \sX&=\bb^{-1}\Big(\sn \Xi_4(\Xi_4+\f12 \sX)+\sX\sn(\f12 \Xi_4-k_{\bN\bN})-2\sn \chih\c \chih\\
&+\sn \N(\Phi, \bp \Phi)-\chih\c \sn \sX\Big)
\end{align*}
as stated in (\ref{12.5.1.21}). (\ref{12.5.2.21}) can be obtained by differentiating (\ref{12.5.1.21}) and using (\ref{cmu2}). 

Next we compute by using (\ref{8.31.5.19}) and (\ref{lb})
\begin{equation*}
L(\bb\Lb \tir)=-\bb k_{\bN\bN} L \tir=-\bb k_{\bN\bN}=L\bb.
\end{equation*}
This is (\ref{8.11.8.22}). We directly derive by using (\ref{8.11.8.22}) that
\begin{align*}
L\big(\bb(\Lb \tir-\f12 \tir \tr\chib)\big)& =-\f12(L \tr\chib\c \bb\tir+\bb \tr\chib+L \bb \tir \tr\chib)+L\bb\\
&=-\f12(L \tr\chib+\tir^{-1}\tr\chib-k_{\bN\bN}\tr\chib) \bb\tir+ L\bb.
\end{align*}
This implies (\ref{8.11.7.22}) in view of (\ref{mub}).

Using the definition of $\mu$, we derive
\begin{align*}
\Lb(\tr\chi-\frac{2}{\tir})&=\mu-k_{\bN\bN}\tr\chi+2\tir^{-2}\Lb \tir-(\f12\tr\chi\tr\chib-k_{\bN\bN}\tr\chi)\\
&=\mu-k_{\bN\bN}\tr\chi+2\tir^{-1}(\Lb \log \tir-\f12\tr\chib)-\f12(\tr\chi-\frac{2}{\tir})\tr\chib+k_{\bN\bN}\tr\chi.
\end{align*}
This gives (\ref{8.11.9.22}). (\ref{8.13.5.22}) follows similarly.

Using (\ref{7.03.1.19}), we have (\ref{4.9.1.24}).
(\ref{8.13.3.22}) then follows by using the definition of $\mho$.
  
(\ref{1.27.6.24}) follows by using (\ref{lb}) and (\ref{cmu_2}).
 (\ref{8.13.4.22}) can be obtained by using (\ref{4.17.1.24}), (\ref{dze}) and  (\ref{mub}).
\end{proof} 

\begin{lemma}\label{10.8.1.23}
 It holds that
\begin{align}
(L+ \tr\chi)\Lb \sX&=(\f12 \sX+\Xi_4+k_{\bN\bN}-\f12\Xi_4)\Lb\Xi_4+2(\zb-\zeta)\sn\sX- \delta L \sX\nn\\
&+ \Lb(\f12 \Xi_L-k_{\bN\bN})\tr\chi+\Lb\Big(-|\chih|^2+\N(\Phi, \bp\Phi)\Big).\label{8.31.4.19}
\end{align}
\end{lemma}
\begin{proof}
Due to (\ref{3.8.3.21}) and (\ref{7.04.8.19}), we write
\begin{equation}
L(\tr\chi+\Xi_L)+\f12 \tr\chi(\tr\chi+\Xi_L)=\f12 \Xi_4^2+(\f12 \Xi_L-k_{\bN\bN})\tr\chi-|\chih|^2+\N(\Phi, \bp\Phi),\label{8.31.6.19-}
\end{equation}
where, on the right-hand side, a term of $\Xi_4\c (\f12 \Xi_L-k_{\bN\bN})$ is included in $\N(\Phi, \bp\Phi)$. 
It follows by using (\ref{8.31.5.19}) and (\ref{8.31.6.19-}) that
\begin{align*}
L\big(\bb\Lb \sX\big)&=[L, \bb \Lb]\sX+\bb \Lb L \sX\\
&=2\bb(\zb-\zeta)\sn \sX+\bb \Lb\big(-\f12 \sX\tr\chi+(\f12\Xi_L-k_{\bN\bN})\tr\chi-|\chih|^2+\N(\Phi, \bp\Phi)\big)\\
&-\bb \delta L \sX +\bb\Lb\Xi_4 \c\Xi_4.
\end{align*}
This gives
\begin{align*}
L(\bb\Lb \sX)+\f12\tr\chi (\bb \Lb \sX)&=2\bb (\zb-\zeta)\sn \sX-\f12 \bb \sX\Lb \tr\chi-\bb \delta L \sX\\
&+\bb \Lb\big((\f12 \Xi_L-k_{\bN\bN})\tr\chi\big)+\bb \Lb\big(-|\chih|^2+\N(\Phi, \bp\Phi)\big)+\bb\Lb\Xi_4 \c\Xi_4.
\end{align*}
Hence
\begin{align*}
L(\bb \Lb \sX)+(\tr\chi+\f12\Xi_4)(\bb \Lb\sX)&=\bb(\f12 \sX+\Xi_4) \Lb \Xi_4+2\bb(\zb-\zeta)\sn\sX-\bb \delta L \sX\\
&+\bb \Lb\big((\f12 \Xi_L-k_{\bN\bN})\tr\chi\big)+\bb \Lb\big(-|\chih|^2+\N(\Phi, \bp\Phi)\big).
\end{align*}
Applying (\ref{lb}) followed with simplification, we can obtain (\ref{8.31.4.19}) from the above identity.
\end{proof}

\begin{lemma}\label{12.21.3.23}
\begin{align}\label{12.21.1.23}
\begin{split}
L\big(\bb\Lb(\tir\Lb \varrho)\big)-&\bb \Lb(\tir k_{\bN\bN}\Lb\varrho)=-\bb\Lb(\tir \bA_b \Lb \varrho)+2\bb(\zb-\ze) \tir\sn\Lb \varrho+\bb\Lb(\tir \sD\varrho)-\bb\Lb(\tir\Box_\bg \varrho)\\
&\qquad\qquad\qquad-\bb\Lb(\tir \hb L \varrho)+2\bb\Lb (\tir\zb\c \sn \varrho)-\bb k_{\bN\bN}L(\tir\Lb \varrho),
\end{split}
\end{align}
\begin{align}\label{12.21.2.23}
\begin{split}
[\sn_\Lb, \sD]\varrho&=-2\chib_{AC}\sn_C \sn_A \varrho-\sdiv\chib\sn\varrho-2\sdiv(\ze+\zb)\bN\varrho-2(\ze+\zb)(\sn_A\bN \varrho+\sn_\bN\sn \varrho)\\
&+\big(\ud \bA(\chi+\chib)+\bR_{AC\Lb A}\big)\sn_C \varrho.
\end{split}
\end{align}
\end{lemma}
Indeed, (\ref{12.21.1.23}) follows by using (\ref{8.31.5.19}) and (\ref{6.30.2.19}). (\ref{12.21.2.23}) can be derived by using (\ref{7.03.1.19}).

\section{Rotation vector fields}\label{rotation}
The set of rotation vector fields are one of the most important vector-fields in our paper.
 It was constructed in \cite{shock_demetrios}. One can find an extensive set of geometric calculations for  rotation vector-fields under the set-up of \cite{shock_demetrios}. Note that the spacetime metric in \cite{shock_demetrios} is conformal to ours. In Section \ref{append}, for completeness, we give the proofs for the main properties, Proposition \ref{2.19.4.22} and Proposition \ref{3.22.6.21}. 

Recall from (\ref{5.14.2.23}) the definition of the rotation vector field
\begin{equation*}
{{}\rp{a}\Omega}^\nu={{}\rp{a}\O}^\mu\Pi_\mu^\nu, \quad a=1,2,3.
\end{equation*}
Let $\hat \bN=c^{-1}\bN$. We denote $\la\rp{a}=c\l {}\rp{a}\O, \bN\r$. This implies
\begin{equation}\label{3.19.1.21}
\Omega\rp{a}+\la\rp{a}\hat \bN={}\rp{a}\O.
\end{equation}
Here $\la\rp{a}$ is a scalar function for a fixed $a$. We can identify $\la\rp{a}=\la^a$. $\la$ then is a $\Sigma$-tangent vector field.
\begin{proposition}\label{2.19.4.22}
Let $\ckc=\cb(u,0)$, $\tir=t+(\ckc)^{-1}u$, and $y^k= \bN^k-\frac{x^k}{\tir}$. There hold
\begin{equation}\label{3.18.2.21}
c\la\rp{a}=\ud\ep_{alk} x^l y^k.
\end{equation}
\begin{align}
\hN\la\rp{a}&=-\Omega\rp{a}\log (c \bb),\label{3.18.1.21}\\
\bT\la\rp{a}&=c{}\rp{a}\Omega \log \bb, \label{3.22.4.21}\\
L\la\rp{a}&=-c {}\rp{a}\Omega \log c.\label{3.22.5.21}
\end{align}
\begin{equation}\label{2.10.2.22}
c\hat e_A(\la\rp{a})=(c\thetac_A^k-\frac{1}{\tir} \delta_A^k){}\rp{a}\Omega_k-\ud\ep_{akj}y^k\hat e_A^j
\end{equation}
where $\hat e_A=c^{-1}e_A$.
\end{proposition}

Let ${}\rp{X}\pi$ denote the deformation tensor of $X$ in the spacetime. Next we give the deformation tensor of $\Omega\rp{a}$, which is denoted by ${}\rp{a}\pi$ for short.
\begin{proposition}\label{3.22.6.21}
Let ${}\rp{a}\ss:=\tr{}\rp{a}\sl{\pi}$. There hold
\begin{align*}
{}\rp{a}\pi_{LL}&=0, \quad {}\rp{a}\pi_{L\Lb}=-2(\zeta+\zb)^A {}\rp{a}\Omega_A=-\f12 {}\rp{a}\pi_{\Lb \Lb},\\
{}\rp{a}\pi_{AB}&=-2c^{-1}\la \hat \theta_{AB}+\f12 {}\rp{a}\ss\ga_{AB}\\
&=-2c^{-1}\la\rp{a} \hat\theta_{AB}-2({}\rp{a}\Omega \log c+c^{-1}\la\rp{a} L\log c) \ga_{AB}\\
&+(c^{-1} \la\rp{a} \p_C v_D \delta^{CD}-c^{-1}\la\rp{a} \tr\chi)\ga_{AB},\\
{}\rp{a}\pi_{\bN A}&=c^{-1}\la\rp{a}\sn_A\log \bb -\sn_A(c^{-1}\la\rp{a}),\\
{}\rp{a}\pi_{\bT A}&=-c^{-1}\la\rp{a} \sn\log \bb+v^j \tensor{\ud\ep}{^a_j_l} {e_A}^l c^{-2}-{}\rp{a}O(v^l) {e_A}_l 
\end{align*}
where $\O$ is decomposed in (\ref{3.19.1.21}).
For future reference, we denote ${}\rp{a}\pi^-_{\bT A}=v^j \tensor{\ud\ep}{^a_j_l} {e_A}^l c^{-2}$.

\begin{align}
 {}\rp{a}\ss&=-4({}\rp{a}\Omega \log c+c^{-1}\la\rp{a} L\log c)+2(c^{-1}\la \p_c v_D\delta^{CD}-c^{-1}\la \tr\chi)\label{5.02.2.21}\\
\sn_\Lb {}\rp{a}\Omega_A&=\chib_{AB}{}\rp{a}\Omega_B+{}\rp{a}\pi_{A\Lb}\label{5.6.01.21}\\
\sn_L{}\rp{a}\Omega_A&=\chi_{AB}{}\rp{a}\Omega_B+{}\rp{a}\pi_{A L}\label{5.6.02.21}\\
(\sn_k {}\rp{i}\Omega)_n&=\Pi_{kk'}\tensor{\ud\ep}{^{ik'm}}\Pi_{mn}-\la\rp{i}\thetac_{kn}\nn\\
&-{}\rp{i}\Omega\log c \Pi_{kn}-\sn_k \log c {}\rp{i}\Omega_{n}+\sn_{n}\log c{}\rp{i}\Omega_k
\label{5.6.03.21}\\
\displaybreak[0]
\snc_l {}\rp{a}\Omega^k&=\Pi_l^n\Pi_k^m \ud\ep_{anm}-\la\rp{a}\thetac_{kl}\label{2.10.1.22}\\
\sDc {}\rp{a}\Omega_i&=-\f12 r^{-1}\tr\thetac {}\rp{a}\Omega_i-\f12 \tr\thetac \tensor{\ud\ep}{^a_m_n}\Pi^n_i {y'}^m+\hat \thetac_{Ai}\hat \bN^m \tensor{\ud\ep}{^a_m_n}e_A^n-\snc^A(\la\rp{a} \thetac_{Ai})\label{6.8.1.22}
\end{align}
where  ${y'}^k=c^{-1}\bN^k-\frac{x^k}{r}$.
\begin{equation}\label{5.13.10.21}
\begin{split}
&[L, {}\rp{a}\Omega]^A={}\rp{a}\pi_L^A, \quad [\Lb, {}\rp{a}\Omega]^A={}\rp{a}\pi_\Lb^A,\\
&[\Lb, {}\rp{a}\Omega]=-2{}\rp{a}\Omega^A(\zeta_A+\zb_A)\bN+{}\rp{a}\pi_\Lb^A e_A, \quad [L, {}\rp{a}\Omega]={}\rp{a}\pi_L^A e_A\\
&[\bb \Lb, {}\rp{a}\Omega]=-2\bb {}\rp{a}\Omega^A \zb_A \bN-\bb{}\rp{a}\Omega^A \ze_A L+\bb {}\rp{a}\pi_\Lb^A e_A. 
\end{split}
\end{equation}
\begin{equation}\label{2.10.3.22}
\begin{split}
\snc_B \snc_A(\la\rp{a})&=\snc_B (\thetac_A^k) {}\rp{a}\Omega_k+(\thetac_A^k-(c\tir)^{-1}\delta_A^k)(\Pi_B^n\Pi_k^m \ud\ep_{anm}-\la\rp{a} \thetac_{Bk})\\
&+(\thetac_B^k-(c\tir)^{-1}\delta_B^k)\Pi_A^n\Pi_k^m \ud\ep_{anm}
\end{split}
\end{equation}
 \begin{equation}\label{12.17.1.23}
 \sDc\la\rp{a}={}\rp{a}\Omega\tr\thetac-(\thetac-c^{-1}\tir^{-1}\Pic)\la\rp{a}\thetac.
 \end{equation}
\end{proposition}

\begin{proposition}\label{7.9.1.24}
With ${y'}^k=c^{-1}\bN^k-\frac{x^k}{r}$, there hold
\begin{align}
&\sum_{i=1}^3 {}\rp{i}\Omega^a {}\rp{i}\Omega^b=r^2(\delta_{cd}-{y'}^c {y'}^d) \Pi_c^a\Pi_d^b\label{5.22.1.22}\\
&|\l c^{-1}\bN, \p_r\r_e-1|\les r^{-2}|\la|^2\label{10.16.3.22}\\
&|y'|\approx  r^{-1}|\la|\label{12.20.3.21}\\
&\bN(r-u)=c-\bb^{-1}+c {y'}^k \frac{x^k}{r}\label{12.20.2.21}\\
&L r=v(\hat\bN)+{y'}^k(\frac{c x^k}{r}-v^k)+c.\label{8.17.1.22}
\end{align}
\end{proposition}
\begin{proof}
(\ref{5.22.1.22}) follows by direct checking.

Next we prove
\begin{equation}\label{3.20.2.21}
c^2\l \bN, \p_r\r_g^2+r^{-2}\sum_{a=1}^3(\la^a)^2=1.
\end{equation}
To see it, we decompose
\begin{equation*}
c^{-1}\bN=\l c^{-1}\bN, \p_r\r_e \p_r+r^{-2}\sum_{i=1}^3\l c^{-1}\bN,{}\rp{i}\O\r_e {}\rp{i}\O.
\end{equation*}
Note
\begin{align*}
\l c^{-1}\bN, \bN\r_g&=\l\p_r, c^{-1}\bN\r_e\l\p_r, \bN\r_g+r^{-2}\l c^{-1}\bN, {}\rp{i}\O\r_e\l {}\rp{i}\O, \bN\r_g\\
&=c\l \p_r, \bN\r_g^2+c^{-1}r^{-2}\sum_{i=1}^3 (\la^i)^2,
\end{align*}
which gives (\ref{3.20.2.21}). By using (\ref{3.20.2.21}), (\ref{10.16.3.22}) is proved.

By using (\ref{3.19.1.21}).
\begin{align*}
\l \bN, \p_r\r_g&=c^{-2}\sum_{k=1}^3\bN^k \frac{x^k}{r}=c^{-1}\sum_{k=1}^3({y'}^k+\frac{x^k}{r})\frac{x^k}{r}\\
&=c^{-1}\sum_{k=1}^3{y'}^k \frac{x^k}{r}+c^{-1}.
\end{align*}
Then by direct decomposition, in view of (\ref{3.18.2.21}), we also have
\begin{equation*}
|y'|^2
=|(\l c^{-1}\bN, \p_r\r_\delta-1)\p_r|^2+\frac{1}{r^2}\sum_{i=1}^3 (\la^i)^2.
\end{equation*}
Using the above identity and (\ref{10.16.3.22}), we obtained (\ref{12.20.3.21}).

Next we derive
\begin{equation}\label{12.20.1.21}
\bN r=c {y'}^k \p_k r+c\frac{x^k}{r}\p_k r= c {y'}^k \frac{x^k}{r}+c,
\end{equation}
and hence (\ref{12.20.2.21}) is proved.

Using (\ref{12.20.1.21}), by direct calculation, we have
\begin{align*}
L r&=\bT r+\bN r=r^{-1} v^i x^i+c+c {y'}^k \frac{x^k}{r}\\
&=v(\hat\bN)+{y'}^k(\frac{c x^k}{r}-v^k)+c
\end{align*}
which gives (\ref{8.17.1.22}).
\end{proof}

\section{Bootstrap assumptions and preliminary estimates }\label{BA_decay}
We adopt the metric induced by the null geodesic flow  for the spacetime,
\begin{equation*}
-2\bb d u dt+\bb^2 du^2+\ga_{AB}(d \omega^A+\b^A du)(d \omega^B+\b^B du).
\end{equation*}
Here $\bb \bN=\p_u-\b_A \frac{\p}{\p\omega^A}$ with the shift $\b(0)=0$. With the help of (\ref{3.19.2}) we can obtain
\begin{equation}\label{1.22.4.22}
[L, \b]=-\bb(\zb-\zeta).
\end{equation}
  
 For a scalar function $f$,
\begin{equation}\label{12.10.1.23}
\p_u\int_{S_{t,u}} f d\mu_\ga=\int_{S_{t,u}} \bb (\bN+\tr\theta) f d\mu_\ga
\end{equation}
and 
\begin{equation*}
\p_u f=\bb\bN f-\b(f).
\end{equation*}
 If $f$ vanishes as $u=u_*$,  there holds
\begin{equation}\label{3.28.2.21}
f(u,t,\omega)= -\int^{u_*}_u \{\bb \bN f-\b (f) \}du'.
\end{equation}

For an $S_{t,u}$-tangent tensor or scalar $F$,  define the $L^p$ norm on $S_{t,u}$ by 
\begin{equation*}
\|F\|^p_{L^p(S_{t,u})}=\int_{S_{t,u}} |F|^p d\mu_\ga; \quad \|F\|^p_{L_\omega^p}=\int_{S_{t,u}}|F|^p(\omega) d\mu_{\mathbb S^2}
\end{equation*}
where  $|\cdot|$ is  taken with respect to $\ga$ if $F$ is a tensor.  $d\mu_{\mathbb S^2}$ in the latter is often written as $d\omega$ for short.

We  introduce the following convention:
For a tensor or scalar field $F$, which can be decomposed as
\begin{equation*}
F=F_1+F_2
\end{equation*}
with $F_1$  and $F_2$ verifying the estimates
 \begin{equation}\label{7.8.1.24}
 |F_1|\les f_1(\Delta_0)\l t\r^{\a_1}, \quad  \|F_2\|_{\sN}\les f_2(\Delta_0)\l t\r^{\a_2},
\end{equation}
where $\sN$ denotes some norm, $f_1$ and $f_2$ are some functions of $\Delta_0$, we write
\begin{equation*}
F=O(f_1(\Delta_0)\l t\r^{\a_1})+O(f_2(\Delta_0)\l t\r^{\a_2})_\sN \mbox{ or } F-O(f_1(\Delta_0)\l t\r^{\a_1})=O(f_2(\Delta_0)\l t\r^{\a_2})_\sN.
\end{equation*}
We may also write $O(f_1(\Delta_0)\l t\r^{\a_1})$ as $f_1(\Delta_0)O(\l t\r^{\a_1})$.

We may abuse the symbolic notation when no confusion occurs: with
 \begin{equation*}
 F=F_1+O(F_2)_\sN \mbox{ or } F=F_1+O(\|F_2\|_\sN)_{\sN} 
 \end{equation*}
 we mean $F=F_1+G$ with $\|G\|_{\sN}\les \|F_2\|_{\sN}$.  
Moreover, we may keep $F_1$ as is, and symbolically write $F_2$ which is bounded as in (\ref{7.8.1.24}), i.e.
\begin{equation*}
F=F_1+O(f_2(\Delta_0)\l t\r^{\a_2})_\sN.
\end{equation*}
\begin{definition}\label{3.24.2.24}
For scalar functions or $3$-vectorfields $f$, we denote the weighted energy for $f$ on $\Sigma_t\cap \{u_1\le u\le u_*\}$ by  $W_{m, -\a}[f](t, [u_1, u_*])$  with $m=1,2$ and $\a\ge 0$
\begin{align*}
&W_{m,-\a}[f](t, [u_1, u_*])=\int_{u_1}^{u_*} \int_{S_{t, u}} \aaa^{-\a}\bb\tir^m \left(|L(v_t^{\f12}f )|^2+|\sn f|^2v_{t}\right) d\omega du, 
\end{align*}
where on each $S_{t,u}$, $v_t:=\sqrt{|\ga|/|\zga|}$ and $\aaa=\log (\tir+3)$. 

We also define the standard energy for $f$
\begin{equation*}
E_{-\a}[f](t, [u_1, u_*])=\f12\int_{\Sigma_t\cap \{u_1\le u\le u_*\}}\aaa^{-\a}(|\nab_g f|^2+|\bT f|^2)d\mu_g.
\end{equation*}
We will hide the interval $[u_1, u_*]$ if there occurs no confusion. 
We denote by $\H_u^t$  the null cone $\H_u$ with temporal parameter $t'\in[0,t]$ and let $\D_{u_1}^{t_1}=\{0\le t\le t_1, u_1\le u\le u_*\}$.  

We define the weighted energy fluxes for $f$ by 
\begin{align*}
&F_{m, -\a}[f](\H_u^t)=\int_{0}^{t} \int_{S_{t', u}} \tir^m \aaa^{-\a}|L(v_{t'}^{\f12}f )|^2 d\omega dt',\, m=1,2
\end{align*}
and the standard energy flux 
\begin{equation*}
F_{0,-\a}[f](\H_u^t)=\int_{\H_u^t} \aaa^{-\a}\{|\sn f|^2+|L f|^2\}d\mu_\ga dt'.
\end{equation*}
We introduce the following notations for various combinations of weighted energies:
\begin{align*}
   WL_{m, -\a}[f](t)&=\int_{\Sigma_t}\aaa^{-\a}\tir^{m}\{\tir^{-2}|f|^2+|\sn f|^2+|(L+h)f|^2\}d\mu_g\\
   I_{m, -\a}[f](\D^{t_1}_{u_1})&=\int_{\D^{t_1}_{u_1}} \tir^{m-1}\aaa^{-\a}(2-m+\a\frac{\tir}{\tir+3}\aaa^{-1})(|\sn\psi|^2+|\psi|^2\tir^{-2})d\mu_g dt \\
   WFIL_{m, -\a}[f](\D_{u_1}^{t_1})&=WL_{m, -\a}[f](t_1)+F_{m, -\a}[f](\H_{u_1}^{t_1})+I_{m, -\a}[f](\D^{t_1}_{u_1})\\
   WFL_{m,-\a}[f](\D_{u_1}^{t_1})&=WL_{m, -\a}[f](t_1)+F_{m, -\a}[f](\H_{u_1}^{t_1})
\end{align*} 
where $m=1,2$ and $\a\ge 0$. 
\end{definition}
In application, we may drop the standard volume elements $d\mu_\ga$, $d\mu_g$ and $d\mu_g dt$ if integrating on $S_{t,u}$, $\Sigma$, and the spacetime region $\D_{u_1}^{t_1}$ respectively. If $\a=0$ in the above definitions, we will drop the corresponding subscript of $-\a$. 

\subsection{Initial slice}\label{9.9.1.22}
We give the result on the initial slice $\Sigma_0$.  
\begin{proposition}\label{12.21.1.21}
Assuming (\ref{1.12.1.22}), it holds that
\begin{equation}
|\bp\Phi(0)|\les 1.\label{9.25.1.22}
\end{equation}
At $\Sigma_0$, 
 we have the following geometric properties 
\begin{equation}\label{12.20.4.21}
\la, r-u, \bb^{-1}-c, \hat \theta, \beta, \tr\thetac-\frac{2}{u}, \sta{\circ}K-\frac{1}{u^2}=0, |S_{0,u}|_\gac=4\pi u^2.
\end{equation}
 Assuming (\ref{1.12.1.22})-(\ref{9.22.1.22}), there is a small bound $\La_0=C\ve$, with the constant $C>0$ merely depending on $\A_0$ and $c_*$, such that the following estimates hold at $t=0$ for $X\in \{\Omega, S\}$
\begin{align}
&\|X\Phi\|^2_{L^2(\Sigma_0)}+\sum_{m=1}^3W_2[X^m \Phi](0)\le \La_0^2, \|X^{\le 3}\Phi\|_{L_u^pL_\omega^4}+\|X^{\le 2}\Phi\|_{L_u^p L_\omega^\infty}\le \La_0^{\f12+\frac{1}{p}},\, p=2, \infty\label{5.14.1.23}\\
&\sum_{l=0}^1 \sum_{m=l}^2E[\Omega^{1+\le m-l}\bT^l\varrho](0)\le \La_0^2\label{9.29.4.23}
\end{align}
Moreover, there hold at $\Sigma_0$ with $Y=e_A, L, \Lb$ 
\begin{align}
&| X^{\le 2}\Phi, \bA_b|\le \La_0^\f12, \quad |v_A|\le \La_0\label{9.18.1.23}\\
&\| X^{\le 2}\Phi, \bA_b\|_{L_u^2 L_\omega^\infty}\le \La_0,\, |c_*-c|\le\La_0^\f12,\, \f12 c_*<c<\frac{3}{2}c_*, \,|c-\ckc|\le \La_0 \label{9.30.1.23}\\
&|\sn_Y\bA_{g,1}, \sn^{\le 1}\bA_g, \sn \bA_b, \sn_S^{\le 1}\ze|\le \La_0,\, \|\sn_\Omega^{1+\le 1}\bA_b, (\tir\sn)^{\le 1}\sF, \sn_X^{\le 1}\ze, (1-\vs^-(X^l))\sn_X^l[\Lb\Phi]),\nn\\
& (1-\vs^-(X^l))\sn_X^l[L\Phi], \tir\sn_X^{\le 1}\sn_Y \bA_{g,1}\|_{L_\omega^4}\le \La_0, l=1,2\label{12.22.1.21}\\
\displaybreak[0]
&\left\{\begin{array}{lll}
\|\Sc(X^{1+\le l}\Phi), \Ac(X^{1+\le l}\Phi)\|_{L_\omega^4}+|\Sc(X^{\le l}\Phi), \Ac(X^{\le l}\Phi)|\le\La_0^\f12,\\
 \|\Sc(X^{\le l}\Omega\Phi),\Ac(X^{\le l}\Omega\Phi), \tir\Sc(X^{\le l} L\Phi), \tir\Ac(X^{\le l}L \Phi)\|_{L^2_\Sigma}\le \La_0, l\le 3.
\end{array}\right.\label{10.22.3.23}\\
& \|X^{\le 2}(\tir\sn_\Lb \bA_{g,1}, \tir\sn [\Lb\Phi]), \tir \sn_\Lb \sn_X^{\le 2}\bA_{g,1}, \sn_X^{\le 2}\bA, \sn_X^3 \bAn, \sn_X^{\le 2} \ze, \nn\\
&\qquad\qquad \qquad\qquad\qquad(\tir\sn)^{1+\le 2}(\tr\chi,\ze), \tir(\tir\sn)^{\le 2}\sF\|_{L^2_{\Sigma_0}}\le \La_0\label{1.12.3.22}\\
&\|X^{\le 1}\Box_\bg\Phi, \Box_\bg\Omega\Phi\|_{L^2_{\Sigma_0}}\le\La_0. \label{1.12.4.22}\\
&\|\sta{\Omega^m, \Lb}\varrho, \Omega^m\bN\varrho(0)\|_{L^2_{\Sigma_0}}+\|\Omega^m \fB(0)\|_{L^2_{\Sigma_0}}+|\Omega \bT\varrho(0), \Omega^{1+\le 1}\varrho(0)|\les \ve, m=1,2,3\label{10.1.9.23'}\\
&\|\Omega^{1+\le 1}\fB(0), \Omega^{1+\le 2}\varrho(0)\|_{L_\omega^4}\les \ve\label{3.21.6.24'}\\
&-\ve^\f12\les \varrho(0)\les \ve\label{9.30.15.23}\\
&\|\Lb \varrho(0)\|_{L_u^2 L_\omega^\infty}\les \ve^\frac{1}{4}, \quad q_0\les \ve^\f12\label{12.6.2.23}
\end{align}
where, for a scalar function $f$,  $\sta{\Omega^m, \Lb}f=Y^{m+1}f$ with only one of vector field $Y=\Lb$ and the remaining vector fields are all $\in\{\Omega\}$. 
\begin{remark}\label{11.22.1.23}
Recall $\bb_0=\inf_{\Sigma_0}c^{-1}$. From (\ref{9.30.1.23}), we see that $\bb_0>0$ and is comparable to $c_*$. 
\end{remark}
\end{proposition}
The complete proof of Proposition \ref{12.21.1.21} can be found in Section \ref{9.25.2.22}. It is presented in the later part of the paper mainly because the proof relies on geometric decompositions, comparisons, and commutation formulas that will be derived while controlling the energy propagation. We emphasize that the formulas we apply to the initial slice and the analysis carried out on the initial slice are independent of the temporal evolution of the compressible Euler flow. The results in  Proposition \ref{12.21.1.21} will be used as the initial data for the propagation of energies and geometry, and for setting up the bootstrap argument in the sequel.

\subsection{Bootstrap assumptions}\label{BAsec}
In this section, we make a set of bootstrap assumptions. Due to our small-large regime for amplitude,
it is crucial to determine both the decay rates and appropriate bounds for each quantity involved.
 Table (\ref{5.24.1.23}) provides us with the proper expectation of the $L^\infty_u$ bounds for all the quantities. The decay rate primarily is related to the energy hierarchy and the geometric structure of the quantity.  We make assumptions on the quantities that have smallness property, and  then use the assumptions to prove that all the quantities with vanishing $\al$-value are bounded with certain decay for large $t$, while improving our assumptions.

\begin{assumption}\label{5.13.11.21+}
Let $X\in\{\Omega, S\}$, and we fix the convention that $\l t\r=t+2$. Let $0<\delta<1$ be fixed. With  $0<\Delta_0<1$ a small constant comparable to $\La_0$, to be chosen,
assume that $T_*>0$ is the maximal life span such that there hold for $0<t<T_*$,
 \begin{equation}\label{3.12.1.21}\tag{\bf BA}\left\{
 \begin{array}{lll}
 \|X\Phi\|^2_{L^2_\Sigma}\le \Delta_0^2\l t\r^{2\delta}\\
 E[\Omega^m\bT\varrho](t)+F_0[\Omega^m\bT\varrho](\H_u^t)\le \Delta_0^2\l t\r^{2\delta},\,  m=1,2,\\
  E[\Omega^m\varrho](t)+F_0[\Omega^m\varrho](\H_u^t)\le \Delta_0^2\l t\r^{\max(0, m-2)2\delta}, \, 1\le m\le 3,\\
W_2[X^m\Phi](t)+F_2[X^m\Phi](\H_u^t)\le \Delta_0^2 \l t\r^{2\delta}, \, 1\le m\le 2,\\
W_1[X^3\Phi](t)+F_1[X^3\Phi](\H_u^t)\le \Delta_0^2 \l t\r^{2\delta},
\end{array}\right.
 \end{equation}
  where $L^2_\Sigma$ is the short-hand notation for $L^2$-norm on $\Sigma_t=\{t'=t,u\in [u_0, u_*]\}$.
  
  Moreover, for ease of analysis, we give the auxiliary bootstrap assumptions, with $X=\tir \sn, \sn_S$
\begin{equation}\label{L4BA1}\tag{L4BA1}
\begin{array}{lll}
&\|(\tir\sn)[L\Phi], X^{\le 1}\bA_{g,1}\|_{L_\omega^4}+\Delta_0^\f12 \|S^{\le 1}([L\Phi]), \tir^{-1}\varrho\|_{L^4_\omega}\le \Delta_0 \l t\r^{-2+\delta},\\
&\|X(\tir \sn)[L\Phi], \tir \sn X[L\Phi], \sn_X^2\bA_{g,1}\|_{L_\omega^4}+\Delta_0^\f12\|S^2([L\Phi])\|_{L^4_\omega}\le \Delta_0 \l t\r^{-\frac{7}{4}+\delta}\\
& \|\tir X^{\le 1}\sn_\Lb \bA_{g,1}, X^{\le 1}(\tir\sn)[\Lb\Phi]\|_{L_\omega^4}\le \Delta_0 \l t\r^{-1+\delta}
\end{array}
\end{equation}
\begin{align}
&
\| \l t\r^{-\frac{l}{2}}X^{l+\le 1} Y (\bA^\natural), \l t\r^{-\frac{l}{2}-1} X^{l+\le 2}Y\Phi\|_{L^2_\Sigma}\le \Delta_0 \l t\r^{-2+\delta}, Y=\sn_L, \sn, l=0,1 \label{L2BA2}\tag{L2BA2}\\
&\|X^{\le 2}(\sn_\Lb \bA_{g,1}, \sn[\Lb\Phi])\|_{L^2_\Sigma}\le \Delta_0\l t\r^{-1+\delta}\label{LbBA2}\tag{LbBA2}\\
&\|\tir^\frac{3}{4} \sn(\bA_{g,2}, \bA_b), \bA_{g,2}\|_{L_\omega^4}+\Delta_0^\f12\|\bA_b\|_{L_\omega^4}\le\Delta_0 \l t\r^{-2+\delta}\label{L4conn}\tag{L4conn}\\
\displaybreak[0]
&\|\bb^{-\f12}\Big(\bA_b, \bA_{g,2}, \l t \r^{-\f12(m-1)(1-\vs^+(X^m))}X^{m}(\bA_b, \bA_{g,2})\Big)\|_{L^2_\Sigma}\nn\\
&\qquad\qquad\qquad\qquad\qquad\qquad\qquad\le \Delta_0\l t\r^{-1+\delta}, \,m=1, 2\label{L2conndrv}\tag{L2conndrv}\\
&\|\bb^{-\f12}\tir^3 \sn^3\tr\chi\|_{L^2_\Sigma}\le \l t\r^\delta\Delta_0, \|\tir^2\sn^2\tr\chi\|_{L_\omega^4}\le \l t\r^{-1+\delta}\Delta_0 \label{ConnH}\tag{ConnH}\\
 &|\tir \ze|\le\l t\r^{\delta}\Delta_0, \|\sn \ze\|_{L_\omega^4}\le \l t\r^{-2+\delta}\Delta_0, \|(\tir\sn)^{\le 1}\sn \ze\|_{L^2_u L_\omega^2}\le \Delta_0\l t\r^{-2+\delta}\label{zeh}\tag{ZetaH}\\
&\|\Box_\bg \Omega\Phi\|_{L^2_\Sigma}\le \Delta_0\l t\r^{-2+2\delta}\label{wave_ass}
\end{align}
\begin{footnote}
{The more important part of the bootstrap assumptions is for $t>\kappa$ with $\kappa>0$ a large fixed constant. For $0<t\le \kappa$, we can actually drop the various power of $\l t\r$ in the assumptions. In this case, by letting $\ve_0>0$ be sufficiently small and running the same energy arguments, we can improve $\Delta_0$ on the right-hand to be a constant multiple of $\La_0+\Delta_0^\frac{5}{4}$, with $\Delta_0\approx \La_0$ properly chosen. }   
\end{footnote}
 as well as the following assumptions  
 \begin{align}
&-\varrho<-\tir\Lb \varrho\le \M_0 \big(1+\f12\wp\log(\frac{\l t\r}{2})\big)^{-1} \label{6.5.1.21}\\
& |y'|\le \Delta_0^\f12 \l t\r^{2\delta-1}, \quad |\sn \la|\le \Delta_0 \label{1.25.1.22}\\
& |\tr\chi-\frac{2}{\tir}, \chih|\le \Delta_0^\f12\l t\r^{2\delta-\frac{7}{4}}\label{3.7.1.21-}\\
&\frac{\bb_0}{20}\le \bb\label{6.20.2.21}
\end{align}
\end{assumption}
where $y'$ is defined in Proposition \ref{7.9.1.24}; in (\ref{6.5.1.21}) $\M_0>1$ is a constant $\approx \A_0$, to be specified; and $\bb_0$ has been defined in Remark \ref{11.22.1.23}. \begin{footnote}{ We can assume  $\Delta_0^\f12 \M_0\les 1$, which is achievable with $\ve_0$ sufficiently small. This fact will be used without explicit mentioning.  }\end{footnote}

If $T_*<\infty$, we will be able to improve both inequalities in (\ref{6.5.1.21}), and improve the right hand-sides of all other inequalities to $<$, which gives the contradiction to the definition of $T_*$. For $t<T_*$ and $X\in \{\Omega, S\}$, (\ref{3.12.1.21}) will be improved to
\begin{align*}
 &\|X\Phi\|^2_{L^2_\Sigma}\les(\La_0^2+\Delta_0^\frac{5}{2}) \log \l t\r^{\M}\\
 &E[\Omega^m\varrho](t)+F_0[\Omega^m\varrho](\H_u^t)\les \La_0^2+\Delta_0^\frac{5}{2}, m=1,2\\
&E[\Omega^3\varrho](t)+F_0[\Omega^3\varrho](\H_u^t)\les(\La_0^2+\Delta_0^\frac{5}{2})\log \l t\r^{2\M^2-\M}\\
 &E[\Omega^m\bT\varrho](t)+F_0[\Omega^m\bT\varrho](\H_u^t)\les \log \l t\r^{\M^2}(\La_0^2+\Delta_0^\frac{5}{2}),  m=1,2\\
&W_2[X^m\Phi](t)+F_2[X^m\Phi](\H_u^t)\les\log \l t\r^{2\M^2}(\La_0^2+\Delta_0^\frac{5}{2}), m=1, 2\\
&W_1[X^3\Phi](t)+F_1[X^3\Phi](\H_u^t)\les\log \l t\r^{2\M^2-\M}(\La_0^2+\Delta_0^\frac{5}{2})
\end{align*}
where the constant $\M$ is a constant multiple of $\M_0$ and $\M \ge 15$. We will show that the right-hand side of (\ref{L4BA1})-(\ref{zeh}) hold with the factor $\l t\r^{\delta}$ on the right-hand side replaced by $(\log \l t\r)^{\M^2+1}$ and $\Delta_0$ replaced by $\La_0+\Delta_0^\frac{5}{4}$, and  (\ref{wave_ass}) will be improved to (\ref{wave_ass'}). The details will be presented in Section \ref{Improve}. We will also improve (\ref{6.5.1.21})-(\ref{6.20.2.21}) to
\begin{align}
&-\varrho+C\ve \l t\r^{-\frac{3}{4}+\delta} <-\tir\Lb \varrho<\ti C(\mathfrak{C}+1)\left(1+\f12\wp \log (\f12\l t\r)\right)^{-1} \label{6.5.1.21_+}\\
&|y'|\les  (\La_0+\Delta_0^\frac{5}{4})\l t\r^{-1+\delta},\quad | \sn\la|\les (\La_0+\Delta_0^\frac{5}{4})\l t\r^{-\frac{3}{4}+\delta}\label{9.14.1.24}\\
\displaybreak[0]
& |\bA_b, \chih|\les (\La_0^\f12+\Delta_0^\frac{5}{4}) \log\l t\r^{\M^2}\l t\r^{-\frac{7}{4}}\label{3.7.1.21}\\
&\bb> \frac{1}{4} \bb_0\label{6.20.2.21+}
\end{align}
where $C, \ti C>0$ are universal constants, and the constant $\mathfrak{C}\approx \A_0$. Hence with $\M_0=2\ti C(\mathfrak{C}+1)$, (\ref{6.5.1.21}) is improved. 

\subsection{Basic decay properties}
We first give basic analytic properties by using (\ref{L4conn}) and (\ref{3.7.1.21-}).
\begin{itemize}
\item{Volume comparison:}
On $S_{t,u}$, with $v_t:=\sqrt{|\ga|/|\zga|}$ and $\zga$ is the canonical round metric on $\mathbb S^2$, 
\begin{equation}
 v_t\approx \tir^2,\quad r\les \tir \approx \l t\r\label{9.29.3.22}
\end{equation}
which will be used frequently without mentioning explicitly. We shall improve the second comparison formula to (\ref{comp1}) shortly. 

\item{Comparison with the standard round metric:} Relative to the transport coordinate $\omega=(\omega_1, \omega_2)$ on each $S_{t,u}$  \begin{equation}\label{11.28.1.23}
|\tir^{-2}\ga_{ab}-\gaz_{ab}|\les\Delta_0^\f12,\quad \|\p_\omega(\tir^{-2}\ga_{ab}-\gaz_{ab})\|_{L_\omega^4}\les \Delta_0,\quad a, b=1,2.
\end{equation}
\item{Sobolev embedding on spheres:} There hold for $S_{t,u}$ tangent tensor $F$ or scalar functions
\begin{equation*}
\|F\|_{L_\omega^\infty}\les \|(\tir \sn)^{\le 1} F\|_{L_\omega^4};\quad \|F\|_{L_\omega^p}\les \|(\tir \sn)^{\le 1} F\|_{L_\omega^2},\quad 2\le p<\infty.
\end{equation*}
\item{\Poincare inequality on spheres:} There holds for scalar functions $f$ 
\begin{equation*}
\|\Osc(f)\|_{L^p_\omega}\les \|\tir\snc f\|_{L^p_\omega}, \quad 1<p<\infty,
\end{equation*}
where $\Osc(f)=f-\bar f$.
\item
 For $F$ vanishing at $u=u_*$ we obtain
\begin{equation}\label{10.22.2.22}
\sup_{[u,u_*]}\|F\|_{L_\omega^p}\les \|\bb \bN F\|_{L_u^1[u, u_*] L_\omega^p},\quad 2\le p<\infty.
\end{equation}
\end{itemize}
Indeed, using the assumption in (\ref{3.7.1.21-}) and (\ref{12.20.4.21}), the first estimate in (\ref{9.29.3.22}) follows by using $L v_t=\tr\chi v_t$.
Using the definition of $\tir$ in Section \ref{10.25.5.23}, we have $\tir\approx \l t\r$ due to (\ref{9.30.1.23}). The second estimate in (\ref{9.29.3.22}) then follows by using $r(t,u)\le r(t,u_*)=t+u_*$. 

Recall the construction of the acoustical coordinates by the geodesic flow in (\ref{6.29.2.19}). Relative to the transport coordinates $t, \omega_1, \omega_2$, we can verify that
\begin{equation}\label{11.28.2.23}
\frac{d\ga_{ab}}{dt}=2\chi_{ab},\quad a, b=1,2
\end{equation}
along null geodesics $\Upsilon_{\omega, u}(t)$.
 The first estimate in (\ref{11.28.1.23}) can be obtained by using (\ref{3.7.1.21-}); the second one follows by using (\ref{3.7.1.21-}) and the bounds on $\bA_b, \bA_{g,2}$ in (\ref{L4conn}). (See details in \cite[Proposition 5.3]{Wangrough}). The remaining two sets of inequalities can be obtained by (\ref{11.28.1.23}).

Finally we prove (\ref{10.22.2.22}). 
Note for an $S_{t,u}$-tangent tensor field $F$, applying (\ref{12.10.1.23}) to $|F|^p v_t^{-1}$ leads to
\begin{equation}\label{11.12.5.23}
\p_u (\int_{S_{t,u}}|F|^p v_t^{-1})= \int_{S_{t,u}}\bb(\bN+\tr\theta)(|F|^p v_t^{-1}), \quad p\ge 2.
\end{equation}
In view of $\bN v_t=\tr\theta v_t$, we obtain $\p_u (\|F\|^p_{L_\omega^p})=\int_{S_{t,u}}\bb p\bN F \c F|F|^{p-2} d\mu_{\mathbb S^2}$. Integrating it in $[u,u_*]$ gives (\ref{10.22.2.22}).
 $\hfill\square$

We also recall the result from \cite[Section 5]{Wangrough}.
\begin{lemma}[Transport lemma]\label{tsp2}
For any $S_{t,u}$-tangent tensor field $F$ satisfying
\begin{equation*}
\sn_L F+m {\emph\tr}\chi F= W
\end{equation*}
with a constant $m$, there holds
\begin{equation*}
v_t^m F(t)=v_{0}^m F(0)+\int_0^t v_{t'}^m Wdt'.
\end{equation*}
Similarly, for the transport equation
\begin{equation*}
\sn_L F+\frac{2m}{\tir} F=G\c F+W
\end{equation*}
with a constant $m$, if $\|G\|_{L_\omega^\infty L_t^1(\H_u)}\le C$, then there holds
\begin{equation*}
\tir^{2m} |F(t)|\les\tir^{2m}(0)|F(0)|+\int_0^t \tir^{2m} |W| d\tt.
\end{equation*}
The same result holds when $\frac{2}{\tir}$ in the transport equation is replaced by ${\emph\tr} \chi$.
The above integrals are taken along null geodesics on $\H_u$.
\end{lemma}

Next we derive basic decay properties.
\begin{lemma}[Decay lemma]\label{5.13.11.21}
Under the assumptions (\ref{3.12.1.21})-(\ref{6.20.2.21}), with $X\in \{\tir\sn, \sn_S\}$, there hold the following decay estimates 
\begin{enumerate}
\item 
\begin{itemize}
\item With $p=2, \infty$ and $l=0,1$, there hold
\begin{align}
&\left\{\begin{array}{lll}
\|\tir^2(X^l L \Phi, X^l\sn \Phi)\|_{L_u^p L_\omega^\infty}\les \Delta_0^{\f12+\frac{1}{p}} \l t\r^{\delta+\frac{l}{4}},\\
 \tir^2 |X^l\bA_{g,1}|\les \Delta_0\l t\r^{\delta+\frac{l}{4}}, \Delta_0^{\f12-\frac{1}{p}}\|\bA_b\|_{L_u^p L_\omega^\infty}+|\bA_{g,2}|\les \l t\r^{-\frac{7}{4}+\delta}\Delta_0,\\ 
 \tir(|\ze|+|\tir (\sn[\Lb\Phi], \sn_\Lb \bA_{g,1})|)+ \tir^2|\sn \tr\chi|\les \Delta_0  \l t\r^\delta
\end{array}\right.
 \label{3.6.2.21}\\
 &\left\{\begin{array}{ll}
\tir^3\Delta_0^\f12 |X^l L[L\Phi]|+\tir^3|X^l(\sn, \sn_L)\bA_{g,1},\tir^{-1}(\tir\sn)^{l}(\bA_b, \bA_{g,2}), X^l\sn[L\Phi]|+\\
 +\tir^3|X^l\sn(\sn\Phi^\dagger, L \Phi^\dagger)|=O(\Delta_0 \l t\r^{\delta+\frac{l}{4}})_{L_\omega^4},\\
  \|\sn \ze, X^{\le 1}\sn_\Lb \bA_{g,1}, X^{\le 1}\sn[\Lb\Phi]\|_{L_\omega^4}\les \l t\r^{-2+\delta}\Delta_0
  \end{array}\right.
    \label{3.11.3.21}\\
&-\l t\r^{-1+\delta}\Delta_0^\f12\les -\tir\Lb \varrho\le \M_0 \big(1+\f12\wp\log(\frac{\l t\r}{2})\big)^{-1} \label{11.11.2.23}\\
&\frac{1}{4}c^{-1}(0)<\bb \les\log \l t\r+1\label{11.13.3.23}\\
&\tir \bb|\Lb \Phi, \bT \Phi|+\tir\bb |k_{\bN\bN}|\les 1\label{3.11.4.21}\\
&\|\bb\tir\Lb \varrho\|_{L^2_u L_t^\infty L_\omega^\infty}\les \Delta_0^\frac{1}{4}\label{12.19.1.23}
\displaybreak[0]
\end{align}
\begin{align}
& |\tir^{-1}\b|\les \Delta_0\label{3.28.3.21}
\end{align}
\begin{align}&\left\{
\begin{array}{lll}
|v_A|\les \Delta_0 \l t\r^{-1+\delta}, \,\, |v, \varrho|\les \Delta_0^\f12 \l t\r^{-1+\delta}\\
|c-\ckc|+|c-c_*|\les \Delta_0^\f12.
\end{array}\right.\label{6.24.1.21}\\
 &\|v\|_{L^2_u L_\omega^2}\les\sum_{X=S,\Omega}\|X\Phi\|_{L^2_u L_\omega^2}, \quad\|v_A\|_{L_\omega^2}\les \l t\r^{-2+2\delta}\sum_{X=\Omega, S}\|X\Phi\|_{L^2_\Sigma}.\label{9.18.3.23}
\end{align}

\item There hold the following derivative estimates for $\fB$ or $\Lb\Phi$
\begin{equation}\label{6.22.1.21}
\tir(\Sc([L\Lb \Phi]), \Sc([\Lb L\Phi]), L\Lb\Phi, \Lb L\Phi)=O(1)(\fB). 
\end{equation}
\begin{align}
&X(\Sc([L\Lb\Phi]), \Sc([\Lb L\Phi]))\nn\\
&=\left\{\begin{array}{ll}
\vs(X)(\l t\r^{-1}[\Lb\Phi]+O(\l t\r^{-2+\delta}\Delta_0)_{L^2_\Sigma})+(1-\vs(X))O(\l t\r^{-1+\delta}\Delta_0)_{L^2_\Sigma}\\
\vs(X)(\l t\r^{-1}[\Lb\Phi]+O(\l t\r^{-\frac{11}{4}+\delta}\Delta_0^\f12)_{L_\omega^4})+(1-\vs(X))O(\l t\r^{-2+\delta}\Delta_0)_{L_\omega^4}.
\end{array}
\label{6.7.4.23}\right.
\end{align}
\begin{equation}\label{8.23.1.23}\left\{
\begin{array}{lll}
X\fB-\vs(X)[\Lb \Phi]=(1-\vs(X))O(\Delta_0\l t\r^{-1+\delta}), (1-\vs(X))O(\Delta_0\l t\r^\delta)_{L^2_\Sigma}\\
X\Lb \Phi,\Lb X\Phi
=\left\{\begin{array}{lll}
&[\Lb \Phi]+(1-\vs(X))\Delta_0\l t\r^{-1+\delta},\\
 &[\Lb\Phi]+(1-\vs(X))O(\Delta_0\l t\r^\delta)_{L^2_\Sigma},
\end{array}\right.\\
X^l\fB=\left\{\begin{array}{lll}
\vs^-(X^l)([\Lb \Phi]+S^{l-1}[L\Phi]+O(\l t\r^{-\frac{7}{4}+\delta}\Delta_0)_{L_\omega^4})+(1-\vs^-(X^l))O(\Delta_0\l t\r^{-1+\delta})_{L_\omega^4}\\
\vs^-(X^l)([\Lb\Phi]+ O(\l t\r^{-1+\delta}\Delta_0)_{L_\Sigma^2} )+(1-\vs^-(X^l))O(\Delta_0\l t\r^\delta)_{L^2_\Sigma}, l=2
\end{array}\right.
\end{array}\right.
\end{equation}
where we assume $\Lb \tir=O(1)$ for the estimate of $\sn_\Lb X\Phi$. 
\begin{equation}\label{8.23.1.23'}
X^l\Lb \Phi=\left\{\begin{array}{lll}
&[\Lb \Phi]+(1-\vs^-(X^l))O(\Delta_0\l t\r^{-1+\delta})_{L_\omega^4}+\vs^-(X^l)O(\l t\r^{-\frac{7}{4}+\delta}\Delta_0^\f12)_{L_\omega^4};\\
& [\Lb\Phi]+\vs^-(X^l) O(\l t\r^{-1+\delta}\Delta_0)_{L^2_\Sigma}+(1-\vs^-(X^l))O(\Delta_0 \l t\r^\delta)_{L^2_\Sigma}, l=2
\end{array}\right.  
\end{equation}
\end{itemize}
\item
\begin{itemize}
\item
For $\wt\N(\Phi,\bp\Phi)$ defined in Proposition \ref{geonul_5.23_23}, there hold
\begin{align}
&|\bb X^l\wt\N(\Phi, \bp\Phi)|\les \Delta_0^\f12 \tir^{-3+\frac{l}{4}+\delta},\quad\|\bb X^{l}\wt\N(\Phi, \bp \Phi)\|_{L_u^2 L_\omega^\infty}\les \Delta_0\l t\r^{-3+\frac{l}{4}+\delta}, l=0,1\label{4.3.3.21}
\end{align}
\begin{equation}\label{3.29.1.23}\left\{
\begin{array}{lll}
\|X^l\wt\N(\Phi, \bp\Phi)\|_{L_\omega^4}\les \Delta_0^\f12\l t\r^{-3+(2-\max(l-1,0))\delta+\frac{1}{4}\max(l-1, 0)}, l=0, 1,2\\
\|\tir X^{l+1} \wt\N(\Phi, \bp\Phi)\|_{L_u^2 L_\omega^4}\les \l t\r^{-2+\frac{l}{4}+(2-l)\delta}\Delta_0, l=0,1\\
\| X^{\le 2}\wt\N(\Phi, \bp \Phi)\|_{L^2_\Sigma}\les \l t\r^{-2+2\delta}\Delta_0.
\end{array}\right.
\end{equation}
\item There hold for $l=0,1$ and $0\le m\le l$ that
\begin{align}\label{3.16.1.22}
\begin{split}
& \|X^{l-m}\sn_L X^m \bA_b, \tir^{-1}X^{l-m}\sn_LX^m(\tir\tr\chi)\|_{L_\omega^4}\les \l t\r^{-3+\frac{l}{4}+\delta}\Delta_0^{\f12+\f12(1-\vs(X^l))},\\
& \|X^{l-m}\sn_L X^m\bA_{g,2}\|_{L_\omega^4}\les \l t\r^{-3+\frac{l}{4}+\delta}\Delta_0, \Delta_0^\f12\sn_L \bA_b, \sn_L\bA_{g,2}=O(\l t\r^{-\frac{11}{4}+\delta}\Delta_0),\\
 &L \tr\chi, L\tr\chib= O(\l t\r^{-2})\\
&\|X^{l-m}\sn_LX^m(\tir \ze)\|_{L_\omega^4}\les \l t\r^{-1+\delta}\Delta_0,  (\sn_L+\tir^{-1}) \ze=\l t\r^{-2+\delta}\Delta_0.\\
&\|\bb^{-\f12}X^2 \ze\|_{L^2_\Sigma}\les \Delta_0\l t\r^\delta
 \end{split}
\end{align}
\begin{align}
&\Delta_0^\f12 (\l t\r^{1-\frac{l}{4}}\|X^l\widetilde{L(\Xi_4)}\|_{L_\omega^4}+\l t\r^\frac{3}{4}|\widetilde{L(\Xi_4)}|)
+\|\tir^{-\f12 l-(1-l)\delta} X^{l+\le 1}\widetilde{L(\Xi_4)}\|_{L^2_\Sigma}\les\l t\r^{-2+\delta}\Delta_0.\label{6.8.4.23}
\end{align}

\end{itemize}

\item
\begin{align}\label{1.29.4.22}
\begin{split}
&\l t\r^{-\frac{l}{4}}\Big(X^l(\bR_{ABLC},\widehat{\bR_{A4B4}}, \widehat{\bR_{A3B4}}, \ga^{AB}\bR_{A43B}-\varpi)\Big)\\
&\qquad\qquad\qquad\qquad\qquad=\left\{\begin{array}{lll}
O(\l t\r^{-3+\delta}\Delta_0)_{L_\omega^4}, \, l=0,1, \\
O(\l t\r^{-\frac{11}{4}+\delta}\Delta_0),\, l=0
\end{array}\right.\\
 & X^{l}\bR_{AB43}=O(\l t\r^{-3+2\delta}\Delta_0^2)_{L_\omega^4}, l=0,1; O(\l t\r^{-3+2\delta}\Delta_0^2), l=0\\
& X^l \bR_{AB\Lb C}=O(\l t\r^{-2+\delta}\Delta_0), \mbox{ if }  l=0;\, O(\l t\r^{-2+\delta}\Delta_0)_{L_\omega^4}, \mbox{ if } l=0,1. 
\end{split}
\end{align}
\begin{align}
&\begin{array}{lll}
\|X^{l+\le 1}( \bR_{ABCL}, \widehat{\bR_{4A4B}}, \widehat{\bR_{A3B4}}, \ga^{AB}\bR_{A43B}-\varpi)\|_{L^2_\Sigma}\les \l t\r^{-2+\f12 l+\delta}\Delta_0, l=0,1\\
\|X^{\le 2}\bR_{ABC\Lb}\|_{L^2_\Sigma}\les \Delta_0\l t\r^{-1+\delta}, \|X^{\le 2}\bR_{AB43}\|_{L^2_\Sigma}\les \l t\r^{-2+2\delta}\Delta_0^2 
\end{array}\label{10.15.2.22}
\end{align}
where $\varpi=(\bA+\frac{1}{\tir}+\fB)\fB$ symbolically.
\item
\begin{align}\label{10.1.1.22}
\begin{split}
 &|c^2\stc{K}-\tir^{-2}|\les \l t\r^{-\frac{11}{4}+\delta}\Delta_0^\f12, \|\tir\sn \stc{K} \|_{L_\omega^4}\les  \l t\r^{-\frac{11}{4}+\delta}\Delta_0,\\
 &\|X^l(K-\tir^{-2})\|_{L^4_\omega}\les \l t\r^{-3+\frac{l}{4}+\delta}\Delta_0^\f12, l=0,1; |K-\tir^{-2}|\les \l t\r^{-\frac{11}{4}+\delta}\Delta_0^\f12.
 \end{split}
 \end{align}
\item 
 \begin{align}
&\tir^2 \sn(\bb \fB), \tir^2 \bb\sn_\Lb\sn \varrho=O(\Delta_0)_{L_\omega^4} \label{1.27.2.24}\\
&\tir (\tir\sn)^{1+\le 1}(\bb \fB)=O(\Delta_0)_{L_u^2 L_\omega^2}\label{1.27.4.24}\\
&\tir\|\sn \bb\|_{L_\omega^4}\les \log \l t\r\Delta_0\label{1.27.5.24}\\
&\begin{array}{lll}
\|\tir^2\log \l t\r^{-1}\bb \sn^2\log \bb, \tir^2\sn^2\bb, \log \l t\r\tir\sn_L(\tir\sn)^{1+\le 1}\bb\|_{L_u^2 L_\omega^2}\les \log \l t\r\Delta_0,\\ \|\sn\bb\|_{L_u^2 L_\omega^\infty}+\|\tir\sn^2 \bb\|_{L_u^2 L_\omega^4}\les \Delta_0 \l t\r^{-1+\f12\delta}\log \l t\r,\, \|\sn_S(\tir \sn)\bb\|_{L_u^2 L_\omega^4}\les \Delta_0.
\end{array}
\label{2.9.2.24}
\end{align} 
\item
\begin{equation}
\|\bb^3 \Lb^2\varrho\|_{L_\omega^4}\les \l t\r^{-1},\quad \|\bb^3(\Lb^2\Phi,\sn_\Lb(\Lb\Phi^\dagger), \bT\bT\Phi)\|_{L_\omega^4}\les \l t\r^{-1}.\label{3.6.1.21}
\end{equation}
\end{enumerate}
\end{lemma}
\begin{remark}
The left-hand inequality in (\ref{11.13.3.23}) gives the estimate (\ref{6.20.2.21+}). 
(6) is crucial for bounding the general standard energy. The proof is postponed to (\ref{1.29.2.22}).

\end{remark}
\begin{remark}
We have from (\ref{3.6.2.21}) that
\begin{equation}\label{9.4.1.22}
|\bA|\les \l t\r^{-\frac{7}{4}+\delta}\Delta_0^\f12,
\end{equation}
which is a better estimate than (\ref{3.7.1.21-}). We can drop the assumption (\ref{3.7.1.21-}) and use the above instead. 
\end{remark}
\begin{remark}\label{5.11.3.24}
The second estimate in the (\ref{8.23.1.23}) can be further improved by using (5), with $O(\l t\r^\delta\Delta_0)_{L^2_\Sigma}$ replaced by $O(\Delta_0\log \l t\r)_{L^2_\Sigma}$ in the last term.  The last term in  the last line can also by improved by using (5) with with $O(\l t\r^\delta\Delta_0)_{L^2_\Sigma}$ replaced by $O\Big(\Delta_0(\log \l t\r)^2\Big)_{L^2_\Sigma}$. Similarly, the last term in the first line of (\ref{6.7.4.23}) can be improved to $O(\Delta_0\log \l t\r\l t\r^{-1})_{L^2_\Sigma}$. 
\end{remark}
\begin{remark}

Using (\ref{3.6.2.21}), (\ref{3.11.3.21}), (\ref{1.27.2.24}) (\ref{8.23.1.23}) and (\ref{3.16.1.22}), we summarize that, with  $\Phi^0=\varrho,  \Phi^i=v^i$, 
\begin{equation} \label{8.23.2.23}
\begin{split}
&X\Lb\Phi^\mu-(\min(\mu,1)+\vs(X))\fB=O(\l t\r^{-1+\delta}\Delta_0)(1-\vs(X)), O(\l t\r^{-1}\log \l t\r\Delta_0)_{L_\omega^4}(1-\vs(X))\bb^{-1}, \\
 &X\bA=O(\l t\r^{-\frac{7}{4}-\frac{\vs(X)}{4}+\delta}\Delta_0^{\f12+\frac{1}{p}})_{L_u^p L_\omega^4},\, X=\tir \sn, \sn_S, p=2,\infty.
 \end{split}
\end{equation}

\end{remark}
\begin{proof}[Proof of (\ref{3.6.2.21})-(\ref{9.18.3.23})] The $p=2$ case in (\ref{3.6.2.21}) can be obtained by using (\ref{L2BA2}) and (\ref{L2conndrv}) by Sobolev embedding on spheres. The $p=\infty $ estimates in (\ref{3.6.2.21}) follow by using (\ref{L4BA1}), (\ref{L4conn}), (\ref{zeh}), (\ref{ConnH}), (\ref{5.23.1.23}) and Sobolev embedding on spheres. 
We note the estimate for $c$ in (\ref{6.24.1.21}) is obtained by integrating the schematic identity $L\log c=L\varrho$ along null geodesics and using (\ref{3.6.2.21}) and the estimate for $c_*-c$ and $c-\ckc$ at $t=0$ in (\ref{9.30.1.23}). We will frequently use the bound of $|c, c^{-1}|\les 1$ without mentioning it explicitly. Moreover from (\ref{5.14.1.23}) and the bound on $L\Phi$ in (\ref{3.6.2.21}), integrating along null geodesics, we obtain $|\Phi|\les \Delta_0^\f12\l t\r^{-1+\delta}$ in (\ref{6.24.1.21}). (\ref{11.11.2.23}) follows as a consequence of the estimate of $\varrho$ and (\ref{6.5.1.21}). 

The estimates of $[L\Phi]$, $\bA_g$ and $\bA_b$ in (\ref{3.11.3.21}) follow from using (\ref{L4BA1}) and (\ref{L4conn}). Similarly, in view of Proposition \ref{6.7.1.23} (1), we also obtain 
\begin{align*}
&\sn\Phi^\dagger=[\sn\Phi], \eta;\quad Z\Phi^\dagger=[Z\Phi], [\sn\Phi], \, Z=L, \Lb. 
\end{align*}
Thus,  also noting $\eh\in \bA_{g,1}$ and $\tr\eta= [L\Phi]$, we conclude the estimates for $\sn\Phi^+$ and $L\Phi^+$ in (\ref{3.11.3.21}). The last three estimates have been given in (\ref{L4BA1}) and (\ref{zeh}). The proof of (\ref{3.11.3.21}) is complete.

 Due to (\ref{3.22.1.21}) and (\ref{lb})
\begin{equation*}
L\log \bb=-k_{\bN\bN}=-\f12\wp \Lb \varrho+[L\Phi].
\end{equation*}
Integrating along null geodesics $\Upsilon_{\omega,u}$,  with the help of $\tir[L\Phi], \varrho=O(\l t\r^{-1+\delta}\Delta_0^\f12)$ in (\ref{3.6.2.21}) and (\ref{6.24.1.21}),  and the assumption (\ref{6.5.1.21}), yields
\begin{equation*}
\log \bb-\log\bb(0)\ge\int_0^t -\f12\wp[\Lb\varrho]_++[L\Phi]> -C\Delta_0^\f12, 
\end{equation*}
where the constant $C>0$ depends only on the $\delta$ in Assumption \ref{5.13.11.21+}, and $[f]_+$ stands for the positive part of the scalar function $f$. With $\Delta_0$ sufficiently small, we obtain
\begin{equation*}
\bb>c^{-1}(0)e^{-C\Delta_0^\f12}>\frac{c^{-1}(0)}{4}
\end{equation*}
which gives the left-hand inequality (\ref{11.13.3.23}). 
Moreover  
\begin{equation*}
|\log \bb-\log c^{-1}(u,\omega, 0)|\les \int_0^t \{\M_0 \big(1+\f12\wp\log(\frac{\l t'\r}{2})\big)^{-1}\l t'\r^{-1}+\Delta_0^\f12\l t'\r^{-2+\delta} \}
\end{equation*}
which gives the upper bound of
\begin{equation}\label{7.9.2.24}
\frac{1}{4}c^{-1}(0)<\bb\les(1+\f12\wp\log(\f12 \l t\r))^{\frac{2\M_0}{\wp}} 
\end{equation}
 by using (\ref{11.11.2.23}) and (\ref{9.30.1.23}).

 
Now consider (\ref{3.11.4.21}). 
We derive by (\ref{1.14.4.22}) and $L v_t=2 h v_t$ that
\begin{align}\label{11.28.4.23}
L(\bb\tir \Lb \varrho)&=\tir\bb\{ -\Box_\bg \varrho+\sD \varrho-\hb L \varrho+2\zb\c \sn \varrho+\bA_b\Lb\varrho\}.
\end{align}

Noting that due to (\ref{7.04.9.19}),
\begin{equation*}
[\Lb \Phi]=\Lb \varrho+[L\Phi],
\end{equation*}
in view of Proposition \ref{geonul_5.23_23}, we have
\begin{equation}\label{8.21.1.23}
\bb \tir ([\Box_\bg \Phi], \N(\Phi, \bp\Phi))=\bb\tir \{\Lb \varrho[L \Phi]+([\sn \Phi]+\eh)([\sn\Phi]+\eh)+[L\Phi]^2\}.
\end{equation}
Moreover,
noting that as an immediate consequence of  (\ref{3.6.2.21}), (\ref{k1}) and $\tr\chi+\tr\chib=-2\tr k$, we derive schematically
\begin{align*}
\bb \tir \hb [L\Phi]&=\bb \tir (\Lb\varrho+h+ [L\Phi])[L\Phi]\\
&=\bb \tir (\Lb\varrho+h)[L\Phi]+\bb O(\Delta_0\l t\r^{-3+2\delta}).
\end{align*}
It follows by combining the above two calculations  and using (\ref{3.6.2.21}) that
\begin{equation}\label{8.21.4.23}
\bb \tir ([\Box_\bg \Phi], \N(\Phi, \bp\Phi),\hb[L\Phi])=O(1)(\tir\bb  \Lb \varrho+\bb\tir h)[L\Phi]+ \bb O(\l t\r^{-3+2\delta}\Delta_0).
\end{equation}
Therefore substituting the above two estimates and (\ref{8.21.4.23}) to (\ref{11.28.4.23}) with the help of (\ref{3.6.2.21}), we deduce
 \begin{align}\label{11.28.5.23}
L(\bb\tir \Lb \varrho)=O(1)\tir\bb \Lb \varrho([L\Phi]+\bA_b)+O(1)\bb[L\Phi]+\bb O(\l t\r^{-\frac{7}{4}+\delta}\Delta_0).
\end{align} 

Using  (\ref{3.6.2.21}), (\ref{7.9.2.24}) and the transport lemma to integrate along any null geodesic $\Upsilon_{\omega, u}$, we obtain
\begin{equation*}
|\bb \tir \Lb \varrho|\les 1+\Delta_0^\f12\les 1.
\end{equation*}
Using  $[\Lb \Phi]=\Lb \varrho+[L \Phi]$ and (\ref{3.6.2.21}) again, we conclude
\begin{equation}\label{8.21.1.23+}
|\bb \tir [\Lb \Phi]|\les 1.
\end{equation}
Due to Proposition \ref{6.7.1.23}, schematically,\begin{footnote}{The sums below are schematic sums, meaning each term in the sum appears when taking $|\cdot|$.}\end{footnote}
$$\Lb\Phi=[\Lb\Phi]+[\sn\Phi]=[\Lb\Phi]+O(\l t\r^{-2+\delta}\Delta_0).$$ Thus, in view of (\ref{8.21.1.23+}), the first estimate in (\ref{3.11.4.21}) is proved. The estimates for $\bb \tir(|\bT \Phi|+|k_{\bN\bN})$ in (\ref{3.11.4.21}) thus follow by using (\ref{3.6.2.21}) and (\ref{3.22.1.21}). The proof of (\ref{3.11.4.21}) is complete.

Using the last estimate in (\ref{3.11.4.21}) and (\ref{lb}), following the same way in the above we have 
$$|\bb-\bb(0)|\les \log \l t\r.$$
Thus the other estimate in (\ref{11.13.3.23}) follows in view of $\bb(0)=c^{-1}(0)$.

Substituting (\ref{8.21.1.23+}) to (\ref{11.28.5.23}), and using (\ref{3.6.2.21}) leads to 
\begin{equation*}
\bb^{-1}L(\bb\tir \Lb \varrho)=O(\l t\r^{-\frac{7}{4}+\delta}\Delta_0^\f12).
\end{equation*}
Integrating it along the null geodesic $\Upsilon_{\omega, u}$ yields
\begin{equation*}
\bb\tir \Lb \varrho=c^{-1} \ckc^{-1} u\Lb \varrho(0,u,\omega)+O(\Delta_0^\f12).
\end{equation*}    
(\ref{12.19.1.23}) follows as a consequence of (\ref{12.6.2.23}) and $c^{-1}, \ckc^{-1}=O(1)$ (due to (\ref{9.30.1.23})).



By using  the estimate for $\ud\bA$ in (\ref{3.6.2.21}), (\ref{11.13.3.23}), (\ref{1.22.4.22}) and transport lemma, we have (\ref{3.28.3.21}).

Next we will improve the decay estimate of $v_A$ in (\ref{6.24.1.21}).
Note, due to $\sn_\bN \bN=-\sn\log\bb$, $\sn \bN=\theta$, and (\ref{7.04.11.19}), 
\begin{equation}\label{12.11.1.23}
\sn_\bN (v_A)=v_\bN\sn\log \bb +(\sn v)_\bN,\, \, \sn (v_A)=(\sn v)_A-v_\bN \theta,\, \, \sn(v_\bN)=(\sn v)_\bN+v_A\c \theta.
\end{equation}
Using (\ref{5.17.1.21}), we write $\tr\theta=\frac{2}{\tir}+\fB+\bA$. Thus $\tir \theta=O(1)$ due to (\ref{3.11.4.21}) and (\ref{3.6.2.21}).
Using (\ref{3.6.2.21}) and also $|v,\varrho|\les \l t\r^{-1+\delta}\Delta_0^\f12$, we have 
\begin{equation}\label{11.30.1.23}
|\sn_\bN(v_A)|\les\l t\r^{-2+\delta}\Delta_0+|v_\bN||\sn\log \bb|,\quad |\sn (v_A)|\les \l t\r^{-2+\delta}\Delta_0^\f12,\quad  |\sn (v_\bN)|\les |[\sn v]|+\l t\r^{-1}|v_A|
\end{equation}
Applying (\ref{3.28.2.21}) to $|v_A|^2_\ga$ with the help of (\ref{3.28.3.21}) and the first two estimates in the above,
\begin{align*}
|v_A|_\ga^2&\les\int_u^{u_*} (\bb|\sn_\bN (v_A)|+\Delta_0\l t\r|\sn(v_A)|)|v_A| du'\\
&\les \int_u^{u_*} \l t\r^{-1+\delta}\Delta_0 |v_A| du'
\end{align*}
which gives $|v_A|\les \l t\r^{-1+\delta}\Delta_0$, as desired.  The proof of (\ref{6.24.1.21}) is completed.

To see the decay of $v_\bN$ in (\ref{9.18.3.23}), we recast (\ref{9.29.6.23}) symbolically  
\begin{equation}\label{11.21.1.23}
\tir^{-1}\overline{c^{-2} v_\bN}=\overline{\tr\eta+c^{-2}\bA_b v_\bN}.
\end{equation}
 It follows by  using \Poincare  inequality and Sobolev embedding on spheres that
    \begin{align*}
    \|c^{-2}v_\bN-\overline{c^{-2}v_\bN}\|_{L_\omega^2}\les \|\tir \sn (c^{-2}v_{\bN})\|_{L_\omega^2}. 
    \end{align*}

  By the above inequality, (\ref{7.04.9.19}) and (\ref{11.30.1.23}), we infer from (\ref{11.21.1.23}) that
 \begin{align*}
  \|v_\bN\|_{L_\omega^2}&\les |\overline{ c^{-2}v_\bN}|+\|[\tir\sn \Phi], v_A\|_{L^2_\omega}\\
  &\les\|[S\Phi], \tir\bA_b v_\bN, [\tir\sn\Phi], v_A\|_{L_\omega^2}.
 \end{align*}
 Using  $\tir \bA=O(\Delta_0^\f12 \l t\r^{-\frac{3}{4}+2\delta})$,  we then derive
 \begin{align*}
  \| v_\bN\|_{L_\omega^2}\les \sum_{X=S, \Omega}\|[X\Phi]\|_{L_\omega^2}+\|v_A\|_{L_\omega^2}.
 \end{align*} 
 Using $\sn_\bN(v_A)=[\sn_A v]+\ud\bA v_\bN$, applying (\ref{10.22.2.22}) to $v_A$ gives
 \begin{align}\label{8.7.1.24}
 \|v_A\|_{L_\omega^p}&\les \|\bb\sn_\bN (v_A)\|_{L_u^1 L^p_\omega}\les \|\bb\ud\bA v_\bN\|_{L_u^1 L_\omega^p}+\|\bb[\sn v]\|_{L^1_u L_\omega^p}, 2\le p<\infty
\end{align}
Combining the above two estimates and using (\ref{3.6.2.21}), we deduce
\begin{align*}
\|v\|_{L^2_u L_\omega^2}\les \sum_{X=S,\Omega}\|[X\Phi]\|_{L^2_u L_\omega^2}, \|v_A\|_{L_\omega^2}\les \l t\r^{-1+\delta}\log \l t\r\sum_{X=S, \Omega}\|X\Phi\|_{L^2_u L_\omega^2}
\end{align*}
which gives (\ref{9.18.3.23}).
\end{proof}
\begin{proof}[Proof of (\ref{6.22.1.21})-(\ref{10.1.1.22})] 
Note the following symbolic formula due to (\ref{5.23.1.23}) and (\ref{11.30.2.23})
\begin{equation}\label{8.23.5.23}
[L\Lb\Phi]=L[\Lb\Phi]+\zb [\sn\Phi]+L\log c[\Lb\Phi], \quad \Lb[L\Phi]=[\Lb L\Phi]+\sn\log \bb[\sn\Phi]+\Lb\log c[L\Phi].
\end{equation}
Using (\ref{3.19.2}), we have
\begin{equation*}
[\Lb L\Phi]=[L\Lb\Phi]-2(\zb-\ze)[\sn\Phi]+2k_{\bN\bN}[\bN\Phi].
\end{equation*}
Due to $[\Lb\Phi]=\Lb \varrho+[L\Phi]$   
\begin{equation*}
L[\Lb\Phi]=L\Lb \varrho+L[L\Phi], 
\end{equation*}
where for $L\Lb \varrho$ we can apply (\ref{11.28.5.23}) and (\ref{lb}).  Combining the above formulas, the estimates of the scalar components in (\ref{6.22.1.21}) in the first line follow by using (\ref{3.11.4.21}) and (\ref{3.6.2.21}). 

Similar to the above, we can obtain the following formula,  with $l=0,1,2$,
\begin{align}\label{8.26.1.23}
\begin{split}
&X^l[L\Lb\Phi],  X^l L[\Lb \Phi], X^l[\Lb L\Phi], X^l \Lb[L\Phi]\\
&=X^l\Big((\tir^{-1}+[\Lb\Phi])\fB\Big)+X^l\big(\bA_b\c(\fB+[L\Phi])+\sD\varrho+L[L\Phi]+\N(\Phi,\bp\Phi)\\
&+\ud\bA\bA_{g,1}\big)
\end{split}
\end{align}
where the last term only appears in $X^l[\Lb L\Phi]$ and $X^l \Lb[L\Phi]$, we included $\fB[L\Phi]$ into $\N(\Phi,\bp\Phi)$, 
and we used (\ref{6.30.2.19}), $h=\bA_b+\frac{1}{\tir}$, and $\hb=-(\bA_b+\frac{1}{\tir})+\fB$.

Moreover using (\ref{5.23.1.23}), (\ref{11.30.2.23}) and (\ref{3.19.2}), we obtain
\begin{align*}
L\Lb\Phi^\|=(\sn_L+L\log c)[\sn\Phi]+[\Lb\Phi]\bA_{g,1},\quad \Lb L\Phi^\|=L \Lb\Phi^\|-2(\zb-\ze)\sn\Phi^\|+2 k_{\bN\bN}\bN\Phi^\|.
\end{align*}
Thus due to (\ref{3.6.2.21}) and $[\Lb\Phi], k_{\bN\bN}=O(\bb^{-1}\l t\r^{-1})$, we obtain 
\begin{equation*}
L\Lb\Phi^\|, \Lb L\Phi^\|= O(\l t\r^{-3+2\delta}\Delta_0^\f12).
\end{equation*}
The above pointwise decay is much better than $\l t\r^{-1}\fB$, the proof of  (\ref{6.22.1.21}) is thus complete. 
  (\ref{4.3.3.21}) follows by using Proposition \ref{geonul_5.23_23}, (\ref{6.22.1.21}), (\ref{3.6.2.21}) and (\ref{3.11.4.21}). 

Next we consider the lower order estimates in (\ref{8.23.1.23}). The first estimate in (\ref{8.23.1.23}) follows as a consequence of (\ref{6.22.1.21}) by noting $\fB=[\Lb\Phi]+[\bar\bp\Phi]$, (\ref{3.6.2.21}) and (\ref{LbBA2}). For the second estimate, the case of $X=S$ can be proved by applying  (\ref{6.22.1.21}), (\ref{3.6.2.21}) and the assumption $\Lb \tir=O(1)$. If $X=\tir \sn$, we use (\ref{5.23.1.23}) to derive
\begin{equation}\label{4.5.5.24}
\begin{split}
&X\Lb \Phi^\|=(X\log c+\sn_X)[\sn \Phi]+X\bN^A[\Lb \Phi],\quad [X \Lb\Phi]=(X\log c+X)[\Lb\Phi]+X\bN^A[\sn\Phi]\\
&\Lb X\Phi^\|=(\Lb\log c+\sn_\Lb)(X \Phi^\|)+\Lb \bN^A[X\Phi], \quad [\Lb X\Phi]=(\Lb\log c+\Lb)[X\Phi]+\Lb \bN^A \c X\Phi_A
\end{split}
\end{equation}
where due to (\ref{11.30.2.23}) and (\ref{3.6.2.21}), $X\bN^A=O(1)$ and  $\Lb \bN^A=\ud\bA$. Noting $\tr\eta=[L\Phi]$, if $\Lb \tir =O(1)$, by  (\ref{3.6.2.21}), (\ref{LbBA2}) and (\ref{L2BA2})  we have
\begin{equation*}
\tir\sn\Lb \Phi,\sn_\Lb (\tir\sn)\Phi=O(1)[\Lb\Phi]+O(\l t\r^{-1+\delta}\Delta_0), \quad O(1)[\Lb\Phi]+O(\l t\r^\delta\Delta_0)_{L^2_\Sigma}.
\end{equation*}
Thus we completed the second estimate in (\ref{8.23.1.23}). 

Using the second estimate in (\ref{8.23.1.23}) and  Proposition \ref{geonul_5.23_23}, estimates upto the first order  $X$-derivative on $\wt\N(\Phi, \bp\Phi)$ in (\ref{3.29.1.23}) can be obtained by using (\ref{3.6.2.21}), (\ref{3.11.4.21}), (\ref{L4BA1}) and (\ref{L2BA2}). 

Next, we will proceed to establish (\ref{6.8.4.23})-(\ref{10.15.2.22}) and (\ref{3.16.1.22}) successively. We will first prove the estimates containing no $X$-derivative, then apply the first order $X$-derivative.
  
In view of (\ref{3.20.1.22}), we write
\begin{equation}\label{12.1.2.23}
X^l\widetilde{L (\Xi_4)}=X^l(\sD, L^2)\varrho+X^l\Big(L\varrho(\bA+\tir^{-1})+\bA_{g,1}^2+\N(\Phi, \bp\Phi)\Big)
\end{equation}
where we included $\fB\c L\varrho$ into $\N(\Phi, \bp\Phi)$ already. 
By using (\ref{3.6.2.21}), (\ref{3.11.3.21}), (\ref{4.3.3.21}) and the lowest order estimates in (\ref{3.29.1.23}) we have
\begin{align}
\widetilde{L (\Xi_4)}&=\bAn\c \bA+O(\l t\r^{-1})[L\Phi]+L^2\varrho+\sD\varrho+\N(\Phi, \bp\Phi)\nn\\
&=O(\l t\r^{-\frac{11}{4}+\delta}\Delta_0^\f12), O(\l t\r^{-3+\delta}\Delta_0^\f12)_{L_\omega^4}, O(\l t\r^{-2+\delta}\Delta_0)_{L^2_\Sigma}. \label{4.5.2.24}
\end{align}
This gives the lowest order estimates in (\ref{6.8.4.23}). 
The $l=0$ case in  (\ref{1.29.4.22})  and (\ref{10.15.2.22}) follow by using (\ref{1.21.2.22})-(\ref{4.17.1.24}), (\ref{3.6.2.21}), (\ref{3.11.3.21}), (\ref{L2BA2}) and (\ref{LbBA2}). 

Now we are ready to show the $l=0$ case in (\ref{3.16.1.22}). In view of (\ref{6.3.1.23}),
using the lowest order estimates in (\ref{6.8.4.23}) and (\ref{4.3.3.21}), (\ref{3.6.2.21}) and (\ref{3.11.3.21}), we infer
\begin{equation}\label{2.1.1.22}
L(\tr\chi-\frac{2}{\tir})+\frac{2}{\tir}(\tr\chi-\frac{2}{\tir})+\f12 (\tr\chi-\frac{2}{\tir})^2=O(\Delta_0^\f12 \l t\r^{-3+\delta})_{L_\omega^4}, O(\Delta_0^{\f12} \l t\r^{-\frac{11}{4}+\delta}).
\end{equation}
This gives 
$$\sn_S \bA_b,\sn_L(\tir\tr\chi)=O(\Delta_0^\f12 \l t\r^{-2+\delta})_{L_\omega^4}, O(\Delta_0^{\f12} \l t\r^{-\frac{7}{4}+\delta}),$$ in view of (\ref{L4conn}) and (\ref{3.6.2.21}).
Since $\tr\chib+\tr\chi=-2\tr k=\fB$,  using the above estimate of $L\bA_b$ and (\ref{6.22.1.21}), we derive  $\sn_L \tr\chi, \sn_L \tr\chib=O(\l t\r^{-2})$.

Using the transport equations (\ref{s2}), the $l=0$ estimate in (\ref{1.29.4.22}), (\ref{3.11.3.21}) and $\tir|k_{\bN\bN}, \tr\chi|\les 1$,
we have
$$ \sn_L \chih= O(\l t\r^{-3+\delta}\Delta_0)_{L_\omega^4}, O(\l t\r^{-\frac{11}{4}+\delta}\Delta_0).$$
 The lower order estimates in the first three lines (\ref{3.16.1.22}) are complete. 
 
Using (\ref{1.27.6.24}) and (\ref{3.6.2.21}), we have 
\begin{equation*}
\sn_L(\tir\sn\log \bb)=O(\l t\r^{-1+\delta}\Delta_0)
\end{equation*}
Due to $\ze+\zb=\sn\log \bb$ and (\ref{3.6.2.21}), this gives the lower order estimate for $\ze$ in (\ref{3.16.1.22}).

We then use (\ref{12.1.2.23}), the lower order estimates in (\ref{3.16.1.22}), the estimates in (\ref{3.29.1.23}) that contains only one $X$-derivative,  and (\ref{3.6.2.21}), (\ref{3.11.3.21}), (\ref{L2BA2}) and (\ref{L2conndrv}) to derive
\begin{equation}\label{4.5.3.24}
\sn_X\widetilde{L(\Xi_4)}=O(\l t\r^{-\frac{11}{4}+\delta}\Delta_0^\f12)_{L_\omega^4}, O(\l t\r^{-2+2\delta}\Delta_0)_{L^2_\Sigma}.
\end{equation}
Thus, the estimates of the first order $X$-derivative in (\ref{1.29.4.22}) and (\ref{10.15.2.22}) can be obtained by using the lower order estimate in (\ref{3.16.1.22}), (\ref{3.6.2.21}), (\ref{3.11.3.21}), the first estimate in (\ref{8.23.1.23}), (\ref{L2BA2}) and (\ref{LbBA2}). 

Differentiating (\ref{6.3.1.23}) and (\ref{s2}) by $X$,  using these proved curvature estimates, (\ref{3.6.2.21}) and (\ref{3.11.3.21}),  (\ref{4.5.3.24}), estimates in (\ref{3.29.1.23}) upto the first order in-$X$, we can obtain the $l=1$ and $m=0$ estimates for $\sn_L \bA_b, \sn_L \bA_{g,2}$ in (\ref{3.16.1.22}).

 Using (\ref{cmu2}), the lower order curvature estimate in (3) and (\ref{3.6.2.21}), we bound
\begin{equation}\label{4.5.4.24}
|[\sn_L, \tir\sn]F|\les \l t\r^{-\frac{3}{4}+\delta}\Delta_0^\f12 |\sn F|+\l t\r^{-\frac{7}{4}+\delta}\Delta_0 |F|. 
\end{equation}
Applying it to $F=\bA_b, \bA_{g,2}$ with the help of (\ref{3.11.3.21}), it shows that these commutators are negligible terms. Thus we obtained the estimates with   $l=1, m=1$  in the first two lines of (\ref{3.16.1.22}) from their counter parts with $l=1, m=0$. 

Differentiating (\ref{1.27.6.24}) by $X$, using (\ref{L4BA1}), the lower order estimate in (\ref{3.16.1.22}) and (\ref{3.6.2.21}), we can obtain $l=1, m=0$ case in the last line of (\ref{3.16.1.22}). Applying (\ref{4.5.4.24}) to $F=\tir\ze$, in view of (\ref{zeh}), it gives
\begin{equation*}
\|[\sn_L, \tir\sn](\tir\ze)\|_{L_\omega^4}\les \l t\r^{-\frac{7}{4}+2\delta}\Delta_0^\frac{3}{2}. 
\end{equation*} 
 The case $l=1, m=1$ follows as a consequence of the case $l=1$ and $m=0$ and the above commutator estimate. Combining (\ref{zeh}) and  the proved $L^4_\omega$ estimates for $\ze$, the last line in (\ref{3.16.1.22}) drops as a consequence. 
  
Next, we show (\ref{6.7.4.23}) by letting $l=1$ in (\ref{8.26.1.23}). Using the first estimate in (\ref{8.23.1.23}), the lower order estimate in (\ref{3.16.1.22}) and $X\N(\Phi,\bp\Phi)$ estimate in (\ref{3.29.1.23}), we obtain (\ref{6.7.4.23}) by using (\ref{3.6.2.21}), (\ref{3.11.3.21}), (\ref{L2BA2}) and (\ref{L2conndrv}). 

Note by using (\ref{1.21.2.22})-(\ref{4.17.1.24}) again, the top order estimate in (\ref{10.15.2.22}) can be derived by using (\ref{L2BA2}), (\ref{LbBA2}), (\ref{L2conndrv})  and (\ref{3.6.2.21}).
 
To prove the top order estimates in (\ref{3.29.1.23}) and (\ref{6.8.4.23}), we need to prove the $l=2$ case in (\ref{8.23.1.23}) and (\ref{8.23.1.23'}). For the case $X^2=S^2$, we set $X=S$ in (\ref{8.26.1.23}), noting that we can drop the last term since we only need the formula for $S L[\Lb\Phi]$. Using the first estimate in (\ref{8.23.1.23}), the lower order estimate in (\ref{3.16.1.22})
 and the first order $X$-derivative estimates in (\ref{3.29.1.23}), (\ref{3.6.2.21}), (\ref{3.11.3.21}) and (\ref{L2BA2}), we deduce
\begin{align*}
SL[\Lb\Phi]=O(\l t\r^{-1})\fB+S L[L\Phi]+O(\l t\r^{-\frac{11}{4}+\delta}\Delta_0)_{L_\omega^4}, O(\l t\r^{-1}\fB)+O(\l t\r^{-2+\delta}\Delta_0)_{L^2_\Sigma}.   
\end{align*}
For the case $\vs^-(X^2)=0$, by setting $X=\tir \sn$ in (\ref{8.26.1.23}), we derive by using the first estimate in (\ref{8.23.1.23}), the $\tir\sn$-derivative estimates in (\ref{3.29.1.23}), (\ref{3.6.2.21}), (\ref{3.11.3.21}), (\ref{L2BA2}) and (\ref{LbBA2})
\begin{align*}
\tir \sn L[\Lb\Phi]=O(\l t\r^{-2+\delta}\Delta_0)_{L_\omega^4}, O(\l t\r^{-1+\delta}\Delta_0)_{L^2_\Sigma}. 
\end{align*}
Applying (\ref{4.5.4.24}) to  $F=[\Lb\Phi]$ gives
\begin{equation*}
\sn_L (\tir \sn)[\Lb\Phi]=O(\l t\r^{-2+\delta}\Delta_0)_{L_\omega^4}, O(\l t\r^{-1+\delta}\Delta_0)_{L^2_\Sigma}. 
\end{equation*}
Finally, due to (\ref{3.11.3.21}) and (\ref{LbBA2})
$$(\tir\sn)^2[\Lb\Phi]=O(\l t\r^{-2+\delta}\Delta_0)_{L_\omega^4}, O(\l t\r^{-1+\delta}\Delta_0)_{L^2_\Sigma}.$$
Therefore we summarize the above estimates into the last estimate in (\ref{8.23.1.23}).

To prove (\ref{8.23.1.23'}), we derive by using (\ref{5.23.1.23}) and (\ref{11.30.2.23})
\begin{align*}
[X_2 X_1 \Lb \Phi]=X_2[X_1 \Lb\Phi]+ O(1)X_1 \Lb\Phi, (X_2 X_1\Lb\Phi)^\|=\sn_{X_2}(X_1 \Lb\Phi^\|)+O(1) X_1 \Lb\Phi.
\end{align*}
If $X_1=S$, using the second estimate of (\ref{8.23.1.23}), the above formulas can be simplified to be 
\begin{align*}
[X_2 S \Lb \Phi]=O([\Lb\Phi])+X_2[S \Lb\Phi], (X_2 S\Lb\Phi)^\|=\sn_{X_2}(S\Lb \Phi^\|)+O([\Lb\Phi]). 
\end{align*}
Differentiating (\ref{4.5.5.24}) with the help of (\ref{8.23.1.23}), (\ref{3.11.3.21}) and (\ref{L2BA2}) leads to
\begin{align*}
X_2[S \Lb\Phi]&=X_2((S\log c+S)(\fB))+X_2(\bA_{g,1}\tir [\sn\Phi])\\
&=\left\{\begin{array}{lll}
\vs(X_2)([\Lb\Phi]+O(\l t\r^{-\frac{7}{4}+\delta}\Delta_0^\f12)_{L_\omega^4})+(1-\vs(X_2))O(\l t\r^{-1+\delta}\Delta_0)_{L_\omega^4}\\
\vs(X_2) ([\Lb\Phi]+O(\l t\r^{-1+\delta}\Delta_0)_{L^2_\Sigma})+(1-\vs(X_2)) O(\l t\r^\delta\Delta_0)_{L^2_\Sigma}
\end{array}\right.\\ 
\sn_{X_2}(S\Lb \Phi^\|)&=\sn_{X_2}\big((S\log c+\sn_S)[\sn\Phi]\big)+\sn_{X_2}(\bA_{g,1}\tir \fB)\\
&=O(\l t\r^{-\frac{7}{4}+\delta}\Delta_0)_{L_\omega^4}, O(\l t\r^{-1+\delta}\Delta_0)_{L^2_\Sigma}. 
\end{align*} 
Combining the above estimates, we obtain the case that $X_1=S$ in (\ref{8.23.1.23'}). If $X_2=S, X_1=\tir \sn$, using (\ref{cmu2})
\begin{align*}
[S, \tir \sn]\Lb\Phi=O(\bA)\tir^2 \sn \Lb\Phi+O(1)(\zb, \tir \bR_{AC4B})\c \tir \sn\Lb\Phi.
\end{align*}
Hence, using (\ref{8.23.1.23}), (\ref{3.11.3.21}) and the lower order estimate of (3) we obtain
\begin{align*}
[S, \tir \sn]\Lb\Phi=O(\tir \bA)[\Lb\Phi]+O(\l t\r^{-\frac{7}{4}+2\delta}\Delta_0),  O(\tir \bA)[\Lb\Phi]+O(\l t\r^{-\frac{3}{4}+2\delta}\Delta_0)_{L^2_\Sigma}. 
\end{align*}
Using the proved result for $X_2 X_1=\tir \sn S$, combining the above commutator estimate with $\tir [\Lb\Phi]=O(\bb^{-1})$,  we proved (\ref{8.23.1.23'}) for $X_2 X_1=S\tir \sn$. 

If $X_2=X_1=\tir \sn$, we derive
\begin{align*}
[(\tir\sn)^2\Lb \Phi]=\tir \sn[\tir \sn \Lb\Phi]+ O(1)\tir \sn \Lb\Phi, ((\tir \sn)^2\Lb\Phi)^\|={\tir \sn}(\tir\sn \Lb\Phi^\|)+O(1) \tir \sn \Lb\Phi.
\end{align*}
Differentiating the first line in (\ref{4.5.5.24}) by using (\ref{3.6.2.21}), (\ref{3.11.3.21}), (\ref{L2conndrv}), (\ref{LbBA2}) and (\ref{L2BA2}), we can obtain 
\begin{align*}
(\tir \sn)^2\Lb\Phi&=O([\Lb\Phi])+O(\l t\r^{-1+\delta}\Delta_0)_{L_\omega^4}, O[\Lb\Phi]+O(\l t\r^\delta\Delta_0)_{L^2_\Sigma}. 
\end{align*}
Hence (\ref{8.23.1.23'}) is proved. 

Next we consider the second order $X$-derivative estimate in (\ref{3.29.1.23}). Using Proposition \ref{geonul_5.23_23}, (\ref{8.23.1.23}), (\ref{8.23.1.23'}) we obtain the estimates by using (\ref{3.6.2.21}), (\ref{3.11.3.21}) and (\ref{L2BA2}). Using the last estimate in (\ref{3.29.1.23}), setting $l=2$ in (\ref{12.1.2.23}), we obtain the top order estimate in (\ref{6.8.4.23}) by using (\ref{L2BA2}), (\ref{L2conndrv}) and (\ref{3.6.2.21}). 


Now we consider (4). For the term $K$, we recall from \cite[(2.4.1) in Page 47, (3.1.2d) in Page 56]{CK} the formulas of the Gaussian curvature under the conformal change of metric $\ga=c^{-2}\gac$
\begin{equation*}
2\Kc-(\tr\thetac)^2+|\thetac|^2=\stackrel{\circ}R-2\stackrel{\circ}R_{\bN\bN}=0;\quad K=c^2\Kc+\sD\log c.
\end{equation*}
Since $\tr\thetac=\frac{2}{c\tir}+\bA_b$, the pointwise estimate of (4)  follows by using (\ref{3.6.2.21}); and the estimate of $\sn\Kc$ can be obtained by using (\ref{3.6.2.21}) and (\ref{3.11.3.21}). 

In view of the above two formulas, we also have
\begin{equation}\label{3.31.4.22}
K-\frac{1}{4}c^2 (\tr\thetac)^2=-\f12|\hat\theta|^2+\sD \log c.
\end{equation}
Due to  (\ref{1.6.1.21}) and using $\tr\eta=[L\Phi]$
\begin{equation*}
c\tr\stackrel{\circ}\theta=\tr\chi+[L\Phi]
\end{equation*}
Combining the above two lines
\begin{equation}\label{10.1.1.23}
K-\tir^{-2}=(\bA_b+[L \Phi]+\tir^{-1})(\bA_b+[L\Phi])+\sD \log c+\bA^2_{g,2}.
\end{equation}
Using (\ref{3.6.2.21}), (\ref{3.11.3.21}) and (\ref{L4BA1}) we conclude the estimates of $K$ and $X(K-\frac{1}{\tir^2})$ in (4).
\end{proof}
\begin{proof}[Proof of (5)] We first prove (\ref{1.27.2.24}). It suffices to regard $\fB=\Xi_4$ in the proof since $\fB=\Xi_4+[L\Phi]$ with $[L\Phi]$ verifying much better estimates.

 It is straightforward to derive
\begin{equation*}
\bb^{-1}\tir^{-1}L (\tir \bb \Xi_4)=\wt{L \Xi_4}+(h-\tir^{-1})\Xi_4. 
\end{equation*}
Using (\ref{cmu2}), we deduce
\begin{align*}
L \sn(\tir \bb \Xi_4)+\chi\c \sn (\tir \bb \Xi_4)=\sn L(\tir \bb \Xi_4)=\sn(\bb \tir \wt{L \Xi_4}+\bb \tir\bA_b \Xi_4).  
\end{align*}
Using transport Lemma, Proposition \ref{12.21.1.21},  (\ref{3.20.1.22}), (\ref{6.8.4.23}),  (\ref{3.11.3.21}) and (\ref{3.6.2.21}), we obtain
\begin{align*}
\|\tir^2\sn(\bb \Xi_4)\|_{L_\omega^4}&\les \La_0+\int_0^t \tir^2\|\sn(\bb \wt{L \Xi_4}+\bb \bA_b \Xi_4)\|_{L_\omega^4}\\
&\les \La_0+\int_0^t \Delta_0 \l t'\r^{-\frac{7}{4}+\delta}\log \l t'\r dt'\les \Delta_0,
\end{align*}
as desired in (\ref{1.27.2.24}). The lower order estimate in (\ref{1.27.4.24}) can be obtained similarly.  
Using (\ref{4.6.1.24}), we derive
\begin{align*}
|\bb \Lb \sn \varrho|\les |\sn(\bb \fB)|+\l t\r^{-1}|\bb\sn\varrho|+|\bb \ud \bA[L\Phi]|+|\bb\sn [L\Phi]|.
\end{align*}
It follows by using (\ref{3.6.2.21}) and (\ref{3.11.3.21}) that
\begin{equation*}
\|\bb \Lb \sn \varrho\|_{L_\omega^4}\les \|\sn(\bb \fB)\|_{L_\omega^4}+\l t\r^{-3+2\delta}\Delta_0 \log \l t\r\les  \l t\r^{-2}\Delta_0,
\end{equation*}
where we used the first estimate in (\ref{1.27.2.24}) to derive the last estimate. Hence (\ref{1.27.2.24}) is proved. 

To see (\ref{1.27.5.24}), we derive by (\ref{1.27.6.24}) that 
\begin{align}\label{1.31.2.24}
\sn_L\sn \bb+\f12\tr\chi\sn\bb=\sn(\bb k_{\bN\bN})-\chih\c \sn \bb,
\end{align}
which, in view of (\ref{1.27.2.24}) and Proposition \ref{12.21.1.21}, gives
\begin{align*}
\|\tir\sn \bb\|_{L_\omega^4}&\les \La_0+\int_0^t \tir\|\sn(\bb k_{\bN\bN})\|_{L_\omega^4}\les\La_0+\log \l t\r\Delta_0,
\end{align*}
as stated.

Next we prove the second order estimate in (\ref{1.27.4.24}). Due to $\fB=\Lb\varrho+[L\Phi]$ and the fact that $[L\Phi]$ verifies much better estimate, it suffices to consider $\fB=\Lb \varrho$ in the proof. Note
\begin{align*}
L(\tir\sn)^2(\bb \tir \Lb \varrho)=[(\tir \sn)^2, L](\bb \tir \Lb \varrho)+(\tir \sn)^2 L(\bb\tir \Lb \varrho).
\end{align*}
 Differentiating (\ref{11.28.4.23}) gives 
\begin{align*}
(\tir \sn)^2 L(\bb\tir \Lb \varrho)&=(\tir \sn)^2(\tir \bb \bA_b\Lb \varrho)+(\tir\sn)^2(\tir \bb (\sD\varrho+\Box_\bg \varrho))+(\tir \sn)^2\big(\tir \bb(\hb L\varrho+\bA_{g,1}^2)\big) \\
&=O(\l t\r^{-\frac{3}{2}+\delta}\Delta_0\log \l t\r^\f12)_{L_u^2 L_\omega^2},
\end{align*}
where we used (\ref{L2conndrv}), (\ref{L2BA2}), (\ref{LbBA2}), (\ref{3.6.2.21}) and (\ref{3.29.1.23}). We bound the commutator by using (\ref{cmu_2}), (\ref{3.6.2.21}) and (\ref{3.11.3.21})
\begin{align*}
\tir \sn[\tir \sn, L]f&=\tir \sn(\tir\bA \sn f)=\tir^2 \sn\bA \sn f+\bA \tir^2 \sn^2 f\\
&=O(\l t\r^{-\frac{11}{4}+\delta}\Delta_0^\f12)_{L_\omega^4}\tir^2 \sn f+O(\l t\r^{-\frac{7}{4}+\delta}\Delta_0^\f12)\tir^2 \sn^2 f.
\end{align*}
We then use  Sobolev embedding on spheres to derive
\begin{align*}
\|\tir \sn[\tir \sn, L]f\|_{L_u^2 L_\omega^2}\les \l t\r^{-\frac{7}{4}+\delta}\Delta_0^\f12\|(\tir\sn)^{1+\le 1} f\|_{L_u^2 L_\omega^2}.  
\end{align*}
Combining the above estimate with (\ref{4.5.4.24}) gives
\begin{equation}\label{4.6.2.24}
\|[(\tir \sn)^2, \sn_L]f\|_{L_u^2 L_\omega^2}\les \l t\r^{-\frac{7}{4}+\delta}\Delta_0^\f12\|(\tir\sn)^{1+\le 1}f\|_{L_u^2 L_\omega^2}. 
\end{equation}
Hence applying the above estimate to $f=\bb\tir \Lb \varrho$ leads to
\begin{equation*}
\|\sn_L(\tir\sn)^2(\bb \tir \Lb \varrho)\|_{L_u^2 L_\omega^2}\les \l t\r^{-\frac{7}{4}+\delta}\Delta_0^\f12\|(\tir\sn)^{1+\le 1}(\bb\tir\Lb \varrho)\|_{L_u^2 L_\omega^2}+\l t\r^{-\frac{3}{2}+\delta}\Delta_0\log \l t\r^\f12. 
\end{equation*}
Integrating the above estimate along the null cone by using the transport lemma  and Proposition \ref{12.21.1.21} yields 
$$
(\tir \sn)^2(\bb \Lb \varrho)=O(\l t\r^{-1}\Delta_0)_{L_u^2 L_\omega^2}, 
$$
which is the top order estimate in (\ref{1.27.4.24}).

Next applying (\ref{4.6.2.24}) to $f=\bb$ gives 
\begin{align*}
\|[(\tir\sn)^2, \sn_L]\bb\|_{L_u^2 L_\omega^2}\les \l t\r^{-\frac{7}{4}+\delta}\Delta_0^\f12\|(\tir\sn)^{1+\le 1}\bb\|_{L_u^2 L_\omega^2}. 
\end{align*}
Using (\ref{lb}) and (\ref{1.27.4.24}), we derive
\begin{equation*}
\|\sn_L(\tir\sn)^2 \bb\|_{L^2_u L_\omega^2}\les \|(\tir\sn)^2(\bb k_{\bN\bN})\|_{L_u^2 L_\omega^2}+\l t\r^{-\frac{7}{4}+\delta}\Delta_0^\f12\|(\tir\sn)^{1+\le 1}\bb\|_{L_u^2 L_\omega^2}.
\end{equation*}
Integrating along the null cone by using the transport lemma and using (\ref{1.27.4.24}) and (\ref{1.27.5.24}),  we obtain
\begin{equation*}
\|(\tir \sn)^2\bb\|_{L^2_u L_\omega^2}\les \log \l t\r\Delta_0, \|\sn_L(\tir\sn)^2 \bb\|_{L^2_u L_\omega^2}\les \l t\r^{-1}\Delta_0. 
\end{equation*}
Due to Sobolev embedding, (\ref{zeh}), (\ref{1.27.5.24}) and the above estimate
\begin{align*}
\|\sn\bb\|_{L_u^2 L_\omega^\infty}&\les \|(\tir\sn)^{\le 1} \sn \bb\|_{L_u^2 L_\omega^4}\les \|\tir^2\sn^3 \bb\|_{L_u^2 L_\omega^2}^\f12 \|\tir \sn^2 \bb\|^\f12_{L_u^2 L_\omega^2}+\l t\r^{-1}\log \l t\r\Delta_0\\
&\les\log \l t\r\l t\r^{-1+\f12\delta}\Delta_0.
\end{align*}
The remaining estimates in (\ref{2.9.2.24}) can be  derived by using the proved part, (\ref{4.5.4.24}) and (\ref{1.27.5.24}). Here we omit the detailed checking. 
\end{proof}
\subsection{Comparison results for vector-fields and frames}\label{framecmp}
\begin{lemma}\label{3.19.5.24}
Under the bootstrap assumptions (\ref{3.12.1.21})-(\ref{1.25.1.22}), we have

\begin{enumerate}
\item
There exists a constant $C>0$ such that,
\begin{equation}\label{comp1}
\tir\approx \l t\r\les \ckc\tir-C\Delta_0^\f12 t\le r\le t+u_*.
\end{equation}

\item For an $S_{t,u}$ tangent 1-form $V$, there holds
\begin{equation}\label{9.30.2.22}
\sum_{a=1}^3(V_j {}\rp{a}\Omega^j)^2=r^2 |V|^2_e=r^2c^{-2}|V|^2\approx\tir^2 |V|^2.
\end{equation}
\item We have $|y'|\les \l t\r^{-1+\delta}\Delta_0$, which improves the first assumption in (\ref{1.25.2.22}). 
\end{enumerate}
\end{lemma}
\begin{proof}
(\ref{9.29.3.22}) has given the last estimate in (\ref{comp1}). To see the remaining inequality for $r$, we first make an auxiliary bootstrap assumption that  $r\ge \l t\r^\f12$ for all $t>2$.

With $r(0,u)=u$, using (\ref{8.17.1.22}), integrating along null geodesics yields
\begin{equation}\label{9.30.1.22}
|r-u-\ckc t|\les\int_0^t \{|c-\ckc|+|v(\hN)|+|y'|\} dt'
\end{equation}
where we used $\frac{c x^k}{r}-v^k=O(1)$.  Substituting (\ref{6.24.1.21}) and (\ref{1.25.1.22}) to the above, we infer
\begin{equation*}
r-\ckc\tir=O(\Delta_0^\f12) \l t\r.  
\end{equation*}
This leads to $r>\ckc t-C\l t\r \Delta_0^\f12 +u$. Since $\ckc>\f12 c_*$  due to (\ref{9.30.1.23}), $r\ge \ckc\tir-C\Delta_0^\f12\l t\r\ges \tir$ since $\Delta_0$ is sufficiently small. It implies $r>\l t\r^\f12$, which improves the assumption. This gives (\ref{comp1}). (2) follows as a consequence of (\ref{5.22.1.22}) and (\ref{comp1}).

From (2) and (\ref{3.6.2.21}), we have $|\Omega \log c|\les \l t\r^{-1+\delta}\Delta_0$. Using (\ref{3.22.5.21}) and $\la(0)=0$ in Proposition \ref{12.21.1.21}, we obtain $\la=O(\l t\r^\delta)\Delta_0$. 
With the help of this estimate, using (\ref{12.20.3.21}) and $r\approx \l t\r$, we obtain $|y'|\les \l t\r^{-1+\delta}\Delta_0$ as stated in (3). 

\end{proof}
\begin{corollary}
Under the bootstrap assumptions (\ref{3.12.1.21})-(\ref{1.25.1.22}), the following estimates hold
\begin{align}
&|\pio_{\bN A}|\les \Delta_0, \quad |\pio_{L A}|\les \Delta_0^\f12\label{1.25.3.22}\\
&\|\pio_{AB}|\les \Delta_0 \l t\r^{-\frac{3}{4}+2\delta} \quad |\pio_{L \Lb}, \pio_{\Lb \Lb}|\les \Delta_0 \l t\r^{\delta}\label{1.25.4.22}\\
&\sn_Y\Omega_A, \tir\sn_Y L =O(1), Y=e_A,\, L, \sn_\Lb \Omega_A=O(1), \tir\sn_\Lb L=O(\l t\r^\delta\Delta_0) \label{1.25.2.22}
\end{align}
\end{corollary}
\begin{proof}
We rewrite the nontrivial components in Proposition \ref{3.22.6.21} symbolically below
\begin{align}\label{8.23.11.23}
\begin{split}
{}\rp{a}\pi_{L\Lb}, {}\rp{a}\pi_{\Lb\Lb}=\sn\log \bb\c \Omega,\quad &{}\rp{a}\pi_{AB}=\la(\tr\chi+\bAn+\bA_{g,2})+\Omega \log c\\
{}\rp{a}\pi_{\bN A}-c^{-1}\la\rp{a}\sn_A\log \bb=\sn(c^{-1}\la), &{}\rp{a}\pi_{\bT A}+c^{-1}\la\rp{a}\sn_A\log \bb=v^*_A+(\la\bA_{g,1}+\bAn)\Omega
\end{split}
\end{align}
where $v^*_A=v^j \tensor{\ud\ep}{^a_j^l} {e_A}_l$ and we dropped the factor of $c^m$ in the above formulas.

We then obtain the estimates in (\ref{1.25.3.22}) and (\ref{1.25.4.22})  by using  the bound for $\sn\la$ in (\ref{1.25.1.22}), $\tir\tr\chi\approx 1$,  (\ref{3.6.2.21}), (\ref{6.24.1.21}) and (\ref{9.30.2.22}).

Due to the bounds of $\bA_b$ and $\bA_{g,2}$ in (\ref{3.6.2.21}), we have $|\tir\chi, \tir\chib|\les 1$. 
The estimates of $\sn_Y \Omega_A$  in (\ref{1.25.2.22}) then can be obtained by using  Proposition \ref{3.22.6.21}, (\ref{3.6.2.21}) and (\ref{1.25.3.22}), 
and the estimates of $\sn_Y L$ is obtained by using  Proposition \ref{6.7con}, (\ref{3.6.2.21}) and $\bb\tir k_{\bN\bN}=O(1)$.  
\end{proof}
Using (\ref{9.30.2.22}) and the following result, all the estimates with $X=\tir \sn$ in Lemma \ref{5.13.11.21} hold with $\tir \sn$ replaced by $\Omega$.
\begin{lemma}
 Let $F$ be an $S$-tangent tensor, and $X=\Omega, S, \tir \Lb$. If $\Lb \tir=O(1)$, there holds
 \begin{align}\label{9.8.2.22}
 \sn^{\le 1}_X(F\c \Omega)\approx\sn_X^{\le 1}(\tir F_A).
 \end{align}
\end{lemma}
\begin{proof}
 Using (\ref{9.30.2.22}) and (\ref{1.25.2.22}), we have
\begin{equation*}
|\sn_X(\tir^{-1}{}\rp{b}\Omega^A)|\les 1.
\end{equation*}
We then conclude
\begin{align*}
\sn_X(F_A \Omega^A)=\sn_X(F_A\tir) \tir^{-1}\Omega^A+\tir F_A \c O(1)
\end{align*}
which gives (\ref{9.8.2.22}) in view of (\ref{9.30.2.22}).


\end{proof}
 

\begin{lemma}\label{6.30.4.23}
Let $p=2,$ or $\infty$. With $n=2, 3$, $X_2, \cdots, X_n\in\{\Omega, S\}$, there hold 
\begin{align}
\Sc(X^n v)&=[X^n v]+O(1)\sn_X^{\le n-2}(X_1 v^\|)+O(\l t\r^{-\frac{3}{4}+\delta}\Delta_0^{1-\f12\vs^+(X_n\cdots X_2)})X^{\le n-2}[X_1 v]\nn\\
&+\max(n-2,0)\{O(\l t\r^{-\frac{3}{4}+\delta}\Delta_0^{\f12+\frac{1}{p}})_{L_u^p L_\omega^4}X_1 v^\|\},\nn\\
\Sc(X^n v)&=[X^n v]+O(1)\sn_X^{\le n-2}(X_1 v^\|)+O(\l t\r^{-1+\delta}\Delta_0)_{L_u^2 L_\omega^2}X^{\le n-2}[X_1 v]\nn\\
&+\max(n-2,0)\{O(\l t\r^{-\frac{3}{4}+\delta}\Delta_0^{\f12+\frac{1}{p}})_{L_u^p L_\omega^4}X_1 v^\|\},\label{1.10.1.23}\\
\Ac(X_n\cdots X_1 v)&=(X^n v)^\|+O(1)\sn_X^{\le n-2}[X_1 v]+O(\l t\r^{-\frac{3}{4}+\delta}\Delta_0^{1-\f12\vs^+(X_n\cdots X_2)})X^{\le n-2}(X_1 v^\|)\nn\\
&+\max(n-2,0)\{O(\l t\r^{-\frac{3}{4}+\delta}\Delta_0^{\f12+\frac{1}{p}})_{L_u^p L_\omega^4}[X_1 v]\}.\nn
\end{align}
Here
$X^{\le n-2}$ represents differentiation up to the $(n-2)$-th order with by $X_i$. These $X_i$-derivatives appear among the terms $X_n, X_{n-1}, \ldots, X_2$ on the left-hand side.
 
 With $n=4$, if $\tir^{-1}\sn_X^2 \Omega=O(1)+O(\l t\r^{\delta-\frac{3}{4}}\Delta_0^{\f12+\frac{1}{p}})_{L_u^p L_\omega^4}$,  we have 
\begin{align*}
\Sc(X^n v)
&=[X^n v]+O(1)\sn_X^{\le n-2}(X_1 v^\|)+O(\l t\r^{-\frac{3}{4}+\delta}\Delta_0^{\f12+\frac{1}{p}})_{L_u^p L_\omega^4}\sn_X(X_1 v^\|)\\
&+O(\l t\r^{-\frac{3}{4}+\delta}\Delta_0^{1-\f12
\vs^+(X_n\cdots X_2)})X^{1+\le 1}[X_1 v]+O(\l t\r^{-\f12+\delta}\Delta_0)_{L_u^2 L_\omega^2}X_1 v^\|\\
&+O(\l t\r^{-\frac{3}{4}+\delta}\Delta_0^{\f12+\frac{1}{p}})_{L_u^p L_\omega^4}[X_1 v]\\
\displaybreak[0]
\Ac(X_n\cdots X_1 F)
&=(X^n F)^\|+O(1)\sn_X^{\le n-2}([X_1 v])+O(\l t\r^{-\frac{3}{4}+\delta}\Delta_0^{\f12+\frac{1}{p}})_{L_u^p L_\omega^4}\sn_X[X_1 v]\\
&+O(\l t\r^{-\frac{3}{4}+\delta}\Delta_0^{1-\f12 \vs^+(X_n\cdots X_2)})X^{1+\le 1}(X_1 v^\|)+O(\l t\r^{-\f12+\delta}\Delta_0)_{L_u^2 L_\omega^2}[X_1 v]\\
&+O(\l t\r^{-\frac{3}{4}+\delta}\Delta_0^{\f12+\frac{1}{p}})_{L_u^p L_\omega^4}X_1 v^\|
\end{align*}
The factor of $O(\l t\r^{-\frac{3}{4}+\delta}\Delta_0^{\f12+\frac{1}{p}})_{L_u^p L_\omega^4}$ in the above can be improved to be $ O(\l t\r^{-\frac{3}{4}+\delta}\Delta_0)_{L_\omega^4}$ if $X_i=\Omega, i=2\cdots n$. 
\end{lemma}
\begin{remark}
 Note if $X_1=L, \Lb$, $X_1 v^\|=[\sn\Phi]$ and if $X_1=\Omega$, $X_1 v^\|=\eta(\Omega)$. This fact will be used frequently together with the above result. 
 The assumption in $n=4$ case will be checked in Lemma \ref{3.17.2.22} . Before that the $n=4$ result will not be used. 
\end{remark}
\begin{remark}
For application, we can replace the $\Omega$ vectors appeared in $X_n\cdots X_2$ by $\tir \sn$ without changing the result. 
\end{remark}
\begin{proof} 
Let $\sn_X^n(X\bN)$ stand for applying $\sn_X^n$ to either the angular or radial components of  $ X\bN$. $\sn_X^{n}[\bN\cdots\bN]$ stands for the terms of 
 $$ \sum_{a_1+\cdots+a_l=n-l}\sn_X^{a_1}(X\bN)\sn_X^{a_2}(X\bN)\cdots \sn_X^{a_l}(X\bN),\, n\ge l\ge 1.$$  And let $\sn_X^{n}[\bN\cdots\bN]^*$ stand for the terms of 
 $$\sum_{a_1+\cdots+a_l=n-l}\sn_X^{a_1}X\log c\sn_X^{a_2} X\bN\cdots \sn_X^{a_l}X\bN, \, n\ge l\ge 1.$$ 

In view of Proposition \ref{6.7.1.23} (3), it suffices to show 
\begin{align}\label{6.30.1.23}
\begin{split}
&|\sn_X[\bN\cdots \bN]|\les 1,\\
 &\sn_X^{l+ 1}[\bN\cdots \bN]=\left\{\begin{array}{lll}
O(1)+O(\l t\r^{-\frac{3}{4}+\delta}\Delta_0^\f12 )_{L_\omega^4}, O(1)+O(\l t\r^{-\frac{3}{4}+\delta}\Delta_0)_{L_u^2 L_\omega^4}, l=1\\
O(1)+O(\l t\r^{-\f12+\delta}\Delta_0)_{L_u^2 L_\omega^2}, l\le 2
\end{array}\right. \\ 
 &\sn_X^{l+1}[\bN\cdots \bN]^*=\left\{\begin{array}{lll}
 O(\l t\r^{-\frac{3}{4}+\delta}\Delta_0^{1-\f12\vs^+(X_n\cdots X_2)}), \, l\le 1\\
 O(\l t\r^{-\frac{3}{4}+\delta}\Delta_0^\f12 )_{L_\omega^4}, O(\l t\r^{-\frac{3}{4}+\delta}\Delta_0 )_{L_u^2 L_\omega^4}, O(\l t\r^{\delta-1}\Delta_0)_{L^2_u L^2_\omega}\, l\le 2
 \end{array}\right.
\end{split}
\end{align}
If $X_i=\Omega, i=2,\cdots, n$, the $O(\l t\r^{-\frac{3}{4}+\delta}\Delta_0^\f12)_{L_\omega^4}$ appeared in the above can be replaced by $O(\l t\r^{-\frac{3}{4}+\delta}\Delta_0)_{L_\omega^4}$.

Note $\Omega\bN=c\thetac(\Omega), S\bN_A=\tir\zb, X\bN^i\c \bN_i=X\log c$. 
 We can write 
 \begin{align}\label{6.21.1.24}
\sn_X^l(X\bN)&=\sn_X^{l}(c\thetac(\Omega))+\sn_X^{l}(\tir\zb) +X^l X\log c.
\end{align}
The lower order estimates in (\ref{6.30.1.23}) follows by using Lemma \ref{5.13.11.21} (1), (\ref{L2conndrv}), (\ref{L2BA2}) and (\ref{1.25.2.22}). 
  For the $n=4$ case in Lemma \ref{6.30.4.23}, we have to use the case that $l=2$ in (\ref{6.30.1.23}). To achieve these estimates, we have to use the assumption on $\sn_X^2\Omega$ in Lemma \ref{6.30.4.23}. Using (\ref{9.8.2.22}), (\ref{3.6.2.21}),  (\ref{3.11.3.21}), (\ref{L2BA2}) and (\ref{L2conndrv}) we can obtain the higher order estimates in (\ref{6.30.1.23}).

We next write with the help of (\ref{5.23.1.23}) that
\begin{align*}
	\Sc(X^n v)&=[X^n v]+\sum_{a=0}^{n-2}\sn_X^{a}(X_1 v^\|)\sn_X^{n-1-a}[\bN\cdots\bN]+\sum_{a=0}^{n-2}\sn_X^{n-1-a}[\bN\cdots\bN]^* X^a[X_1 v], \\
\Ac(X_n\cdots X_1 v)&=(X^n v)^\|+ \sum_{a=0}^{n-2}\sn_X^{a}(X_1 v^\|)\sn_X^{n-1-a}[\bN\cdots\bN]^*+\sum_{a=0}^{n-2}\sn_X^{n-1-a}[\bN\cdots\bN] X^a[X_1 v]
\end{align*}
where the union of the vector fields in $X^{n-1-a}$ and $X^a$ are $X^{n-1}=X_n\cdots X_2$.  
Using the above formula,  $n=2,3$ case in (\ref{1.10.1.23}) follows by using (\ref{6.30.1.23}) and (\ref{3.6.2.21}). The $n=4$ case follows in the same way.  
\end{proof}

\subsection{Decay of the transversal derivatives of the connection coefficients}\label{trans_control}
In this subsection, under the assumptions (\ref{3.12.1.21})-(\ref{1.25.1.22}), we control higher order transversal derivatives of $\Phi$ and of connection coefficients. We first give the following result which is on the mixed derivatives of $\Phi$ by tangential and transversal directions.
\begin{lemma}\label{7.6.6.23}
Let (\ref{3.12.1.21})-(\ref{1.25.1.22}) hold. If $\tir^{-1}\sn_X^2 \Omega=O(1)+O(\l t\r^{\delta-\frac{3}{4}}\Delta_0^{\f12+\frac{1}{p}})_{L_u^p L_\omega^4}$,  there hold for $X\in \{\Omega, S\}$ 
\begin{align*}
X^3([\Lb \Phi], k_{\bN\bN})&=\vs^-(X^{3})(O(1)[\Lb \Phi]+O(\l t\r^{-\f12+\delta}\Delta_0)_{L_\Sigma^2})\\
&+(1-\vs^-(X^{3}))O(\l t\r^{\delta}\Delta_0)_{L^2_\Sigma},\\
X^3\Lb v&=O(1)\fB+O(1)\vs^-(X^3)O(\l t\r^{-\f12+\delta}\Delta_0)_{L_\Sigma^2}\\
&+(1-\vs^-(X^3))O(\l t\r^{\delta}\Delta_0)_{L^2_\Sigma},
\end{align*}
\end{lemma}
\begin{proof}
If $X^3=X^2 \Omega$, we apply (\ref{LbBA2}), (\ref{3.6.2.21}) and the assumption on $\tir^{-1}\sn_X^2 \Omega$ to obtain the first inequality in Lemma \ref{7.6.6.23}. For the case that $X^3=X^2 S$, by letting $l=2$ in (\ref{8.26.1.23}), the first estimate in Lemma \ref{7.6.6.23} can be obtained by using (\ref{8.23.1.23}), (\ref{L2conndrv}), (\ref{L2BA2}), (\ref{3.6.2.21}), Lemma \ref{5.13.11.21} (5) and (\ref{3.29.1.23}).
For the second estimate in this lemma,  using the first formula in Lemma \ref{6.30.4.23} with $n=4$,  $\Lb v^\|=[\sn \Phi],$ (\ref{L2BA2}), Lemma \ref{5.13.11.21} (1), we infer with $X_1=\Lb$ 
\begin{align*}
[X^4 v]&=X^3[\Lb v]+O(1)\sn_X^{\le 2}[\sn \Phi]+O(\l t\r^{-\frac{3}{4}+\delta}\Delta_0)_{L_u^2 L_\omega^4}\fB\\
&+O(\l t\r^{-\frac{3}{4}+\delta}\Delta_0^{1-\f12
\vs^+(X_4\cdots X_2)})X^{1+\le 1}\fB\\
&+O(\l t\r^{-\frac{3}{4}+\delta}\Delta_0^\f12)_{L_\omega^4}\sn_X^{\le 1}([\sn\Phi])+O(\l t\r^{-\f12+\delta}\Delta_0)_{L_u^2 L_\omega^2}[\sn\Phi]\\
&=X^3[\Lb v]+O(\l t\r^{-\frac{1}{2}+\delta}\Delta_0)_{L^2_\Sigma}+O(\l t\r^{-\frac{3}{4}+\delta}\Delta_0^\f12)\vs^-(X_4\cdots X_2)\fB. 
\end{align*}
 \begin{align*}
(X^4 v)^\|&=\sn_{X_4}\cdots \sn_{X_2} [\sn \Phi]+O(1)\sn_X^{\le 2}\fB+O(\l t\r^{-\frac{3}{4}+\delta}\Delta_0^\f12)_{L_\omega^4}[\sn \Phi]\\
&+O(\l t\r^{-\frac{3}{4}+\delta}\Delta_0^{1-\f12 \vs^+(X_4\cdots X_2)})X^{1+\le 1}[\sn \Phi]\\
&+O(\l t\r^{-\frac{3}{4}+\delta}\Delta_0^{\f12+\frac{1}{p}})_{L_u^p L_\omega^4}\sn_X^{\le 1}\fB+O(\l t\r^{-\f12+\delta}\Delta_0)_{L_u^2 L_\omega^2}\fB\\
&=\vs^-(X_4 \cdots X_2)\fB+(1-\vs^-(X_4\cdots X_2))O(\l t\r^\delta\Delta_0)_{L^2_\Sigma}+O(\l t\r^{-\f12+\delta}\Delta_0)_{L^2_\Sigma}.
\end{align*} 
Substituting the first estimate in Lemma \ref{7.6.6.23} to the bound for $[X^3 \Lb v]$ in the above, combining the above two estimates, we conclude the second estimate in Lemma \ref{7.6.6.23}.  
\end{proof}

Next we will consider the transversal derivatives of $\Lb(\tr\chi-\frac{2}{\tir})$. Related to it, we need to control $\tir\mho$, which is not necessarily close to $0$ for large $t$ when $\bb$ is no longer close to $c^{-1}$. We first prove the initial bound of $\mho$, which relies on a delicate cancelation.  
\begin{lemma}
As a corollary of Proposition \ref{12.21.1.21}, 
\begin{equation}\label{mhoini}
\bb\tir \mho(0)=O(\La_0^\f12),\quad O(\La_0)_{L^2_u L_\omega^\infty}.
\end{equation}
\end{lemma}
\begin{proof}
At $t=0$, due to $\bb^{-1}=c$,
\begin{align*}
\Lb \tir-1&=-2\bN (\ckk c^{-1} u)=\ckk c^{-2}2\bN(\ckk c) u-\ckk c^{-1} 2\bN u\\
&=2\tir\bN\log (\ckk c)(0,u,\omega) -2\ckk c^{-1}c.
\end{align*}
Hence
\begin{equation}\label{7.11.5.22}
\bb \Lb \tir(t,u,\omega)=2 \tir c^{-1}\bN \log \ckk c(0,u,\omega)-2\ckk c^{-1}+c^{-1}(0,u,\omega).
\end{equation}
 Note $|c-\ckk c|\les \La_0$ by \Poincare inequality. We have  at $t=0$,
\begin{equation}\label{7.11.1.22}
|\ckk c^{-1} c-1|\les \La_0.
\end{equation}
Next we claim
\begin{equation}\label{7.11.3.22}
c^{-1}\tir (\bN\log \ckk c- \bN \log c)(0,u,\omega)=O(\La_0).
\end{equation}
Indeed, note $\bN\log c=(\wp-1)\bN \varrho$. Using (\ref{9.25.1.22}), we have $\bN\log c, \tr\thetac=O(1)$ at $t=0$.
Using these estimates, at $\Sigma_0$, by applying 
\begin{equation*}
\bb \bN \bar f=\overline{\bb \bN f}+\overline{\bb c\tr\thetac(f-\bar f)}
\end{equation*}
to $f=c$ and using (\ref{7.11.1.22}), we obtain at $t=0$ that
\begin{equation*}
 \overline{\bN \log c}- \bN \log \ckk c=O(\La_0).
\end{equation*}
   Applying \Poincare inequality to have $\Osc(\bN\log c)=O(\La_0)$ due to (\ref{12.22.1.21}), and noting $c^{-1}\tir=O(1)$  at $t=0$, we obtain (\ref{7.11.3.22}).
Using (\ref{7.11.5.22})-(\ref{7.11.3.22}), we infer
\begin{equation*}
\bb \Lb \tir=2c^{-1}\tir \bN\log c -c^{-1}+O(\La_0), \quad \mbox{ at }t=0.
\end{equation*}

Note also by using (\ref{k1}),
\begin{equation}\label{8.12.1.22}
\tr\chi+\tr\chib=-4 \bT\log c+2\p_{\hat A} v_{\hat B} \delta^{\hat A\hat B}.
\end{equation}
At $t=0$
\begin{align*}
\bb\Lb \tir-\f12\bb\tir \tr\chib&=2\bb\tir \bN \log c-c^{-1}-\f12 \tir\bb \tr\chib+O(\La_0)\\
&=2\bb\tir(L \log c+\f12 \tr\eta)+\f12\tir\bb \tr\chi-c^{-1}+O(\La_0).
\end{align*}
Thus at $t=0$, it follows by using (\ref{9.18.1.23}), (\ref{9.30.1.23}) and $\tr\eta=[L\Phi]$ that
\begin{equation*}
\bb (\Lb \tir-\f12\tir \tr\chib)=O(\La_0^\f12); \quad \bb(\Lb\tir-\f12\tir \tr\chib)=O(\La_0)_{L^2_u L_\omega^\infty}, 
\end{equation*}
as stated.
\end{proof}
In the sequel, we control $\mho$ and the transversal derivatives of connection coefficients. 
We will show that  $\tir \mho$ is $O(1)$ for $t>0$ and verifies various derivative estimates.
\begin{proposition}\label{7.22.2.22}
Let (\ref{3.12.1.21})-(\ref{1.25.1.22}) hold. With $X, X_1, X_2\in\{\Omega, S\}$, there hold for $0<t<T_*$ and $u_0\le u\le u_*$
\begin{align}
&\begin{array}{lll}
\|\mu- \tr\chi k_{\bN\bN}-\varpi\|_{L_\omega^4}\les \l t\r^{-2+\delta}\log\l t\r\Delta_0\\
\|\Omega^l(\ud\mu-k_{\bN\bN}\tr\chib-\varpi)\|_{L_\omega^4}\les \Delta_0\l t\r^{-3+\frac{l}{4}+\delta},\, l=0,1
\end{array}\label{3.26.2.21}\\
\displaybreak[0]
&\begin{array}{lll}
\|X^{\le 1}(\mu- \tr\chi k_{\bN\bN}-\varpi)\|_{L^2_u L_\omega^2}\les \l t\r^{-2+\delta}\Delta_0\log \l t\r,\\
\|X^{\le 1}(\ud\mu-k_{\bN\bN}\tr\chib-\varpi)\|_{L_\Sigma^2}\les \Delta_0\l t\r^{-2+\delta}
\end{array}
\label{9.8.1.22}\\
\displaybreak[0]
&\|\tir^{2-2\delta-}\sn_\Omega^{\le 1}(\mu- \tr\chi k_{\bN\bN}-\varpi)\|_{L_u^2 L_t^\infty L_\omega^2}\les\Delta_0\label{1.7.1.23}\\
&\|\tir\sn_\Lb \bA_g\|_{L_\omega^4}+\|\bb^{-\f12}\sn_X^{\le 1} \sn_\Lb \bA_{g}\|_{L^2_\Sigma}\les \l t\r^{-1+2\delta}\Delta_0\label{7.22.1.22}\\
& \Lb\bA_b-O(\tir^{-1})(\fB+\mho)-O(\l t\r^{-1})\bA_b=\left\{\begin{array}{lll}
O(\l t\r^{-2+\delta}\log \l t\r\Delta_0)_{L_\omega^4}\\
 O(\l t\r^{-2+\delta}\log \l t\r\Delta_0)_{L^2_u L_\omega^2}
 \end{array}\right.
 \label{8.3.5.23}\\
&\sn_X^l(\Lb \bA_b)=\sn_X^l\big(\tir^{-1}( k_{\bN\bN}+\mho)+\varpi\big)+O(\l t\r^{-2+\delta}\log \l t\r\Delta_0)_{L^2_u L_\omega^2},\quad l=0,1\label{8.12.3.22}\\
&\bb(\Lb \tir-1)(t, u, \omega)=c^{-1}(0)(\Lb \tir(0)-1)(u, \omega)\label{7.12.1.22}\\
&\tir\mho-2=O(\bb^{-1})\label{11.13.1.23}\\
&\|\bb^{-\f12}\Omega\mho\|_{L^2_\Sigma}\les \Delta_0\log\l t\r,\, \|\Omega \mho\|_{L_\omega^4}\les \Delta_0\l t\r^{-1}\log \l t\r \label{9.12.3.22}\\
&S\mho-O(1)(\mho+\fB)=O(\l t\r^{-2+\delta}\Delta_0)_{L_\omega^4}, O(\l t\r^{-2+\delta}\Delta_0)_{L_u^2 L_\omega^2}\label{8.29.3.23}\\
&X_2 X_1 \mho=\vs^-(X^2)O(1)(\mho+\fB)+O(\l t\r^{-1-\vs^-(X^2)+\delta}\Delta_0)_{L^2_u L_\omega^2}.\label{8.13.7.22}
\end{align}
\end{proposition}
\begin{proof} 

 (\ref{3.26.2.21}) can be obtained by applying  Lemma \ref{5.13.11.21} (1) and (\ref{1.27.5.24})  to (\ref{8.13.4.22}). And (\ref{9.8.1.22}) can be obtained similarly by, in addition, using (\ref{L2conndrv}), (\ref{zeh}), (\ref{4.17.1.24}), (\ref{3.16.1.22}) and (\ref{L2BA2}).

For simplicity, we only present the proof for the first order estimate  in (\ref{1.7.1.23}). The lower order estimate follows similarly and in an easier way.   
Using (\ref{lb}) and (\ref{tran1}), we derive
\begin{align*}
L\sD \log \bb +  \tr\chi \sD \log \bb&= -\sD k_{\bN\bN}-2\chih\c \sn^2 \log \bb-\sn_A\chi_{AC}\sn_C \log \bb\nn
\\&+(\tr\chi \zb_C-\chi_{AC}\zb_{A}-\delta^{AB}\bR_{CA4B}) \sn_C \log \bb.
\end{align*}
Applying (\ref{cmu_2}) to the above transport equation gives
\begin{align*}
L(\tir \sn)(\tir^2 \sD\log \bb)&=[L, \tir\sn](\tir^2 \sD \log \bb)+\tir^3\sn\Big(\bA\c \sn^2\log \bb+\sD k_{\bN\bN}\\
&+(\sn\bA+\chi \zb+\bR_{CA4B})\sn\log \bb\Big).
\end{align*}
Due to (\ref{4.5.4.24}) and (\ref{zeh})
\begin{align*}
\|[L, \tir\sn](\tir^2 \sD \log \bb)\|_{L^2_u L_\omega^2}&\les \l t\r^{-\frac{7}{4}+\delta}\Delta_0^\f12 \|\tir^2(\tir\sn)^{\le 1}\sD\log \bb\|_{L_u^2 L_\omega^2}\les \l t\r^{-\frac{7}{4}+2\delta}\Delta_0^\frac{3}{2}.
\end{align*}
Moreover, due to (\ref{zeh}), (\ref{L2conndrv}), (\ref{L2BA2}), (\ref{3.6.2.21}) and (\ref{3.11.3.21})
\begin{equation*}
\tir^3\sn\Big(\bA\c \sn^2\log \bb, (\sn \bA, \chi\zb)\sn\log \bb\Big)=O(\l t\r^{-\frac{3}{2}+2\delta}\Delta_0^\frac{3}{2})_{L_u^2 L_\omega^2}
\end{equation*}
Using the above two estimates,  (3) in Lemma \ref{5.13.11.21} and (\ref{LbBA2}), we have
\begin{align*}
\|\sn_L\Big(\tir \sn(\tir^2 \sD\log \bb)\Big)\|_{L_u^2 L_\omega^2}\les \l t\r^{-1+\delta}\Delta_0.
\end{align*}
Noting $\ze=\sn\log \bb+k_{A\bN}$ and using (\ref{dze}), (2) in Lemma \ref{5.13.11.21}, (\ref{L2BA2}), (\ref{L2conndrv}), (\ref{3.16.1.22}), (\ref{3.6.2.21}), (\ref{3.11.3.21}), and (\ref{4.17.1.24})  we derive
\begin{align*}
\sn_L(\tir^2(\tir \sn)(\mu-\tr\chi k_{\bN\bN}-\varpi)&=\sn_L\Big(\tir \sn(\tir^2 \sD\log \bb)\Big)+\sn_L\Big(\tir^2 (\tir \sn)(\sn \bA_{g,1}+\ud\bA^2+\bA_g^2)\Big)\\
&=O(\l t\r^{-1+2\delta}\Delta_0)_{L_u^2 L_\omega^2}+\sn_L\Big(\tir \sn(\tir^2 \sD\log \bb)\Big).
\end{align*} 
Integrating the above estimate along null cones gives
\begin{equation*}
\|\tir^{2-2\delta-} (\tir \sn)(\mu-\tr\chi k_{\bN\bN}-\varpi)\|_{L_u^2 L_t^\infty L_\omega^2}\les \Delta_0.  
\end{equation*}
Thus in view of (\ref{9.8.2.22}) the proof of (\ref{1.7.1.23}) is completed.

Next we show 
\begin{equation}
\tir\|\sn_\Lb \chih\|_{L_\omega^4}+\|\bb^{-\f12}\sn_X^{\le 1}\sn_\Lb \chih\|_{L^2_\Sigma}\les \Delta_0\l t\r^{-1+2\delta}.\label{1.28.4.22}
\end{equation}

To see (\ref{1.28.4.22}), in view of (\ref{3chi}), we symbolically write that
\begin{equation}\label{8.22.4.23}
\sn_X^l\sn_\Lb \chih=\sn_X^l\{(\tir^{-1}+\fB+[L\Phi]+\bA_b)(\chih+\eh)+\sn\hot\ze+\ze\hot\ze+\widehat{\bR_{A43B}}\}
\end{equation}
where $l=0,1$ and $X=S, \Omega$. 

Using $\tir|\fB|\les 1$ and (\ref{3.6.2.21}), with $l=0$ in the above, we bound
\begin{align*}
|\sn_\Lb \chih|&\les \l t\r^{-1}(|\chih|+|\eh|)+|\sn \zeta|+|\zeta\c \zeta|+|\widehat{\bR_{A43B}}|.
\end{align*}
By using (\ref{3.6.2.21}), (\ref{3.11.3.21}) and (\ref{1.29.4.22}), we obtain the first estimate in (\ref{1.28.4.22}).

We now prove the second estimate in (\ref{1.28.4.22}). It suffices to consider $l=1$. With $l=1$ in (\ref{8.22.4.23}),
 using (1) in Lemma \ref{5.13.11.21}, (\ref{L2BA2}), (\ref{3.16.1.22}) and (\ref{10.15.2.22}), we obtain
 \begin{align*}
\|\bb^{-\f12}\sn_X\sn_\Lb \chih\|_{L^2_\Sigma}&\les \l t\r^{-1}\|\bb^{-\f12}\sn_X^{\le 1} (\chih+\eh)\|_{L^2_\Sigma}+\Delta_0^\f12\l t\r^{-1+\delta}\|\tir(\chih+\eh)\|_{L_u^2 L_\omega^4}\\
&+\|\bb^{-\f12}\sn_X (\sn\ze, \ze\c \ze, \widehat{\bR_{A43B}})\|_{L^2_\Sigma}\\
&\les \l t\r^{-1+2\delta}\Delta_0.
\end{align*}
This proves the second estimate in (\ref{1.28.4.22}). The remaining estimates in (\ref{7.22.1.22}) is from (\ref{L4BA1}) and (\ref{LbBA2}).

To see the $\Lb \bA_b$ estimate in (\ref{8.3.5.23}), recall from (\ref{8.11.9.22}) that
\begin{equation}\label{10.6.1.23}
\Lb \bA_b=2\tir^{-1}(k_{\bN\bN}+\mho)+\varpi+\mu-\tr\chi k_{\bN\bN}-\varpi+(k_{\bN\bN}-\f12\tr\chib) \bA_b.
\end{equation}
The estimate of $\Lb \bA_b$ follows by using $k_{\bN\bN}, \tr\chib=O(\tir^{-1})$, $|\varpi|\les \l t\r^{-1}|\fB|$, (\ref{3.26.2.21}), (\ref{9.8.1.22}) and (\ref{3.6.2.21}). 

Differentiating (\ref{10.6.1.23}), we use (\ref{8.23.2.23}) and $\tr\chib=-\tr\chi+\fB$ from (\ref{8.12.1.22}) to derive
\begin{align}\label{8.23.4.23}
\begin{split}
\sn_X \Lb (\tr\chi-\frac{2}{\tir})&=\sn_X\big(\tir^{-1} (k_{\bN\bN}+\mho)+\varpi\big)+\sn_X(\mu-k_{\bN\bN}\tr\chi-\varpi)\\
&+O(\l t\r^{-1})\sn_X^{\le 1}\bA_b+O(\l t\r^{-1+\delta}\Delta_0)\bA_b.
\end{split}
\end{align}
Using (\ref{9.8.1.22}) and (\ref{L2conndrv}), we then obtain (\ref{8.12.3.22}) for $\bA_b=\tr\chi-\frac{2}{\tir}$. Due to (\ref{5.17.1.21}), (\ref{6.7.4.23}) and $\tr\eta=[L\Phi]$, the estimate for $\sn_X \Lb (\tr\thetac-\frac{2}{\ckr})$ follows  from using the above estimate.
Hence the proof of (\ref{8.12.3.22}) is complete.

(\ref{7.12.1.22}) follows immediately from (\ref{8.11.8.22}). By the definition of $\mho$ and using (\ref{7.12.1.22}), we write
\begin{align*}
\tir \mho=\Lb \tir-1-\f12\tir\tr\chib+1=O(\bb^{-1})+1+\f12\tir(\tr\chi+2\tr k)=2+\tr k+O(\bb^{-1})+\tir \bA_b
\end{align*} 
(\ref{11.13.1.23}) follows by using (\ref{k1}), (\ref{3.6.2.21}) and (\ref{3.11.4.21}). 

Next consider (\ref{9.12.3.22}) and (\ref{8.29.3.23}). 
In view of  (\ref{8.11.7.22}) and (\ref{8.13.3.22}), we write schematically
\begin{align}\label{9.8.4.22}
\begin{split}
\Omega \mho&=(\tir^{-1}+\mho+\fB+\bA_b) \ud\bA\c \Omega+\Omega(\bA_b+\fB)\\
S\mho&=\tir(\ud\mu-k_{\bN\bN}\tr\chib-\varpi)+(\tir k_{\bN\bN}-1) \mho-k_{\bN\bN}+\tir\varpi+\tir \tr\chib \bA_b.
\end{split}
\end{align}
Then (\ref{9.12.3.22}) can be obtained by using (\ref{11.13.1.23}), Lemma \ref{5.13.11.21} (1) and (5), (\ref{L2conndrv}). Again noting that  $|\varpi|\les \l t\r^{-1}|\fB|$, (\ref{8.29.3.23}) can be obtained by the lower order estimate in (\ref{3.26.2.21}), (\ref{9.8.1.22}) and Lemma \ref{5.13.11.21} (1). 

At last we consider (\ref{8.13.7.22}).
Using (\ref{9.8.4.22}), we derive
\begin{align*}
X\Omega\mho&=X(\tir^{-1}+\mho+\fB+\bA_b)\ud\bA\c\Omega+(\tir^{-1}+\mho+\fB+\bA_b)\sn_X(\ud\bA\c \Omega)+X\Omega(\bA_b+\fB)\\
XS\mho&=X(k_{\bN\bN}\tir)\mho+(k_{\bN\bN}\tir -1) X\mho+X (k_{\bN\bN}+\tir \varpi)+X\left(\tir(\ud\mu-\tr\chib k_{\bN\bN}-\varpi)\right)+X(\tir \tr\chib \bA_b).
\end{align*}
Hence we obtain by using (\ref{8.23.2.23}), (\ref{3.6.2.21})  and (\ref{11.13.1.23})-(\ref{8.29.3.23})
\begin{align*}
\|\bb^{-\f12}X\Omega \mho\|_{L^2_\Sigma}&\les \|\bb^{-\f12}\sn_X^{\le 1}\ud\bA\|_{L^2_\Sigma}+\|\bb^{-\f12}X\Omega \bA_b\|_{L^2_\Sigma}+\|\bb^{-\f12}X\Omega \fB\|_{L^2_\Sigma}\\
&+\l t\r^{-1+\delta}\Delta_0^\f12\|\tir\ud\bA\|_{L_u^2 L_\omega^4};\\
XS\mho&=X(k_{\bN\bN}\tir)\mho+O(1)X\mho+X(k_{\bN\bN}+\tir \varpi)+X(\tir(\ud\mu-\tr\chib k_{\bN\bN}-\varpi))\\
&+(O(\Delta_0^\f12\log \l t\r)_{L_\omega^4}+X^{\le 1})\bA_b.
\end{align*}
Using (\ref{9.8.2.22}), (\ref{9.8.1.22}), (\ref{9.12.3.22}), (\ref{8.29.3.23}), (\ref{3.11.3.21}), (\ref{8.23.2.23}), (\ref{1.27.5.24}), (\ref{3.16.1.22}), (\ref{LbBA2}) and (\ref{L2conndrv}), we obtain
\begin{align*}
&\|\bb^{-\f12}X\Omega\mho\|_{L^2_\Sigma}\les \Delta_0\l t\r^{\delta}\\
& XS\mho=\vs(X)(\mho+k_{\bN\bN})+(1-\vs(X))O(\Delta_0 \l t\r^{-1+\delta})_{L_u^2 L_\omega^2}+O(\l t\r^{-2+\delta}\Delta_0)_{L_u^2L_\omega^2}.
\end{align*}
  Thus (\ref{8.13.7.22}) is proved.
\end{proof}

Next we control $\Lb\tir$ and its derivatives.  We do not expect to be close to the Minkowskian value that $\Lb\tir=-1$. The following estimates suffice for our purpose. 
\begin{proposition}\label{1.8.1.22}
Let (\ref{3.12.1.21})-(\ref{1.25.1.22}) hold. We have the following estimates
\begin{align}
&\Lb \tir-1=\bb^{-1}O(1)\label{10.9.5.22}\\
&L \Lb \tir=k_{\bN\bN}(\Lb \tir-1)\label{7.16.3.22}\\
&\sn L \Lb \tir=\{\sn k_{\bN\bN}-\tir^{-1}(\zeta-k_{A\bN})\}O(1)\label{7.20.1.22}\\
&S^2\Lb \tir=O(1)S^{\le 1}(\tir k_{\bN\bN})\label{9.4.1.23}\\
&XL \Lb \tir-\vs(X)\fB=(1-\vs(X))O(\l t\r^\delta\Delta_0)_{L^2_\Sigma}\label{9.26.1.23}
\end{align}
where $X\in\{\Omega, S\}$ in (\ref{9.26.1.23}).
\end{proposition}
\begin{remark}
We will frequently use (\ref{10.9.5.22}) together with various estimates, such as (\ref{8.23.1.23'}) and the second estimate in (\ref{1.25.2.22}). Together with $L\tir=1$, we use it to have $\bN \tir=O(1)$, which will be frequently used as well.    
\end{remark}
\begin{proof}
(\ref{10.9.5.22}) is a direct consequence of (\ref{7.12.1.22}). (\ref{7.16.3.22}) follows from (\ref{3.19.2}) immediately.


Differentiating (\ref{7.16.3.22}) by using (\ref{7.03.1.19}) yields
\begin{align*}
\sn L\Lb \tir&=\sn k_{\bN\bN}(\Lb \tir-1)+k_{\bN\bN}\sn\Lb \tir\\
&=\sn k_{\bN\bN}(\Lb \tir-1)+k_{\bN\bN}(\zeta-k_{A\bN})\c(1-\Lb\tir)
\end{align*}
which gives (\ref{7.20.1.22}).

By using (\ref{7.16.3.22}) and (\ref{6.22.1.21}), we calculate
\begin{align*}
S^2 \Lb \tir&=S\big(\tir k_{\bN\bN}(\Lb \tir-1)\big)=S(\tir k_{\bN\bN})(\Lb \tir -1)+\tir k_{\bN\bN}S \Lb \tir\\
&=(S(\tir k_{\bN\bN})+\tir^2 k^2_{\bN\bN})(\Lb \tir-1).
\end{align*}
This gives (\ref{9.4.1.23}) in view of (\ref{10.9.5.22}).

By using (\ref{8.23.1.23}), (\ref{7.20.1.22}) and (\ref{9.4.1.23}), we derive 
\begin{align*}
XL\Lb \tir
&=\vs(X) \fB+(1-\vs(X))O(\l t\r^\delta\Delta_0)_{L^2_\Sigma}.
\end{align*}
This gives (\ref{9.26.1.23}).
\end{proof}
From now on, we frequently use (\ref{10.9.5.22}) without explicitly mentioning. Next, we prove (6) in Lemma \ref{5.13.11.21} by establishing the following crucial derivative estimates on $\Lb^2\varrho$. 
\begin{proposition}\label{6.27.1.24}
 Let (\ref{3.12.1.21})-(\ref{1.25.1.22}) hold. We  have the following estimates
\begin{align}
&\bb^3\Lb^2\varrho=O(\l t\r^{-1})_{L_\omega^4}, \quad\bb^3(\sn_\Lb(\Lb\Phi^\dagger), \Lb^2\Phi, \bT\bT\Phi)=O(\l t\r^{-1})_{L_\omega^4}\label{1.29.2.22}\\
&\bb\tir L\Lb\fB=O(\bb^{-2}\fB)_{L_\omega^4}\label{11.16.1.23}\\
&\|\tir\Omega(\bb^3 \Lb^2 \varrho)\|_{L_u^2 L_\omega^2}\les \Delta_0\label{1.27.3.24}
\end{align}
\end{proposition}
\begin{proof}
We first consider (\ref{1.29.2.22}) by using Lemma \ref{12.21.3.23}. Recall from (\ref{12.21.1.23}) that
\begin{equation*}
L(\bb\Lb(\tir \Lb \varrho))-\bb \Lb(\tir k_{\bN\bN}\Lb \varrho)=-\bb\Lb(\bA_b \tir \Lb \varrho)-\bb\Lb(\hb L \varrho \tir)+\bb F_1,
\end{equation*}
with 
\begin{align}\label{4.9.2.24}
F_1=\tir (\zb-\ze)\sn \Lb \varrho+\Lb(\zb\c \sn \varrho \tir)+\Lb(\tir \sD \varrho)+(\tir\Lb)^{\le 1} \N(\Phi, \bp\Phi) +k_{\bN\bN} L(\tir \Lb \varrho)
\end{align}
where we also used Proposition \ref{geonul_5.23_23} in the above; and the first two terms of $F_1$ can be written as $\tir \ud \bA \sn\Lb \varrho+\tir \bA_{g,1}\sn_\Lb \bA_{g,1}=O(\l t\r^{-2+2\delta}\Delta_0^2)$ due to (\ref{3.6.2.21}). 
  
We further derive due to the fact $\tr\chi+\tr\chib=-2\tr k=[L\Phi]+\Lb\varrho$ and $k_{\bN\bN}=\f12\wp\Lb\varrho+[L\Phi]$,
\begin{align}
L(\bb \Lb(\tir \Lb \varrho))&=2\bb\Lb(\Lb \varrho \tir)k_{\bN\bN}+\bb\Lb(\tir \Lb \varrho)(\bA_b+[L\Phi])+\bb \Lb\varrho(\fB \Lb \tir+\tir \Lb[L\Phi])+\Lb \bA_b\c \bb \tir \Lb \varrho\nn\\
&+\bb\Lb\big(\tir \tr\chi L\varrho\big)+\bb F_1.\label{2.10.1.24}
\end{align}
Using (\ref{8.3.5.23}), (\ref{3.6.2.21}), (\ref{3.11.3.21}) and (\ref{11.13.1.23}), we derive
\begin{equation}\label{1.12.1.23}
\Lb(\tir \tr\chi), \tir\Lb \bA_b=O(\l t\r^{-1})+O(\l t\r^{-1+\delta}\log \l t\r\Delta_0)_{L_\omega^4}. 
\end{equation}
Also using (\ref{6.22.1.21}) and $\bb \tir \fB=O(1)$, we have
\begin{align}\label{12.21.7.23}
L(\bb \Lb(\tir \Lb \varrho))-2\bb\Lb(\Lb \varrho \tir)k_{\bN\bN}&=\bb\Lb(\tir \Lb \varrho)(\bA_b+[L\Phi])+\bb F_1\\
&+O(\l t\r^{-2})+O(\l t\r^{-2+\delta}\log \l t\r\Delta_0)_{L_\omega^4}\nn.
\end{align}
Next we will show
\begin{equation}\label{12.21.5.23}
\bb F_1=[L\Phi]\c\bb\Lb( \tir  \Lb\varrho)+\bb \fB O(\tir^{-1})+  \bb O(\l t\r^{-2+2\delta}\Delta_0)_{L_\omega^4}.
\end{equation}
Indeed, in view of Proposition \ref{geonul_5.23_23} and $[\Lb\Phi]=\Lb \varrho+[L\Phi]$, we derive for $\Lb \Box_\bg \varrho$,
\begin{align}
\Lb \N(\Phi, \p \Phi)&=\Lb([\Lb \Phi][L\Phi])+\Lb([\sn \Phi][\sn\Phi]+|\eh|^2)\nn\\
&=\Lb (\Lb \varrho+[L\Phi])[L\Phi]+\fB\sn_\Lb[L\Phi]+[\sn \Phi]\sn_\Lb [\sn\Phi]+\sn_\Lb \eh\c\eh.\label{2.10.7.24}
\end{align}
Using (\ref{3.6.2.21}) and (\ref{6.22.1.21}), we derive
\begin{align}\label{12.21.4.23}
\Lb\N(\Phi, \p\Phi)=\Lb^2\varrho [L\Phi]+\fB O(\tir^{-1}\fB)+O(\l t\r^{-4+2\delta}\Delta_0^\frac{3}{2}).
\end{align}
Using (\ref{12.21.2.23}), $\bb \fB=O(\tir^{-1})$, (\ref{3.6.2.21}), (\ref{3.11.3.21}) and (3) in Lemma \ref{5.13.11.21}, we have
\begin{align*}
[\Lb, \sD]\varrho=O(\l t\r^{-3+2\delta}\Delta_0)_{L_\omega^4}. 
\end{align*}
Noting $\|\sD\Lb \varrho\|_{L_\omega^4}\les \l t\r^{-3+\delta}\Delta_0$ by (\ref{L4BA1}), 
\begin{equation}\label{7.13.1.24}
\|\Lb\sD\varrho\|_{L^4_\omega}\les \l t\r^{-3+2\delta}\Delta_0.
\end{equation}
Note the first two terms in (\ref{4.9.2.24}) can be included into the last error in (\ref{12.21.5.23}). Thus (\ref{12.21.5.23}) is proved. 

Integrating (\ref{12.21.7.23}) with the help of (\ref{12.21.5.23}) by using the transport lemma, (\ref{lb}) and (\ref{3.6.2.21}), we derive
\begin{align}\label{12.21.6.23}
\|\bb^3\Lb(\tir \Lb \varrho)\|_{L_\omega^4}\les\|\bb^3\Lb(\tir \Lb \varrho)\|_{L_\omega^4(S_{0,u})}+\int_0^t \l t\r^{-2+2\delta+} \les 1.
\end{align}
We thus obtain 
\begin{equation*}
\|\bb^3 \Lb^2 \varrho\|_{L_\omega^4}\les \l t\r^{-1} 
\end{equation*}
as stated in (\ref{1.29.2.22}).

Next we note
\begin{equation}\label{1.28.6.22}
|\Lb (\Lb \Phi)^\dagger, \Lb\Lb\Phi|\les |\Lb\Lb \varrho|+(\tir^{-1}+|\ud \bA|)|\fB|+\l t\r^{-2+\delta}\Delta_0.
\end{equation}
Indeed, we only need to consider $\Phi=v^i$. In view of Proposition \ref{6.7.1.23}
\begin{equation*}
[\Lb \Lb v]=(\Lb\log c+\Lb) [\Lb v]+\ud\bA\Lb v^\|,  \quad\Lb \Lb v^\|=(\Lb\log c+\sn_\Lb) [\sn \Phi]+[\Lb v]\ud\bA
\end{equation*}
Using $[\Lb v]=\Lb\varrho+[L\Phi]=\fB$, $\Lb v^\|=[\sn\Phi]$ and (\ref{6.22.1.21}), we have
\begin{equation*}
|\sn_\Lb (\Lb v^\dagger)|+|\Lb^2 v^i|\les |\Lb \Lb \varrho|+|\sn_\Lb [\sn \Phi]|+\l t\r^{-1}|\fB|+|\ud\bA||\fB|.
\end{equation*}
(\ref{1.28.6.22}) then follows by applying Lemma \ref{5.13.11.21} (1) to the right-hand side of the above inequality. 
Hence $$\bb^3(|\sn_\Lb (\Lb v^\dagger)|+|\Lb^2 v^i|)\les\l t\r^{-1}$$
as stated in (\ref{1.29.2.22}). The last estimate in (\ref{1.29.2.22}) can then be obtained by using Lemma \ref{5.13.11.21} (1) and the proved estimates in (\ref{1.29.2.22}).

(\ref{11.16.1.23}) follows as a consequence of (\ref{12.21.6.23}), (\ref{12.21.7.23}) and (\ref{12.21.5.23}). 

Next we consider (\ref{1.27.3.24}). We derive by using (\ref{2.10.1.24}), (\ref{lb}) and (\ref{cmu_2}) that
\begin{align}\label{2.10.3.24}
\sn_L \sn(\bb^3 \Lb(\tir \Lb \varrho))+\f12\tr\chi\sn(\bb^3\Lb(\tir \Lb \varrho))=-\chih\sn(\bb^3 \Lb(\tir \Lb \varrho))+\sn F_3
\end{align}
where
\begin{align*}
F_3&=\bb^3 \Lb (\tir \Lb \varrho)(\bA_b+[L\Phi])+\bb^3 \Lb \varrho(\fB \Lb \tir+\tir\Lb[L\Phi])+\bb^2\Lb \bA_b\c \bb\tir \Lb \varrho\\
&+\bb^3 \Lb(\tir \tr\chi L\varrho)+\bb^3 F_1.
\end{align*}
Next we prove 
\begin{equation}\label{2.10.2.24}
\|\tir \sn F_3\|_{L_u^2 L_\omega^2}\les \|\tir \sn\big(\bb^3 \Lb(\tir \Lb \varrho)\big)\|_{L^2_u L_\omega^2}\l t\r^{-\frac{7}{4}+\delta}\Delta_0^\f12+ \Delta_0 \l t\r^{-\frac{7}{4}+\delta}.
\end{equation}
(\ref{1.27.3.24}) will follow as a consequence of (\ref{2.10.2.24}) by applying the transport lemma to (\ref{2.10.3.24}). Using (\ref{12.21.6.23}) and (\ref{3.11.3.21}) gives
\begin{align*}
&\|\tir \sn(\bA_b+[L\Phi]) \bb^3 \Lb (\tir \Lb \varrho)\|_{L_u^2 L_\omega^2}\les \l t\r^{-\frac{7}{4}+\delta}\Delta_0.
\end{align*}
Using (\ref{4.9.1.24}), Lemma \ref{5.13.11.21} (5) and (\ref{6.7.4.23}), we obtain
\begin{align*}
\|\tir \sn\big(\bb^3\Lb \varrho&(\fB\Lb \tir+\tir \Lb[L\Phi])\big)\|_{L_u^2 L_\omega^2}\\
&\les \|\tir\sn(\bb \Lb \varrho)\|_{L_u^2 L_\omega^2}\l t\r^{-1}\log \l t\r+\|\tir \sn(\bb^2(\fB\Lb \tir+\tir\Lb[L\Phi]))\|_{L_u^2 L_\omega^2}\l t\r^{-1}\\
&\les  \l t\r^{-2+\delta}\log \l t\r^\frac{3}{2}\Delta_0.
\end{align*}
Note
\begin{align*}
\|\tir \sn(\bb^2 \Lb \bA_b\c \bb \tir \Lb \varrho)\|_{L_u^2 L_\omega^2}&\les \|\tir \sn(\bb^2\Lb \bA_b)\|_{L_u^2 L_\omega^2}+\|\bb^2 \Lb \bA_b\c \tir^2 \sn(\bb\Lb \varrho)\|_{L_u^2 L_\omega^2}.
\end{align*}
Using Proposition \ref{7.22.2.22}, (\ref{3.6.2.21}), (\ref{3.11.3.21}) and (\ref{LbBA2}), we derive 
\begin{equation}\label{2.10.4.24}
\|\tir \sn\Lb \bA_b\|_{L^2_u L_\omega^2}\les \l t\r^{-2+\delta}\log \l t\r\Delta_0, \Lb \bA_b=O(\l t\r^{-2})+O(\l t\r^{-2+\delta}\log \l t\r\Delta_0)_{L_\omega^4}. 
\end{equation}
We then use Lemma \ref{5.13.11.21} (5) and the above estimate to derive
\begin{align*}
&\|\tir \sn(\bb^2 \Lb \bA_b\c \bb \tir \Lb \varrho)\|_{L_u^2 L_\omega^2}\les \l t\r^{-2+\delta}\Delta_0\log \l t\r^3.
\end{align*}
Next we deduce
\begin{align*}
\tir \sn(\bb^3\Lb(\tir \tr\chi L\varrho))&=\tir\sn(\bb^3\Lb(\tir \bA_b) L\varrho)+\tir^2\sn(\bb^3\tr\chi\Lb L\varrho)\\
&=\tir \sn(\bb^3(\Lb \tir \bA_b+\tir \Lb \bA_b)L \varrho)+\tir^2 \sn(\bb^3\tr\chi \Lb L \varrho).
\end{align*}
It follows by using (\ref{4.9.1.24}), Lemma \ref{5.13.11.21} (5) and (\ref{3.11.3.21}) that
\begin{align*}
\|\tir\sn(\bb^3\Lb \tir \bA_b L\varrho)\|_{L_u^2 L_\omega^2}&\les \|\tir\sn(\bb^3 \Lb \tir)\|_{L_u^2 L_\omega^2} \l t\r^{-\frac{15}{4}+2\delta}\Delta_0+\|\bb^3\tir \sn(\bA_b L\varrho)\|_{L_u^2 L_\omega^2}\\
&\les \l t\r^{-\frac{15}{4}+3\delta}\Delta_0\log \l t\r^3.
\end{align*}
Using (\ref{2.10.4.24}), (\ref{8.3.5.23}), (\ref{L2BA2}), (\ref{3.6.2.21}) and (\ref{3.11.3.21}), we derive
\begin{align*}
\|\tir \sn(\bb^3 \tir \Lb \bA_b L \varrho)\|_{L_u^2 L_\omega^2}\les \l t\r^{-3+3\delta}\Delta_0\log \l t\r^4. 
\end{align*}
Using (\ref{6.7.4.23}), (\ref{1.27.5.24}) and (\ref{L2conndrv}), we obtain
\begin{align*}
\|\tir^2 \sn(\bb^3 \tr\chi \Lb L \varrho)\|_{L_u^2 L_\omega^2}&\les \l t\r\|\sn(\bb^3\tr\chi) \fB\|_{L_u^2 L_\omega^2}+\|\bb^3 \tir \sn(\Lb L\varrho)\|_{L_u^2 L_\omega^2}\\
&\les\log \l t\r^3 \Delta_0 \l t\r^{-2+\delta}. 
\end{align*}
In view of (\ref{3.6.2.21}), we summarize the above estimates as
\begin{equation}\label{2.10.5.24}
\|\tir \sn F_3\|_{L_u^2 L_\omega^2}\les \|\tir \sn\big(\bb^3 \Lb(\tir \Lb \varrho)\big)\|_{L^2_u L_\omega^2}\l t\r^{-\frac{7}{4}+\delta}\Delta_0^\f12+\|\tir \sn(\bb^3 F_1)\|_{L_u^2 L_\omega^2}+\l t\r^{-\frac{7}{4}+\delta}\Delta_0. 
\end{equation}
It remains to show
\begin{align}\label{2.10.6.24}
\|\tir \sn(\bb^3 F_1)\|_{L_u^2 L_\omega^2}\les\|\tir \sn\big(\bb^3 \Lb(\tir \Lb \varrho)\big)\|_{L^2_u L_\omega^2}\l t\r^{-2+\delta}\Delta_0^\f12+\l t\r^{-2+2\delta}\log \l t\r^3 \Delta_0. 
\end{align}
(\ref{2.10.2.24}) will follow as a consequence. 

We first derive by using Lemma \ref{5.13.11.21} (1), (3) and (5), (\ref{7.13.1.24}), (\ref{LbBA2}), (\ref{4.9.1.24}) and (\ref{7.03.1.19}) 
\begin{align*}
\tir \sn(\bb^3 \Lb (\bA_{g,1}^2\tir)), \tir\sn(\bb^3\Lb(\tir \sD\varrho)), \tir \sn(\bb^3 k_{\bN\bN}L(\tir \Lb\varrho))=O(\l t\r^{-2+2\delta}\log \l t\r^3\Delta_0)_{L_u^2 L_\omega^2}.
\end{align*}
Using (\ref{LbBA2}), (\ref{zeh}) and (\ref{3.6.2.21}), we have
\begin{align*}
\|\tir^2 \sn(\bb^3\ud\bA\sn \Lb \varrho)\|_{L^2_u L_\omega^2}\les \l t\r^{-2+2\delta}\log \l t\r^3\Delta_0^2
\end{align*}
and the following estimate in view of (\ref{1.27.5.24})
\begin{align*}
\|\tir^2 \sn(\bb^3 \Lb \N(\Phi, \bp\Phi))\|_{L^2_u L_\omega^2}&\les\|\tir^2 \bb^2\Lb \N(\Phi, \bp\Phi)\|_{L_u^2 L_\omega^4}\l t\r^{-1}\log \l t\r\Delta_0+\|\bb^3\tir^2\sn\Lb \N(\Phi, \bp\Phi)\|_{L_u^2 L_\omega^2}.  
\end{align*}
Using (\ref{12.21.4.23}) and (\ref{1.29.2.22}), the first term on the right-hand side is bounded by 
$$\l t\r^{-2}\log \l t\r\Delta_0+\l t\r^{-2+\delta}\log \l t\r^3\Delta_0^\frac{3}{2}.$$ 

Differentiating (\ref{2.10.7.24}) by $\sn$, by using (\ref{1.29.2.22}), Sobolev embedding on spheres, (\ref{L2BA2}), (\ref{LbBA2}), (\ref{3.6.2.21}), (\ref{3.11.3.21}), (\ref{1.27.4.24}), (\ref{1.27.5.24}), (\ref{6.22.1.21}) and (\ref{6.7.4.23}) we deduce 
\begin{align}
&\|\bb^3\tir\sn\Lb \N(\Phi,\bp\Phi)\|_{L_u^2 L_\omega^2}\nn\\
&\les\|\bb^3 \Omega(\Lb^2\varrho)\|_{L_u^2 L_\omega^2}\|[L\Phi]\|_{L_x^\infty}+\|\bb^3 \Lb^2\varrho\|_{L_u^\infty L_\omega^4}\|\Omega[L\Phi]\|_{L_u^2 L_\omega^4}+\|\bb^3\Omega(\Lb[L\Phi])\|_{L_u^2 L_\omega^2}\|[L\Phi]\|_{L^\infty_x}\nn\\
&+\|\Omega[L\Phi]\|_{L_u^2 L_\omega^2}\|\bb^3\Lb[L\Phi]\|_{L_x^\infty}+\|\bb\Omega\fB\|_{L_u^2 L_\omega^2}\|\bb^2\Lb[L\Phi]\|_{L^\infty_x}+\|\bb\fB\|_{L_x^\infty}\|\bb^2 \Omega(\Lb[L\Phi])\|_{L_u^2 L_\omega^2}\nn\\
&+\|\bb^3\Omega(\sn_\Lb[\sn\Phi])\|_{L_u^2 L_\omega^2}\|[\sn\Phi]\|_{L_x^\infty}+\|\sn_\Omega[\sn\Phi]\|_{L_u^2 L_\omega^4}\|\bb^3\sn_\Lb[\sn\Phi]\|_{L_\omega^4}\nn\\
&+\|\bb^3\sn_\Omega\eh\|_{L_u^2 L_\omega^4}\|\sn_\Lb \eh\|_{L_\omega^4}+\|\sn_\Omega\sn_\Lb \eh\|_{L_u^2 L_\omega^2}\|\bb^3\eh\|_{L^\infty_x}\nn\\
&\les\|\Omega(\bb^3 \Lb^2\varrho)\|_{L_u^2 L_\omega^2}\l t\r^{-2+\delta}\Delta_0^\f12 +\l t\r^{-3+\delta}\Delta_0(1+\l t\r^\delta\Delta_0^\f12)\log \l t\r^3.\label{2.10.9.24}
\end{align}

Moreover, we directly derive
\begin{align*}
\tir \sn(\bb^3 \Lb(\tir \Lb\varrho))=\tir^2 \sn(\bb^3\Lb^2\varrho)+\tir\sn(\bb^3 \Lb \tir \Lb\varrho).
\end{align*}
For the last term, we bound by using (\ref{4.9.1.24}), (\ref{1.27.4.24}) and (\ref{1.27.5.24}) that 
\begin{align*}
\|\tir \sn(\bb^3 \Lb\tir \Lb \varrho)\|_{L_u^2 L_\omega^2}&\les \|\tir \sn(\bb^3\Lb\varrho)\|_{L^2_u L_\omega^2}+\log \l t\r^2 \l t\r^{-1}\|\tir \sn \Lb \tir\|_{L_u^2 L_\omega^2}\\
&\les \log \l t\r^2\l t\r^{-1+\delta} \Delta_0.
\end{align*}
Thus we have the comparison estimate that
\begin{equation}\label{2.10.8.24}
\tir^2\sn(\bb^3\Lb^2\varrho)=\tir \sn(\bb^3 \Lb(\tir \Lb\varrho))+O(\log \l t\r^2\l t\r^{-1+\delta} \Delta_0)_{L^2_u L_\omega^2}.
\end{equation}
Substituting (\ref{2.10.8.24}) into (\ref{2.10.9.24}), we conclude (\ref{2.10.6.24}). The proof of (\ref{1.27.3.24}) is completed. 
\end{proof}

\subsection{Sobolev inequalities}
We give the trace inequalities and Sobolev inequalities under Assumption \ref{5.13.11.21+}.
\begin{lemma}[Sobolev embedding]
\begin{align}\label{1.31.1.24}
\|\bb^\a F\|_{L_\omega^4}\les \|\bb^\a (\tir \sn)F\|_{L_\omega^2}^\f12 \|\bb^\a F\|^\f12_{L_\omega^2}+\|\bb^\a F\|_{L^2_\omega}+|\a|\log \l t\r\Delta_0\|\bb^{\a-1}F\|_{L_\omega^2}, \a\in {\mathbb R}
\end{align}
\end{lemma}
\begin{proof}
We apply the standard Sobolev embedding to $|\bb^\a F|^2$ and (\ref{1.27.5.24}) to obtain 
\begin{align*}
\||\bb^\a F|^2\|_{L_\omega^2}&\les \|(\tir \sn)^{\le 1}(|\bb^\a F|^2)\|_{L_\omega^1}\\
&\les\|\bb^\a (\tir \sn)^{\le 1} F\|_{L_\omega^2}\|\bb^\a F\|_{L_\omega^2}+\|\tir \sn(\bb^{2\a}) |F|^2\|_{L_\omega^1}\\
&\les\|\bb^\a (\tir \sn)^{\le 1} F\|_{L_\omega^2}\|\bb^\a F\|_{L_\omega^2}+|\a|\|\tir\sn \bb\|_{L_\omega^4}\|\bb^{2\a-1} |F|^2\|_{L_\omega^\frac{4}{3}}\\
&\les \|\bb^\a(\tir \sn)^{\le 1}F\|_{L_\omega^2}\|\bb^\a F\|_{L_\omega^2}+|\a|\log \l t\r\Delta_0\|\bb^{\a-1} F\|_{L_\omega^2}\|\bb^\a F\|_{L_\omega^4}
\end{align*}
The consequence drops out by using Cauchy-Schwarz.
\end{proof}
\begin{lemma}[Trace inequalities]
For scalar functions $f$, there hold the following inequalities
\begin{align}
&\|f \tir^{-\f12-\a}\|_{L^2(\H_u^t)}\les (\int_{S_{0,u}}f^2 \big)^\f12+ F_1[f]^\f12(\H_u^t),\quad \mbox{ if } \a>0\label{1.14.1.22}\\
&\int_S \bb^{n} f^2 \tir d\mu_\omega+\int_0^t \int_S \bb^n f^2 d\mu_\omega dt\nn\\
&\qquad\qquad\qquad\les\int_{S_{0,u}} \bb^n f^2 \tir d\mu_\omega+\|\bb^\frac{n}{2}\tir(\sn_L f+h f)\|_{L_t^2 L_\omega^2}^2,\, n=0,1\label{1.13.1.21}
\end{align}
where $S$ stands for $S_{t,u}$ and $d\mu_\omega=d\mu_{{\mathbb S}^2}$.
\end{lemma}
\begin{proof}
We first prove (\ref{1.14.1.22}). Since $
L(f^2 v_t)=2v_t  (Lf+hf)f,
$
integrating along $\H_u^t$ gives
\begin{equation}\label{1.14.2.22}
\int_{S_{t,u}}f^2 v_td\mu_\omega-\int_{S_{0,u}}f^2 v_0d\mu_\omega=\int_{\H_u^t} 2 (Lf +hf ) fd\mu_\ga dt.
\end{equation}
Dividing (\ref{1.14.2.22}) by $\tir^{-1-2\a}$ gives
\begin{align*}
&\|\tir^{-\f12-\a} f\|_{L^2(\H_u^t)}^2\\
&\le \int_0^t \l t'\r^{-1-2\a}dt'\{\int_{S_{0,u}}f^2 v_0 d\mu_\omega+\int_{\H_u^{t'}}|\tir^\f12(L f+hf)| |f \tir^{-\f12}|d\mu_\ga dt''\}\\
&\les  \int_{S_{0,u}}f^2 v_0 d\mu_\omega+\|\tir^\f12 (Lf+hf)\|_{L^2(\H_u^t)}\|f\tir^{-\f12-\a}\|_{L^2(\H_u^t)}.
\end{align*}
Due to Cauchy-Schwarz inequality, we have
\begin{align*}
\|\tir^{-\f12-\a} f\|_{L^2(\H_u^t)}^2\les\int_{S_{0,u}}f^2 v_0 d\mu_{{\mathbb S}^2}+\|\tir^\f12 (Lf+hf)\|_{L^2(\H_u^t)}^2,
\end{align*}
as desired in (\ref{1.14.1.22}). 

Next we prove (\ref{1.13.1.21}). We first derive
\begin{align*}
\int_S \bb^n f^2 \tir d\mu_\omega-\int_{S_{0,u}}\bb^n f^2 \tir d\mu_\omega&=\int_0^t \int_S (2\tir \sn_L f \c f +f^2)\bb^n d\mu_\omega dt'-\int_S n\bb^n k_{\bN\bN} \tir f^2 d\mu_\omega\\
&=\int_0^t\int_S \{2\tir (\sn_L f+h f)f+(1-2\tir h-n\tir k_{\bN\bN}) f^2\}\bb^n d\mu_\omega.
\end{align*}
Hence
\begin{align*}
\int_S \bb^n f^2\tir d\mu_\omega&+\int_0^t\int_S \bb^n f^2 d\mu_\omega dt=\int_{S_{0,u}} f^2 \tir \bb^n d\mu_\omega\\
\displaybreak[0]
&+\int_0^t\int_S \{2\tir (\sn_L f+h f)f+(2(1-\tir h)-n\tir k_{\bN\bN}) f^2\}\bb^nd\mu_\omega dt.
\end{align*}
Using (\ref{3.6.2.21}) and $|\bb\tir k_{\bN\bN}|\les 1$, Cauchy-schwarz inequality, and by induction on $n$, we have
\begin{align*}
\int_S\bb^n f^2\tir d\mu_\omega+\int_0^t\int_S\bb^n f^2 d\mu_\omega dt&\les \int_{S_{0,u}} \bb^n f^2 \tir d\mu_\omega+\int_0^t\int_S \bb^n\{2\tir (\sn_L f+h f)f\}d\mu_\omega dt\\
&+n\int_0^t \int_S \bb^{n-1}f^2 d\mu_\omega dt\\
&\les  \int_{S_{0,u}}\bb^n f^2 \tir d\mu_\omega+\|\bb^{\frac{n}{2}}\tir(\sn_L f+h f)\|^2_{L^2_t L^2_\omega},
\end{align*}
where we used $\bb>\frac{\bb_0}{4}$ to conclude the last line.   
\end{proof}

\begin{lemma}[Sobolev inequalities]
For any smooth function $f$ and constants verifying the relation $\ga_0'+2\ga_2=2\ga$ and $\ga_0', \ga_2, \ga\ge 0$, we have, for all $  u\in [u_0, u_*]$ and $0<t< T_*$,
\begin{align}
\int_{S_{t, u}}|\tir^\ga f|^4 \tir^{-2}
&\les \int_{S_{0,u}}|\tir^{\ga} f|^4 \tir^{-2} \nn \\
& \quad \, + \int_{\H_u^t}\tir^{\ga_0'}|L (\tir^\ga f)|^2 \tir^{-2}  \cdot \sum_{l\le 1}\int_{\H_u^t}\tir^{2\ga_2}|\Omega^l  f|^2\tir^{-2}.\label{6.24.12.18}
\end{align}
For $S_{t,u}$-tangent tensor fields $F$, there hold
\begin{align}
\int_{S_{t, u}}|\tir^\ga F|^4 \tir^{-2}
&\les \int_{S_{t,u_*}}|\tir^\ga F|^4 \tir^{-2} + \int_{\Sigma_t\cap[u, u_*]}\tir^{\ga_0'}|\sn_\bN(\tir^\ga F)|^2 \tir^{-2} \nn\\
&\cdot \sum_{l\le 1}\int_{\Sigma_t\cap[u, u_*]}\tir^{2\ga_2}(|\sn_\Omega^l F|^2+\log \l t\r^2 \Delta_0^2 |\bb^{-1} F|^2)\tir^{-2}\label{7.21.3.21}\\
\displaybreak[0]
\|\tir F\|^2_{L_\omega^4}&\les\|\tir F\|^2_{L_\omega^4(S_{t, u_*})}+\left(\|\tir^{-1} F\|_{L^2_{\Sigma_t}\cap[u, u_*]}+\|\sn_\bN F\|_{L^2_{\Sigma_t}\cap[u, u_*]}\right)\nn\\
&\times \|\sn^{\le 1}_\Omega F, \Delta_0 \log \l t\r\bb^{-1}F\|_{L^2_{\Sigma_t}\cap[u, u_*]}.\label{9.20.1.23} 
\end{align}
\end{lemma}
\begin{proof}
To obtain (\ref{6.24.12.18}), we repeat the proof of \cite[Lemma 1 in Sect 2.1]{KWY} with the help of $|\Omega f|\approx \tir |\sn f|$ (due to (\ref{9.30.2.22})). Note that (\ref{10.9.5.22}) and $L\tir=1$ implies $\bN\tir=O(1)$. Thus (\ref{7.21.3.21}) can be proved in the same way as for (\ref{6.24.12.18}) by applying (\ref{11.12.5.23}) to $\p_u(\int_{S_{t,u}}|\tir^\ga F|^4 v_t^{-1})$ and using (\ref{1.27.5.24}). Applying (\ref{7.21.3.21}) with $\ga=\ga_2=1$ and $\ga_0'=0$ yields (\ref{9.20.1.23}).
\end{proof}

\section{Energy inequalities}\label{mul_1}
In this section, we derive the fundamental energy inequalities in Proposition \ref{10.10.3.22} and Proposition \ref{MA2} under Assumption \ref{5.13.11.21+}. To prove these results, we take advantage of the bound (\ref{11.11.2.23}) for $\Lb \varrho$.  
\subsection{Standard energy inequalities}
\begin{proposition}[Standard energy estimate]\label{10.10.3.22}
Assume (\ref{3.12.1.21})-(\ref{1.25.1.22}) hold. For $0< t_1< T_*$ and $u_0\le u_1\le u_*$ and scalar functions $\psi$, the following standard energy estimate holds with $\a\ge 0$ 
\begin{equation}\label{6.21.2.21}
\begin{split}
&\int_{\D_{u_1}^{t_1}}\{-\wp[\Lb \varrho]_{-}(|\Lb \psi|^2+|L \psi|^2)\aaa^{-\a}+\aaa^{-\a-1}\frac{\a}{\tir+3}(|\Lb\psi|^2+|\sn\psi|^2)\}\\
&+E_{-\a}[\psi](t_1)+F_{0,-\a}[\psi](\H_{u_1}^{t_1})\\
&\les\int_{\D_{u_1}^{t_1}} \aaa^{-\a}\left|\bT \psi \Box_\bg \psi\right|+E_{-\a}[\psi](0)+F_{0,-\a}[\psi](\H_{u_*}^{t_1}),
\end{split}
\end{equation}
where $\aaa=\log(\tir+3)$.
\end{proposition}
\begin{proof}
 Denote the deformation tensor of $\bT$ by $\piT$. It is direct to obtain
 \begin{equation}\label{7.13.1.22}
 \begin{matrix}
&\piT_{LL}=-2k_{\bN\bN}, & \piT_{\Lb \Lb}=-2 k_{\bN\bN}, & \piT_{L\Lb}= 2k_{\bN\bN}, \\
& \piT_{AB}=-2 k_{AB}, & \piT_{AL}=-2k_{A\bN}, & \piT_{A\Lb}=2k_{A\bN}.
\end{matrix}
\end{equation}
With $Q_{\ga\b}=\p_\ga\psi\p_\b\psi-\f12\bg_{\ga\b}\bd^\mu \psi\p_\mu\psi$ and $\P_\ga=Q_{\ga\b}\bT^\b$,
applying divergence theorem to $\bd^\ga\P_\ga$ in the region of $\D_{u_1}^{t_1}$, we obtain the standard energy identity
\begin{equation*}
\int_{\D_{u_1}^{t_1}}\bd^\b(\aaa^{-\a}Q_{\ga\b}[\psi]\bT^\ga)=E_{-\a}[\psi](0)+F_{0,-\a}[\psi](\H_{u_*}^{t_1})-\left(E_{-\a}[\psi](t_1)+F_{0,-\a}[\psi](\H_{u_1}^{t_1})\right)
\end{equation*}
which implies 
\begin{align*}
\int_{\D_{u_1}^{t_1}}&\aaa^{-\a}\{\f12 Q_{\ga\b}[\psi]\piT^{\ga\b}+\Box_\bg \psi\bT \psi+\a(\f12 L\log \aaa \P_{\Lb}+\f12\Lb \log \aaa \P_L)\}\\
&+E_{-\a}[\psi](t_1)+F_{0,-\a}[\psi](\H_{u_1}^{t_1})=E_{-\a}[\psi](0)+F_{0,-\a}[\psi](\H_{u_*}^{t_1}).
\end{align*}
 Hence 
\begin{align*}
 \int_{\D_{u_1}^{t_1}}&\frac{1}{8}\aaa^{-\a}(Q_{\Lb \Lb} \piT_{LL}+Q_{LL}\piT_{\Lb \Lb})+E_{-\a}[\psi](t_1)+F_{0,-\a}[\psi](\H_{u_1}^{t_1})\\
&=E_{-\a}[\psi](0)+F_{0,-\a}[\psi](\H_{u_*}^{t_1})+\int_{\D_{u_1}^{t_1}}-\aaa^{-\a}\a(\f12 L\log\aaa \P_{\Lb}+\f12\Lb \log \aaa \P_L)\\
&-\aaa^{-\a}(\frac{1}{4} Q_{L\Lb}\piT_{L\Lb}+\f12 Q_{A\nu}\piT^{A\nu}+\Box_\bg \psi \bT \psi).
\end{align*}
It is direct to check 
\begin{align*}
\P_L=\f12(|L\psi|^2+|\sn\psi|^2), \P_\Lb=\f12(|\Lb\psi|^2+|\sn\psi|^2).
\end{align*}
Combining (\ref{3.22.1.21}), (\ref{3.6.2.21}) and (\ref{11.11.2.23}), we infer for $k_{\bN\bN}$,
\begin{align}\label{4.12.1.24}
2k_{\bN\bN}&=\wp[\Lb \varrho]_{-}+\wp[\Lb \varrho]_++O(\l t\r^{-2+\delta}\Delta_0^\f12)=\wp[\Lb \varrho]_{-}+O(\l t\r^{-2+\delta}\Delta_0^\f12).
\end{align}
Substituting the above identity to the energy identity implies
\begin{align}\label{7.13.2.22}
&\int_{\D_{u_1}^{t_1}}\{-\frac{1}{8}\wp[\Lb \varrho]_{-}(|\Lb \psi|^2+|L \psi|^2)\aaa^{-\a}+\frac{1}{4}\aaa^{-\a-1}\frac{\a}{\tir+3}(|\Lb\psi|^2+|\sn\psi|^2)\}\nn\\
&+E_{-\a}[\psi](t_1)+F_{0,-\a}[\psi](\H_{u_1}^{t_1})\\
&\les \int_{\D_{u_1}^{t_1}}|k_{\bN \bN} Q_{L \Lb}, k^{AB}Q_{AB}|\aaa^{-\a}+\int_0^{t_1}\l t'\r^{-2+\delta} \Delta_0^\f12 E_{-\a}(t')dt'+\left|\int_{\D_{u_1}^{t_1}} \aaa^{-\a}\Box_\bg \psi\bT \psi\right|\nn\\
&+\a\int_{\D_{u_1}^{t_1}}\aaa^{-\a-1}\{(\tir+3)^{-1}|\Lb\tir|(|L\psi|^2+|\sn\psi|^2)\} +E_{-\a}[\psi](0)+F_{0,-\a}[\psi](\H_{u_*}^{t_1}),\nn
\end{align}
where we used (\ref{7.13.1.22}), the rough estimate $\int_{\Sigma_t}\aaa^{-\a} |Q_{AL}, Q_{A\Lb}|\les E_{-\a}(t)$ and used (\ref{3.6.2.21}) for the bound of $k_{A\bN}$.  Noting $\bb \aaa^{-1}(\tir+3)^{-1}|\Lb \tir|\les 1$,
$$
\a\int_{\D_{u_1}^{t_1}}\aaa^{-\a-1}\{(\tir+3)^{-1}|\Lb\tir|(|L\psi|^2+|\sn\psi|^2)\}\les \a\int_{u_1}^{u_*} F_{0,-\a}[\psi](\H_u^{t_1})du.$$
 This term can be treated by Gronwall's inequality. 

Further, since
$$Q_{L\Lb}=|\sn \psi|^2,\quad Q_{AB}=\sn_A\psi\c \sn_B\psi-\f12\ga_{AB}(-L\psi\Lb\psi+|\sn\psi|^2),$$
by noting the symbolic formula $\tr\sl{k}=\fB+[L\Phi]$ (due to (\ref{k1})) and using (\ref{3.6.2.21}) and (\ref{3.11.4.21}), we derive that
\begin{align*}
&k^{AB}Q_{AB}=\f12\tr\sl{k} L \psi\Lb \psi+\hk_{AB}\sn_A \psi \sn_B\psi \\
&|\bb k_{\bN\bN} Q_{L \Lb}, \bb k_{AB}Q_{AB}| \les \tir^{-1}(|\Lb \psi||L\psi|+|\sn \psi|^2).
\end{align*}
With the help of Cauchy-Schwarz inequality, we bound
\begin{align*}
\int_{\D_{u_1}^{t_1}}\aaa^{-\a}|k_{\bN\bN} Q_{L\Lb}, k^{AB}Q_{AB}|&\les \int_{u_1}^{u_*} F_{0,-\a}[\psi](\H_u^{t_1})du+\int_{\D_{u_1}^{t_1}} \l t'\r^{-2}\aaa^{-\a}|\Lb \psi|^2\\
    &\les \int_{u_1}^{u_*} F_{0,-\a}[\psi](\H_u^{t_1})du+\int_0^{t_1} \l t'\r^{-2}E_{-\a}[\psi](t') dt'
\end{align*}
which can be easily absorbed by the left-hand side of (\ref{7.13.2.22}) by virtue of Gronwall's inequality. Substituting the above estimate into (\ref{7.13.2.22}) and using Gronwall's inequality imply (\ref{6.21.2.21}).
\end{proof}

\subsection{Weighted energy inequalities}
In this subsection, we give the estimate for weighted energies based on the method given in \cite[Section 7]{Wangrough}.
We recall the following standard result.
\begin{lemma}[Identities of integration by parts]
Let $\D_{u}^{t} = \{0\le t'\le t, u\le u'\le u_*\}$. There holds for smooth scalar functions $f$ that 
\begin{equation}\label{8.14.1.20}
\int_{\D_u^t} (Xf+f \bd^\mu X_\mu) d\mu_g dt'= \int_{\Sigma_t} f X(t)+\int_{\H_u^t} f X(u)-\int_{\Sigma_0} fX(t)-\int_{\H_{u_*}^t} f X(u).
\end{equation}
The following formulas can be obtained by using Proposition \ref{6.7con}
\begin{equation*}
\bd^\mu L_\mu=\tr\chi-k_{\bN\bN}, \quad \bd^\mu \Lb_\mu=\tr\chib-k_{\bN\bN}.
\end{equation*}
\end{lemma}

Next we derive the weighted energy identity for scalar functions $\psi$. Rewriting (\ref{6.30.1.19}) gives
\begin{equation} \label{wave2}
\begin{split}
\Box_{\bg} \psi& =\sD \psi+2\zeta\c \sn \psi-\Lb(L\psi+ h\psi)
-(\hb-k_{\bN\bN})(L\psi+h\psi)\\
& \quad \, +(\Lb h+h \hb- h k_{\bN\bN})\psi.
\end{split}
\end{equation}
With $\a\ge 0$, we multiply (\ref{wave2}) by $(L\psi+h \psi) \tir^m \aaa^{-\a}$ and integrate over $\D_{u_1}^{t_1}$ to obtain
\begin{align*}
&\int_{\D_{u_1}^{t_1}}\Box_\bg \psi\c (L\psi+ h \psi) \tir^m \aaa^{-\a}\bb d\mu_\ga du dt\\
& = \int_{\D_{u_1}^{t_1}} \left(\sD \psi+2\zeta\c \sn \psi\right)
(L\psi +  h \psi) \tir^m \aaa^{-\a}\bb  d\mu_\ga du dt \nn \displaybreak[0]\\
& - \int_{\D_{u_1}^{t_1}} \left(\frac{1}{2} \Lb( (L\psi+h \psi)^2 )
+(\hb-k_{\bN\bN}) (L\psi+ h\psi)^2 \right) \tir^m \aaa^{-\a}\bb d\mu_{\ga}du dt \nn \displaybreak[0]\\
& + \int_{\D_{u_1}^{t_1}} \left(\f12 \mu   - h k_{\bN\bN}\right)\psi
(L\psi+h \psi) \tir^m \aaa^{-\a}\bb d\mu_{\ga} du dt\nn,
\end{align*}
where we used $2(\Lb h + h\hb)=\mu$. 
 Using  (\ref{cmu_2}),  $\chi =\chih + h \ga$ and (\ref{lb}), by integrating in $\D_{u_1}^{t_1}$  and rearranging the terms, we obtain that
\begin{align*}
& \f12 \int   \big\{L(\aaa^{-\a}\bb \tir^m|\sn \psi|^2)+(4h+k_{\bN\bN}+\frac{\a}{\aaa(\tir+3)}-\frac{m}{\tir}) \aaa^{-\a}\bb\tir^m |\sn \psi|^2\nn\\
& +\aaa^{-\a}\bb(\Lb(\tir^m |L\psi+h\psi|^2 )+\tir^m(2\hb -m\tir^{-1}\Lb\tir)|L\psi+h\psi|^2) \big\} d\mu_\ga du dt \nn \displaybreak[0]\\
&=\int \tir^m\aaa^{-\a}\bb\left\{(2\zeta-\sn\log\bb)\c \sn \psi (L\psi+ h\psi) +{k}_{\bN\bN} |L\psi+ h\psi|^2
-\chih\c\sn\psi\c\sn\psi \right\}d\mu_\ga du dt  \nn \displaybreak[0]\\
&+\int \tir^m \aaa^{-\a}\bb\psi\left\{\Big(\f12 \mu- h k_{\bN\bN} \Big) (L\psi+ h \psi)-\sn  \psi \sn h\right\} d\mu_{\ga}dudt\\
&-\int \Box_\bg \psi (L\psi+ h\psi)\tir^m \aaa^{-\a}\bb d\mu_\ga du dt\nn.
\end{align*}
Using (\ref{8.14.1.20}), we derive
\begin{lemma}\label{prop:MA}
Let (\ref{3.12.1.21})-(\ref{1.25.1.22}) hold. For $0<t_1<T_*$ and $u_0\le u_1\le u_*$,  with $m=1,2$, there holds for scalar functions $\psi$ 
\begin{align}\label{B_1}\tag{\bf {MA}}
&2\int_0^{t_1} \int_{S_{t,u_1}} \tir^m \aaa^{-\a}|L(v_t^{\f12}\psi)|^2 d\omega dt+\int_{u_1}^{u_*} \int_{S_{t_1, u}} \aaa^{-\a}\bb \tir^m \left(|L(v_t^{\f12}\psi)|^2+|\sn\psi|^2v_{t}\right) d\omega du\nn\\
&+\int_{\D_{u_1}^{t_1}}  \tir^m\Big(\aaa^{-\a}\bb\big((\a\frac{\tir}{\aaa(\tir+3)}-m)\Lb \tir\c\tir^{-1}- k_{\bN\bN}
\big)|L(v_t^{\f12}\psi)|^2\nn\\&\qquad\qquad\qquad+(\frac{2}{\tir}-\frac{m}{\tir}+k_{\bN\bN}+\frac{\a}{\aaa(\tir+3)})\aaa^{-\a}\bb|\sn\psi|^2v_t \Big) d\omega du dt\nn \displaybreak[0]\\
& =\int_{\D_{u_1}^{t_1}}-2 \Box_\bg \psi(L\psi+h \psi)\tir^m\aaa^{-\a}\bb d\mu_\ga du dt+\int_{u_1}^{u_*} \int_{S_{0, u}}\aaa^{-\a} \bb \tir^m \left(|L(v_t^{\f12}\psi)|^2+|\sn\psi|^2v_{t}\right) d\omega du\nn\\
&+ 2\int_0^{t_1} \int_{S_{t, u_*}} \aaa^{-\a}\tir^m |L(v_t^{\f12}\psi)|^2 d\omega d t +\mbox{error}(m) \nn 
\end{align}
where
\begin{align*}
\mbox{error}(m)&=-\int_{\D_{u_1}^{t_1}}\aaa^{-\a}\bb\tir^m\{\big(\tr\chi-\frac{2}{\tir}\big)|\sn \psi|^2 +2\chih\c\sn\psi\c\sn\psi \}d\mu_\ga du dt\\
&+2\Big(\int_{\D_{u_1}^{t_1}} \tir^m\aaa^{-\a}\bb\left\{(2\zeta-\sn\log \bb)\c \sn \psi(L\psi+h\psi)
\right\}d\mu_\ga dudt  \nn \displaybreak[0]\\
&+\int_{\D_{u_1}^{t_1}} \tir^m \aaa^{-\a}\bb \psi\left\{\Big(\f12 \mu- h k_{\bN\bN}\Big) (L\psi+h \psi)-\sn \psi \sn  h\right\} d\mu_\ga dudt\Big).
\end{align*}
\end{lemma}
The above result is obtained by direct calculation with the help of (\ref{8.14.1.20}), $L v_t=2 h v_t$ and 
\begin{equation*}
|L\psi+ h\psi|^2 v_t= |L(v_t^{\f12} \psi)|^2.
\end{equation*}

Now we recall the symbol $\bA$ includes the terms $\chih, k_{\bN A}$, and $ \tr\chi-\frac{2}{\tir}$,  and
let the symbol $\J[\psi]$ denotes any quadratic terms with factors being either $\sn \psi$ or $L\psi+h \psi$. $\mbox{error}(m)$ is written as
\begin{equation}\label{3.13.4.21}
\begin{split}
\mbox{error}(m)&=\int\tir^m\aaa^{-\a}\bb\big(\bA \J[\psi]+\ud \bA\c\sn\psi (L+h)\psi\big)d\mu_\ga du dt \\
&+2\int\tir^m \aaa^{-\a}\bb\psi\left\{\Big(\f12 \mu- h k_{\bN\bN}\Big)(L\psi+ h \psi)-\sn  \psi \sn  h\right\} d\mu_\ga dudt.
\end{split}
\end{equation}
We first consider the second line on the left-hand side of (\ref{B_1}). By using (\ref{4.12.1.24}) and (\ref{11.13.3.23}) we derive
\begin{align*}
\bb \tir k_{\bN\bN}&=\f12\wp\bb\tir \Lb \varrho+O(\Delta_0^\f12 \l t\r^{-1+2\delta})=\f12\wp \bb\tir ([\Lb\varrho]_{-})+O(\Delta_0^\f12 \l t\r^{-1+2\delta}).
\end{align*}
The contribution of the last error term is
 $$|\int_{\D_{u_1}^{t_1}} \tir^{m-1} \aaa^{-\a} O(\Delta_0^\f12\l t\r^{-1+2\delta}) |L(v_t^{\f12}\psi)|^2 d\omega du dt|\les \Delta_0^\f12 \int_{u_1}^{u_*}F_{m, -\a}[\psi](\H_u^{t_1})du.$$ 
  The other term can be dropped since $[\Lb\varrho]_{-}<0$. 
Moreover using (\ref{10.9.5.22}), we bound
\begin{align*}
|\int_{\D_{u_1}^{t_1}}  \tir^m \aaa^{-\a-1}\frac{(\a\frac{\tir}{\tir+3}-m\aaa)\bb\Lb \tir}{\tir}|L(v_t^{\f12}\psi)|^2 d\omega du dt|\les \int_{u_1}^{u_*} F_{m, -\a}[\psi](\H_u^{t_1}) du.
\end{align*}
Hence we derive that
\begin{align*}
&2\int_0^{t_1} \int_{S_{t,u_1}} \tir^m |L(v_t^{\f12}\psi)|^2\aaa^{-\a} d\omega dt+\int_{\D_u^t}2 \Box_\bg \psi(L\psi+ h \psi)\tir^m \aaa^{-\a}\bb d\mu_\ga du dt\nn\\
&\qquad +\int_{\D^{t_1}_{u_1}} \tir^{m-1}
  \Big(2-m+\frac{\tir}{\tir+3} \a\aaa^{-1}+\tir k_{\bN\bN}\Big)|\sn\psi|^2v_t \aaa^{-\a}\bb d\omega du dt \nn\\
  \displaybreak[0]
&\qquad +\int_{u_1}^{u_*} \int_{S_{t_1, u}} \aaa^{-\a}\bb \tir^m \left(|L(v_{t_1}^{\f12}\psi)|^2+|\sn\psi|^2v_{t_1}\right) d\omega du \\
& \les \int_0^{t_1} \int_{S_{t, u_*}} \tir^m \aaa^{-\a}|L(v_t^{\f12}\psi)|^2 d\omega d t+\int_{u_1}^{u_*} F_{m, -\a}[\psi](\H_u^{t_1}) du +|\mbox{error}(m)|.\nn
\end{align*}
We now consider the last line in $\mbox{error}(m)$, which is recast below,
\begin{equation*}
\Er(m)=\int_{\D^{t_1}_{u_1}} \tir^m \aaa^{-\a}\bb\psi\left\{\Big(\f12 \mu- h k_{\bN\bN}-\f12\varpi+\f12\varpi\Big) (L\psi+h \psi)-\sn \psi \sn h\right\} d\mu_\ga dudt.
\end{equation*}
Due to $\tir \varpi=O(\fB)$ and Cauchy-Schwarz inequality, 
\begin{align*}
\int_{\D^{t_1}_{u_1}} &\tir^m \aaa^{-\a}\bb|\psi||\varpi||L\psi+h \psi|d\mu_\ga dudt\\
&\les \int_{u_1}^{u_*}F_{m, -\a}[\psi](\H_u^{t_1}) du+\int_{\D^{t_1}_{u_1}} \l t\r^{-2}|\tir^{\f12 m}\aaa^{-\frac{\a}{2}}\psi|^2 du dt d\omega.
\end{align*}
Using (\ref{1.7.1.23}), $\|\tir^{2-3\delta}(\f12\mu-h k_{\bN\bN}-\f12\varpi)\|_{L_u^2 L_t^\infty L_\omega^4}\les \Delta_0$, we derive
\begin{align*}
|\int_{\D_{u_1}^{t_1}}&\tir^m \aaa^{-\a}\bb\psi(\f12 \mu- h k_{\bN\bN}-\f12\varpi\Big) (L\psi+h \psi) d\mu_\ga du dt|\\
&\qquad\les \left(\int_0^{t_1}\l t'\r^{-1-2\delta}\|\tir^{\f12m}\aaa^{-\f12\a}\psi\|^2_{L_u^2 L_\omega^4} dt'\right)^\f12\|\tir^{1+\f12m}\aaa^{-\f12 \a}(L \psi+h\psi)\|_{L_u^\infty L_t^2 L_\omega^2}\\
&\qquad\cdot\|\bb\tir^{\f12+\delta+1}(\mu-\tr\chi k_{\bN\bN}-\varpi)\|_{L_u^2 L_t^\infty L_\omega^4}\\
&\qquad\les \Delta_0 \sup_{u\in[u_1, u_*]} F_{m, -\a}[\psi]^\f12(\H_u^{t_1})(\int_0^{t_1} \l t'\r^{-1-2\delta} WL_{m, -\a}[\psi](t') dt')^\f12
\end{align*}
where we applied $\|\tir^{\f12m}\aaa^{-\f12\a}\psi\|_{L_u^2 L_\omega^4}\les  WL_{m, -\a}[\psi]^\f12(t)$ (due to Sobolev embedding). 

Next, using (\ref{3.11.3.21}) for $\|\sn\bA_b\|_{L_\omega^4}$, we bound the other term in $\Er(m)$,
\begin{align*}
\int_{\D^{t_1}_{u_1}}\tir^m \aaa^{-\a}\bb|\psi \sn \psi \sn h| d\mu_\ga du dt\les \int_0^{t_1} \Delta_0 \l t\r^{-\frac{7}{4}+\delta}\|\tir^{\f12m}\aaa^{-\f12 \a}\sn \psi\|_{L^2(\Sigma_t)}\|\bb^\f12 \aaa^{-\f12\a}\tir^{\f12m}\psi\|_{L_u^2 L_\omega^4}.
\end{align*}
Therefore, by Sobolev embedding on spheres, we conclude that
\begin{align*}
|\Er(m)|&\les\int_0^{t_1} \l t'\r^{-\frac{7}{4}+2\delta} WL_{m, -\a}[\psi](t')+\int_{u_1}^{u_*}F_{m, -\a}[\psi](\H_u^{t_1}) du\\
&+ \Delta_0 \sup_{u\in[u_1, u_*]} F_{m, -\a}[\psi]^\f12(\H_u^{t_1})(\int_0^{t_1} \l t'\r^{-1-\delta} WL_{m, -\a}[\psi](t') dt')^\f12.
\end{align*}
In view of (\ref{9.4.1.22}), we derive
$$|\int_{\D_{u_1}^{t_1}}\tir^m \bA \J[\psi]\aaa^{-\a}\bb d\mu_\ga du dt|\les \Delta_0^\f12 \int_0^{t_1}\l t\r^{-\frac{7}{4}+\delta}W_{m, -\a}[\psi](t) dt.$$
Using (\ref{3.6.2.21}), we have
\begin{align*}
|\int_{\D_{u_1}^{t_1}}&\tir^m \ud\bA \sn \psi (L+h)\psi \aaa^{-\a}\bb d\mu_\ga du dt|\\
&\les \| \l t'\r^{\f12+\frac{\delta}{2}}\bb^\f12\ud \bA\|_{L_u^2 L_t^\infty L_\omega^\infty}\sup_u F_{m, -\a}[\psi]^\f12(\H_u^{t_1})(\int_0^{t_1} \l t'\r^{-1-\delta} W_{m, -\a}[\psi](t') dt')^\f12\\
&\les \Delta_0 \sup_{u\in[u_1, u_*]} F_{m, -\a}[\psi]^\f12(\H_u^{t_1})(\int_0^{t_1} \l t'\r^{-1-\delta} W_{m, -\a}[\psi](t') dt')^\f12.
\end{align*}
 Thus we have treated all the terms in $\mbox{error}(m)$.
For the remaining term in (\ref{B_1}), using (\ref{4.12.1.24}) we have
\begin{align}\label{11.14.1.23}
&\int_{\D_{u_1}^{t_1}}\aaa^{-\a}\bb  k_{\bN\bN}|\sn \psi|^2 \tir^m du d\mu_\ga dt\nn\\
&= \int_0^{t_1} \int_{\Sigma_t\cap[u_1,u_*]}\Big(\f12\wp[\Lb \varrho]_-+O(\Delta_0^\f12\l t\r^{-2+\delta})\Big)|\sn\psi|^2\tir^m\aaa^{-\a}\bb  du d\mu_\ga dt.
\end{align}

Finally, to bound $\|\tir^{\f12 m-1}\psi\|_{L^2_\Sigma}$, we note due to (\ref{8.14.1.20})
\begin{align}\label{4.12.3.24}
&\int_{\Sigma\cap[u_1, u_*]}|\psi|^2\aaa^{-\a}\tir^{m-2}(t_1)+\int_{\D_{u_1}^{t_1}}|\psi|^2\aaa^{-\a} (2-m+\a\frac{\tir}{\tir+3}\aaa^{-1})\tir^{m-3}\\
  &=\int_{\Sigma_0\cap[u_1, u_*]}\tir^{m-2}\aaa^{-\a}|\psi|^2+\int_{\D_{u_1}^{t_1}} \aaa^{-\a}\{2(L\psi+h\psi)\psi \tir^{m-2}-k_{\bN\bN}|\psi|^2 \tir^{m-2}\},\nn
\end{align}
where the last term can be similarly treated as (\ref{11.14.1.23}). For the second term, note
\begin{align*} 
|\int_{\D_{u_1}^{t_1}} \aaa^{-\a}(L\psi+h\psi)\psi \tir^{m-2}|&\les\int_{u_1}^{u_*} F_{m, -\a}[\psi]^\f12(\H_{u}^{t_1})\|\aaa^{-\frac{\a}{2}}\bb\tir^{\f12 m-2}\psi\|_{L^2(\H_u^{t_1})} du\\
&\les \int_{u_1}^{u_*}F_{m, -\a}[\psi](\H_{u}^{t_1})du+\int_{\D_{u_1}^{t_1}}\l t' \r^{-2}\log \l t'\r \psi^2 \aaa^{-\a}\bb\tir^{m-2} du d\mu_\ga dt'.   
\end{align*}
The above two terms can be treated by Gronwall's inequality. 

 Adding a positive constant multiple of (\ref{4.12.3.24}) to the energy inequality (\ref{B_1}), summarizing the above control on error terms, we conclude
\begin{proposition}\label{MA2}
Assume (\ref{3.12.1.21})-(\ref{1.25.1.22}) hold. For $0<t_1<T_*$, $u_0\le u_1\le u_*$, with $m=1,2$ and $\a\ge 0$,  there holds the weighted energy inequality for smooth scalar functions $\psi$,  
\begin{align*}
&\int_0^{t_1} \int_{S_{t,u_1}} \aaa^{-\a}\tir^m |L(v_t^{\f12}\psi)|^2 d\omega dt+\int_{u_1}^{u_*} \int_{S_{t_1, u}} \aaa^{-\a}\bb \tir^m \left(|L(v_t^{\f12}\psi)|^2+(|\sn\psi|^2+\tir^{-2}|\psi|^2)v_{t}\right) d\omega du \\
\displaybreak[0]
&\qquad +\int_{\D^{t_1}_{u_1}}  \tir^{m-1}
  (2-m+\a\frac{\tir}{\aaa(\tir+3)})(|\sn\psi|^2+|\psi|^2\tir^{-2})\aaa^{-\a}d\mu_g dt \\
&\qquad \les \int_{\D_{u_1}^{t_1}}|\Box_\bg \psi(L\psi+ h \psi)|\tir^m \aaa^{-\a}\bb d\mu_\ga du dt+\int^{u_*}_{u_1} F_{m, -\a}[\psi](\H_u^{t_1}) du\\
&\qquad+\int_0^{t_1} (-\wp[\Lb\varrho]_{-}+\l t\r^{-1-\delta})WL_{m, -\a}[\psi](t, [u_1,u_*])dt \\
&\qquad+\int_{u_1}^{u_*} \int_{S_{0, u}} \aaa^{-\a}\bb\tir^m \left(|L(v_t^{\f12}\psi)|^2+(|\sn\psi|^2+|\psi|^2 \tir^{-2})v_{t}\right) d\omega du \\
&\qquad+ \int_0^{t_1} \int_{S_{t, u_*}} \tir^m \aaa^{-\a}|L(v_t^{\f12}\psi)|^2 d\omega d t.
\end{align*}
\end{proposition}
We give an inequality of Gronwall's type below, which is adapted from \cite[Lemma 2.4]{KWY}.
\begin{lemma}\label{1.6.4.18}
Let $f_0(\iota)$ be an nonnegative non-increasing function. Let $f$ and $g$ be two nonnegative continuous functions such that for any $t>0, 2\le \iota\le u_*$ there holds 
\begin{equation}\label{1.6.1.18}
\begin{split}
f(t,\iota)+g(t,\iota)&\les f_0(\iota)+\int_0^t f(t', \iota) (\l t'\r^{-a}+\M_0\l t'\r^{-1}(\f12\wp\log (\frac{\l t'\r}{2})+1)^{-1})d t'\\
&+\int_{\iota}^{u_*}  g(t,\iota') d\iota',
\end{split}
\end{equation}
with the constant $a>1$.
If $f(t, \iota)$ is non-increasing with respect to $\iota$ and $g(t, \iota)$ is non-decreasing with respect to $t$, then
\begin{equation*}
f(t, \iota)+g(t, \iota)\les f_0(\iota)(\f12\wp\log (\frac{\l t\r}{2})+1)^\M, \quad\forall t>0
\end{equation*}
where we take $\M=C \M_0\wp^{-1}$ and $\M \ge 15$, with $C>2$ a constant.
\end{lemma}
\begin{proof}
We first fix $u_1\le \iota\le u_*$. Applying Gronwall's inequality to $f(t, \cdot)+g(t,\cdot)$ gives
\begin{equation*}
f(t,\iota)+g(t,\iota)\les f_0(u_1)+\int_0^t f(t', u_1) (\l t'\r^{-a}+\M_0\l t'\r^{-1}(\f12\wp\log (\frac{\l t'\r}{2})+1)^{-1})d t', 
\end{equation*}
which certainly holds true for $\iota=u_1$, i.e. 
\begin{align*}
f(t,u_1)+g(t,u_1)&\les f_0(u_1)+\int_0^t f(t', u_1) (\l t'\r^{-a}+\M_0\l t'\r^{-1}\big(\f12\wp\log (\frac{\l t'\r}{2})+1\big)^{-1})d t'
\end{align*}
Applying Gronwall's inequality to $f(\cdot, u_1)+g(\cdot, u_1)$,  in view of (\ref{1.6.1.18}), we have 
\begin{equation*}
f(t, u_1)+g(t, u_1)\les f_0(u_1)(\f12\wp\log(\frac{\l t\r}{2})+1)^\M,
\end{equation*}
as stated. 
\end{proof}

\section{Elliptic estimates}\label{Elliptic}
In this section, under the assumptions (\ref{3.12.1.21})-(\ref{1.25.1.22}), we first give a set of basic elliptic estimates on $S_{t,u}$ in Lemma \ref{hodge_1}. Then in Lemma \ref{error_prod} (and Lemma \ref{8.28.9.23}) we provide basic product estimates that will be frequently used for nonlinear estimates. Using these two types of estimates, we derive in Lemma \ref{3.31.3.22} higher order elliptic estimates for Hodge-operators on spheres $S_{t,u}$.  The results in this section will be crucially used to control connection coefficients.
\begin{lemma}\label{hodge_1}
  Assume (\ref{3.12.1.21})-(\ref{1.25.1.22}) hold. Let $\D_1$ be the operator that takes any $1$-form $F$ on $S_{t,u}$ into the pair of functions $(\sl\div F, \sl\curl F)$. Let $\D_2$ be the operator that takes a $2$-covariant, symmetric, traceless tensor $F$ to the $1$-form $\sl{\div} F$. 
  The kernel of both $\D_1$ and $\D_2$ in $L^2_\omega$ is trivial.  With $\D$ either $\D_1$ or $\D_2$, 
 there hold for $2\le p<\infty$ and the tensor-fields $F$ in the domain of $\D$ and scalar functions $f$
\begin{align}
&\|\sn F\|_{L_\omega^p(S_{t,u})}+\|\tir^{-1} F\|_{L_\omega^p(S_{t,u})}\les \|\D F\|_{L_\omega^p(S_{t,u})},\label{4.4.1.21}\\
\displaybreak[0]
&\|\snc^2 f, \tir^{-1}\snc f\|_{L^p_\omega}+\|\sn^2 f, \tir^{-1}\sn f\|_{L_\omega^p}\les \|\sD f\|_{L^p_\omega}\label{11.19.1.23}\\
&\|\snc^3 f\|_{L_\omega^2}+\tir^{-1}\|\snc^2 f\|_{L^2_\omega}\les \|\snc\sDc f\|_{L_\omega^2}\label{7.16.3.24} 
\end{align}
Denote by $\snc$, $\stc{\D}_1$, $\stc{\D}_2$ the operators similarly defined as $\sn$, $\D_1$ and $\D_2$, with respect to the metric $\gac$ instead of $\ga$. The above results hold the same.   
\end{lemma}
The proof of the above result mainly relies on (\ref{10.1.1.22}) and (\ref{L4BA1}). We first adapt the uniformization theorem \cite[Lemma 2.3.2]{CK} into the following result.
\begin{lemma}[Uniformization Lemma]\label{11.3.1.23}
Let $S=S_{t,u}$ be fixed. There exists a conformal transformation of the metric such that with a conformal factor $\tau$, $\gac=e^{2\tau}(\ckc \tir)^2\zga$, with $\zga$ a standard round metric such that $K(\zga)=1$ on $S$.  Moreover the conformal factor $\tau$ can be chosen such that, 
\begin{align*}
 |\tau|+|\zsn \tau| +\|\zsn{}^2\tau\|_{L^p_\omega}\les 1, \quad 1<p<\infty 
\end{align*}
where $\zsn$ denotes the  Levi-civita connection of $\zga$. 
\end{lemma}
Indeed, due to the estimate (\ref{10.1.1.22}) and $c-\ckc=O(\Delta_0^\f12)$ from (\ref{6.24.1.21}), we can obtain 
\begin{equation}\label{11.1.3.23}
|\stc{K}-\frac{1}{(\ckc\tir)^2}|\les \l t\r^{-2}\Delta_0^\f12.
\end{equation}
Note by using  Gauss-Bonnet theorem 
\begin{align*}
\int_S(\stc{K}-\frac{1}{(\ckc\tir)^2}) d\mu_\gac+\ckc^{-2}\tir^{-2}\int_S 1 d\mu_{\gac}=4\pi. 
\end{align*}
Due to $|S|_\gac\approx |S|_\ga=O(\l t\r^2)$, thus we infer
\begin{align*}
\ckc^{-2}\tir^{-2}\int_S 1 d\mu_{\gac}=4\pi+O(\Delta_0^\f12).
\end{align*}
Consequently,
$$|S|_{\gac}\c \stc{K}=4\pi+O(\Delta_0^\f12).$$ Thus we conclude Lemma \ref{11.3.1.23} by using \cite[Lemma 2.3.2]{CK}.

\begin{proof}[Proof of Lemma \ref{hodge_1}]
 We first control the difference between  $\gac$ and $\zga$.  On a fixed sphere $S_{t,u}$, we set $\tzga:=(\ckc \tir)^2\zga$, where $\ckc \tir$ is a constant on $S_{t,u}$. It is straightforward to check, by applying \cite[Page 47, (2.4.1)]{CK}, that 
\begin{align}\label{11.2.8.23}
\sta{(0)}\sD\tau-1=-(\ckc \tir)^2\stc{K}e^{2\tau}.
\end{align}
Differentiating the above identity and using B\"{o}chner indentity gives
\begin{equation}\label{dsdtau}
\sta{(0)}\sD\zsn\tau=\zsn\tau(1-2(\ckc \tir)^2\stc{K} e^{2\tau})-(\ckc \tir)^2\zsn\stc{K}e^{2\tau}.
\end{equation}
As the first step, we will improve the estimates of $\tau$ in Lemma \ref{11.3.1.23} to
\begin{equation}\label{11.1.4.23}
\|\zsn\tau\|_{W^{2,p}_\omega}\les \Delta_0,\, |\tau|\les \Delta_0^\f12, \, 2\le p\le 4,
\end{equation}
where $W^{2,p}_\omega$ stands for the standard $W^{2,p}$ norm on $({\mathbb S}^2, \zga)$. 

We make an auxiliary bootstrap assumption that
\begin{equation}\label{11.2.7.23}
|\tau|<1.
\end{equation}
Note at $t=0$, on the sphere $S_{0,u}$, that $\tau=0$.  We will improve (\ref{11.2.7.23}) to $\tau=O(\Delta_0^\f12)$ for all $t<T_*$, as stated in the second estimate in (\ref{11.1.4.23}).
   
Let $2<p\le 4$ in the sequel unless specified otherly. Note, for scalar functions $f$,  $|\zsn f|\les |\tir\sn f|$  due to (\ref{11.2.7.23}). Applying the elliptic estimates on $({\mathbb S}^2, \zga)$, we derive by using (\ref{dsdtau}) and (\ref{10.1.1.22}) that
\begin{align}\label{11.2.11.23}
\begin{split}
\|\zsn\tau\|_{W^{2,p}_\omega}&\les\|\zsn \tau(1-2e^{2\tau})\|_{L^p_\omega}+\|(\ckc \tir)^2\stc{K}-1\|_{L^\infty_\omega}\|e^{2\tau}\zsn \tau\|_{L^p_\omega}+(\ckc \tir)^2\|\zsn\stc{K}\|_{L_\omega^p}\\
&\les \|\zsn\tau (1-2e^{2\tau}) \|_{L^p_\omega}+\Delta_0^\f12\|e^{2\tau}\zsn\tau\|_{L_\omega^p}+\l t\r^{-\frac{3}{4}+\delta}\Delta_0\\
&\les\|\tau \c \zsn\tau\|_{L^p_\omega}+\Delta_0^\f12 \|\zsn\tau\|_{L_\omega^p}+\Delta_0\l t\r^{-\frac{3}{4}+\delta},
\end{split}
\end{align}
 where we used (\ref{11.2.7.23}) to derive the last inequality. Next we use Sobolev embedding to bound 
 \begin{equation}\label{11.2.10.23}
 \|\tau\|_{L^\infty_\omega}\les \|\zsn\tau\|_{L^p_\omega} +\|\tau\|_{L^2_\omega}.
 \end{equation} 
It remains to bound $\|\tau\|_{L^2_\omega}$. For this purpose, we first recast (\ref{11.2.8.23}) to be
\begin{equation*}
(\ckc \tir)^2 \sDc\tau=e^{-2\tau}\sta{(0)}\sD\tau=e^{-2\tau}-1+(1-(\ckc\tir)^2\stc{K}).
\end{equation*}
It follows by integrating by parts the above identity on spheres $S_{t,u}$ that
\begin{equation*}
\int_{S_{t,u}} -(\ckc \tir)^2 |\snc\tau|^2=\int_{S_{t,u}} \left(e^{-2\tau}-1+(1-(\ckc\tir)^2\stc{K})\right)\tau.
\end{equation*}
Using (\ref{11.2.7.23}), (\ref{11.1.3.23}) and Cauchy-Schwarz inequality, we derive
\begin{align}\label{11.2.9.23}
\|\tau\|^2_{L^2_\omega}\les \|\ckc \tir\snc \tau\|^2_{L_\omega^2}+\Delta_0.
\end{align}
Due to (\ref{11.2.7.23}),  $|\cdot|_\gac\approx |\cdot|_{\tzga}$. Hence, combining (\ref{11.2.10.23}) and (\ref{11.2.9.23}), we obtain 
\begin{align*}
\|\tau \|_{L^\infty_\omega}\les \|\zsn \tau\|_{L^p_\omega}+\Delta_0^\f12. 
\end{align*}
 Substituting the above estimate to (\ref{11.2.11.23}), we obtain $\|\zsn\tau\|_{L_\omega^p}\les \Delta_0$. Then the first estimate in (\ref{11.1.4.23}) follows by using (\ref{11.2.11.23}). The second one follows as its consequence by using the above estimate. (\ref{11.2.7.23}) is therefore improved. 
  
 Due to $$\tensor{\Ga(\gac)}{^a_{bc}}-\tensor{\Ga(\tzga)}{^a_{bc}}=\p_b \tau \delta^a_c+\p_c\tau \delta^a_b-\p_d \tau{\tzga}^{da}\tzga_{bc},$$  it follows from (\ref{11.1.4.23}) that 
\begin{footnote}{Remark that  $\Ga(\gac)-\Ga(\tzga)$, regarded as $(1,2)$ tensor, could be evaluated with respect to metrics of different scalings. To be consistent, we evaluate the tensor and its $\zsn$-derivatives  by $\zga$.}\end{footnote}
\begin{align}\label{11.2.4.23}
\||\zsn{}^{\le 2}(\Ga(\gac)-\Ga(\tzga))|_{\zga}\|_{L^p_\omega}\les\Delta_0,\,\, 2\le p\le 4. 
\end{align}
 
Using  $\Ga(\ga)-\Ga(\gac)\approx\p_\omega\log c$ and the estimate $\|(\tir\sn)^{\le 2}\sn\log c\|_{L^4_\omega}\les \l t\r^{-\frac{7}{4}+\delta}\Delta_0$ in (\ref{L4BA1}),
we have
\begin{align}\label{11.3.2.23}
\||\zsn{}^{\le 2}(\Ga(\ga)-\Ga(\tzga))|_{\zga}\|_{L^p_\omega}\les \Delta_0,\,\, 2\le p\le 4. 
\end{align}
Using the above estimate, repeating the proof for \cite[Corollary 2.2.2.1, Lemma 2.3.1]{CK}, we obtain (\ref{11.19.1.23}), the following first order estimate for the Hodge operators relative to $\ga$,
 \begin{align*}
 &\|\sn F\|_{L^p_\omega}+\tir^{-1}\|F\|_{L^p_\omega}\les \|\D F\|_{L^p_\omega}, \, 2\le p<\infty,
 \end{align*}
and the same $L^p$ estimate as above holds with $\sn$ and $\D$ replaced by $\snc$ and $\stc{\D}$ which are same operators defined with respect to the metric $\gac$. (\ref{7.16.3.24}) can be obtained by applying B\"ochner indentity, integrating by parts on spheres and using the estimates of $K$ in Lemma \ref{5.13.11.21} (3). 
\end{proof}
\begin{lemma}\label{error_prod}
Under the assumptions (\ref{3.12.1.21})-(\ref{1.25.1.22}), there hold for $G, F$ either scalar functions or $S_{t,u}$-tangent tensor-fields,
\begin{align}
&\|G F\|_{L_\omega^p}\les \|G\|_{L_\omega^4}\|\sn_\Omega^{\le 1}F\|_{L_\omega^p}, \, 2\le p\le 4\label{error_prod_1}\\
&\sn_X(\fB\c F)= O(\l t\r^{-1})\big((\sn_X+\vs(X))F\big)+(1-\vs(X))O(\l t\r^{-1+\delta}\Delta_0)\c F\label{8.24.5.23}\\
&\sn_X(\chi\c F)=O(\l t\r^{-1})\big((\sn_X+\vs(X))F\big)+(1-\vs(X))O(\l t\r^{-\frac{7}{4}+\delta}\Delta_0)_{L^4_\omega}\c F\label{8.24.4.23}\\
&\sn_X\sn_{\tir \Lb}(\bA_b F)=X(\bA_b \sn_{\tir\Lb}F)+\sn_X F\big(\fB+\mho+\bA_b+O(\l t\r^{-1+\delta}\log \l t\r\Delta_0)_{L_\omega^4}\big)\nn\\
&\qquad\qquad\quad\quad+\big((\sn_X+\vs(X))(\fB+\mho+\tir \varpi)+X^{\le 1}\bA_b+O(\l t\r^{-1+\delta}\log \l t\r\Delta_0)_{L^2_u L_\omega^2}\big) \c F.\label{8.31.4.23}
\end{align}
\end{lemma}
\begin{proof}
(\ref{error_prod_1}) is a direct consequence of Sobolev embedding. (\ref{8.24.5.23}) can be obtained by using  (\ref{8.23.1.23}). 
(\ref{8.24.4.23}) is obtained by using (\ref{8.23.2.23}). 
 
To see (\ref{8.31.4.23}), we directly compute
\begin{align*}
\sn_X\sn_{\tir \Lb}(\bA_b F)&=\sn_X (\sn_{\tir \Lb}\bA_b F+\bA_b \sn_{\tir\Lb}F)\\
&=\sn_X F \sn_{\tir \Lb}\bA_b+X(\tir \Lb\bA_b) F+\sn_X(\bA_b \sn_{\tir\Lb}F).
\end{align*}
Using Proposition \ref{7.22.2.22} and (\ref{8.23.2.23}) to treat the first two terms
\begin{align*}
\sn_X F\sn_{\tir \Lb}\bA_b&=\sn_X F(\fB+\mho+\bA_b+O(\l t\r^{-1+\delta}\log \l t\r\Delta_0)_{L_\omega^4})\\
X(\tir \Lb\bA_b) F&=((\sn_X+\vs(X)) (\fB+\mho+\tir \varpi)+X^{\le 1}\bA_b+O(\l t\r^{-1+\delta}\log \l t\r\Delta_0)_{L^2_u L_\omega^2}) \c F.
\end{align*}
(\ref{8.31.4.23}) follows by combining the above calculations.
\end{proof}

Next we use Lemma \ref{hodge_1} to obtain the following structure equations and elliptic estimates under the assumption of (\ref{3.12.1.21})-(\ref{1.25.1.22}).
\begin{lemma}\label{3.31.3.22}
Let $X\in \{\Omega, S\}$. 

(1) There holds the symbolic formula 
 \begin{equation}\label{9.21.2.22}
[\D, \sn_X] F=\sn X\c \sn F+(\fR+\tir^{-2})\c F\star X.
\end{equation}
where $\fR=K-\frac{1}{\tir^2}, \bR_{LABC}, \chi\bA_{g,1},$ and $F\star X=\vs(X) \tir F+(1-\vs(X))X\c F$.

(2) There holds with $2\le p\le 4$ the following elliptic estimate
\begin{align*}
\|\tir\sn \sn_X^\ell F\|_{L^p_\omega}&+\|\sn_X^\ell F\|_{L^p_\omega}\les \|\tir\sn_X^{\le \ell}\D F\|_{L^p_\omega}, \ell=0,1. 
\end{align*}

(3)
Let $\D=\D_1$ or $\D_2$ and $X^2=X_2 X_1$, with $X_1, X_2\in \{S,\Omega\}$. We have for $\a=0, \f12$
\begin{align}
\|\bb^\a\tir&\sn\sn_{X}^2F\|_{L^2_u L_\omega^2}+\|\bb^\a\sn_X^2 F\|_{L^2_u L_\omega^2}\label{9.14.10.22}\\
&\les \|\bb^\a\tir\{\sn_X^2 \D F, \l t\r^{-1}\sn_\Omega^{\le 1}\sn_{X_1}F, \sn_{X_2}^{\le 1}(\sn F, \tir^{-1}F)\}\|_{L^2_u L_\omega^2}\nn\\
&+\l t\r^{-\frac{7}{4}+\delta}\Delta_0^\f12 \|\bb^\a\tir (\tir\sn F, F)\|_{L^2_u L_\omega^4}+\a\log \l t\r\Delta_0\|\bb^{\a-1}\sn_X^2 F\|_{L_u^2 L_\omega^4}.\nn
\end{align}
\end{lemma}

\begin{proof}
We claim for an $S$-tangent 1-tensor $F$, there holds
\begin{align*}
&[\sl{\div}, \sn_X]F=\sn_A X^B \c \sn_B F_A+(\fR +\tir^{-2}) F\star X\\
&[\sl{\curl},\sn_X]F=\ep^{AC} \sn_A X^B \sn_B F_C+(\fR +\tir^{-2}) F\star X\nn;
\end{align*}
 and for a symmetric traceless $S$ tangent 2-tensor $F$, there holds
\begin{equation*}
\sn_A \sn_X F_{AB}=\sn_X \sl{\div}F_B +\sn_A X^D \sn_D F_{AB}+ (\fR+\tir^{-2}) F\star X.
\end{equation*}
Indeed, if $X=S$,
using (\ref{6.22.16.19}) and  (\ref{2.18.1.22}), we can write
\begin{align*}
\D \sn_X F&=\sn_X \D F+\sn X\c \sn F+\big(\bR_{ACXB}+\tir\chi\c \bA_{g,1}\big) F.
\end{align*}
If $X=\Omega$, we use (\ref{8.10.1.21}). Hence in view of the definition of $\fR$, we obtain symbolically the above result for $S$-tangent 1-tensor. The identity for $S$-tangent symmetric traceless 2-tensor follows similarly. We then write them into (1) schematically.

Note, due to (1), (3) and (4) in Lemma \ref{5.13.11.21} and (\ref{1.25.2.22})
\begin{equation}\label{11.16.6.23}
\|\tir^2\fR\|_{L_\omega^\infty}\les \l t\r^{-\frac{3}{4}+\delta}\Delta_0^\f12,\quad \sn X=O(1).
\end{equation}
Therefore, 
\begin{equation}\label{9.1.5.23}
\|\tir \fR F\|_{L_\omega^p}\les \l t\r^{-\frac{7}{4}+\delta}\Delta_0^\f12\|F\|_{L_\omega^p}.
\end{equation}

Combining (\ref{11.16.6.23}) with (\ref{9.1.5.23}), we obtain the following inequality for an $S$-tangent tensor $F$
\begin{equation}\label{9.14.9.22}
\|\D \sn_X F\|_{L_\omega^p}\les \|\sn_X \D F\|_{L_\omega^p}+\|\sn F, \tir^{-1}F\|_{L_\omega^p}+ \l t\r^{-\frac{7}{4}+\delta}\Delta_0^\f12\|F\|_{L_\omega^p}.
\end{equation}
  We then apply (\ref{4.4.1.21}) in Lemma \ref{hodge_1} to $F$ and $\sn_X F$. 
  The elliptic estimate of (2) follows.
 
 Next we prove (3). Applying (\ref{9.21.2.22}) to $F_1=\sn_{X_1} F$ and $X=X_2$, also in view of (\ref{9.21.2.22}) with $X=X_1$ therein, we write
\begin{align}\label{9.17.1.22}
\begin{split}
\D\sn_{X_2} F_1&=\sn_{X_2} \D F_1+\sn {X_2}\c\sn F_1+(\fR+\tir^{-2})\c F_1\star X_2\\
\sn_{X_2} \D F_1&=\sn_X^2 \D F+\sn_{X_2}(\sn_A X_1\c \sn F_A)+\sn_{X_2}\{(\fR+\tir^{-2})\c F\star X_1\}.
\end{split}
\end{align}

We infer from Lemma \ref{5.13.11.21}, (\ref{5.6.03.21}), (\ref{1.25.1.22}), $L\la=-\Omega c=O(\l t\r^{-1+\delta}\Delta_0^\f12)$, (\ref{1.25.2.22}) and the basic calculation (\ref{10.17.1.23}) and (\ref{11.5.3.23}) that
\begin{align*}
\begin{split}
\sn_{X_2}\sn X_1&=\sn_{X_2}(\tir\chi)+ \sn_{X_2}\sn\Omega=O(1)+O(\l t\r^{-\frac{3}{4}+\delta}\Delta_0^\f12)_{L_\omega^4},\\
\sn_{X_2}\fR&=O(\l t\r^{-\frac{11}{4}+\delta}\Delta_0^\f12)_{L^4_\omega}, \quad \tir^{-1}(\sn^{\le 1}_{X_2} X_1, \sn_{X_1}^{\le 1}X_2)=O(1),\\
 \Delta_0^\f12\tir\fR&=O(\Delta_0\l t\r^{-\frac{7}{4}+\delta}).
\end{split}
\end{align*}
Therefore, with $F_1=\sn_{X_1}F$,  we derive in view of (\ref{9.17.1.22}) that
\begin{align*}
\begin{split}
\D \sn_{X_2} F_1&=\sn_X^2\D F+ O(\l t\r^{-1})\sn^{\le 1}_\Omega F_1+O(1)\sn_{X_2}^{\le 1}( \sn F, \tir^{-1}F)\\
&+O(\l t\r^{-\frac{7}{4}+\delta}\Delta_0^\f12)_{L_\omega^4} (\tir\sn F,  F).
\end{split}
\end{align*}
  Applying Lemma \ref{hodge_1} to $\bb^\a \sn_X^2 F$ with the help of (\ref{1.27.5.24}), we conclude (\ref{9.14.10.22}). 
  
\end{proof}
\section{Comparison estimates on connection coefficients}\label{low_ricci}
Connection coefficients $\bA_b, \chih, \zeta$ are not defined directly by derivatives of $\Phi$. Nevertheless, our assumptions are only on initial derivatives of $\Phi$. In this section, we provide comparison estimates in Proposition \ref{11.4.1.22} for connection coefficients $\bA_b, \chih, \zeta$. In particular  $\eh$ is a quantity defined via $\sn\Phi(\in \al_\f12)$, in Proposition \ref{8.12.1.23} and Proposition \ref{10.4.4.23}, we will derive various derivative estimates of $\eh$ which are crucial for justifying $\eh\in \bA_{g,1}$. In Proposition \ref{7.15.5.22}, we summarize a set of important decay estimates derived from Assumption \ref{5.13.11.21+} and geometric comparison. These estimates will be frequently used in later sections.
\subsection{Preliminary results}
We first provide a set of commutator estimates under the assumptions (\ref{3.12.1.21})-(\ref{1.25.1.22}). 
\begin{lemma}[Commutation lemma]\label{6.9.3.23}
\begin{enumerate}
	\item 
	With $X_1, X_2\in \{\Omega, S\}$, for a scalar function $f$,  there holds
	\begin{equation}\label{7.17.6.21}
	\tir[X_1, L]f=\left\{\begin{array}{lll}
O(1)X_1 f\\
O(1)\left((1-\vs(X_1))\Delta_0^\f12 \l t\r^{-\frac{3}{4}+\delta}+\vs(X_1)\right)X_1 f
\end{array}\right.
	\end{equation}
\begin{footnote}{If $X_1$ is a rotation field on the left-hand side, the $X_1$ on the right-hand side stands for all the rotation fields $\Omega$. We adopt this convention without explicitly mentioning  whenever there is no confusion occurs.}\end{footnote}
where the second line is the improved estimate under the assumption that
\begin{equation}\label{7.16.2.24} 
\pio_{LA}=O(\l t\r^{-\frac{3}{4}+\delta}\Delta_0^\f12).
\end{equation}
	If 
\begin{equation}\label{7.16.1.24}
\sn_{X_2}(\pio_{\Lb A}, \pio_{L A})=O(\l t\r^{-\frac{3}{4}+\delta}\Delta_0^\f12)_{L_\omega^4},
\end{equation} 
then
	\begin{equation}\label{7.17.7.21}
		X_2 [X_1, L]f, [X_2, L]X_1 f= O(\l t\r^{-1})X_2^{\le 1}X_1 f+(1-\vs(X_1))O(\l t\r^{-\frac{3}{4}+\delta}\Delta_0^\f12)_{L_\omega^4}\sn f
	\end{equation}
	 
	\item There hold for $S$-tangent tensor fields \begin{footnote}{We always regard scalar functions as a special subclass of $S$ tangent tensors.}\end{footnote} $F$ and $G$
	\begin{align}
		&[\sn_S, \sn_\Omega]F, [\sn_S, \tir \sn]F=\left\{\begin{array}{lll}
O(\Delta_0^\f12) \sn_\Omega F+ O(\l t\r^{-\frac{3}{4}+\delta}\Delta_0) F\\
O(\l t\r^{-\frac{3}{4}+\delta}\Delta_0^\f12)\sn^{\le 1}_\Omega F
\end{array}\right.\label{3.21.1.23}
\end{align}
where to obtain the second line for $[\sn_S, \sn_\Omega]F$, we assumed (\ref{7.16.2.24}).
 
 With $\a=-\f12, 0$,  there hold
\begin{align}\label{3.20.2.23}
\begin{split}
		&\|\bb^\a (\sn_S [\sn_S, \sn] G, [\sn_S, \sn]\sn_S G)\|_{L^2_\Sigma}\les \|\bb^\a(\sn \sn_S G, \sn_S^{\le 1} \sn G)\|_{L^2_\Sigma}\\
&\qquad\qquad\qquad+\l t\r^{-\frac{7}{4}+\delta}\log \l t\r^{\f12}\Delta_0(\|\bb^\a\sn_S^{\le 1} G\|_{L^2_\Sigma}+\|\bb^\a\tir(\tir\sn)^{\le 1} G\|_{L_u^2 L_\omega^4});\\
&\|\bb^\a(\sn_\Omega[\sn_\Omega, \sn_L]G, [\sn_\Omega, \sn_L]\sn_\Omega G)\|_{L^2_\Sigma}\les \l t\r^{-\frac{7}{4}+\delta}\log \l t\r^{\f12}\Delta_0^\f12\|\bb^\a(\sn^{1+\le 1}_\Omega G, \Delta_0^\f12 G)\|_{L^2_\Sigma}.
\end{split}
	\end{align}
where for obtaining the second estimate, we assumed (\ref{7.16.2.24}) and (\ref{7.16.1.24}).
	\item With $X_1=L, \sn$, it holds for $S$-tangent tensor $F$ that
	\begin{equation}\label{7.4.1.21}
\begin{split}
	[\sn_\bN, \sn_{X_1}] F&= O(\l t\r^{-2+\delta}\Delta_0) F+O(\l t\r^{-1})\{(\log \l t\r\Delta_0)_{L_\omega^4}^{1-\vs(X_1)}\bb^{-1}\sn_\bN F\\
&+(\l t\r^{\delta}\Delta_0\vs(X_1)+(1-\vs(X_1)))\sn F\}.
\end{split}
	\end{equation}
\end{enumerate}
\end{lemma}
\begin{proof}
Consider (1) first. 
It is straightforward to have, in view of (\ref{5.13.10.21}), that
	\begin{equation*}
		[X_1, L]f=\tir^{-1}X_1 f, \mbox{ if }X_1=S; \quad	[X_1, L]f=\pio_{LA} \sn_A  f, \mbox{ if }X_1=\Omega.
	\end{equation*}
	This implies the first line in (\ref{7.17.6.21}) by using (\ref{1.25.3.22}) and $|v|\les \l t\r^{-1}$, as well as the second estimate by using the assumption of $\pio_{AL}$.
	
	To show (\ref{7.17.7.21}) we derive  by using (\ref{1.25.3.22})
	\begin{align*}
		X_2[X_1, L]f&= O(\l t\r^{-1})X_2^{\le 1}(X_1 f), \mbox{ if } X_1=S\\
		X_2[X_1, L]f&=X_2(\pio_{LA} \sn_A f)\\
&=O(\Delta_0^\f12\l t\r^{-1})X_2^{\le 1} \Omega f+O(\l t\r^{-\frac{3}{4}+\delta}\Delta_0^\f12)_{L_\omega^4}\sn f, \mbox{ if } X_1=\Omega.
	\end{align*}
We then	summarize the above two cases into (\ref{7.17.7.21}).

	For (2), using (\ref{5.13.10.21}) and (\ref{cmu2}), we deduce
	\begin{align}\label{6.4.1.24}
		[\sn_S, \sn_\Omega]F&=\tir\pio_{LA}\sn F+\tir(\chi\c \zb+\bR_{AC4 B})F\c \Omega.
	\end{align}
Recall by using (\ref{1.29.4.22}), (\ref{3.16.1.22}) and (\ref{3.11.3.21})
	\begin{align}\label{8.22.7.23}
		\|\sn_X^l(\tir \bR_{AC4B}, \tir\chi\c \zb)\|_{L_\omega^4}\les \l t\r^{-2+\frac{l}{4}+\delta}\Delta_0, \,  l=0,1; \,  |\l t \r^{-\frac{1}{4}}\bR_{AC4B}, \chi\c \zb|\les \l t\r^{-3+\delta}\Delta_0.
	\end{align}
Using (\ref{3.6.2.21}), (\ref{cmu2}), (\ref{1.25.3.22}), the above pointwise estimate, we can obtain the first estimate and the estimate of $[\sn_S, \tir \sn]F$. The refined estimate follows similarly by using the assumption on $\pio_{LA}$.
 
 In view of (\ref{cmu2}), we write for $S$-tangent tensor field $G$ and $F$
	\begin{align*}
		\sn_S^l[\sn_S, \sn]G&=\sn_S^l\{\tir (\chi\sn G+(\chi\c \zb+\bR_{AC4B})\c G)\}.
	\end{align*}
With $l=1$, or when $l=0$ with $G$ replaced by $\sn_S G$, we obtain the first set of estimates in (\ref{3.20.2.23}) from the above formula by using (\ref{8.22.7.23}), (\ref{8.24.4.23}) and Sobolev embedding on spheres.

To see the second estimate, by using (\ref{6.4.1.24}), we derive
\begin{align*}
\sn_\Omega [\sn_L, \sn_\Omega]G=\sn_\Omega(\pio_{LA}\sn G+(\chi\c \zb+\bR_{AC4 B})G\c \Omega)\\
[\sn_L, \sn_\Omega]\sn_\Omega G=\pio_{LA}\sn\sn_\Omega G+(\chi\c \zb+\bR_{AC4 B})\sn_\Omega G\c \Omega. 
\end{align*}
The second set of estimates of (\ref{3.20.2.23}) can be obtained by applying the assumptions on $\sn^{\le 1}_{X_2}\pio_{L A}$, (\ref{8.22.7.23}), Sobolev embedding. 

	Recall from (\ref{7.04.5.21}), (\ref{cmu2}) and (\ref{7.03.1.19}), for $S_{t,u}$-tangent tensor field $F$,
	\begin{align}\label{7.04.2.21}
		\begin{split}
			[\sn_\bN, \sn_L]F_A&=(\zb^A-\zeta^A)\sn_A F-k_{\bN\bN} \sn_\bN F_A +\bR_{AB L\Lb} F_B+\ze\c\zb\c F\\
			[\sn_\bN, \sn_B]F_A&=-\theta_{BC}\sn_C F_A+(\zeta_B+\zb_B)\sn_\bN F_A+\bR_{AC\bN B}F_C+(\chi+\theta)\c (\zeta-\zb)F.
		\end{split}
	\end{align}
	Hence using (\ref{8.22.7.23}), (\ref{zeh}), (\ref{1.27.5.24}) and (\ref{3.6.2.21})  we summarize the above two lines into (\ref{7.4.1.21}).
\end{proof}

\subsection{Estimates of $\la$}
$\la$ and its derivatives are important quantities for controlling rotation fields, whose deformation tensor $\pio$ and its derivatives are crucial terms in energy estimates. In  Proposition \ref{10.16.1.22} we derive control on $\la$, which improves the bootstrap assumption (\ref{1.25.1.22}).

\begin{proposition}\label{10.16.1.22}
Let $X=\Omega, S$ and $\la$ stand for any of $\la\rp{a}$. Under the assumptions (\ref{3.12.1.21})-(\ref{1.25.1.22}), we have the following estimates
\begin{align}
&\l t\r^{-1}|\la|+|L\la|+\l t\r^{-1}|\Lb \la|\les \l t\r^{-1+\delta}\Delta_0,\label{2.4.4.22}\\
&\|\sn \la\|_{L_\omega^4}\les \l t\r^{-1+\delta}\Delta_0, \|(\tir\cir\sn)^2 \la\|_{L_\omega^4}+|\tir\sn\la|\les \Delta_0\l t \r^{\frac{1}{4}+\delta}, \label{5.8.1.21}
\end{align}
\begin{equation}\label{2.11.3.22}\left\{
\begin{array}{lll}
\sn_Y \sn\la, \sn Y\la=O(\l t\r^{-2+\delta}\Delta_0)_{L_\omega^4}, Y=L; O(\l t\r^{-\frac{7}{4}+\delta}\Delta_0)_{L_\omega^4}, Y=\sn\\
\| LY\la, Y L\la\|_{L_\omega^4}\les \l t\r^{-2+\max(-\vs(Y), 0)+\delta}\Delta_0, \quad Y=L\,\Lb, \sn
\end{array}\right.
\end{equation}
\begin{align}
&|\sn_S \sn\la, \sn S\la|\les \l t\r^{-\frac{3}{4}+\delta}\Delta_0\label{12.19.2.23}\\
&\|(\tir\snc)^{l+\le 1}\sn \la\|_{L^2_u L_\omega^2}\les \l t\r^{-1+\delta+\frac{l}{2}}\Delta_0,\, l=0,1\label{2.25.1.22}\\
&\|\sn^3\la\|_{L_\omega^2}\les \l t\r^{-2+\delta}\Delta_0 \label{6.10.1.24}\\
&\|\sn_L\sn \la, \sn L\la\|_{L^2_\Sigma}+\|\tir\sn^2\la\|_{L^2_u L_\omega^2}\les \l t\r^{\delta-1}\Delta_0\label{8.24.6.23}\\
&\|\bb^{-\f12}\sn_X \sn_L \sn\la\|_{L^2_\Sigma}\les \l t\r^{-1+\delta}\Delta_0,\label{2.19.5.22}\\
&\|\bb^{-\f12}\sn_X^n\sn\la\|_{L^2_\Sigma}\les \l t\r^{\delta+\f12(n-1)(1-\vs^+(X^n))}\Delta_0,\quad 1\le n\le 2\label{7.18.7.22}\\
&\|X^{l}Y\la\|_{L_u^2 L_\omega^2}\les \l t\r^{\delta-\max(\vs(Y),0)}\Delta_0, Y=L, \Lb,\, l\le 2\label{9.21.1.22}
\end{align} 
\end{proposition}
\begin{remark}
The second estimate in (\ref{5.8.1.21}) improves the second estimate in the assumption (\ref{1.25.1.22}).
\end{remark}
  The lower order result below will be proved as a consequence of (\ref{2.4.4.22})-(\ref{2.11.3.22}) in Proposition \ref{10.16.1.22}. It will be used to prove (\ref{2.19.5.22}) and (\ref{7.18.7.22}) in Proposition \ref{10.16.1.22}.
\begin{lemma}\label{3.17.2.22} 
 With $X=S, \Omega$, using the results (\ref{2.4.4.22})-(\ref{2.11.3.22}), there hold for $p=2$ or $\infty$
\begin{equation}\label{8.23.9.23}
\begin{split}
 &\sn_X\pio_{AL}, \sn_X\pio_{A\Lb}=O(\l t\r^{-\frac{3}{4}+\delta}\Delta_0^{\f12+\frac{1}{p}})_{L^p_u L_\omega^4}, \pio_{AL}=O(\l t\r^{-\frac{3}{4}+\delta}\Delta_0^\f12)\\
 &\sn_S\pio_{AL}, \l t\r^{-\delta}\sn_S\pio_{A\Lb}=O(\l t\r^{-1+\delta})_{L_\omega^4}.
 \end{split}
\end{equation}
Consequently, for $Y=L, \Lb, \sn$, $p=2$ or $\infty$
\begin{align} 
\sn_X(\sn_Y\Omega_B)-O(1)&=\left\{\begin{array}{lll}
 O(\Delta_0^{\f12+\frac{1}{p}}\l t\r^{\delta-\frac{3}{4}})_{L_u^p L^4_\omega},\mbox{ if }Y\neq\Lb\\
O(\Delta_0^\f12 \l t\r^{\delta-\vs(X)(1-\delta)})_{L_\omega^4}.
\end{array}\right.\label{4.22.4.22}\\
\sn_\Omega\sn \Omega-O(1)&=O(\l t\r^{\delta-1}\Delta_0)_{L_\omega^4}\label{7.11.1.24}
\end{align}

It follows by further assuming (\ref{7.18.7.22}) and (\ref{9.21.1.22}) that
\begin{equation}\label{3.28.1.24}
\sn_X^2(\pio_{A L}, \pio_{A\Lb})=O(\l t\r^{\delta-1+\f12(1-\vs^+(X^2))+\vs^+(X^2)\delta}\Delta_0)_{L^2_u L_\omega^2}
\end{equation}
and consequently
\begin{align}\label{9.1.1.21}
\begin{split}
&\sn_X^2 \sn_Y {\Omega^A}=O(1)+O(\l t\r^{-\f12+\delta+\f12\max(-\vs(Y), 0)}\Delta_0)_{L^2_u L_\omega^2}, Y=\Lb, L, e_1, e_2\\
&\sn_\Omega^2 \sn \Omega^A=O(1)+O(\l t\r^{-\frac{3}{4}+\delta}\Delta_0)_{L_\omega^4}
\end{split}
\end{align}

\end{lemma}
\begin{remark}
In order to apply the above lemma, the assumption (\ref{8.23.9.23}) will be proved after (\ref{2.4.4.22})-(\ref{2.11.3.22}) are established.Then (\ref{4.22.4.22}) will be used to prove  (\ref{2.19.5.22}) and afterwards. (\ref{9.1.1.21}) will be used for showing the highest order $L^2_\Sigma$ estimate in Lemma \ref{8.28.9.23} and later sections.    
\end{remark}
\begin{proof}[Proof of Lemma \ref{3.17.2.22}]
The estimates in (\ref{2.4.4.22}) and (\ref{5.8.1.21}) imply 
  \begin{align}
&|\la, S\la|\les \l t\r^\delta\Delta_0, |\Omega\la|\les \l t\r^{\delta+\frac{1}{4}}\Delta_0, \|\Omega \la\|_{L_\omega^4}\les \l t\r^\delta\Delta_0.\label{8.23.10.23}
\end{align}
Now we prove (\ref{8.23.9.23}). Differentiating the last line in (\ref{8.23.11.23}), we write, with $v^*$ the short-hand notation for ${}\rp{a}v^*=v^i \tensor{\ud\ep}{^{a}_i_j}e_A^j$,
\begin{equation}\label{3.28.2.24}
\begin{split}
&\sn_X\pio_{A\bN}-\sn_X(\ud\bA\la)=\sn_X\sn(c^{-1}\la),\\
& \sn_X\pio_{\bT A}+\sn_X(\ud\bA\la)=\sn_X(c^{-2} v^*)+\sn_X((\la\bA_{g,1}+\bAn)\Omega).
\end{split}
\end{equation}
Note, due to  the estimate of $\|\sn_X \ud\bA\|_{L^4_\omega}$ in (\ref{3.11.3.21}), (\ref{3.16.1.22}) and using (\ref{8.23.10.23}), 
\begin{align*}
&\|\sn_X (\ud \bA \la)\|_{L_\omega^4}\les \l t\r^{-1+2\delta}\Delta_0. 
\end{align*} 
(\ref{8.23.9.23}) follows by using $|\sn_X(c^{-2} v^*)|\les|X^{\le 1}v|=O(\l t\r^{-1+\delta}\Delta_0^\f12), O(\l t\r^{-1+\delta}\Delta_0)_{L_u^2 L_\omega^\infty}$ (due to (\ref{3.28.3.24}) and Lemma \ref{5.13.11.21}), (\ref{2.11.3.22}), (\ref{3.11.3.21}), (\ref{3.6.2.21}) and the above estimate. 

With $Y=L,\Lb, \sn$, in view of (\ref{5.6.01.21})-(\ref{5.6.03.21}) and $\tr\eta=[L\Phi]$, we write with $l=1,2$
\begin{align}
\sn^l_X \sn_Y {}\rp{a}\Omega&=\sn_X^l\Big(\chi\c {}\rp{a}\Omega+\max(-\vs(Y), 0)(\chib\c {}\rp{a}\Omega+{}\rp{a}\pi_{A\Lb})+\la([L\Phi]+\chi+\eh)c^{-1}+\sn \varrho\c{}\rp{a}\Omega\nn\\
&+\Pi\c \tensor{\ud\ep}{^{(a)}_{ij}}\c\Pi+{}\rp{a}\pi_{AL}\Big).\label{1.13.1.23}
\end{align}
Note, due to (\ref{11.5.3.23}) and (\ref{1.25.2.22}), $\sn_X \tensor{\ud\ep}{^{(a)}_{AB}}, \tir^{-1}\sn_X\Omega=O(1)$.  Using (\ref{8.23.2.23}), (\ref{3.6.2.21}), (\ref{8.23.10.23}) and (\ref{8.23.9.23}), we obtain (\ref{4.22.4.22}). 

Similarly, in view of (\ref{5.6.03.21}), we apply (\ref{8.23.10.23}), (\ref{3.6.2.21}) and (\ref{3.11.3.21}) to derive
\begin{equation*}
\sn_\Omega \sn\Omega=O(1)+\sn_\Omega(\la \thetac)+\sn_\Omega(\bA_{g,1}\c \Omega)=O(\l t\r^{-1+\delta}\Delta_0)_{L_\omega^4}
\end{equation*}
as stated in (\ref{7.11.1.24}).

Using (\ref{8.23.10.23}), (\ref{2.11.3.22}) and the $n=1$ case in (\ref{7.18.7.22}) we derive by using (\ref{3.16.1.22}) and (\ref{zeh}) that
\begin{equation*}
\sn_X^2(\ud \bA\la)=O(\l t\r^{-1+2\delta}\Delta_0^2)_{L_u^2 L_\omega^2}.
\end{equation*}
(\ref{3.28.1.24}) follows by using the above estimate, (\ref{8.23.10.23}), (\ref{3.28.2.24}), (\ref{7.18.7.22}), (\ref{L2BA2}) and (\ref{3.28.3.24}). 

The first estimate of (\ref{9.1.1.21}) can be obtained by using (\ref{3.28.1.24}), (\ref{8.23.10.23}), (\ref{1.13.1.23}), (\ref{7.18.7.22}), (\ref{9.21.1.22}), (\ref{8.23.1.23}), (\ref{3.6.2.21}), (\ref{11.5.3.23}), (\ref{L2BA2}) and  (\ref{L2conndrv}). Similarly, the second estimate can be derived in view of  (\ref{5.6.03.21}), (\ref{5.8.1.21}), (\ref{3.11.3.21}), (\ref{ConnH}), (\ref{L2conndrv}) and (\ref{10.17.1.23}). Hence the proof of Lemma \ref{3.17.2.22} is complete.
\end{proof}
\begin{proof}[Proof of Proposition \ref{10.16.1.22}] (\ref{2.4.4.22}) can be obtained by using (\ref{3.18.1.21})-(\ref{3.22.5.21}), with the help of (\ref{3.6.2.21}).
 Using (\ref{11.19.1.23}), (\ref{12.17.1.23}), (\ref{3.6.2.21}), (\ref{3.11.3.21}) and (\ref{2.4.4.22}) we bound with $2\le p\le 4$ that
 \begin{align}
 \|\snc^2\la\|_{L^p_\omega}&+\tir^{-1}\|\snc \la\|_{L_\omega^p}\les \|\cir\sD \la\|_{L_\omega^p}\nn\\
 &\les \|\tir\sn\bA_b, \tir^{-1}\Omega \log c, \tir^{-1} (\bA_b+\hat\theta)\la\|_{L^p_\omega}\les \l t\r^{-\frac{7}{4}+\delta}\Delta_0\label{1.7.3.23}.
 \end{align}
This gives (\ref{5.8.1.21}) with the help of Sobolev embedding on spheres.
Moreover, for higher order estimates, we apply (\ref{7.16.3.24}) to derive
\begin{align}
\|\snc^3\la\|_{L^2_\omega}&+\tir^{-1}\|\snc^2\la\|_{L^2_\omega}\les\|\snc \sDc \la\|_{L^2_\omega}\nn\\
&\les \|\tir \sn^2 \bA_b, \tir^{-1}\sn \Omega(c^{-1}), \tir^{-1} \sn((\bA_b+\hat\theta)\la)\|_{L^2_\omega}+\l t\r^{-\frac{19}{4}+3\delta}\Delta_0^\frac{5}{2},\label{1.7.2.23}
\end{align}
where the last error, which  arises from $\sn\bA\c \bA\la =O(\l t\r^{-\frac{19}{4}+3\delta}\Delta_0^\frac{5}{2})_{L_\omega^2}$ (due to (\ref{3.11.3.21}) and (\ref{8.23.10.23})), is negligible. 
 Due to $\Ga(\gac)-\Ga(\ga)=\p_\omega\log c$, using (\ref{5.8.1.21}), (\ref{L2conndrv}) and (\ref{3.6.2.21}), taking $L^2_u$ of the above estimate leads to 
\begin{align*}
\|\sn^3\la\|_{L_u^2 L^2_\omega}+\l t\r^{-1}\|\sn^2\la\|_{L_u^2 L^2_\omega}&\les \|\tir \sn^2 \bA_b,\tir^{-1}\sn\Omega(c^{-1}), \tir^{-1} \sn((\bA_b+\hat\theta)\la)\|_{L_u^2 L^2_\omega}\les \l t\r^{-\frac{5}{2}+\delta}\Delta_0,
\end{align*}
as stated in (\ref{2.25.1.22}).

(\ref{6.10.1.24}) can be obtained by using (\ref{5.8.1.21}), (\ref{1.7.2.23}), (\ref{ConnH}), (\ref{3.6.2.21}) and (\ref{3.11.3.21}).

Next we prove (\ref{2.11.3.22}). The $Y=\sn$ case in the first line follows from (\ref{5.8.1.21}) directly. It suffices to consider $Y=L$.
Using (\ref{3.22.5.21}) and (\ref{cmu_2}), symbolically, we write
\begin{equation}\label{8.24.7.23}
\sn_L \sn\la+\chi\sn \la, \sn L\la=\sn\Omega c.
\end{equation}
Using (\ref{3.6.2.21}), (\ref{L4BA1}) and (\ref{5.8.1.21}), we obtained the first line in (\ref{2.11.3.22}). We also obtain (\ref{12.19.2.23}) by using (\ref{3.6.2.21}) and the last estimate in (\ref{5.8.1.21}).

 Using (\ref{3.19.2}), (\ref{3.22.5.21}) and (\ref{3.18.1.21}), we bound
\begin{align*}
|L\Lb \la, \Lb L\la|&\le |\Lb L\la|+|[L, \Lb]\la|\les |\Lb \Omega c|+|\ud\bA||\sn\la|+|k_{\bN\bN}||\Omega\rp{a}\log (c \bb)|.
\end{align*}
Using (\ref{3.6.2.21}), (\ref{3.11.3.21}), $k_{\bN\bN}=O(\l t\r^{-1})$ and (\ref{5.8.1.21}), we obtain the estimate for $Y=\Lb$ in (\ref{2.11.3.22}). If $Y=L$, using (\ref{3.22.5.21}) and (\ref{L4BA1}), we have $LL\la=O(\l t\r^{-2+\delta})_{L_\omega^4}$. If $Y=\sn$, the result has been obtained in the first set of estimates in (\ref{2.11.3.22}).

Using (\ref{8.24.7.23}), (\ref{3.6.2.21}) and (\ref{L2BA2}), we obtain (\ref{8.24.6.23}) by also combining with the estimate (\ref{2.25.1.22}).
 This also gives the $n=1$ case in (\ref{7.18.7.22}).

Next we prove (\ref{2.19.5.22}), for which we will use (\ref{4.22.4.22}). Applying $X$-derivative to (\ref{8.24.7.23}) yields
\begin{align*}
\sn_X \sn_L\sn\la&=\sn_X\big(\chi\sn \la+\sn\Omega c\big).
\end{align*}
Note that, due to (\ref{3.11.3.21}), $\sn_X\big(\bA_{g,1}\c \Omega \log c\big)=O(\l t\r^{-3+2\delta}\Delta_0^2)_{L_u^2 L_\omega^2}$.
Then the leading term in the above can be bounded by using (\ref{L2BA2}), (\ref{1.25.2.22}) and (\ref{4.22.4.22})  
\begin{equation}\label{2.20.3.22}
\|\sn_X \sn \Omega c\|_{L^2_\Sigma}\les \Delta_0 \l t\r^{\delta-1}.
\end{equation}
Using (\ref{8.24.4.23}), (\ref{8.24.6.23}) and (\ref{5.8.1.21}), we have
\begin{align*}
\sn_X\big(\chi\c\sn\la\big)&=O(\l t\r^{-1}) \sn_X^{\le 1}\sn\la +O(\l t\r^{-\frac{7}{4}+\delta}\Delta_0^\f12)_{L_\omega^4} \sn\la\\
&=O(\l t\r^{-2+\delta}\Delta_0)_{L_u^2 L_\omega^2}.
\end{align*}
 Hence, (\ref{2.19.5.22}) is proved.

Next, we prove the $n=2$ case in (\ref{7.18.7.22}). It is direct to write
\begin{equation*}
\sn_L \sn_\Omega \sn\la=\sn_\Omega \sn_L\sn\la+[L, \Omega]\sn\la.
\end{equation*}
 Using (\ref{3.21.1.23}) and (\ref{2.25.1.22}), we derive
\begin{align}
\|[S, \Omega]\sn \la\|_{L^2_u L_\omega^2}&\les\Delta_0^\f12\|\tir\sn^2\la\|_{L^2_u L_\omega^2}+\l t\r^{-\frac{3}{4}+\delta}\Delta_0^\f12\|\sn\la\|_{L_u^2 L_\omega^2}\nn\\
&\les \l t\r^{-1+\delta}\Delta_0^\frac{3}{2}.\label{8.23.7.23}
 \end{align}
 Thus using (\ref{2.19.5.22}), we conclude
 $$
 \|\bb^{-\f12}(\sn_S\sn_\Omega\sn\la, \sn_\Omega \sn_S \sn\la)\|_{L^2_\Sigma}\les \l t\r^\delta\Delta_0.
 $$
The case of $X^2=S^2$ can be obtained by (\ref{2.19.5.22}) and (\ref{8.24.6.23}). 
By using (\ref{2.25.1.22}) and (\ref{9.8.2.22}) we have
\begin{equation*}
\|\bb^{-\frac{1}{2}}\sn_\Omega^2 \sn\la\|_{L^2_\Sigma}\les \l t\r^{\delta+\f12}\Delta_0.
\end{equation*} 
Thus (\ref{7.18.7.22}) is proved.

Finally, using Proposition \ref{2.19.4.22}, (\ref{9.21.1.22}) is a straightforward consequence of the estimates for $\ud\bA$ in (\ref{3.16.1.22}), (\ref{L2BA2}), (\ref{3.6.2.21}) and (\ref{4.22.4.22}).
\end{proof}
Using (\ref{9.8.2.22}), (\ref{4.22.4.22}) and (\ref{9.1.1.21}), the covariant derivatives $\tir\sn$ appeared in Assumption \ref{5.13.11.21} and Lemma \ref{5.13.11.21} can be replaced by $\sn_\Omega$.
For future reference, we give a set of product estimates below. 
\begin{lemma}\label{8.28.9.23}
Let $X\in \Omega, S$. Under the assumptions (\ref{3.12.1.21})-(\ref{6.5.1.21}), with $i=1,2$, we have 
\begin{align}
\sn_X^{\le 1}(\bA\c \bA_{g,i})&=O(\l t\r^{-\frac{7}{2}-\frac{2-i}{4}+2\delta}\Delta_0^\frac{3}{2})_{L^4_\omega}, \sn_X^{\le 1}(\bA\c (\bAn, \bA_g))= O(\l t\r^{-\frac{15}{4}+2\delta}\Delta_0^\frac{3}{2})_{L_u^2 L_\omega^2}\label{8.3.4.23}\\
\sn_X^{\le 1}(\bA_g\c\bA_g)&= O(\l t\r^{-\frac{15}{4}+2\delta}\Delta_0^2)_{L_u^2 L_\omega^2}, O(\l t\r^{-\frac{7}{2}+2\delta}\Delta_0^2)_{L_\omega^4}\label{8.3.3.23}\\
\sn_X^{\le 1}(\bA_{g,1}^2)&= O(\l t\r^{-3+2\delta}\Delta_0^2)_{L^2_\Sigma}, O(\l t\r^{-4+2\delta}\Delta_0^2)_{L_\omega^4}\label{8.28.1.23}\\
\sn_X^2(\bA\c \bA_{g,i})&=O(\l t\r^{-\frac{13}{4}-\frac{2-i}{4}+2\delta}\Delta_0^\frac{3}{2})_{L^2_u L_\omega^2}\label{8.28.2.23}\\
\sn_X^{\le 2}(\bA\c \bA)&=O(\l t\r^{-\frac{13}{4}+2\delta}\Delta_0^\frac{3}{2})_{L^2_u L_\omega^2}, \sn_X^{\le 2}(\bAn\c\bAn)=O(\l t\r^{-4+2\delta}\Delta_0^\frac{3}{2})_{L^2_u L_\omega^2}\label{2.27.1.24}\\
\sn_X^{\le 2}(\bA_{g,1}\c \bA_{g,1})&=  O(\l t\r^{-4+2\delta}\Delta_0^\frac{3}{2})_{L^2_u L_\omega^2}, O(\l t\r^{-\frac{15}{4}+2\delta}\Delta_0^2)_{L_\omega^4}\label{12.22.4.23}\\
\displaybreak[0]
\sn_X^3(\bA_{g,1}\c\bA_{g,1})&=O(\l t\r^{-\frac{5}{2}+2\delta}\Delta_0^2)_{L^2_\Sigma}\label{12.22.2.23}\\
 \sn_X^{\le 2}(\ud\bA\c \bA_{g,1})&=O(\l t\r^{-3+2\delta}\Delta_0^2)_{L^2_u L_\omega^2}, \sn_X^2(\ud \bA^2)=O(\l t\r^{-2+2\delta}\Delta_0^2)_{L_u^2 L_\omega^2}\label{1.30.1.24}\\
 \sn_X^{\le 2}(\ud \bA\c \bA)&=O(\l t\r^{-\frac{5}{2}+2\delta}\Delta_0^\frac{3}{2})_{L_u^2 L_\omega^2}.\label{2.29.6.24}
\end{align}
\end{lemma}
\begin{proof} 
Using (\ref{9.8.2.22}), (\ref{4.22.4.22}) and (\ref{9.1.1.21}), we summarize from (\ref{3.11.3.21}), (\ref{L2BA2}), (\ref{3.16.1.22}) and (\ref{L2conndrv}) that
\begin{align*}
\sn_X^{\le 1} \bA&=O(\l t\r^{-\frac{7}{4}+\delta}\Delta_0^\f12)_{L_\omega^4}, O(\l t\r^{-2+\delta}\Delta_0)_{L^2_u L_\omega^2} \\
\sn_X^{\le 1}\bA_g&=O(\l t\r^{-\frac{7}{4}+\delta}\Delta_0)_{L_\omega^4}, O(\l t\r^{-2+\delta}\Delta_0)_{L^2_u L_\omega^2}\\
\sn_X^{\le 1+l}\bA_{g,1}&=O(\l t\r^{-2+\frac{l}{4}+\delta}\Delta_0)_{L_\omega^4}, O(\l t\r^{-1+\delta}\Delta_0)_{L^2_\Sigma},\, l=0,1\\
\sn_X^3\bA_{g,1}&=O(\l t\r^{-\f12+\delta}\Delta_0)_{L^2_\Sigma}\\
\sn_X^{2}\bA_b, & \sn_X^{2}\bA_{g,2}= O(\l t\r^{-\frac{3}{2}+\delta}\Delta_0)_{L^2_u L_\omega^2}\\
\sn_X^{\le 2}\bAn&=O(\l t\r^{-1+\delta}\Delta_0)_{L^2_\Sigma}, \sn_X^{\le 1}\bAn=O(\l t\r^{-2+\delta}\Delta_0^\f12)_{L_\omega^4}\\
\sn_X^{\le 2}\ud\bA&=O(\l t\r^{-1+\delta}\Delta_0)_{L_u^2 L_\omega^2}, \sn_X^{\le 1}\ud \bA=O(\l t\r^{-1+\delta}\Delta_0)_{L_\omega^4}. 
\end{align*}
Lemma \ref{8.28.9.23} can be obtained by using the above estimates with the help of (\ref{3.6.2.21}) and (\ref{3.11.3.21}).
\end{proof}

\subsection{Estimates of null forms} \begin{footnote}{
 From now on, until the completion of the top order weighted energy estimates, all estimates will be established under the assumptions (\ref{3.12.1.21})-(\ref{6.5.1.21}), which may not be explicitly stated repetitively.}
\end{footnote}
Next we derive estimates of null forms.

\begin{proposition}\label{6.24.10.23}
Assume  (\ref{3.12.1.21})-(\ref{6.5.1.21}) hold. For the null form $\N(\Phi,\bp\Phi)$ and $\wt \N(\Phi, \bp\Phi)$ in Proposition \ref{geonul_5.23_23}, with $X\in\{\Omega, S\}$, the following estimates hold with $l=0,1$ and $n\le 2$
\begin{equation}\label{6.7.3.23}
\begin{split}
X^n\N(\Phi, \bp \Phi)&=X^n([L\Phi]\c([\Lb\Phi]+[L\Phi]))+O(\l t\r^{-\frac{15}{4}+2\delta}\Delta_0^2)_{L_\omega^4}, \\
X^{n+l}\N(\Phi, \bp \Phi)&=X^{n+l}([L\Phi]\c ([\Lb\Phi]+[L\Phi]))+O(\l t\r^{-3+\f12 l+2\delta}\Delta_0^2)_{L_\Sigma^2}\\
X^{n+l}\sC(\varrho, v)&=X^{n+l}(L\Phi[\Lb\Phi]+[L\Phi]\Lb\Phi)+O(\l t\r^{-3+\f12 l +2\delta} \Delta_0^\frac{3}{2})_{L^2_\Sigma}.
\end{split}
\end{equation}
\begin{align}
X^{n}\N(\Phi, \bp \Phi)&=[\Lb\Phi]X^{n}[L\Phi]+\sum_{a+b=n, b\ge 1}\vs^-(X^b) X^a[L\Phi][\Lb\Phi]\nn\\
&+O\big(\l t\r^{-3+2\delta}\Delta_0^{1+\f12(1-\vs^-(X^{n}))}\big)_{L_\omega^4}\label{10.22.5.22}\\
X^{n+l}\N(\Phi, \bp \Phi)&=O(1)X^{\le n+l}[L\Phi]\c[\Lb\Phi]+\vs^-(X^{n+l})O(\l t\r^{-\frac{15}{4}+\frac{l}{4}+2\delta})\Delta_0^{\frac{3}{2}})_{L_u^2 L_\omega^2}\nn\\
&+(1-\vs^-(X^{n+l}))O(\bb^{-\f12 l}\l t\r^{-2+2\delta}\Delta_0^\frac{3}{2})_{L^2_\Sigma}.\label{6.6.1.23}\\
X^{n+l} \wt \N(\Phi, \bp\Phi), & X^{n+l}\sC(\varrho, v)=O(1)(X^{\le n+l}L\Phi+X^{\le n+l}[L\Phi])\c [\Lb\Phi]\nn\\
\displaybreak[0]
&\qquad\qquad\qquad+\vs^-(X^{n+l})O(\l t\r^{-\frac{15}{4}+\frac{l}{4} +2\delta}\Delta_0^{\frac{3}{2}})_{L_u^2 L_\omega^2}\nn\\
&\qquad\qquad\qquad+(1-\vs^-(X^{n+l}))O(\bb^{-\f12 l}\l t\r^{-2+2\delta}\Delta_0^\frac{3}{2})_{L^2_\Sigma}\label{9.16.1.23}
\end{align}
where $O(\bb^{-\f12 }\l t\r^{-2+2\delta}\Delta_0^\frac{3}{2})_{L^2_\Sigma}=O(\l t\r^{-3+2\delta}\Delta_0^\frac{3}{2})_{L^2_u L_\omega^2}$. 
(\ref{10.22.5.22}) and (\ref{6.6.1.23}) also hold for $[X^m \sG(\varrho,v)]$, with  $[X^{\le m}L\Phi]$ added to the term $X^{\le m}[L\Phi]$ on the right-hand side.
\end{proposition}

\begin{proof}
In view of Proposition \ref{geonul_5.23_23}, symbolically we write 
\begin{align*}
X^m\N(\Phi,\bp\Phi)&=X^m([\Lb\Phi][L\Phi]+[L\Phi]^2)+X^m(\bA_{g,1}^2).
\end{align*}
Using (\ref{8.28.1.23}), (\ref{12.22.4.23}) and (\ref{12.22.2.23})  to treat the last term,
we obtain the first two lines in (\ref{6.7.3.23}). In view of (\ref{5.22.2.23}) by repeating the proof for the second line in (\ref{6.7.3.23}) we obtain with $m=n+l$
\begin{align*}
X^m \sC(\varrho, v)&=X^m(L\Phi\c [\Lb\Phi]+[L\Phi]\Lb\Phi+\sn\varrho\c\sn v)\\
&=X^m(L\Phi\c [\Lb\Phi]+[L\Phi]\Lb\Phi)+O(\l t\r^{-3+\f12 l+2\delta}\Delta_0^\frac{3}{2})_{L^2_\Sigma}, \quad n\le 2.
\end{align*}
This gives the last line in (\ref{6.7.3.23}).

To see (\ref{10.22.5.22}), using (\ref{3.6.2.21}) and (\ref{3.11.3.21}), we obtain for $n\le 2$,
\begin{align*}
X^n([L\Phi]^2)&=[L\Phi]X^n[L\Phi]+O\big(\l t\r^{-4+2\delta}\Delta_0^{1+\f12(1-\vs^-(X^n))}\big)_{L_\omega^4}\\
&=O\big(\l t\r^{-\frac{15}{4}+2\delta}\Delta_0^{1+\f12(1-\vs^-(X^n))}\big)_{L_\omega^4}. 
\end{align*}
 Taking $L^4_\omega$ norm with the help of (\ref{8.23.1.23})  and  (\ref{3.6.2.21}) yields
\begin{align*}
&X^n([L\Phi][\Lb\Phi])\\
&= [\Lb\Phi]X^n[L\Phi]+\sum_{a+b=n, b\ge 1}\vs^-(X^b) X^a[L\Phi]([\Lb\Phi]+S^{b-1}[L\Phi]+O(\l t\r^{-1+\delta}\Delta_0)_{L_\omega^4})\\
&=[\Lb\Phi]X^n[L\Phi]+\sum_{a+b=n, b\ge 1}\vs^-(X^b)X^a[L\Phi]([\Lb\Phi]+S^{b-1}[L\Phi])+O(\l t\r^{-3+2\delta}\Delta_0^\frac{3}{2})_{L_\omega^4}. 
\end{align*}
Note that the term $\vs^-(X^b) X^a[L\Phi]  S^{b-1}[L\Phi]$ can be treated similarly as for $X^n([L\Phi]^2)$. 
Thus, combining the above two estimates,   we obtain (\ref{10.22.5.22}).

Next we prove $L^2_\Sigma$ type estimates. With $m=0,1,2,3$, we write
\begin{equation}\label{8.24.1.23}
X^m([L\Phi][\Lb\Phi])=[L\Phi]X^{m}[\Lb \Phi]+[\Lb \Phi]X^{m}[L\Phi]+X^{m-1}[L\Phi]X [\Lb\Phi]+X[L\Phi]X^{m-1}[\Lb\Phi].
\end{equation}
The $0$-order case in (\ref{6.6.1.23}) holds trivially. 
Using (\ref{8.23.1.23}), symbolically, we write
\begin{align*}
X^{m-1}[L\Phi]X[\Lb\Phi]&=X^{m-1}[L\Phi]\left(\vs(X)[\Lb\Phi]+(1-\vs(X))\l t\r^{-1+\delta}\Delta_0\right).
\end{align*}
Taking $L^2_\Sigma$ norm with the help of (\ref{L2BA2}) gives
\begin{align*}
X^{m-1}[L\Phi]X[\Lb\Phi]&=\vs(X)X^{m-1}[L\Phi][\Lb\Phi]+(1-\vs(X))O(\l t\r^{-2+2\delta}\Delta_0^2)_{L^2_\Sigma}.
\end{align*}
  
Applying Lemma \ref{7.6.6.23} to the first term in (\ref{8.24.1.23}) and using (\ref{3.6.2.21}), we have with $m=n+l$
\begin{align}
[L\Phi]X^{n+l}[\Lb\Phi]&=\Big(\vs^-(X^{n+l})([\Lb \Phi]+O(\l t\r^{-1+\f12 l+\delta}\Delta_0)_{L^2_\Sigma})\nn\\
&+(1-\vs^-(X^{n+l}))O(\l t\r^\delta\Delta_0)_{L^2_\Sigma}\Big)[L\Phi]\nn\\
&=\vs^-(X^{n+l})([\Lb\Phi][L\Phi]+O(\l t\r^{-3+\f12 l+2\delta}\Delta_0^\frac{3}{2})_{L^2_\Sigma})\nn\\
&+(1-\vs^-(X^{n+l}))O(\l t\r^{2\delta-2}\Delta_0^\frac{3}{2})_{L^2_\Sigma}.\label{9.16.2.23}
\end{align}
Similarly, it follows by using (\ref{8.23.1.23}) and (\ref{3.11.3.21}) with the help of Sobolev embedding on spheres that
\begin{align*}
X[L\Phi] X^{m-1}[\Lb\Phi]&=\vs^-(X^{m-1})([\Lb\Phi]X[L\Phi]+O(\l t\r^{-\frac{15}{4}+2\delta}\Delta_0^\frac{3}{2})_{L^2_u L_\omega^2})\nn\\
&+(1-\vs^-(X^{m-1}))O(\bb^{-\f12 l}\l t\r^{2\delta-2}\Delta_0^\frac{3}{2})_{L^2_\Sigma}.
\end{align*}
Substituting the above estimates to (\ref{8.24.1.23}) implies (\ref{6.6.1.23}).

To show (\ref{9.16.1.23}), we will rely on the last estimate in (\ref{6.7.3.23}).
 Note using Lemma \ref{7.6.6.23}, (\ref{8.23.1.23}), (\ref{8.23.1.23'}) and (\ref{L2BA2})
\begin{align*}
[L\Phi]X^m\Lb\Phi&=[L\Phi](O(1)[\Lb\Phi]+O(1)\vs^-(X^{n+l})O(\l t\r^{-1+\delta+\f12l}\Delta_0)_{L_\Sigma^2}\\
&+(1-\vs^-(X^{n+l}))O(\l t\r^{\delta}\Delta_0)_{L^2_\Sigma}),
\end{align*}
where $[L\Phi]$ can be bounded by using (\ref{3.6.2.21}); $L\Phi X^m[\Lb\Phi]$ can be bounded as above, and 
\begin{align*}
&X^{m-1}[L\Phi] X\Lb\Phi, X^{m-1}\Lb\Phi X[L\Phi]\\
&=[\Lb\Phi]X^{\le m-1}[L\Phi]+(1-\vs^-(X^{m-1}))O(\bb^{-\f12l}\l t\r^{-2+2\delta}\Delta_0^\frac{3}{2})_{L^2_\Sigma}\\
&+\vs^-(X^{m-1})O(\l t\r^{-\frac{15}{4} +2\delta}\Delta_0^{\frac{3}{2}})_{L_u^2 L_\omega^2}.
\end{align*}
Summarizing the above calculations gives (\ref{9.16.1.23}). To avoid repetition, we omit the proof for the estimates (\ref{10.22.5.22}) and (\ref{6.6.1.23}) for $[X^n \sG(\varrho,v)]$. They can be bounded similarly, with the additional term $[X^{\le n}L\Phi]$ added to $X^{\le n}[L\Phi]$ on the right-hand side.



\end{proof}

Using Lemma \ref{6.30.4.23}, we give a useful result for future reference.
\begin{lemma}\label{10.10.3.23}
With $X^n=X_n\cdots X_1$, $X_i\in \{S, \Omega\}$, there hold for scalar functions $f$,  $v$ and $\Phi$ 
\begin{align}
X^nL f&=LX^n f+O(\l t\r^{-1})X_n^{\le 1}\cdots X_1 f\label{8.30.4.23}\\
&+\max(n-1, 0)O(\l t\r^{-\frac{3}{4}+\delta}\Delta_0^\f12)_{L_\omega^4}\sn f, \, n=1,2,\nn\\
\Sc(X^n Lv)&=[LX^n v]+O(\l t\r^{-1})[X_n^{\le 1}\cdots X_1 v]+O(\l t\r^{-\frac{3}{4}+\delta}\Delta_0^{1-\f12\vs^+(X^n)})X^{\le n-1}[L v]\nn\\
&+O(1)\sn_X^{\le n-1}[\sn \Phi]+\max(n-1, 0)O(\l t\r^{-\frac{3}{4}+\delta}\Delta_0^\f12)_{L_\omega^4}[\sn \Phi],\,\, n=1,2\label{6.30.2.23}
\end{align}
\begin{align}
[X^3, L]\Phi&=\l t\r^{-1}X_3^{\le 1} X_2 X_1\Phi+O(\l t\r^{-\frac{5}{2}+2\delta}\Delta_0^\frac{3}{2})_{L_u^2 L_\omega^2}\label{7.1.1.23}\\ 
\|\sn_X^n[\tir \sn \Phi], \sn_X^n\eta(\Omega)\|_{L^2_\Sigma}&\les \|X^{\le n}\Omega\Phi\|_{L^2_\Sigma}+\l t\r^{-\frac{3}{4}+\frac{\max(n-2,0)}{4}+2\delta}\log \l t\r^\f12\Delta_0^\frac{3}{2}, n=1,2,3. \label{7.25.5.23}
\end{align}
\end{lemma}
\begin{proof}
 We will frequently use the estimates of $\sn_X^l \pio_{LA}$ in Lemma \ref{3.17.2.22}.
In view of (\ref{8.23.9.23}), we can apply (\ref{7.17.7.21}) to scalar functions $f$ to obtain 
\begin{align*}
[X_2 X_1, L] f&= O(\l t\r^{-1})X_2^{\le 1}X_1 f+O(\l t\r^{-\frac{3}{4}+\delta}\Delta_0^\f12)_{L_\omega^4}\sn f.
\end{align*} 
 This gives (\ref{8.30.4.23}) for the $n=2$ case. We apply (\ref{7.17.6.21}) instead if $n=1$. 

We next write by using Lemma \ref{6.30.4.23} with $X_1=L$,
\begin{align*}
\Sc(X^n Lv)&=[X^n Lv]+O(1)\sn_X^{\le n-1}[\sn \Phi]+O(\l t\r^{-\frac{3}{4}+\delta}\Delta_0^{1-\f12\vs^+(X^n)})X^{\le n-1}[L v]\nn\\
&+\max(n-1,0)\{O(\l t\r^{-\frac{3}{4}+\delta}\Delta_0^\f12)_{L_\omega^4}[\sn \Phi]\}.
\end{align*}
 Substituting (\ref{8.30.4.23}) to the above symbolic formula gives (\ref{6.30.2.23}).
 
 For (\ref{7.1.1.23}), we first write for scalars $f$
\begin{equation*}
[X_3 X_2 X_1, L]f=X_3 X_2[X_1, L]f+X_3[X_2, L]X_1 f+[X_3, L] X_2 X_1 f.
\end{equation*}
Using Lemma \ref{6.9.3.23} (1), we bound
\begin{align*}
[X^3,L]f= X_3 X_2[X_1, L]f+\l t\r^{-1}X_3^{\le 1}X_2 X_1 f+O(\l t\r^{-\frac{3}{4}+\delta}\Delta_0^\f12)_{L_\omega^4}\sn X_1 f.
\end{align*}
Note that due to (\ref{5.13.10.21}) if $X_1=\Omega$, 
\begin{align*}
X^2[\Omega, L]f=X^2(\pio_{LA} \sn f)=\sn_X^2\pio_{LA}\sn f+\pio_{LA}\sn_X^2\sn f+\sn_X\pio_{LA}\sn_X \sn f.
\end{align*}
 Thus using (\ref{3.6.2.21}), (\ref{3.11.3.21}), (\ref{3.12.1.21}) and Lemma \ref{3.17.2.22}, we infer
\begin{align*}
\|X^2[\Omega, L]\Phi\|_{L^2_u L_\omega^2}&\les \|\sn_X^2\pio_{LA}\|_{L^2_u L_\omega^2}\|\sn \Phi\|_{L_\omega^\infty}+\| \sn_X\sn \Phi\|_{L_u^2 L_\omega^4}\|\sn_X \pio_{LA}\|_{L_\omega^4}\\
&+\|\pio_{AL}\|_{L_\omega^\infty}\|\sn_X^2 \sn \Phi\|_{L^2_u L_\omega^2}\\
&\les \l t\r^{-\frac{5}{2}+2\delta}\Delta_0^\frac{3}{2}.
\end{align*}
Combining  $X^2[S,L]\Phi=X^2 L\Phi$, we deduce by using (\ref{3.11.3.21}) that 
\begin{align*}
[X^3,L] \Phi&=\l t\r^{-1}X_3^{\le 1}X_2 X_1 \Phi+O(\l t\r^{-\frac{5}{2}+2\delta}\Delta_0^\frac{3}{2})_{L_u^2 L_\omega^2}.
\end{align*}
This is (\ref{7.1.1.23}). 

Next we prove the following estimates.
\begin{align}\label{7.25.5.23'}
\begin{split}
\|\sn_X^n[\tir \sn \Phi]\|_{L^2_\Sigma}&\les \|X^{\le n}\Omega\Phi\|_{L^2_\Sigma}+\|\sn_X^{\le n-1}\eta(\Omega)\|_{L^2_\Sigma}+\l t\r^{-\frac{3}{4}+2\delta}\log \l t\r^\f12\Delta_0^\frac{3}{2};\\
\|\sn_X^n(\eta(\Omega))\|_{L^2_\Sigma}&\les \|X^{\le n}\Omega\Phi\|_{L^2_\Sigma}+\|\sn_X^{\le n-1}[\Omega v]\|_{L^2_\Sigma}+\l t\r^{-\frac{3}{4}+2\delta}\log \l t\r^\f12\Delta_0^\frac{3}{2}
\end{split}
\end{align}
(\ref{7.25.5.23}) follows directly from using the above estimate. 

 Note due to Lemma \ref{3.17.2.22}, we can bound with $1\le n\le 2$ 
\begin{align*}
\sn_X^n[\tir \sn \Phi]=O(1) X^{\le n}[\Omega \Phi]+\max(n-1, 0)O(\Delta_0^\f12\l t\r^{\delta-\frac{3}{4}})_{L_\omega^4}[\Omega\Phi].
\end{align*} 
 Note $\|\tir[\Omega\Phi]\|_{L^2_u L_\omega^4}\les \l t\r^\delta\Delta_0$, due to (\ref{3.11.3.21}). Using this estimate, if $\Phi=\varrho$ in the above, then we can obtain (\ref{7.25.5.23}) in this case;
if $\Phi=v$, by using Lemma \ref{6.30.4.23}, with $X_1=\Omega$ therein, we write
\begin{align*}
X^n[\Omega v]&=[X^n\Omega v]+O(1)\sn_X^{\le n-1}(\eta(\Omega))+O(\l t\r^{-\frac{3}{4}+\delta}\Delta_0^{1-\f12\vs^+(X^n)})X^{\le n-1}[\Omega v]\\
&+\max(n-1,0)O(\l t\r^{-\frac{3}{4}+\delta}\Delta_0^\f12)_{L_\omega^4}(\eta(\Omega)).
\end{align*}
 Then we conclude the first line in (\ref{7.25.5.23'}) by using (\ref{L2BA2}) for $n\le 2$. The second line in (\ref{7.25.5.23'}) and  the $n=3$ case in (\ref{7.25.5.23}) can be proved similarly by using Lemma \ref{6.30.4.23}, with details omitted here. 
 
\end{proof}

\subsection{Derivative control on $\tr\chi$, $\chih$ and $\zeta$}
In this subsection, we give the comparison estimates of the connection coefficients in terms of the derivatives of $\Phi$. 
We first give the lowest order comparison estimates.
\begin{proposition}\label{7.13.4.22}
Under the assumptions (\ref{3.12.1.21})-(\ref{6.5.1.21}), there hold
 \begin{align}
& \|\tir^2(\tr\chi-\frac{2}{\tir})(t)\|_{L^p_\omega}\les\|(\tr\chi-\frac{2}{\tir})(0)\|_{L^p_\omega}+\Delta_0^2\nn\\
&\qquad\qquad\qquad\,\,+ \|\tir [\bar \bp\Phi],\tir^2 LL\varrho, \tir^2 \sD\varrho \|_{L_t^1 L_\omega^p}, \,\quad 2\le p\le 4.\label{7.10.5.22}\\
&\|\la\|_{L_\omega^4}\les \int_0^t \|\Omega c\|_{L_\omega^4}\label{8.25.3.23}\\
&\|\bb^{-\f12}\sn \la\|_{L^2_\Sigma}\les \int_0^t \|\Omega^2 c\|_{L^2_u L_\omega^2} dt\label{9.18.2.23}\\
&\begin{array}{lll}
\|SL\la\|_{L^2_\Sigma}\les \|S\Omega c\|_{L^2_\Sigma}\\
\|\sn_L(\tir \sn \la), \Omega L\la\|_{L^2_\Sigma}\les \|\Omega^2 c\|_{L^2_\Sigma}+\l t\r^{-1+\delta}\Delta_0^\frac{3}{2}\log \l t\r^\f12 
\end{array}\label{9.21.2.23}\\
&\|L\Lb\la\|_{L^2_u L_\omega^2}\les \|\Omega k_{\bN\bN}, L\Omega c\|_{L^2_u L_\omega^2}+\|(\pio_{AL}, S(c))\sn\log (\bb c)\|_{L^2_u L_\omega^2}\label{2.15.1.24}\\  
&\|\bb^{-\frac{1}{2}}\sn^2 \la \|_{L^2_\Sigma}\les \|\bb^{-\frac{1}{2}}(\tir \sn\bA_b, \tir^{-1}\Omega \log c)\|_{L^2_\Sigma}+\l t\r^{-2+2\delta}\Delta_0^2\label{9.21.3.23}.
\end{align}
\end{proposition}
\begin{proof}
Due to (\ref{6.3.1.23}) and
\begin{align*}
L\big(\tir&(\tr\chi-\frac{2}{\tir})\big)+\f12\tir \tr\chi(\tr\chi-\frac{2}{\tir})\nn\\
&=\tir\big(L(\tr\chi-\frac{2}{\tir})+\frac{2}{\tir}(\tr\chi-\frac{2}{\tir})+\f12(\tr\chi-\frac{2}{\tir})^2\big),
\end{align*}
we bound
\begin{align}\label{7.10.7.22}
|L(v_t^\f12 \tir (\tr\chi-\frac{2}{\tir}))|&\les v_t^\f12 \tir(|\tr\chi, \fB, \zb||[\bar \bp \Phi]|+|LL\varrho, \sD\varrho|+|\chih\c \chih|),
\end{align}
which leads to
\begin{align*}
\|L(v_t^\f12 \tir (\tr\chi-\frac{2}{\tir}))\|_{L^1_t L_\omega^p}&\les \|\tir ([\bar \bp\Phi], \tir\sD \varrho, \tir LL \varrho, \tir\chih\c \chih)\|_{L_t^1 L_\omega^p}.
\end{align*}
Using (1) in Lemma \ref{5.13.11.21}, we then obtain (\ref{7.10.5.22}). 

Revisiting the proof of $\|\la\|_{L_\omega^4}$, we obtain (\ref{8.25.3.23}).  Due to (\ref{cmu_2}) and (\ref{3.22.5.21})
\begin{equation*}
\sn_L \sn\la=-\chi\sn \la+\sn\Omega c.
\end{equation*}
(\ref{9.18.2.23}) then follows by the transport lemma.  The  estimates of $X L\la$ in (\ref{9.21.2.23}) follow by using (\ref{3.22.5.21}), the  one for $\sn_L(\tir \sn \la)$ relies on the transport equation and noting $\|\tir\bA \sn\la\|_{L^2_u L_\omega^2}\les \l t\r^{-2+2\delta}\Delta_0^\frac{3}{2}$ (due to (\ref{5.8.1.21}) and (\ref{3.11.3.21})). (\ref{2.15.1.24}) can be obtained by using Proposition \ref{2.19.4.22} and (\ref{5.13.10.21}).   
 (\ref{9.21.3.23}) can be obtained from (\ref{1.7.3.23}), (\ref{L2conndrv}) and (\ref{2.4.4.22}).
\end{proof}

\begin{proposition}\label{11.4.1.22}
Under the assumptions (\ref{3.12.1.21})-(\ref{6.5.1.21}), with $X, X_1, X_2\in \{\Omega, S\}$, $X^2=X_2 X_1$, and in the region $u_0\le u_1\le u\le u_*$ there hold

(1)
With $2\le p\le 4$, 
\begin{align}
& \|\sn\sn_X^\ell\chih, \tir^{-1}\sn_X^\ell\chih)\|_{L_\omega^p}\les \|\sn_X^{\le\ell}(\sn\tr\chi, \tir^{-1} \zb, \tensor{\bR}{_{C4}^C_B})\|_{L_\omega^p} +\Delta_0^2 \l t\r^{-\frac{15}{4}+2\delta}\label{8.25.2.23}, \ell=0,1\\
 &\|\tir\sn\sn_{X}^2\chih\|_{L^2_u L_\omega^2}+\|\sn_X^2 \chih\|_{L^2_u L_\omega^2}\nn\\
&\qquad\qquad\quad\les \|\bb^{-\f12}\{\sn_X^2(\sn\tr\chi, \bR_{B4BA}, \tir^{-1}\zb), \l t\r^{-1}\sn_\Omega^{\le 1}\sn_{X_1}\chih, \sn_{X_2}^{\le 1}(\sn \chih, \tir^{-1}\chih)\}\|_{L^2_\Sigma}\nn\\
&\qquad\qquad\quad+\l t\r^{-\frac{5}{2}+2\delta}\Delta_0^\frac{3}{2}\label{2.25.2.24}\\
&\|\tir^2(\tir\sn)^l\wt{\tr\chi}(t)\|_{L_u^2 L_\omega^2}\les\|\tir^2(\tir\sn)^l\Big(\sD\varrho, \tir^{-1}[L\Phi], [L\Phi][\Lb\Phi]\Big)\|_{L_t^1 L_u^2 L_\omega^2(\D_{u_1}^t)}\nn\\
&\qquad\qquad\qquad\qquad\quad+\La_0+\Delta_0^\frac{3}{2},\quad l=1,2\label{7.21.2.22+}\\
&\|\tir^2(\tir\sn) \wt{\tr\chi}\|_{L_\omega^4}\les\|\tir^2(\tir\sn)(\sD\varrho,\tir^{-1}[L\Phi], [L\Phi][\Lb\Phi])\|_{L_t^1 L_\omega^4(\H_u^t)}+ \La_0+\Delta_0^\frac{3}{2}\label{3.16.14.24}\\
&\|\tir^{-2}\sn_L( \tir^3\sn \tr\chi)\|_{L^2_u L_\omega^2}\les\|\tir\sn\widetilde{L(\Xi_4)}, \sn[L\Phi], \tir\sn([L\Phi][\Lb\Phi])\|_{L^2_u L_\omega^2}+\l t\r^{-\frac{15}{4}+2\delta}\Delta_0^\frac{3}{2}.\label{8.5.1.22}
\end{align}
(2)With $l=0,1$ below, there hold
\begin{align}
&\sn_L(\bb^l \sn\log \bb)+\f12\tr\chi \bb^l \sn \log \bb=-\bb^l\chih\c\sn\log \bb-\sn(\bb^l k_{\bN\bN}). \label{7.17.2.24}\\
&\sn_L (\tir\sn)^m(\bb^l\sn\log \bb)+\chi (\tir \sn)^m(\bb^l\sn\log\bb)+(\tir\sn)^m \sn (\bb^l k_{\bN\bN})\nn\\
&=\sum_{a=1}^m(\tir \sn)^{m-a}\big((\frac{1}{\tir}-\f12\tr\chi-\chih)(\tir \sn)^a (\bb^l \sn\log \bb)\big)\nn\\
&-\sum_{a=0}^{m-1}(\tir \sn)^{m-a}(\f12\tr\chi+\chih)(\tir\sn)^a (\bb^l\sn\log \bb)\label{12.17.3.23}\\
&+\tir\sum_{a=0}^{m-1}(\tir \sn)^a(\bR_{AC4B}+\chi\c \zb) (\tir \sn)^{m-1-a}(\bb^l\sn \log \bb)\big),\, m=1,2,\nn
\end{align}
where the constant coefficients on the right-hand side are symbolic. 

We have the following estimates
\begin{align} \label{10.29.1.22}
\begin{split}
\|\tir (\tir \sn)^m(\bb^l\sn\log \bb)\|_{L^2_u L_\omega^2}&\les \|(\tir \sn)^{m+1}(\bb^l k_{\bN\bN})\|_{L^1_t L^2_u L_\omega^2(\D_{u_1}^t)}+\Delta_0^2+\La_0, m=0,1,2 \\
\|\tir(\tir \sn)^{m-1} (\bb^l\sn\log \bb)\|_{L^4_\omega}&\les \|(\tir \sn)^m(\bb^l k_{\bN\bN})\|_{L_t^1 L^4_\omega(\H_u^t)}+\Delta_0^2+\La_0,\, m=1,2.
\end{split}
\end{align}
(3)\begin{align}\label{8.5.1.22+}\left\{
\begin{array}{lll}
\|\sn_X^n(\sn L \fB, \sn_L \sn\fB)\|_{L^2_\Sigma}
&\les \|\sn_X^n\sn(\sD\varrho, L[L \Phi]), \tir^{-1}\sn_X^{\le n}\sn([L\Phi]+\fB), \bb^{-1}\tir^{-1}\sn_X^{\le n}\sn\bA_b\|_{L^2_\Sigma}\\
&+\l t\r^{-2}(\log \l t\r)^2\Delta_0^\frac{3}{2},
 \qquad n\le 1.\\
\|\bb^{-\f12}\sn_X^\ell\sn_L\zeta\|_{L^2_\Sigma}&\les \|\bb^{-\f12}(\sn_X^\ell \sn k_{\bN\bN}, \sn_X^\ell \sn_L \zb)\|_{L^2_\Sigma}
+\|\tir^{-1}\bb^{-\f12}\sn_X^\ell\sn \log \bb\|_{L^2_\Sigma}\\
&+\l t\r^{-2+\f12\max(\ell-1,0)+2\delta}\Delta_0^2, \qquad\ell\le 2.
\end{array}\right.
\end{align}
There also hold the following comparison formulas 
\begin{align}
X^\ell L\big(\tir^2 (\tr\chi&-\frac{2}{\tir})\big)=O\big(\l t\r^{-\frac{3}{4}+\delta}\Delta_0^\f12\big)\Big(X^{\le \ell}\big(\tir(\tr\chi-\frac{2}{\tir})\big)
+\sn_X^{\le\ell}(\tir\chih)\Delta_0^\f12\Big)\nn\\
&\quad\qquad+X^\ell(\tir^2\widetilde{L(\Xi_4)})+X^{\le \ell}(\tir[L\Phi])+X^\ell(\tir^2\N(\Phi, \bp \Phi))\label{6.8.3.23}, \ell\le 1\\
\sn_X^{\ell}\sn_\Lb\chih&=\sn_X^{\ell}((\bA+\tir^{-1})(\chih+\eh)+\fB\chih+\sn_X^{\ell}\sn\ze+\sn_X^{\ell}(\ze\c\ze)+\sn_X^{\ell}\widehat{\bR_{4A3B}}, \ell\le 1\label{8.3.1.23} \\
\sn_X^{\ell}\sn_S \chih=&\tir\sn_X^{\le\ell}(\fB\chih)+O(1)\chih+\sn_X^{\le \ell} \widehat{\bR_{4A4B}}
+O(\l t\r^{-\frac{11}{4}+\frac{\ell-1}{2}+2\delta}\Delta_0^\frac{3}{2})_{L_u^2 L_\omega^2},\quad \ell=1,2\label{8.25.4.23}
\end{align}
where $\widehat{\bR_{4A3B}}=\bR_{4A3B}-\f12\ga_{AB}\ga^{CD}\bR_{4C3D}$, and $\hat F$ in the above represents the traceless part of the $S$-tangent symmetric 2-tensor $F$.  
\end{proposition}

\subsection{Proof of Proposition \ref{11.4.1.22}} 
 
\noindent{\it Proof of (\ref{8.25.2.23})-(\ref{8.5.1.22}).} 
We first derive from (\ref{1.26.1.23}) that
\begin{align*}
\sn_L\Big((\tir \sn)^l \tir^3\sn\wt{\tr\chi}\Big)&=(\tir\sn)^l\Big(\tir^3\bA \sn\wt{\tr\chi}+\tir^3\big(\sn \chih\c \chih+\sn(\sD \varrho+\hb L \varrho+2\zb^A \sn_A \varrho)\\
&+ [L\Phi] \sn \fB+ \tr\chi \sn[L \Phi]+\sn \N(\Phi, \bp\Phi)\big)\Big)\\
&+l[\sn_L, \tir \sn](\tir^3 \sn \wt{\tr\chi}),
\end{align*}
where $l=0,1$ and the last term can be treated by (\ref{3.21.1.23}). 
Applying transport lemma and using (\ref{3.6.2.21}), we integrate the above identity along the null cone $\H_u^t$, with $u_1\le u\le u_*$, to obtain with $2\le p\le 4$ that 
\begin{align}
&\|(\tir\sn)^l\tir^3\sn\wt{\tr\chi}\|_{L^p_\omega(S_{t,u})}\nn\\
&\les \|(\tir \sn)^l(\tir^3\sn\wt{\tr\chi})\|_{L^p_\omega(S_{0,u})}+\|\tir^3(\tir \sn)^l\sn(\sD\varrho+\hb L\varrho+\zb\c \sn\varrho+\N(\Phi,\bp\Phi))\nn\\
&+\tir^3(\tir \sn)^l([L\Phi] \sn \fB+\tr\chi\sn [L\Phi])\|_{L_t^1 L_\omega^p}+\|\tir^3(\tir\sn)^l\big(\sn\chih\c \chih\big)\|_{L_t^1 L_\omega^p}\nn\\
&+l\int_0^t (\l t'\r^{-\frac{7}{4}+\delta}\Delta_0^\f12+\|(\tir \sn)^l\bA\|_{L_\omega^4})\|\tir^3(\tir \sn)^{\le 1} \sn\wt{\tr\chi}\|_{L_\omega^p}\label{10.25.6.23}
\end{align}
where the last line is obtained by applying  (\ref{error_prod_1}) to $(G, F)=(\tir \sn\bA, \sn\wt{\tr\chi})$. 

It follows from (\ref{6.7.3.23}), (\ref{L2BA2}) and (\ref{3.6.2.21}) that
\begin{align}\label{6.9.5.23}
\|\tir^3 (\tir \sn)^l\sn\N(\Phi, \bp\Phi)\|_{L_t^1 L_u^2 L_\omega^2}&\les \Delta_0^\frac{3}{2}+\|\tir (\tir \sn)^{l+1}([S\Phi]\fB)\|_{L_t^1 L_u^2 L_\omega^2}\\
&\les\Delta_0^\frac{3}{2}+\|\tir (\tir \sn)^{l+1}([S\Phi][\Lb\Phi])\|_{L_t^1 L_u^2 L_\omega^2}\nn
\end{align}
Note that $\hb=\fB-h$, and that $(\tir\sn)^l\sn (\fB L\varrho)$ and $(\tir \sn)^l([L\Phi] \sn \fB)$ in (\ref{10.25.6.23}) can be included in the last term of (\ref{6.9.5.23}). Therefore we drop such error terms in (\ref{10.25.6.23}).
For the other error terms in (\ref{10.25.6.23}), we apply (\ref{3.6.2.21}), (\ref{3.11.3.21}), (\ref{L2conndrv}), (\ref{L2BA2}) and Sobolev embedding on spheres to derive
\begin{align}\label{7.17.1.24}
\|(\tir \sn)^l\big(\sn\chih\c \chih\big), (\tir \sn)^l(\sn(\bA_b[L\Phi]),\sn(\bA_{g,1}^2), \bA_b \sn [L\Phi])\|_{L_u^2 L_\omega^2}\les\l t\r^{-\frac{19}{4}+\frac{l}{2}+2\delta}\Delta_0^\frac{3}{2}.
\end{align}
By taking $L_t^1$ norm of the above line, then substituting the bound $\Delta_0^\frac{3}{2}$ and (\ref{6.9.5.23}) to (\ref{10.25.6.23}), we derive (\ref{7.21.2.22+}) with the help of Minkowski inequality and Proposition \ref{12.21.1.21}. 

Similarly, with $l=0$ in (\ref{10.25.6.23}), again, we include the error $\sn(\fB L\varrho)$ into the last term of the following estimate obtained from (\ref{6.7.3.23})
\begin{align*}
\|\tir^3 \sn\N(\Phi, \bp\Phi)\|_{L_t^1 L_\omega^4}&\les \Delta_0^\frac{3}{2}+\|\tir \Omega([S\Phi]\fB)\|_{L_t^1 L_\omega^4}\\
&\les \Delta_0^\frac{3}{2}+\|\tir \Omega([S\Phi][\Lb\Phi])\|_{L_t^1 L_\omega^4},
\end{align*}
where the last line is derived by using (\ref{3.6.2.21}) and (\ref{3.11.3.21}). For other error terms, noting
\begin{align*}
\|\sn\chih\c \chih, \sn(\bA_b[L\Phi]), \sn(\bA_{g,1}^2), \bA_b \sn [L\Phi]\|_{L_\omega^4}\les\l t\r^{-\frac{9}{2}+2\delta}\Delta_0^\frac{3}{2},
\end{align*}
which are due to (\ref{3.6.2.21}) and (\ref{3.11.3.21}), we thus can obtain (\ref{3.16.14.24}). 

Using (\ref{9.15.3.22}) instead of (\ref{1.26.1.23}), the quadratic error estimates in (\ref{7.17.1.24}), we can similarly obtain (\ref{8.5.1.22}).

 We will rely on (\ref{dchi}) to bound the term $\|(\tir\sn)^{\le 1}(\sn_X \chih)\|_{L^p_\omega}$. We first apply (2) in Lemma \ref{3.31.3.22} to $\sn_X^\ell\chih$ to derive, 
\begin{align}\label{12.12.2.21}
\|\sn(\sn_X^\ell\chih)\|_{L_\omega^p}+\tir^{-1}\|\sn_X^{\ell}\chih\|_{L_\omega^p}&\les \|\sn_X^{\le\ell} \sl{\div}\chih\|_{L_\omega^p}, \, \ell=0,1.
\end{align}
 
In view of (\ref{dchi}), using (\ref{8.3.4.23}), we bound
\begin{align}
 \|\sn_X^{\le\ell} \sl{\div}\chih\|_{L_\omega^p}&\les \| \sn_X^{\le\ell}\sn\tr\chi\|_{L_\omega^p}+\|\sn_X^{\le\ell}(\tir^{-1}\zb, \tensor{\bR}{_{C4}^C_B})\|_{L_\omega^p}+\|\sn_X^{\le\ell}(\zb(\chih, \bA_b))\|_{L_\omega^p}\nn\\
 &\les \|\sn_X^{\le \ell}(\sn\tr\chi, \tir^{-1} \zb, \tensor{\bR}{_{C4}^C_B})\|_{L_\omega^p}+\l t\r^{-\frac{15}{4}+2\delta}\Delta_0^\frac{3}{2}\label{12.12.3.21}.
\end{align}
 Substituting the above estimate to (\ref{12.12.2.21}), we conclude (\ref{8.25.2.23}).
 
 Next we apply (\ref{9.14.10.22}) to $F=\chih$ to obtain
 \begin{align*}
 \|\tir\sn\sn_{X}^2\chih\|_{L^2_u L_\omega^2}+\|\sn_X^2 \chih\|_{L^2_u L_\omega^2}&\les \|\bb^{-\f12}\{\sn_X^2 \sdiv \chih, \l t\r^{-1}\sn_\Omega^{\le 1}\sn_{X_1}\chih, \sn_{X_2}^{\le 1}(\sn \chih, \tir^{-1}\chih)\}\|_{L^2_\Sigma}\nn\\
&+\l t\r^{-\frac{7}{4}+\delta}\Delta_0^\f12 \|\tir (\tir\sn \chih, \chih)\|_{L^2_u L_\omega^4}
 \end{align*}
 where the  terms in the last line is bounded by $\l t\r^{-\frac{5}{2}+2\delta}\Delta_0^\frac{3}{2}$ due to (\ref{3.11.3.21}).
  
 Using (\ref{dchi}) and (\ref{8.28.2.23}), we have
\begin{align*}
\|\sn_X^2\sdiv \chih\|_{L_u^2 L_\omega^2}&\les \|\sn_X^2(\sn\tr\chi+\tir^{-1}\zb+(\bA_b+\bA_{g,2})\zb+\bR_{B4BA})\|_{L_u^2 L_\omega^2}\\ 
&\les \|\sn_X^2 (\sn\tr\chi, \tir^{-1}\zb, \bR_{B4BA})\|_{L^2_\Sigma}+\l t\r^{-\frac{7}{2}+2\delta}\Delta_0^\frac{3}{2}.
\end{align*}

Combining the above two estimates, we conclude (\ref{2.25.2.24}).
 \\

\noindent{\it Proof of (\ref{7.17.2.24})-(\ref{10.29.1.22}).} Using (\ref{lb}), (\ref{1.27.6.24}) and (\ref{cmu2}), it is straightforward to derive  (\ref{7.17.2.24}) and (\ref{12.17.3.23}).
Symbolically, we rewrite (\ref{12.17.3.23}) as 
\begin{align*}
&\sn_L(\tir \sn)^m(\bb^l \sn\log \bb)+\chi (\tir \sn)^m(\bb^l \sn\log \bb)\\
&=-(\tir\sn)^m \sn (\bb^l k_{\bN\bN})+(\bA_b+\bA_{g,2}) (\tir \sn)^m(\bb^l\sn\log \bb)+\bb^l \sn\log \bb ((\tir \sn)^m \chih+(\tir \sn)^m\tr\chi)\\
&+(\tir \sn)\bA\c (\tir \sn)^{m-1}(\bb^l \sn\log \bb)+\tir\sum_{a=0}^{m-1}(\tir \sn)^a(\bR_{AC4B}+\chi\c \zb) (\tir \sn)^{m-1-a}(\bb^l\sn \log \bb).
\end{align*}
The second term on the right-hand side is negligible when applying the transport lemma. Due to (\ref{L2conndrv}) and (\ref{3.6.2.21})
\begin{align*}
\|\bb^l \sn\log \bb  ((\tir \sn)^m \chih+(\tir \sn)^m\tr\chi)\|_{L_u^2 L_\omega^2}\les \l t\r^{2\delta-1-\frac{3}{2}}\log\l t\r\Delta_0^2.
\end{align*}
Using (\ref{3.11.3.21}) and Sobolev embedding, we derive
\begin{align*}
\|(\tir\sn)^{m-1}&(\bb^l \sn\log \bb)\tir \sn \bA\|_{L_u^2 L_\omega^2}\les \|\tir \sn\bA\|_{L_\omega^4}\|(\tir\sn)^{m-1}(\bb^l \sn\log \bb)\|_{L_u^2 L_\omega^4}\\
&\les \l t\r^{-\frac{7}{4}-1+\delta}\Delta_0 \|\tir(\tir \sn)^m(\bb^l \sn\log \bb)\|_{L_u^2 L_\omega^2}^\f12 \|\tir(\tir \sn)^{m-1}(\bb^l \sn\log \bb)\|_{L_u^2 L_\omega^2}^\f12.
\end{align*}
Using (\ref{zeh}), (\ref{1.29.4.22}) and (\ref{10.15.2.22}), we have
\begin{align*}
\sum_{a=0}^{m-1}&\|\tir (\tir \sn)^a (\bR_{AC4B}+\chi\c \zb)(\tir \sn)^{m-1-a}(\bb^l\sn\log \bb)\|_{L_u^2 L_\omega^2}\\
&\les \|\tir(\tir \sn)^{m-1}(\bR_{AC4B}+\chi\c \zb)\|_{L_u^2 L_\omega^2}\|\bb^l \sn\log \bb\|_{L^\infty_\omega}\\
&+\|\tir(\bR_{AC4B}+\chi\c \zb)\|_{L_\omega^4}\|(\tir \sn)^{m-1}(\bb^l \sn\log \bb)\|_{L_u^2 L_\omega^4}\\
&\les \l t\r^{-3+2\delta}\log \l t\r\Delta_0^2 +\l t\r^{-2+\delta}\Delta_0\|(\tir \sn)^{m-1}(\bb^l \sn\log \bb)\|_{L_u^2 L_\omega^4}\\
&\les  \l t\r^{-3+3\delta}\log \l t\r\Delta_0^2.
\end{align*}
where the last term was  further treated by using Sobolev embedding on spheres. 

We thus summarize the above calculations by using transport lemma with $m=\f12$ therein, 
\begin{align*}
\|\tir (\tir \sn)^m(\bb^l \sn\log \bb)\|_{L_u^2 L_\omega^2}&\les\|\tir (\tir \sn)^m(\bb^l \sn\log \bb)(0)\|_{L_u^2 L_\omega^2}+\|(\tir\sn)^{m+1} (\bb^l k_{\bN\bN})\|_{L_t^1 L_u^2 L_\omega^2}\\
&+\int_0^t \l t\r^{-\frac{7}{4}+\delta}\Delta_0\|\tir(\tir \sn)^{m}(\bb^l \sn\log \bb)\|_{L_u^2 L_\omega^2}+\Delta_0^2. 
\end{align*}
Using (\ref{7.17.2.24}), and the transport lemma, we have the lowest order estimates of (\ref{10.29.1.22}). Using the above inequality and transport lemma we obtain
\begin{align*}
\|\tir (\tir\sn)^m&(\bb^l \sn\log \bb)\|_{L_u^2 L_\omega^2}\les\|\tir (\tir \sn)^m(\bb^l \sn\log \bb)(0)\|_{L_u^2 L_\omega^2}+\|(\tir\sn)^{m+1} (\bb^l k_{\bN\bN})\|_{L_t^1 L_u^2 L_\omega^2}+\Delta_0^2,
\end{align*}
as stated in (\ref{10.29.1.22}) in view of Proposition \ref{12.21.1.21}. 

Next we consider the second estimate in (\ref{10.29.1.22}) with $m=2$. With $m=1$ in (\ref{12.17.3.23}), we rewrite the right-hand side as
\begin{align*}
(\chih, \bA_b)(\tir \sn)(\bb^l \sn\log \bb)+\tir \bb^l\sn\log \bb (\sn\chih, \sn\tr\chi)+\tir (\bR_{AC4B}+\chi\c \zb)\bb^l \sn\log \bb. 
\end{align*}
Again the first term is negligible when applying the transport lemma. Using (\ref{zeh}), (\ref{3.11.3.21}) and (\ref{1.29.4.22}), we bound
\begin{align*}
&\|\tir \bb^l \sn\log \bb(\sn\chih, \sn\tr\chi)\|_{L_\omega^4}\les \l t\r^{-\frac{7}{4}-1+2\delta}\Delta_0^2\log \l t\r;\\
&\|\tir (\bR_{AC4B}+\chi\c \zb)\bb^l \sn\log \bb\|_{L_\omega^4}\les \l t\r^{-2-1+2\delta}\log \l t\r\Delta_0^2. 
\end{align*}
Thus, applying the transport lemma to (\ref{12.17.3.23}) in this case, with the help of Proposition \ref{12.21.1.21}, leads to 
\begin{align*}
\|\tir (\tir\sn)&(\bb^l \sn\log \bb)\|_{L_\omega^4}\les \La_0+\Delta_0^2+\|(\tir\sn)^2 (\bb^l k_{\bN\bN})\|_{L_t^1 L_\omega^4},
\end{align*}
as desired. The proof for (2) is complete.

\noindent{\it Proof of (\ref{8.5.1.22+})-(\ref{8.25.4.23}).}
Next, using (\ref{3.20.1.22}) and (\ref{cmu_2}), we decompose
\begin{align*}
\sn_X^n(\sn L (k_{\bN\bN}), (\sn_L+\chi) \sn k_{\bN\bN})&=\sn_X^n\{\sn(\widetilde{L(\Xi_4)})+\sn L[L\Phi]+\sn\big((h-k_{\bN\bN}) \Xi_4\big)\}\\
&=\sn_X^n\{\sn(\sD \varrho-\hb L \varrho+2\zb^A \sn_A \varrho-\Box_\bg \varrho+(L+(h-k_{\bN\bN}))L \varrho)\\
&+\sn(\fB(\bA_b+\fB))+\tir^{-1}\sn\fB +\sn L[L\Phi]\}.
\end{align*}
By using (\ref{6.6.1.23}) with $n=2, l=0$, (\ref{L2BA2}), (\ref{LbBA2}), (\ref{L2conndrv}), (\ref{3.6.2.21}), (\ref{3.11.3.21}) and $\bb\fB=O(\l t\r^{-1})$, the right-hand side in the above is bounded by 
\begin{align*}
 &\|\sn_X^n\sn(\sD\varrho, L[L \Phi]), \tir^{-1}\sn_X^{\le n}\sn([L\Phi]+\fB),\sn_X^n(\fB\sn  \fB),  \bb^{-1}\tir^{-1}\sn_X^{\le n}\sn\bA_b\|_{L^2_\Sigma}\\
 &+\l t\r^{-3+2\delta}\log \l t\r^\f12\Delta_0^\frac{3}{2}.
\end{align*}
In view of (\ref{8.23.2.23}), we obtain $\|\Omega \fB \c \sn\fB \bb^2\|_{L_\omega^2}\les \l t\r^{-3}(\log \l t\r)^2 \Delta_0^2$ and hence
\begin{equation*}
\sn_X^n(\fB \sn \fB)=\l t\r^{-1}\sn_X^{\le n}\sn\fB+O(\Delta_0^2\l t\r^{-3}(\log \l t\r)^2 \bb^{-2})_{L_\omega^2}
\end{equation*}
with the help of $\bb \fB=O(\l t\r^{-1})$. Summarizing the above estimates, we obtain the bound of $\|\sn_X^n \sn L\fB\|_{L^2_\Sigma}$ in (\ref{8.5.1.22+}).

Using (\ref{8.24.4.23}) and (\ref{3.11.3.21}), we derive
\begin{align*}
\sn_X (\chi\c \sn k_{\bN\bN})&=O(\l t\r^{-1})(\sn_X+\vs(X))\sn k_{\bN\bN}+(1-\vs(X))O(\l t\r^{-\frac{7}{4}+\delta}\Delta_0)_{L_\omega^4}\sn k_{\bN\bN}\\
&=O(\l t\r^{-1})(\sn_X+\vs(X))\sn k_{\bN\bN}+O(\l t\r^{-\frac{15}{4}+2\delta}\Delta_0^2)_{L_u^2 L_\omega^2}.
\end{align*}
Hence, in view of $\fB=k_{\bN\bN}+[L\Phi]$,  we obtain the  estimate of $\|\sn_X^n \sn_L \sn\fB\|_{L^2_\Sigma}$  in (\ref{8.5.1.22+}).

Next we consider the estimate for $\sn_X^\ell\sn_L\zeta$ in (\ref{8.5.1.22+}). 
Symbolically, due to (\ref{lb}) and (\ref{cmu_2}),  
 \begin{equation}\label{12.8.5.21}
-\sn k_{\bN\bN}=\sn L\log \bb=(\sn_L+\chi) \sn \log \bb.
\end{equation}
Due to $\zeta+\zb=\sn\log \bb$ and (\ref{12.8.5.21}), we obtain
\begin{equation}\label{6.23.6.24}
\sn_X^\ell\sn_L(\zeta+\zb)=\sn_X^\ell \sn_L \sn\log \bb=-\sn_X^\ell (\sn k_{\bN\bN}+\chi\c \sn \log \bb).
\end{equation}
In view of the above formula, using (\ref{L2conndrv}), (\ref{3.6.2.21}), (\ref{3.11.3.21}) and (\ref{3.16.1.22}) we infer
\begin{align*}
\|\bb^{-\f12}\sn_X^\ell\sn_L\zeta\|_{L^2_\Sigma}&\les \|\bb^{-\f12}(\sn_X^\ell \sn k_{\bN\bN}, \sn_X^\ell \sn_L \zb)\|_{L^2_\Sigma}+\|\bb^{-\f12}\sn_X^\ell(\chi\c\sn \log \bb)\|_{L^2_\Sigma}\\
&\les \|\bb^{-\f12}(\sn_X^\ell \sn k_{\bN\bN}, \sn_X^\ell \sn_L \zb)\|_{L^2_\Sigma}+\|\bb^{-\f12}(\tir^{-1}\sn_X^{\le \ell}\sn\log \bb, \sn_X^\ell\chi\c \sn\log \bb)\|_{L^2_\Sigma}\\
&+\l t\r^{-2+\frac{\max(\ell-1,0)}{4}+2\delta}\Delta_0^2.
\end{align*}
Using (\ref{L2conndrv}) and (\ref{3.6.2.21}), we further obtain
\begin{align*}
\sn_X^\ell\chi\c \sn\log \bb&=\sn\log \bb( O(\l t\r^{-1})+O(\l t\r^{-2+\frac{\max(\ell-1,0)}{2}+\delta})_{L_u^2 L_\omega^2})\\
&=O(\l t\r^{-1})\sn \log \bb+O(\l t\r^{-3+\f12\max(\ell-1,0)+2\delta}\Delta_0^\frac{3}{2})_{L^2_u L_\omega^2}. 
\end{align*}
Hence we obtain the second estimate in (\ref{8.5.1.22+}).

Using (\ref{6.3.1.23}), for $\ell\le 1$, we derive
\begin{align*}
X^\ell L(\tir^2 (\tr\chi-\frac{2}{\tir}))&=X^{\ell}(\tir^2(\tr\chi-\frac{2}{\tir})(\bA_b+[L\Phi]))
+X^\ell (\tir^2(\chih\c \chih))+X^\ell(\tir^2\widetilde{L(\Xi_4)})\\
&+X^\ell(\tir[L\Phi])+X^\ell(\tir^2\N(\Phi, \bp \Phi)).
\end{align*}
This gives (\ref{6.8.3.23}) in view of (\ref{3.6.2.21}).

 (\ref{8.3.1.23}) follows by using (\ref{3chi}). 
It follows by using (\ref{s2}) that
\begin{equation}
\sn_X^{\ell}\sn_L(\tir\chih)=\sn_X^{\ell }(\tir (\bA_b+\fB)\chih)+\sn_X^{\ell}(\tir\widehat{\bR_{4A4B}})\label{6.8.2.23}
\end{equation}
 and thus (\ref{8.25.4.23}) can be obtained by using (\ref{6.8.2.23}), (\ref{8.3.4.23}) and (\ref{8.28.2.23}).

 The proof of Proposition \ref{11.4.1.22} is complete. $\hfill\square$

\subsection{Estimates of $\eta$}\label{12.6.1.23}
In this subsection, we will use the equations (\ref{6.1.1.23}) and (\ref{12.16.1.23}) to control $\eh$ and  $\sdiv(\eta(\Omega))$. In particular, we can prove that $\sn_\Lb \eh$ verifies better estimates than $\sn_\Lb \bA_{g,1}$. We also gain better estimates for $\sdiv(\eta(\Omega))$ and its derivatives, which will be used in energy estimates. 
   
\begin{proposition}\label{8.12.1.23}
There holds with $X\in\{\Omega, S\}$ that
\begin{align}
\|(\tir\sn)^{\le 1}\sn_X^\ell\eh\|_{L^p_\omega}\les \|\sn_X^\ell(\Omega[L\Phi]), \sn_X^{\le \ell}\ep\|_{L_\omega^p}+\l t\r^{-\frac{11}{4}+2\delta}\Delta_0^\frac{3}{2}, \label{8.2.1.23}
\end{align}
where $X=\Omega, S$, $\ell=0,1$ and $2\le p\le 4$.

Moreover, there holds the identity 
\begin{align}\label{6.26.2.23}
-\sdiv(\eta({}\rp{a}\Omega))&=\tr\thetac \ep({}\rp{a}\Omega)-{}\rp{a}\Omega[L\Phi]+{}\rp{a}\la \thetac\c \eta+\snc_A\log c \Omega v^A.
\end{align}
We have the following estimates for $\sdiv(\eta({}\rp{a}\Omega))$, 
\begin{align}\label{8.2.2.23}
\begin{split}
&X^l\sdiv\big(\eta(\Omega)\big)=\sn_X^{\le l}\zb+X^l\Omega[L\Phi]+O(\l t\r^{-3+\frac{l}{4}+2\delta}\Delta_0^\frac{3}{2})_{L_\omega^4},\,  l=0,1\\
&X^l \sdiv\big(\eta(\Omega)\big)=\tir^{-1}\sn_X^{\le l}[\Omega v]+X^l(\Omega[L\Phi])+O(\l t\r^{-\frac{11-(l-1)}{4}+2\delta}\Delta_0^\frac{3}{2})_{L_u^2 L_\omega^2}, l=1,2
\end{split}
\end{align}
\end{proposition}
\begin{proof}
Using Lemma \ref{3.31.3.22}, we derive with $X=\Omega, S$
\begin{align*}
\|\tir^{-1}(\tir \sn)^{\le 1}\sn_X^\ell \eh\|_{L_\omega^4}\les\|X^{\le \ell}\D_2 \eh\|_{L^4_\omega}.
\end{align*}
(\ref{8.2.1.23}) follows by substituting (\ref{6.1.1.23}) into the right-hand side and using (\ref{8.24.4.23}) and (\ref{8.3.4.23}). 

Next we derive by using (\ref{6.1.1.23}) and (\ref{2.10.1.22})
\begin{align*} 
-\snc_A({}\rp{a}\Omega v^A)&=-\snc\Omega^A\c\p_A v^\|+\f12{}\rp{a}\Omega\tr\eta+(\stc\sdiv\eh)({}\rp{a}\Omega)\\
&=-{}\rp{a}\la \thetac\c \eta+{}\rp{a}\Omega[L\Phi]-\tr\thetac\ep({}\rp{a}\Omega).
\end{align*}
This gives (\ref{6.26.2.23}) in view of (\ref{1.27.1.22}). 

The first estimate of (\ref{8.2.2.23}) is obtained by using (\ref{3.6.2.21}), (\ref{3.11.3.21}), (\ref{3.16.1.22}), (\ref{8.3.4.23}) and Proposition \ref{10.16.1.22}. 
Using Proposition \ref{10.16.1.22}, Lemma \ref{5.13.11.21} (1) and (2), (\ref{L2conndrv}) and (\ref{L2BA2}), we have
\begin{align*}
\|X^l(\bAn \la \thetac)\|_{L^2_u L_\omega^2}\les \l t\r^{-3+2\delta}\Delta_0^\frac{3}{2}.
\end{align*}
The second estimate follows by using  (\ref{4.22.4.22}), (\ref{6.26.2.23}), (\ref{8.3.4.23}), (\ref{8.28.2.23}) and the above estimate.
\end{proof}

\begin{proposition}\label{10.4.4.23} 
Let $X_i\in \{\Omega, S\}$. There hold 
\begin{align}
\|\sn_X^n (\eh({}\rp{a}\Omega))\|_{L^2_\Sigma} &\les  \|X^{\le n} \Omega v\|_{L^2_\Sigma}+\|\sn_X^{n}(\tr\eta\c \Omega)\|_{L^2_\Sigma}\nn\\
&+\l t\r^{-\frac{3}{4}+\frac{\max(n-2,0)}{4}+2\delta}\log \l t\r^\f12\Delta_0^\frac{3}{2}, n\le 3,
\label{8.9.5.23}\\
\|\sn_X^n \sn_\bN \eh\|_{L^2_\Sigma}&\les \|\sn_X^n (\snc \hot \ep, \ud\bA\c \ep, \fB \bA_g)\|_{L^2_\Sigma}+\|\tir^{-1}\sn_X^{\le n}\eh\|_{L^2_\Sigma}\nn\\
&+\l t\r^{-\frac{5}{2}+2\delta}\log \l t\r^\f12\Delta_0^\frac{3}{2}, n\le 2\label{1.30.2.24}\\
\|\sn_X^n \sn_\bN \eh\|_{L^4_\omega}&\les \|\sn_X^n (\snc \hot \ep, \ud\bA\c \ep, \fB  \bA_g), \tir^{-1}\sn_X^{\le n}\eh\|_{L_\omega^4}+\l t\r^{-\frac{15}{4}+2\delta}\Delta_0^\frac{3}{2}, n\le 1.\label{2.1.2.24}
\end{align}  
\end{proposition}

\begin{proof}
(\ref{8.9.5.23}) is a consequence of (\ref{7.25.5.23}).
(\ref{1.30.2.24}) follows by using (\ref{12.16.1.23}), (\ref{8.3.4.23}) and (\ref{8.28.2.23}). (\ref{2.1.2.24}) can be derived similarly by using (\ref{8.3.4.23}).

\end{proof}

\subsection{Summary}

For ease of future reference, we summarize Assumption \ref{5.13.11.21+}, the decay estimates obtained in Lemma \ref{5.13.11.21}, and their consequences, which can be derived by applying Lemma \ref{6.30.4.23}, Lemma \ref{3.17.2.22}, Lemma \ref{6.9.3.23}, Proposition \ref{11.4.1.22}, and (\ref{7.22.1.22}), into one comprehensive result. This summary will be frequently used and will be complemented by the more precise estimates provided in Proposition \ref{11.4.1.22}, Proposition \ref{8.12.1.23}, and Proposition \ref{10.4.4.23}.
\begin{proposition}\label{7.15.5.22}
Under the assumptions of (\ref{3.12.1.21})-(\ref{6.5.1.21}), with $X_i\in \{\Omega, S\}$, there hold
\begin{align*}
&\Delta_0^\f12\l t\r^\frac{l}{4}|\sn_X^l (L\Phi, \sn\Phi)|+|\bA_{g,1}|+\l t\r^{-\frac{1+3l}{4}}( |\sn_X^{l}\bA_{g,2}|+\Delta_0^\f12 |\sn_X^l\bA_b|)\les \Delta_0\l t\r^{-2+\delta}, l\le 1\\
 &\tir(|\sn_S^{\le1}\ze|+|\tir (\sn[\Lb\Phi], \sn_\Lb \bA_{g,1})|)+ \tir^2|\sn \tr\chi|\les \Delta_0  \l t\r^\delta\\
&\|\l t\r^{-\frac{l}{4}}(\sn_X^{l+\le 1} \bA_{g,1}, \sn_X^l\sn_\Omega[L\Phi]), \l t\r^{-\frac{1}{4}}(\Delta_0^{\f12\vs(X)}\sn_X \bA_b, \sn_X\bA_{g,2}),\bA_{g,2}, \Delta_0^\f12\bA_b\|_{L_\omega^4}\\
&\qquad\qquad\qquad\qquad\qquad\qquad\quad\les \l t\r^{-2+\delta}\Delta_0, l=0,1\\
&\|\sn_X^2 \bA_{g,2}, \Delta_0^{\f12\vs^-(X^2)} \sn_X^2 \bA_b\|_{L_\omega^4}\les \l t\r^{\delta-1-\frac{3}{4}\vs^+(X^2)}\Delta_0\\
&\|\l t\r^{-\delta}\sn_\Lb \bA_{g,2}, \sn_X^{\le 1}\sn_\Lb \bA_{g,1}, \l t\r^{-1}(1-\vs^-(X^l))X^l\fB\|_{L_\omega^4}\les \l t\r^{-2+\delta}\Delta_0, l\le 2\\
& \|\sn_X^{\le 1} \ze, \vs^+(X^2)\sn_X^2\ze\|_{L_\omega^4}\les\l t\r^{-1+\delta}\Delta_0\\ 
&\|X^{2+l} \Phi, \Sc(X^{2+l}\Phi), \Ac(X^{2+l}\Phi)\|_{L^4_\omega}\les \l t\r^{-1+\frac{l}{4}+\delta}\Delta_0^\f12, l=0,1\\
 \displaybreak[0]
 &\|\Sc(X^l\Phi,\tir X^{l-1}L\Phi,\tir\sn_X^{l-1}\sn\Phi), \Ac(X^l\Phi, \tir X^{l-1}L\Phi,\tir \sn_X^{l-1}\sn\Phi)\|_{L^2_\Sigma}\\
 &\qquad\qquad\qquad\qquad\les \Delta_0\l t\r^{\f12\max(l-3,0)+\delta}, 1\le l\le 4\\
&\|\bb^{-\f12}\Big(\bA_{g,2}, \bA_b, \l t\r^{(-\f12+\f12\vs^+(X^{l+1})) l}\sn_X^{l+1}(\bA_{g,2}, \bA_b)\Big), \l t\r^{-\f12 l}\sn_X^{l+\le 2}\bAn\|_{L^2_\Sigma}\\
&+\l t\r^{-1}\sum_{\vs^-(X^i)=0}\| \sn_X^i \fB\|_{L^2_\Sigma}+\l t\r^{-l\delta}\|\bb^{-\f12 l}\sn_X^{\le 1}\sn_\Lb\bA_{g, 1+l}\|_{L^2_\Sigma}\les \Delta_0\l t\r^{-1+\delta} \\
&\qquad\qquad\qquad\qquad l=0,1, i=1,2,3\\
&\|\sn_X^{\le 2} \ze, \tir\sn_X^{\le 2}\sn_L \ze\|_{L_u^2 L_\omega^2}+\|\tir(\sn_\Omega^2 \sn\bA_b, \sn\sn_\Omega^2\bA_{g,2})\|_{L_u^2 L_\omega^2}\les \l t\r^{-1+\delta}\Delta_0\\
&\|\bb^{-\frac{1}{2}}\tir\Omega^{1+\le 1}(\bb^{-1}\sn(\tr\chi+\Xi_4))\|_{L^2_\Sigma}+\|\tir^2\sn(\bb^{-1}(\tr\chi+\Xi_4))\|_{L_\omega^4}\les \l t\r^\delta\Delta_0.
\end{align*}
\end{proposition}
We remark that the last line is a consequence of (\ref{L2BA2}), (\ref{L2conndrv}), (\ref{LbBA2}) and Lemma \ref{5.13.11.21} (1), (5). We omit the detailed checking for the remaining estimates for simplicity. 

\section{Error estimates}
\subsection{Preliminaries}\label{error_manu} 
\begin{proposition}[Decomposition of commutators]
Let
\begin{equation}\label{8.28.11.23}
\fm{X}:=-{}\rp{X}\pi_{L \Lb}, \quad {}\rp{X}\pih^{\mu\nu}:={}\rp{X}\pi^{\mu\nu}-\f12 \bg^{\mu\nu}\fm{X}.
\end{equation}
 It holds
\begin{equation}\label{5.02.3.21_1}
 \Er_1(f,X):=[\Box_\bg, X]f=\f12 \fm{X}\Box_\bg f+\sP[X, f]
 \end{equation}
where
\begin{align*}
&\sP[X, f]:={}\rp{X}\pih^{\a\b} \bd^2_{\a\b} f+(\J[X]^\nu -\f12 \bd^\nu{}\rp{X}\ss)\p_\nu f,\\
&{}\rp{X}\ss=\ga^{AB}{}\rp{X}\pi_{AB}, \quad\J[X]^\nu:=\bd_\mu \pih^{\mu\nu}.
\end{align*}

With $X_i\in \{\Omega, S\}$, $i=1,2,3$ and denote the higher order commutators for the function $X_n \cdots X_1 f$
\begin{equation}\label{10.30.2.21}
\Er_n(f, X_n\cdots X_1):=\Box_\bg(X_n\cdots X_1 f)-X_n\cdots X_1\Box_\bg f.
\end{equation}
We have
\begin{equation*}
\Er_2(f, X_2X_1)=\Er_1(X_1f, X_2)+X_2(\Er_1(f, X_1)).
\end{equation*}
With $f=\Phi$, denoting the terms in the above identity that have a factor of $\Box_\bg X^{\le n-1}\Phi$ by $I[X^2\Phi]$, 
\begin{equation}\label{3.29.2.23}
	I[X^2\Phi]:= X_2(\f12\fm{X_1}\Box_\bg \Phi)+\f12\fm{X_2}\Box_\bg X_1 \Phi,
\end{equation}
 we write
\begin{align}\label{5.02.4.21}
	\Er_2(\Phi, X_2 X_1)&= X_2(\sP[X_1, \Phi])+\sP[X_2, X_1 \Phi]+I[X^2\Phi].
\end{align}
With $n=3$ and $f=\Phi$ in (\ref{10.30.2.21}), using the above convention, we write
\begin{align}\label{3.29.2.23'}
	\Er_3(\Phi, X_3 X_2 X_1)&=X_3 X_2(\sP[X_1, \Phi])+X_3(\sP[X_2, X_1 \Phi])+\sP[X_3, X_2 X_1\Phi]+I[X^3\Phi],
\end{align}
where
\begin{equation*}
	2I[X^3\Phi]:=X_3 X_2(\fm{X_1}\Box_\bg \Phi)+X_3(\fm{X_2}\Box_\bg X_1 \Phi)+\fm{X_3}\Box_\bg X_2 X_1 \Phi.
\end{equation*}
\end{proposition}
\begin{proof}
 Note it is direct to have \begin{footnote}{ See the first line in \cite[Page 169]{shock_demetrios}.}\end{footnote}
 \begin{align}
 [\Box_\bg, X]f&=\bd^\nu({}\rp{X}\pi_{\mu\nu}\p^\mu f)-\f12\bd_\mu(\bg^{\a\b}{}\rp{X}\pi_{\a\b})\p^\mu f\label{10.12.1.23}\\
 &=\bd^\nu\big(({}\rp{X}\pih_{\mu\nu}+\f12 \bg_{\mu\nu}\fm{X})\p^\mu f\big)-\f12\bd_\mu(\bg^{\a\b}{}\rp{X}\pi_{\a\b} )\p^\mu f\nn\\
 &=\bd^\nu\big({}\rp{X}\pih_{\mu\nu}\p^\mu f)+\f12 \fm{X}\Box_\bg f+\f12 \bd_\mu(\fm{X}-\bg^{\a\b}{}\rp{X}\pi_{\a\b}) \p^\mu f.\nn
 \end{align}
 The consequence drops out due to $\fm{X}-\bg^{\a\b}{}\rp{X}\pi_{\a\b}=-\ga^{AB}{}\rp{X}\pi_{AB}=-{}\rp{X}\ss$. The remaining parts are direct checking.
 \end{proof}
  Remark that the decomposition of the commutator is different from the classical treatment. The classical treatment is to decompose ${}\rp{X}\pi$ into the trace and traceless part.
   We instead decompose it in (\ref{8.28.11.23}). In this way, the non-small quantity $\fm{S}$ is paired with $\Box_\bg  f$, with the latter satisfying the desired smallness due to geometric null condition modulo controllable errors in our setting. Thanks to this treatment, for scalars $f$, the large transversal terms in the commutator are combined into a good term. Nevertheless applying the commutation formula to velocity, for $X=\Omega$, we will see more large terms requiring further treatment. 

\subsubsection{Deformation tensors ${}\rp{X}\pih$ and the current $\bJ[X]$}
We have given the components of $\pio$ in Proposition \ref{3.22.6.21}. Now using Proposition \ref{6.7con}, we can directly obtain the basic properties of $\piS$.
\begin{proposition}\label{8.18.3.21}
The components of $\piS$ are listed below
\begin{equation}\label{6.23.4.21}
\begin{split}
&\piS_{AB}=2\tir\chi_{AB} \quad\quad\,\,\piS_{\Lb A}=2\tir(\zeta_A-\zb_A),\\
&\piS_{LA}=0=\piS_{LL} \quad \piS_{L\Lb}=-\fm{S}=-2+2\tir k_{\bN\bN},\\
&\f12\piS_{\Lb \Lb}=-2 \tir k_{\bN\bN}-2\Lb \tir=-\tir(2 k_{\bN\bN}+2\mho+\tr\chib).
\end{split}
\end{equation}
With ${}\rp{X}\ss=\ga^{AB}{}\rp{X}\pi_{AB}$, there holds
 \begin{equation*}
 {}\rp{S}\ss=2 \tir \tr\chi.
 \end{equation*}
 It holds for $X={}\rp{a}\Omega$ and $S$ that
$${}\rp{X}\pih_{L\Lb}=0={}\rp{X}\pih_{LL}.$$
 Other components of ${}\rp{X}\pih$ relative to the null tetrad are the same as components as ${}\rp{X}\pi,$ except that
\begin{equation}\label{7.15.1.22}
\begin{split}
{}\rp{a}\pih_{AB}&=-2c^{-1}\la\hat\theta_{AB}+\f12 {}\rp{a}\ss \ga_{AB}-\f12\ga^{AB}\fm{{}\rp{a}\Omega}\\
{}\rp{S}\pih_{AB}&=2\tir \chih_{AB}+\f12 {}\rp{S}\ss\ga_{AB}-\f12\ga^{AB}\fm{S}.
\end{split}
\end{equation}
 \end{proposition}
We omit the details of calculation, since they are rather straight forward.
Next we give the following result on the terms in commutators.
 \begin{proposition}\label{11.12.2.22}
Let 
$$\ckk\J[X]_A=\J[X]_A,\, \ckk\J[X]_L=\J[X]_L+{}\rp{X}\sl{\pih}\c \chi, \ckk\J[X]_\Lb=\J[X]_\Lb+{}\rp{X}\sl{\pih}\c \chib,$$ 
where ${}\rp{X}\sl{\pih}$ denotes the fully angular components of ${}\rp{X}\pih$. There hold
\begin{equation}\label{6.28.3.21}
\ckk\J[S]_L=0,\quad \ckk\J[\Omega]_L=\sn_A {}\rp{a}\pih_{AL}+2\zeta_A\c {}\rp{a}\pih_{AL}.
\end{equation}

We denote $\bJ[X]_\mu=\ckk \J[X]_\mu-\f12\p_\mu {}\rp{X}\ss$, and ${}\rp{a}\bJ:=\bJ[{}\rp{a}\Omega]$ for short.
 It holds that
\begin{equation}\label{5.18.3.21}
\begin{split}
\sP[X, f]&=\bJ[X]^\mu\p_\mu f+{}\rp{X}\pih^{AB}\sn_A \sn_B f\\
  &+\frac{1}{4}{}\rp{X}\pih_{\Lb\Lb}\bd^2_{LL}f-\f12({}\rp{X}\pih_{L A}\bd^2_{\Lb A}f+{}\rp{X}\pih_{\Lb A}\bd^2_{LA} f).
  \end{split}
\end{equation}
where, importantly, ${}\rp{S}\pih_{LA}=0$. 
\end{proposition}
We remark that it is important the term $\bJ[X]_L$ which is multiplied to $\Lb f$ is either $0$ or small. 
\begin{proof}
 With 
$I[X, \psi]^\mu:={}\rp{X}\pih^{\mu\nu}\p_\nu \psi,$  we derive
\begin{align}\label{5.02.3.21}
\bd_\a I[X, f]^\a&={}\rp{X}\pih^{\a\b} \bd^2_{\a\b} f+\J[X]^\la \p_\la f; \sP[X, f]=\bd_\a I[X, f]^\a-\f12\bd^\nu{}\rp{X}\ss\p_\nu f. 
\end{align}
To see (\ref{5.18.3.21}), for a symmetric two tensor-field $F$, we directly derive by using Proposition \ref{6.7con}
\begin{equation}\label{5.13.1.21}
\begin{split}
\bd_A F_{BL}&=\sn_A F_{BL}-F(\bd_A \Pi_{\mu'}^\mu e_B^{\mu'}, L)-F(\bd_A L, e_B)\\
&=\sn_A F_{BL}-\f12 F_{\mu\nu}\bd_A(L_{\mu'} \Lb^\mu+\Lb_{\mu'}L^\mu) L^\nu e_B^{\mu'}-F(\bd_AL, e_B)\\
&=\sn_A F_{BL}-\f12 F_{L\Lb}\chi_{AB}-\f12 F_{LL}\chib_{AB}-F(e_B, \chi_{AC}e_C-k_{A\bN}L).
\end{split}
\end{equation}
Thus
\begin{equation*}
\bd^A F_{AL}=\sn_A F_{AL}-\f12 F_{L\Lb}\tr\chi-\f12 F_{LL} \tr\chib-F_{AC}\chi_{AC}+k_{A\bN} F_{AL}.
\end{equation*}
We can similarly obtain
\begin{align*}
\bd^L F_{LL}&=-\f12 \bd_\Lb F_{LL}=-\f12\Lb (F_{LL})+F_{\bd_\Lb L, L}\\
\displaybreak[0]
&=-\f12 \Lb (F_{LL})+F(2\zeta_A e_A+k_{\bN\bN}L, L),\\
\bd^\Lb F_{\Lb L}&=-\f12 \bd_L F_{\Lb L}=-\f12\left(L(F_{\Lb L})-F(\bd_L \Lb, L)-F(\Lb, \bd_L L)\right)\\
&=-\f12 (L( F_{\Lb L})-F(2\zb_A e_A+k_{\bN\bN}\Lb, L)+F_{\Lb L}k_{\bN\bN}).
\end{align*}
We derive by using   ${}\rp{X}\pih_{L\Lb}=0={}\rp{X}\pih_{LL}$ in Proposition \ref{8.18.3.21} that
\begin{equation}\label{6.28.2.21}
\begin{split}
&\bd^\Lb {}\rp{X}\pih_{\Lb L}=\zb_A {}\rp{X}\pih_{AL}, \quad \bd^L {}\rp{X}\pih_{LL}=2\zeta_A {}\rp{X}\pih_{AL},\\
&\bd^A {}\rp{X}\pih_{AL}=\sn_A {}\rp{X}\pih_{AL}-{}\rp{X}\pih_{AC}\chi_{AC}+k_{A\bN} {}\rp{X}\pih_{A L}.
\end{split}
\end{equation}
We further derive by ${}\rp{S}\pih_{AL}=0$ in Proposition \ref{8.18.3.21} that
\begin{align*}
\bd^A {}\rp{S}\pih_{AL}&=-{}\rp{S}\sl{\pih}\c \chi,\quad  \bd^L {}\rp{S}\pih_{LL}=0, \quad \bd^\Lb {}\rp{S}\pih_{\Lb L}=0,
\end{align*}
which implies the first identity in (\ref{6.28.3.21}). We can obtain the second identity in  (\ref{6.28.3.21}) by adding up the formulas in (\ref{6.28.2.21}).

 We now consider (\ref{5.18.3.21}).
Combining (\ref{5.02.3.21}) with
$$
\bd^2_{AB}f=\sn_A\sn_B f-\f12 \chi_{AB} \Lb f-\f12 \chi_{AB}L f,
$$
 since ${}\rp{X}\pih_{L\Lb}=0={}\rp{X}\pih_{LL}$ in Proposition \ref{8.18.3.21} we derive
\begin{align*}
\bd_\a I[X, f]^\a&=\J[X]_A\c \sn_A f+{}\rp{X}\pih^{AB}\sn_A\sn_B f+({}\rp{X}\pih\c \bd^2 f)_{s\neq 0, -2}\\
&-\f12\{(\J[X]_L+\chib\c {}\rp{X}\sl{\pih})\Lb f+(\J[X]_\Lb+\chi\c {}\rp{X}\sl{\pih})L f\},
\end{align*}
where we denote the signature of  $\bd^2 f(X_1, X_2)$ relative to the null tetrad  by $s:=\vs(X_1)+\vs(X_2)$.

Further decomposing $({}\rp{X}\pih\c \bd^2 f)_{s\neq 0, -2}$ gives (\ref{5.18.3.21}).
\end{proof}

Next we derive other components of $\bJ[X]$ or $\ckk\J[X]$.
\begin{lemma}\label{3.23.2.23}
\begin{align}
{}\rp{a}\bJ_B&={}\rp{a}\wt{\eth_B}-\frac{3}{4}(\tr\chi{}\rp{a}\pih_{\Lb B}+\tr\chib {}\rp{a}\pih_{LB})-\f12 \chih_{AB}{}\rp{a}\pih_{A\Lb}-\f12 \chibh_{AB}{}\rp{a}\pih_{AL}\nn\\
&+(\zeta_A+\zb_A){}\rp{a}\pih_{AB}+k_{\bN\bN}{}\rp{a}\pih_{\bT B}\label{6.28.5.21}\\
{}\rp{a}\ckk\J_\Lb&=\sn^A {}\rp{a}\pih_{A\Lb}-\f12 L({}\rp{a}\pih_{\Lb \Lb})-(\f12 \tr\chi-k_{\bN\bN}) {}\rp{a}\pih_{\Lb \Lb}\nn\\
&+(\zeta_A+3\zb_A){}\rp{a}\pih_{A\Lb}-(\zeta-k_{\bN A}){}\rp{a}\pih_{LA}\label{6.28.1.21}
\end{align}
where \begin{footnote}{$\J[{}\rp{a}\Omega]$ is denoted by ${}\rp{a}\J$ for short.}\end{footnote}
\begin{align}
{}\rp{a}\wt{\eth_B}&:=\rp{a}\eth_B-\f12(\sn_\Lb {}\rp{a}\pih_{LB}+\sn_L {}\rp{a}\pih_{\Lb B})\nn\\
{}\rp{a}\eth_B&:=\sn_A{}\rp{a}\pih_{AB}-\f12\sn_B{}\rp{a}\ss=-2\sn_A(c^{-1}\la\hat\theta_{AB})-\f12\sn_B\fm{{}\rp{a}\Omega}\label{7.15.2.22}\\
\fm{\Omega}&=2{}\rp{a}\Omega\log \bb\label{8.29.1.23}\\
{}\rp{a}\pih_{AB}&=-2 c^{-1}\la\rp{a} \hat\theta_{AB}+(-{}\rp{a}\Omega(\log \bb+2\log c)+c^{-1}\la\rp{a} \p_c v_D \delta^{CD}-c^{-1}\la\rp{a} \tr\chi\nn\\
&-2c^{-1}\la\rp{a} L \log c)\ga_{AB}\label{2.4.2.22}.
\end{align}
and other components of ${}\rp{a}\pih$ (except ${}\rp{a}\pih_{L\Lb}$) appeared above are the same as the corresponding  ${}\rp{a}\pi$ components in Proposition \ref{3.22.6.21}.
\begin{align}
\ckk\J[S]_\Lb&=\sn^A \piSh_{A\Lb}+(\zeta_A+3\zb_A)\piSh_{A\Lb}\nn\\
&-\f12 (L+\tr\chi-2k_{\bN\bN})(\piSh_{\Lb\Lb})\label{5.01.2.21}\\
\bJ[S]_B&=\eth[S]_B-\f12 \chih_{AB}\piSh_{\Lb A}-\frac{3}{4}\tr\chi \piSh_{B\Lb} \nn\\
&+(\zeta_A+\zb_A) \piSh_{AB}-\f12(\sn_L \piSh_{\Lb B}-k_{\bN\bN} \piSh_{\Lb B}),
\label{5.01.4.21}
\end{align}
where
\begin{equation}\label{7.15.3.22}
\begin{split}
\eth[S]_B:&=\sn_A\piSh_{AB}-\f12\sn_B {}\rp{S}\ss=2\tir \sn_A\chih_{AB}-\f12\sn_B\fm{S}\\
&=\tir (2\bR_{C4CB}+k_{B\bN}\tr\chi+\sn_B \tr\chi+\sn_B k_{\bN\bN}-2\chih_{BC}k_{C\bN})
\end{split}
\end{equation}
with the second line obtained from (\ref{dchi}).
\end{lemma}
\begin{remark}
In $\bJ[S]_\Lb$ which will be given in Proposition \ref{error_terms}, the highest order term $\sdiv \ze$ is cancelled. 
\end{remark}
\begin{proof}
 (\ref{7.15.2.22}) and (\ref{7.15.3.22}) can be obtained in view of (\ref{7.15.1.22}) directly.
 Using Proposition \ref{6.7con}, with $F_{L\Lb}=F_{LL}=0$, we directly compute
\begin{align*}
\bd_A F_{B\Lb}&=\sn_A F_{B\Lb}-F(\f12\chi_{AB}\Lb+\f12\chib_{AB} L, \Lb)-F(e_B, \chib_{AC} e_C+k_{A\bN}\Lb)\\
&=\sn_A F_{B\Lb}-\f12\chi_{AB} F_{\Lb\Lb}-\chib_{AC} F_{BC}-k_{A\bN} F_{B\Lb},\\
\bd_\Lb F_{L\Lb}&=\Lb (F_{L\Lb})-F(\bd_\Lb L, \Lb)-F(L, \bd_\Lb \Lb)\\
&=\Lb(F_{L\Lb})-F(2\zeta_A e_A+k_{\bN\bN}L, \Lb)-F(L, (-2\zeta_A+2k_{\bN A})e_A -k_{\bN\bN} \Lb)\\
&=-2\zeta_A F_{A\Lb}+2(\zeta_A-k_{\bN A}) F_{L A},\\
\displaybreak[0]
\bd_L F_{\Lb \Lb}&=L(F_{\Lb \Lb})-2 F(\bd_L \Lb, \Lb)\\
&=L(F_{\Lb \Lb})-2 F(2\zb_A e_A+k_{\bN\bN}\Lb ,\Lb)\\
&=L(F_{\Lb \Lb})-4 F_{A\Lb} \zb_A -2 k_{\bN\bN} F_{\Lb \Lb},\\
\bd_A F_{BC}&=\sn_A F_{BC}-F(\f12 \chi_{AB} \Lb+\f12\chib_{AB}L, e_C)-F(e_B, \f12 \chi_{AC}\Lb+\f12 \chib_{AC}L)\\
&=\sn_A F_{BC}-\f12 \chi_{AB}F_{\Lb C}-\f12 \chib_{AB}F_{L C}-\f12 \chi_{AC} F_{B \Lb}-\f12 \chib_{AC} F_{BL},\\
\bd_\Lb F_{LB}&=\sn_\Lb F_{LB}-F(\bd_\Lb L, e_B)-F(L,  -\zeta_B \Lb+(\zeta_B-k_{\bN B})L)\\
&=\sn_\Lb F_{LB}-2\zeta_A F_{AB}-k_{\bN\bN} F_{L B}\\
\bd_L F_{\Lb B}&=\sn_L F_{\Lb B}-F(\bd_L \Lb, e_B)-F(\Lb, \zb_B L)\\
&=\sn_L F_{\Lb B}-2\zb_A F_{AB}-k_{\bN\bN} F_{\Lb B}.
\end{align*}
We can obtain (\ref{6.28.5.21})-(\ref{5.01.4.21}) by using the above calculations and Proposition \ref{3.22.6.21} and Proposition \ref{8.18.3.21}.
\end{proof}

\subsection{Preliminary estimates of ${}\rp{X}\pi$ and $\bJ[X]$}\label{Jformulas} In the sequel, we provide estimates of  ${}\rp{X}\pi$ and $\bJ[X]$ for completing  the first order energy estimates.
  
Symbolically, we rewrite the following formulas in Proposition \ref{8.18.3.21}, Proposition \ref{3.22.6.21} and (\ref{2.4.2.22}) below with the help of Lemma \ref{dcom_s}. 
\begin{align}
&\piSh_{\Lb A}=2\tir \ud\bA, \quad \piSh_{AB}-(1+\tir k_{\bN\bN})\ga_{AB}=2\tir (\bA_{g,2}+\bA_b)\nn\\
&\f12\piSh_{\Lb\Lb}-2-2(\wp-2)\tir\bT\varrho=-2\tir \mho+\tir\bA_b+\tir [L\Phi]\nn\\
&{}\rp{S}\ss=4+2\tir(\tr\chi-\frac{2}{\tir})\nn\\
&{}\rp{a}\pih_{\Lb \Lb}=\ud\bA\c \Omega^A; \quad {}\rp{a}\pih_{\bN A}=-c^{-1}\sn\la^a+c^{-1}\la^a\c \ud\bA\nn\\
&{}\rp{a}\pih_{L A}={}\rp{a}\pih_b+{}\rp{a}\pih_{L A}^+,\label{7.16.2.22}\\
&{}\rp{a}\pih_b=c^{-2}v^i \tensor{\ud\ep}{^a_i_l}e_A^l, \quad {}\rp{a}\pih_{L A}^+=-c^{-1}\sn \la^a +c^{-1}\eta({}\rp{a}\Omega)+c^{-1}\la^a\c \bA_{g,1}\nn\\
&{}\rp{a}\pih_{AB}=c^{-1}\tir^{-1}\la^a+c^{-1}\la^a \c\bA+{}\rp{a}\Omega (\log c+\log \bb)\nn\\
&{}\rp{a}\ss=\Omega \varrho+c^{-1}([L\Phi]+\bA_b+\tir^{-1}) \la^a\nn
\end{align}
For convenience we set
\begin{equation}\label{syb_1}
\begin{split}
{}\rp{a}\ud\pih&={}\rp{a}\pih_{AB}, {}\rp{a}\pih_{\Lb \Lb}; \quad {}\rp{a}\pih_b={}\rp{a}\pih^-_{L A};  {}\rp{a}\pih^+= {}\rp{a}\pih^+_{AL}, {}\rp{a}\pih_{\bN A}\\
\piShb&=\piSh_{\Lb A}, \piSh^\sharp=\piSh_{AB}, \piSh_{\Lb\Lb} 
\end{split}
\end{equation}

In view of (\ref{6.28.3.21}), Lemma \ref{3.23.2.23}, (\ref{7.16.2.22}) and (\ref{12.17.1.23}), we give the symbolic formulas for $\ckk \J[\Omega], \ckk \J[S]$ and derivatives of ${}\rp{a}\ss$ and ${}\rp{S}\ss$.
\begin{proposition}\label{error_terms}
 \begin{footnote}{
Regarding ${}\rp{Y}\pih_{A\bT}$, ${}\rp{Y}\pih_{A\bN}$ as $S_{t,u}$ tangent 1-tensors, $\sn_X$ is understood as the operator $\bd_X$ applied to the $S_{t,u}$ tangent tensors. $\sn_X{}\rp{Y}\pih_{AB}$ is understood in the same way.}
\end{footnote}
\begin{align}
{}\rp{a}\ckk\J_L&=(\ud\bA+\sn_A{}){}\rp{a}\pih_{AL}\nn\\
{}\rp{a}\ckk\J_\Lb&=\sn_A{}\rp{a}\pih_{A \Lb}+\ud\bA {}\rp{a}\pih_{A\bN}+\bA_{g,1}{}\rp{a}\pih_{LA}+(L+\tr\chi-2k_{\bN\bN}){}\rp{a}\pih_{\Lb\Lb}\nn\\
{}\rp{a}\bJ_A&={}\rp{a}\wt{\eth}_A+(\fB+\bA_{g,1}){}\rp{a}\pih_b+{}\rp{a}\pih_{AB}\ud\bA+(\tir^{-1}+\bA+\Lb \varrho)({}\rp{a}\pih_{\bN \bA}+{}\rp{a}\pih^+_{\bT A})\nn\\
\displaybreak[0]
{}\rp{a}\wt\eth&=-2\sn_A(c^{-1}\la_a \hat\theta_{AB})+\sn{}\rp{a}\Omega\log c+c^{-1}\hat\theta \sn \la_a+\f12 \tr\thetac \sn\la_a+\ud \bA^2\c{}\rp{a}\Omega\nn\\
&-\{\sn_\bT{}\rp{a}\pih_b+\sn_\bT(c^{-1}\la_a\c \bA_{g,1}+c^{-1}\eta({}\rp{a}\Omega))+\sn_L(c^{-1}\la_a \ud \bA)\}\nn\\
\sn_A {}\rp{a}\pih_{A\bN}&=c^{-1}(\sn \varrho+\sdiv)(\la_a\c \ud\bA+\sn\la_a)\nn\\
\sn_A{}\rp{a}\pih^+_{A\bT}&=c^{-1}(\sn \varrho+\sdiv)\big(-\la_a \c \ud\bA+\la_a\bA_{g,1}+\eta({}\rp{a}\Omega)\big)\nn\\
\sD{}\la_a&=c^2\big({}\rp{a}\Omega (\tr\thetac) +(\bA_{g,2}+\bA_b)\la_a \c \thetac\big)\nn\\
Y{}\rp{a}\ss&=Y{}\rp{a}\Omega \varrho+c^{-1}(Y \varrho+Y)\big(([L\Phi]+\tir^{-1}+\bA_b)\la^a\big), \,Y=L, \Lb\nn\\
\bJ[S]_\Lb&=\tir(\sn\bA_{g,1}+\varpi+ k_{\bN\bN}\tr\chi+\tr\chi\mho+\ud\bA \ud\bA+\bA_g^2)+2(L+\tr\chi-2k_{\bN\bN})(\tir k_{\bN\bN}+\Lb \tir)\nn\\
\bJ[S]_A&=\eth[S]_A+\tir(\chih+k_{\bN\bN}+\bA_b)\ud\bA+\sn_S^{\le 1}\ud\bA\nn\\
\eth[S]_A&=\tir(\bR_{B4BA}+k_{A\bN}\c \chi+\sn \bA_b+\sn k_{\bN\bN})\nn\\
L{}\rp{S}\ss&=
\tir(\bA^2_{g,2}+(h, \fB, \bA_{g,1}, \eh) ([\bar\bp\Phi]+\eh)+\sD\varrho+LL\varrho+\tr\chi(\tr\chi-\frac{2}{\tir}))\nn\\
\Lb{}\rp{S}\ss&=2\tir\tr\chi \mho+2\tir\mu.\nn
\end{align}

\end{proposition}
\begin{proof}
$L{}\rp{S}\ss$ follows by combining (\ref{6.3.1.23}) and (\ref{3.20.1.22}). The last line follows by using (\ref{8.13.5.22}). Most of the above formulas have been either obtained or can be derived straightforwardly. The derivation of formulas of ${}\rp{a}\wt{\eth}$  and $\bJ[S]_\Lb$ involves delicate cancellations. In what follows we will only give detailed calculations for these two. 

We first derive the formula of ${}\rp{a}\wt{\eth}$.  It is direct to compute
\begin{align*}
\sn_\Lb{}\rp{a}\pih_{LB}&=\sn_\Lb {}\rp{a}\pih_b+\sn_\Lb\big(-c^{-1}\sn \la_a+c^{-1}\la_a \c \bA_{g,1}+c^{-1}\eta({}\rp{a}\Omega)\big)\\
\sn_L{}\rp{a}\pih_{\Lb B}&=\sn_L{}\rp{a}\pih_b+\sn_L\big(c^{-1}\sn\la_a+c^{-1}\la_a \c \bA_{g,1}+c^{-1}\eta({}\rp{a}\Omega)+c^{-1}\la\ud \bA\big).
\end{align*}
Combining them gives
\begin{align*}
\f12(\sn_\Lb{}\rp{a}\pih_{L B}+\sn_L{}\rp{a}\pih_{\Lb B})&=\sn_\bN(c^{-1}\sn\la_a)+\sn_\bT{}\rp{a}\pih_b+\sn_\bT(c^{-1}\la_a\c \bA_{g,1}\\
&+c^{-1}\eta({}\rp{a}\Omega))+\sn_L(c^{-1}\la \ud \bA).
\end{align*}
Now we substitute (\ref{7.15.2.22}) and the above identity to the definition of ${}\rp{a}\wt{\eth}$ to obtain 
\begin{align*}
{}\rp{a}\wt\eth&=-2\sn_A(c^{-1}\la \hat\theta_{AB})-\sn_B{}\rp{a}\Omega\log \bb-\sn_\bN(c^{-1} \sn\la_a)\\
&-\{\sn_\bT{}\rp{a}\pih_b+\sn_\bT(c^{-1}\la_a\c \bA_{g,1}+c^{-1}\eta({}\rp{a}\Omega))+\sn_L(c^{-1}\la \ud \bA)\}.
\end{align*} 
Using (\ref{7.04.2.21}), (\ref{5.17.1.21}) and Proposition \ref{2.19.4.22} we have
\begin{align*}
\sn_\bN(c^{-1}\sn\la_a)&=c^{-1}\sn \bN \la_a +c^{-1}[\sn_\bN, \sn]\la_a+\bN(c^{-1})\sn\la_a\\
&=-\sn\Omega\log(\bb c)-\sn\log c\Omega\log (\bb c)-c^{-1}\hat \theta \sn\la_a-\f12\tr\thetac\sn\la_a-\ud \bA^2\c \Omega.
\end{align*}
Note that the first term in the above can be canceled by the second term on the right-hand side in ${}\rp{a}\wt\eth$. Thus 
\begin{align*}
{}\rp{a}\wt\eth&=-2\sn_A(c^{-1}\la \hat\theta_{AB})+\sn\Omega\log c+c^{-1}\hat\theta \sn \la+\f12 \tr\thetac \sn\la_a+\ud \bA^2\c\Omega\\
&-\{\sn_\bT{}\rp{a}\pih_b+\sn_\bT(c^{-1}\la_a\c \bA_{g,1}+c^{-1}\eta({}\rp{a}\Omega))+\sn_L(c^{-1}\la \ud \bA)\}
\end{align*}
as stated.

Using Proposition \ref{8.18.3.21}, the formula of $\Lb {}\rp{X}\ss$ in Proposition \ref{error_terms}, (\ref{dze}) and (\ref{4.17.1.24}), we have 
\begin{align*}
\bJ[S]_\Lb&=\sn_A{}\rp{S}\pih_{A\Lb}-\f12\Lb{}\rp{S}\ss+(\ze_A+3\zb_A) {}\rp{S}\pih_{A\Lb}-\f12(L+\tr\chi-2k_{\bN\bN})({}\rp{S}\pih_{\Lb\Lb})\\
&=2\tir(\sdiv\ze-\f12\mu)-2\tir \sdiv \zb-\tir \tr\chi\mho+\tir\ud\bA^2-\f12(L+\tr\chi-2k_{\bN\bN})({}\rp{S}\pih_{\Lb\Lb})\\
&=2\tir (-k_{\bN\bN}\tr\chi-\bA_g^2)+\tir \varpi+\tir \sn\bA_{g,1}-\tir\tr\chi \mho+\tir\ud\bA^2-\f12(L+\tr\chi-2k_{\bN\bN})({}\rp{S}\pih_{\Lb\Lb})
\end{align*}
as desired.
\end{proof}
Next, we control the deformation tensors $\piSh$ and $\pioh$, $\J[\Omega], \J[S]$ and derivatives of ${}\rp{a}\ss, {}\rp{S}\ss$.
\begin{proposition}\label{5.24.1.21}
 For the nontrivial components of $\piSh$, there hold
\begin{equation}\label{5.25.2.21}
|\piShb|\les \Delta_0\l t\r^{\delta},\quad |\piSh^\sharp|\les 1
\end{equation}
\begin{equation}\label{7.5.1.21}
\begin{split}
&\sn_\Omega\piSh_{AB}, \sn_\Omega \piSh_{\Lb\Lb}, \sn_\Omega\piSh_{\Lb A}=O(\Delta_0 \l t\r^{\delta})_{L_\omega^4}\\
&\| (\sn_X^l (\piSh_{\Lb\Lb}+4\tir \mho, \piSh_{AB})-\vs^-(X^l)O(\tir\fB)\|_{L^2_u L_\omega^2}\les \l t\r^{\delta} \Delta_0,\, l=1,2\\
&\|\sn_X^l \piSh_{\Lb A}\|_{L^2_u L_\omega^2}\les \l t\r^\delta\Delta_0, l=1,2
\end{split}
\end{equation}
\begin{equation}\label{7.5.2.21}
\sn_S\piSh^\sharp=\tir([\Lb \Phi]+\mho)+O(\l t\r^{-\frac{3}{4}+\delta}\Delta_0^\f12)_{L_\omega^4}, \sn_S\piShb=O(\Delta_0 \l t\r^{\delta})
\end{equation}
\begin{equation}\label{10.30.1.21}
\begin{split}
& X(\fm{S})=\vs(X)\tir\Lb \varrho+(1-\vs(X))O(\Delta_0\log \l t\r\bb^{-1})_{L_\omega^4}\\
&\|\bb^{-\f12}X^{\le 2}\fm{\Omega}\|_{L^2_\Sigma}+\sum_{\vs^-(X^n)=0}\|X^n \fm{S}\|_{L^2_\Sigma}\les \l t\r^{1+\delta} \Delta_0, \quad n=1,2
\end{split}
\end{equation}
For the nontrivial components of $\pioh$, and ${}\rp{\Omega}\ss$, with $Y=L, \Lb, e_A$, there hold the following estimates 
\begin{equation}\label{5.21.1.21}\left\{
\begin{array}{lll}
\piohb=O(\l t\r^{\delta}\Delta_0), \pioh_{A\bN}, \pioh^+_{LA}, \Delta_0^\f12\pioh_{LA}=O(\l t\r^{-\frac{3}{4}+\delta}\Delta_0),\\
\pioh^+_{LA}, \Delta_0^\f12\pioh_{LA}=O(\l t\r^{-1+\delta}\Delta_0)_{L_\omega^4}, \pioh_{\bN A}=O(\l t\r^{-1+2\delta} \Delta_0)_{L_\omega^4}\\
\sn_X\pioh_b=O(\Delta_0^\f12\l t\r^{-1+\delta}),\sn_X\piohb=O(\l t\r^{\delta}\Delta_0)_{L_\omega^4}, \\
\sn_S\pioh_{A\bN}=O(\l t\r^{-1+2\delta}\Delta_0)_{L_\omega^4}, \sn_\Omega\pioh_{A\bN}=O(\l t\r^{-\frac{3}{4}+\delta}\Delta_0)_{L_\omega^4}\\
\sn_X\pioh_{LA}=O(\l t\r^{-\frac{3}{4}-\frac{1}{4}\vs(X)+\delta}\Delta_0^\f12)_{L^4_\omega} 
\end{array}\right.
\end{equation}
\begin{equation}\label{3.25.1.22}
\begin{split}
&\sn_X^l\piohb=O(\l t\r^\delta\Delta_0)_{L^2_u L_\omega^2}, \pioh_{LA}, \l t\r^{-\delta}\pioh_{\bN A}=O(\l t\r^{\delta-1}\Delta_0)_{L_u^2 L_\omega^2}\\
 &\|\sn_X^l\pioh_{LA}, \l t\r^{-\delta\vs^+(X^l)}\sn_X^l\pioh_{\bN A}\|_{L_u^2 L_\omega^2}\les \l t\r^{\delta-1+\f12(l-1)(1-\vs^+(X^l))}\Delta_0, l= 1, 2
 \end{split}
\end{equation}
\begin{equation}\label{10.8.2.23}
\sn_X^{l}\pioh_{L A}=O(1)\sn_X^{l}(\eta(\Omega), {}\rp{a}v^*, \sn\la)+O(\Delta_0^2\l t\r^{-2+\frac{\max(l-1,0)}{4}+2\delta})_{L^4_\omega}, l\le 2
\end{equation}
\begin{equation}\label{8.30.2.23}
{}\rp{a}\wt\eth=\fB+O(\l t\r^{-1+2\delta}\Delta_0^2)+O(\l t\r^{-2+2\delta}\Delta_0)_{L_\omega^4}
\end{equation}
\begin{equation}\label{7.21.1.22}
\begin{split}
{}\rp{a}\ckk\J_L&=c^2\Omega \tr\thetac+\sdiv\big(\eta(\Omega)\big)+\sdiv{}\rp{a}\pih_b+O(\l t\r^{-2+2\delta}\Delta_0^2)_{L^4_\omega}\\
Y{}\rp{a}\ss&=Y \Omega \varrho+\bA_{g,1}+\max(-\vs(Y),0) \ud\bA+O(\tir^{-2})\la\\
&+O(\l t\r^{-3+\frac{1}{4}+\frac{3}{4}\max(-\vs(Y),0)+2\delta}\Delta_0^\frac{3}{2}\log \l t\r^{\max(-\vs(Y),0)})_{L^4_\omega}.
\end{split}
\end{equation}
\begin{equation}\label{7.3.1.22}\left\{
\begin{array}{lll}
Y{}\rp{a}\ss=O(\l t\r^{-2-\min(\vs(Y),0)+\delta}\Delta_0)_{L_\omega^4}, O(\l t\r^{-2-\min(\vs(Y),0)+\delta}\Delta_0)_{L_u^2 L_\omega^2}\\
L{}\rp{a}\ss=O(\l t\r^{-\frac{7}{4}+2\delta}\Delta_0),\quad \|L{}\rp{S}\ss\|_{L_u^2 L_\omega^4}\les \Delta_0\l t\r^{-2+\delta}\\
|L{}\rp{S}\ss|\les |1, \tir[\Lb \Phi], \tir\bA_{g,1}||[\bar\bp\Phi]+\eh|+\tir|\sD\varrho, LL\varrho|+|\tr\chi-\frac{2}{\tir}, \tir \bA_{g,2}^2|\\
\qquad\quad=O(\l t\r^{-\frac{7}{4}+\delta}\Delta_0^\f12).
\end{array}\right.
\end{equation}
\begin{equation}\label{5.25.1.21}
\begin{split}
& |\bJ[\Omega]_L|\les \Delta_0 \l t\r^{-1+\delta}, \|\bJ[\Omega]_L\|_{L_\omega^4}\les \Delta_0 \l t\r^{-\frac{7}{4}+\delta},  \|\bJ[\Omega]_\Lb\|_{L_\omega^4}\les \Delta_0 \l t\r^{-1+\delta};\\
&\|\bb^{-\f12}\bJ[\Omega]_L\|_{L^2_\Sigma}\les \Delta_0\l t\r^{-1+2\delta}, \|\bb^{-\f12}\bJ[\Omega]_\Lb\|_{L^2_\Sigma}\les \Delta_0\l t
\r^\delta\\
&\bJ[\Omega]_B=O(\fB)+O(\l t\r^{-1+2\delta}\Delta_0)+O(\Delta_0\l t\r^{-2+2\delta})_{L_\omega^4}
\end{split}
\end{equation}

\begin{align}
&\ckk\J[S]_L=0, \quad\|\bJ[S]_A\|_{L_\omega^4}+\|\bJ[S]_A\|_{L^2_u L_\omega^2}\les\Delta_0\l t\r^{-1+\delta},\label{6.29.3.21}\\
& \bJ[S]_\Lb=O(\tir^{-1})+O(\l t\r^{-1+\delta}\log \l t\r\Delta_0^\f12)_{L_\omega^4} \bb^{-1}\label{6.29.4.21}%
\end{align}

\end{proposition}

\begin{proof} We will frequently use (\ref{comp1}), (\ref{1.25.2.22}), Lemma \ref{3.17.2.22} and Proposition \ref{7.15.5.22} in the sequel.

In view of (\ref{7.16.2.22}), the estimates in (\ref{7.5.1.21}) and (\ref{7.5.2.21}) can be obtained by using (\ref{8.23.1.23}), (\ref{3.16.1.22}), Proposition \ref{7.15.5.22} and Proposition \ref{7.22.2.22}. Recall from (\ref{6.23.4.21}) and (\ref{8.29.1.23}) that $\fm{S}=2-2\tir k_{\bN\bN}$ and $\fm{\Omega}=2\Omega\log \bb$. (\ref{10.30.1.21}) follows by using Proposition \ref{7.15.5.22} and (\ref{8.23.2.23}). 
 
Next with $X\in \{S, {}\rp{a}\Omega\}$, due to (\ref{7.16.2.22}) we have
\begin{align}\label{7.18.5.22}
\begin{split}
&X(\pioh_{\Lb\Lb})=\sn_X( \ud\bA\c \Omega^A)\\
&\sn_X\pioh_{\bN A}=c^{-1}(X\varrho+\sn_X)(\la^a\c \ud\bA+\sn\la^a)\\
&\sn_Y \pioh_{\bT A}=\sn_Y \pioh_b+c^{-1}(Y\varrho+\sn_Y)(-\la^a \c \ud \bA+\eta(\Omega)+\la \c \bA_{g,1}), Y=\Omega, S, \tir\Lb\\
&\sn_X\pioh_{AB}=c^{-1}(X\varrho+\sn_X)(\tir^{-1}\la)+ c^{-1}(X \varrho+\sn_X)(\la \c \bA)+X\Omega  (\log c+\log \bb).
\end{split}
\end{align}
Note that it is easier to treat the multiplier $X\varrho\cdot$ which appeared simultaneously with $\sn_X$ in the above. We will ignore them in analysis.

We also can obtain by using (\ref{3.6.2.21}), (\ref{6.24.1.21}) and (\ref{3.6.1.22})-(\ref{6.28.6.21}) that 
\begin{equation}\label{7.18.1.22}
\pih_b=O(\l t\r^{-1+\delta}\Delta_0^\f12), \sn{}\rp{a}\pih_b, \sn_L{}\rp{a}\pih_b=O(\Delta_0^\f12 \l t\r^{-2+\delta}), \sdiv v^*=O(\tir^{-1} v_A).
\end{equation}
This gives the estimates of $\pih_b$ in (\ref{5.21.1.21}). Noting the decomposition of $\pioh_{LA}=\pih_b+\pioh^+_{LA}$ in (\ref{7.16.2.22}), the estimates on $\pioh_{LA}$ follow by combining the estimates on each part. 
  
 With $X=\Omega, S$ and  $l\le 2$, using (\ref{zeh}), Proposition \ref{10.16.1.22} and Proposition \ref{7.15.5.22} we have 
\begin{equation}\label{dpio}
\begin{split}
&\sn_X\sn\la= O(\l t\r^{-1+\delta+\frac{1-\vs(X)}{4}}\Delta_0)_{L_\omega^4}, \sn_X(\la \ud \bA)=O(\l t\r^{-1+2\delta}\Delta_0^2)_{L_\omega^4}\\
&\sn_X^l(\la \ud \bA)=O(\l t\r^{-1+2\delta}\Delta_0^2)_{L_u^2 L_\omega^2},\, \sn_X^l\sn \la=O(\l t\r^{-1+\delta+\f12\max(l-1,0)(1-\vs^+(X^l))}\Delta_0)_{L_u^2 L_\omega^2}\\
&\sn_X^l (\ud\bA\c \Omega)=O(\Delta_0 \l t\r^\delta)_{L_u^2 L_\omega^2},\, \sn_X^l(\bAn\c\Omega)=O(\Delta_0\l t\r^{\delta})_{L^2_\Sigma}\\
&\sn_X^l(\la \bA)=O(\l t\r^{-2+\f12\max(l-1,0)(1-\vs^+(X^l))+2\delta}\Delta_0^\frac{3}{2})_{L_u^2 L_\omega^2},\, \sn_X^l(\la \bA_{g,1})=O(\l t\r^{-1+2\delta}\Delta_0^2)_{L^2_\Sigma}.
\end{split}
\end{equation}
Using the $L^4_\omega$ estimate in the above, the remaining estimates in (\ref{5.21.1.21}) can be obtained by using (\ref{6.24.1.21}), Proposition \ref{7.15.5.22} and Proposition \ref{10.16.1.22} in view of (\ref{7.18.5.22}).

Moreover, using (\ref{3.28.3.24}) and  
 (\ref{L2BA2}),   also in view of (\ref{9.18.3.23}),
we have the rough bound
\begin{equation*}
\|\bb^{-\f12}\sn_X^{\le 2} {}\rp{a}v^*_A\|_{L^2_\Sigma}\les \l t\r^\delta\Delta_0.
\end{equation*} 
 Hence
\begin{equation*}
\|\bb^{-\f12} \sn_X^{\le 2}\pioh_b\|_{L^2_\Sigma}\les \l t\r^\delta\Delta_0.
\end{equation*}
We combine the above estimates with the $L^2$ estimates in (\ref{dpio}) to conclude (\ref{3.25.1.22}).

We further check by using Proposition \ref{10.16.1.22}, (\ref{3.6.2.21}) and (\ref{3.11.3.21}) that
\begin{equation}\label{2.7.1.24}
\|\sn_X^l(\bA_{g,1}\la)\|_{L_\omega^4}\les \l t\r^{-2+\frac{\max(l-1,0)}{4}+2\delta}\Delta_0^2,\,  l\le 2. 
\end{equation}
 (\ref{10.8.2.23}) follows immediately. 

Next we prove (\ref{7.21.1.22}). Rewrite ${}\rp{a}\ckk\J_L$ 
 symbolically
\begin{align*}
{}\rp{a}\ckk \J_L &=(\ud \bA+\sn_A){}\rp{a}\pih_{AL}^++(\ud\bA+\sdiv){}\rp{a}\pih_b.
\end{align*}
Using (\ref{5.21.1.21}) and (\ref{3.6.2.21}) we also obtain
\begin{equation*}
\ud \bA(\pioh_b, \pioh^+_{AL})=O(\l t\r^{-2+2\delta}\Delta_0^\frac{3}{2})_{L_\omega^4}, O(\l t\r^{-\frac{7}{4}+2\delta}\Delta_0^\frac{3}{2}).
\end{equation*}
Thus in view of Proposition \ref{error_terms}, Proposition \ref{7.15.5.22}, Proposition \ref{10.16.1.22}, (\ref{6.22.1.21}) and (\ref{8.3.5.23})  we derive that
\begin{align}\label{3.5.11.24}
\begin{split}
{}\rp{a}\ckk\J_L&-(c\Omega \tr\thetac+c^{-1}\sdiv\big(\eta(\Omega)\big)+\sdiv{}\rp{a}\pih_b)=O(\l t\r^{-2+2\delta}\Delta_0^\frac{3}{2})_{L_\omega^4}, O(\l t\r^{-\frac{7}{4}+2\delta}\Delta_0^\frac{3}{2})\\
Y{}\rp{a}\ss&=Y\Omega \varrho+c^{-1}(Y\varrho+Y)([L\Phi]\la+\tir^{-1}\la+\bA_b\la)\\
&=Y\Omega\varrho+O(1)(\bA_{g,1}+\max(-\vs(Y), 0)\ud \bA)+O(\tir^{-2})\la\\
&+O(\l t\r^{-3+\frac{1}{4}+\frac{3}{4}\max(-\vs(Y),0)+2\delta}\Delta_0^\frac{3}{2}\log \l t\r^{\max(-\vs(Y),0)})_{L_\omega^4},
\end{split}
\end{align}
where $Y=L$ and $\Lb$. Thus (\ref{7.21.1.22}) is proved. We can similarly obtain the pointwise estimate of $L{}\rp{a}\ss$ in (\ref{7.3.1.22}); and the estimates in the first line of (\ref{7.3.1.22}) follows by using Proposition \ref{7.15.5.22} and Proposition \ref{10.16.1.22}.

Note it follows from Proposition \ref{8.12.1.23}, Proposition \ref{7.15.5.22} and Proposition \ref{10.16.1.22} that 
\begin{equation*}
\sdiv(\eta({}\rp{a}\Omega))=O(\l t\r^{-2+\delta}\Delta_0)_{L^4_\omega},  O(\l t\r^{-1+\delta}\Delta_0)_{L^2_\Sigma}, O(\l t\r^{-\frac{7}{4}+\delta}\Delta_0).
\end{equation*}
We can conclude $\ckk \J_L=O(\l t\r^{-1+\delta}\Delta_0)$ by using (\ref{ConnH}), $|\sdiv v^*|\les \tir^{-1}|v_A|$, (\ref{3.5.11.24}) and the above estimate. Combining with the pointwise bound of $L{}\rp{a}\ss$ in (\ref{7.3.1.22}), we obtain the pointwise bound of $\bJ[\Omega]_L$ in (\ref{5.25.1.21}).

Next we consider the other estimates of $\bJ[\Omega]_L$ and the estimate of $\bJ[\Omega]_\Lb$ in (\ref{5.25.1.21}). Recall from Proposition \ref{error_terms} that 
\begin{equation*}
 {}\rp{a}\ckk\J_{\Lb}=\sn^A{}\rp{a}\pih_{A\Lb}+\ud\bA {}\rp{a}\pih_{A\bN}+\bA_{g,1}{}\rp{a}\pih_{LA}+(\tir^{-1}+k_{\bN\bN}+\bA_b+L){}\rp{a}\pih_{\Lb\Lb}.
 \end{equation*}
 In view of the formulas of $\sn_A \pioh_{L A}$ and $\sn_A \pioh_{\bN A}$ in Proposition \ref{error_terms}, using (\ref{2.7.1.24}), (\ref{8.2.2.23}), Proposition \ref{10.16.1.22} for the estimates of $\la$ and Proposition \ref{7.15.5.22}, we obtain
\begin{align*}
\sn_A\pioh^+_{L A}=O(\l t\r^{-\frac{7}{4}+\delta}\Delta_0)_{L^4_\omega},  O(\l t\r^{-2+\delta}\Delta_0)_{L^2_u L_\omega^2}\\
\sn_A\pih_{\bN A}=O(\l t\r^{-\frac{7}{4}+\delta}\Delta_0)_{L^4_\omega},  O(\l t\r^{-2+2\delta}\Delta_0)_{L^2_u L_\omega^2}.
\end{align*} 
From the estimate $\sn_S^{\le 1}\pioh_{\Lb \Lb}$ in (\ref{5.21.1.21}), we obtain
\begin{align*}
&(\tir^{-1}+k_{\bN\bN}+\bA_b+L) \pioh_{\Lb \Lb}=O(\l t\r^{-1+\delta}\Delta_0)_{L_\omega^4}, O(\l t\r^{\delta-1}\Delta_0)_{L_u^2 L_\omega^2}.
\end{align*}
Thus we conclude  by using the above estimates, (\ref{7.18.1.22}) and (\ref{6.24.1.21}) that
$$
{}\rp{a}\ckk \J_Y=O(\l t\r^{-\frac{7}{4}+\frac{3}{4}\max(-\vs(Y),0)+\delta}\Delta_0)_{L_\omega^4}, {}\rp{a}\ckk \J_Y=O(\l t\r^{-2+\max(-\vs(Y),0)(1-\delta)+2\delta}\Delta_0)_{L_u^2 L_\omega^2}.
$$
Combining with the first line of (\ref{7.3.1.22}), we conclude the estimates in the first two lines of (\ref{5.25.1.21}).

Next we use the formula of ${}\rp{a}\bJ_A$ in Proposition \ref{error_terms},  (\ref{3.6.2.21}) and (\ref{5.21.1.21}) to obtain
\begin{align}\label{7.17.3.22}
\bJ[\Omega]_A&=O(\l t\r^{-2+\delta}\Delta_0^\f12)+O(\l t\r^{-1+2\delta}\Delta_0)+{}\rp{a}\wt{\eth}_A.
\end{align}
To see (\ref{8.30.2.23}), we first write
\begin{align*}
{}\rp{a}\wt\eth_A&=\sdiv(c^{-1}\la \bA_{g,2})+\sn\Omega\log c+\ud\bA^2\Omega+\sn_L(c^{-1}\la\ud\bA)+(\bA+\tir^{-1}) \sn\la\\
&+\sn_\bT\pioh_b+\sn_\bT(c^{-1}\eta(\Omega))+\sn_\bT(c^{-1}\bA_{g,1}\la).
\end{align*}
Using (\ref{6.22.1.21}), (\ref{2.1.2.24}) for the estimates of $\sn_\Lb (\eta(\Omega))$,  (\ref{6.28.6.21}),  Proposition \ref{7.15.5.22}, and Proposition \ref{10.16.1.22}, we can obtain (\ref{8.30.2.23}). 
Substituting the estimates of ${}\rp{a}\wt\eth$ into (\ref{7.17.3.22}), we conclude the last estimate in (\ref{5.25.1.21}).

Since symbolically ${}\rp{S}\ss=1+\tir\bA_b$, the estimates of $L{}\rp{S}\ss$ in (\ref{7.3.1.22}) follow in view of Proposition \ref{error_terms} and by using Proposition \ref{7.15.5.22}.

The first identity in (\ref{6.29.3.21}) can be found in (\ref{6.28.3.21}). 
Using Proposition \ref{7.15.5.22}, Lemma \ref{5.13.11.21} (3) and Proposition \ref{error_terms}, we estimate that
\begin{align*}
&\|\eth[S]\|_{L_\omega^4}+\|\eth[S]\|_{L^2_u L_\omega^2}\les \l t\r^{-1+\delta}\Delta_0.
\end{align*}
Thus noting that $\bJ[S]_A=\eth[S]_A+O(1)\sn_S^{\le 1}\ud \bA$, we completed the proof of (\ref{6.29.3.21}) by using (\ref{3.16.1.22}).

Finally, we prove (\ref{6.29.4.21}). 
Due to (\ref{10.9.5.22}), it follows by using (\ref{6.22.1.21}), (\ref{3.16.1.22}) and (\ref{8.29.3.23}) 
\begin{align*}
(L+\tr\chi-2 k_{\bN\bN})(\tir k_{\bN\bN}+\Lb \tir)&=O(\tir^{-1})+L(\tir k_{\bN\bN}+\Lb \tir)\\
&=O(\tir^{-1})+ L\Big(\tir(\mho+\bA_b+\fB)\Big)\\
&=O(\tir^{-1})+O(1)([\Lb\Phi]+[L\Phi]+\mho)+O(\l t\r^{-2+\delta}\Delta_0^\f12)_{L_\omega^4}.
\end{align*}
In view of Proposition \ref{error_terms}, (\ref{6.29.4.21}) follows by using the above estimate, (\ref{1.27.5.24}) and Proposition \ref{7.15.5.22}.   
 
\end{proof}
\subsection{Comparison lemma}\label{comp_sec}
To derive the energy estimates, we need the following comparison result.
\begin{lemma}\label{comp}
There hold the following comparison results for scalar functions $f$, 
\begin{align}\label{7.04.7.21}
\tir \sn_Y \sn f =O(1)Y^{\le 1}\Omega f,  \, Y=\Omega, S, \tir \bN, \tir e_A.
\end{align}
\begin{equation}\label{7.03.4.21}\left\{
\begin{array}{lll}
|\bd^2_{LL}f|\les |L(Lf)|+|k_{\bN\bN}| |Lf| \\
\tir |\bd^2_{A L} f|\les |L \Omega f|+|\sn f|+|\tir \zb||Lf|\\
\tir|\bd^2_{A \Lb} f|\les |\Lb\Omega f|+|\sn f|+|\tir \ze||\bN f|+|\bA_{g,1}\c S f|\\
\tir|\bd^2_{A\Lb}f|\les |\Omega\Lb f|+|\sn f|+\tir|\bA_{g,1}\Lb f| 
\end{array}\right.
\end{equation}
With $X=\Omega, S$,  there also hold the following comparison results for the higher order terms,
\begin{align}\label{7.31.4.22}
\begin{split}
\tir^2\sn_X \sn^2 f&= O(1)X^{\le 1}\Omega^{1+\le 1} f+O(\l t\r^{-\frac{3}{4}+\delta}\Delta_0^\f12)_{L_\omega^4}\Omega f\\
\tir^2\sn_X (\bd^2_{LA}, \sn_L\sn_A) f&=O(1)(SX^{\le 1}\Omega f+X^{\le 1}\Omega f)\\
&+O(\l t\r^{-\frac{3}{4}+\delta}\Delta_0^\f12)_{L_\omega^4}\Omega f+O(1)\tir^2\sum_{l=0}^1\sn_X^l \bA_{g,1}\sn_X^{1-l}Lf\\
\tir^2\sn_X (\bd^2_{\Lb A} f)&=O(\tir)(X^{\le 1}\Lb \Omega f)+X^{\le 1}\Omega f+O(\l t\r^{\delta-\vs(X)(1-\delta)}\Delta_0^\f12)_{L_\omega^4}\Omega f\\
&+O(1)\tir^2\sum_{l=0}^1(\sn_X^l \ud\bA\sn_X^{1-l}  \bN f+ \sn_X^l \bA_{g,1}\sn_X^{1-l}Lf)\\
\tir X(LLf, \ell\bd^2_{LL} f)&=O(1)(L+\f12\tr\chi)XS^{\le 1}f+O(\tir^{-1})X^{\le 1} S f+O(\tir^{-1})Xf\\
&+\ell(1-\vs(X))O(\l t\r^{-1}\log \l t\r\Delta_0)_{L_\omega^4}\bb^{-1}Sf,\quad \ell=0,1.
\end{split}
\end{align}
\begin{align}
S\Lb \Omega f,	\tir\sn_S \sn_\Lb \sn f&=O(1)(\tir\sD \Omega f+\tir\Box_\bg \Omega f+\Lb \Omega f+L \Omega f)+O(\l t\r^{-1+\delta}\Delta_0)\sn\Omega f\nn\\
&+(O(1)+O(\l t\r^{-1+2\delta}\Delta_0^\f12)_{L_\omega^4})\sn f.\label{7.11.7.21}
\end{align}
\end{lemma}
\begin{proof}
 (\ref{7.04.7.21}) is due to (\ref{9.8.2.22}).  We derive by using Proposition \ref{6.7con} that
\begin{align}\label{7.03.3.21}
\begin{split}
\bd^2_{AB}f&=\sn_A\sn_B f-\f12 \chi_{AB} \Lb f-\f12 \chib_{AB}L f, \quad \bd^2_{AL}f=\sn_L\sn_A f+k_{A\bN}L f\\
\bd^2_{A\Lb}f&=\sn_\Lb\sn_A f+k_{A\bN} L f-2\ze \bN f, \quad \bd^2_{LL}f=L(L f)+k_{\bN\bN}L f\\
\bd^2_{\Lb A}f&=\sn_A \Lb f-\chib \c\sn f-k_{\bN A} \Lb f.
\end{split}
\end{align}
Thus the first and the last lines in (\ref{7.03.4.21}) follow directly, and the other two can be obtained by using (\ref{9.8.2.22}).

For the higher order comparison, we claim
\begin{align}
	\tir \sn_X \sn^2 f&=O(1)(\sn_X^{\le 1}\sn(\Omega f)+\sn_X^{\le 1}\sn f) +O(\l t\r^{-\frac{3}{4}+\delta}\Delta_0^\f12)_{L_\omega^4}\sn f\label{7.03.5.21}\\
	\tir\sn_X \sn^2 f&=O(1)(\tir^{-1}X\Omega^{1+\le 1}f+\sn \Omega\rp{\le 1} f)+O(\l t\r^{-\frac{3}{4}+\delta}\Delta_0^\f12)_{L_\omega^4}\sn f\label{7.5.3.21}\\
	\tir\sn_X \sn_\Lb\sn f&= O(1)(X^{\le 1}\Lb\Omega f+\sn_X^{\le 1}\sn f)+O(\l t\r^{\delta-\vs(X)(1-\delta)}\Delta_0^\f12)_{L_\omega^4}\sn f\label{7.04.3.21}\\
	\tir\sn_X \sn_L \sn f&=O(1)(LX\Omega f+ \sn_S \sn f+\sn_X^{\le 1}\sn f)+O(\l t\r^{-\frac{3}{4}+\delta}\Delta_0^\f12)_{L_\omega^4} \sn f\label{7.04.11.21}\\
	\tir XLLf&=O(1)(L+\f12\tr\chi)XSf+(O(\tir^{-1})+ O((1-\vs(X))\l t\r^{-\frac{7}{4}+\delta}\Delta_0^\f12))X^{\le 1} S f\label{7.5.5.21}
\end{align}

In view of (\ref{9.30.2.22}), to prove (\ref{7.03.5.21})-(\ref{7.04.11.21}), by using (\ref{4.22.4.22}) and (\ref{1.25.2.22}),  with $\sn_Y=\sn_A, \sn_\Lb, \sn_L$ we deduce
\begin{align}
\Omega^B \sn_X \sn_Y \sn_B f&=\sn_X(\Omega^B \sn_Y \sn_B f)-\sn_X \Omega^B \sn_Y \sn_B f\nn\\
&=\sn_X\sn_Y \Omega f-\sn_X(\sn_Y \Omega^B \sn_B f)-\sn_X \Omega^B \sn_Y \sn_B f\nn\\
&=XY\Omega f+O(1)(\sn_X \sn f+\tir\sn_Y \sn f)+\sn f \c (O(1)\nn\\
&+\min(\vs(Y)+1,1) O(\Delta_0^\f12\l t\r^{\delta-\frac{3}{4}})_{L^4_\omega}\nn\\
&+\max(-\vs(Y),0)O(\l t\r^{\delta-\vs(X)(1-\delta)}\Delta_0^\f12)_{L_\omega^4}).\label{9.1.1.23}
\end{align}
Then we apply (\ref{7.04.7.21}) to the right-hand side to obtain (\ref{7.03.5.21})-(\ref{7.04.11.21}). Here to obtain (\ref{7.04.11.21}), we also applied (\ref{7.17.6.21}). (\ref{7.5.3.21}) gives the first estimate in (\ref{7.31.4.22}).

Now by using (\ref{7.17.6.21}) we derive
\begin{align*}
&\tir XLLf=SX(\tir^{-1}Sf)+\tir|[X,L]Lf\\
\displaybreak[0]
&=O(\tir^{-1})X^{\le 1}Sf+(L+\f12\tr\chi)XSf+ O((1-\vs(X))\l t\r^{-\frac{3}{4}+\delta}\Delta_0^\f12)X L f+\vs(X) XLf \\
&=(L+\f12\tr\chi)XSf+(O(\tir^{-1})+ O((1-\vs(X))\l t\r^{-\frac{7}{4}+\delta}\Delta_0^\f12))X^{\le 1} S f
\end{align*}
as desired in (\ref{7.5.5.21}).

Using (\ref{8.23.2.23}), (\ref{6.22.1.21}) and (\ref{7.17.6.21}), we further compute
\begin{align*}
\tir X(k_{\bN\bN} Lf)&=\tir X(k_{\bN\bN}) Lf+\tir k_{\bN\bN}(LX f+[X, L]f)\nn\\
&=O(1)\{(L+\f12\tr\chi)Xf+\tir^{-1}Xf\}+(1-\vs(X))O(\l t\r^{-1}\log \l t\r\Delta_0)_{L_\omega^4}\bb^{-1} Sf.
\end{align*}
In view of (\ref{7.03.3.21})  we write
 \begin{align*}
 \sn_X \bd^2_{LA}f&=\sn_X \sn_L \sn_A f+\sn_X(k_{A\bN} Lf)\\
 \sn_X \bd^2_{\Lb A} f&=\sn_X\sn_\Lb \sn_A f+\sum_{l=0}^1(\sn_X^{1-l}\ud\bA \sn_X^l\bN f +\sn_X^{1-l}\bA_{g,1}\sn_X^l Lf).
 \end{align*}
 Using the above three calculations together with (\ref{7.04.3.21})-(\ref{7.5.5.21}), also applying (\ref{7.04.7.21}), we obtain all the remaining estimates in (\ref{7.31.4.22}).

Next we consider $\sn_S \sn_\Lb \sn_A f$.
Using (\ref{9.1.1.23}), we derive
\begin{equation}\label{7.11.6.21}
\tir\Omega^A \sn_L \sn_\Lb \sn_A f=S(\Lb \Omega f)+O(1)(\sn_S^{\le 1}\sn_A f+\tir\sn_\Lb \sn_A f)+O(\l t\r^{-1+2\delta}\Delta_0^\f12)_{L_\omega^4}\sn f.
\end{equation}
Applying (\ref{6.30.2.19}) to $\Omega f$ yields
\begin{equation*}
L \Lb \Omega f=\sD \Omega f-\Box_\bg \Omega f-(h-k_{\bN\bN})\Lb \Omega f-\hb L \Omega f+2\zb^A \sn_A \Omega f.
\end{equation*}
Using $|\tir h, \tir \hb, \tir k_{\bN\bN}|\les 1$ and $|\zb|\les \Delta_0 \l t\r^{-2+\delta}$, we derive from (\ref{7.11.6.21}) and (\ref{7.04.7.21}) that
\begin{align*}
S\Lb \Omega f, \tir\Omega^A\sn_L \sn_\Lb \sn_A f&=\tir\sD \Omega f+\tir\Box_\bg \Omega f+O(\Lb \Omega f+L \Omega f)+O(\l t\r^{-1+\delta}\Delta_0)\sn\Omega f\\
&+\Big(O(1)+O(\l t\r^{-1+2\delta}\Delta_0^\f12)_{L_\omega^4}\Big)\sn f.
\end{align*}
 This gives (\ref{7.11.7.21}).
\end{proof}

Next we compare $\bN \Omega^l v$ and $\Omega^l \bN v$ with $\bN \Omega^l \varrho$.  
\begin{lemma}\label{2.9.3.23}
With $l=1,2$ we have the symbolic formulas
\begin{align}
\Omega^l\Lb \varrho&=\Omega^l k_{\bN\bN}+\Omega^l[L\Phi]\label{2.9.1.23}\\
{}\rp{a}\Omega \hat \bN v^j&=-\sk_{\hat\bN\hat\bN} \thetac_{{}\rp{a}\Omega}^i+{}\rp{a}\Omega(\hat\bN v^i\Pic_i^{j'})\Pic_{j'}^j+\stc{k}_{\hat\bN A}\c\thetac_{{}\rp{a}\Omega}^A\hat \bN^j+{}\rp{a}\Omega\stc{k}_{\hat \bN\hat\bN} \hat\bN^j\label{5.26.1.23}\\
[\Omega \bN v]&=[\Omega L\Phi]+[\sn \Phi]\thetac(\Omega)+\fB \Omega \varrho+\Omega \Lb \varrho\label{2.9.2.23}\\
X \bN v^\|&=\fB\thetac(X)+\sn_X([\sn v])+X \log c [\sn \Phi], \quad X=\Omega, e_A\label{9.16.3.23}\\
\Lb \Omega v&=\Omega \Lb v-2\Omega^A(\ze+\zb)_A\c \bN v+\pio_{\Lb}^A \sn_A v\label{2.8.1.24}
\end{align}
where $\thetac_\Omega$ denotes the $S_{t,u}$ tangent 1-form $\thetac(\Omega, \cdot)$.
\end{lemma}
\begin{corollary}\label{6.29.1.23}
\begin{align*}
&\|[\Omega\Lb v]\|_{L^2_\Sigma}\les \|\Omega\Lb\varrho\|_{L^2_\Sigma}+\l t\r^{-1}W_2[\Omega \Phi]^\f12(t)+\|[\sn \Phi]\|_{L^2_\Sigma}+\l t\r^{-\frac{3}{4}+\delta}\Delta_0\|[L \Phi]\|_{L^2_\Sigma}\\
&\sn\Lb v-(O(\tir^{-1})\fB+\sn \Lb \varrho) =O(\l t\r^{-2+\delta}\Delta_0)_{L^2_\Sigma}, O(\l t\r^{-3+\delta}\Delta_0)_{L_\omega^4}\\
&\tir^{-1}(\Lb\Omega v-O(\bb^{-1}\ud \bA))-O(1)\sn\Lb v=O(\l t\r^{-1+2\delta}\Delta_0^\f12)_{L_\omega^4}\sn v.
\end{align*}
\end{corollary}
\begin{proof}[Proof of Lemma \ref{2.9.3.23} and Corollary \ref{6.29.1.23}]
 By using (\ref{7.04.7.19}) and (\ref{7.04.8.19}),
\begin{equation}\label{7.26.4.21}
\bN \varrho=L \varrho-\wp^{-1}(k_{\bN\bN} +L\varrho+L v_\bN).
\end{equation}
This gives (\ref{2.9.1.23}).
(\ref{5.26.1.23}) is obtained similar to (\ref{5.23.1.23}). (\ref{2.9.2.23})  and (\ref{9.16.3.23}) are direct consequences of (\ref{5.26.1.23}) and (\ref{2.9.1.23}), where we used the fact that
\begin{align*}
 \Omega[L v]=[\Omega Lv]+[\sn \Phi]\thetac(\Omega).
\end{align*}
(\ref{2.8.1.24}) is a direct consequence of (\ref{5.13.10.21}). Corollary \ref{6.29.1.23} is a consequence of (\ref{6.30.2.23}), (\ref{2.9.2.23}) and (\ref{9.16.3.23}). The last estimate follows by using (\ref{2.8.1.24}) and (\ref{5.21.1.21}). 
\end{proof}
As a consequence of Lemma \ref{comp}, we have the following result.
\begin{corollary}\label{9.2.5.23}
With $l=0,1$, there hold
\begin{align}
&\|\tir(\sn^2 X^{\le 1}\Phi, \bd^2_{LA}X^{\le 1}\Phi, \bd^2_{LL}X^{\le 1}\Phi), \bd^2_{\Lb A}\Omega^{\le 1}\varrho, \tir^{-1}\Omega \Lb \Omega^{\le 1} \varrho\nn\\
&\sn_S^{1-a}\bd^2_{\Lb A}S^a\varrho, \tir^{-1}\Omega \Lb S\varrho\|_{L^2_\Sigma}\les \l t\r^{-1+\delta}\Delta_0, \, a=0,1,\label{9.1.3.23}\\
&\|\sn_X^l \sn^2\Phi, \sn_X^l \bd^2_{LA}\Phi, (\log \l t\r)^{-l(1-\vs(X^l))} X^l(\bd^2_{LL}\Phi), \l t\r^{-1}\sn_X^l\bd^2_{\Lb A}\varrho\|_{L^2_\Sigma}\les \l t\r^{-2+\delta}\Delta_0,\label{8.15.3.21}\\
&\|\tir \sn_X\bd^2_{\Lb A}\varrho, \tir\bd^2_{\Lb A}X\varrho\|_{L^2_u L_\omega^4}\les \l t\r^{-1+\delta}\Delta_0\label{12.25.1.23}\\
&\|\sn_{X_2}^{\le 1}\sn {X_1}^{\le 1}\Phi, \sn_{X_2}^{\le 1}L X_1^{\le 1}\Phi, \tir^{-1} X^{1+\le 1}\Phi\|_{L^2_\Sigma}\les \l t\r^{-1+\delta}\Delta_0\label{9.2.3.23}\\
&\|\tir (\sn_X^{l_1} \sn X^{l_2}\Phi, \sn_X^{l_1} LX^{l_2} \Phi, \bd^2_{LL} \Phi, \bd^2_{LA}\Phi,\sn^2\Phi)\|_{L_u^2 L^4_\omega}\nn\\
&\qquad\qquad\qquad\qquad\les    \l t\r^{-2+\delta}\Delta_0,\, 0\le l_1+l_2\le 1\label{9.2.4.23}\\
&\bd^2_{\Lb A}\Phi^\mu=\min(1,\mu)O(\tir^{-1}\fB)+\sn\Lb \varrho+O(\l t\r^{-2+\delta}\Delta_0)_{L^2_\Sigma}\label{9.18.4.23}\\
&\bd^2_{\Lb A}\Phi=O(\tir^{-1})\fB+ O(\l t\r^{-2+\delta}\Delta_0)\label{12.25.3.23}\\
&\bd^2_{\Lb A} X\Phi^\mu, \sn_X\bd^2_{\Lb A}\Phi^\mu-\min(\mu, 1)O(\l t\r^{-1}\fB)\nn\\\
&\quad=O(\l t\r^{-1+\delta}\Delta_0)_{L^2_\Sigma},  O(\l t\r^{-2+\delta}\Delta_0)_{L_u^2 L_\omega^4}\label{9.5.5.23}
\end{align}
where we assumed 
\begin{equation}\label{3.3.2.24}
\|\Box_\bg\Omega\Phi\|_{L_u^2 L_\omega^4}\les \l t\r^{-\frac{11}{4}+\delta}\Delta_0
\end{equation}
 to obtain $L_u^2 L_\omega^4$ estimates in  (\ref{12.25.1.23}) and (\ref{9.5.5.23}) for the case $X=S$. 
 
 \begin{equation}\label{9.5.2.23}
 \sn \Omega \bN v-O(\l t\r^{-1})\Omega^{\le 2}\fB=O(\l t\r^{-\frac{11}{4}+2\delta}\Delta_0)_{L_u^2 L_\omega^4}, O(\l t\r^{-\frac{7}{4}+2\delta}\Delta_0)_{L^2_\Sigma}.
 \end{equation}
 With $f=v, \varrho$
 \begin{equation}\label{3.10.4.24}
 \Omega\Lb S f-O(1) \Omega \Lb f, S\Lb \Omega f-\Lb \Omega f=O(\l t\r^{-\frac{7}{4}+\delta}\Delta_0)_{L_u^2 L_\omega^4}, O(\l t\r^{-1+2\delta}\Delta_0)_{L^2_\Sigma}.
 \end{equation}
\end{corollary}
\begin{remark}
Note (\ref{12.25.1.23}) and (\ref{9.5.5.23}) will be used for deriving the second order energy estimates. Before that, the assumption (\ref{3.3.2.24}) will have been verified in (\ref{8.23.1.21}).  
\end{remark}
\begin{proof}

 Note by using Lemma \ref{comp} and Lemma \ref{5.13.11.21}, we derive 
 \begin{align}\label{9.2.1.23}
 \begin{split}
 &\|\tir \sn^2 f\|_{L^2_\Sigma}\les\|\sn\Omega^{\le 1} f\|_{L^2_\Sigma}\\
&\|\tir \bd_{LL}^2 f\|_{L^2_\Sigma}\les \l t\r^{-1} W_2[Sf]^\f12(t)+\|Lf\|_{ L_\Sigma^2}\\
&\|\tir \bd_{\Lb A}^2 f\|_{L^2_\Sigma}\les E[\Omega f]^\f12 (t)+\|\sn f\|_{L^2_\Sigma}+\|\tir\ud\bA\c \bN f\|_{L^2_\Sigma}+\Delta_0 \l t\r^{-2+\delta}\|S f\|_{L^2_\Sigma}\\
&\|\tir \bd_{LA}^2 f\|_{L^2_\Sigma}\les \|L \Omega f\|_{L^2_\Sigma}+\|\sn f\|_{L^2_\Sigma}+\Delta_0 \l t\r^{-1+\delta}\|Lf\|_{L^2_\Sigma}\\
&\tir \bd_{\Lb A}^2 f=\tir \sn_A \Lb f+O(1)\sn f+\tir\bA_{g,1}\Lb f.
\end{split}
\end{align}
 The estimates of $\sn^2 X^{\le 1}\Phi, \bd^2_{LA}X^{\le 1}\Phi, \bd^2_{LL}X^{\le 1}\Phi, \bd^2_{\Lb A}\Omega^{\le 1}\varrho$ in (\ref{9.1.3.23}) follow by using (\ref{3.12.1.21}), (\ref{3.6.2.21}) and the above estimates. The estimate of $\Omega \Lb \Omega^{\le 1} \varrho$  in (\ref{9.1.3.23}) can be obtained by using (\ref{LbBA2}), (\ref{3.6.2.21}) and (\ref{4.22.4.22}). 
The $l=0$ case in (\ref{8.15.3.21}) follows easily by using (\ref{9.2.1.23}) and (\ref{3.12.1.21}).
The estimates of $\sn_X \sn^2\Phi, \sn_X \bd^2_{LA}\Phi, X(\bd^2_{LL}\Phi)$  in (\ref{8.15.3.21}) can be obtained by using (\ref{7.31.4.22}), (\ref{3.12.1.21}), (\ref{3.6.2.21}) and (\ref{3.11.3.21}). Applying (\ref{7.31.4.22})  to $\sn_X(\bd^2_{\Lb A}\varrho)$, also using Lemma \ref{5.13.11.21} (1) and (2),  implies
\begin{align}\label{4.22.1.24}
\tir\sn_X(\bd^2_{\Lb A} \varrho)&=X^{\le 1}\Lb \Omega \varrho+O(\tir^{-1})X^{\le 1}\Omega\varrho+\sum_{l=0}^1\sn_X^l \ud \bA \sn_X^{1-l}\bN \varrho+O(\l t\r^{2\delta-2}\Delta_0^\frac{3}{2})_{L_\omega^4}.
\end{align}  
Due to Lemma \ref{5.13.11.21} (5) and (\ref{8.23.2.23}), 
\begin{align}\label{7.28.1.24}
\sum_{l=0}^1\sn_X^l \ud \bA \sn_X^{1-l}\bN \Phi=O(\l t\r^{\f12\delta-2}\Delta_0\log \l t\r^2)_{L_u^2 L_\omega^4}.    
\end{align}
Applying (\ref{7.11.7.21}) to $f=\varrho$, we obtain $S\Lb \Omega \varrho=O(\l t\r^{-1+2\delta}\Delta_0)_{L^2_\Sigma}$ by using (\ref{wave_ass}), (\ref{3.12.1.21}) and (\ref{3.11.3.21}). Substituting this estimate, (\ref{7.28.1.24}) and the $L^2_\Sigma$ bound of $\Omega^{\le 1}\Lb \Omega \varrho$ in (\ref{9.1.3.23}) to (\ref{4.22.1.24}) we can obtain the last estimate in (\ref{8.15.3.21}). The first estimate of (\ref{12.25.1.23}) can be obtained by using (\ref{3.3.2.24})(for $X=S$), (\ref{7.28.1.24}), (\ref{4.22.1.24}), (\ref{4.22.4.22}), (\ref{3.6.2.21}) and (\ref{L4BA1}).

(\ref{9.2.3.23}) and (\ref{9.2.4.23}) can be obtained by using Proposition \ref{7.15.5.22}, (\ref{7.03.4.21}) and Sobolev embedding on spheres.

The case $\Phi=\varrho$ in (\ref{9.18.4.23}) can be obtained by applying the last estimate in (\ref{9.2.1.23}) to $f=\varrho$ and using (\ref{L2BA2}). 
 If $\Phi=v$, (\ref{9.18.4.23}) follows from using the last lines in (\ref{9.2.1.23}) and Corollary \ref{6.29.1.23}, (\ref{LbBA2}) and $\sn v=O(\l t\r^{-1+\delta}\Delta_0)_{L^2_\Sigma}$ due to (\ref{3.12.1.21}). (\ref{12.25.3.23}) follows similarly by using (\ref{3.6.2.21})

Next we consider the estimates in the second line of (\ref{9.1.3.23}). The case $a=0$ has been given in the last estimate of  (\ref{8.15.3.21}). For the rest of the estimates, we derive by using (\ref{6.30.1.19}) that
\begin{align*}
\tir\sn\Lb Sf=\tir (\tir \sn)(\sD f-\Box_\bg f+(\mho+k_{\bN\bN}) Lf-h\Lb f+2\zeta^A\sn_A f).
\end{align*}

Using Proposition \ref{7.22.2.22} for treating $\mho$ and its derivatives, also using (\ref{7.31.4.22}) for treating $\Omega \sD f$, we have
\begin{align}
\tir \sn \Lb Sf&=\tir (\tir \sn) (\sD f)+\tir(\tir \sn)\Box_\bg f+\tir(\tir \sn)((\mho+k_{\bN\bN}) L f)+\tir(\tir \sn)(h\Lb f)\nn\\
&+(\tir\sn)^2 f\c \ze+\tir\sn \ze\tir\sn f\nn\\
&=\l t\r^{-1}O(1)\Omega^{1+\le 2}f +\tir(\tir \sn) \Box_\bg f+O(1)\Omega Lf+\tir \sn (\ze, \mho, k_{\bN\bN})\c\sum_{X=S, \Omega} Xf\nn\\
&+\tir h\c\tir \sn\Lb f+\tir^2\sn h\c\Lb f+\tir^2\sn^2f\c \ze,\label{8.22.2.21}
\end{align}
where in the above we used the pointwise estimates $\mho, \bb \fB=O(\l t\r^{-1})$. 

Using  Proposition \ref{7.22.2.22} for $\Omega\mho$ and (\ref{zeh}), applying (\ref{8.22.2.21}) to $f=\varrho$ and $v^i$, also using Sobolev embedding and (\ref{3.12.1.21}), (\ref{3.29.1.23}) and Proposition \ref{7.15.5.22} gives
\begin{align}
\|\tir \sn\Lb S f-\tir h\tir \sn\Lb f\|_{L^2_u L_\omega^4}&\les\l t\r^{-1}\|\Omega^{1+\le 2}f\|_{L_u^2 L_\omega^4}+\l t\r^{2\delta-2}\Delta_0^2+\|\tir \Omega \Box_\bg f\|_{L^2_u L_\omega^4}\nn\\
&+\|\bb^{-1}\Omega h, \Omega L f\|_{L^2_u L_\omega^4}\nn\\
&\les \l t\r^{-\frac{7}{4}+\delta}\Delta_0.\label{9.17.1.23}
\end{align}
Using Lemma \ref{5.13.11.21} (5) and (\ref{9.12.3.22}), we have 
$$ \sn (\ze, \mho, k_{\bN\bN})=O(\l t\r^{-2}\log \l t\r^2\Delta_0)_{L_u^2 L_\omega^2}.
$$  
Taking $L^2_\Sigma$ norm instead, with the help of the above estimates, using  (\ref{3.29.1.23}) and Proposition \ref{7.15.5.22}, we have for $f=v, \varrho$
\begin{equation}\label{3.10.5.24}
\|\sn\Lb S f-h\tir \sn\Lb f\|_{L^2_\Sigma}\les \l t\r^{-2+2\delta}\Delta_0.
\end{equation}
Combining  the above estimate with (\ref{LbBA2}), we obtain $\|\Omega \Lb S\varrho\|_{L^2_\Sigma}\les \l t\r^{\delta}\Delta_0$, as stated in (\ref{9.1.3.23}).  

In view of (\ref{9.17.1.23}), applying the last line in (\ref{9.2.1.23}) to $f=v^i, \varrho$ and using (\ref{9.2.4.23}), (\ref{3.11.3.21}) and (\ref{6.22.1.21}) we obtain
\begin{align}\label{3.3.1.24}
\begin{split}
\|\bd^2_{A\Lb} Sf-\tir h\sn\Lb f\|_{L^2_u L_\omega^4}&\les\l t\r^{-1}\|\sn S\Phi\|_{L^2_u L_\omega^4}+\|\bA_{g,1}\|_{L_\omega^4}\|\Lb S\Phi\|_{L^\infty_x}+\l t\r^{-\frac{11}{4}+\delta}\Delta_0\\
&\les \l t\r^{-\frac{11}{4}+\delta}\Delta_0.
\end{split}
\end{align}
Taking $L^2_\Sigma$ norm instead, we have for $f=v, \varrho$,
$$
\|\bd^2_{A \Lb} S f-\tir h\sn\Lb f\|_{L^2_\Sigma}\les  \l t\r^{-2+2\delta}\Delta_0.
$$
This gives $\|\bd^2_{A \Lb} S\varrho\|_{L^2_\Sigma}\les \l t\r^{-1+\delta}\Delta_0$ in view of (\ref{LbBA2}) and $h\approx\tir^{-1}$, as stated in (\ref{9.1.3.23}). Thus the proof of (\ref{9.1.3.23}) is completed. We also obtain (\ref{12.25.1.23}) for $\bd_{A\Lb}^2 S\varrho$ as a consequence of (\ref{3.3.1.24}) and (\ref{L4BA1}). Applying (\ref{9.2.1.23}) to $f=\Omega \varrho$, and using (\ref{L4BA1}), (\ref{4.22.4.22}) and (\ref{3.6.2.21}), we have 
$\tir \bd^2_{\Lb A}\Omega \varrho=O(\l t\r^{-1+\delta}\Delta_0)_{L_u^2 L_\omega^4}$, as desired in (\ref{12.25.1.23}). The proof of (\ref{12.25.1.23}) is completed.
 
The $\mu=0$ case in (\ref{9.5.5.23}) has been given in (\ref{9.1.3.23})-(\ref{12.25.1.23}). It suffices to consider it for $\Phi^i=v^i$.  
Using  Corollary \ref{6.29.1.23},  we infer from (\ref{3.3.1.24}) and the above $L^2_\Sigma$ comparison result that 
\begin{align*}
\tir\bd^2_{A\Lb}(S \Phi^\mu)&=\min(\mu, 1)O(1) \fB+\tir \sn\Lb \varrho +O(\l t\r^{-\frac{7}{4}+\delta}\Delta_0)_{L_u^2 L_\omega^4}\\
&=\min(\mu, 1)O(1) \fB+O(\l t\r^{\delta}\Delta_0)_{L^2_\Sigma}.
\end{align*}
This gives the first estimate of (\ref{9.5.5.23}) if $X=S$.

Next, for $X=S$, we consider the other estimate in (\ref{9.5.5.23}). Using (\ref{7.11.7.21}), with $f=\varrho$ or $v^i$, we infer
\begin{align*}
\|(\tir^{-1}S\Lb \Omega f,& \sn_S \sn_\Lb \sn f)-O(\tir^{-1}) \Lb \Omega f\|_{L^2_\Sigma}\\
&\les\|\sD \Omega f\|_{L^2_\Sigma}+\|\Box_\bg \Omega f\|_{L^2_\Sigma}+\l t\r^{-1}\|L \Omega f\|_{L^2_\Sigma}\nn\\
&\quad+\Delta_0\l t\r^{-2+\delta}\|\sn\Omega f\|_{L^2_\Sigma}+\l t\r^{-1}\|\sn f\|_{L^2_\Sigma}+\l t\r^{-2+2\delta}\Delta_0^\frac{1}{2}\|\bb^\f12\tir\sn f\|_{L^2_u L_\omega^4}\\
&\les\l t\r^{-2} W_2[\Omega^{1+\le 1}f](t)+\|\Box_\bg \Omega f\|_{L^2_\Sigma}+\l t\r^{-1}\|\sn f\|_{L^2_\Sigma}+\l t\r^{-3+3\delta}\Delta_0^\frac{3}{2}\log \l t\r^\f12\\
&\les \l t\r^{-2+2\delta}\Delta_0
\end{align*}
where we used (\ref{3.12.1.21}) and (\ref{wave_ass}).

Taking $L_u^2 L_\omega^4$ norm by using (\ref{9.2.4.23}), similar to (\ref{9.17.1.23}) we have for $f=\varrho, v$ that
\begin{align*}
\|(\tir^{-1}S\Lb \Omega f,& \sn_S \sn_\Lb \sn f)-O(\tir^{-1}) \Lb \Omega f\|_{L^2_u L_\omega^4}\\
&\les \l t\r^{-\frac{11}{4}+\delta}\Delta_0+\|\Box_\bg \Omega f\|_{L_u^2 L_\omega^4}\les \l t\r^{-\frac{11}{4}+\delta}\Delta_0
\end{align*}
where we used Proposition \ref{7.15.5.22} and the assumption (\ref{3.3.2.24}) to derive the last estimate. 
We conclude
\begin{equation*}
(S\Lb \Omega f, \tir \sn_S \sn_\Lb \sn f)-O(1) \Lb \Omega f=O(\l t\r^{-1+2\delta}\Delta_0)_{L^2_\Sigma}, O(\l t\r^{-\frac{7}{4}+\delta}\Delta_0)_{L_u^2 L_\omega^4}.
\end{equation*}
Using the above formula, (\ref{9.17.1.23}) and (\ref{3.10.5.24}) gives (\ref{3.10.4.24}).
Using again the above formula, with the help of Proposition \ref{7.15.5.22}, Lemma \ref{5.13.11.21} (5) and (\ref{8.23.2.23}), we derive from (\ref{7.31.4.22}) for $X=S$,
\begin{align}\label{8.16.5.21}
\begin{split}
\|\tir\sn_X& \bd^2_{\Lb A}f-O(1)\Lb\Omega f\|_{L^2_\Sigma}\\
&\les \l t\r^{-1+\delta}\Delta_0+ \l t\r^{-1}\|X^{\le 1}\Omega \Phi\|_{L^2_\Sigma}+\|\tir\fB\sn_X \ud\bA\|_{L^2_\Sigma}+\|\tir\sn_X\bN \Phi\c \ud\bA\|_{L^2_\Sigma}\\
&+\Delta_0 \l t\r^{-2+\delta}\|\tir\sn_{X}L \Phi\|_{L^2_\Sigma}+\|\sn_S  \bA_{g,1}\|_{L^\infty_x}\|S\Phi\|_{L^2_\Sigma}+\Delta_0^\frac{3}{2} \l t\r^{-2+3\delta}\log \l t\r^\f12\\
&\les \Delta_0(\log \l t\r)^2
\end{split}
\end{align}
Taking $L_u^2 L_\omega^4$ norm instead, we can similarly obtain for $f=v, \varrho$ that 
\begin{equation*}
\|\tir\sn_X \bd^2_{\Lb A}f-O(1)\Lb\Omega f\|_{L^2_u L_\omega^4}\les \l t\r^{-1+\delta}\Delta_0
\end{equation*}
 which gives
$$
\sn_S \bd^2_{\Lb A}f- O(\tir^{-1})\Lb \Omega f=O(\Delta_0(\log \l t\r)^2\l t\r^{-1})_{L^2_\Sigma}, O(\l t\r^{-2+\delta}\Delta_0)_{L_u^2 L_\omega^4}.
$$
For $f=v$, using (\ref{2.8.1.24}), Corollary \ref{6.29.1.23} and (\ref{5.21.1.21}) we have from the above estimate that
\begin{equation*}
\sn_S \bd^2_{\Lb A}v-O(1)(\tir^{-1}\fB+\sn \Lb \varrho+\tir^{-1}\bb^{-1}\ud\bA)=O(\l t\r^{-1}(\log \l t\r)^2\Delta_0)_{L^2_\Sigma},  O(\l t\r^{-2+\delta}\Delta_0)_{L_u^2 L_\omega^4}.
\end{equation*}
Using Proposition \ref{7.15.5.22},  the proof of (\ref{9.5.5.23}) for $X=S$ is completed. 
  
Next we consider the case that $X=\Omega$.  
 It follows by using Proposition \ref{7.15.5.22} and (\ref{8.23.1.23}) that
\begin{align}\label{3.3.3.24}
\sum_{l=0}^l\sn_\Omega^l \bA_{g,1}(\sn_\Omega^{1-l}\Lb\Phi, \Lb\sn_\Omega^{1-l}\Phi)=O(\l t\r^{-2+2\delta}\Delta_0)_{L^2_\Sigma}, O(\l t\r^{-3+2\delta}\Delta_0)_{L_u^2 L_\omega^4}.
\end{align}
Using (\ref{7.03.3.21}), we compute
\begin{align*}
\bd^2_{\Lb A}{}\rp{a}\Omega v&=\sn_A \Lb {}\rp{a}\Omega v-\chib\sn{}\rp{a}\Omega v+\bA_{g,1}\Lb{}\rp{a}\Omega v\\
&=\sn_A({}\rp{a}\Omega \bN v)+\sn_AL{}\rp{a}\Omega v+\sn_A([\bN, {}\rp{a}\Omega]v)+\chib\sn{}\rp{a}\Omega v+\bA_{g,1}\Lb{}\rp{a}\Omega v,
\end{align*}
where the last line is a schematic identity.

Using (\ref{5.13.10.21}) the commutator term is written explicitly below 
\begin{align*}
\sn_A([\bN, \Omega]v)&=\sn_A(\pioh_{\bN B}\sn_B v+\Omega^B(\zb+\ze)_B \bN v)\nn\\
&=O(\l t\r^{-1}\log \l t\r^2\Delta_0)_{L^2_\Sigma}, O(\l t\r^{-2+\delta}\Delta_0)_{L_u^2 L_\omega^4},
\end{align*}
where the estimate follows by using Proposition \ref{7.15.5.22}, Lemma \ref{5.13.11.21} (5) and (\ref{5.21.1.21}).
Hence  using the above estimate, Proposition \ref{7.15.5.22} and (\ref{3.3.3.24}), we infer 
\begin{align}\label{9.5.1.23}
\bd^2_{\Lb A}{}\rp{a}\Omega v-\sn_A({}\rp{a}\Omega \bN v)=O(\l t\r^{-1}\log \l t\r^2\Delta_0)_{L^2_\Sigma}, O(\l t\r^{-2+\delta}\Delta_0)_{L_u^2 L_\omega^4}
\end{align}

Using Lemma \ref{6.30.4.23} and noting that $\bN v^\|=[\sn v]$ and $[\bN v]=\fB$ we derive
\begin{align}\label{3.3.5.24}
\begin{split}
\tir \sn \Omega \bN v^\|&=\tir\sn\sn_\Omega[\sn v]+O(1)\sn_\Omega^{\le 1}\fB+O(\l t\r^{-\frac{3}{4}+\delta}\Delta_0)\Omega^{\le 1}[\sn v]\\
&+O(\l t\r^{-\frac{3}{4}+\delta}\Delta_0)_{L_\omega^4}\fB,\\
[\tir \sn \Omega \bN v]&=\tir \sn\Omega\fB+O(1)\sn_\Omega^{\le 1}[\sn v]+O(\l t\r^{-\frac{3}{4}+\delta}\Delta_0)\Omega^{\le 1}\fB\\
&+O(\l t\r^{-\frac{3}{4}+\delta}\Delta_0)_{L_\omega^4}[\sn v].
\end{split}
\end{align}

Hence using Proposition \ref{7.15.5.22}, we infer
\begin{equation*}
\sn \Omega \bN v-O(\l t\r^{-1})\Omega^{\le 2}\fB=O(\l t\r^{-\frac{11}{4}+2\delta}\Delta_0)_{L_u^2 L_\omega^4}, O(\l t\r^{-\frac{7}{4}+2\delta}\Delta_0)_{L^2_\Sigma},
\end{equation*}
as stated in (\ref{9.5.2.23}). Substituting the above estimate into (\ref{9.5.1.23}) with the help of Proposition \ref{7.15.5.22} for $\Omega^{1+\le 1}\fB$, we complete the proof of the first estimate in (\ref{9.5.5.23}) for $\Phi^i=v^i$.

 Using (\ref{7.03.3.21}), similar to (\ref{9.5.1.23}), we derive
\begin{align*}
\sn_{{}\rp{a}\Omega}\bd^2_{\Lb A} v&=\sn_{{}\rp{a}\Omega}(\sn_A\Lb v+\chib\sn v+\bA_{g,1}\Lb v)
\end{align*}
Using Proposition \ref{7.15.5.22}, (\ref{8.24.4.23}) and (\ref{8.24.5.23}) we obtain
\begin{align*}
\sn_{{}\rp{a}\Omega}(\chib\c \sn v)=O(\l t\r^{-2+2\delta}\Delta_0)_{L^2_\Sigma}, O(\l t\r^{-3+2\delta}\Delta_0)_{L_u^2 L_\omega^4}.
\end{align*}
Applying (\ref{3.3.3.24}), the above estimate and Proposition \ref{7.15.5.22} yields
\begin{align*}
\sn_{{}\rp{a}\Omega}\bd^2_{\Lb A} v+2\sn_{{}\rp{a}\Omega}(\sn_A\bN v)=O(\l t\r^{-2+2\delta}\Delta_0)_{L^2_\Sigma}, O(\l t\r^{-3+2\delta}\Delta_0)_{L^2_u L_\omega^4}. 
\end{align*}

Similar to (\ref{9.5.2.23}), we have
\begin{equation*}
 \sn_\Omega\sn \bN v-O(\l t\r^{-1})\Omega^{\le 2}\fB=O(\l t\r^{-\frac{11}{4}+\delta}\Delta_0)_{L_u^2 L_\omega^4}, O(\l t\r^{-\frac{7}{4}+\delta}\Delta_0)_{L^2_\Sigma}.
\end{equation*}
Combining the above two estimates, using Proposition \ref{7.15.5.22} for $\Omega^{1+\le 1}\fB$, we proved the second estimate in (\ref{9.5.5.23}) for $X=\Omega$ and $\Phi=v$. 
\end{proof}

\section{Lower order energies}\label{low_energy}
In this section, under the assumptions (\ref{3.12.1.21})-(\ref{6.5.1.21}), we provide the lower order energy estimates by commuting one vector field $X\in\{\Omega, S\}$ with the wave operator. We begin by deriving the base-order commutator estimates.

Recast the formula for $\sP[X, f]$ from (\ref{5.18.3.21}),  
\begin{align}
	\sP[{}\rp{a}\Omega, f]-{}\rp{a}\bJ^\a\p_\a f+\f12 {}\rp{a}\pih_{LA}\bd^2_{\Lb A} f&=\piohb\c(\sn^2 f, \bd^2_{LL}f)+\pioh_{\Lb A} \bd^2_{LA}f\label{6.23.1.23}\\
\sP[S, f]-\bJ[S]^\a\p_\a f&=\piShb\bd^2_{LA}f+(\piSh_{AB}\sn^2 f+ \piSh_{\Lb \Lb}\bd^2_{LL}f) \label{6.24.2.23}.
	\end{align}
For future reference, for treating the second term in the line of (\ref{6.23.1.23}), we write symbolically
\begin{equation}\label{6.24.1.23}
\bJ[X_1]^\a \c \sn_X^l \p_\a f:=\bJ[X_1]_L \sn_X^l \Lb f+\bJ[X_1]_\Lb \sn_X^l L f+\bJ[X_1]_A \sn_X^l \sn f
\end{equation}
where $l=0, 1$ and $X\in \{\Omega, S\}$.  

Using (\ref{5.25.1.21})-(\ref{6.29.4.21}), we write
\begin{align}
{}\rp{a}\bJ^\a \c \sn_X^l \p_\a f&={}\rp{a}\bJ_L \sn_X^l\Lb f+( O(\l t\r^{-1+2\delta}\Delta_0)+O(\fB))\sn_X^l \sn f\nn\\
&+O(\l t\r^{-1+\delta}\Delta_0)_{L_\omega^4}(\sn_X^l Lf, \l t\r^{-1+\delta} \sn_X^l \sn f)\label{6.24.3.23}\\
\bJ[S]^\a \c \sn_X^l \p_\a f&=\frac{1}{4}L{}\rp{S}\ss\c \sn_X^l \Lb f+O(\l t\r^{-1+\delta}\Delta_0)_{L_\omega^4}\sn_X^l \sn f\nn\\
&+(O(\tir^{-1})+O(\l t\r^{-1+\delta}\log \l t\r\Delta_0^\f12)_{L_\omega^4}\bb^{-1})\sn_X^l Lf.\label{6.24.4.23}
\end{align}
\begin{lemma}\label{9.18.5.23}
Let $l=0,1$. There hold that
\begin{align}
 \| {}\rp{a}\pih_{LA}(\bd^2_{\Lb A} X^l\Phi, \sn_X^l \bd^2_{\Lb A}\Phi)\|_{L^2_\Sigma}&\les\l t\r^{-1} \sum_{X=\Omega, S}\|X\Phi\|_{L^2_u L_\omega^2}+\l t\r^{-2+2\delta}\log \l t\r^\f12\Delta_0^\frac{5}{4}\label{9.18.6.23}\\
  \|\bb \tir\sn_X\pioh_{LA} \bd^2_{\Lb A}\Phi\|_{L^2_u L_\omega^2}&\les \l t\r^{-\frac{7}{4}+\delta}\Delta_0^\frac{5}{4}+\l t\r^{-1}\sum_{X=\Omega, S}\|X\Phi, X\eta(\Omega)\|_{L^2_u L_\omega^2}\label{9.18.6.23+}
\end{align}
\end{lemma}
\begin{proof}
Using (\ref{12.25.3.23}), (\ref{9.5.5.23}), we deduce for $l=0,1$,
 \begin{align*}
 {}\rp{a}\pih_{LA}&(\bd^2_{\Lb A} X^l\Phi, \sn_X\bd^2_{\Lb A} \Phi)={}\rp{a}\pih_{LA}(O(\tir^{-1})\fB+O(\l t\r^{-2+\delta}\Delta_0)_{L^2_u L_\omega^4}).
 \end{align*}

In view of (\ref{10.8.2.23}), (\ref{5.8.1.21}), (\ref{9.18.3.23}) and (\ref{12.19.1.23}), we derive
\begin{align*}
\|{}\rp{a}\pih_{LA}\fB\|_{L^2_\Sigma}&\les \|(\pioh_b, \sn \la, \eta(\Omega))\fB\|_{L^2_\Sigma}+\l t\r^{-2+2\delta}\Delta_0^2\\
&\les\sum_{X=\Omega, S}\|X\Phi\|_{L^2_u L_\omega^2}+\l t\r^{-1+\delta}\Delta_0^\frac{5}{4}.
\end{align*} 
 Hence we conclude (\ref{9.18.6.23}) with the help of (\ref{5.21.1.21}). 
  
Next we apply (\ref{9.18.4.23}) to write
 \begin{align*}
 \sn_X\pioh_{LA}\bd^2_{\Lb A}\Phi=\sn_X\pioh_{LA} (O(\l t\r^{-1})\fB+O(\l t\r^{-2+\delta}\Delta_0)). 
 \end{align*}
 We again use (\ref{10.8.2.23}) to treat $\sn_X\pioh_{LA}$.
Due to (\ref{6.24.2.21}), (\ref{6.28.6.21}) and (\ref{9.18.3.23}), $\|X v^*\|_{L^2_u L_\omega^2}\les \sum_{X=S, \Omega}\|X\Phi\|_{L^2_u L_\omega^2}$. Also due to Proposition \ref{10.16.1.22}
\begin{align*}
\|\bb\tir \sn_X\pioh_{LA}\fB\|_{L_u^2 L_\omega^2}&\les\|\sn_X(\sn\la)\|_{L_\omega^2}\|\bb \tir \fB \|_{L_u^2 L_\omega^\infty}+\|\sn_X(\pioh_b, \eta(\Omega)) \|_{L^2_u L_\omega^2}+\l t\r^{-2+2\delta}\Delta_0^2 \\
&\les \l t\r^{-\frac{3}{4}+\delta}\Delta_0^\frac{5}{4}+\sum_{X=\Omega, S}\|X\Phi, X(\eta(\Omega))\|_{L^2_u L_\omega^2}.
\end{align*}
Using (\ref{3.25.1.22}) to treat the remaining term, we then conclude (\ref{9.18.6.23+}).
\end{proof}
 
 \begin{proposition}\label{9.5.8.23} Under the assumptions (\ref{3.12.1.21})-(\ref{6.5.1.21}),
  if $f=X^{\le 1}\Phi$ with $X\in \{\Omega, S\}$, it holds
\begin{align}\label{9.2.6.23}
&\|\sP[{}\rp{a}\Omega, f]\|_{L^2_\Sigma}\les\l t\r^{-\frac{7}{4}+\delta}\Delta_0^\frac{5}{4}+\|\fB \sn f\|_{L^2_\Sigma}+\l t\r^{-2}\sum_{X=\Omega, S}\|X\Phi\|_{L^2_\Sigma}.
\end{align}
There hold for all $f=X^{\le 1}\Phi$
\begin{align}
 \|\sP[S, f]-\frac{1}{4}L{}\rp{S}\ss \Lb f\|_{L^2_\Sigma}&\les  \l t\r^{-2} \left(W_2[\Omega f]^\f12(t)+W_2[Sf]^\f12(t)\right)+ \l t\r^{-1} \|L f\|_{L^2_\Sigma}\nn\\
 &+\l t\r^{-1+\delta}\Delta_0^\f12\|\tir(\log \l t\r\bb^{-\f12} L f, \bb^\f12 \sn f)\|_{L_u^2 L_\omega^4}\label{8.21.2.21},\\
  \|\bb\tir^2 L{}\rp{S}\ss\Lb f\|_{L_u^1 L_t^2 L_\omega^2}&\les \Delta_0^\frac{5}{4}.\label{12.20.1.23}
\end{align}
\end{proposition}
\begin{remark}
Using (\ref{9.2.4.23}), the last line in (\ref{8.21.2.21}) can be bounded by $\l t\r^{-2+2\delta}\log \l t\r\Delta_0^\frac{3}{2}$. This fact may be used without mentioning explicitly.   
\end{remark}

\begin{proof}
 By using (\ref{5.21.1.21}), (\ref{6.23.1.23}) and (\ref{6.24.3.23}), we derive symbolically
\begin{align}
\sP[\Omega, f]&+\f12 {}\rp{a}\pih_{LA}\bd^2_{\Lb A} f+\f12\Lb f \bJ[\Omega]_L\nn\\
&= O(\l t\r^{\delta}\Delta_0)(\sn^2 f+\bd^2_{LL}f) +O(\l t\r^{-\frac{3}{4}+\delta}\Delta_0^\f12)\bd^2_{LA} f\nn\\
&+O(\l t\r^{-1+\delta}\Delta_0)_{L_\omega^4}( L f+\l t\r^{-1+\delta}\sn f)\nn\\
&+O(\fB+\l t\r^{-1+2\delta}\Delta_0)\sn f.\label{2.14.3.24}
\end{align} 
For $f=\Phi, \Omega\Phi, S\Phi$ we have from (\ref{8.23.1.23}) that
\begin{equation*}
\Lb X^l\Phi=O(1)\fB+ l(1-\vs(X^l))O(\l t\r^{-1+\delta}\Delta_0).
\end{equation*}
Then using (\ref{12.19.1.23}) and (\ref{5.25.1.21}), we infer
\begin{align*}
\|\Lb f \bJ[\Omega]_L\|_{L^2_\Sigma}&\les \|\bb\fB\|_{L_u^2 L_\omega^\infty}\|\tir\bb^{-\f12} \bJ[\Omega]_L\|_{L^\infty_u L_\omega^2}+\l t\r^{-1+\delta}\Delta_0 \|\bJ[\Omega]_L\|_{L^2_\Sigma}\\ 
&\les \Delta_0^\frac{5}{4}\l t\r^{-\frac{7}{4}+\delta}. 
\end{align*}
Applying Corollary \ref{9.2.5.23} and (\ref{12.19.1.23}) to the right-hand of (\ref{2.14.3.24}), for $f=X^{\le 1} \Phi$, we bound
\begin{align*}
\|\sP[\Omega, f]+\f12 {}\rp{a}\pih_{LA}\bd^2_{\Lb A} f\|_{L^2_\Sigma}&\les \Delta_0^\frac{5}{4}\l t\r^{-\frac{7}{4}+\delta}+\|\fB \sn f\|_{L^2_\Sigma}.
\end{align*}
Applying (\ref{9.18.6.23}) to the left-hand side, we obtain (\ref{9.2.6.23}).
 

Next we prove (\ref{8.21.2.21}).
 By using (\ref{5.25.2.21}), (\ref{6.24.2.23}) and (\ref{6.24.4.23})  we estimate
 \begin{align}\label{3.5.2.24}
 \begin{split}
 \sP[S, f]&=\frac{1}{4}L{}\rp{S}\ss \Lb f+O(1)(\sn^2 f+\bd^2_{LL} f+\l t\r^{-1}Lf)+O(\Delta_0 \l t\r^\delta)\bd^2_{LA}f\\
 &+O(\l t\r^{-1+\delta}\Delta_0^\f12)_{L_\omega^4}(\Delta_0^\f12\sn f+\log \l t\r\bb^{-1} Lf).
 \end{split}
 \end{align}
 Applying (\ref{9.2.1.23}), (\ref{9.1.3.23}) to the right-hand side, we obtain (\ref{8.21.2.21}).
 
 Using (\ref{7.3.1.22}), (\ref{8.23.1.23}) and (\ref{12.19.1.23}), we obtain 
 \begin{align*}
 \|\bb\tir^2 L{}\rp{S}\ss\Lb f\|_{L_u^1 L_t^2 L_\omega^2}&\les \|\bb\fB\tir\|_{L^2_u L_t^\infty L_\omega^\infty}\|\tir L{}\rp{S}\ss\|_{L_u^2 L_t^2 L_\omega^2}\\
 &+\Delta_0\|\l t'\r^{-1+\delta}\bb \tir^2 L{}\rp{S}\ss\|_{L_u^2 L_t^2 L_\omega^2}\les \Delta_0^\frac{5}{4},
 \end{align*}
as stated.
 \end{proof}

Next we give the first order energy estimates by commuting the vector-fields $X=\Omega, S$ with the  wave equations for $\Phi$ and applying energy inequalities in Section \ref{mul_1}. 
\begin{proposition}[The first order energy estimates]\label{1steng}
Under the assumptions (\ref{3.12.1.21})-(\ref{6.5.1.21}), we have the following results.
\begin{enumerate}
\item With $X=S, \Omega$, there hold 
\begin{equation}\label{8.21.4.21}
\begin{split}
&E[\Omega\varrho](t)+ F_0[\Omega \varrho](\H_u^t)\les \La_0^2+\Delta_0^\frac{5}{2}\\
&WFIL_2[X\Phi](\D_u^t)\les (\log \l t\r)^\M(\La_0^2+\Delta_0^\frac{5}{2})
\end{split}
\end{equation}
\begin{align}
&\|\bb^{-\f12}v, \Sc(X^l\Phi),\Ac(X^l\Phi), \tir^2(\bd^2_{LL}\Phi, \bd^2_{LA}\Phi, \sn^2\Phi)\|_{L^2_\Sigma}\nn\\
&\qquad\qquad\qquad\qquad\qquad\les (\log \l t\r)^{\f12\M}(\La_0+\Delta_0^\frac{5}{4}), \quad l=1,2\label{10.10.2.23}\\
&\|\bb^\f12X^{l} \Big(\Box_\bg \Phi,\, \N(\Phi,\bp\Phi)\Big) \|_{L^2_\Sigma}\les \l t\r^{-2}\log \l t\r^{\f12\M+l}(\La_0+\Delta_0^\frac{5}{4}), \, l=0,1\label{8.26.4.21}\\
&\|\Er_1(\Phi, X)\|_{L^2_\Sigma}\les \l t\r^{-2}\Big(\log\l t\r^{\f12\M+1}(\La_0+\l t\r^{\delta}\Delta_0^\frac{5}{4})+\l t\r^{(1-\vs(X))(\frac{1}{4}+\delta)}\Delta_0^\frac{5}{4}\Big)\label{8.26.3.21}\\
&\begin{array}{lll}
\| \Box_\bg X\Phi\|_{L^2_\Sigma}\les \l t\r^{-2}\Big(\log\l t\r^{\f12\M+1}(\La_0+\l t\r^\delta\Delta_0^\frac{5}{4})+\l t\r^{(1-\vs(X))(\frac{1}{4}+\delta)}\Delta_0^\frac{5}{4}\Big)\\
\|\Box_\bg \Omega\Phi\|_{L_u^2 L_\omega^4}\les \l t\r^{-\frac{11}{4}+\delta}\Delta_0
 \end{array}
 \label{8.23.1.21}\\
&\|\bb^{-\f12}L{}\rp{S}\ss\|_{L^2_\Sigma}\les \l t\r^{-1}\log \l t\r^{\f12\M+1}(\La_0+\Delta_0^\frac{5}{4}).\label{8.21.1.22}
\end{align}
\item There hold 
\begin{align}\label{7.13.5.22}
\begin{split}
&\|\tir\sn_X^{\le 1}\bA_{g,1}\|_{L^2_\Sigma}\les \log \l t\r^{\f12\M}(\La_0+\Delta_0^\frac{5}{4}),  \|\sn_\Lb \eh\|_{L^2_\Sigma}\les\l t\r^{-2}(\Delta_0^\frac{5}{4}\l t\r^\delta+\log \l t\r^{\f12\M+1}\La_0)\\
&\|\tir(\log \l t\r)^{-1}(\sn k_{\bN\bN}, \sn\Lb\varrho), \tir \sn_\Lb (\sn \varrho, \zb), (\log \l t\r)^{-2}\bb^{-\f12}\sn_S^{\le 1}\ze\|_{L^2_\Sigma}\les \La_0+\Delta_0^\frac{5}{4}. 
\end{split}
\end{align}
\begin{align}\label{10.11.3.23}
\begin{split}
&\|\bR_{ABCL}, (\log \l t\r)^{-1}(\l t\r^{-1}\bR_{ABC\Lb}, \bR_{AB43})\|_{L^2_\Sigma}\les\l t\r^{-2} \log \l t\r^{\f12\M}(\La_0+\Delta_0^\frac{5}{4}) \\
&\|\widehat{\bR_{A4B4}}, \widehat{\bR_{A3B4}}\|_{L^2_\Sigma}\les \l t\r^{-2} (\log \l t\r^{\f12\M}\La_0+\l t\r^\delta\Delta_0^\frac{5}{4}\Big).
\end{split}
\end{align}
\begin{align}
&\|\bb^{-\f12}\sn_S^{\le 1}(\tir\bA_b)\|_{L^2_\Sigma}\les\log\l t\r^{\f12\M+1}(\Delta_0^\frac{5}{4}+\La_0),\, \|\tir^2\bb^{-1}\sn_S^{\le 1}\bA_{g,2}\|_{L_u^2 L_\omega^2}\les \l t\r^\delta(\La_0+\Delta_0^\frac{5}{4})\label{8.8.6.22}\\
&\l t\r^2\Big\|\zb, \sn\varrho, \eh\Big\|_{L_\omega^4}\les \log \l t\r^{\frac{\M}{4}+\f12}(\La_0+\Delta_0^\frac{5}{4})\label{10.11.2.23}\\
&\begin{array}{lll}
S[\Lb \Phi], \tir\Lb[L\Phi]=O(1)[\Lb \Phi]+O(\l t\r^{-1}\log \l t\r^{\f12\M+1}(\Delta_0^\frac{5}{4}+\La_0))_{L^2_\Sigma}\\
X\Lb\Phi^\mu=O(1)(\vs(X)+\mu)[\Lb\Phi]\\
\qquad\qquad+O\Big((\vs(X)\l t\r^{-1}\log \l t\r^{\f12\M}+(1-\vs(X)))\log \l t\r(\Delta_0^\frac{5}{4}+\La_0)\Big)_{L^2_\Sigma}
\end{array}\label{8.2.1.22}\\
&\|\Sc(\Omega \Lb \Phi), \log \l t\r^{-1}\Sc(\Lb\Omega\Phi)\|_{L^2_\Sigma}\les\log \l t\r(\Delta_0^\frac{5}{4}+\La_0) \label{9.20.4.22}
\displaybreak[0]
\end{align}
\begin{equation}\label{11.26.2.23}\left\{
\begin{array}{lll}
\|\la\|_{L_\omega^4}\les (\log \l t\r)^{\frac{\M}{4}+\frac{3}{2}}(\Delta_0^\frac{5}{4}+\La_0);\, \|\bb^{-\f12}\sn \la\|_{L^2_\Sigma}\les\log \l t\r^{\f12\M+1}(\La_0+\Delta_0^\frac{5}{4})\\
\|\bb^{-\f12}\log \l t\r^{-1}\sn_S\sn\la, \sn S\la, L S^{\le 1}\la\|_{L^2_\Sigma}\les (\log \l t\r)^{\f12\M}(\Delta_0^\frac{5}{4}+\La_0),\\
\|\bb^{-\f12}(L\Lb\la, \log \l t\r^{-1}\tir^{-1}\Lb \la)\|_{L^2_\Sigma}\les(\La_0+\Delta_0^\frac{5}{4})\log \l t\r.
\end{array}\right.
\end{equation}
\begin{equation}
\begin{array}{lll}
\Big\|\log \l t\r^{-2}\sn_S^{\le 1}(\pioh_{A\bN}), \log \l t\r^{-1}\sn_S^{\le 1} \pioh^+_{L A},\\
\qquad\qquad\qquad\qquad\sn_S^{\le 1} \pioh_b\Big\|_{L^2_u L_\omega^2}\les \l t\r^{-1} (\Delta_0^\frac{5}{4}+\La_0)\log \l t\r^{\f12\M}\\
\|\sn_S^{\le 1} \piohb\|_{L^2_u L_\omega^2}\les (\Delta_0^\frac{5}{4}+\La_0)(\log \l t\r)^2
\end{array}\label{7.26.2.22}
\end{equation}
\begin{align}
\|\bJ[\Omega]_L\|_{L_u^2 L_\omega^2}&\les \l t\r^{-2}\int_0^t \l t'\r^{-1}W_2[\Omega^2\Phi]^\f12(t')+\l t\r^{-2}(\log \l t\r^{\f12\M+1}\La_0+\l t \r^{2\delta}\Delta_0^\frac{5}{4})\label{5.7.2.24} 
\end{align}
\end{enumerate}
\end{proposition}
We remark that the estimate of $\sn_S^{\le 1}\bA_{g,2}$ in (\ref{8.8.6.22}) will be further improved. The above estimates improved some of the lower order estimates in the assumptions (\ref{3.12.1.21})-(\ref{6.5.1.21}). The $X=\Omega$ case in (\ref{8.23.1.21}) will be further improved later.  We will summarize the improvements of the bootstrap assumptions, after completing the top order energy.

We divide the proof into two parts: for (1) we mainly apply energy regime given in Section \ref{mul_1}; for (2) we apply the results in (1) together with the geometric estimates given in Section \ref{low_ricci} and Section \ref{Jformulas}.  
\begin{proof}
Let us decompose $\Er_1(X, \Phi)=\wt \Er_1(X, \Phi)-\vs(X)\f12 \bJ[X]_L\c \Lb\Phi$, where  
$$\wt\Er_1(X,\Phi)=\fm{X}\Box_\bg\Phi+ \sP[X,\Phi]+\vs(X)\f12\bJ[X]_L\Lb \Phi.$$ 
 In view of (\ref{5.02.3.21_1}), applying (\ref{9.2.6.23}), noting $\fm{\Omega}=\ud\bA\c \Omega$ due to Proposition \ref{3.22.6.21}, we derive by using (\ref{3.6.2.21}) and (\ref{3.29.1.23}) that
\begin{equation}\label{12.21.1.22}
\begin{split}
\|\Er_1(\Omega, \Phi)\|_{L^2_\Sigma}&\les \|\fm{\Omega}\Box_\bg \Phi\|_{L^2_\Sigma}+\|\sP[\Omega,\Phi]\|_{L^2_\Sigma}\\
&\les \l t\r^{-\frac{7}{4}+\delta}\Delta_0^\frac{5}{4}+\l t\r^{-2}\sum_{X=S, \Omega}\|X\Phi\|_{L^2_\Sigma}
\end{split}
\end{equation}
Next we bound $\|\wt\Er_1(S, \Phi)\|_{L^2_\Sigma}$. By using $|\fm{S}|\les 1$, (\ref{8.21.2.21}) and (\ref{5.02.3.21_1})  we derive
\begin{align}
\|\wt\Er_1(S, \Phi)\|_{L^2_\Sigma}&\les \|\fm{S}\Box_\bg \Phi\|_{L^2_\Sigma}+\|\sP[S,\Phi]-\frac{1}{4} L{}\rp{S}\ss \Lb \Phi\|_{L^2_\Sigma}\nn\\
&\les\|\Box_\bg \Phi\|_{L^2_\Sigma}+\l t\r^{-2+2\delta}\log \l t\r\Delta_0^\frac{3}{2}+ \l t\r^{-2}\sum_{X=\Omega, S}WL_2[X \Phi]^\f12(t). \label{12.21.2.22}
\end{align}
Using Proposition \ref{6.24.10.23} with $n=1$, Lemma \ref{10.10.3.23} and Proposition \ref{7.15.5.22}  we deduce
\begin{align}
\|X^{\le 1}\Box_\bg \Phi\|_{L^2_\Sigma}&\les \l t\r^{-1}\| X^{\le 1}L\Phi\|_{L^2_\Sigma}+(1-\vs(X))\l t\r^{-2+2\delta}\Delta_0^\frac{3}{2}+\l t\r^{-\frac{11}{4}+2\delta}\log \l t\r^\f12\Delta_0^{\frac{3}{2}}\label{12.21.3.22}\\
&\les \l t\r^{-2}\Big((1-\vs(X))\l t\r^{2\delta}\Delta_0^\frac{3}{2}+\l t\r^{-\frac{3}{4}+2\delta}\Delta_0^\frac{3}{2}\log \l t\r^\f12\Big)\nn\\
&+\l t\r^{-2}\sum_{X=S, \Omega}WL_2[X\Phi]^\f12(t)\nn.
\end{align}
With the help of the above estimates on $\wt\Er_1(X, \Phi)$ and $\|\bb X \Box_\bg \Phi\|_{L^2_\Sigma}$, we estimate the error integral for $\sum_{X=S, \Omega}WFIL_2[X\Phi](\D_{u_1}^{t_1})$ below
\begin{align*}
&\int_{\D_{u_1}^{t_1}}\{\sum_{X=\Omega, S}|\bb^\f12\tir (\wt\Er_1(X, \Phi)+X\Box_\bg\Phi)||\bb^{-\f12}\tir (L+\f12\tr\chi)X\Phi|\}d\mu_g dt\nn\\
&\les\sum_{X=S,\Omega}(\int_{u_1}^{u_*}F_2[X\Phi](\H_u^{t_1}) du+\int_0^t \|\bb^\f12\tir (\wt\Er_1(X, \Phi)+X\Box_\bg\Phi)\|^2_{L^2_\Sigma}dt')\nn\\
&\les \Delta_0^\frac{3}{2}+\int_0^{t_1}  \l t\r^{-2+}\sum_{X=S, \Omega}WL_2[X\Phi](t')d t'+\sum_{X=S,\Omega}\int_{u_1}^{u_*}F_2[X\Phi](\H_u^{t_1}) du.
\end{align*}
Noting that $F_2[X\Phi](\H_{u_*}^{t_1})=0$, using (\ref{12.20.1.23}) in Proposition \ref{9.5.8.23} to treat the remaining term in $\Er_1(S,\Phi)$, using the above error integral, applying Proposition \ref{MA2} to $\psi=X\Phi$ and using Proposition \ref{12.21.1.21} for the initial data, we obtain 
\begin{align*}
&\sum_{X=\Omega, S}WFIL_2[X\Phi](\D_{u_1}^{t_1})\nn\\
&\les\Delta_0^\frac{5}{2}+\La_0^2+\sum_{X=\Omega, S}\{\int_0^{t_1} \sup_{\Sigma_{t'}}(-[\Lb\varrho]_{-}+\l t'\r^{-1-\delta})WL_2[X\Phi](t') dt'\nn\\
&+\int_{\D_{u_1}^{t_1}} \big| \Box_\bg X\Phi\c (L+h)X\Phi\big|\tir^2 \bb d\mu_{\ga} du dt+\int_{u_1}^{u_*} F_2[X\Phi](\H_u^{t_1}) du\}\nn\\
\displaybreak[0]
&\les\La_0^2+\Delta_0^\frac{5}{2}+\int_0^{t_1}  \l t'\r^{-2+}\sum_{X=S, \Omega}WL_2[X\Phi](t')+\sum_{X=S,\Omega}\int_u^{u_*}F_2[X\Phi](\H_{u'}^{t_1}) du'\\
&+\sum_{X=\Omega, S}\int_0^{t_1} (\sup_{\Sigma_{t'}}(-[\Lb\varrho]_{-})+\l t'\r^{-1-\delta})WL_2[X\Phi](t') dt'.\nn
\end{align*}

Using Lemma \ref{1.6.4.18}, we then conclude
\begin{align*}
\sum_{X=\Omega,S}WFL_2[X\Phi](\D_{u_1}^{t_1})\les (\Delta_0^\frac{5}{2}+\La_0^2)(\f12\wp\log(\frac{\l t_1\r}{2})+1)^\M.
\end{align*}
Next, consider the error integral for $E[\Omega\varrho](t)$. Substituting the above estimates to (\ref{12.21.1.22}) and (\ref{12.21.3.22}), we derive with the help of Cauchy-Schwarz inequality that
\begin{align}
&\int_{\D^{t_1}_{u_1}}|\Er_1(\Omega, \varrho)+\Omega \Box_\bg \varrho||\bd_\bT \Omega\varrho|\label{8.10.2.23}\\
\displaybreak[0]
&\les \int_0^{t_1} \l t'\r^{1+\delta}\|\Er_1(\Omega, \varrho)+\Omega \Box_\bg\varrho\|^2_{L^2_\Sigma}(t') dt'+\int_0^{t_1} \l t'\r^{-1-\delta}E[\Omega\varrho](t').\nn\\
&\les \Delta_0^\frac{5}{2}+\La_0^2+\int_0^{t_1} \l t'\r^{-1-\delta}E[\Omega\varrho](t').\nn
\end{align}
Using Proposition \ref{12.21.1.21} for the initial data, applying Proposition \ref{10.10.3.22} to $\Omega\varrho$, we bound 
\begin{equation}\label{8.11.2.23}
\begin{split}
E[\Omega\varrho](t_1)+F_0[\Omega\varrho](\H_{u_1}^{t_1})&\les\int_{\D_{u_1}^{t_1}}|\Box_\bg \Omega\varrho \bT \Omega\varrho|+E[\Omega\varrho](0)+F_0[\Omega\varrho](\H_{u_*}^{t_1})\\
&\les \Delta_0^\frac{5}{2}+\La_0^2+\int_0^{t_1} \l t'\r^{-1-\delta}E[\Omega\varrho](t')
\end{split}
\end{equation}
where we used the fact that $F_0[\Omega\varrho](\H_{u_*}^{t_1})=0$. Using the standard Gronwall's inequality, 
 we conclude
\begin{align*}
E[\Omega \varrho](t)+F_0[\Omega\varrho](\H_u^t)\les \La_0^2+\Delta_0^\frac{5}{2}.
\end{align*}
We can obtain (\ref{10.10.2.23}) from the obtained bound on $\sum_{X=S, \Omega}WFIL_2[X\Phi](\D_{u_1}^{t_1})$ by using Lemma \ref{6.30.4.23} and (\ref{9.18.3.23}).

Combining (\ref{7.10.5.22}) with (\ref{10.10.2.23}) implies
\begin{equation*}
\|\bb^{-\f12}\tir(\tr\chi-\frac{2}{\tir})\|_{L^2_\Sigma}\les \log \l t\r^{\f12\M+1}(\La_0+\Delta_0^\frac{5}{4})
\end{equation*}
which is the lower order estimate for $\bA_b$ in (\ref{8.8.6.22}).

Using  the above estimate, it follows by using (\ref{7.3.1.22}), Lemma \ref{comp}, (\ref{10.10.2.23}) and Proposition \ref{7.15.5.22} that
\begin{align*}
\|\bb^{-\f12}L {}\rp{S}\ss\|_{L^2_\Sigma}&\les \|\bb^{-\f12}(\tir (\sD\varrho, LL\varrho), \tr\chi-\frac{2}{\tir}, \bar\bp\Phi, \eh, \bA_{g,2}^2\tir)\|_{L^2_\Sigma}\nn\\
&\les\l t\r^{-1} \log \l t\r^{\f12\M+1}(\La_0+\Delta_0^\frac{5}{4}).
\end{align*}
(\ref{8.21.1.22}) is proved. (\ref{8.26.4.21}) follows from (\ref{8.21.4.21}), (\ref{10.10.2.23}), Proposition \ref{6.24.10.23} and (\ref{8.23.2.23}). 
 
 Refining (\ref{8.21.2.21}) by using $\bd^2_{LA}\Phi=O(\l t\r^{-2}\log \l t\r^{\f12\M}(\Delta_0^\frac{5}{4}+\La_0))_{L^2_\Sigma}$  and (\ref{8.21.4.21}), we can derive due to $\fm{S}=O(1)$,
\begin{align*}
\|\Er_1(S, \Phi)\|_{L^2_\Sigma}\les \l t\r^{-2} \log \l t\r^{\f12\M+1}(\La_0+\Delta_0^\frac{5}{4}+\l t\r^\delta\Delta_0^\frac{3}{2}).
\end{align*}
The case $X=\Omega$ of (\ref{8.26.3.21}) follows by substituting the lower order estimates in (\ref{8.21.4.21}) to (\ref{12.21.1.22}). 
In view of $\Box_\bg X \Phi=X \Box_\bg \Phi+\Er_1[\Phi, X]$,  using (\ref{8.26.4.21}) and (\ref{8.26.3.21}), we can obtain the first estimate in (\ref{8.23.1.21}). 

Using (\ref{9.2.4.23}), Proposition \ref{7.15.5.22},  and (\ref{2.14.3.24}),  we obtain
\begin{equation*}
\sP[\Omega, \Phi]+\f12 {}\rp{a}\pih_{LA}\bd^2_{\Lb A} \Phi+\f12\Lb \Phi\bJ[\Omega]_L=O(\l t\r^{-3+3\delta}\Delta_0)_{L_u^2 L_\omega^4}.
\end{equation*}
Moreover, due to (\ref{5.21.1.21}) and (\ref{12.25.3.23})
\begin{align*}
{}\rp{a}\pih_{LA}\bd^2_{\Lb A} \Phi=O(\l t\r^{-3+2\delta}\Delta_0)_{L_u^2 L_\omega^4}.
\end{align*}
Using (\ref{5.25.1.21}) and $\bb\Lb\Phi=O(\l t\r^{-1})$, we derive
\begin{equation*}
\bb\bJ[\Omega]_L \Lb\Phi=O(\l t\r^{-\frac{11}{4}+\delta}\Delta_0)_{L_\omega^4}.
\end{equation*}
Thus we have 
\begin{equation*}
\|\sP[\Omega, \Phi]\|_{L_u^2 L_\omega^4}\les \l t\r^{-\frac{11}{4}+\delta}\Delta_0. 
\end{equation*}
Due to (\ref{4.3.3.21}) and (\ref{3.6.2.21})
\begin{equation*}
\fm{\Omega}\Box_\bg \Phi=O(\l t\r^{-3+2\delta})_{L_u^2 L_\omega^\infty}.
\end{equation*}
Hence 
$$
[\Box_\bg, \Omega]\Phi=O(\l t\r^{-\frac{11}{4}+\delta}\Delta_0)_{L_u^2 L_\omega^4}.
$$
Combining the above estimate with the second estimate in (\ref{3.29.1.23}) with $l=0$, we conclude the second estimate in (\ref{8.23.1.21}).

Next we consider the estimates in (2) by using the results in (1). 

We first prove (\ref{7.13.5.22}). Since $\bA_{g,1}$ are some specific combinations of components of $\sn\Phi^\dagger$, we conclude the same estimates for $\|\sn_X^{\le 1}\bA_{g,1}\|_{L^2_\Sigma}$ due to (\ref{10.10.2.23}).
We next derive by using (\ref{8.21.4.21}), (\ref{5.13.10.21}), Lemma \ref{5.13.11.21} (5), (\ref{5.21.1.21}) and (\ref{12.19.1.23}) 
\begin{align*}
\|\tir\sn \Lb \varrho\|_{L^2_\Sigma}&\les E[\Omega\varrho]^\f12(t)+\|[\Omega, \Lb]\varrho\|_{L^2_\Sigma}\\
&\les \La_0+\log \l t\r \Delta_0^\frac{5}{4}.  
\end{align*}
Noting that $k_{\bN\bN}=\Lb\varrho+[L\Phi]$,
 we can then obtain the estimate of $\|\sn k_{\bN\bN}\|_{L^2_\Sigma}$ in (\ref{7.13.5.22}) by also using (\ref{10.10.2.23}).
 
Noting that $\ze=\sn\log \bb-\zb$,  using (\ref{10.29.1.22}), (\ref{8.5.1.22+}), the first estimate of (\ref{7.13.5.22}) and the bound on $\|\sn k_{\bN\bN}\|_{L^2_\Sigma}$, we obtain
\begin{equation*}
\|\bb^{-\f12}\sn_S^l \sn \log \bb\|_{L^2_\Sigma}\les \log \l t\r^{2}(\Delta_0^\frac{5}{4}+\La_0), \, l=0,1.
\end{equation*}
 We then obtain the estimate of $\ze$ and $\sn_L\ze$ by using the  above estimate and the first estimate of (\ref{7.13.5.22}).
 
 Recall from (\ref{5.23.1.23}) that
 \begin{align*}
 \sn_\Lb [\sn\Phi]=[\sn_\Lb \sn\Phi]+\sn\Phi^A \Lb \bN^A+\Lb\log c[\sn\Phi]
 \end{align*}
 Using (\ref{9.8.2.22}), the fact that $\Lb \bN^A=\ud \bA$, (\ref{3.6.2.21}), (\ref{8.21.4.21}) and (\ref{10.10.2.23}), we have $\sn_\Lb [\sn\Phi]=O(\l t\r^{-1}(\Delta_0^\frac{5}{4}+\La_0))_{L^2_\Sigma}.$ Thus the estimates in the second line in (\ref{7.13.5.22}) are proved. 

Using (\ref{12.19.1.23}) and (\ref{3.11.3.21}), we derive
\begin{align}\label{5.13.2.24}
\|\chih\c \fB\|_{L^2_\Sigma}\les \Delta_0^\frac{5}{4}\l t\r^{-2+\delta}
\end{align}
Due to (\ref{10.10.2.23}) and (\ref{1.27.5.24})
\begin{align*}
\|\ud \bA \c \ep\|_{L^2_\Sigma}\les\|\tir\bb \ud\bA\|_{L_\omega^4}\|\bb^{-\f12}\ep\|_{L_u^2 L_\omega^4}\les \log \l t\r^{\f12\M+1}\Delta_0^2\l t\r^{-2}.
\end{align*}
Using the above estimate and substituting the first estimate in (\ref{7.13.5.22}) into (\ref{1.30.2.24}), we obtain
\begin{align*}
\|\sn_\Lb \eh\|_{L^2_\Sigma}\les \l t\r^{-2}(\Delta_0^\frac{5}{4}\l t\r^\delta+\log \l t\r^{\f12\M+1}\La_0).
\end{align*}
This is the estimate for $\sn_\Lb \eh$ in (\ref{7.13.5.22}).

Using (\ref{5.13.2.24}), (\ref{1.30.2.22}), (\ref{5.21.3.23}) (\ref{3.11.3.21}), (\ref{3.6.2.21}) and (\ref{10.10.2.23}), we deduce 
\begin{equation}\label{5.13.3.24}
\|\widehat{\bR_{4A4B}}, \widehat{\bR_{3A4B}}\|_{L^2_\Sigma}\les \l t\r^{-2}(\log \l t\r^{\f12\M}\La_0+\l t\r^\delta\Delta_0^\frac{5}{4}). 
\end{equation}
Similarly in addition by using (\ref{7.13.5.22}), (\ref{1.27.5.24}) and  (\ref{1.21.2.22})-(\ref{1.27.2.22}), 
we obtain the remaining estimates (\ref{10.11.3.23}).

Next we consider the estimates $\bA_{g,2}$ in (\ref{8.8.6.22}).
Using (\ref{s2}), (\ref{lb}) and (\ref{5.13.3.24}), we have
\begin{align*}
\|\bb^{-1}\chih \tir^2\|_{L_u^2 L^2_\omega}&\les\La_0+\|\bb^{-1} \tir^2\widehat{\bR_{4A4B}}\|_{L_t^1 L_u^2 L_\omega^2}\les \l t\r^\delta(\La_0+\Delta_0^\frac{5}{4}).
\end{align*}
Thus the lower order estimate of $\bA_{g,2}$ in (\ref{8.8.6.22}) is proved. 

Substituting the above estimate of $\bA_{g,2}$ and (\ref{5.13.3.24}) we obtain 
\begin{equation*}
\|\bb^{-1}\tir^2\sn_S\chih\|_{L_u^2 L_\omega^2}\les\l t\r^\delta(\La_0+\Delta_0^\frac{5}{4})
\end{equation*}
This gives the higher order estimate for $\bA_{g,2}$ in (\ref{8.8.6.22}).

Next we prove the first estimate in (\ref{8.8.6.22}). It only remains to prove
\begin{align}
&\|\bb^{-\f12}\Big(S(\tir(\tr\chi-\frac{2}{\tir}))\Big)\|_{L^2_\Sigma}\les \log \l t\r^{\f12\M+1}(\La_0+\Delta_0^\frac{3}{2}).\label{7.10.6.22}
\end{align}
Using (\ref{7.10.7.22}) and Proposition \ref{7.15.5.22} for $\chih$, we derive by using (\ref{10.10.2.23}) 
\begin{align*}
\|L\Big(v_t^\f12\tir(\tr\chi-\frac{2}{\tir})\Big)\|_{L^2_u L_\omega^2}\les\l t\r^{-1} \log\l t\r^{\f12\M}(\La_0+\Delta_0^\frac{5}{4})
\end{align*}
which gives (\ref{7.10.6.22}), in view of $L v_t=\tr\chi v_t$. Consequently, the estimates of $\bA_b$ in (\ref{8.8.6.22}) are obtained. 

Next using (\ref{9.20.1.23}), (\ref{10.10.2.23}) and (\ref{7.13.5.22}), we obtain the estimates of $\zb, \eh, \sn\varrho$ in (\ref{10.11.2.23}).
 
Using (\ref{8.26.1.23}) with $l=0$, (\ref{10.10.2.23}), (\ref{7.13.5.22}), (\ref{8.8.6.22}), (\ref{1.27.5.24}) and (\ref{3.6.2.21}), we can obtain the first line in  (\ref{8.2.1.22}). 
Using Corollary \ref{6.29.1.23}, (\ref{7.13.5.22}), (\ref{8.21.4.21}) and (\ref{10.10.2.23}), we can obtain estimate of $[\Omega\Lb\Phi], [\Lb\Omega\Phi]$ in (\ref{9.20.4.22}). Using these estimates, we further obtain the estimates  $\Omega[\Lb\Phi]$ and $\Lb[\Omega\Phi]$  by using (\ref{5.23.1.23}), (\ref{3.6.2.21}) and (\ref{10.10.2.23}). Hence the proof of (\ref{9.20.4.22}) is complete. Using (\ref{4.5.5.24}), (\ref{10.10.2.23}), the proved estimates in (\ref{8.2.1.22}) and (\ref{9.20.4.22}), we obtain the second estimate in (\ref{8.2.1.22}). 

Using Proposition \ref{7.13.4.22},  (\ref{10.10.2.23}), (\ref{7.13.5.22}), (\ref{8.8.6.22}) and (\ref{10.11.2.23}), we can obtain (\ref{11.26.2.23}). Here to obtain the estimate of $\Lb\la$, we used Proposition \ref{2.19.4.22}.

We can obtain (\ref{7.26.2.22}) by using (\ref{7.16.2.22}), $|\sn_X v^*|\les |X^{\le 1} v|$ together with (\ref{10.10.2.23}) for $v$ and $X\Phi$, Lemma \ref{5.13.11.21} (5), (\ref{3.6.2.21}) and (\ref{11.26.2.23}).

Due to (\ref{3.6.1.22}), $\sdiv\pioh_b=O(\l t\r^{-1}v)=O(\l t\r^{-2}\log \l t\r^{\f12\M}(\Delta_0^\frac{5}{4}+\La_0))_{L^2_u L_\omega^2}$. 
  In view of (\ref{7.21.1.22}), by using this estimate, (\ref{7.21.2.22+}), Proposition \ref{8.12.1.23},  Proposition \ref{7.15.5.22}, (\ref{10.10.2.23}), (\ref{11.26.2.23}) for $\la$, and (\ref{7.13.5.22}) for $\bA_{g,1}$, we can obtain (\ref{5.7.2.24}).
\end{proof}
\section{The second order energy estimates}\label{2nd}
In this section, under the assumptions (\ref{3.12.1.21})-(\ref{6.5.1.21}), we establish in Proposition \ref{8.29.8.21} the set of second order weighted energy estimates, part of which will be further improved in Proposition \ref{3.14.4.24} when the top order energies are controlled.  
\subsection{Geometric estimates}
In this subsection, we establish more geometric estimates which will be used for deriving the second order energy estimates.

\begin{lemma}[Differentiation formulas]\label{3.23.1.23}
Let $l=0,1,2$ below.
\begin{equation*}
\begin{split}
\sn_X^l{}\rp{a}\ckk\J_L&=\sn_X^l (\sdiv{}\rp{a}\pih_b+\ud\bA{}\rp{a}\pih_b)+\sn_X^l\big(\ud\bA {}\rp{a}\pih^+_{L A}\big)\\
&+\sn_X^l\sdiv{}\rp{a}\pih^+_{AL}
\end{split}
\end{equation*}
where
\begin{align*}
\sn_X^l\sdiv{}\rp{a}\pih_b&=-\f12\sn_X^l(c^{-1}\tr\thetac{}\rp{a}v^\sharp)+\sn_X^l(c^{-1}{}\rp{a}v^*\snc\log c)\\
&=-\sn^l_X(\ckr^{-1}{}\rp{a}v^\sharp)+\sn_X^l(\bA_b{}\rp{a}v^\sharp)+\sn_X^l(\sn\varrho \c {}\rp{a}v^*)
\end{align*}
where we dropped the factor of $c^{-1}$ in the last line before being differentiated by $X^l$ since there is no difference in the treatment. 
\begin{align}\label{xdjo}
\sn_X^l {}\rp{a}\ckk\J_\Lb&=\sn_X^l\sn_A{}\rp{a}\pih_{A\Lb}+\sn_X^l(\ud \bA{}\rp{a}\pih_{A\bN}+\bA_{g,1}{}\rp{a}\pih_{AL})\nn\\
&+\sn_X^l \Big((L+\tr\chi-2k_{\bN\bN}){}\rp{a}\pih_{\Lb \Lb}\Big)\nn\\
\sn_X^l{}\rp{a}\bJ_B&=\sn_X^l{}\rp{a}\wt{\eth}_A+\sn_X^l\big((\fB+\bA_{g,1}){}\rp{a}\pih_b\big)+\sn_X^l(\ud \bA {}\rp{a}\pih_{AB})\nn\\
&+\sn_X^l\big((\tir^{-1}+\bA+\Lb \varrho)({}\rp{a}\pih_{\bN \bA}+{}\rp{a}\pih^+_{\bT A})\big)\\
\displaybreak[0]
\sn_X^l\sn_A{}\rp{a}\pih_{A\bN}&=c^{-1}(X\varrho+\sn_X)^l\circ(\sn \varrho+\sdiv)(\la^a\c \ud\bA+\sn \la^a)\nn\\
X^l\sD{}\la_a&=c^2(X\varrho+X)^l\Big({}\rp{a}\Omega \tr\thetac +\bA\la_a \c \thetac\Big)\nn\\
\sn_X^l \sn_A {}\rp{a}\pih^+_{A\bT}&=c^{-1}(X\varrho+\sn_X)^l\circ(\sn\varrho+\sdiv)\Big(-\la^a \c\ud\bA+\la_a \c \bA_{g,1}+\eta(\Omega)\Big)\nn\\
X^l Y{}\rp{a}\ss&=X^l Y\Omega \varrho+c^{-1}(X\varrho+X)^l\{(Y \varrho+Y)(([L\Phi]+\tir^{-1}+\bA_b)\la)\}\nn
\end{align}
where $ Y=L ,\, \Lb$ and 
$
{}\rp{a}\pih_{YA}={}\rp{a}\pih_b+{}\rp{a}\pih^+_{YA}, 
$
and ${}\rp{a}\wt{\eth}_A$ can be found in Proposition \ref{error_terms}.
\begin{align}\label{xdjs}
\sn_X^l\bJ[S]_\Lb&=\sn_X^l\Big(\tir(\sn\bA_{g,1}+\varpi+ \tr\chi(k_{\bN\bN}+\mho)+\ud \bA\ud \bA+\bA^2_g)\nn\\
&+2(L+\tr\chi-2k_{\bN\bN})(\tir k_{\bN\bN}+\Lb \tir)\Big)\nn\\
\sn_X^l\bJ[S]_B&=\sn_X^l\eth[S]_B+\sn_X^l\big(\tir(\chih+k_{\bN\bN}+\bA_b)\ud\bA\big)+\sn_X^l(\sn_S^{\le 1} \ud\bA)\nn\\
\displaybreak[0]
\sn_X^l \eth[S]_B&=\sn_X^l\Big(\tir(\bR_{B4BA}+k_{A\bN}\c \chi+\sn \bA_b+\sn k_{\bN\bN})\Big)\\
X^l L{}\rp{S}\ss&=X^l\{ \tir(\bA^2_{g,2}+(h, \fB, \bA_{g,1}) ([\bar\bp\Phi]+\eh)+\sD\varrho+LL\varrho+\tr\chi(\tr\chi-\frac{2}{\tir}))\}\nn
\end{align}
where the definition of ${}\rp{a}v^\sharp$ can be found in Section \ref{geocal}.
\end{lemma}

The first identity is derived by using  Proposition \ref{error_terms}, (\ref{3.6.1.22}) and (\ref{1.27.1.22}). The remaining identities follow by direct  differentiating the corresponding formulas in Proposition \ref{error_terms}.

 \begin{proposition}
  Under the assumptions (\ref{3.12.1.21})-(\ref{6.5.1.21}), we have
\begin{align}\label{1.26.2.23}
\begin{split}
\|\tir^\frac{3}{2}(\tir \sn)^2 \wt{\tr\chi}\|_{L_\omega^2}+\|(\tir \sn)^2 \wt{\tr\chi}\|_{L^2(\H_u^t)}&\les\|L \Omega^{1+\le 1}\Phi, \sn_\Omega^{\le 1}[\sn\Phi]\|_{L^2(\H_u^t)}+F_0[\Omega^{1+\le 2}\varrho](\H_u^t)\\
&+\Delta_0^\frac{3}{2}+\La_0\\
\tir^\frac{3}{2}(\tir \sn)^2 \Big(\tr\chi, \wt{\tr\chi}\Big)&=O(\l t\r^\delta\Delta_0)_{L_\omega^2}
\end{split}
\end{align}
\begin{align}
&\tir^\f12\sn_\Omega^2 \pioh_{A L}=O(\l t\r^\delta\Delta_0^\f12)_{L_\omega^2},\quad  \tir^\f12\sn_\Omega^2\sn\la=O(\l t\r^\delta\Delta_0)_{L_\omega^2}\label{2.24.5.24} 
\end{align}
\end{proposition}
\begin{proof}
We first derive by using (\ref{10.25.6.23}), Proposition \ref{6.24.10.23}  and Proposition \ref{7.15.5.22} that
\begin{align*}
&\|\tir^\frac{3}{2}(\tir\sn)^2 \wt{\tr\chi}\|_{L_\omega^2}+\|\tir^3 \sn^2 \wt{\tr\chi}\|_{L_t^2 L_\omega^2(\H_u^t)}\\
&\les \|\sn^2\wt{\tr\chi}\|_{L^2_\omega(S_{0,u})}+\|\tir^4 \sn^2 \sD\varrho, \tir^3 \sn^2[L\Phi]\|_{L_t^2 L_\omega^2(\H_u^t)}+\|\tir^4 \sn(\sn\bA\c \bA)\|_{L^2_t L_\omega^2}\\
&+\|\tir^4 \sn^2(\N(\Phi, \bp\Phi)+\bA_{g,1}^2), \tir^4 (\sn \bA)^2 \|_{L_t^2 L_\omega^2}\\
&\les \|\sn^2\wt{\tr\chi}\|_{L^2_\omega(S_{0,u})}+\|\tir^4 \sn^2 \sD\varrho, \tir^3 \sn^2[L\Phi]\|_{L_t^2 L_\omega^2(\H_u^t)}+\Delta_0^\frac{3}{2}.
\end{align*} 
The first estimate in
(\ref{1.26.2.23})  then follows by using Proposition \ref{12.21.1.21}, Lemma \ref{10.10.3.23} and Lemma \ref{3.17.2.22}. 

Using Lemma \ref{3.17.2.22}, Proposition \ref{7.15.5.22} and (\ref{3.12.1.21}), we also have
\begin{equation*}
\|\tir^4 \sn^2 \sD\varrho, \tir^3 \sn^2[L\Phi]\|_{L_t^2L_\omega^2}\les \Delta_0+F_0[\Omega^{1+\le 2}\varrho](\H_u^t)\les \l t\r^\delta\Delta_0
\end{equation*}
The second estimate follows as a consequence. 

Note by using (\ref{1.7.2.23}), we bound 
\begin{align*}
\|\sn_\Omega^2 \sn \la\|_{L_\omega^2}\les \|\tir^3 \sn^2\bA_b\|_{L_\omega^2}+\|\Omega^2\log c, \Omega(\bA \la)\|_{L_\omega^2} +\l t\r^{-\frac{11}{4}+3\delta}\Delta_0^\frac{5}{2}.
\end{align*}
Applying  Proposition \ref{7.15.5.22} and Proposition \ref{10.16.1.22}, we can obtain
\begin{equation*}
\|\tir^\f12\sn_\Omega^2 \sn\la\|_{L_\omega^2}\les \|\tir^\frac{7}{2} \sn^2 \bA_b\|_{L_\omega^2}+\l t\r^{-\f12+\delta}\Delta_0.
\end{equation*}
Using the second estimate in (\ref{1.26.2.23}), we then can obtain the second estimate in (\ref{2.24.5.24}). Then in view of (\ref{7.16.2.22}), using (\ref{3.28.3.24}) and Proposition \ref{7.15.5.22},  we can derive the first estimate of (\ref{2.24.5.24}). 
\end{proof}
\begin{proposition}\label{djoest}
With $X=\Omega, S$,  under the assumptions (\ref{3.12.1.21})-(\ref{6.5.1.21}), we have
\begin{align}
&X\bJ[{}\rp{a}\Omega]_L=O(\tir^{-1}X^{\le 1}v)+c^2X\Omega \tr\thetac+X\sdiv(\eta({}\rp{a}\Omega))\label{7.11.5.21}\\
\displaybreak[0]
&\qquad\qquad\quad+XL\Omega \varrho+ O\Big(\l t\r^{-2}(\log \l t\r^{\f12\M+1}\La_0+ \l t\r^{2\delta}\Delta_0^\frac{5}{4})\Big)_{L^2_u L_\omega^2}\nn\\
&\qquad\qquad\quad=O(\l t\r^{-\frac{3}{2}+\delta-\f12\vs(X)(1-2\delta)}\Delta_0)_{L_u^2 L_\omega^2}\nn\\
&\sn_X{}\rp{a}\bJ_\Lb= O(\l t\r^{-1}v) +O(\l t\r^{\delta-1}\Delta_0)_{L^2_u L_\omega^2}\label{9.5.16.23}\\
&\sn_X {}\rp{a}\bJ_B=\fB\c O(1)+O(\l t\r^{-1+\delta}\Delta_0)_{L^2_u L_\omega^2}\label{2.24.1.22}\\
&\sn_X\bJ[S]_\Lb=\vs(X)O(\tir^{-1})+O(\l t\r^{-1+\delta}\Delta_0)_{L^2_u L_\omega^2}\label{10.19.1.23}\\
&\|\bb^{-\f12}\sn_X\bJ[S]_A\|_{L^2_\Sigma}\les \l t\r^{\delta}\Delta_0, \label{7.12.1.21}\\
&\|\bb^{-1} X L {}\rp{S}\ss\|_{L^2_\Sigma}\les\l t\r^{-1}\left(W_2[\Omega^2 \Phi]^\f12(t)+W_2[XS\Phi]^\f12(t)\right)\label{7.9.6.22}\\
& \qquad\qquad\qquad+(1-\vs(X))\|\bb^{-1}X\bA_b\|_{L^2_\Sigma}+\l t\r^{-1}\log \l t\r^{\f12\M+1}(\La_0+\Delta_0^\frac{5}{4})\nn
\end{align}
\end{proposition}
   Both $X\bJ[{}\rp{a}\Omega]_L$ and $XL{}\rp{S}\ss$ will be paired to $\Lb \Phi$ in energy estimates, which need more careful treatments, since $\Lb\Phi$ does not have smallness.
   
    In order to prove Proposition \ref{djoest}, we further give results on derivatives of $\pioh$ and $\piSh$. We will frequently use (\ref{1.25.2.22}) and (\ref{4.22.4.22}) in the sequel without mentioning explicitly. 
\begin{lemma}
With $X\in \{\Omega, S\}$,  under the assumptions (\ref{3.12.1.21})-(\ref{6.5.1.21}), we have
\begin{align}
&\sn_X\sn_\bT{}\rp{a}\pih_b=\fB O(1)+ \Big((1-\vs(X))+\vs(X) \l t\r^{-1}\log \l t\r^{\f12\M}\Big)O(\log \l t\r(\Delta_0^\frac{5}{4}+\La_0))_{L^2_\Sigma}\label{3.11.1.22}
\end{align}
\end{lemma}
\begin{proof}
 We first prove with $X\in \{\Omega, S\}$, 
\begin{align}
\sn_X \sn_L {}\rp{a}\pih_b&= O(\l t\r^{-1}\log \l t\r^{\f12\M}(\Delta_0^\frac{5}{4}+\La_0))_{L^2_\Sigma}\label{9.3.2.23}\\
\sn_X\sn_\Lb{}\rp{a}\pih_b
&=((1-\vs(X))+\vs(X) \l t\r^{-1}\log \l t\r^{\f12\M})O(\log \l t\r(\Delta_0^\frac{5}{4}+\La_0))_{L^2_\Sigma}\nn\\
&+ O(1)\fB\label{9.3.3.23}
\end{align}
By the convention given in Section \ref{geocal}, for $Y=L, \Lb$, we have
\begin{equation}\label{9.3.4.23}
X Yv\sta{a}\wedge\Pi=O(1)XYv,\quad Yv\sta{a}\wedge\hN=[\sn v]. 
\end{equation} 
Using the above identity and (\ref{7.19.6.22}), (\ref{9.3.2.23}) follows by using (\ref{3.28.3.24}), Proposition \ref{7.15.5.22}, Proposition \ref{1steng} and Lemma \ref{5.13.11.21} (1) and (5). If $Y=\Lb$ in (\ref{7.19.6.22}) and (\ref{9.3.4.23}), applying (\ref{8.2.1.22}), (\ref{9.20.4.22}) and (\ref{9.16.3.23}), and treating the lower order errors similar to (\ref{9.3.2.23}), we can obtain (\ref{9.3.3.23}).
Finally, we obtain (\ref{3.11.1.22}) by using  $2\bT=L+\Lb$ and combining (\ref{9.3.2.23})  with (\ref{9.3.3.23}). 
\end{proof}
  
 Using Proposition \ref{7.15.5.22} and Proposition \ref{10.16.1.22}, we have the following result
 \begin{lemma}
 Under the assumptions (\ref{3.12.1.21})-(\ref{6.5.1.21}), there hold for $X\in\{\Omega, S\}$ that
\begin{align}\label{7.19.2.22}
\begin{split}
&\|\bb^{-\f12}\sn_X\sn_L\sn \la\|_{L^2_\Sigma}\les \l t\r^{\delta-1}\Delta_0,\\
&X \Lb(\tr\eta\c\Omega)=O(1)[\Lb\Phi]+\vs(X)O(\l t\r^{\delta-1}\Delta_0)_{L^2_\Sigma}+(1-\vs(X))O(\l t\r^{\delta}\Delta_0)_{L^2_\Sigma}\\
&\|\sn_X\sn_Y\eh(\Omega)\|_{L^2_\Sigma}\les \l t\r^{-1+\delta+\max(-\vs(Y),0)}\Delta_0, Y=L, \Lb, e_A\\
&\|\sn_X\sn_Y(\tr\eta\c\Omega)\|_{L^2_\Sigma}\les \l t\r^{-1+\delta}\Delta_0,\, Y=L,  e_A\\
&\|\sn_X(\sn_L,\sn) (\la \c \ud\bA)\|_{L_u^2 L_\omega^2}\les \l t\r^{-2+2\delta}\Delta_0\\
&\|\bb^{-\f12}\sn_X\sn_Y (\la \c \bA_{g,1})\|_{L^2_\Sigma}\les \l t\r^{-2+\max(-\vs(Y), 0)+2\delta}\Delta_0^2 \\
&\|\bb^{-\f12}\sn_X\sn_Y (\la \c \bA_g)\|_{L^2_\Sigma}\les \l t\r^{-\frac{3}{2}+(\frac{1}{2}+\delta)\max(-\vs(Y), 0)+2\delta}\Delta_0^2,\, Y=L, \Lb, e_A.
\end{split}
\end{align}
\end{lemma}
\begin{proof} 
 The first line in (\ref{7.19.2.22}) is (\ref{2.19.5.22}). The second estimate follows by using (\ref{6.7.4.23}) and (\ref{6.22.1.21}). For the third estimate, we derive by (\ref{4.22.4.22}), (\ref{1.25.2.22}) and Proposition \ref{7.15.5.22}
\begin{align*}
\sn_X\sn_\Lb(\eh\c \Omega)&=\sn_X \sn_\Lb \eh \c O(\tir)+ O(\l t\r^{\delta-\vs(X)(1-\delta)}\Delta_0^\f12)_{L_\omega^4}\eh+O(\tir)\sn_\Lb \eh+O(1)\sn_X^{\le 1}\eh\\
&=O(\l t\r^{\delta}\Delta_0)_{L^2_\Sigma}
\end{align*}
and 
\begin{equation*}
\sn_X(\sn_L, \sn)(\eh\c \Omega)=O(\l t\r^{-1+\delta}\Delta_0)_{L^2_\Sigma}.
\end{equation*}
The estimate of $\sn_X(\sn_L, \sn)(\tr\eta\c \Omega)$ follows similarly by noting $\tr\eta=[L\Phi]$. 

 The estimate of $\sn_X^l(\la\c\ud\bA)$ has been given in (\ref{dpio}).  Similar to the last line in (\ref{dpio}), the last two estimates can be obtained by using  Proposition \ref{7.22.2.22}, Proposition \ref{7.15.5.22} and Proposition \ref{10.16.1.22}.
\end{proof}

\begin{proof}[Proof of Proposition \ref{djoest}]
We first claim  
\begin{align}\label{9.4.7.23}
\begin{split}
&\|\bb^{-\f12}\{\sn_X\left(\ud\bA (\pioh^+_{LA},\pioh_b)\right), \sn_X (\ud \bA \pioh_{A\bN})\}\|_{L^2_\Sigma}\les \l t\r^{-1+2\delta}\Delta_0^\frac{3}{2}\\
&\|\sn_X\left(\bA ({}\rp{a}v^\sharp, \rp{a}v^*)\right)\|_{L^2_\Sigma}\les \l t\r^{-2+2\delta}\Delta_0^\frac{3}{2}\\
&\|\tir^{-1+\frac{\vs(X)}{4}}\sn_X\left(\ud\bA (\pioh^+_{LA},\pioh_b)\right), \Delta_0^\f12\sn_X\left(\bA ({}\rp{a}v^\sharp, \rp{a}v^*)\right)\|_{L^2_\omega}\les \l t\r^{-\frac{11}{4}+2\delta}\Delta_0^\frac{3}{2}\\
 &\|\bb^{-\f12}\{\sn_X(\bA_{g,1}\pioh^+_{AL}), \l t\r^{-\delta}\sn_X(\bA_{g,1}\pioh_{A\bN})\}\|_{L^2_\Sigma}\les \l t\r^{-2+2\delta}\Delta_0^2.
 \end{split}
\end{align}

Indeed for $F={}\rp{a}v^\sharp,\pioh_b$, using Lemma \ref{5.13.11.21} (1) and (\ref{3.28.3.24}) 
\begin{align*}
|v|\les \l t \r^{-1+\delta}\Delta_0^\f12, |\sn_X F|\les |X^{\le 1} v|. 
\end{align*}
Using the above estimate, (\ref{10.10.2.23}) and Proposition \ref{7.15.5.22}, we obtain
\begin{align*}
\|\sn_X(\bA\c F)\|_{L^2_\Sigma}&\les \|\sn_X\bA\|_{L^2_\Sigma}\|F\|_{L_x^\infty}+\|\bA\|_{L_\omega^4}\|\bb^\f12\tir\sn_X F\|_{L_u^2 L^4_\omega}\\
&\les \l t\r^{-2+2\delta}\Delta_0^\frac{3}{2}
\end{align*}
as desired in (\ref{9.4.7.23}). The third estimate follows similarly by also using  (\ref{5.21.1.21}). The last estimate follows by using (\ref{5.21.1.21}), (\ref{3.25.1.22}), (\ref{3.6.2.21}) and (\ref{3.11.3.21}). The first estimate follows similarly by using (\ref{3.6.2.21}) and Lemma \ref{5.13.11.21} (5) for the estimate of $\ud \bA$ and $\sn_X \ud\bA$ 

We first consider (\ref{7.11.5.21}) and (\ref{9.5.16.23}).
In view of (\ref{xdjo}), applying (\ref{9.4.7.23}) leads to  
\begin{align*}
\sn_X\sdiv{}\rp{a}\pih_b&=-\sn_X(\ckr^{-1}{}\rp{a}v^\sharp)+\sn_X(\bA_b{}\rp{a}v^\sharp)+\sn_X(\bA_{g,1} \c \pioh_b)\\
&=-\sn_X(\ckr^{-1}{}\rp{a}v^\sharp)+\left\{\begin{array}{lll}
O(\l t\r^{-2+2\delta}\Delta_0^\frac{3}{2})_{L^2_\Sigma}\\
O(\l t\r^{-\frac{11}{4}+2\delta}\Delta_0)_{L_\omega^2}
\end{array}\right.
\end{align*}

Using (\ref{7.25.1.22}) and (\ref{9.4.7.23}),  with $Y=L$ or $\Lb$, we summarize the above estimates to obtain
\begin{equation}\label{7.20.7.22}
\begin{split}
\sn_X{}\rp{a}\ckk\J_Y&= -\sn_X(\ckr^{-1}{}\rp{a}v^\sharp)+(\sn_X{}\rp{a}\ckk\J_Y)_g+O(\l t\r^{-2+2\delta}\Delta_0^\frac{3}{2})_{L^2_u L_\omega^2}\\
\sn_X{}\rp{a}\ckk \J_L&=-\sn_X(\ckr^{-1}{}\rp{a}v^\sharp)+(\sn_X{}\rp{a}\ckk\J_L)_g+O(\l t\r^{\frac{-7-\vs(X)}{4}+2\delta}\Delta_0)_{L_\omega^2}
\end{split}
\end{equation}
with
\begin{align*}
(\sn_X{}\rp{a}\ckk\J_L)_g&=X \sn_A {}\rp{a}\pih^+_{AL}\\
 (\sn_X{}\rp{a}\ckk\J_\Lb)_g&=(\sn_X(\sdiv\sn\la+\sdiv\eta(\Omega))+\sn_X\left(\ud\bA{}\rp{a}\pih_{A\bN}+\sn(\ud \bA\c \la)+\bA_{g,1}{}\rp{a}\pih_{AL}\right)\\
 &+\sn_X \big((L+\tr\chi+k_{\bN\bN}){}\rp{a}\pih_{\Lb \Lb}\big).
\end{align*}
where we used  (\ref{7.16.2.22}) and Proposition \ref{error_terms} to derive the second symbolic formula.

We will prove
\begin{align}
&(\sn_X{}\rp{a}\ckk\J_L)_g-\{c^2X(\Omega \tr\thetac)+X\sdiv\eta({}\rp{a}\Omega)\nn\}\\
&\qquad\qquad\qquad\qquad=O(\l t\r^{-3+2\delta}\Delta_0^\frac{3}{2})_{L^2_u L_\omega^2}, O(\l t\r^{-\frac{11}{4}+2\delta}\Delta_0^2)_{L_\omega^2}\label{7.22.6.22}\\
&\sn_X \big((L+\tr\chi+k_{\bN\bN}){}\rp{a}\pih_{\Lb \Lb}\big)=O(\l t\r^{\delta-1}\Delta_0)_{L^2_u L_\omega^2}\label{9.4.2.23}\\
&\sn_X{}\rp{a}\ckk\J_\Lb=X(\tir^{-1}v^\sharp)+c^2X(\Omega \tr\thetac)+X\sdiv\eta({}\rp{a}\Omega)+O(\l t\r^{\delta-1}\Delta_0)_{L^2_u L_\omega^2}\nn\\
&\qquad\qquad= O(\l t\r^{-1}(X^{\le 1}v)) +O(\l t\r^{-1+\delta}\Delta_0)_{L_u^2 L_\omega^2}\label{7.22.8.22}\\
&XL{}\rp{a}\ss-X(\tir^{-1}+L) \Omega \varrho-\tir^{-2}X^{\le 1}\la=O(\l t\r^{-3+2\delta}\Delta_0^\frac{3}{2})_{L^2_u L^2_\omega}, O(\l t\r^{-\frac{11}{4}+2\delta}\Delta_0^\frac{3}{2})_{L^\infty_u L^2_\omega}\label{7.22.7.22}
\end{align}
We first prove (\ref{7.22.6.22}). Using (\ref{xdjo}) and the definition of $(\sn_X{}\rp{a}\ckk\J_L)_g$ in (\ref{7.20.7.22}), we write symbolically
\begin{align*}
(\sn_X{}\rp{a}\ckk\J_L)_g&=\sn_X\Big(\sD \la+\sdiv\eta({}\rp{a}\Omega)+\sn(\la \c \bA_{g,1})\Big)+\sn_X\big(\sn\varrho(\la \bA_{g,1}+\sn\la+\eta(\Omega))\big)\\
&+X\varrho(\sn \varrho+\sdiv)(\sn\la+\eta(\Omega)+\la \bA_{g,1})\\
&=\sn_X\Big(\sD \la+\sdiv\eta({}\rp{a}\Omega)\Big)+\left\{\begin{array}{lll}
O(\l t\r^{-3+2\delta}\Delta_0^\frac{3}{2})_{L^2_u L_\omega^2}\\
O(\l t\r^{-\frac{11}{4}+2\delta}\Delta_0^\frac{3}{2})_{L^\infty_u L_\omega^2}\end{array}\right.
\end{align*}
where we used Proposition \ref{7.15.5.22}, (\ref{dpio}) and Proposition \ref{10.16.1.22} to obtain the last estimate.
 
Note from Proposition \ref{error_terms}
\begin{align*}
X\sD \la^a&=c^2(X\log c\c+X)(\Omega\tr\thetac+\bA\la \thetac)
\end{align*}

Next we show
\begin{align}\label{7.22.5.22}
\begin{split}
&\|X(\bA\la \thetac)\|_{L^2_\Sigma}\les\l t\r^{-3+2\delta}\Delta_0^\frac{3}{2}.
\end{split}
\end{align}
Indeed we write symbolically $\thetac=\bA_g+\bA_b+\tir^{-1}$ due to  (\ref{5.17.1.21}). Therefore, using (\ref{dpio}) and Proposition \ref{7.15.5.22}
\begin{align*}
\|X(\bA\la \thetac)\|_{L^2_u L_\omega^2}&\les\|X^{\le 1}(\bA\la)\|_{L^2_u L_\omega^2}\l t\r^{-1}+\| X\bA\|_{L_u^2 L^4_\omega}\|\la\bA\|_{L_u^\infty L_\omega^4}\\
&\les \l t\r^{-3+2\delta} \Delta_0^\frac{3}{2},
\end{align*}
as desired. Thus by (\ref{7.22.5.22}), Proposition \ref{7.15.5.22} and $\la=O(\l t\r^\delta\Delta_0)$ we have
\begin{equation*}
X\sD\la\rp{a}=c^2X\Omega\tr\thetac+ O(\l t\r^{-3+2\delta}\Delta_0^\frac{3}{2})_{L^2_u L_\omega^2}.
\end{equation*}
Similarly, we have $X\sD\la\rp{a}=c^2X\Omega\tr\thetac+O(\l t\r^{-\frac{11}{4}+2\delta}\Delta_0^\frac{3}{2})_{L_\omega^2}$. 
Hence we conclude (\ref{7.22.6.22}).
 
Next we prove (\ref{9.4.2.23}). Using Lemma \ref{5.13.11.21} (2) and (5),  (\ref{8.23.2.23}) and (\ref{8.24.4.23}), we infer
\begin{align}
&\sn_X\big((\tr\chi+k_{\bN\bN}){}\rp{a}\pih_{\Lb\Lb})=O(\l t\r^{-1+\delta}\Delta_0)_{L^2_u L_\omega^2}.\label{9.4.6.23}
\end{align}
 Combining (\ref{9.4.6.23}) in the above with (\ref{3.25.1.22}) yields (\ref{9.4.2.23}). 
 
Using Proposition \ref{7.15.5.22} and (\ref{4.22.4.22}), we have
 $$\|\bb^{-\f12}X \Omega \tr\thetac\|_{L^2_\Sigma}\les\Delta_0 \l t\r^{\delta-\f12-\f12\vs(X)}.$$
 Moreover using  (\ref{9.4.7.23}), (\ref{dpio}), Lemma \ref{5.13.11.21} (5)  and (\ref{3.25.1.22}), we have the following error estimates
 \begin{align*}
 \sn_X(\ud\bA{}\rp{a}\pih_{A\bN}),\sn_X\sn(\ud \bA\la), \sn_X(\bA_{g,1}{}\rp{a}\pih_{AL})&=O(\l t\r^{-2+2\delta}\Delta_0^\frac{3}{2})_{L^2_u L_\omega^2}
 \end{align*}
   Then using (\ref{8.2.2.23}), (\ref{10.10.2.23}), (\ref{7.20.7.22})-(\ref{9.4.2.23}) and the above two lines of estimates we obtain the following estimate
\begin{align}\label{9.4.4.23}
&\|(\sn_X{}\rp{a}\ckk\J_Y)_g\|_{L^2_u L_\omega^2}\les \Delta_0\left(\l t\r^{\delta-\frac{3}{2}-\f12\vs(X)(1-2\delta)}+\max(-\vs(Y),0) \l t\r^{-1+\delta}\right), Y=L, \Lb.
\end{align}
Using (\ref{3.28.3.24}) we obtain 
\begin{align}
\sn_X(\ckr^{-1}{}\rp{a}v^\sharp)&=O(\tir^{-1})(|v|+|Xv|)\label{9.4.3.23}.
\end{align}
Combining (\ref{7.20.7.22}) with (\ref{9.4.4.23}) and (\ref{9.4.3.23}), we obtain the second line in (\ref{7.22.8.22}). Similarly using the above error estimates, (\ref{9.4.2.23}) and the estimate of $X\sD\la^a$ gives the first line in (\ref{7.22.8.22}).

Next we prove 
\begin{equation}
\|\bb^{-\f12}X \Lb{}\rp{a}\ss\|_{L^2_\Sigma}\les \l t\r^{\delta}\Delta_0\label{7.10.11.22}.
\end{equation}
 For convenience, by neglecting  multi-linear terms which have better decay properties, we recast the formula from (\ref{xdjo})
\begin{equation*}
X\Lb {}\rp{a}\ss=X\Lb \Omega \varrho+X\Lb\big(([L\Phi]+\tir^{-1}+\bA_b)\la\rp{a}\big).
\end{equation*}
Using (\ref{8.31.4.23}), (\ref{6.7.4.23}), Proposition \ref{7.22.2.22} and Proposition \ref{10.16.1.22}, we estimate
\begin{equation*}
\|\bb^{-\f12}X\Lb\big(([L\Phi]+\bA_b)\la\big)\|_{L^2_\Sigma}\les \l t\r^{-1+2\delta}\log \l t\r\Delta_0,\, \|X\Lb(\tir^{-1}\la)\|_{L_u^2 L_\omega^2}\les \l t\r^{-1+\delta}\Delta_0,
\end{equation*}
where we also applied (\ref{7.16.3.22}) and (\ref{4.9.1.24}) to derive the second estimate. Using (\ref{LbBA2}) and Lemma \ref{3.17.2.22}, $\|X\Lb\Omega\varrho\|_{L^2_\Sigma}\les \l t\r^{\delta}\Delta_0$. Combining the three estimates gives (\ref{7.10.11.22}).  In view of  (\ref{7.10.11.22}), (\ref{7.22.8.22}) and (\ref{10.10.2.23}), we conclude (\ref{9.5.16.23}).

Next we prove (\ref{7.22.7.22}).   By (\ref{xdjo}), we first schematically write
$$
XL{}\rp{a}\ss=XL\Omega \varrho+XL\big(([L\Phi]+\bA_b+\tir^{-1})\la\big),$$
since other terms are of lower order and have better decay properties. 

Recall $[L\Phi]+\bA_b=\bA$ and recall from the last line in (\ref{dpio}) that
\begin{align*}
&\|\bb^{-\frac{1}{2}}\sn_X\sn_L (\la \c \bA)\|_{L^2_\Sigma}\les \l t\r^{-2+2\delta}\Delta_0^\frac{3}{2}.
\end{align*}
Using Proposition \ref{7.15.5.22} and Proposition \ref{10.16.1.22}, we derive
\begin{equation*}
\|\sn_X \sn_L(\la \c \bA)\|_{L_\omega^2}\les \l t\r^{-\frac{11}{4}+2\delta}\Delta_0^\frac{3}{2}.
\end{equation*}
Note that, due to (\ref{3.22.5.21}), $L\la=c \Omega \log c$. Combining this fact with the above estimates, we obtain (\ref{7.22.7.22}).

We have from Proposition \ref{1steng} 
\begin{equation*}
\|\tir^{-1}X^{\le 1}\Omega \varrho, (\log \l t\r)^{-1}\tir^{-2}X^{\le 1}\la\|_{L^2_u L_\omega^2}\les \l t\r^{-2}\log \l t\r^{\f12\M}(\Delta_0^\frac{5}{4}+\La_0).
\end{equation*}
Substituting the above estimate to (\ref{7.22.7.22}) leads to 
$$
XL{}\rp{a}\ss=XL\Omega\varrho+O( \l t\r^{-2}\log \l t\r^{\f12\M+1}(\Delta_0^\frac{5}{4}+\La_0))_{L_u^2 L_\omega^2}.
$$
Combining the above estimate with (\ref{7.20.7.22}), (\ref{7.22.6.22}), (\ref{9.4.4.23}), (\ref{9.4.3.23}) and (\ref{10.10.2.23}) gives (\ref{7.11.5.21}). 

Next we prove (\ref{2.24.1.22}). In the formula for $\sn_X{}\rp{a}\bJ_B$ in (\ref{xdjo}), we note that the terms in $\sn_X{}\rp{a}\wt\eth_B$ can be treated by using (\ref{7.19.2.22}), (\ref{3.11.1.22}), Proposition \ref{7.15.5.22}, Lemma \ref{5.13.11.21} (5) and Proposition \ref{10.16.1.22}
\begin{align}
\sn_X {}\rp{a}\wt\eth_B&=\sn_X(\sdiv(c^{-1}\la \bA_{g,2})+\sn\Omega\log c)+\sn_X \sn_\bT (\pioh_b+c^{-1}\la\bA_{g,1}+c^{-1}\eta(\Omega))\nn\\
&+\sn_X\big(\bA\sn\la+c^{-1}\tir^{-1}\sn \la+c^{-1}\ud \bA^2\c \Omega+\sn_L(c^{-1}\la \ud \bA)\big)\nn\\
&=O(\l t\r^{-1+\delta}\Delta_0)_{L^2_u L_\omega^2}+O(\fB).\label{2.23.5.22}
\end{align}

Moreover using (\ref{5.21.1.21}), (\ref{3.25.1.22}), (\ref{3.11.3.21}), Lemma \ref{5.13.11.21} (5) and (\ref{8.24.5.23}), we obtain
\begin{equation*}
\|\sn_X\big((\fB+\bA+\tir^{-1})(\pioh_{\bN A}, \pioh^+_{\bT A})\big), \sn_X(\ud \bA\pioh_{AB})\|_{L^2_u L_\omega^2}\les \l t\r^{\delta-1}\Delta_0.
\end{equation*}
Applying (\ref{8.24.5.23}) to $F=\pioh_b$ and using (\ref{3.28.3.24}) and Lemma \ref{5.13.11.21} (1) yields $\sn_X(\fB \pioh_b)=O(\l t\r^{-2+\delta}\Delta_0)_{L_u^2 L_\omega^2}$.
Combining (\ref{9.4.7.23}) and the above two estimates, we conclude
\begin{equation*}
\sn_X {}\rp{a}\bJ_B=O(\fB)+O(\l t\r^{\delta-1}\Delta_0)_{L^2_u L_\omega^2}
\end{equation*}
as desired in (\ref{2.24.1.22}).

Next consider $\sn_X \bJ[S]_\Lb$. By using (\ref{L2BA2}), (\ref{8.3.4.23}), Lemma \ref{5.13.11.21} (5) and Sobolev embedding on spheres we bound
\begin{equation*}
\|\sn_X(\tir \sn\bA_{g,1}+\tir\ud\bA^2+\tir \bA_g^2)\|_{L^2_u L_\omega^2}\les \l t\r^{-1+\delta}\Delta_0.
\end{equation*}
In view of (\ref{6.23.4.21}), we write
\begin{align*}
\sn_X((L+\tr\chi-2k_{\bN\bN}) \piSh_{\Lb\Lb})=\sn_X\big((L+\tr\chi-2k_{\bN\bN})(\tir k_{\bN\bN}+\Lb\tir)\big).
\end{align*}
We use (\ref{6.7.4.23}) to have
\begin{equation}
XL(\tir k_{\bN\bN})=\vs(X) [\Lb \Phi]+O(\l t\r^\delta\Delta_0)_{L^2_\Sigma}.\label{9.8.11.23}
\end{equation}
Hence using (\ref{9.26.1.23}) and (\ref{9.8.11.23}), we derive 
\begin{align*}
X\big((L&+\tr\chi-2k_{\bN\bN})(\tir k_{\bN\bN}+\Lb \tir)\big)\\
&=X(L\Lb \tir+L(\tir k_{\bN \bN}))+X\left((\tr\chi-2k_{\bN\bN})(\tir k_{\bN\bN}+\Lb \tir)\right)\\
&=\vs(X)\fB+O(\l t\r^\delta\Delta_0)_{L^2_\Sigma}+X((\tr\chi-2k_{\bN\bN})(\tir k_{\bN\bN}+\Lb \tir)).
\end{align*}
For the last term in the above, using (\ref{L2conndrv}), (\ref{8.23.1.23}), Proposition \ref{7.22.2.22} for estimate of $\mho$, Lemma \ref{5.13.11.21} and the first estimate in (\ref{9.12.3.22}), we infer
\begin{align*}
X\big((\bA+\tir^{-1}+\fB)(\tir\fB+1+\tir\mho)\big)&=O(1)(X(\tir^{-1})+X(\bA)+\vs(X)(\fB+\mho)+X\fB+X\mho)\\
&=\vs(X)(\tir^{-1}+\mho+\fB)+O(\l t\r^{-1+\delta}\Delta_0)_{L^2_u L_\omega^2}.
\end{align*}
Note that the estimate of $\sn_X\big(\tir (\varpi+\tr\chi (\mho+\fB))\big)$ can be included in the above.
We combine the above two estimates to obtain
\begin{align*}
X\big((L&+\tr\chi-2k_{\bN\bN})(\tir k_{\bN\bN}+\Lb \tir)\big)+\sn_X\big(\tir (\varpi+\tr\chi (\mho+\fB))\big)\\
&=\vs(X)(\tir^{-1}+\mho+\fB)+O(\l t\r^{-1+\delta}\Delta_0)_{L^2_u L_\omega^2}.
\end{align*}
Hence we conclude in view of (\ref{xdjs})
\begin{align}
&\sn_X  \bJ[S]_\Lb=\vs(X)(\tir^{-1}+\mho+\fB)+O(\l t\r^{\delta-1}\Delta_0)_{L^2_u L_\omega^2}\label{7.24.2.22}.
\end{align}
We thus obtained (\ref{10.19.1.23}) by noting $\mho, \bb \fB=O(\tir^{-1})$.

Next we control $\sn_X \eth[S]_B$.
In view of (\ref{1.21.2.22}) and (\ref{cmu_2}), we symbolically write
\begin{equation*}
 \bR_{B4BA}=(\sn_L, \sn) \bA_{g,1}+(\chi+\bp\Phi)\bA_{g,1}.
\end{equation*}
Hence we can recast $\sn_X \eth[S]_B$ as
\begin{equation*}
\sn_X \eth[S]_B=\sn_X\Big(\tir (\sn_L, \sn) \bA_{g,1}+\tir(\chi+\bp\Phi)\bA_{g,1}+\tir (\sn \bA_b+\sn k_{\bN\bN})\Big).
\end{equation*}
Using (\ref{8.23.2.23}), (\ref{8.3.4.23}) and (\ref{7.13.5.22})  we derive
\begin{align*}
\sn_X(\tir(\chi+\bp\Phi)\bA_{g,1})
&=O(1)\sn_X^{\le 1}\bA_{g,1}+O(\l t\r^{-1+\delta}\Delta_0^\frac{3}{2})_{L^2_\Sigma}=O(\l t\r^{-1+\delta}\Delta_0)_{L^2_\Sigma}
\end{align*}
and recall from Proposition \ref{7.15.5.22} that
\begin{equation*}
\|\sn_X(\tir (\sn_L, \sn) \bA_{g,1}, \tir \sn \bA_b, \tir \sn k_{\bN\bN})\|_{L^2_\Sigma}\les \l t\r^{\delta}\Delta_0.
\end{equation*}
We conclude
\begin{equation*}
\|\sn_X \eth[S]_B\|_{L^2_\Sigma}\les \l t\r^{\delta}\Delta_0.
\end{equation*}
Next, recall from (\ref{xdjs}) that
\begin{align*}
\sn_X \bJ[S]_B
&=\sn_X^l\big(\tir(\chih+k_{\bN\bN}+\bA_b)\ud\bA\big)+\sn_X^l(\sn_S^{\le 1} \ud\bA)+\sn_X\eth[S]_B.
\end{align*}
 By using (\ref{8.23.2.23}), Lemma \ref{5.13.11.21} (5) and Proposition \ref{7.15.5.22}, we derive
\begin{align*}
\sn_X^l&\big(\tir(\chih+k_{\bN\bN}+\bA_b)\ud\bA\big)+\sn_X^l(\sn_S^{\le 1} \ud\bA)=O(\l t\r^{-1+\delta}\Delta_0)_{L^2_u L_\omega^2}, l=0,1.
\end{align*}
We thus conclude
\begin{equation*}
\|\bb^{-\f12}\sn_X  \bJ[S]_B\|_{L^2_\Sigma}\les \l t\r^{\delta}\Delta_0,
\end{equation*}
as desired in (\ref{7.12.1.21}).


By using (\ref{xdjs}), Lemma \ref{10.10.3.23}, (\ref{8.23.2.23}), (\ref{8.3.4.23}), Lemma \ref{comp} and Proposition \ref{1steng}, we bound
\begin{align*}
&\|\bb^{-1}X L{}\rp{S}\ss\|_{L^2_\Sigma}\\
&\les \|\bb^{-1}X^{\le 1}([\bar \bp\Phi], \eh)\|_{L^2_\Sigma}+\|\bb^{-1}X^{\le 1} (\tir \sD \varrho, \tir LL\varrho)\|_{L^2_\Sigma}+\|\bb^{-1}X^{\le 1}\big((\tr\chi-\frac{2}{\tir}), \tir\bA_{g,2}^2\big)\|_{L^2_\Sigma}\\
&+\|\bb^{-1}X(\tir \fB \bAn)\|_{L^2_\Sigma}+\l t\r^{-\frac{7}{4}+2\delta}\Delta_0^\frac{3}{2}\\
&\les  \l t\r^{-1}(W_2[\Omega^2 \Phi]^\f12(t)+W_2[XS\Phi]^\f12(t))+\|\bb^{-1}X^{\le 1}(\tr\chi-\frac{2}{\tir})\|_{L^2_\Sigma}+\l t\r^{-1}\log \l t\r^{\f12\M+1}(\La_0+\Delta_0^\frac{5}{4})\\
&\les \l t\r^{-1}\log \l t\r^{\f12\M+1}(\La_0+\Delta_0^\frac{5}{4})+\l t\r^{-1}\left(W_2[\Omega^2 \Phi]^\f12(t)+W_2[XS\Phi]^\f12(t)\right)\\
&+(1-\vs(X))\|\bb^{-1}X \bA_b\|_{L^2_\Sigma}.
\end{align*}
This gives (\ref{7.9.6.22}).
\end{proof}

\subsubsection{The second order error estimates}
In this subsection we control the error generated from commuting two tangential vector fields with the wave operator. 
We first give an important result. 
\begin{lemma}
 Under the assumptions (\ref{3.12.1.21})-(\ref{6.5.1.21}), there holds
\begin{equation}\label{2.13.1.24}
\bb\Lb^2 \varrho=O(\log \l t\r\l t\r^{-1+\f12\delta}\Delta_0)_{L_u^2 L^\infty_\omega}+\l t\r^{-1}.
\end{equation}
\end{lemma}
\begin{proof}
 We first use \Poincare inequality and Sobolev embedding to derive
\begin{align}
\|\Osc(\bb\Lb^2 \varrho)\|_{L_u^2 L^\infty_\omega}&\les \|\Omega(\bb\Lb^2\varrho)\|_{L_u^2 L_\omega^4}\nn\\
&\les \|\Omega^2(\bb \Lb^2\varrho)\|_{L_u^2 L_\omega^2}^\f12\|\Omega(\bb\Lb^2\varrho)\|_{L_u^2 L_\omega^2}^\f12.\label{2.11.2.24}
\end{align} 
It follows by using (\ref{1.27.3.24}), (\ref{1.29.2.22}) and (\ref{1.27.5.24}) that
\begin{equation}\label{2.11.1.24}
\bb^2 \tir\Omega(\bb\Lb^2\varrho)=O(\log \l t\r\Delta_0)_{L_u^2 L_\omega^2}. 
\end{equation}
For the higher order term, we derive by using (\ref{5.13.10.21}) 
\begin{align*}
[\Omega^2, \bb\Lb]\Lb\varrho&=\Omega ([\bb\Lb, \Omega]\Lb \varrho)+[\Omega, \bb\Lb]\Omega\Lb \varrho\\
&=\sum_{a=0}^1\Omega^{1-a}(\bb\Omega^A \zb_A \bN+\bb\Omega^A \ze_A L+\bb\pioh_{A\Lb}\sn)\Omega^a \Lb \varrho.
\end{align*}
Using (\ref{LbBA2}),  (\ref{3.12.1.21}), (\ref{3.11.3.21}), (\ref{6.22.1.21}), (\ref{8.23.1.23}), (\ref{1.27.5.24}), (\ref{1.29.2.22}), (\ref{5.21.1.21}) and (\ref{2.11.1.24}), we have
\begin{align*}
\Omega^2(\bb \Lb^2\varrho)=\bb\Lb\Omega^2\Lb \varrho+O(\l t\r^{-2+2\delta}\Delta_0\log \l t\r^\f12)_{L_u^2 L_\omega^2}.
\end{align*}
Hence we conclude by using (\ref{3.12.1.21})
\begin{align*}
\|\Omega^2(\bb \Lb^2\varrho)\|_{L^2_u L_\omega^2}&\les\log \l t\r^\f12\l t\r^{-1}\{ E[\Omega^2\Lb \varrho]^\f12(t)+\l t\r^{-1+2\delta}\Delta_0\}\\
&\les \log \l t\r^\f12\l t\r^{-1+\delta}\Delta_0.
\end{align*}
In view of (\ref{1.29.2.22}), substituting the above estimate and (\ref{2.11.1.24}) into (\ref{2.11.2.24}) gives  (\ref{2.13.1.24}). 
\end{proof}
\begin{proposition}\label{8.24.3.21}
 Under the assumptions (\ref{3.12.1.21})-(\ref{6.5.1.21}), there hold with $X_2=S$ or $\Omega$,
\begin{align}
\|\sP[\Omega, X_2\Phi]\|_{L^2_\Sigma}&\les\l t\r^{-2}(\l t\r^{\frac{1}{4}+\delta}\Delta_0^\frac{5}{4}+\log \l t\r^{\f12\M}\La_0);\label{7.6.1.22}
\end{align}
and for $f=\varrho, v$, 
 \begin{align}
&\|X_2(\sP[\Omega, f])-(1-\vs(X_2))X_2\Omega \wt{\tr\chi} \Lb f\|_{L^2_\Sigma}\nn\\
&\les\l t\r^{-2}\sum_{\vs(X)\le \vs(X_2)}W_2[X\Omega\Phi]^\f12(t)+\l t\r^{-\frac{7}{4}+\delta}\Delta_0^\frac{5}{4}+\l t\r^{-2} \log \l t\r^{\f12\M+1}\La_0\label{8.24.1.21},\\
&\|X_2(\sP[S, f])\|_{L_\Sigma^2}\nn\\
&\les \l t\r^{-2}(W_2[\Omega^2\Phi]^\f12(t)+W_2[X_2 S\Phi]^\f12(t))+\l t\r^{-1}(1-\vs(X))\|\bb^{-1}X\bA_b\|_{L^2_\Sigma}\nn\\
&+\l t\r^{-2 }(\log \l t\r^{\f12\M+1}\La_0+\l t\r^{2\delta}\log \l t\r^\f12\Delta_0^\frac{5}{4}),\label{8.7.1.21}\\
&\|\sP[\Omega, \bT\varrho]\|_{L^2_\Sigma}\les \l t\r^{-2} \int_0^t \l t'\r^{-1}W_2[\Omega^2 \Phi]^\f12(t') dt'+\l t\r^{-2}\log \l t\r^{\f12\M+1}\La_0+\l t\r^{-\frac{7}{4}+2\delta}\Delta_0^\frac{5}{4}. \label{2.13.5.24}
\end{align}
\end{proposition}

\begin{proof}

(\ref{7.6.1.22}) is a direct consequence of (\ref{9.2.6.23}) and (\ref{8.21.4.21}).

If $X_1=\Omega$, it follows by differentiating (\ref{6.23.1.23}) by $X_2$ and using (\ref{5.21.1.21}) that
\begin{align}\label{3.4.6.24}
&X_2(\sP[{}\rp{a}\Omega, f]-\bJ[{}\rp{a}\Omega]^\a \p_\a f)+\f12 X_2({}\rp{a}\pih_{LA}\bd^2_{\Lb A}f)\nn\\
&=\sn_{X_2}\big(\piohb(\bd^2_{LL}f+\sn^2 f)+{}\rp{a}\pih_{\Lb A}\bd^2_{LA}f\big)\\
&=O(\l t\r^\delta)\Big(\Delta_0\sn_{X_2}(\bd^2_{LL} f+\sn^2 f)+\Delta_0^\f12 \l t\r^{-\frac{3}{4}} \sn_{X_2}\bd^2_{LA}f\Big)\nn\\
&+O(\Delta_0\l t\r^\delta)_{L_\omega^4}(\bd^2_{LL} f+\sn^2 f)+O(\l t\r^{\delta-\frac{3}{4}}\Delta_0^\f12)_{L_\omega^4}\bd^2_{LA} f.\nn
\end{align}

With $f=\Phi$, using (\ref{5.21.1.21}), (\ref{8.15.3.21}) and (\ref{9.2.4.23}), the right-hand side is bounded by $O(\l t\r^{-2+3\delta}\Delta_0^\frac{3}{2})_{L^2_\Sigma}$.
Next we consider the left-hand side of (\ref{3.4.6.24}). Combining Lemma \ref{9.18.5.23} with Proposition \ref{1steng} gives for $f=\varrho, v$
\begin{align}\label{2.14.2.24}
\|X_2(\pioh_{LA} \bd^2_{\Lb A} f)\|_{L^2_\Sigma}\les\l t\r^{-2}\log \l t\r^{\f12\M}(\La_0+\Delta_0^\frac{5}{4})+\l t\r^{-\frac{7}{4}+\delta}\Delta_0^\frac{5}{4}.
\end{align} 
Using (\ref{6.24.3.23}), (\ref{8.23.1.23}), (\ref{12.19.1.23}), (\ref{5.25.1.21}), (\ref{10.10.2.23}) and (\ref{9.2.4.23}), we derive
\begin{align}
\|{}\rp{a}\bJ^\a \c \sn_{X_2} \p_\a f\|_{L^2_\Sigma}&\les \l t\r^{-1+\delta}\Delta_0\|\bb^\f12\sn_{X_2}\sn_L f\|_{L^2_u L_\omega^4}+\l t\r^{-1}\|\bb^{-1}\sn_{X_2}\sn f\|_{L^2_\Sigma}\nn\\
&+\|\fB\c {}\rp{a}\bJ_L\|_{L^2_\Sigma}+\l t\r^{-1+\delta}\Delta_0\|{}\rp{a}\bJ_L\|_{L^2_\Sigma}+\l t\r^{-2+3\delta}\Delta_0^\frac{3}{2}\nn\\
&\les \l t\r^{-2}(\l t\r^{\frac{1}{4}+\delta}\Delta_0^\frac{5}{4}+\log \l t\r^{\f12\M}\La_0).\label{7.10.2.21}
\end{align}
 
 Using (\ref{7.11.5.21})-(\ref{2.24.1.22}), we infer
 \begin{align}
 \sn_{X_2}{}\rp{a}\bJ^\a \bd_\a f
&=\big(O(\tir^{-1})X^{\le 1}\Phi+cX_2(\Omega \tr\thetac)+X_2\sdiv\eta(\Omega)+X_2L\Omega \varrho\big)\c \Lb f\nn\\
 &+O(1)\fB\c\sn f+\bb^{-1}O(\l t\r^{-3}(\l t\r^{2\delta}\Delta_0^\frac{5}{4}+\log \l t\r^{\f12\M+1}\La_0))_{L^2_u L_\omega^2}\nn\\
 &+O(\l t\r^{\delta-1}\Delta_0)_{L_u^2 L_\omega^2}\c (\sn f, Lf) \label{6.24.7.23}
 \end{align}
 
 For $f=\Phi$, using Proposition \ref{1steng}, we have
\begin{equation*}
\|\tir^{-1}X^{\le 1} \Phi\c \Lb f\|_{L^2_\Sigma}+\|\fB\sn f\|_{L_\Sigma^2}+ \l t\r^{-3}\|S f\|_{L^2_\Sigma}\les \l t\r^{-2}\log \l t\r^{\f12\M}(\Delta_0^\frac{5}{4}+\La_0).
\end{equation*} 
 By (\ref{3.6.2.21}), $|\sn f, Lf|\les \l t\r^{-2+\delta}\Delta_0^\f12$, then all the terms in the second and the third line in (\ref{6.24.7.23}) are $O\big(\l t\r^{-2}(\l t\r^{2\delta}\log \l t\r^\f12\Delta_0^\frac{5}{4}+\log \l t\r^{\f12\M+1}\La_0)\big)_{L^2_\Sigma}$. 

Using (\ref{1.6.1.21}), we write by using (\ref{10.10.2.23}), (\ref{L2conndrv}), (\ref{3.6.2.21}) and (\ref{3.16.1.22}) that
\begin{align*}
cX_2 \Omega\tr\thetac&=X_2\Omega(\tr\chi+L\varrho+\tr\eta)-X_2(\Omega c \tr\thetac)-X_2 c \Omega \tr\thetac\\
&=X_2\Omega(\tr\chi+[L\Phi])+O(\l t\r^{-1}\log \l t\r^{\f12\M}(\Delta_0^\frac{5}{4}+\La_0))_{L^2_\Sigma}
\end{align*}
and the first term on the right-hand side can be written as $X_2\Omega\wt{\tr\chi}$ alternatively. 
  
Hence we bound for $f=\varrho, v$,  with the help of Lemma \ref{10.10.3.23} and (\ref{10.10.2.23})
\begin{align}
X_2(\sP[{}\rp{a}\Omega, f])&=(X_2(\Omega \wt{\tr\chi})+ c(X_2\sdiv\eta(\Omega)+X_2\Omega [L\Phi]+X_2L\Omega \varrho)\big)\c \Lb f\nn\\
&+O(\l t\r^{-\frac{7}{4}+\delta}\Delta_0^\frac{5}{4}+\l t\r^{-2} \log \l t\r^{\f12\M+1}\La_0)_{L^2_\Sigma}\nn\\
&=(X_2(\Omega \wt{\tr\chi})+ c(X_2\sdiv\eta(\Omega)+[LX_2\Omega \Phi])\big)\c \Lb f\nn\\
&+O(\l t\r^{-\frac{7}{4}+\delta}\Delta_0^\frac{5}{4}+\l t\r^{-2} \log \l t\r^{\f12\M+1}\La_0)_{L^2_\Sigma}
\label{9.5.12.23}
\end{align}
where the first term on the right-hand side can be written as $X_2\Omega\tr\chi$ alternatively. 
For the case $X_2=\Omega$, the first term on the right-hand side will be treated by using (\ref{1.26.2.23}). Now we consider the case $X_2=S$.
Using (\ref{3.20.1.22}), (\ref{7.17.7.21}), (\ref{6.6.1.23}), Lemma \ref{comp} and Proposition \ref{1steng}, we bound
\begin{align*}
\|\sn\wt{L\Xi_4}\|_{L^2_\Sigma}&\les \l t\r^{-3}(W_2[S\Omega \Phi]^\f12(t)+W_2[\Omega^2\Phi]^\f12(t))+\|\sn((h, \fB)[L\Phi]), \sn(\zb\sn \varrho)\|_{L^2_\Sigma}\\
&+\l t\r^{-3}(\log \l t\r^{\f12\M}\La_0+\l t\r^{2\delta}\Delta_0^\frac{5}{4}).
\end{align*}
In view of (\ref{8.5.1.22}), applying (\ref{8.24.5.23}), (\ref{8.24.4.23}) and (\ref{8.28.1.23}) to  the right-hand side of the above estimate and also using Proposition \ref{1steng}, we infer
\begin{align}\label{2.16.1.24}
\|\sn\wt{L\Xi_4}\|_{L^2_\Sigma}\les 
\l t\r^{-3}(W_2[S\Omega \Phi]^\f12(t)+W_2[\Omega^2\Phi]^\f12(t))+\l t\r^{-3}(\log \l t\r^{\f12\M}\La_0+\l t\r^{2\delta}\Delta_0^\frac{5}{4})
\end{align}
and
\begin{align*}
\begin{split}
\|\tir^{-2}\sn_L(\tir^3\sn \tr\chi)\|_{L^2_u L_\omega^2}&\les \l t\r^{-3}\Big(W_2[S\Omega \Phi]^\f12(t)+W_2[\Omega^2\Phi]^\f12(t)\Big)\\
&+\l t\r^{-3}(\log \l t\r^{\f12\M}\La_0+\l t\r^{2\delta}\Delta_0^\frac{5}{4}).
\end{split}
\end{align*}
Note due to (\ref{12.19.1.23}) and (\ref{3.11.3.21})
\begin{align*}
\|\Omega \tr\chi \Lb \Phi\|_{L^2_\Sigma}\les \Delta_0^\frac{5}{4}\l t\r^{-\frac{7}{4}+\delta}. 
\end{align*}
Hence 
\begin{align}\label{10.20.2.23}
\|S\Omega\tr\chi \Lb \Phi\|_{L^2_\Sigma}&\les \l t\r^{-2}\Big(W_2[S\Omega \Phi]^\f12(t)+W_2[\Omega^2\Phi]^\f12(t)\Big)\\
&+\l t\r^{-2}\log \l t\r^{\f12\M}\La_0+\Delta_0^\frac{5}{4}\l t\r^{-\frac{7}{4}+\delta}. \nn
\end{align}

Next, using (\ref{8.2.2.23}), (\ref{10.10.2.23}) and Lemma \ref{10.10.3.23} we bound
 \begin{align*}
 \|X_2 \sdiv\eta(\Omega)\|_{L^2_\Sigma}&\les \|\sn_{X_2}^{\le 1}(\tir^{-1}[\Omega v], \Omega[L\Phi])\|_{L^2_\Sigma}+\l t\r^{-\frac{7}{4}+2\delta}\Delta_0^\frac{3}{2}\log \l t\r^\f12\nn\\
 &\les \|L X_2^{\le 1}\Omega \Phi\|_{L^2_\Sigma}+\l t\r^{-1}\log \l t\r^{\f12\M}(\Delta_0^\frac{5}{4}+\La_0).
 \end{align*}

 Using Proposition \ref{1steng}, we bound 
 $$\|LX_2\Omega \Phi\|_{L^2_\Sigma}\les \l t\r^{-1}\left(W_2[X_2 \Omega\Phi]^\f12(t)+(\log\l t\r)^{\f12\M}(\La_0+\Delta_0^\frac{5}{4})\right).$$
Hence we combine the above two estimates with (\ref{9.5.12.23}) and (\ref{10.20.2.23}) to conclude for $f=\Phi$,
\begin{align*}
\|X_2(\sP[\Omega, f])&-(1-\vs(X_2))X_2\Omega \wt{\tr\chi} \Lb f\|_{L^2_\Sigma}\\
&\les \l t\r^{-2}\sum_{\vs(X)\le \vs(X_2)}W_2[X\Omega\Phi]^\f12(t)+\l t\r^{-\frac{7}{4}+\delta}\Delta_0^\frac{5}{4}+\l t\r^{-2} \log \l t\r^{\f12\M+1}\La_0.
\end{align*}
This is (\ref{8.24.1.21}).

Next we prove (\ref{8.7.1.21}) in view of (\ref{6.24.2.23}).  By using   (\ref{5.25.2.21})-(\ref{7.5.2.21}) we derive
\begin{align*}
X_2\big(\piSh^\sharp(\sn^2 f, \bd^2_{LL} f)\big)&=\sn_{X_2}\rp{\le 1}(\sn^2 f, \bd^2_{LL}f)+O(\l t\r^{\delta}\Delta_0^\f12)_{L_\omega^4}\big(\sn^2 f, \bd^2_{LL}f\big),\\
 X_2(\piShb\bd^2_{L A} f)&=O(\Delta_0\l t\r^{\delta})\sn_{X_2}\bd^2_{LA}f+O(\l t\r^{\delta}\Delta_0)_{L_\omega^4}\bd^2_{LA}f.
\end{align*}
Hence
\begin{align}\label{3.5.3.24}
\begin{split}
X_2(\sP[S, f]-\bJ[S]^\a\p_\a f)&=\sn_{X_2}\rp{\le 1}(\sn^2 f, \bd^2_{LL}f)+O(\l t\r^{\delta}\Delta_0^\f12)_{L_\omega^4}\big(\sn^2 f, \bd^2_{LL}f, \bd^2_{LA}f\big)\\
&+O(\Delta_0\l t\r^{\delta})\sn_{X_2}\bd^2_{LA}f
\end{split}
\end{align}
Let $f=\varrho, v$ in the sequel. Then using (\ref{8.15.3.21}), (\ref{9.2.3.23}) and (\ref{10.10.2.23}), we deduce 
\begin{equation}\label{5.7.1.24}
\|\sn^2 \Phi, \bd^2_{LL}\Phi, \bd^2_{LA}\Phi, \tir^{-1}(X L \Phi, (L, \sn)X\Phi, \sn_X \sn \Phi)\|_{L_u^2 L_\omega^4}\les \l t\r^{-3+\frac{1}{2}\delta}\Delta_0(\log \l t\r)^{\frac{\M}{4}}.
\end{equation}
 Combining the above estimates and using Corollary \ref{9.2.5.23} and (\ref{10.10.2.23}) again, we have
\begin{align}
\|X_2(\sP[S, f]&-\bJ[S]^\a\p_\a f)\|_{L^2_\Sigma}\nn\\
&\les \l t\r^{-2+2\delta}\Delta_0^\frac{3}{2}+\|\sn_{X_2}^{\le 1}(\sn^2 f, \bd^2_{LL}f)\|_{L^2_\Sigma}\nn\\
&\les \l t\r^{-2}(\log \l t\r^{\f12\M}\La_0+\l t\r^{2\delta}\Delta_0^\frac{5}{4})+\|\sn_{X_2}(\sn^2 f, \bd^2_{LL}f)\|_{L^2_\Sigma}.\label{9.5.14.23}
\end{align}
Using (\ref{6.24.4.23}), (\ref{8.23.1.23}),  (\ref{8.21.1.22}), (\ref{10.10.2.23}) and (\ref{5.7.1.24})  we bound
\begin{align}
\|\bJ[S]^\a \sn_{X_2} \sn_\a f \|_{L_\Sigma^2}&\les \|\fB L{}\rp{S}\ss\|_{L^2_\Sigma}+\l t\r^{-1+\delta}\Delta_0\|L{}\rp{S}\ss\|_{L^2_\Sigma}+\l t\r^{-1}\|\sn_{X_2}Lf\|_{L^2_\Sigma}\nn\\
&+\l t\r^{-1+\delta}\Delta_0^\f12\|\bb^\f12\tir(\sn_{X_2}\sn f, \bb^{-1}\log \l t\r\sn_{X_2}Lf)\|_{L_u^2 L_\omega^4}\nn\\
&\les \l t\r^{-2}\log \l t\r^{\f12\M+1}(\La_0+\l t\r^{\frac{3}{2}\delta}\Delta_0^\frac{5}{4}).\label{9.5.15.23}
\end{align}

Next we consider $\sn_{X_2}\bJ[S]^\a\p_\a f$. Note that by using (\ref{7.12.1.21}) and (\ref{3.6.2.21}) 
\begin{align*}
\|\sn_X \bJ[S]_A \sn_A f\|_{L_\Sigma^2}&\les  \Delta_0 \l t\r^{\delta}\|\sn f\|_{L^\infty_x}\log \l t\r^\f12\les \l t\r^{-2+2\delta}\log \l t\r^\f12\Delta_0^\frac{3}{2}.
\end{align*}
Using (\ref{10.19.1.23}), (\ref{3.6.2.21}) and Proposition \ref{1steng}, we derive
\begin{align*}
\|\sn_{X_2} \bJ[S]_\Lb Lf \|_{L^2_\Sigma}&\les \l t\r^{-1}\| L f\|_{L^2_\Sigma}+\l t\r^{-1+\delta}\Delta_0\log \l t\r^\frac{1}{2}\|Lf\|_{L^\infty_x}\\
&\les \l t\r^{-2}\log \l t\r^{\f12\M}\La_0+\l t\r^{-2+2\delta}\log\l t\r^\f12\Delta_0^\frac{5}{4}.
\end{align*}

At last using  (\ref{7.9.6.22}) we bound
\begin{align*}
\|X_2L{}\rp{S}\ss \Lb f\|_{L^2_\Sigma}&\les\l t\r^{-1}\|\bb^{-1} X_2^{\le 1} L{}\rp{S}\ss\|_{L^2_\Sigma}\\
&\les \l t\r^{-2}\left(W_2[\Omega^2 \Phi]^\f12(t)+W_2[X_2S\Phi]^\f12(t)\right)+\l t\r^{-2}\log \l t\r^{\f12\M+1}(\La_0+\Delta_0^\frac{5}{4})\\
&+(1-\vs(X))\l t\r^{-1}\|\bb^{-1}X\bA_b\|_{L^2_\Sigma}.
\end{align*} 
We then combine the above estimates with (\ref{9.5.14.23}) and (\ref{9.5.15.23}) and use Lemma \ref{comp} to obtain
\begin{align*}
\|X_2(\sP[S, f])\|_{L_\Sigma^2}&\les \l t\r^{-2}\left(W_2[\Omega^2 \Phi]^\f12(t)+W_2[X_2 S\Phi]^\f12(t)\right)+(1-\vs(X))\l t\r^{-1}\|\bb^{-1}X\bA_b\|_{L^2_\Sigma}\\
&+\l t\r^{-2 }(\log \l t\r^{\f12\M+1}\La_0+\l t\r^{2\delta}\log \l t\r^\f12\Delta_0^\frac{5}{4}).
\end{align*}
 Thus (\ref{8.7.1.21}) is proved.
 
 Next we prove (\ref{2.13.5.24}). Applying (\ref{2.14.3.24}) to $f=\bT\varrho$, we write
 \begin{align}\label{2.15.6.24}
 \begin{split}
\sP[\Omega, \bT\varrho]&+\f12 {}\rp{a}\pih_{LA}\bd^2_{\Lb A} \bT \varrho+\f12\Lb \bT\varrho \bJ[\Omega]_L\\
&= O(\l t\r^{\delta}\Delta_0)(\sn^2 \bT \varrho+\bd^2_{LL}\bT \varrho) +O(\l t\r^{-\frac{3}{4}+\delta}\Delta_0^\f12)\bd^2_{LA} \bT\varrho\\
&+O(\l t\r^{-1+\delta}\Delta_0)_{L_\omega^4}(L \bT\varrho+\l t\r^{-1+\delta}\sn \bT \varrho)+O(\fB+\l t\r^{-1+2\delta}\Delta_0)\sn\bT\varrho.
\end{split}
 \end{align}
 Applying  (\ref{7.03.4.21}) to $f=\bT \varrho$ gives
 \begin{align}\label{2.15.2.24}
 \begin{split}
 &|\bd^2_{LL}\bT\varrho|\les |LL\bT \varrho|+|k_{\bN\bN}||L\bT\varrho|\\
 &\tir|\bd^2_{LA}\bT \varrho|\les |L \Omega \bT\varrho|+|\sn \bT\varrho|+|\tir \zb||L\bT \varrho|\\
 &\tir|\bd^2_{\Lb A}\bT\varrho|\les |\Omega\Lb \bT\varrho|+|\sn \bT \varrho|+\tir|\bA_{g,1}\Lb \bT \varrho|.
 \end{split}
 \end{align}
To bound the above quantities in $L^2_\Sigma$,  
using (\ref{8.23.1.23}), also using (\ref{6.22.1.21}) and (\ref{LbBA2}) we have
\begin{align}\label{2.15.5.24}
\begin{split}
&\|\tir \bd^2_{LA}\bT \varrho\|_{L^2_\Sigma}\les\Delta_0\l t\r^{-1+\delta}\\
&\|\bb^\f12(\tir \bd^2_{LL}\bT\varrho-O(\tir^{-1})\fB)\|_{L^2_u L_\omega^2}\les\l t\r^{-3+\delta}\Delta_0.
\end{split}
\end{align}
We also claim 
\begin{equation}\label{2.15.4.24}
\|\tir\bd^2_{\Lb A}\bT\varrho\|_{L^2_\Sigma}\les \l t\r^\delta\Delta_0
\end{equation}
which is obtained by assuming
\begin{align}\label{2.15.3.24}
\Omega\Lb \bT\varrho&=\Lb \Omega \bT\varrho\nn\\
&+\left\{\begin{array}{lll}
O(\Delta_0 \log \l t\r)_{L^2_\Sigma}, \mbox{ or }\\
\bb^{-2}\Omega \bb O(\l t\r^{-1})_{L_\omega^4}+O(\l t\r^{-1}\log \l t\r^{\f12\M}(\La_0+\Delta_0^\frac{5}{4}))_{L^2_\Sigma}
 \end{array}\right.
\end{align}
combined with (\ref{3.12.1.21}), (\ref{LbBA2}), (\ref{3.6.2.21}) and (\ref{1.29.2.22}).

Substituting (\ref{2.15.5.24}) and (\ref{2.15.4.24}) to (\ref{2.15.6.24}), also using Proposition \ref{7.15.5.22},  (\ref{8.23.2.23}), (\ref{12.19.1.23}) and (\ref{3.25.1.22})  gives
\begin{align*}
\|\sP[\Omega, \bT\varrho]&+\f12\Lb \bT\varrho \bJ[\Omega]_L\|_{L^2_\Sigma}\les \l t\r^{-2+2\delta}\log \l t\r\Delta_0^\frac{5}{4}+\l t\r^{-\frac{7}{4}+2\delta}\Delta_0^\frac{3}{2}.
\end{align*}
To treat the second term on the left-hand side, by using  (\ref{5.7.2.24}), (\ref{5.25.1.21}) and (\ref{2.13.1.24}), we deduce
\begin{align*}
\|\bJ[\Omega]_L \Lb \bT\varrho\|_{L^2_\Sigma}&\les \|\bb^{-\frac{1}{2}}\bJ[\Omega]_L\|_{L_u^2 L_\omega^2}+\|\bb^{-\frac{1}{2}}\bJ[\Omega]_L\tir\|_{L_\omega^2}\log \l t\r\l t\r^{-1+\f12\delta}\Delta_0\\
&\les \l t\r^{-2}\int_0^t \l t'\r^{-1}W_2[\Omega^2\Phi]^\f12(t')+\l t\r^{-2}\log \l t\r^{\f12\M+1}\La_0\\
&+\l t\r^{-\frac{7}{4}+2\delta}\Delta_0^\frac{3}{2}. 
\end{align*}
Therefore,  we conclude
\begin{align*}
\|\sP[\Omega, \bT\varrho]\|_{L^2_\Sigma}&\les \l t\r^{-2}\int_0^t \l t'\r^{-1}W_2[\Omega^2\Phi]^\f12(t')+\l t\r^{-2}\log \l t\r^{\f12\M+1}\La_0+\l t\r^{-\frac{7}{4}+2\delta}\Delta_0^\frac{5}{4}. 
\end{align*}
It remains to prove (\ref{2.15.3.24}). 
We first write
\begin{align}\label{2.16.3.24}
\Omega \Lb \bT \varrho=\Lb \Omega \bT\varrho-\bb^{-1}\Omega \bb \Lb \bT \varrho+\bb^{-1}[\Omega, \bb\Lb]\bT\varrho.
\end{align}
We recall from (\ref{5.13.10.21}) for the last term on the right-hand side, 
and deduce by using  (\ref{5.21.1.21}), (\ref{6.22.1.21}), (\ref{3.11.3.21}), (\ref{1.29.2.22}) and Proposition \ref{1steng} that
\begin{align*}
&\|\tir\bb^{-\f12}[\bb \Lb, \Omega]\bT \varrho\|_{L^2_u L_\omega^2}\\
&\les \|\bb^\f12 \tir\Omega^A \zb_A \bN \bT\varrho\|_{L_u^2 L_\omega^2}+\|\bb^\f12 \Omega^A \ze_A L\bT \varrho\tir\|_{L^2_u L_\omega^2}+\|\bb^\f12 \pioh_{\Lb A}\sn_A \bT\varrho\tir \|_{L^2_u L_\omega^2}\\
&\les \l t\r^{-1}\log \l t\r^{\f12\M}(\La_0+\Delta_0^\frac{5}{4})+\Delta_0^\frac{5}{4}\log \l t\r^2\l t\r^{-1}+\l t\r^{-2+2\delta}\log \l t\r^\f12\Delta_0^\frac{3}{2}\\
&\les \l t\r^{-1}\log \l t\r^{\f12\M}(\La_0+\Delta_0^\frac{5}{4})
\end{align*}
where we used the fact that $\M>4$. It follows by using the above estimate, (\ref{1.29.2.22}) and (5) in Lemma \ref{5.13.11.21} that
\begin{align*}
\|\Omega\Lb \bT\varrho\|_{L^2_\Sigma}&\les\|\Lb \Omega \bT \varrho\|_{L^2_\Sigma}+|\bb^{-\f12}\tir \Omega \bb \Lb \bT\varrho\|_{L^2_u L_\omega^2}+\|\tir\bb^{-\f12}[\bb \Lb, \Omega]\bT \varrho\|_{L^2_u L_\omega^2}\\
&\les \|\Lb\Omega\bT \varrho\|_{L^2_\Sigma}+\log \l t\r\Delta_0.  
\end{align*}
Without bounding the second term on the right-hand of (\ref{2.16.3.24}) completely, applying merely the estimate (\ref{1.29.2.22}) instead, we can get the alternative bound stated in (\ref{2.15.3.24}). 
\end{proof}
\begin{lemma} 
There holds 
\begin{align}
\Box_\bg \bT \varrho&=\bT \Box_\bg \varrho+\bN \varrho (\div_g k)_\bN+\sn \varrho (\div_g k)_A+k^{\bN \bN} (\bN^2\varrho+\ud\bA\sn\varrho)\label{2.18.11.24}\\
&+k^{AB}(\sn^2_{AB}\varrho+\theta_{AB}\bN\varrho)+k^{A\bN}(\sn(\bN\varrho)+\theta\c \sn\varrho)\nn\\
&+(\Delta_g \varrho+\bT \varrho \Tr k)\Tr k+\nab^i \bT \varrho\p_j \varrho+\bT \varrho \Box_\bg \varrho\nn
\end{align}
where \begin{align}\label{2.18.10.24}
\begin{split}
(\div_g k)_A&=\fB \bA_{g,1}+\sn \bT \varrho+\sn \bAn+(\tr\theta+\bA_g) \bA_{g,1}+\sn_\bN \bA_{g,1}+(\fB+\bAn) \ud\bA\\
(\div_g k)_\bN&=\fB^2+[\sn\Phi]^2+\bN\fB+\sn\bA_{g,1}+\tr\theta\fB+(\ud\bA+\bA_g)\bA_{g,1}.
\end{split}
\end{align} 
\end{lemma}
\begin{proof}
Using (\ref{10.12.1.23}), we have 
\begin{equation}\label{2.18.5.24}
\Box_\bg \bT f=\bT \Box_\bg f+\bd_i\big((-k^{ij}+\f12 \Tr k\bg^{ij} )\p_j f\big)-\f12 \Tr k\Box_\bg f.
\end{equation}
 Hence
\begin{equation}\label{2.18.6.24}
\Box_\bg \bT \varrho=\bT \Box_\bg \varrho-(\div k^j-\f12\bd^j \Tr k) \p_j \varrho-(k^{ij}-\f12 \Tr k \bg^{ij})\bd_i \p_j \varrho-\f12\Tr k \Box_\bg \varrho.
\end{equation} 
Using (\ref{2.18.4.24}), (\ref{4.23.1.19}), Proposition \ref{6.7.1.23} (2) and (3), 
\begin{align}\label{2.18.7.24}
(\div_g k)_A=\fB\c \bA_{g,1}+\p_A \bT \varrho+\Delta v^\|; (\div_g k)_\bN=\fB^2+[\sn\Phi]^2+\p_\bN \bT \varrho+[\Delta v].
\end{align}
Using \cite[P86 (4.2.9)]{CK}, we write
\begin{align*}
[\Delta v]&=[\sD v]+\tr\theta[\bN v]+[\bN\bN v]+\sn\log \bb [\sn v]\\
\Delta v^\|&=\sD v^\|+\tr\theta\bN v^\|+\bN^2 v^\|+\sn\log \bb\sn v^\|.
\end{align*}
Moreover we have $\bd^i \p_i \varrho=\Delta_g \varrho+\bT \varrho \Tr k$ and schematically $\p_i \varrho \nab^i\Tr k=\nab^i\bT\varrho\p_i \varrho$.

In view of (\ref{5.23.1.23}), substituting the above identities to (\ref{2.18.7.24}) yields (\ref{2.18.10.24}).
We further rewrite (\ref{2.18.6.24}) symbolically into (\ref{2.18.11.24}). 
\end{proof}

\subsection{The  second order energy estimates}\label{2ndengpr}
In this subsection, we use the results in Proposition \ref{9.5.8.23} and Proposition \ref{8.24.3.21} to prove the following second order energy estimates.
\begin{proposition}[The preliminary second order energy estimates]\label{8.29.8.21} 
Let $X_1, X_2\in \{\Omega, S\}$ and $u_0\le u\le u_*$ and $0<t<T_*$. Under the assumptions of (\ref{3.12.1.21})-(\ref{6.5.1.21}), we have
\begin{equation}\label{8.25.2.21}
\begin{split}
&\sup_{u_0\le u\le u_*} F_0[\Omega^2\varrho](\H_u^t)+ E[\Omega^2\varrho](t)\les \La_0^2+\Delta_0^\frac{5}{2}\\
&WFIL_2[X_2 X_1 \Phi](\D_u^t)\les (\log \l t\r)^\M(\La_0^2+\Delta_0^\frac{5}{2}+(1-\vs^+(X^2))F_0[\Omega^3\varrho]\Delta_0^\frac{1}{2})\\
&WFIL_2[\Omega^2\Phi](\D_{u}^{t})\les (\Delta_0^\frac{5}{2}\l t\r^{2\delta}+\La_0^2)(\log \l t\r)^\M.\\
&WFIL_1[\Omega^2\Phi](\D_u^t)\les (\log \l t\r)^\M(\La_0^2+\Delta_0^\frac{5}{2}).
\end{split}
\end{equation}
\begin{align}
&\|X^3\Phi, \Sc(X^3\Phi), \Ac(X^3\Phi)\|_{L^2_\Sigma}\nn\\
&\qquad\qquad\qquad\les (\log \l t\r)^{\frac{\M}{2}} (\La_0+\l t\r^{\delta(1-\vs^+(X^3))}\Delta_0^\frac{5}{4}), \quad X\in \{\Omega, S\}\label{8.29.9.21}\\
&\|X^2\left(\Box_\bg \Phi, \N(\Phi, \bp\Phi)\right)\|_{L^2_\Sigma}\les \l t\r^{-2} \log \l t\r^{\f12\M}(\La_0+\Delta_0^\frac{5}{4})\nn\\
&\qquad\qquad\qquad\qquad\qquad\qquad+(1-\vs^-(X^2))\l t\r^{-2+2\delta}\Delta_0^\frac{3}{2} \label{8.30.3.21}
\end{align}
\begin{equation}\label{8.30.3.21+}
\|\Box_\bg (X^2 \Phi), X\Box_\bg X\Phi\|_{L^2_\Sigma}\les\left\{\begin{array}{lll}
\l t\r^{-\frac{3}{2}+\delta}\Delta_0^\frac{5}{4}+\l t\r^{-\frac{3}{2}}\La_0, X^2=\Omega^2\\
\l t\r^{-2}(\log \l t\r^{\f12\M+1}\La_0+\Delta_0^\frac{5}{4}\l t\r^{\frac{1}{4}+2\delta}), X^2=\Omega S, S\Omega\\
\l t\r^{-2}(\log \l t\r^{\f12\M+1}\La_0+\l t\r^{2\delta}\log \l t\r^\f12\Delta_0^\frac{5}{4}), X^2=S^2.
\end{array}\right.
\end{equation}
\begin{align}
E[\Omega\bT \varrho](t)+\|\Omega\Lb \bT\varrho\|_{L^2_\Sigma}^2+F_0[\Omega\bT \varrho](t)\les\log \l t\r^{\M+7}(\La_0^2+\Delta_0^\frac{5}{2})\label{2.13.3.24}
\end{align}
\begin{align}\label{3.9.8.24}
\Omega\Box_\bg \bT\varrho, \Box_\bg \Omega\bT \varrho=O\Big(\l t\r^{-1}\log \l t\r^{\f12(\M+7)}(\La_0+\Delta_0^\frac{5}{4})\Big)_{L^2_\Sigma} 
\end{align}
\begin{align}
&\|\Omega \wt{\tr\chi}, \sn_\Omega^{\le 1} \chih\|_{L_u^2 L_\omega^2}\les  \l t\r^{-2}\log \l t\r^{\f12\M+1}(\La_0+\l t\r^\delta\Delta_0^\frac{5}{4})\label{2.18.2.24}\\
&\|\bb^{-\f12}\Omega\mho\|_{L^2_\Sigma}(\les \log \l t\r)^2(\La_0+\Delta_0^\frac{5}{4})\label{3.10.9.24}\\
& \|\bb^{-1} X L {}\rp{S}\ss\|_{L^2_\Sigma}\les\l t\r^{-1}\log \l t\r^{\f12\M+1}(\La_0+\l t\r^{(1-\vs(X))\delta}\Delta_0^\frac{5}{4})\label{3.5.5.24}\\
&\|\Omega^a(\bb^{-\a}, \log \bb), (S \Omega^a,\Omega S\Omega^{a-1}, \Omega^a S) (\bb^{-\a}, \log \bb) \|_{L_u^2 L_\omega^2}\nn\\
&\qquad\qquad\qquad\les \log \l t\r(\La_0+\log \l t\r^{\frac{5}{2}(a-1)}\Delta_0^\frac{5}{4}). a=1,2, \a\ge-\f12\label{2.20.2.24}\\
&\|\tir\Omega^2 (\Lb \varrho, \fB)\|_{L^2_u L_\omega^2}\les \log \l t\r(\La_0+\log \l t\r^\frac{5}{2}\Delta_0^\frac{5}{4})\label{2.14.1.24}\\
&\|\Omega \Lb \Omega \varrho\|_{L^2_\Sigma}\les \La_0+\log \l t\r^2\Delta_0^\frac{5}{4}\label{2.22.3.24}\\
&\| \Omega (\Lb[L\Phi], L\fB), \Lb \Omega[L\Phi], L\Omega\fB\|_{L^2_\Sigma}\les \l t\r^{-1}\log \l t\r(\Delta_0^\frac{5}{4}+\La_0).\label{2.19.1.24}
\end{align}
\begin{align}\label{3.12.4.24}
(X\Lb\Phi^\mu,& \Lb X\Phi^\mu)-(\mu+\vs(X))\fB-(1-\vs(X))X\fB\nn\\
&=\left\{\begin{array}{lll}
(1-\vs(X))O(\l t\r^{-1}\log \l t\r(\La_0+\log \l t\r^3\Delta_0^\frac{5}{4}))_{L_u^2 L_\omega^4}\\
(1-\vs(X))(\l t\r^{-1}\log \l t\r\Delta_0)_{L_\omega^4}
\end{array}\right.
\end{align}
\end{proposition}
We will divide the proof into three parts, with (\ref{2.18.2.24})-(\ref{3.12.4.24}) as the second part, (\ref{2.13.3.24})-(\ref{3.9.8.24}) as the last part and the rest of the estimates as the first part. 
\begin{proof}[Proof of (\ref{8.25.2.21})-(\ref{8.30.3.21+})]
Using (\ref{9.16.1.23}), Lemma \ref{5.13.11.21} (1) and Lemma \ref{10.10.3.23}, we derive
\begin{align}
&\|X^2 \Box_\bg\Phi\|_{L^2_\Sigma}\nn\\
&\les \l t\r^{-1}\|X^2 L\Phi, X^2[L\Phi]\|_{L^2_\Sigma}+\l t\r^{-2}(\log \l t\r)^{\f12\M}(\La_0+\Delta_0^\frac{5}{4})+(1-\vs^-(X^2))\l t\r^{-2+2\delta}\Delta_0^\frac{3}{2}
\nn\\
&\les \l t\r^{-2}W_2[X^2\Phi]^\f12(t)+\l t\r^{-2}(\log \l t\r)^{\f12\M}(\La_0+\Delta_0^\frac{5}{4})+(1-\vs^-(X^2))\l t\r^{-2+2\delta}\Delta_0^\frac{3}{2},\label{7.9.4.22}
\end{align}
where we used Proposition \ref{1steng} and to bound the lower order terms.

Recall from (\ref{5.02.4.21}) 
\begin{equation}\label{6.24.9.23}
\Er_2(\Phi, X_2 X_1)=I[X_2 X_1\Phi]+X_2(\sP[X_1, \Phi])+\sP[X_2, X_1\Phi]
\end{equation}
where, with the convention given in (\ref{3.29.2.23}),
\begin{equation*}
I[X_2 X_1\Phi]=\fm{X_2}\Box_\bg X_1\Phi+X_2(\fm{X_1}\Box_\bg \Phi).
\end{equation*}
Recall that $\fm{\Omega}=\ud\bA\c \Omega$ and $\fm{S}= 1+\tir k_{\bN\bN}$. Using Proposition \ref{7.15.5.22}, Lemma \ref{5.13.11.21} (5) and (\ref{10.30.1.21}), we have
\begin{align}\label{6.24.8.23}
\begin{split}
&|\fm{S}|\les 1,\quad |\fm{\Omega}|\les \l t\r^{\delta}\Delta_0\\
&Y\fm{X}= \vs^-(YX)(O(\bb^{-1}))+(1-\vs^-(YX))O (\l t\r^\delta\Delta_0)_{L_\omega^4} 
\end{split}
\end{align}
where $Y, X\in \{\Omega, S\}$. Using the above estimates, Sobolev embedding on spheres, (\ref{8.23.1.21}) and (\ref{8.26.4.21}), we obtain
\begin{align}
\|I[X^2 \Phi]\|_{L^2_\Sigma}&\les(1+(1-\vs^-(X^2))\l t\r^\delta\Delta_0)\sum_{X=S, \Omega}(\|\Box_\bg X\Phi\|_{L^2_\Sigma}+\|X\Box_\bg \Phi, \log \l t\r^{\f12}\Omega^{\le 1}(\Box_\bg \Phi)\|_{L^2_\Sigma})\nn\\
&\les  \l t\r^{-2}(\log \l t\r^{\f12\M+1}(\La_0+\l t\r^\delta\Delta_0^\frac{5}{4})+\l t\r^{(1-\vs^-(X^2))(\frac{1}{4}+\delta)}\Delta_0^\frac{5}{4})\nn\\
&\cdot(1+(1-\vs^-(X^2))\Delta_0 \l t\r^\delta).\label{8.17.4.23}
\end{align}

Next we bound the weighted energies in the order of
(1) $X_2=X_1=\Omega,$ (2) $X_2=\Omega,\, X_1=S$; or  $X_2=S, \, X_1=\Omega$, (3) $X_2=X_1=S$.

\noindent{\bf (1) $X_2=X_1=\Omega$.} 

Using (\ref{1.26.2.23}), we first derive by using Proposition \ref{1steng}
\begin{align}
\|\tir^2 \bb \Omega^2&\wt{\tr\chi} \Lb \Phi\|_{L_u^2 L_t^2 L_\omega^2(\D_{u_1}^{t_1})} \les\|\bb \Lb\Phi\tir\|_{L_u^2 L^\infty_t L_\omega^\infty}\sup_u (F_0[\Omega^3\varrho]^\f12(\H_u^t)+\Delta_0^\frac{3}{2}+\La_0)\nn\\
&+\|(L+\f12\tr\chi)\Omega^{1+\le 1}\Phi, \sn\Omega^{\le 1}\Phi\|_{L_u^2 L^2(\H_u^t)}\nn\\
 &\les \Delta_0^\frac{1}{4}\sup_u F_0[\Omega^3\varrho]^\f12(\H_u^{t_1})+(\int_0^{t_1} \l t'\r^{-2}W_2[\Omega^2\Phi](t')+\int_{\D_{u_1}^{t_1}}|\sn\Omega^{\le 1}\Phi \tir|^2 dt d\omega du)^\f12 +\Delta_0^\frac{5}{4}\nn\\
 &\les \Delta_0^\frac{1}{4}\sup_u F_0[\Omega^3\varrho]^\f12(\H_u^{t_1})+(\int_0^{t_1} \l t'\r^{-2}W_2[\Omega^2\Phi](t'))^\f12+\Delta_0^\frac{5}{4}+\La_0. \label{2.17.1.24}
\end{align}
Similarly, using (\ref{1.26.2.23}) and (\ref{12.19.1.23}), we deduce
\begin{align}
\|\tir^2 \bb \Omega^2\wt{\tr\chi} \Lb \Phi\|_{L_u^2 L_\omega^2}&\les \l t\r^\delta\Delta_0\|\tir^\f12 \bb\Lb\Phi\|_{L_u^2 L_\omega^\infty}\les \l t\r^{-\f12+\delta}\Delta_0^\frac{5}{4}.\label{2.17.2.24}
\end{align}
  Using  (\ref{7.6.1.22}), (\ref{8.24.1.21}) and (\ref{8.17.4.23}), 
   we bound 
\begin{align}
\|\Er_2(\Phi,\Omega^2)-\Omega^2\wt{\tr\chi}\Lb\Phi\|_{L^2_\Sigma}& \les\l t\r^{-2}W_2[\Omega^2\Phi]^\f12(t')+\l t\r^{-\frac{7}{4}+2\delta}\Delta_0^\frac{5}{4}+\l t\r^{-2} (\log \l t\r)^{\frac{\M}{2}+1}\La_0\label{10.22.1.23}.
\end{align}
 Using  $\Box_\bg \Omega^2\Phi=\Omega^2\Box_\bg\Phi+\Er_2(\Phi, \Omega^2)$, the above estimate, (\ref{2.17.1.24}) and (\ref{7.9.4.22}),
   Proposition \ref{MA2} with $n=0$, we obtain for $\Psi=\Omega^2\Phi$
   \begin{align}\label{9.19.1.23}
   \begin{split}
& WFIL_2[\Psi](\D_{u_1}^{t_1})\\
 &\les \La_0^2+\Delta_0^\frac{5}{2}+\int_0^{t_1} \Big(\sup_{\Sigma_{t'}}(-[\Lb\varrho]_{-})+\l t'\r^{-1-\delta}\Big)WL_2[\Psi](t') dt'\\
&+\int_{\D_{u_1}^{t_1}} \big| \Box_\bg \Psi(L\Psi+ h \Psi)\big|\tir^2 \bb d\mu_{\ga} du dt+\int_{u_1}^{u_*} F_2[\Psi](\H_u^{t_1}) du\\
&\les \La_0^2+\Delta_0^\frac{5}{2}+\int_0^{t_1} \Big(\sup_{\Sigma_{t'}}(-[\Lb\varrho]_{-})+\l t'\r^{-1-\delta}\Big)WL_2[\Psi](t') dt'\\
&+\|\bb^\f12 \tir (\Omega^2\Box_\bg\Phi+\Er_2(\Phi, \Omega^2))\|_{L^2_t L_\Sigma^2}^2+\int_{u_1}^{u_*} F_2[\Psi](\H_u^{t_1}) du\\
&\les \La_0^2+\Delta_0^\frac{5}{2}+\int_0^{t_1} \Big(\sup_{\Sigma_{t'}}(-[\Lb\varrho]_{-})+\l t'\r^{-1-\delta}\Big)WL_2[\Psi](t') dt'\\
&+\int_{u_1}^{u_*} F_2[\Psi](\H_u^{t_1}) du+\int_0^t \l t'\r^{-2}W_2[\Omega^2\Phi](t')dt'+ \Delta_0^\frac{1}{2}\sup_{u\in[u_1, u_*]} F_0[\Omega^3\varrho](\H_u^{t_1}).
\end{split}
   \end{align}
   Hence using Gronwall's inequality we obtain from the above energy inequality that
   \begin{equation}\label{11.26.1.23}
WFIL_2[\Omega^2\Phi](\D_{u_1}^{t_1})\les (\Delta_0^\frac{5}{2}+\La_0^2+\Delta_0^\frac{1}{2}\sup_u F_0[\Omega^3\varrho](\H_u^{t_1}))(\log \l t_1\r)^\M.
   \end{equation}
   This gives the case $X^2=\Omega^2$ in the second line of (\ref{8.25.2.21}). It also gives the rough estimate in view of (\ref{3.12.1.21})
   \begin{equation}\label{2.17.3.24}
   WFIL_2[\Omega^2\Phi](\D_{u_1}^{t_1})\les (\Delta_0^\frac{5}{2}\l t_1\r^{2\delta}+\La_0^2)\log \l t_1\r^\M.
   \end{equation}
   as stated in (\ref{8.25.2.21}).
Moreover, repeating the above calculation, we can obtain the last estimate in (\ref{8.25.2.21}).
   
Similarly, using (\ref{6.21.2.21}), (\ref{2.17.2.24}) and (\ref{10.22.1.23}) we have the standard energy estimate 
    \begin{align}\label{10.21.1.23}
    \begin{split} 
   E[\Omega^2\varrho](t_1)+F_0[\Omega^2\varrho]&(\H_{u_1}^{t_1})\les\int_{\D_{u_1}^{t_1}}|\Box_\bg \Omega^2\varrho \bT\Omega^2\varrho| +E[\Omega^2\varrho](0)+F_0[\Omega^2\varrho](\H_{u_*}^{t_1})\\
   &\les \int_{\D_{u_1}^{t_1}}|(\Omega^2\Box_\bg\varrho+\Er_2(\Phi, \Omega^2))  \bT\Omega^2\varrho|+\La_0^2\\
   &\les \Delta_0^\frac{5}{2}+\La_0^2+\int_0^{t_1} \l t'\r^{-1-\a}\|\bT\Omega^2\varrho(t')\|^2_{L^2_\Sigma}+\int_0^{t_1} \l t\r^{-3+\a}W_2[\Omega^2\Phi](t')\\
   &\les \Delta_0^\frac{5}{2}+\La_0^2+\int_0^{t_1} \l t'\r^{-1-\a}\|\bT\Omega^2\varrho(t')\|^2_{L^2_\Sigma}
   \end{split}
    \end{align}
    where we applied (\ref{2.17.3.24}) to obtain the last line. 
    
It follows by using Gronwall's inequality that
    \begin{align*}
E[{\Omega^2}\varrho](t_1)+F_0[{\Omega^2}\varrho](\H_{u_1}^{t_1})\les \La_0^2+\Delta_0^\frac{5}{2}
    \end{align*}
as stated in (\ref{8.25.2.21}).
  
    Substituting (\ref{2.17.3.24}) to (\ref{7.9.4.22}) leads to
\begin{align}\label{8.13.2.21}
\|\Omega^2 \Box_\bg \Phi\|_{L^2_\Sigma}\les \l t\r^{-2}(\log \l t\r)^{\f12\M}(\La_0+\l t\r^{2\delta}\Delta_0^\frac{5}{4})+\l t\r^{-2+2\delta}\Delta_0^\frac{3}{2}.
\end{align}
Substituting (\ref{2.17.3.24}) into  (\ref{10.22.1.23}) and using (\ref{2.17.2.24}), we derive 
\begin{align}\label{2.17.4.24}
\|\Er_2(\Phi, \Omega^2)\|_{L^2_\Sigma}\les \l t\r^{-\frac{3}{2}+\delta}\Delta_0^\frac{5}{4}+\l t\r^{-\frac{3}{2}}\La_0, \|\bb^\f12 \tir \Er_2(\Phi, \Omega^2)\|_{L^2[0,t] L^2_\Sigma}\les \l t\r^\delta\Delta_0^\frac{5}{4}+\log \l t\r\La_0, 
\end{align}
which could be further improved whence the top order energy estimates are obtained. 

\noindent{\bf (2) Case $X_2=\Omega,\, X_1=S$}.
Next we apply (\ref{8.24.1.21}) to $X_2=S$ and $f=\Phi$.  Using (\ref{8.24.1.21}) and using the $X_2 X_1=\Omega^2$ case in  (\ref{8.25.2.21}), we have
\begin{align*}
\|S(\sP[\Omega, \Phi])\|_{L^2_\Sigma}&\les \l t\r^{-2}W_2[S\Omega\Phi]^\f12(t)+\l t\r^{-\frac{7}{4}+\delta}\Delta_0^\frac{5}{4}+\l t\r^{-2} \log \l t\r^{\f12\M+1}\La_0.
\end{align*}
Applying (\ref{8.21.2.21}) to $f=\Omega \Phi$ with the help of (\ref{5.7.1.24}),  using the $X_2X_1=\Omega^2$ case in (\ref{8.25.2.21}), we have
\begin{align*}
\|\sP[S, \Omega \Phi]-\frac{1}{4}L{}\rp{S}\ss \Lb \Omega \Phi\|_{L^2_\Sigma}&\les  \l t\r^{-2} \left(W_2[\Omega^2\Phi]^\f12(t)+W_2[S\Omega\Phi]^\f12(t)\right)\nn\\
 &+ \l t\r^{-1} \|L\Omega\Phi\|_{L^2_\Sigma}+\l t\r^{-2+2\delta}\Delta_0^\frac{3}{2}
\end{align*}
Note that due to (\ref{8.23.2.23}), (\ref{7.3.1.22}) and (\ref{8.21.1.22})
\begin{align}\label{2.18.1.24}
\|L{}\rp{S}\ss \Lb X\Phi\|_{L^2_\Sigma}\les \l t\r^{-2}\log \l t\r^{\f12\M+1}(\La_0+\Delta_0^\frac{5}{4})+(1-\vs(X))\l t\r^{-2+2\delta}\Delta_0^2\log \l t\r^\f12.
 \end{align}
 Combining the above two estimates, applying (\ref{2.17.3.24}) again, and using Proposition \ref{1steng} to treat the lower order term, we obtain 
\begin{align*}
\|\sP[S, \Omega \Phi]\|_{L^2_\Sigma}\les  \l t\r^{-2}W_2[S\Omega\Phi]^\f12(t)+ \l t\r^{-2}\log \l t\r^{\f12\M+1}\La_0+\l t\r^{-2+2\delta}\log \l t\r^\f12\Delta_0^\frac{5}{4}.
\end{align*}

We then combine the above estimates with (\ref{6.24.9.23}) and (\ref{8.17.4.23}) to obtain
\begin{align*}
\|\Er_2(\Phi, S\Omega)\|_{L^2_\Sigma}&\les  \l t\r^{-2}W_2[S\Omega\Phi]^\f12(t)+ \l t\r^{-2}\log \l t\r^{\f12\M+1}\La_0+\l t\r^{-\frac{7}{4}+2\delta}\Delta_0^\frac{5}{4}.
\end{align*}

{\bf Case $X_2=S, \, X_1=\Omega$.} 
Applying (\ref{9.2.6.23}) to $f=S\Phi$ and using (\ref{10.10.2.23}), we deduce  
\begin{align*}
\|\sP[{}\rp{a}\Omega, S\Phi]\|_{L^2_\Sigma}\les \l t\r^{-\frac{7}{4}+\delta}\Delta_0^\frac{5}{4}+\l t\r^{-2}(\log \l t\r)^{\f12\M}(\La_0+\Delta_0^\frac{5}{4}).
\end{align*}
Next applying (\ref{8.7.1.21}) to $f=\Phi$ and $X_2=\Omega$,  by using  $X_2X_1=\Omega^2$ in (\ref{8.25.2.21})  and also  (\ref{10.10.2.23}), we obtain
\begin{align*}
\|\Omega(\sP[S, \Phi])\|_{L^2_\Sigma}&\les \l t\r^{-2}\log \l t\r^{\f12\M+1}\La_0+\l t\r^{-2+2\delta}\log \l t\r^\f12\Delta_0^\frac{5}{4}+\l t\r^{-2}W_2[\Omega S\Phi](t)^\f12\\
&+\|\bb^{-\f12}\Omega\wt{\tr\chi}\|_{L^2_u L_\omega^2}.
\end{align*}
Next we use Lemma \ref{comp}, (\ref{7.21.2.22+}), (\ref{2.17.3.24}), Proposition \ref{7.15.5.22} and (\ref{10.10.2.23}) to treat the last term,
\begin{align*}
\|\Omega\wt{\tr\chi}\|_{L_u^2 L_\omega^2}&\les \l t\r^{-2}(\l t\r^{2\delta}\Delta_0^\frac{3}{2}+\log \l t\r^{\f12\M+1}\La_0)+\l t\r^{-2}\int_0^t \l t'\r^{-1}W_2[\Omega^2\varrho]^\f12(t') dt'\\
&\les \l t\r^{-2}(\l t\r^{2\delta}\Delta_0^\frac{5}{4}+\log \l t\r^{\f12\M+1}\La_0),
\end{align*} 
which will be further improved  in (\ref{2.18.2.24}). 

In the same way as for the other case, we then conclude
\begin{align*}
\|\Er_2(\Phi, \Omega S)\|_{L^2_\Sigma}\les \l t\r^{-2}\log \l t\r^{\f12\M+1}\La_0+\l t\r^{-\frac{7}{4}+2\delta}\Delta_0^\frac{5}{4}+\l t\r^{-2}W_2[\Omega S\Phi]^\f12(t).
\end{align*}
Then using the above estimates of $\Er_2(\Phi, S\Omega)$ and $\Er_2(\Phi, \Omega S)$, we derive from Proposition \ref{MA2} for $\Psi=S\Omega\Phi, \Omega S\Phi$
\begin{align*}
WFIL_2[\Psi](\D_{u_1}^{t_1})&\les \La_0^2+\Delta_0^\frac{5}{2}+\int_0^{t_1} \Big(\sup_{\Sigma_{t'}}(-[\Lb\varrho]_{-})+\l t'\r^{-1-\delta}\Big)WL_2[\Psi](t') dt'\\
&+\int_{\D_{u_1}^{t_1}} \big|\Box_\bg \Psi(L\Psi+ h \Psi)\big|\tir^2 \bb d\mu_{\ga} du dt+\int_{u_1}^{u_*} F_2[\Psi](\H_u^{t_1}) du\\
&\les \La_0^2+\Delta_0^\frac{5}{2}+ \int_0^{t_1}\l t\r^{-2}\log \l t\r W_2[\Psi](t) dt+\int_{u_1}^{u_*}F_2[\Psi](\H_u^{t_1}) du\\
&+\int_0^{t_1} \Big(\sup_{\Sigma_{t'}}(-[\Lb\varrho]_{-})+\l t'\r^{-1-\delta}\Big)WL_2[\Psi](t') dt'.
\end{align*}
Thus by Gronwall's inequality we can obtain (\ref{8.25.2.21}) for $X_2 X_1=S \Omega, \Omega S$. Substituting this result to the above estimates of $\|\Er_2(\Phi, S\Omega), \Er_2(\Phi, \Omega S)\|_{L^2_\Sigma}$  and (\ref{7.9.4.22}), we further conclude for $X^2=S\Omega$ or $\Omega S$,
\begin{align}
&\|\Er_2(\Phi, X^2), \Box_\bg X^2\Phi\|_{L^2_\Sigma}\les \l t\r^{-2}(\log \l t\r)^{\f12\M+1}\La_0+\l t\r^{-\frac{7}{4}+2\delta}\Delta_0^\frac{5}{4},\label{8.3.4.21}\\
&\|X^2 \Box_\bg \Phi\|_{L^2_\Sigma}\les\l t\r^{-2}(\log \l t\r)^{\f12\M} (\La_0+\Delta_0^\frac{5}{4})+\l t\r^{-2+2\delta}\Delta_0^\frac{3}{2}.\label{8.3.2.21}
\end{align}

{\bf (3) $X_2=X_1=S$.}
Using the estimate for $W_2[\Omega^2\Phi](t)$ in (\ref{2.17.3.24}) and $W_2[\Omega S\Phi](t)$ in (\ref{8.25.2.21}), with $X_2=S, f=S\Phi$ in (\ref{8.7.1.21}), (\ref{2.18.1.24}), (\ref{8.21.2.21}) and (\ref{5.7.1.24}), we bound
\begin{align*} 
\| S(\sP[S, \Phi]), \sP[S, S\Phi]\|_{L^2_\Sigma}\les \l t\r^{-2}W_2[S^2\Phi]^\f12(t)+\l t\r^{-2 }(\log \l t\r^{\f12\M+1}\La_0+\l t\r^{2\delta}\log \l t\r^\f12\Delta_0^\frac{5}{4}).
\end{align*}
Combining the above estimates with (\ref{6.24.9.23}) and (\ref{8.17.4.23}) implies
\begin{align}\label{8.29.6.21}
\|\Er_2(\Phi, S^2)\|_{L^2_\Sigma}\les  \l t\r^{-2}(\log \l t\r^{\f12\M+1}\La_0+\l t\r^{2\delta}\log \l t\r^\f12\Delta_0^\frac{5}{4})+\l t\r^{-2}W_2[S^2\Phi]^\f12(t).
\end{align}
Recall from (\ref{7.9.4.22}) that
\begin{equation*}
\|S^2\Box_\bg \Phi\|_{L^2_\Sigma}\les \l t\r^{-2}W_2[S^2\Phi]^\f12(t)+\l t\r^{-2}(\log \l t\r)^{\f12\M}(\La_0+\Delta_0^\frac{5}{4}).
\end{equation*}
 Applying Proposition \ref{MA2} to $\Psi=S^2\Phi$, we derive
\begin{align*}
WFIL_2[\Psi](\D_{u_1}^{t_1})&\les \La_0^2+\Delta_0^\frac{5}{2}+\int_0^{t_1} (\sup_{\Sigma_{t'}}(-[\Lb\varrho]_{-})+\l t'\r^{-1-\delta})WL_2[\Psi](t') dt'\\
&+\int_0^{t_1}  \l t\r^{-2}\log \l t\r W_2[\Psi](t)+\int_{u_1}^{u_*} F_2[\Psi](\H_u^{t_1}) du.
\end{align*}
 Using Gronwall's inequality, we conclude the case $X_1=X_2=S$ in (\ref{8.25.2.21}).

We further conclude
\begin{align}\label{2.18.3.24}
&\|\Er_2(\Phi, S^2), \Box_\bg S^2\Phi\|_{L^2_\Sigma}\les   \l t\r^{-2}(\log \l t\r^{\f12\M+1}\La_0+\l t\r^{2\delta}\log \l t\r^\f12\Delta_0^\frac{5}{4})
\end{align}
and \begin{equation}\label{5.10.1.24}
\|S^2\Box_\bg \Phi\|_{L^2_\Sigma}\les\l t\r^{-2}(\log \l t\r)^{\f12\M}(\La_0+\Delta_0^\frac{5}{4}).
\end{equation}
(\ref{8.29.9.21}) can be obtained by using (\ref{8.25.2.21}), Lemma \ref{6.30.4.23}, (\ref{7.17.7.21}), Proposition \ref{1steng} and Lemma \ref{comp}. (\ref{8.30.3.21}) can be obtained by using (\ref{8.29.9.21}) to refine (\ref{7.9.4.22}).
The first estimate in (\ref{8.30.3.21+}) follows by combining (\ref{8.13.2.21}) and (\ref{2.17.4.24}). The second estimate follows by combining (\ref{8.3.4.21}) and (\ref{8.3.2.21})
Finally, we obtain the last estimate in (\ref{8.30.3.21+}) by combining (\ref{2.18.3.24}) and (\ref{5.10.1.24}).

\end{proof}

 By direct checking in view of (\ref{5.13.10.21}), we have the following formulas.
\begin{lemma}\label{3.25.2.24}
Let $\a> -1$ and $\a\neq0$. For a scalar $f$, there holds the commutation formula
\begin{align}\label{2.23.2.24}
[\Omega^n, \bb^{-\a}\Lb]f&=\sum_{i=0}^{n-1}\{\frac{\a+1}{\a}\Omega^{i}(\Omega(\bb^{-\a})\Lb\Omega^{n-1-i} f)-\frac{1}{\a}\Omega^{i}(\Omega(\bb^{-\a})L\Omega^{n-1-i} f)\nn\\
&+\Omega^i(\bb^{-\a}\pioh_{A\Lb}\sn_A\Omega^{n-i-1}f)\}.
\end{align}
If $\a=0$, we have
\begin{equation}\label{5.10.2.24}
[\Omega^n, \Lb]f=\sum_{i=0}^{n-1}\Omega^i\Big((-\Omega \log \bb \Lb +\Omega \log \bb L -{}\rp{\Omega}\pih_{\Lb A}\sn )\Omega^{n-1-i} f\Big).
\end{equation}
In both cases, when $i=n-1$ for the first terms on the right-hand side, the operator $\Omega^n={}\rp{a_1}\Omega{}\rp{a_2}\Omega\cdots{}\rp{a_n}\Omega$, with $a_i\in \{1,2,3\}, i=1\cdots n$, appearing in $\Omega^n(\bb^{-\a})\Lb f$ or $\Omega^n \log \bb \Lb f$ is exactly the same operator $\Omega^n$ on their left-hand side.   
\end{lemma}

\begin{lemma}\label{6.29.2.24}
Let $0\le a\le 3, a\in {\mathbb Z}$. 
\begin{align}\label{2.24.1.24}
L\Omega^a \log \bb&=-\f12\wp\Lb\Omega^a\varrho+ \sum_{i=0}^{a-1}\Omega^i(\pioh_{A L}\sn_A\Omega^{a-i-1}(\log \bb))+\Omega^a([L\Phi])\nn\\
&-\f12\wp[\Omega^a, \Lb]\varrho,
\end{align}
\begin{align}\label{2.23.3.24}
L\Omega^a (\bb^{-\a})&=\f12\a \wp\bb^{-\a}\Lb\Omega^a\varrho+\sum_{i=0}^{a-1}\Omega^i(\pioh_{A L}\sn_A\Omega^{a-i-1}(\bb^{-\a}))+\Omega^a(\bb^{-\a}[L\Phi])\nn\\
&+\f12\a\wp[\Omega^a, \bb^{-\a}\Lb]\varrho, \quad \a>-1, \a\neq 0.
\end{align}
\end{lemma}
\begin{proof}
The case that $a=0$ can be directly obtained by using (\ref{lb}) and (\ref{5.13.10.21}). If $\a\neq 0$,
\begin{align*}
L\Omega^a(\bb^{-\a})&=[L, \Omega^a](\bb^{-\a})+\Omega^a L(\bb^{-\a})\\
&=\sum_{i=0}^{a-1}\Omega^i\big(\pioh_{A L}\sn_A\Omega^{a-1-i}(\bb^{-\a})\big)+\a\Omega^a(\bb^{-\a}k_{\bN\bN}).
\end{align*}
For the second term on the right-hand side, we derive by using (\ref{3.22.1.21})
\begin{align*}
\Omega^a(\bb^{-\a}k_{\bN\bN})&=\f12\wp\Omega^a(\bb^{-\a}\Lb\varrho)+\Omega^a(\bb^{-\a}[L\Phi])\\
 &=\f12 \wp\bb^{-\a}\Lb \Omega^a \varrho+\f12\wp[\Omega^a, \bb^{-\a}\Lb]\varrho+\Omega^a(\bb^{-\a}[L\Phi]).
\end{align*}
Combining the above two calculations gives (\ref{2.23.3.24}). (\ref{2.24.1.24}) can be obtained similarly.
\end{proof}
\begin{proof}[Proof of (\ref{2.18.2.24})-(\ref{3.12.4.24})] We first prove (\ref{2.20.2.24}). 
The case $a=1$ is much easier. We can obtain by using (\ref{2.23.2.24})-(\ref{2.23.3.24}), (\ref{5.21.1.21}) and Proposition \ref{1steng} 
\begin{equation}\label{2.24.2.24}
\|\Omega(\bb^{-\a}), \Omega\log \bb\|_{L^2_u L_\omega^2}\les \log \l t\r(\La_0+\Delta_0^\frac{5}{4}), \quad\a\ge-\f12.
\end{equation}
Using (\ref{2.23.2.24})-(\ref{2.23.3.24}), we write
\begin{align*}
L \sta{cb}{\Omega^2}(\bb^{-\a})&=\wp(\f12 (\a+1)\Lb \varrho+L\varrho) \sta{cb}{\Omega^2}(\bb^{-\a})+\sum_{Y=L, \Lb}\Omega(\bb^{-\a})(Y \Omega \varrho+\Omega Y \varrho)\\
\displaybreak[0]
&+\f12 \wp\a\bb^{-\a}\Lb\Omega^2\varrho+\sum_{i=0}^1\Omega^i(\bb^{-\a}\pioh_{A\Lb}\sn_A\Omega^{1-i}\varrho)\\
&+\sum_{i=0}^1\Omega^i(\pioh_{A L}\sn_A\Omega^{1-i}(\bb^{-\a}))+\Omega^2(\bb^{-\a}[L\Phi]),
\end{align*}
where we only need to keep the sign of the first term on the right-hand size precise in the above, the values of the coefficients of all other terms do not affect our estimates.  Multiplying the above identity by $\sta{cb}{\Omega^2}(\bb^{-\a})$, using (\ref{6.5.1.21}), ignoring the negative part of the first term on the right-hand side, by transport Lemma, we deduce for $\a\ge -\f12$
\begin{align}\label{2.24.3.24}
&\|\Omega^2(\bb^{-\a})(t)\|_{L_u^2 L_\omega^2}\nn\\
&\les \La_0+\int_0^t\|\bb^{-\a}\Lb\Omega^2\varrho\|_{L^2_u L_\omega^2}+\int_0^t \sum_{Y=\Lb, L}\|\Omega(\bb^{-\a})(Y \Omega \varrho+\Omega Y \varrho)\|_{L_u^2 L_\omega^2}\nn\\
\displaybreak[0]
&+ \int_0^t \sum_{i=0}^1\|\Omega^i(\bb^{-\a}\pioh_{A\Lb}\sn_A\Omega^{1-i}\varrho)\|_{L_u^2 L_\omega^2}+\int_0^t \sum_{i=0}^1\|\Omega^i(\pioh_{A L}\sn_A\Omega^{1-i}(\bb^{-\a}))\|_{L^2_u L_\omega^2}\nn\\
&+\int_0^t \|\Omega^2(\bb^{-\a}[L\Phi])\|_{L_u^2 L_\omega^2}.
\end{align}
By using (\ref{5.21.1.21}), (\ref{3.25.1.22}), (\ref{3.6.2.21}) and (\ref{3.11.3.21}), we deduce
\begin{align}\label{5.11.1.24}
\sum_{i=0}^1\|\Omega^i(\bb^{-\a}\pioh_{A\Lb}\sn_A\Omega^{1-i}\varrho)\|_{L_u^2 L_\omega^2}\les \l t\r^{-3+3\delta}\log \l t\r\Delta_0^\frac{3}{2}.
\end{align}
\begin{align}
\|\Omega^i(\pioh_{A L}&\sn_A\Omega^{1-i}(\bb^{-\a}))\|_{L^2_u L_\omega^2}\nn\\
&\les \l t\r^{-\frac{7}{4}+\delta}\Delta_0^\f12\|\Omega^2(\bb^{-\a})\|_{L_u^2 L_\omega^2}+\l t\r^{-2+\delta}\Delta_0\|\Omega(\bb^{-\a})\|_{L_u^\infty L^\infty_\omega}\nn\\
&\les \l t\r^{-\frac{7}{4}+\delta}\Delta_0^\f12\|\Omega^2(\bb^{-\a})\|_{L_u^2 L_\omega^2}+\l t\r^{-2+3\delta}\Delta_0^2.\label{5.11.2.24}
\end{align}
Substituting the above two estimates to (\ref{2.24.3.24}), also using Lemma \ref{5.13.11.21} (5), (\ref{11.11.2.23}), (\ref{3.11.3.21}), (\ref{2.24.2.24}) and (\ref{8.29.9.21}) to treat the remaining terms, we infer
\begin{align*}
\|\Omega^2(\bb^{-\a})(t)\|_{L_u^2 L_\omega^2}&\les\log \l t\r(\La_0+\Delta_0^\frac{5}{4})+\int_0^t \|\Omega(\bb^{-\a})\|_{L_u^2 L_\omega^4}\l t\r^{-1}\Delta_0^\f12\\
&+\int_0^t \l t\r^{-\frac{7}{4}+\delta}\Delta_0^\f12\|\Omega^2(\bb^{-\a})\|_{L_u^2 L_\omega^2}\\
&\les\log \l t\r(\La_0+\Delta_0^\frac{5}{4})+\int_0^t \|\Omega^2(\bb^{-\a})\|^\f12_{L_u^2 L_\omega^2}\log \l t\r^{\frac{1}{2}}\l t\r^{-1}\Delta_0^\f12 (\La_0+\Delta_0^\frac{5}{4})^\f12\\
&+\int_0^t \l t\r^{-\frac{7}{4}+\delta}\Delta_0^\f12\|\Omega^2(\bb^{-\a})\|_{L_u^2 L_\omega^2}.
\end{align*}
Hence we conclude by Gronwall's inequality that
\begin{align*}
\|\Omega^2(\bb^{-\a})(t)\|_{L_u^2 L_\omega^2}\les  \log \l t\r(\La_0+\Delta_0^\frac{5}{4})+\log \l t\r^\frac{7}{2}\Delta_0(\La_0+\Delta_0^\frac{5}{4}).
\end{align*}
We can similarly obtain the same estimate for $\Omega^2\log \bb$. Moreover,
using (\ref{2.24.2.24}), the above estimate, (\ref{zeh}), (\ref{8.25.2.21}), (\ref{8.21.4.21}) and (\ref{8.29.9.21}), in view of (\ref{2.23.2.24})-(\ref{2.23.3.24}), we can also obtain 
\begin{align*}
\|L\Omega^a (\log \bb,\bb^{-\a})\|_{L^2_u L_\omega^2}\les \l t\r^{-1}\log \l t\r(\La_0+\log \l t\r^{\frac{5}{2}(a-1)}\Delta_0^\frac{5}{4})   
\end{align*}
where we used the obtained estimates (\ref{5.11.1.24}) and (\ref{5.11.2.24}) for the quadratic terms related to $\pioh_{\Lb A}$ and $\pioh_{LA}$.  Using (\ref{7.17.7.21}), we can bound for $f=\log \bb,\bb^{-\a}$ that
\begin{align*}
\|[L, \Omega^a]f \|_{L_u^2 L_\omega^2}&\les \l t\r^{-1}\|\Omega^{1+\le 1} f\|_{L_u^2 L_\omega^2}+\l t\r^{-\frac{3}{4}+\delta}\Delta_0^\f12 \|\sn f\|_{L_u^2 L_\omega^4}\\
&\les  \l t\r^{-1}\log \l t\r(\La_0+\log \l t\r^{\frac{5}{2}(a-1)}\Delta_0^\frac{5}{4}).
 \end{align*}
 Combining the above two estimates, we obtain the remaining estimates in (\ref{2.20.2.24}). 

 (\ref{2.14.1.24}) follows as a consequence of (\ref{2.20.2.24}) in view of (\ref{lb}) and (\ref{8.29.9.21}).  
 To see (\ref{2.22.3.24}), we employ (\ref{5.13.10.21}), (\ref{1.27.5.24}), (\ref{8.25.2.21}) to derive
 \begin{align*}
 \|\Omega\Lb \Omega \varrho\|_{L^2_\Sigma}&\les \|[\Lb, \Omega]\Omega \varrho\|_{L^2_\Sigma}+\|\Lb \Omega^2\varrho\|_{L^2_\Sigma}\\
 &\les \|\bb^{-1}\Omega \bb\bN \Omega\varrho\|_{L^2_\Sigma}+\|{}\rp{a}\pih_{\Lb A}\c \sn \Omega\varrho\|_{L^2_\Sigma}+\La_0+\Delta_0^\frac{5}{4}\\
 &\les \log \l t\r\Delta_0\|\Omega\bN \Omega\varrho\|_{L^2_\Sigma}^\f12 \|\bN\Omega\varrho\|_{L^2_\Sigma}^\f12+\La_0+\Delta_0^\frac{5}{4}
 \end{align*}
 where to obtain the last line we used (\ref{10.10.2.23}) and (\ref{5.21.1.21}) to treat the second term on the right-hand side.  We then use (\ref{8.21.4.21}), (\ref{8.29.9.21}), $2\bN=L-\Lb$ to get the leading part of the first term on the right-hand side absorbed by the left-hand side with the help of Cauchy-Schwartz  inequality,
 \begin{equation*}
 \|\Omega\Lb \Omega \varrho\|_{L^2_\Sigma}\les \log\l t\r^2\Delta_0^3+\La_0+\Delta_0^\frac{5}{4}.
 \end{equation*}
  Thus we can conclude (\ref{2.22.3.24}).

Using (\ref{8.29.9.21}), (\ref{2.14.1.24}), Proposition \ref{6.24.10.23},  (\ref{8.23.2.23}), (\ref{3.11.3.21}) and (\ref{3.6.2.21}), we obtain the  improved estimate 
\begin{equation*}
\|\Omega^{1+\le 1}([L\Phi][\Lb\Phi], \N(\Phi,\bp\Phi))\|_{L^2_\Sigma}\les \l t\r^{-2}(\log \l t\r^{\f12\M}\La_0+\log \l t\r^5\l t\r^\delta\Delta_0^\frac{3}{2}).
\end{equation*}
(\ref{2.18.2.24}) then follows by using the above estimate, (\ref{10.11.3.23}), (\ref{7.21.2.22+}), (\ref{8.25.2.23}) and (\ref{8.29.9.21}). Here we used the fact that $\M\ge 8$.  
(\ref{3.10.9.24}) is a consequence of (\ref{9.8.4.22}), (\ref{3.6.2.21}), (\ref{2.18.2.24}) and (\ref{7.13.5.22}).  (\ref{3.5.5.24}) is obtained by using (\ref{7.9.6.22}) and (\ref{2.18.2.24}).
 
 Next we prove (\ref{2.19.1.24}). Using (\ref{8.26.1.23}), (\ref{8.29.9.21}), (\ref{2.18.2.24}) and Proposition \ref{1steng}
\begin{align*}
\Omega \Lb[L\Phi]&=\Omega\big((\tir^{-1}+\fB)\fB+\bA_b\fB+\sD\varrho+L[L\Phi]+\ud\bA\bA_{g,1}+\N(\Phi, \bp\Phi)\big)\\
&= O(\l t\r^{-1}\log \l t\r(\La_0+\Delta_0^\frac{5}{4}))_{L^2_\Sigma}.
\end{align*}
For the third estimate, using (\ref{5.13.10.21}), (\ref{6.22.1.21}), (\ref{2.20.2.24}), (\ref{5.21.1.21}) and (\ref{10.10.2.23}), we bound the commutator below
\begin{align*}
\|[\Lb, \Omega][L\Phi]\|_{L^2_\Sigma}&\les \|\tir \ud\bA\bN[L\Phi]\|_{L^2_\Sigma}+\|\pioh_{A\Lb}\sn[L\Phi]\|_{L^2_\Sigma}\\
&\les \|\fB\ud \bA\|_{L^2_\Sigma}+\l t\r^{-\frac{3}{4}+\delta}\Delta_0^\f12 \|\sn[L\Phi]\|_{L^2_\Sigma}\\
&\les\l t\r^{-1} \log \l t\r(\La_0+\Delta_0^\frac{5}{4}).
\end{align*}
The other estimates follow similarly, also with the help of (\ref{7.17.6.21}).
 
Finally, similar to (\ref{8.2.1.22}), for  (\ref{3.12.4.24}),  the case of $X=\Omega$ follows by applying Lemma \ref{2.9.3.23}, (\ref{10.10.2.23}), (\ref{8.29.9.21}) and using (\ref{5.21.1.21}), (\ref{3.6.2.21}) (\ref{2.20.2.24}) and (\ref{1.27.5.24}); the case of $X=S$ follows in view of (\ref{8.23.1.23}). 
\end{proof}

Next we give a preliminary result for proving (\ref{2.13.3.24}) and (\ref{3.9.8.24}).
\begin{lemma}
\begin{align}
\Box_\bg \bT\varrho&=\Lb^2\varrho \fB+\fB\Box_\bg \varrho+O(\tir^{-1})\fB^2+O\Big(\l t\r^{-3+2\delta}(\Delta_0^\frac{5}{4}+\La_0)\Big)_{L^2_\Sigma} \label{2.18.9.24}\\
\Box_\bg\bT\varrho&=\Lb^2\varrho\fB+O(\l t\r^{-1})\fB^2+O(\l t\r^{-\frac{15}{4}+\delta}\Delta_0^\f12)\label{3.9.6.24}\\
\Omega \Box_\bg \bT \varrho&=\Lb \Omega \bT\varrho \c \fB +O\Big(\l t\r^{-1}\log\l t\r^\frac{9}{4}(\La_0+ \Delta_0^\frac{5}{4})\Big)_{L^2_\Sigma}.\label{2.16.2.24}
\end{align}
\end{lemma}
\begin{proof}
Using Proposition \ref{7.15.5.22},  we derive in view of (\ref{2.18.10.24}) that
\begin{align}\label{2.19.2.24}
&(\div_g k)_A=O(\l t\r^{-1+\delta}\Delta_0)_{L^2_\Sigma}, O(\l t\r^{-2+\delta}\Delta_0),\\
& (\div_g k)_\bN=\bN \fB+\fB^2+\tr\theta \fB+O(\l t\r^{-3+2\delta}\Delta_0^2). \nn
\end{align}
 Decomposing by using (\ref{3.6.2.21}) and \cite[P86 (4.2.9)]{CK},
 \begin{align}
 \Delta_\bg \varrho&=\bN\bN \varrho+\tr\theta\fB+\sD\varrho+\bA_{g,1}\ud \bA\label{5.11.4.24}\\
 &=\bN\bN \varrho+\tr\theta\fB+\sD\varrho+O(\l t\r^{-3+2\delta}\Delta_0^2)\nn
  \end{align}
we use (\ref{2.18.11.24}), the above three estimates, Proposition \ref{7.15.5.22} and (\ref{10.10.2.23}) to derive that
\begin{align*}
\Box_\bg \bT\varrho=\bT\Box_\bg \varrho+\fB \Box_\bg \varrho+\fB(\bN \fB +\fB^2+\tr\theta\fB)+O(\l t\r^{-3+2\delta}(\Delta_0^\frac{5}{4}+\La_0))_{L^2_\Sigma},
\end{align*}
and we also have the rough bound
\begin{equation*}
\Box_\bg \bT\varrho-(\bT\Box_\bg \varrho+\fB \Box_\bg \varrho+\fB(\bN \fB +O(\l t\r^{-1})\fB))=O(\l t\r^{-4+\frac{1}{4}+\delta}\Delta_0^\f12).
\end{equation*}
 We can similarly derive as in (\ref{12.21.4.23}) by using Proposition \ref{7.15.5.22} 
\begin{align*}
\bT \Box_\bg \varrho-(\Lb^2\varrho[L\Phi]+O(\tir^{-1})\fB^2)=O(\l t\r^{-3+2\delta}\Delta_0^\frac{3}{2})_{L^2_\Sigma}, O(\l t\r^{-4+2\delta}\Delta_0^\frac{3}{2}). 
\end{align*}
Thus we conclude (\ref{2.18.9.24}) and (\ref{3.9.6.24}).

To see (\ref{2.16.2.24}), we first differentiate (\ref{2.18.10.24}) with the help of Lemma \ref{error_prod}, Proposition \ref{7.15.5.22}, (\ref{7.13.5.22})  and (5) in Lemma \ref{5.13.11.21} to obtain
\begin{align}\label{2.19.3.24}
\sn_\Omega (\div_g k)_A&=\Omega(\fB \bA_{g,1})+\sn_\Omega \sn\fB+\sn_\Omega \sn \bAn+\sn_\Omega((\tr\theta+\bA_g) \bA_{g,1})\nn\\
&+\sn_\Omega \sn_\bN \bA_{g,1}+\sn_\Omega\big((\fB+\bAn)\ud\bA\big)=O(\l t\r^{-1+\delta}\Delta_0)_{L^2_\Sigma}\\
\sn_\Omega\big((\div_g k)_\bN\big)&=\Omega(\fB^2+[\sn\Phi]^2+\bN \fB+\sn\bA_{g,1}+\tr\theta\fB+(\bA_g+\ud\bA)\bA_{g,1})\nn\\
&=(\fB+\tr\chi)\Omega\fB+\fB\Omega\bA_b+\Omega(\Lb^2\varrho+\Lb[L\Phi]+L\fB)+\sn_\Omega \sn\bA_{g,1}\nn\\
&+O(\l t\r^{-3+2\delta}\Delta_0^2)_{L^2_u L_\omega^2}\nn\\
&=\Omega(\Lb\bT\varrho)+O(\l t\r^{-1}\log\l t\r(\Delta_0^\frac{5}{4}+\La_0))_{L^2_\Sigma}\nn
\end{align}
where we applied (\ref{9.20.4.22}), (\ref{8.29.9.21}), (\ref{2.18.2.24}) and (\ref{2.19.1.24}) to conclude the second estimate.

Differentiating (\ref{2.18.11.24}) and using (\ref{5.11.4.24}), symbolically, we write
\begin{align}\label{3.7.4.24}
\begin{split}
\Omega \Box_\bg\bT\varrho&=\Omega \bT\Box_\bg \varrho+\Omega(\fB \Box_\bg \varrho)+\Omega(\fB (\div_g k)_\bN)+\Omega(\sn\varrho(\div_g k)_A)+\Omega(\fB\ud \bA \sn\varrho)\\
&+\Omega\Big(\fB(\sn^2\varrho+\theta \fB+\fB^2+\bN \fB)\Big)+\Omega\Big(\bA_{g,1}(\sn\fB+\theta\sn\varrho)\Big).
\end{split}
\end{align}
Here we remark that $\Omega(\fB \bN \fB), \Omega \fB \l t\r^{-1}\fB$ appear in both the third term on the right-hand side and the last line.

Using (\ref{2.19.2.24}), (\ref{2.19.3.24}) and (\ref{9.20.4.22}), we deduce
\begin{align*}
\Omega(\fB (\div_g k)_\bN)&=\fB\c (\Lb\Omega\bT\varrho+[\Omega,\Lb]\bT\varrho)+\fB\c O(\l t\r^{-1}\log\l t\r(\Delta_0^\frac{5}{4}+\Delta_0))_{L^2_\Sigma}\\
&+\Omega \fB(\bN\fB+O(\l t\r^{-1}) \fB)+O\Big(\l t\r^{-2}\log\l t\r(\Delta_0^\frac{5}{4}+\Delta_0)\Big)_{L^2_\Sigma}.
\end{align*}

By using (\ref{1.29.2.22}), (\ref{9.20.4.22}) and (\ref{2.14.1.24}), and writing $\bN\fB=\Lb^2\varrho+O(\tir^{-1}\fB)$, 
\begin{align}
\|\bb^\frac{5}{2} \Omega \fB \Lb^2 \varrho \|_{L^2_\Sigma}&\les \|\Omega\fB\|_{L_u^2 L_\omega^4}
\les \|\Omega^2 \fB\|^\f12_{L_u^2 L_\omega^2}\|\Omega\fB\|^\f12_{L_u^2 L_\omega^2}\nn\\
&\les\l t\r^{-1}\log\l t\r(\La_0+(\log \l t\r)^\frac{5}{2} \Delta_0^\frac{5}{4})^\f12(\La_0+ \Delta_0^\frac{5}{4})^\f12. \label{2.19.4.24}
\end{align}
Hence, also using (\ref{2.15.3.24}), (\ref{1.27.5.24}) and again by (\ref{9.20.4.22}) 
\begin{align*}
&\Omega(\fB (\div_g k)_\bN)\\
&=\fB \Lb\Omega\bT\varrho +\fB \bb^{-2}\Omega\bb O(\l t\r^{-1})_{L_\omega^4}+O(\l t\r^{-1}\log\l t\r^\frac{9}{4}(\La_0+ \Delta_0^\frac{5}{4}))_{L^2_\Sigma}\\
&=\fB\Lb \Omega\bT\varrho+\fB\bb^{-2} O(\l t\r^{-1}\log \l t\r\Delta_0)_{L_\omega^2}+O(\l t\r^{-1}\log\l t\r^\frac{9}{4}(\La_0+ \Delta_0^\frac{5}{4}))_{L^2_\Sigma}\\
&=\fB\Lb \Omega\bT\varrho+O(\l t\r^{-1}\log\l t\r^\frac{9}{4}(\La_0+ \Delta_0^\frac{5}{4}))_{L^2_\Sigma}
\end{align*}
where we treated the second term on the right-hand side by using (\ref{12.19.1.23}) to obtain the last line.  

Using (\ref{2.19.3.24}), (\ref{2.19.2.24}), (\ref{10.10.2.23}) and (\ref{3.6.2.21}) we have
\begin{align*}
\Omega\big((\div_g k)_A\sn_A \varrho\big)=O(\l t\r^{-3+2\delta}\Delta_0^2)_{L^2_\Sigma}.
\end{align*}
It follows by using Lemma \ref{5.13.11.21} (5), (\ref{3.6.2.21}) and (\ref{10.10.2.23}), 
\begin{align*}
\|\bb^\f12\Omega(\fB\ud\bA\bA_{g,1})\|_{L^2_u L_\omega^2}\les \l t\r^{-4+2\delta}\Delta_0^2.
\end{align*}
Similarly, using (\ref{8.29.9.21}), (\ref{3.11.3.21}) and the estimate of $\Omega \fB$ in (\ref{2.19.4.24})
\begin{align*}
\|\bb^\f12\Omega(\fB\sn^2\varrho)\|_{L^2_u L_\omega^2}&\les \|\Omega \fB\|_{L_u^2 L_\omega^4}\|\bb^\f12\sn^2\varrho\|_{L_\omega^4}+\l t\r^{-1}\|\Omega\sn^2 \varrho\|_{L_u^2 L_\omega^2}\\
&\les \l t\r^{-4+\delta}\log\l t\r^\frac{3}{2}(\La_0+(\log \l t\r)^\frac{\M}{2} \Delta_0^\frac{5}{4}).
\end{align*}
Using (\ref{9.20.4.22}) and (\ref{2.18.2.24}), we derive
\begin{align*}
\Omega(\fB^2\theta)=(\bb^{-1}\Omega \fB+\bb^{-2}\Omega \theta) O(\l t\r^{-2})=O(\l t\r^{-2}\log \l t\r(\Delta_0^\frac{5}{4}+\La_0))_{L^2_\Sigma}. 
\end{align*}

Using (\ref{3.6.2.21}), (\ref{LbBA2}) and (\ref{L2BA2}), we have
 \begin{align*}
\|\bb^\f12\Omega(\bA_{g,1}\sn \fB)\|_{L^2_u L_\omega^2}\les \l t\r^{-4+2\delta}\Delta_0^2.
\end{align*}
$\Omega(\bA_{g,1}\theta \sn\varrho)$ is a much better term compared with the above one, which is negligible. 

 Next we use (\ref{9.20.4.22}), (\ref{8.26.4.21}) and (\ref{4.3.3.21}) to bound
\begin{align*}
\|\Omega(\fB\Box_\bg \varrho)\|_{L^2_\Sigma}&\les \|\Omega \fB\|_{L^2_u L_\omega^2}\|\bb^\f12\tir\Box_\bg \varrho\|_{L_\omega^\infty}+\|\bb^{-\f12}\Omega\Box_\bg \varrho\|_{L^2_u L_\omega^2}\\
&\les\l t\r^{-3+\delta}\log\l t\r^{\f12\M+1}(\La_0+\Delta_0^\frac{5}{4}).
\end{align*}
Similar to (\ref{2.10.7.24}), by using (\ref{3.11.3.21}), (\ref{6.22.1.21}), (\ref{2.15.3.24}), (\ref{9.20.4.22}), (\ref{2.19.1.24}), (\ref{1.29.2.22}), (\ref{LbBA2}) and (\ref{8.29.9.21}), also writing $\bT\fB=\Lb^2\varrho+\Lb[L\Phi]+L\fB$ we bound
\begin{align*}
\Omega\bT \Box_\bg \varrho&=\Omega \bT [\Lb \Phi][L\Phi]+\bT[\Lb\Phi]\Omega[L\Phi]+\Omega[\Lb\Phi]\bT[L\Phi]+[\Lb\Phi]\Omega\bT[L\Phi]+\Omega\bT([\sn\Phi]^2)+\Omega\bT(|\eh|^2)\\
&=\Lb \Omega \bT \varrho [L\Phi]+O(\l t\r^{-2}\log\l t\r^{\f12\M+1}(\Delta_0^\frac{5}{4}+\Delta_0))_{L^2_\Sigma}.
\end{align*}
where for the last line, we used (\ref{3.6.2.21}), (\ref{10.10.2.23}) and (\ref{1.27.5.24}). Summarizing the above calculations, we conclude (\ref{2.16.2.24}).
\end{proof}
\begin{proof}[Proof of (\ref{2.13.3.24}) and (\ref{3.9.8.24})]
We first derive the standard energy inequality by Proposition \ref{10.10.3.22}
   \begin{align*}
   E[\Omega\bT\varrho](t_1)+F_0[\Omega\bT\varrho]&(\H_{u_1}^{t_1})\les\int_{\D_{u_1}^{t_1}} |\Box_\bg \Omega\bT\varrho \bT\Omega\bT\varrho| +E[\Omega\bT\varrho](0)+F_0[\Omega\bT\varrho](\H_{u_*}^{t_1})\\
   &\les \int_{\D_{u_1}^{t_1}} |\Omega\Box_\bg\bT\varrho+\Er_1(\bT\varrho, \Omega)|| \bT\Omega\bT\varrho|+\La_0^2.
    \end{align*}
    Using (\ref{1.27.5.24}), (\ref{3.9.6.24}) and (\ref{1.29.2.22}), we infer that 
    \begin{align*}
    \fm{\Omega}\Box_\bg \bT\varrho=O(\l t\r^{-1}\log \l t\r\Delta_0^\frac{5}{4})_{L^2_\Sigma}.
    \end{align*}
It follows by  using (\ref{2.13.5.24}) and (\ref{2.17.3.24})
    \begin{align*}
 \|\sP[\Omega, \bT\varrho]\|_{L^2_\Sigma}\les \l t\r^{-2}\log \l t\r^{\f12\M+1}\La_0+\l t\r^{-\frac{7}{4}+2\delta}\Delta_0^\frac{5}{4}.
    \end{align*}
    Combining the above two estimates, in view of (\ref{5.02.3.21_1}) we have 
\begin{equation}\label{3.9.12.24}
\|\Er_1(\bT\varrho, \Omega)\|_{L^2_\Sigma}\les\l t\r^{-2}\log \l t\r^{\f12\M+1}\La_0+\l t\r^{-1}\log \l t\r\Delta_0^\frac{5}{4}. 
\end{equation}
Applying  the above estimate and (\ref{2.16.2.24}) to the standard energy inequality, we obtain
 \begin{align*}
  & E[\Omega\bT\varrho](t_1)+F_0[\Omega\bT\varrho](\H_{u_1}^{t_1})\\
   &\les |\int_{\D_{u_1}^{t_1}} |(\fB\Lb\Omega \bT\varrho+O(\l t\r^{-1}\log\l t\r^\frac{9}{4}(\La_0+ \Delta_0^\frac{5}{4})) _{L^2_\Sigma})\bT\Omega\bT\varrho|+\La_0^2\\
   &\les \La_0^2+\int_0^t \|\Lb\varrho\|_{L_x^\infty}E[\Omega\bT\varrho](t')dt+\int_0^t\l t'\r^{-1}\log \l t'\r^\frac{9}{4}(\La_0+\Delta_0^\frac{5}{4})E[\Omega \bT\varrho]^\f12(t') dt'. 
    \end{align*}
    By Cauchy-Schwarz and Gronwall's inequality, we conclude
    \begin{align*}
    E[\Omega\bT\varrho](t_1)+F_0[\Omega\bT\varrho](\H_{u_1}^{t_1})\les \log \l t\r^{\M+7}(\La_0^2+\Delta_0^\frac{5}{2}),
    \end{align*} 
    as stated in (\ref{2.13.3.24}). The second estimate of (\ref{2.13.3.24}) is a consequence of the first estimate in (\ref{2.13.3.24}), (\ref{2.15.3.24}) and (\ref{2.20.2.24}). (\ref{3.9.8.24}) follows as a consequence of (\ref{2.13.3.24}), (\ref{2.16.2.24}) and (\ref{3.9.12.24}). 
 \end{proof}

\subsection{Decay estimates: the second order}
In this subsection, we derive a set of important decay estimates by using Proposition \ref{8.29.8.21}. 
\begin{proposition}\label{9.8.6.22}
 Under the assumptions of (\ref{3.12.1.21})-(\ref{6.5.1.21}), there hold the following decay estimates: 
\begin{equation}\label{9.19.5.23}\left\{
\begin{array}{lll}
\|\tir(S\Phi, S^2\Phi, S\Omega\Phi, \Omega S^{\le 1}\Phi)\|_{L^4_\omega}^4\les \log \l t\r^{2\M} (\Delta_0^5+\La_0^2)\\
\|\tir( S\Omega\varrho, \Omega S\varrho, \Omega\varrho)\|_{L^4_\omega}^4\les \log \l t\r^{2\M}(\Delta_0^5+\La_0^4)\\
\|\tir \Omega^2\varrho\|_{L^4_\omega}^4\les(\Delta_0^5 +\La_0^4)\l t\r^{2\delta}\log \l t\r^\M\\
\|\tir \Omega^2\Phi\|_{L_\omega^4}^4\les\l t\r^{2\delta}\log \l t\r^{\M+10}(\Delta_0^5+\La_0^2)
\end{array}\right.
\end{equation}
and with $X=\Omega, S$, 
\begin{align}
&\begin{array}{lll}
\|\sn_X^{\le 2} [L\Phi] \|_{L^2_\Sigma}\les
  \l t\r^{-1}(\Delta_0^\frac{5}{4}+\La_0)\log \l t\r^{\f12\M}\\
 \|\sn_X^{l+1}\bA_{g,1}\|_{L^2_\Sigma}\les\l t\r^{-1}  (\log \l t\r)^{\f12\M}\big(\La_0+\Delta_0^\frac{5}{4}(1+(1-\vs^+(X^{l+1}))\l t\r^{l\delta})\big), l=0,1
 \end{array}
 \label{L2BA2'}\\
&\|\sn_X\sn_L \zeta, \sn_L \sn_X \ze\|_{L^2_u L_\omega^2}\les \l t\r^{-2}\log\l t\r(\log \l t\r^3\Delta_0^\frac{5}{4}+\La_0)\label{2.20.4.24}
\end{align}
\begin{equation}\label{LbBA2'}
\left\{
\begin{array}{lll}
&\|\sn_X\sn\fB, \sn_X\sn_\Lb \bA_{g,1}, (\log\l t\r)^{-\f12}\sn_\Lb \sn_X\bA_{g,1}\|_{L^2_u L_\omega^2}\les\l t\r^{-2}(\log \l t\r^3\Delta_0^\frac{5}{4}+\La_0)\log \l t\r,\\
 &\|\sn_\Lb (\tir\sn)[L\Phi]\|_{L^2_u L_\omega^2}\les \l t\r^{-2}(\Delta_0^\frac{5}{4}+\La_0)(\log \l t\r)^2
\end{array} \right.
\end{equation}
\begin{equation}\label{2.20.5.24}
\|\sn_\Lb (\tir\sn\eh)\|_{L^2_\Sigma}\les\l t\r^{-2}\log \l t\r^{\f12\M+2}(\La_0+\l t\r^\delta\Delta_0^\frac{5}{4}).
\end{equation}

\begin{align}\label{L4conn'}
\left\{\begin{array}{lll}
\|\tir\sn_X \bA_{g,1}, \tir^2\sn[L\Phi]\|_{L_\omega^4}\les\left\{
\begin{array}{lll}
\l t\r^{-1}\l t\r^{\f12\delta}(\Delta_0^\frac{5}{4}+\La_0)\log \l t\r^{\frac{10+\M}{4}}\\
\qquad\qquad\qquad\qquad\mbox{ if }X=\Omega, \bA_{g,1}=[\sn \Phi],\\
\l t\r^{-1}\log \l t\r^\frac{\M+14}{4}(\Delta_0^\frac{5}{4}+\La_0), \mbox{ otherwise}
\end{array}\right.\\
\begin{array}{lll}
 \|S^{\le 1}([L\Phi])\|_{L^4_\omega}\les (\La_0^\f12+\Delta_0^\frac{5}{4})\l t\r^{-2}\log \l t\r^{\f12\M}\\
\|S^{\le 1}\bA_b\|_{L^4_\omega}\les \l t\r^{-2}\l t\r^{\f12\delta}(\Delta_0^\frac{5}{4}+\La_0^\f12)\log \l t\r^{\frac{\M}{4}}.
\end{array} 
\end{array}\right.
\end{align}
\begin{align}
&\| (\log \l t\r)^{-1} \sn_\Lb \bA_{g,1},  \sn k_{\bN\bN}, \sn\fB\|_{L^4_\omega}\les\l t\r^{-2}\log \l t\r^{\f12\M+1}(\Delta_0^\frac{5}{4}
+\La_0)\label{9.14.3.22}
\end{align}
\begin{equation}\label{3.15.1.24}
\|\tir^{-1}\log \l t\r^{\frac{3}{4}}\chih, \sn_\Lb\eh\|_{L_\omega^4}\les \l t\r^{-3}(\log \l t\r)^{\f12\M+1}\l t\r^{\f12\delta}(\Delta_0^\frac{5}{4}+\La_0)
\end{equation}
\begin{align}
\|\bb^l\sn\log\bb\|_{L_\omega^4}&\les \l t\r^{-1}\log \l t\r^{\f12\M+2+l}(\Delta_0^\frac{5}{4}
+\La_0), l=0,1\label{2.21.2.24}\\
\|\sn_L \sn\log \bb\|_{L_\omega^4}&\les \l t\r^{-2}\log \l t\r^{\f12\M+2}(\Delta_0^\frac{5}{4}+\La_0).\label{6.17.1.24}
\end{align}
\begin{equation}\label{8.8.6.22'}
\begin{split}
\|\tir^2 S^{\le 1}\bA_b\|_{L^2_u L_\omega^2}\les\log\l t\r^{\f12\M+1}(\Delta_0^\frac{5}{4}+\La_0), \\
\|\tir^2\sn_S^{\le 1}\bA_{g,2}\|_{L_u^2 L_\omega^2}\les \log \l t\r^{\f12\M}(\La_0+\Delta_0^\frac{5}{4}).
\end{split}
\end{equation}
\begin{equation}\label{L2conndrv'}
\left\{
\begin{array}{lll}
\|\bb^{-\f12} \vs^+(X^{2})\sn_X^2 \bA_{g,2}\|_{L^2_\Sigma}\les \l t\r^{-1} \log \l t\r^{\f12\M+1}(\La_0+\l t\r^{\delta(1-\vs^-(X^2))}\Delta_0^\frac{5}{4})\\
\|\bb^{-\f12} \vs^+(X^2)X^2\bA_b\|_{L^2_\Sigma}\les\l t\r^{-1} \log \l t\r^{\f12\M+1}(\La_0+\l t\r^{\delta(1-\vs^-(X^2))}\Delta_0^\frac{5}{4})
\end{array}\right.
\end{equation}
\begin{align}\label{3.5.7.24}
\begin{split}
&\| \vs^+(X^2)\sn_X^2 \sn\la, \l t\r^{-1}X^2\la\|_{L^2_u L_\omega^2}\les \log \l t\r^{\f12\M+1}\l t\r^{-1}(\La_0+\l t\r^\delta\Delta_0^\frac{5}{4}),\\
&\|X\la\|_{L_\omega^4}\les \l t\r^{\frac{\delta}{2}(1-\vs(X))}\log \l t\r^{\frac{\M}{4}}(\La_0+\Delta_0^\frac{5}{4}), \|\la\|_{L_\omega^4}\les \log \l t\r^{\frac{\M}{4}+\frac{3}{2}}(\La_0+\Delta_0^\frac{5}{4})
\end{split}
\end{align}
\begin{align}\label{3.7.1.24}
\sn_X \big((L+\tr\chi+k_{\bN\bN}){}\rp{a}\pih_{\Lb \Lb}\big)=O(\l t\r^{-1}\log \l t\r(\La_0+\log \l t\r^3\Delta_0^\frac{5}{4}))_{L^2_u L_\omega^2}
\end{align}
\begin{align}\label{3.7.2.24}
X\sdiv\eta(\Omega)=O\Big(\l t\r^{-1}(\log \l t\r)^{\frac{\M}{2}} (\La_0+\Delta_0^\frac{5}{4})\Big)_{L^2_\Sigma}
\end{align}
\begin{align}\label{3.5.8.24}
\|\sn_\Omega\pioh_{LA}, \vs^+(X^2)\sn_X^2\pioh_{LA}\|_{L^2_u L_\omega^2}\les \l t\r^{-1}\log \l t\r^{\f12\M+1}(\La_0+\Delta_0^\frac{5}{4}\l t\r^\delta)
\end{align}
\begin{equation}\label{7.25.2.22}
\begin{split}
&|\bA_{g,1}|+\l t\r^{-1}|v_A|\les \l t\r^{-2}\l t\r^{\f12\delta}(\Delta_0^\frac{5}{4}+\La_0)(\log \l t\r)^{\frac{10+\M}{4}}\\
&|\log \l t\r^\frac{5}{2}\la, S\la|\les  \l t\r^{\f12\delta}(\Delta_0^\frac{5}{4}+\La_0)(\log \l t\r)^{\frac{10+\M}{4}}\\
&|L\Phi, \tir^{-1}(\varrho, v_\bN)|\les \l t\r^{-2}\log \l t\r^{\frac{\M}{2}}(\Delta_0^\frac{5}{4}+\La_0^\f12)
\end{split}
\end{equation}
\end{proposition}
\begin{remark}
The last line in (\ref{9.19.5.23}) can be further improved whence the top order energy estimates are completed. But they already improved the bootstrap assumptions.
\end{remark}
\begin{proof}
The proof is divided into three parts. 

\noindent$\bullet$ {\it Proof of (\ref{9.19.5.23}).} 
To begin with, we deduce for scalar functions $f$ by using $\pio_{AL}=O(\l t\r^{-\frac{3}{4}+\delta}\Delta_0^\f12)$,  
\begin{align}\label{9.19.2.23}
\begin{split}
\int_{\H_u^t}\tir^{-2}|\Omega(Sf)|^2&\les \int_{\H_u^t}\tir^{-2}\{|S\Omega f|^2+|[S, \Omega]f|^2\}\\
&\les \int_{\H_u^t}\tir^{-2}\{|\tir L\Omega f|^2 +|\tir\pio_{AL} \sn f|^2\}\\
&\les \int_{\H_u^t}\{|L\Omega f|^2+\Delta_0\l t\r^{-\frac{3}{2}+2\delta}|\sn f|^2\}\\
&\les \int_{\H_u^t} (|(L+\f12\tr\chi)\Omega f|^2+|\sn f|^2).
\end{split}
\end{align}
Moreover, by applying (\ref{1.13.1.21}) to $f=\Omega X^l\Phi,\, l=0,1,2$, we derive 
\begin{align}\label{10.22.2.23}
\|\sn X^l \Phi\|_{L^2(\H_u^t)}^2\les \int_{S_{0,u}}  |\Omega X^l \Phi|^2 \tir d\mu_\omega+\|(\sn_L+h) \Omega X^l\Phi\|_{L^2(\H_u^t)}^2.
\end{align}
Now we are ready to consider the estimates in (\ref{9.19.5.23}). Applying (\ref{6.24.12.18}) to $Sf$ with $(\ga, \ga_2, \ga_0')=(1,0, 2)$ leads to
\begin{align*}
\int_{S_{t, u}}|\tir Sf|^4 \tir^{-2}&\les \int_{S_{0,u}}|\tir Sf|^4 \tir^{-2} + \int_{\H_u^t}|L (\tir Sf)|^2   \cdot \sum_{l\le 1}\int_{\H_u^t}\tir^{-2}|\Omega^l(Sf)|^2.
\end{align*}
Note that 
\begin{equation*}
|L(\tir f)|\les |\tir (Lf+hf)|+|\tir(h-\frac{1}{\tir})f|\les|\tir(L f+hf)|+|\tir\bA_b f|.
\end{equation*}
Applying (\ref{1.14.1.22}) to treat the second term on the right-hand side, by using (\ref{3.6.2.21}) we derive for $G=Sf$
\begin{equation}\label{9.11.1.22}
\int_{\H_u^t} |L(\tir G)|^2\les F_2[G](\H_u^t)+\Delta_0\|G\|_{L^2(S_{0,u})}^2.
\end{equation}

Hence, using (\ref{9.11.1.22}), (\ref{9.19.2.23}) and (\ref{1.13.1.21}), we obtain
\begin{align}\label{9.19.3.23}\tag{D1}
\begin{split}
\int_{S_{t, u}}|\tir Sf|^4 \tir^{-2}&\les \int_{S_{0,u}}|\tir Sf|^4 \tir^{-2}  +(F_2[Sf](\H_u^t)+\Delta_0\|Sf\|_{L^2(S_{0,u})}^2)\\
&\times\big(\int_{S_{0,u}}|S f|^2\tir d\mu_\omega+\sum_{Y=S, \Omega}\|(L+\f12\tr\chi)Y f\|^2_{L^2(\H_u^t)}+\|\sn f\|_{L^2(\H_u^t)}^2\big).
\end{split}
\end{align}

Applying (\ref{6.24.12.18}) to $\Omega f$ with $(\ga, \ga_2, \ga_0')=(1,0, 2)$, using (\ref{9.11.1.22})  we bound
\begin{align}\label{9.19.4.23}\tag{D2}
\begin{split}
\int_{S_{t, u}}|\tir \Omega f|^4 \tir^{-2}
&\les \int_{S_{0,u}}|\tir\Omega f|^4 \tir^{-2} + \int_{\H_u^t}|L (\tir \Omega f)|^2  \cdot \sum_{l\le 1}\int_{\H_u^t}|\Omega^l \Omega f|^2\tir^{-2}\\
&\les \int_{S_{0,u}}|\tir\Omega f|^4 \tir^{-2}+(F_2[\Omega f](\H_u^t)+\Delta_0\|\Omega f\|_{L^2(S_{0,u})}^2)\int_{\H_u^t}|\sn\Omega^{\le 1} f|^2.
\end{split}
\end{align}
Applying (\ref{9.19.3.23}) to $f=S \Phi, \Phi$, if $f=S\Phi$ we further treat the term $\|\sn S\Phi\|_{L^2(H_u^t)}$ in (\ref{9.19.3.23}) by using (\ref{9.19.2.23}),  and by using (\ref{10.22.2.23}) with $l=0$ for $\|\sn \Phi\|_{L^2(H_u^t)}$. Hence we obtain
\begin{align*}
\int_{S_{t, u}}|\tir Sf|^4 \tir^{-2}&\les \int_{S_{0,u}}|\tir Sf|^4 \tir^{-2}  +(F_2[Sf](\H_u^t)+\Delta_0\|Sf\|_{L^2(S_{0,u})}^2)\\
&\times\big(\int_{S_{0,u}}|Sf, \Omega \Phi|^2\tir d\mu_\omega+\sum_{Y=S, \Omega}\|(L+\f12\tr\chi)(Y f, \Omega \Phi)\|^2_{L^2(\H_u^t)}\big)\\
&\les (\Delta_0^\frac{5}{2}+\La_0)^2\log \l t\r^{2\M}
\end{align*}
where we used Proposition \ref{12.21.1.21} and (\ref{8.25.2.21}) to conclude the last estimate. 

Next we consider $\|S\Omega\Phi\|_{L_\omega^4(S_{t,u})}$. 
Applying (\ref{9.19.3.23}) to $f=\Omega\Phi$, also using (\ref{10.22.2.23}), we obtain with the help of Proposition \ref{12.21.1.21} and (\ref{8.25.2.21}) that 
\begin{align*}
\|\tir S\Omega\Phi\|^4_{L_\omega^4(S_{t,u})}&\les \|\tir S \Omega\Phi\|^4_{L^4_\omega(S_{0,u})}+(F_2[Sf]+\Delta_0\|Sf\|_{L^2(S_{0,u})}^2)\\
&\times\big(\int_{S_{0,u}}|Sf, \Omega^2 \Phi|^2\tir d\mu_\omega+\sum_{Y=S, \Omega}\|(L+\f12\tr\chi)(Y f, \Omega^2 \Phi)\|^2_{L^2(\H_u^t)}\big)\\
&\les (\Delta_0^\frac{5}{2}+\La_0)^2\log \l t\r^{2\M}.
\end{align*}

In particular if $\Phi=\varrho$ in the above, instead of (\ref{10.22.2.23}), we use the boundedness of $F_0[\Omega^2\varrho](\H_u^t)$  to obtain 
\begin{equation*}
\|\tir S\Omega\varrho\|^4_{L_\omega^4(S_{t,u})}\les \|\tir S \Omega\varrho\|^4_{L^4_\omega(S_{0,u})}+(\Delta_0^\frac{5}{2}+\La_0^2)^2\log \l t\r^{2\M}.
\end{equation*}
Finally, we apply (\ref{9.19.4.23}) to $f=\Omega^l\varrho$ with $l=0,1$, by using the boundedness of the fluxes for $\Omega^{1+\le 1}\varrho$ in (\ref{8.21.4.21}) and (\ref{8.25.2.21}), to derive 
\begin{equation*}
\|\tir\Omega^{l+1}\varrho\|^4_{L_\omega^4(S_{t,u})}\les \|\tir \Omega^{1+1}\varrho\|^4_{L^4_\omega(S_{0,u})}+(\Delta_0^\frac{5}{2}+\La_0^2) (\Delta_0^\frac{5}{2}\l t\r^{2l\delta}+\La_0^2)\log \l t\r^\M.
\end{equation*}
 Applying (\ref{9.19.4.23}) to $f= v$ with the help of (\ref{10.22.2.23}) and Proposition \ref{12.21.1.21}, we obtain
\begin{equation*}
\|\tir\Omega v\|^4_{L_\omega^4(S_{t,u})}\les \|\tir\Omega v\|^4_{L^4_\omega(S_{0,u})}+(\Delta_0^\frac{5}{2}+\La_0^2)(\Delta_0^\frac{5}{2}+\La_0)\log \l t\r^{2\M}.
\end{equation*}
We remark that using (\ref{10.22.2.23}) with $l=2$ requires top order weighted energy bound. Thus we will instead use geometric estimate to prove the estimate of $\Omega^2 v$ in (\ref{9.19.5.23}), which will be done after justifying the first set of estimate in (\ref{L4conn'}).

Thus we conclude 
\begin{align*}
&\|\tir(X\Phi, S^2\Phi, S\Omega\Phi)\|_{L^4_\omega}^4\les \La_0^2+\log \l t\r^{2\M} (\Delta_0^\frac{5}{2}+\La_0^2)(\Delta_0^\frac{5}{2}+\La_0)\\
&\|\tir( S\Omega\varrho, \Omega\varrho)\|_{L^4_\omega}^4\les \La_0^4+\log \l t\r^{2\M}(\Delta_0^\frac{5}{2}+\La_0^2)^2\\
&\|\tir \Omega^2\varrho\|_{L^4_\omega}^4\les \La_0^4+(\Delta_0^5+\La_0^4)\l t\r^{2\delta}\log \l t\r^\M.
\end{align*}
Using the above estimates, (\ref{9.19.5.23}) follows by commuting $S$ with $\Omega$ with the help of (\ref{7.17.6.21}).\\ 


\noindent$\bullet${\it Proof of (\ref{L2BA2'})-(\ref{2.20.5.24}).}
 (\ref{L2BA2'}) can be obtained by using (\ref{8.29.9.21}), (\ref{9.8.2.22}) and (\ref{4.22.4.22}).
  
 
 The estimate of $\sn_X\sn\fB$ in (\ref{LbBA2'}) is a consequence of (\ref{2.14.1.24}) and (\ref{7.13.5.22}).  

 We first claim
 \begin{align}\label{5.12.1.24}
 \begin{split}
\|\sn_X(\ud \bA\c\fB)\|_{L^2_\Sigma} &\les\l t\r^{-1} \log \l t\r(\La_0+\log \l t\r^3\Delta_0^\frac{5}{4})\\
 \|\sn_X(\fB \bA_{g,1})\|_{L^2_\Sigma}&\les\l t\r^{-2}\log \l t\r^\frac{\M}{2}(\La_0+\log \l t\r\Delta_0^\frac{5}{4})\\
 \|\sn_X(\fB \bA_g)\|_{L^2_\Sigma}&\les \l t\r^{-2}\log \l t\r^{\f12\M+1}(\La_0+\l t\r^\delta\Delta_0^\frac{5}{4})\\
 \| \sn_X(\chi\c \bA_{g,1})\|_{L^2_\Sigma}&\les \l t\r^{-2}\log \l t\r^{\f12\M}(\La_0+\Delta_0^\frac{5}{4})\\
 \|\sn_X(\chi\c \ud\bA)\|_{L_u^2 L_\omega^2}&\les \l t\r^{-2}\log\l t\r(\La_0+(\log \l t\r)^3\Delta_0^\frac{5}{4})\\
 \|\sn_X(\ud \bA \bA_{g,1})\|_{L^2_u L_\omega^2}&\les \Delta_0^2 \l t\r^{-3+\f12\delta}\log \l t\r^{\frac{1}{4}\M+\frac{5}{2}}
 \end{split}
 \end{align}
 Indeed,
 combining (\ref{2.20.2.24}), (\ref{1.27.5.24}) and (\ref{8.23.2.23}) gives 
 \begin{align*}
 \|\sn_X(\ud \bA\c\fB)\|_{L^2_\Sigma}&\les \l t\r^{-1}\log \l t\r(\La_0+\log \l t\r^3\Delta_0^\frac{5}{4})+\|\bb^{-\f12}\Omega \fB\|_{L_u^2 L_\omega^4}\log \l t\r\Delta_0\\
 &\les \l t\r^{-1}\log \l t\r(\La_0+\log \l t\r^3\Delta_0^\frac{5}{4}).
 \end{align*}
 Similarly using (\ref{7.13.5.22}), (\ref{2.20.2.24}) and (\ref{8.23.2.23})
 \begin{align*}
 \|\sn_X(\fB \bA_{g,1})\|_{L^2_\Sigma}&\les  \l t\r^{-2}\log \l t\r^{\f12\M}(\La_0+\Delta_0^\frac{5}{4})+\|\bb^{\f12}\Omega \fB\|_{L_u^\infty L_\omega^4}\|\tir\bA_{g,1}\|_{L_u^2 L_\omega^4}\\
 &\les \l t\r^{-2}\log \l t\r^\frac{\M}{2}(\La_0+\log \l t\r\Delta_0^\frac{5}{4})
 \end{align*}
 Using Proposition \ref{7.15.5.22}, (\ref{s2}),  (\ref{5.13.2.24}) and (\ref{10.11.3.23}), we deduce
\begin{align*}
\|\sn_X(\fB \bA_g)\|_{L^2_\Sigma}&\les\|\fB \bA_g\|_{L^2_\Sigma}+\|\tir \fB \widehat{\bR_{4A4B}} \|_{L^2_\Sigma}+\|\sn_\Omega(\fB \bA_g)\|_{L^2_\Sigma}\\
&\les \l t\r^{-2} \Big(\log \l t\r^{\f12\M}\La_0+\l t\r^\delta\Delta_0^\frac{5}{4}\Big) +\|\sn_\Omega(\fB \bA_g)\|_{L^2_\Sigma}.
\end{align*}
Using (\ref{2.18.2.24}), (\ref{8.23.2.23}) and (\ref{3.11.3.21}), we infer 
\begin{align*}
\|\sn_\Omega(\fB \bA_g)\|_{L^2_\Sigma}&\les \|\bb^\f12\tir\Omega \fB\|_{L_u^2 L_\omega^4}\|\bA_g\|_{L_\omega^4}+\|\fB \sn_\Omega \bA_g\|_{L^2_\Sigma}\\
& \les \l t\r^{-2}\log \l t\r^{\f12\M+1}(\La_0+\l t\r^\delta\Delta_0^\frac{5}{4}).
\end{align*}
This gives the third estimate. 
We can obtain the fourth estimate in (\ref{5.12.1.24}) follows by using Lemma \ref{error_prod}, (\ref{3.11.3.21}) and (\ref{10.10.2.23}).
 
To see the fifth estimate in (\ref{5.12.1.24}), noting $\ud \bA=\sn\log \bb+\bA_{g,1}$, applying (\ref{8.24.4.23}) to $F=\sn\log \bb$ and using (\ref{7.13.5.22}) and (\ref{2.20.2.24}), we deduce
\begin{align*}
\sn_X(\chi\c \sn\log \bb)&=O(\l t\r^{-1})\sn_X^{\le 1}\sn\log \bb+O(\l t\r^{-\frac{7}{4}+\delta}\Delta_0^\f12)_{L_\omega^4}\sn\log \bb\\
&=O(\l t\r^{-2}\log\l t\r(\La_0+(\log \l t\r)^3\Delta_0^\frac{5}{4}))_{L_u^2 L_\omega^2}.
\end{align*}

The last estimate in (\ref{5.12.1.24}) can be obtained by using Lemma \ref{5.13.11.21} (5), (\ref{7.13.5.22}) and (\ref{10.11.2.23}).
 Hence (\ref{5.12.1.24}) is proved.

We use (\ref{8.5.1.22+}), the proved estimate in (\ref{LbBA2'}), $k_{\bN\bN}=\Lb \varrho+[L\Phi]$ and the proved case for $\bAn$ to obtain for $l=0,1$
\begin{align*}
\|\sn_X^l\sn_L \zeta\|_{L^2_u L_\omega^2}&\les \|\sn_X^l \sn k_{\bN\bN}, \sn_X^l \sn_L \zb\|_{L^2_u L_\omega^2}+\|\sn_X^l\sn \log \bb\|_{L^2_u L_\omega^2}+\l t\r^{-3+2\delta}\Delta_0^\frac{3}{2}.
\end{align*}
Combining the above estimate with the  last estimate in (\ref{5.12.1.24}) and using (\ref{L2BA2'}) and the proved estimate in (\ref{LbBA2'}) yields
\begin{align*}
\|\sn_X^l\sn_L \zeta\|_{L^2_u L_\omega^2}&\les \l t\r^{-2}(\log \l t\r^3 \Delta_0^\frac{5}{4}+\La_0)\log \l t\r, l=1
\end{align*}  
as stated in (\ref{2.20.4.24}). The second estimate in (\ref{2.20.4.24}) can be obtained by further using (\ref{3.21.1.23}).

Next we consider the remaining estimates in (\ref{LbBA2'}). The estimate of $\|\Lb \Omega[L\Phi]\|_{L_u^2 L_\omega^2}$ can be obtained by using (\ref{2.19.1.24}), (\ref{7.13.5.22}) and the commutator estimate
\begin{align*}
\|[\Omega, \Lb][L\Phi]\|_{L^2_u L_\omega^2}&\les \|\ud\bA(\Omega)\bN[L\Phi]\|_{L_u^2 L_\omega^2}+\|\pioh_\Lb^A \sn[L\Phi]\|_{L_u^2 L_\omega^2}\\
&\les \l t\r^{-2}(\log \l t\r)^2(\La_0+\Delta_0^\frac{5}{4}),
\end{align*}
which is due to (\ref{5.13.10.21}), Proposition \ref{1steng} and (\ref{5.21.1.21}). Hence due to (\ref{9.8.2.22}) 
\begin{align}\label{2.20.3.24}
\|\sn_\Lb (\tir \sn[L\Phi])\|_{L_u^2 L_\omega^2}\les \l t\r^{-2}(\log \l t\r)^2(\La_0+\Delta_0^\frac{5}{4}).
\end{align}
Next we prove the second estimate in (\ref{LbBA2'}). 
  We derive by using (\ref{7.03.1.19}) that  
 \begin{align*}
 \sn_X \sn_\Lb \sn\varrho&=\sn_X (\sn \Lb \varrho-\chib \sn\varrho+\ud \bA \fB). 
 \end{align*}
  Using (\ref{5.12.1.24}) and the first estimate in (\ref{LbBA2'}) we then have
  \begin{align*}
  \|\sn_X \sn_\Lb \sn\varrho\|_{L^2_u L_\omega^2}\les\l t\r^{-2}(\log \l t\r^3\Delta_0^\frac{5}{4}+\La_0)\log \l t\r.
  \end{align*}
Next we consider $\bA_{g,1}=\ep$, and derive from (\ref{5.30.3.23}) that
\begin{equation}\label{10.1.8.23}
\sn_\bN \ep=-\theta \ep+\sn \fB+\sn\log (\bb c)(\fB+\eta)+\fB \ep.
\end{equation}
Hence, with $X=\Omega, S$, using the first estimate in (\ref{LbBA2'}), (\ref{1.30.1.24}) and (\ref{5.12.1.24}), we have
\begin{align*}
\|\sn_X\sn_\bN \ep\|_{L^2_u L_\omega^2}&\les\|\sn_X^{\le 1}\Big((\chi+\fB)\c \ep+\sn \fB+\ud \bA (\fB+\bA_{g,1})\Big)\|_{L^2_u L_\omega^2}\\
&\les\l t\r^{-2} \log \l t\r\Big(\La_0+\log \l t\r^3\Delta_0^\frac{5}{4}\Big).
\end{align*}
Using the second line of (\ref{L2BA2'}) with $X^2$ being $XS$ and the above estimate, we get the second estimate in (\ref{LbBA2'}) for $\bA_{g,1}=\zb$.

For the second estimate in (\ref{LbBA2'}), it remains to consider the case $\bA_{g,1}=\eh$. Using (\ref{L2BA2'}) and (\ref{1.30.2.24}), we have 
\begin{align*}
\|\sn_X\sn_\bN \eh\|_{L^2_\Sigma}&\les \|\sn_X(\snc \hot \ep, \ud\bA\c \ep, \fB \bA_g)\|_{L^2_\Sigma}+\|\tir^{-1}\sn_X^{\le 1}\eh\|_{L^2_\Sigma}+\l t\r^{-\frac{5}{2}+2\delta}\Delta_0^\frac{3}{2}\log \l t\r^\f12.
\end{align*}
Using (\ref{5.12.1.24}), Proposition \ref{1steng}, (\ref{1.30.1.24}) and (\ref{L2BA2'}), we conclude
\begin{align}\label{5.13.1.24}
\|\sn_X\sn_\Lb \eh, \sn_X\sn_\bN \eh\|_{L^2_\Sigma}\les\l t\r^{-2}\log \l t\r^{\f12\M+1}(\La_0+\l t\r^\delta\Delta_0^\frac{5}{4}).
\end{align}
Hence the second estimate in (\ref{LbBA2'}) is proved. 

Refining (\ref{7.4.1.21}) with the help of (\ref{1.27.5.24}) and applying it to $U=\bA_{g,1}$, with $X_1=L, \sn$ with the help of Proposition \ref{1steng} leads to
\begin{align*}
\|[\sn_\bN, L]\bA_{g,1}\|_{L^2_\Sigma}&\les\|\bA_{g,1}\|_{L^2_\Sigma}\l t\r^{-2+\delta}\Delta_0+\l t\r^{\delta-1}\Delta_0\|\sn \bA_{g,1}\|_{L^2_\Sigma}+\|\bb^{-\f12}\sn_\bN \bA_{g,1}\|_{L_u^2 L_\omega^4}\\
\|[\sn_\bN, \sn]\bA_{g,1}\|_{L^2_\Sigma}&\les \|\bA_{g,1}\|_{L^2_\Sigma}\l t\r^{-2+\delta}\Delta_0+\l t\r^{-1}(\|\sn \bA_{g,1}\|_{L^2_\Sigma}+\Delta_0\log \l t\r \|\tir\sn_\bN \bA_{g,1}\|_{L_u^2 L^4_\omega}).
\end{align*}
In both cases, with the help of Sobolev embedding on spheres, the second estimate in (\ref{LbBA2'}) and (\ref{7.13.5.22}), we have for $X_1=\sn, L$
\begin{equation*}
\|[\sn_\bN, X_1]\bA_{g,1}\|_{L^2_u L_\omega^2}\les\l t\r^{-2}(\log \l t\r^3\Delta_0^\frac{5}{4}+\La_0)\log \l t\r^\frac{3}{2}.
\end{equation*}
Hence using the second estimate in (\ref{LbBA2'}) and the above estimate,
\begin{align*}
\|\sn_\bN (\sn, L) \bA_{g,1}\|_{L^2_u L_\omega^2}\les\l t\r^{-2}(\log \l t\r^3\Delta_0^\frac{5}{4}+\La_0)\log \l t\r^\frac{3}{2}.
\end{align*}
This gives the third estimate in (\ref{LbBA2'}) with the help of (\ref{L2BA2'}). The proof of (\ref{LbBA2'}) is completed. 

Moreover, due to (\ref{7.4.1.21}), (\ref{7.13.5.22}), (\ref{5.13.1.24}) and  Sobolev embedding
\begin{align*}
\|[\sn_\bN, \sn]\eh\|_{L^2_\Sigma}&\les \l t\r^{-1}(\|\sn\eh\|_{L^2_\Sigma}+\log \l t\r\Delta_0\|\tir\sn_\bN\eh\|_{L_u^2 L_\omega^4})+\l t\r^{-2+\delta}\Delta_0\|\eh\|_{L^2_\Sigma}\\
&\les\l t\r^{-3}\log \l t\r^{\f12\M+2}(\La_0+\l t\r^\delta\Delta_0^\frac{5}{4}).
\end{align*}
Hence, 
\begin{align*}
\|\sn_\Lb (\tir\sn\eh)\|_{L^2_\Sigma}&\les\l t\r^{-2}\log \l t\r^{\f12\M+2}(\La_0+\l t\r^\delta\Delta_0^\frac{5}{4}).
\end{align*}
Thus (\ref{2.20.5.24}) are proved.
 \\

\noindent$\bullet${\it Proof of (\ref{L4conn'})-(\ref{7.25.2.22}).}
Applying (\ref{9.20.1.23}) to $\bA_{g,1}$, also using (\ref{1.25.2.22}), we bound
\begin{align*}
\|\tir \sn_X\bA_{g,1}\|_{L_\omega^4}&\les \|\tir\sn_X \bA_{g,1}\|_{L_\omega^4(S_{t,u_*})}+\|\tir\sn_\bN (\sn, L) \bA_{g,1}, \tir^{-1}\sn_X \bA_{g,1}\|^\f12_{L^2_\Sigma}\\
&\times\|\sn_\Omega^{\le 1} \sn_X \bA_{g,1}, \Delta_0\log \l t\r\bb^{-1}\sn_X\bA_{g,1}\|^\f12_{L^2_\Sigma}.
\end{align*}
The first set of estimate for $\sn_X\bA_{g,1}$ in (\ref{L4conn'}) follows by substituting (\ref{LbBA2'}) and (\ref{L2BA2'}) into the  above estimate. The estimate for $\tir\sn[L\Phi]$ follows similarly with the help of (\ref{2.20.3.24}). Using the estimate of $\|\sn_X\bA_{g,1}\|_{L_\omega^4}$ and (\ref{10.11.2.23}), we can obtain the pointwise bound of $\bA_{g,1}$ in (\ref{7.25.2.22}).

Next we use Lemma \ref{6.30.4.23} to write
\begin{align*}
[\Omega^2 v]&=\Omega(\bA_{g,1}(\Omega))+\Omega v^\|+O(\l t\r^{-\frac{3}{4}+\delta}\Delta_0)[\Omega v]\\
\Omega^2 v^\|&=\sn_\Omega(\eta(\Omega))+[\Omega v]+O(\l t\r^{-\frac{3}{4}+\delta}\Delta_0)\Omega v^\|.
\end{align*}
Using the first set of estimate in (\ref{L4conn'}) and the lower order estimate of (\ref{9.19.5.23}), we have
\begin{equation*}
\|[\Omega^2 v], \Omega^2 v^\|\|_{L_\omega^4}\les \l t\r^{-1}\l t\r^{\f12\delta}\log \l t\r^\frac{10+\M}{4}(\Delta_0^\frac{5}{4}+\La_0^\f12).
\end{equation*}
This gives the last line in (\ref{9.19.5.23}).

The estimate of $S^{\le 1}[L\Phi]$ in the second line of (\ref{L4conn'}) can be obtained by using (\ref{9.19.5.23}), the first line in Lemma \ref{6.30.4.23}, (\ref{3.11.3.21}) and $\|[\sn\Phi]\|_{L_\omega^4}\les \l t\r^{-2}\log \l t\r^{\f12\M}(\La_0+\Delta_0^\frac{3}{2})$ due to (\ref{10.11.2.23}). 

The estimate of $\|\bA_b\|_{L_\omega^4}$ in (\ref{L4conn'}) can be obtained by using (\ref{7.10.5.22}) with $p=4$, (\ref{9.19.5.23}) and Proposition \ref{12.21.1.21}. $\|\sn_L\bA_b\|_{L^4_\omega}$ can be obtained in view of (\ref{7.10.7.22}), $L v_t=\tr\chi v_t$, using the proved estimates in (\ref{9.19.5.23}), (\ref{10.11.2.23}),  (\ref{3.6.2.21}) and (\ref{3.11.3.21}). The proof of  (\ref{L4conn'}) is complete. 

Next, using (\ref{2.13.3.24}), and (\ref{7.13.5.22}), (\ref{8.29.9.21}), (\ref{2.14.1.24}) and (\ref{9.20.1.23}), we obtain 
\begin{equation}\label{9.25.2.23}
\|\sn\Lb \varrho, \sn k_{\bN\bN}, \sn\Xi_4\|_{L_\omega^4}\les \l t\r^{-2}\log \l t\r^\frac{\M+15}{4}(\Delta_0^\frac{5}{4}
+\La_0)
\end{equation}
where we used  $k_{\bN\bN}, \Xi_4=\Lb\varrho+[L\Phi]$ and (\ref{L4conn'}) in the above. The last two estimates in (\ref{9.14.3.22}) are proved in view of $\M \ge 15$. 
(\ref{2.21.2.24}) follows as the consequence of the estimate of $\sn k_{\bN\bN}$, (\ref{10.29.1.22}) and $|\bb k_{\bN\bN}|\les 1$. 
(\ref{6.17.1.24}) can be obtained by using (\ref{2.21.2.24}), (\ref{1.27.6.24}) and the estimate of $\sn k_{\bN\bN}$ in (\ref{9.14.3.22}). 

It remains to prove the first estimate in (\ref{9.14.3.22}). Recall that $\bA_{g,1}=[\sn\Phi], \eh$. We will check each of them separately in the sequel.

We first apply (\ref{7.03.1.19}) to $f=\varrho$ to derive
\begin{equation*}
\sn_\Lb \sn \varrho=\sn \Lb \varrho-\chib\c \sn \varrho+\ud\bA\fB.
\end{equation*}
Using the first estimate in (\ref{9.25.2.23}), (\ref{10.11.2.23}) and (\ref{2.21.2.24}) we obtain 
\begin{equation*}
\|\sn_\Lb \sn \varrho\|_{L_\omega^4}\les \l t\r^{-2}\log \l t\r^{\f12\M+2}(\Delta_0^\frac{5}{4}
+\La_0).
\end{equation*}
 Using (\ref{9.25.2.23}), (\ref{10.11.2.23}), (\ref{10.1.8.23}) and (\ref{2.21.2.24}) we derive
\begin{align*}
\|\sn_\bN \ep\|_{L_\omega^4}\les\l t\r^{-2}\log \l t\r^{\f12\M+2}(\Delta_0^\frac{5}{4}+\La_0).
\end{align*}
Thus, we have completed the case that $\bA_{g,1}=[\sn\Phi]$.

Using (\ref{8.3.1.23}), (\ref{7.13.5.22})-(\ref{8.8.6.22}), (\ref{2.20.2.24}) and (\ref{1.27.5.24}), we bound 
\begin{align}\label{3.30.2.24}
\|\sn_\bN\chih\|_{L^2_\Sigma}&\les\l t\r^{-1} \log\l t\r^{\frac{9}{2}}(\Delta_0^\frac{5}{4}+\La_0)+\|\widehat{\bR_{4A3B}}\|_{L^2_\Sigma}\\
&\les \l t\r^{-1} \log\l t\r^{\frac{9}{2}}(\Delta_0^\frac{5}{4}+\La_0).\nn
\end{align}
Due to (\ref{10.22.2.22}) and $\M \ge 15$, we infer from (\ref{3.30.2.24})
\begin{equation*}
\|\chih\|_{L_\omega^2}\les \l t\r^{-2}  \log\l t\r^{\f12\M}(\Delta_0^\frac{5}{4}+\La_0).
\end{equation*}
Substituting the above estimate to  (\ref{1.30.2.22}), repeating the proof of (\ref{5.13.3.24}), we obtain 
\begin{equation*}
\|\widehat{\bR_{4A4B}}\|_{L^2_\Sigma}\les \l t\r^{-2}\log \l t\r^{\f12\M}(\La_0+\Delta_0^\frac{5}{4}). 
\end{equation*}
We then substitute the above two estimates to (\ref{s2}) to obtain 
\begin{equation*}
\|\sn_S \chih\|_{L_u^2 L^2_\omega}\les  \l t\r^{-2}  \log\l t\r^{\f12\M}(\Delta_0^\frac{5}{4}+\La_0).
\end{equation*}
Thus the second estimate in (\ref{8.8.6.22'}) is proved. The first one is obtained by using the first estimate in (\ref{8.8.6.22}) directly. 

Using (\ref{2.18.2.24}), (\ref{3.30.2.24}), the above estimate of $\|\chih\|_{L_\omega^2}$ and (\ref{9.20.1.23}), we infer
\begin{equation}\label{3.30.1.24}
\|\chih\|_{L_\omega^4}\les \l t\r^{-2}(\log \l t\r)^{\f12\M+\frac{1}{4}}\l t\r^{\f12\delta}(\Delta_0^\frac{5}{4}+\La_0).
\end{equation}
To see the estimate for $\eh$, we use (\ref{2.1.2.24}) to derive
\begin{align}
\|\sn_\bN \eh\|_{L^4_\omega}&\les \|\snc \hot \ep, \ud\bA\c \ep, \tir^{-1}(\ep, \bA_{g,2})\|_{L_\omega^4}+\l t\r^{-\frac{15}{4}+2\delta}\Delta_0^\frac{3}{2}\label{9.25.3.23}
\end{align}
Substituting the first set of estimate in (\ref{L4conn'}) to the above, using (\ref{3.30.1.24}), (\ref{1.27.5.24}), (\ref{10.11.2.23}) and the first estimate in (\ref{7.25.2.22}) we derive 
\begin{align*}
\|\sn_\Lb\eh\|_{L_\omega^4}&\les \l t\r^{-3}(\log \l t\r)^{\f12\M+1}\l t\r^{\f12\delta}(\Delta_0^\frac{5}{4}+\La_0)
\end{align*}
as stated in (\ref{3.15.1.24}).
Hence we proved the estimate for $\bA_{g,1}=\eh$ in (\ref{9.14.3.22}). Thus the proof of (\ref{9.14.3.22}) is complete. 

Next we prove (\ref{L2conndrv'}).
Using (\ref{8.8.6.22'}) and (\ref{8.23.2.23}), we can refine the third estimate in (\ref{5.12.1.24}) to 
\begin{equation*}
\|\sn_X^l(\fB \chih)\|_{L^2_\Sigma}\les\log \l t\r^{\f12\M+1}\l t\r^{-2+\delta(1-\vs(X))}(\La_0+\Delta_0^\frac{5}{4}), \, l=0,1.
\end{equation*}
Using the above estimate, (\ref{1.30.2.22}), (\ref{L2BA2'}) and Proposition \ref{1steng}, we have 
\begin{align}\label{9.27.3.23}
&\|\sn_X^{\le 1}\widehat{\bR_{4A4B}}\|_{L^2_\Sigma}\les \log \l t\r^{\f12\M+1}\l t\r^{-2+\delta(1-\vs(X))}(\La_0+\Delta_0^\frac{5}{4}).
\end{align}
Using the above two estimates, (\ref{8.8.6.22'}) and (\ref{8.25.4.23}),  we derive
\begin{align}\label{9.27.4.23}
\|\bb^{-\f12}\sn_X\sn_L \chih \|_{L^2_\Sigma}\les \log \l t\r^{\f12\M+1}\l t\r^{-2+\delta(1-\vs(X))}(\La_0+\Delta_0^\frac{5}{4}).
\end{align}
Applying (\ref{3.21.1.23}) to $\bA_{g,2}$ and using (\ref{L2conndrv}) we derive
\begin{align*}
\|\bb^{-\f12}[\sn_\Omega, \sn_L]\bA_{g,2}\|_{L^2_\Sigma}&\les \l t\r^{-\frac{11}{4}+2\delta}\Delta_0^\frac{3}{2}.
\end{align*} 
Thus we obtain the first line in (\ref{L2conndrv'}).

Using (\ref{3.20.1.22}) we bound
\begin{align*}
\|X^l(\wt{L \Xi_4})\|_{L^2_\Sigma}&\les \|X^l(\sn^2 \varrho, LL \varrho, \Box_\bg \varrho)\|_{L^2_\Sigma}+\|X^l(\zb\sn\varrho, \big((\chi, \chib, \fB)L \varrho\big))\|_{L^2_\Sigma}.
\end{align*}
 The first part on the right-hand side can be treated by using (\ref{L2BA2'}) and (\ref{8.26.4.21}).
The other part on the right-hand side is of lower order, which can be treated by using (\ref{8.23.2.23}), (\ref{8.24.4.23}), (\ref{2.27.1.24}) and Proposition \ref{1steng} \begin{footnote}{Noting that $\fB L\varrho$ is a term in $\N(\Phi, \bp\Phi)$, thus $X^l(\fB L\varrho)$ can be treated by using (\ref{8.26.4.21}).}\end{footnote}.  We thus have 
\begin{align*}
\|X^{l}\big(\zb\sn\varrho, (\chi, \chib, \fB)L \varrho\big)\|_{L^2_\Sigma}\les \l t\r^{-2}\log \l t\r^{\f12\M+l}(\Delta_0^\frac{5}{4}+\La_0),\, l=0,1.
\end{align*}
Hence, we conclude for $l=0,1$ that
\begin{align}\label{9.24.1.23}
\|X^l\wt{L \Xi_4}\|_{L^2_\Sigma}\les \l t\r^{-2}\log \l t\r^{\f12\M+l}(\l t\r^{l\delta(1-\vs(X^l))}\Delta_0^\frac{5}{4}+\La_0). 
\end{align}
Using (\ref{8.8.6.22'}), (\ref{2.18.2.24}), (\ref{6.8.3.23}), (\ref{3.6.2.21}), (\ref{L2conndrv}), (\ref{9.24.1.23}) and (\ref{8.26.4.21}), we derive
\begin{align*}
\|\bb^{-\f12}\sn_X\sn_S\bA_b\|_{L^2_\Sigma}&\les \|\bb^{-\f12}\sn_X^{\le 1}(\bA_b, \tir \bA_b^2, \tir\wt{L\Xi_4}, [L\Phi], \tir\N(\Phi, \bp\Phi))\|_{L^2_\Sigma}+\l t\r^{-\frac{7}{4}+2\delta}\Delta_0^\frac{3}{2}\\
&\les  \l t\r^{-1}\log \l t\r^{\f12\M+1}(\l t\r^{\delta(1-\vs(X))}\Delta_0^\frac{5}{4}+\La_0).
\end{align*}
 The remaining second order cases in (\ref{L2conndrv'}) can be obtained by using (\ref{3.21.1.23}) and (\ref{2.18.2.24}).
We completed the proof of (\ref{L2conndrv'}).

In (\ref{3.5.7.24}), the estimates of $X^2\la$ can be obtained by using Proposition \ref{7.13.4.22} together with (\ref{2.18.2.24}) and (\ref{10.10.2.23}). In view of (\ref{8.24.7.23}) and (\ref{8.23.7.23}), we can obtain the estimate of $\vs^+(X^2)\sn_X^2 \sn \la$ by using (\ref{8.29.9.21}), (\ref{10.10.2.23}) and (\ref{8.24.4.23}). Using (\ref{8.24.7.23}), (\ref{9.19.5.23}) and (\ref{3.6.2.21}) , we can derive the estimate of $\Omega \la$ in the second line of (\ref{3.5.7.24}). The estimate of $S\la$ therein follows directly by using (\ref{3.22.5.21}) and (\ref{10.11.2.23}). The $L^4$ estimate of $\la$ is from (\ref{11.26.2.23}). 

(\ref{3.7.1.24}) can be obtained by using (\ref{2.20.4.24}), (\ref{L2BA2'}) and (\ref{5.12.1.24}). (\ref{3.7.2.24}) follows by using (\ref{8.2.2.23}), (\ref{10.10.2.23}) and (\ref{8.29.9.21}). 

In view of (\ref{7.16.2.22}), (\ref{3.5.8.24}) can be obtained by using (\ref{3.5.7.24}), (\ref{10.10.2.23}), (\ref{8.29.9.21}), (\ref{3.28.3.24}) and applying Proposition \ref{10.16.1.22} and Proposition \ref{7.15.5.22} for bounding derivatives of $\la \bA_{g,1}$.

 The estimate in  (\ref{7.25.2.22})  for $L\Phi$ follows by Sobolev embedding on spheres and (\ref{9.19.5.23}). The $S\la$ estimate follows as a consequence and the estimate of $\la$ is a consequence of (\ref{3.5.7.24}) and Sobolev embedding.

It only remains to prove smallness bound for $\Phi$ in (\ref{7.25.2.22}). Recall from (\ref{11.21.1.23}) that
\begin{equation*}
\tir^{-1}\overline{c^{-2} v_\bN}=\overline{\tr\eta+c^{-2}\bA_b v_\bN}.
\end{equation*}
Using the $L^4_\omega$ estimates of $\bA_b$ in (\ref{L4conn'}), and of $[L\Phi]$ in (\ref{9.19.5.23}) and  $v_\bN=O(\l t\r^{-1+\delta}\Delta_0^\f12)$, we obtain
\begin{equation*}
\overline{c^{-2} v_\bN}=O(\l t\r^{-1}\log \l t\r^{\f12\M}(\La_0^\f12+\Delta_0^\frac{3}{2})).
\end{equation*} 
Moreover, due to (\ref{8.7.1.24}), (\ref{1.27.5.24}) and (\ref{10.10.2.23}), $v_A=O(\l t\r^{-2+\delta}\log \l t\r(\Delta_0^\frac{5}{4}+\La_0))_{L_\omega^4}$.
Using this estimate and (\ref{10.11.2.23}), we can obtain the estimate of 
\begin{equation}\label{11.12.4.23}
\sn (v_\bN)=[\sn \Phi]+v_A\c \theta=O(\l t\r^{-2} \log \l t\r^{\frac{\M}{4}+\f12}(\La_0+\Delta_0^\frac{5}{4}))_{L_\omega^4}.
\end{equation}
 Due to the above estimate, \Poincare  inequality, Sobolev embedding   we obtain
\begin{equation}\label{11.12.3.23}
v_\bN=O\big(\l t\r^{-1}\log \l t\r^{\f12\M}(\Delta_0^\frac{5}{4}+\La_0^\f12)\big).
\end{equation}

In view of (\ref{9.29.5.23}), (\ref{9.29.6.23}) and $\tr\eta=[L\Phi]$, with integration by parts on spheres,  
\begin{equation}\label{11.12.2.23}
\overline{\bb\bN(c^{-2}v_\bN-\varrho)-\bb\bN(c^{-2}) v_\bN}=\overline{\bb([L\Phi]+\ud \bA v_A)}.
\end{equation}
We directly compute for scalar functions $f$
\begin{equation}\label{11.11.3.23}
\p_u \bar f=\overline{\bb \bN f+\bb c\tr\thetac (f-\bar f)}.
\end{equation}
Applying the above formula to $f=c^{-1}v_{\hN}-\varrho$ and using (\ref{11.12.2.23}) gives
\begin{align*}
\overline{c^{-2}v_\bN-\varrho}&=\int_u^{u_*}\overline{\bb \bN f+\bb c\tr\thetac (f-\bar f)}\\
&=\int_u^{u_*}\{\overline{\bb([L\Phi]+\ud \bA v_A+\bN(c^{-2})v_\bN)}+\overline{\Osc(\bb c\tr\thetac)\Osc(f)}\}\\
&=\int_u^{u_*}\overline{\bb[L\Phi]}+O(\l t\r^{-1} v_\bN)_{L_u^1L_\omega^1}+O(\l t\r^{-2}(\log \l t\r)^{\f12\M+1}\Delta_0^2),
\end{align*}
where the last estimate is obtained by using  (\ref{6.24.1.21}),  (\ref{11.12.4.23}), (\ref{3.11.3.21}), (\ref{1.27.5.24}), \Poincare  inequality. We can control the first two terms on the right-hand side due to (\ref{10.10.2.23}). 
Consequently,
\begin{align*}
\overline{c^{-2}v_\bN-\varrho}=O(\l t\r^{-2}\log \l t\r^{\f12\M+1}(\Delta_0^\frac{5}{4}+\La_0)).
\end{align*}
Substituting (\ref{11.12.3.23}) to the left-hand side leads to
\begin{equation*}
\bar \varrho=O\big(\l t\r^{-1}\log \l t\r^{\f12\M}(\Delta_0^\frac{5}{4}+\La_0^\f12)\big).
\end{equation*}
Using \Poincare  inequality, Sobolev embedding and (\ref{10.11.2.23}), we obtain
\begin{equation*}
|\varrho|\les \l t\r^{-1}\log \l t\r^{\f12\M}(\Delta_0^\frac{5}{4}+\La_0^\f12).
\end{equation*}
Repeating  the proof of the $v_A$ bound in (\ref{6.24.1.21}), we have by using (\ref{11.12.4.23}), (\ref{11.12.3.23}), the first estimate in (\ref{7.25.2.22}), (\ref{L2BA2'}) and (\ref{3.11.3.21}) that 
$$
v_A=O\left( \l t\r^{-1}\l t\r^{\f12\delta}(\Delta_0^\frac{5}{4}+\La_0)(\log \l t\r)^{\frac{10+\M}{4}}\right).
$$
\end{proof}
\section{Top order energy estimates}\label{top_eng}
In this section, we obtain the  weighted energies including the standard energies upto the top order under the bootstrap assumptions (\ref{3.12.1.21})-(\ref{6.5.1.21}), with which we improve these assumptions except (\ref{6.5.1.21}). The improvement of (\ref{6.5.1.21}) will be derived in Section \ref{10.24.1.23}.
\subsection{Commutator estimates}
We first provide a set of preliminary estimates, which will be frequently used throughout this section.  
\begin{lemma}
Let $X\in \{\Omega, S\}$.
There holds
\begin{align}
X^l \fB &=\left\{\begin{array}{lll}
\vs(X^l)O(\fB)+(1-\vs(X^l))O(\l t\r^{-1} \log \l t\r(\Delta_0^\frac{5}{4}+\La_0))_{L_u^2 L_\omega^2},\, l=1\\
O\Big(\l t\r^{-1}(\log \l t\r^3\Delta_0^\frac{5}{4}+\La_0)\log \l t\r\Big)_{L^2_u L_\omega^2}, \vs^-(X^l)=0,\, l=2\\
O(\fB)+O(\l t\r^{-2}\log \l t\r^{\f12\M+1}(\La_0+\Delta_0^\frac{5}{4}))_{L_u^2 L_\omega^2}, \vs^-(X^l)=1,\, l=2.
\end{array}\right.\label{7.29.1.22}
\end{align}
There hold the following estimates for $S$-tangent tensor $F$ that
\begin{align}
X^2(\fB F)&=O(\bb^{-1} \l t\r^{-1})\sn_X^{\le 2}F+O\Big(\Delta_0(\l t\r^{-1}\log \l t\r)\Big)_{L_\omega^4} \bb^{-1} X F\nn\\
&+O(\Delta_0(\log \l t\r)^2)_{L^2_\Sigma})F\label{8.9.4.22}\\
\|\bb^{-\frac{1}{2}}\sn_X^2(\bA \c F)\|_{L^2_\Sigma}&\les \l t\r^{-\frac{7}{4}+\delta}\Delta_0^\f12(\sum_{i=1,2}\|\sn_\Omega^{\le 1} \sn_{X_i} F\|_{L^2_\Sigma}+\|\sn_X^2 F\|_{L^2_\Sigma})\nn\\
&+\l t\r^{-\f12+\delta}\Delta_0\|F\|_{L^\infty_x}\label{8.9.3.22+}
\end{align}
\end{lemma}
Indeed, we can refine (\ref{8.23.1.23}) to derive (\ref{7.29.1.22}) by using (\ref{7.13.5.22}) for the case $l=1$; and for the case $l=2$ by also using (\ref{8.26.1.23}), (\ref{2.19.1.24}), (\ref{LbBA2'}), (\ref{8.29.9.21}) and Proposition \ref{1steng}.  (\ref{8.9.4.22}) follows by using Lemma \ref{5.13.11.21} (5) together with (\ref{8.23.2.23}). (\ref{8.9.3.22+}) can be obtained by using (\ref{L4conn'}), Proposition \ref{7.15.5.22} and Sobolev embedding on spheres.

 Next, for quantities comparable with $E[\Omega^3\varrho](t)$ and $E[\Omega^2\bT\varrho](t)$, we provide a set of comparison estimates by establishing commutator estimates. 
\begin{lemma}\label{3.24.3.24}
Let $0\le \sig\le \frac{3}{2}$ and $\sig\neq 1$. We have
\begin{align}
\|\tir \bb^{-1+\sig}\Omega^3 \Lb \varrho,\, & \tir\Omega^3(\bb^{-1+\sig}\Lb\varrho)\|_{L_u^2 L_\omega^2}\label{2.22.1.24}\\
&\les \|\tir\bb^{-1+\sig}\Lb\Omega^3\varrho\|_{L_u^2 L_\omega^2}+\Delta_0^\frac{3}{2}\log \l t\r^7\nn\\
&+\big(\int_0^t \{\|\bb^{-1+\sig}\Lb\Omega^3\varrho\|_{L_u^2 L_\omega^2}+\|\Omega^3(\bb^{-1+\sig}[L\Phi])\|_{L_u^2 L_\omega^2}\}\nn\\
&+\La_0+\log\l t\r^8\Delta_0^2+\Delta_0^\frac{5}{4}\big)\c (C\log \l t\r^{-1}\M_0+\log \l t\r^{-\frac{3}{2}}\Delta_0^\f12)\nn
\end{align}
\begin{align}\label{2.26.2.24}
\begin{split}
\|\Omega^3(\bb^{-1+\sig})\|_{L_u^2 L_\omega^2}&\les\int_0^t \{\|\bb^{-1+\sig}\Lb\Omega^3\varrho\|_{L_u^2 L_\omega^2}+\|\Omega^3(\bb^{-1+\sig}[L\Phi])\|_{L_u^2 L_\omega^2}\}\\
&+\La_0+\log\l t\r^8\Delta_0^2+\Delta_0^\frac{5}{4}.
\end{split}
\end{align}
\begin{align}
\|\Omega^2\bN \bT\varrho\|_{L^2_\Sigma}&\les \|\Omega^2(\bb^{-1})\|_{L_u^2L_\omega^4}(1+\Delta_0\log \l t\r)+\|\bN\Omega^2\bT\varrho\|_{L^2_\Sigma}\nn\\
&+\Delta_0^2 (\log \l t\r)^{\f12(\M+11)}\label{3.7.10.24}\\
\|\Omega \bN\Omega\bT\varrho\|_{L^2_\Sigma}&\les 
\|\bN\Omega^2\bT\varrho\|_{L^2_\Sigma}+\Delta_0^2\log \l t\r^{\f12(\M+11)}(\La_0+\Delta_0^\frac{5}{4}).\label{3.7.11.24}
\end{align}

\end{lemma}
\begin{proof} Consider (\ref{2.22.1.24}) first. We will only give the proof of the first estimate, since the second one can be obtained in  the same way. 
Similar to the proof of (\ref{2.14.1.24}), we write
\begin{align}\label{2.22.2.24}
\sta{cba}{\Omega^3}\Lb \varrho&=\Lb \sta{cba}{\Omega^3} \varrho+[{}\rp{c}\Omega, \Lb]\sta{ba}{\Omega^2}\varrho+{}\rp{c}\Omega[{}\rp{b}\Omega, \Lb]{}\rp{a}\Omega \varrho+\sta{cb}{\Omega^2}[{}\rp{a}\Omega, \Lb]\varrho
\end{align}
where we denote ${}\rp{a_n}\Omega {}\rp{a_{n-1}}\Omega\cdots{}\rp{a_1}\Omega f$ by $\sta{a_n a_{n-1}\cdots a_1}{\Omega^n} f$ for short.  
  
To control the second term on the right-hand side, we first bound by using (\ref{8.25.2.21})
\begin{align*}
\|\tir\bb\Omega\log \bb \bN \sta{ba}{\Omega^2}\varrho\|_{L^2_u L_\omega^2}&\les \|\bb\Omega\log \bb\|_{L_\omega^4}\|\tir \bN\sta{ba}{\Omega^2}\varrho\|_{L_u^2 L_\omega^4}\\
\displaybreak[0]
&\les \log \l t\r\Delta_0\|\tir\Omega\bN\sta{ba}{\Omega^2}\varrho\|^\f12_{L^2 L_\omega^2}\|\tir\bN\sta{ba}{\Omega^2}\varrho\|^\f12_{L^2_u L_\omega^2}\\
&\les \log \l t\r\Delta_0\big(\|\tir[\Omega, \bN]\sta{ba}{\Omega^2}\varrho\|^\f12_{L^2_u L_\omega^2}+\|\tir \bN\Omega\sta{ba}{\Omega^2}\varrho\|_{L^2_u L_\omega^2}^\f12\Big)(\La_0+\Delta_0^\frac{5}{4})^\f12\\
&\les \log \l t\r\Delta_0\|\tir(\Omega\log \bb\bN\sta{ba}{\Omega^2}\varrho+\pioh_{A\Lb}\sn\sta{ba}{\Omega^2}\varrho)\|^\f12_{L_u^2 L_\omega^2} (\La_0+\Delta_0^\frac{5}{4})^\f12\\
&+\log \l t\r\Delta_0(\La_0+\Delta_0^\frac{5}{4})^\f12 \|\tir \bN\Omega \sta{ba}{\Omega^2} \varrho\|^\f12_{L_u^2 L_\omega^2}.
\end{align*}
Using Cauchy-Schwarz, (\ref{5.21.1.21}), (\ref{3.11.3.21}) and (\ref{8.29.9.21}) this implies
\begin{align*}
\|\tir\bb\Omega\log \bb \bN \sta{ba}{\Omega^2}\varrho\|_{L^2_u L_\omega^2}
&\les (\log \l t\r)^2\Delta_0^2(\La_0+\Delta_0^\frac{5}{4})+\Delta_0^\frac{7}{4}(\La_0+\Delta_0^\frac{5}{4})^\f12\l t\r^{-\frac{7}{4}+2\delta}\\
&+\log \l t\r\Delta_0(\La_0+\Delta_0^\frac{5}{4})^\f12 \|\tir \bN\Omega \sta{ba}{\Omega^2} \varrho\|^\f12_{L_u^2 L_\omega^2}\\
&\les (\log \l t\r)^2\Delta_0^\frac{9}{4}+\log \l t\r\Delta_0(\La_0+\Delta_0^\frac{5}{4})^\f12 \|\tir \bN\Omega \sta{ba}{\Omega^2}\varrho\|^\f12_{L_u^2 L_\omega^2}.
\end{align*}
Using the above estimate, (\ref{5.21.1.21}) and (\ref{8.29.9.21}), we deduce
\begin{align}\label{2.22.5.24}
\|\tir\bb[{}\rp{c}\Omega, \Lb]\sta{ba}{\Omega^2}\varrho\|_{L_u^2 L_\omega^2}&\les\|\tir\bb\Omega\log \bb \bN \sta{ba}{\Omega^2}\varrho\|_{L^2_u L_\omega^2}+\|\tir\bb \pioh_{\Lb A}\sn\sta{ba}{\Omega^2}\varrho\|_{L_u^2 L_\omega^2}\nn\\
&\les (\log \l t\r)^2\Delta_0^\frac{3}{2}+\log \l t\r\Delta_0(\La_0+\Delta_0^\frac{5}{4})^\f12 \|\tir \bN\Omega \sta{ba}{\Omega^2}\varrho\|^\f12_{L_u^2 L_\omega^2}.
\end{align}
Next, using (\ref{5.13.10.21}) again we consider the third term on the right-hand side of (\ref{2.22.2.24}).
\begin{align*}
{}\rp{c}\Omega[{}\rp{b}\Omega, \Lb]{}\rp{a}\Omega\varrho&={}\rp{c}\Omega\big(-2{}\rp{b}\Omega(\bb^{-1}) \bb\bN{}\rp{a}\Omega\varrho+{}\rp{b}\pih_{\Lb A}\sn_A{}\rp{a}\Omega\varrho\big)\\
&=-2\big(\sta{cb}{\Omega^2}(\bb^{-1})\bb \bN{}\rp{a}\Omega \varrho+{}\rp{b}\Omega(\bb^{-1}) {}\rp{c}\Omega (\bb\bN {}\rp{a}\Omega \varrho)\big)\\
&+{}\rp{c}\Omega({}\rp{b}\pih_{\Lb A}\sn_A{}\rp{a}\Omega\varrho).
\end{align*}
Using (\ref{5.21.1.21}) and Proposition \ref{7.15.5.22}, the last term in the above is bounded by
 $$
 \|\tir \Omega(\pioh_{\Lb A}\sn\Omega \varrho)\|_{L^2_u L_\omega^2}\les \l t\r^{-\frac{7}{4}+2\delta}\Delta_0^\frac{3}{2}.$$
 
 To treat the second term on the right-hand side, we use Lemma \ref{5.13.11.21} (5), (\ref{3.6.2.21}), (\ref{10.10.2.23}) and (\ref{8.29.9.21}) and Sobolev embedding on spheres to derive 
\begin{align}\label{2.24.6.24}
&\|\bb\tir{}\rp{b}\Omega(\bb^{-1}) {}\rp{c}\Omega(\bb\bN{}\rp{a}\Omega \varrho)\|_{L^2_u L_\omega^2}\nn\\
&\les \|\bb\tir\Omega(\bb^{-1}) \Omega(\bb L{}\rp{a}\Omega\varrho)\|_{L^2_u L_\omega^2}+\|\Omega\bb\|_{L_\omega^4}\|\bb^{-1}\tir{}\rp{c}\Omega(\bb\Lb{}\rp{a}\Omega\varrho)\|_{L_u^2 L_\omega^4}\nn\\
&\les\l t\r^{-1+\delta}\Delta_0\log \l t\r^{\f12\M+2}(\Delta_0^\frac{5}{4}+\La_0)+\Delta_0\log \l t\r\|\bb^{-1}\tir{}\rp{c}\Omega(\bb\Lb{}\rp{a}\Omega\varrho)\|_{L_u^2 L_\omega^4}.
\end{align}
For the last term in the above, we have
\begin{align*}
\|\bb^{-1}\tir{}\rp{c}\Omega(\bb\Lb{}\rp{a}\Omega\varrho)\|_{L_u^2 L_\omega^4}\les \|\Lb \sta{ca}{\Omega^2}\varrho\tir\|_{L^2_u L_\omega^4}+\|\bb^{-1}\tir[{}\rp{c}\Omega, \bb\Lb]{}\rp{a}\Omega\varrho\|_{L^2_u L_\omega^4}.
\end{align*}
For the last term, we proceed by using (\ref{5.13.10.21}) to derive
\begin{align*}
\bb^{-1}(\Omega[\Omega, \bb\Lb]\varrho, [\Omega, \bb\Lb]\Omega\varrho)&=\bb^{-1}\Omega\{\bb(\pioh_{A\Lb} \sn \varrho+\zb(\Omega)\bN \varrho+\ze(\Omega)L \varrho)\}\\
&+\pioh_{A\Lb}\sn\Omega \varrho+\zb(\Omega)\bN\Omega\varrho+\ze(\Omega) L\Omega \varrho.
\end{align*}
Then it follows by using (\ref{5.21.1.21}), (\ref{3.6.2.21}), (\ref{3.11.3.21}), Lemma \ref{5.13.11.21} (5), (\ref{10.10.2.23}) and (\ref{8.29.9.21}) that
\begin{align}\label{2.24.4.24}
\|\bb^{-1}\tir[\Omega, \bb\Lb]\Omega\varrho, \bb^{-1}\tir[\Omega^2, \bb\Lb]\varrho\|_{L_u^2 L_\omega^4}\les\l t\r^{-1}(\l t\r^{2\delta}\Delta_0^\frac{5}{4}+\log \l t\r^{\f12\M+1}\La_0).
\end{align}
Therefore
\begin{align*}
\|\bb^{-1}\tir{}\rp{c}\Omega(\bb\Lb{}\rp{a}\Omega\varrho)\|_{L_u^2 L_\omega^4}\les \|\Lb \sta{ca}{\Omega^2}\varrho\|_{L^2_u L_\omega^4}+\l t\r^{-1}(\l t\r^{2\delta}\Delta_0^\frac{5}{4}+\log \l t\r^{\f12\M+1}\La_0).
\end{align*}
Substituting the above estimate to (\ref{2.24.6.24}) yields
\begin{align*}
\|\bb\tir{}\rp{b}\Omega(\bb^{-1}) {}\rp{c}\Omega(\bb\bN{}\rp{a}\Omega \varrho)\|_{L^2_u L_\omega^2}&\les 
\Delta_0 \log \l t\r  \{\|\tir\Lb \sta{ca}{\Omega^2}\varrho\|_{L^2_u L_\omega^4}\\
&+\l t\r^{-1}(\l t\r^{2\delta}(\Delta_0^\frac{5}{4}+\La_0)+\log \l t\r^{\f12\M+1}\La_0)\}.
\end{align*}
By Sobolev embedding on spheres, (\ref{2.22.5.24}), Cauchy-Schwarz and (\ref{8.25.2.21}), we have
\begin{align}
\|\tir\Lb \Omega^2\varrho\|_{L^2_u L_\omega^4}&\les\|\tir\Omega\Lb\Omega^2\varrho\|_{L_u^2 L_\omega^2}^\f12\|\tir\Lb\Omega^2\varrho\|^\f12_{L_u^2 L_\omega^2}\label{2.24.7.24}\\
&\les (\|\tir\Lb\Omega^3\varrho\|_{L_u^2 L_\omega^2}+ \|\tir[\Omega, \Lb]\Omega^2\varrho\|_{L_u^2 L_\omega^2})^\f12(\Delta_0^\frac{5}{4}+\La_0)^\f12\nn\\
&\les (\|\tir\Lb\Omega^3\varrho\|_{L_u^2 L_\omega^2}+(\log \l t \r)^2\Delta_0^\frac{3}{2})^\f12(\Delta_0^\frac{5}{4}+\La_0)^\f12\nn.
\end{align}
Hence
\begin{align*}
\|\tir\bb{}\rp{b}\Omega(\bb^{-1})& {}\rp{c}\Omega(\bb\bN{}\rp{a}\Omega \varrho)\|_{L^2_u L_\omega^2}\\
&\les 
\Delta_0 \log \l t\r  \{ (\|\tir\Lb\Omega^3\varrho\|_{L_u^2 L_\omega^2}+(\log \l t \r)^2\Delta_0^\frac{3}{2})^\f12(\Delta_0^\frac{5}{4}+\La_0)^\f12\\
&+\l t\r^{-1}(\l t\r^{2\delta}(\Delta_0^\frac{5}{4}+\La_0)+\log \l t\r^{\f12\M+1}\La_0)\}.
\end{align*}
  
Let $0\le \sig\le \frac{3}{2}$ and $\sig\neq 1$. Summing over $a, b, c=1,2, 3$, also using (\ref{1.27.2.24}) and Cauchy-Schwartz, we conclude
\begin{align*}
\|\tir\bb^{-1+\sig}\Omega[\Omega, \Lb]\Omega\varrho\|_{L_u^2 L_\omega^2}&\les \|\bb^{\sig-1}\Omega^2(\bb^{-1})\|_{L_u^2 L_\omega^4}\|\tir \bb\bN\Omega\varrho\|_{L_\omega^4}+\ve'\|\bb^{-1}\tir\Lb\Omega^3\varrho\|_{L_u^2 L_\omega^2}\\
&+\log\l t\r^3\Delta_0 (\La_0+\Delta_0^\frac{5}{4})+\|\bb^\f12\tir{}\Omega({}\rp{b}\pi_{\Lb A}\sn_A\Omega\varrho)\|_{L_u^2 L_\omega^2}\\
&\les  \|\bb^{\sig-1}\Omega^2(\bb^{-1})\|_{L_u^2 L_\omega^4}\Delta_0+\ve'\|\bb^{-1}\tir\Lb\Omega^3\varrho\|_{L_u^2 L_\omega^2}\\
&+\log\l t\r^3\Delta_0 (\La_0+\Delta_0^\frac{5}{4}),
\end{align*}
where $\ve'>0$ is a small constant. 

It is direct to compute
\begin{align*}
\|\bb^{\sig-1}\Omega^2(\bb^{-1})\|_{L_u^2 L_\omega^4}&\les \|\bb^{-1}\Omega^2(\bb^{-1+\sig})\|_{L_u^2 L_\omega^4}+\|\bb^{\sig-2}(\Omega(\log \bb))^2\|_{L_u^2 L_\omega^4}\\
&\les \|\Omega^2(\bb^{-1+\sig})\|_{L_u^2 L_\omega^4}+\Delta_0 \log \l t\r \|\Omega(\bb^{-1+\sig})\bb^{\sig-2}\|_{L_u^2 L_\omega^\infty}.
\end{align*}
With the help of Sobolev embedding on spheres and (\ref{1.27.5.24}), we conclude
\begin{align*}
\|\tir\bb^{-1+\sig}\Omega[\Omega, \Lb]\Omega\varrho\|_{L_u^2 L_\omega^2}&\les \Delta_0\log \l t\r \|\Omega^2(\bb^{-1+\sig})\|_{L_u^2 L_\omega^4}+\ve'\|\bb^{-1}\tir\Lb\Omega^3\varrho\|_{L_u^2 L_\omega^2}\\
&+\log\l t\r^3\Delta_0 (\La_0+\Delta_0^\frac{5}{4}).
\end{align*}
As the consequence of the above estimate, (\ref{2.22.5.24}) and (\ref{2.22.3.24}),  we have 
\begin{align}\label{12.25.1.24}
\begin{split}
\|\bb^{-1+\sig}\tir\Omega\Lb \Omega\varrho\|_{L_u^2 L_\omega^4}&\les(\|\bb^{-1+\sig}\tir[\Omega^2, \Lb]\Omega\varrho\|^\f12_{L_u^2 L_\omega^2}+\|\bb^{-1+\sig}\tir\Lb\Omega^3\varrho\|^\f12_{L_u^2 L_\omega^2})\|\bb^{-1+\sig}\tir\Omega\Lb\Omega\varrho\|^\f12_{L_u^2 L_\omega^2} \\
&+\log \l t\r\Delta_0 \|\bb^{-1+\sig}\tir\Omega\Lb \Omega\varrho\|_{L_u^2L_\omega^2}\\
&\les(\|\tir\Lb\Omega^3\varrho\|^\f12_{L_u^2 L_\omega^2}\log \l t\r^\frac{1}{4} +\Delta_0^\f12\log \l t\r^\f12 \|\Omega^2(\bb^{-1+\sig})\|^\f12_{L_u^2 L_\omega^4})\\
&\cdot
(\La_0+\log \l t\r^2 \Delta_0^\frac{5}{4})^\f12+\log \l t\r^3 \Delta_0^\frac{3}{2}.
\end{split}
\end{align}
We will refer to this estimate shortly. 

 For the last term in (\ref{2.22.2.24}),  we write, in view of (\ref{5.13.10.21}), that
\begin{align*}
&\bb^{-1+\sig}\sta{cb}{\Omega^2}[{}\rp{a}\Omega, \Lb]\varrho\\
&=\frac{2}{-1+\sig}\{\sta{cba}{\Omega^3}(\bb^{-1+\sig})\bN\varrho+\bb^{-1+\sig}\Big(\Omega^2(\bb^{-1+\sig})\Omega(\bb^{1-\sig}\bN \varrho)+\Omega(\bb^{-1+\sig})\Omega^2(\bb^{1-\sig}\bN \varrho)\Big)\}\\
&+\bb^{-1+\sig}\Omega^2(\pioh_{\Lb A}\sn\Lb \varrho).
\end{align*}


Now we treat the second and the last term on the right-hand side. Note by using (\ref{2.20.2.24}) and (\ref{1.27.2.24})
\begin{align}
\label{12.24.1.24}
\begin{split}
\|\Omega^2(\bb^{-1+\sig})\|_{L_u^2 L_\omega^4}\|\tir\Omega(\bb^{1-\sig} \bN\varrho)\bb^{-1+\sig}\|_{L_u^\infty L_\omega^4}&\les \Delta_0\|\Omega^2 (\bb^{-1+\sig})\|_{L_u^2 L_\omega^4}\log \l t\r\\
&\les \Delta_0 \|\Omega^3(\bb^{-1+\sig})\|_{L_u^2 L_\omega^2}^\f12\log \l t\r^\frac{3}{2}(\log\l t\r^\frac{5}{2}\Delta_0^\frac{5}{4}+\La_0)^\f12.
\end{split}
\end{align}

Using (\ref{3.25.1.22}), (\ref{5.21.1.21}), (\ref{3.6.2.21}) and (\ref{LbBA2}), we have
\begin{align*}
\|\tir\Omega^2(\pioh_{\Lb A}\sn\Lb \varrho)\|_{L_u^2 L_\omega^2}&\les\|\tir\sn_\Omega^2 \pioh_{\Lb A}\|_{L_u^2 L_\omega^2}\|\sn\Lb \varrho\|_{L^\infty_\omega}+\|\tir\sn_\Omega\pioh_{\Lb A}\|_{L_\omega^4}\|\sn_\Omega \sn \Lb \varrho\|_{L_u^2 L_\omega^4}\\
&+\|\pioh_{\Lb A}\|_{L^\infty_\omega}\|\tir\sn_\Omega^2 \sn \Lb \varrho\|_{L_u^2 L_\omega^2}\\
&\les \l t\r^{2\delta-\frac{3}{2}}\Delta_0^2+\l t\r^{-\frac{7}{4}+2\delta}\Delta_0^\frac{3}{2}.
\end{align*}

Using (\ref{1.27.5.24}) and Proposition \ref{7.15.5.22}, we bound
\begin{align*}
\|\bb^{-1+\sig}&\tir\Omega(\bb^{-1+\sig})\Omega^2(\bb^{1-\sig}\bN \varrho)\|_{L_u^2 L_\omega^2}\\
&\les \|\tir \Omega(\log \bb)\bb^{-1+\sig}\Lb\Omega^2\varrho, \tir\bb^{-2+2\sig}\Omega(\log \bb)[\Omega^2,\bb^{1-\sig}\Lb]\varrho\|_{L_u^2 L_\omega^2}\\
&+\|\bb^{-1+\sig}\tir\Omega(\bb^{-1+\sig})\Omega^2(\bb^{1-\sig}L \varrho)\|_{L_u^2 L_\omega^2}\\
&\les \Delta_0 \log \l t\r \|\tir\bb^{-2+\sig}\Lb \Omega^2\varrho, \tir\bb^{-3+2\sig} [\Omega^2, \bb^{1-\sig}\Lb]\varrho\|_{L_u^2 L_\omega^4}+\l t\r^{-1+3\delta}\Delta_0^2.
\end{align*}
Next if $0<\sigma\le \frac{3}{2}$, we apply (\ref{2.23.2.24}) to $\a=\sig-1$ to derive
\begin{align*}
[\Omega^2, \bb^{1-\sig}\Lb]\varrho&=\sum_{i=0}^1\{\frac{\sig}{\sig-1}\Omega^i(\Omega(\bb^{1-\sig})\Lb \Omega^{1-i}\varrho)-\frac{1}{\sig-1}\Omega^i(\Omega(\bb^{1-\sig})L\Omega^{1-i}\varrho)\\
&+\Omega^i(\bb^{1-\sig}\pioh_{A\Lb} \sn_A \Omega^{1-i}\varrho)\}.
\end{align*}
Using the above identity, similar to the estimate of (\ref{2.24.4.24}), for $0\le \sig\le \frac{3}{2}$ we bound
\begin{align*}
\|\bb^{-3+2\sig}\tir&[\Omega^2, \bb^{1-\sig}\Lb]\varrho\|_{L_u^2 L_\omega^4}\\
&\les\|\tir\bb^{-3+2\sig}\Omega^2(\bb^{1-\sig})\c \Lb\varrho\|_{L_u^2 L_\omega^4}+\|\bb^{-2+\sig}\tir\Omega \log\bb(\Omega\Lb \varrho, \Lb\Omega\varrho)\|_{L_u^2 L_\omega^4}\\
&+\l t\r^{-1+\delta}(\l t\r^{2\delta}\Delta_0^\frac{5}{4}+\log \l t\r^{\f12\M+1}\La_0).
\end{align*}
It follows by using the fact that $$\bb^{-2+2\sig}\Omega^2(\bb^{1-\sig})=-(\Omega^2(\bb^{-1+\sig})+\Omega(\bb^{1-\sig})\Omega(\bb^{-2(1-\sig)}))$$ and Lemma \ref{5.13.11.21} (5), and $\bb\tir |\Lb\varrho|\les 1$ that
\begin{align*}
\|\bb^{-3+2\sig}\tir[\Omega^2, \bb^{1-\sig}\Lb]\varrho\|_{L_u^2 L_\omega^4}&\les \|\Omega^2(\bb^{-1+\sig})\|_{L_u^2 L_\omega^4}+ \|\Omega(\bb^{-1+\sig})\|_{L_u^2 L_\omega^\infty}\\
&+\l t\r^{-1+\delta}(\l t\r^{2\delta}\Delta_0^\frac{5}{4}+\log \l t\r^{\f12\M+1}\La_0).
\end{align*}
 Hence, in view of Sobolev embedding, summarizing the above estimates implies
\begin{align*}
\|\bb^{-1+\sig}\tir\Omega(\bb^{-1+\sig})&\Omega^2(\bb^{1-\sig}\bN \varrho)\|_{L_u^2 L_\omega^2}\\
&\les \Delta_0\log \l t\r\big(\|\tir\bb^{-2+\sig}\Lb \Omega^2\varrho\|_{L_u^2 L_\omega^4}+\|\Omega^{1+\le 1}(\bb^{-1+\sig})\|_{L_u^2 L_\omega^4}\\
&+\l t\r^{-1+\delta}(\l t\r^{2\delta}\Delta_0^\frac{5}{4}+\log \l t\r^{\f12\M+1}\La_0)\big).
\end{align*}
Thus  also using (\ref{2.24.7.24}) together with (\ref{1.27.5.24}) and (\ref{8.25.2.21}), we derive
\begin{align}\label{6.15.2.24}
\|\bb^{-1+\sig}\tir\sta{cb}{\Omega^2}&[{}\rp{a}\Omega, \Lb]\varrho+\frac{2\tir}{1-\sig}\sta{cba}{\Omega^3}(\bb^{-1+\sig}) \bN\varrho)\|_{L_u^2 L_\omega^2}\nn\\
\displaybreak[0]
&\les \Delta_0 \|\Omega^3(\bb^{-1+\sig})\|_{L_u^2 L_\omega^2}^\f12\log \l t\r^\frac{3}{2}(\log\l t\r^\frac{5}{2}\Delta_0^\frac{5}{4}+\La_0)^\f12\nn\\
&+ \Delta_0\log \l t\r\big(\|\tir\Lb \Omega^2\varrho\|_{L_u^2 L_\omega^4}+\log \l t\r\Delta_0^\f12\big)\\
&\les \Delta_0 \|\Omega^3(\bb^{-1+\sig})\|_{L_u^2 L_\omega^2}^\f12\log \l t\r^\frac{3}{2}(\log\l t\r^\frac{5}{2}\Delta_0^\frac{5}{4}+\La_0)^\f12\nn\\
&+\ve'\|\bb^{-1+\sig}\tir\Lb\Omega^3\varrho\|_{L_u^2 L_\omega^2}+\log\l t\r^3\Delta_0^\frac{3}{2}.\nn
\end{align}
Summarizing the above calculations implies
\begin{align}\label{2.25.1.24}
&\|\tir\Big(\bb^{-1+\sig}(\sta{cba}{\Omega^3}\Lb \varrho-\Lb \sta{cba}{\Omega^3} \varrho)+\frac{2}{1-\sig}\sta{cba}{\Omega^3}(\bb^{-1+\sig}) \bN \varrho\Big)\|_{L_u^2 L_\omega^2}\nn\\
&\les\ve'\|\bb^{-1+\sig}\tir \Lb\Omega^3\varrho\|_{L^2_u L_\omega^2}+\Delta_0 \|\Omega^3(\bb^{-1+\sig})\|_{L_u^2 L_\omega^2}^\f12\log \l t\r^\frac{3}{2}(\log\l t\r^\frac{5}{2}\Delta_0^\frac{5}{4}+\La_0)^\f12+\log\l t\r^3\Delta_0^\frac{3}{2}.
\end{align}

Next we bound $\|\Omega^3(\bb^{-1+\sig})\|_{L_u^2 L_\omega^2}$ using (\ref{2.23.2.24})-(\ref{2.23.3.24}). We derive by using (\ref{lb}) and (\ref{3.22.1.21}) that
\begin{align*}
L\sta{cba}{\Omega^3}(\bb^{-1+\sig})&=\sta{cba}{\Omega^3} L(\bb^{-1+\sig})+[L, \sta{cba}{\Omega^3}](\bb^{-1+\sig})\\
&=\sta{cba}{\Omega^3}(\bb^{-1+\sig}(1-\sig) k_{\bN\bN})+[L, \sta{cba}{\Omega^3}](\bb^{-1+\sig})\\
&=\sta{cba}{\Omega^3} (\bb^{-1+\sig}(1-\sig)(\f12\wp\Lb \varrho+[L\Phi]))+[L, {\Omega^3}](\bb^{-1+\sig})\\
&=(1-\sig)\big(\f12\wp([\sta{cba}{\Omega^3},\bb^{-1+\sig}\Lb]\varrho+\f12\wp\bb^{-1+\sig}\Lb \sta{cba}{\Omega^3} \varrho)\\
&+\Omega^3 (\bb^{-1+\sig}[L\Phi])\big)+[L, \Omega^3](\bb^{-1+\sig}).
\end{align*}
For the first commutator on the right-hand side, we apply (\ref{2.23.2.24}) to derive the detailed expansion
\begin{align*}
L\sta{cba}{\Omega^3}(\bb^{-1+\sig})&= \f12\wp(2-\sig) \sta{cba}{\Omega^3}(\bb^{-1+\sig})\Lb \varrho+\f12 \wp \bb^{-1+\sig}\Lb \sta{cba}{\Omega^3} \varrho+\Omega^3 (\bb^{-1+\sig}[L\Phi])\\
&+[L, \Omega^3](\bb^{-1+\sig})+\sum_{i=0}^2\sum_{a=0, a<2}^i\Omega^{1+a}(\bb^{-1+\sig})\Omega^{i-a}(\Lb \Omega^{2-i} \varrho)\\
&+\sum_{i=0}^2\Omega^i\{\Omega(\bb^{-1+\sig})L\Omega^{2-i}\varrho+\bb^{-1+\sig}\pioh_{A\Lb}\sn_A\Omega^{2-i}\varrho\}
\end{align*}
where it is only necessary to keep the coefficient of the first term on the right-hand side precise.
Noting that the sign of the first term on the right-hand side of  is $+$, in the same way as for deriving (\ref{2.24.3.24}), using (\ref{6.5.1.21}), ignoring the negative part of $\Lb \varrho$ due to the $+$ sign, multiplying the above identity by $\sta{cba}{\Omega^3}(\bb^{-1+\sig})$ we derive  by applying the transport lemma and (\ref{5.13.10.21}) that
\begin{align*}
&\|\Omega^3(\bb^{-1+\sig})(t)\|_{L_u^2 L_\omega^2}\\
&\les \|\Omega^3(\bb^{-1+\sig})(0)\|_{L^2_u L_\omega^2}+\f12\wp\int_0^t\|\bb^{-1+\sig}\Lb \Omega^3 \varrho\|_{L_u^2 L_\omega^2}+\int_0^t  \|\Omega^3(\bb^{-1+\sig}[L\Phi])\|_{L_u^2 L_\omega^2}\\
&+\int_0^t \{ \sum_{i=0}^1\|\Omega^i(\Omega(\bb^{-1+\sig})\Lb\Omega^{2-i}\varrho)\|_{L^2_u L_\omega^2}+\sum_{a=0}^1\Omega^{1+a}(\bb^{-1+\sig})\Omega^{2-a}(\Lb \varrho)\}\\
 &+\int_0^t \sum_{i=0}^2 \|\Omega^i\{\pioh_{AL}\sn_A \Omega^{2-i}(\bb^{-1+\sig})+\Omega(\bb^{-1+\sig})L\Omega^{2-i}\varrho+\bb^{-1+\sig}\pioh_{A\Lb}\sn_A\Omega^{2-i}\varrho\}\|_{L_u^2 L_\omega^2}. 
\end{align*}
Next we treat the error terms on the right-hand side. Using (\ref{2.24.5.24}), (\ref{1.27.5.24}), (\ref{zeh}), (\ref{5.21.1.21}) and (\ref{7.11.1.24}), we deduce
\begin{align*}
&\sum_{i=0}^2 \|\Omega^i\{\pioh_{AL}\sn_A \Omega^{2-i}(\bb^{-1+\sig})\}\|_{L_u^2 L_\omega^2}\\
&\les \|\sn_\Omega^2 \pioh_{A L}\sn(\bb^{-1+\sig})\|_{L^2_u L_\omega^2}+\l t\r^{-1}\|\sn_\Omega\pioh_{AL}\|_{L_\omega^4}\|\Omega^2(\bb^{-1+\sig})\|_{L_u^2 L_\omega^4}\\
&+\|\pioh_{AL}\|_{L_\omega^\infty}\|\sn_\Omega^2 \sn(\bb^{-1+\sig}), \sn_\Omega \sn \Omega(\bb^{-1+\sig}), \sn\Omega^2(\bb^{-1+\sig})\|_{L_u^2 L_\omega^2}\\
&\les\l t\r^{-\f12+\delta}\Delta_0^\f12 \|\sn(\bb^{-1+\sig})\|_{L_u^2 L_\omega^\infty}+\l t\r^{-\frac{3}{4}-1+\delta}\Delta_0^\f12 \|\Omega^2(\bb^{-1+\sig})\|_{L^2_u L_\omega^4}\\
&+\l t\r^{-\frac{3}{4}-1+\delta}\Delta_0^\f12(\|\Omega^{1+\le 2}(\bb^{-1+\sig})\|_{L_u^2 L_\omega^2}+\Delta_0\l t\r^{\delta}\|\sn(\bb^{-1+\sig})\|_{L_\omega^4})\\
\displaybreak[0]
&\les \l t\r^{-\frac{3}{2}+\delta}(\Delta_0^\frac{3}{2}+\Delta_0^\f12\|\Omega^3(\bb^{-1+\sig})\|_{L_u^2 L_\omega^2}^\f12\|\Omega^2(\bb^{-1+\sig})\|_{L_u^2 L_\omega^2}^\f12)\\
&\les \l t\r^{-\frac{3}{2}+\delta}\Big(\Delta_0\|\Omega^3(\bb^{-1+\sig})\|_{L_u^2 L_\omega^2}+\log \l t\r(\La_0+\log \l t\r^\frac{5}{2} \Delta_0^\frac{5}{4})\Big)
\end{align*}
where we applied (\ref{2.20.2.24}) and Cauchy-Schwarz to obtain the last line. Using (\ref{3.25.1.22}), (\ref{5.21.1.21}), Proposition \ref{7.15.5.22} and (\ref{4.22.4.22}), we infer
\begin{align*}
&\sum_{i=0}^2\|\Omega^i\{\bb^{-1+\sig}\pioh_{A\Lb}\sn \Omega^{2-i}\varrho\}\|_{L_u^2 L_\omega^2}\\
&\les\|\sn_\Omega^2 \pioh_{A\Lb}\|_{L_u^2 L_\omega^2}\|\bb^{-1+\sig}\sn\varrho\|_{L^\infty_x}+\|\sn_\Omega\pioh_{\Lb A}\|_{L_u^2 L_\omega^4}\|\sn(\bb^{-1+\sig}\Omega\varrho), \sn_\Omega(\bb^{-1+\sig}\sn\varrho)\|_{L_\omega^4}\\
\displaybreak[0]
&+\|\pioh_{A\Lb}\|_{L_\omega^\infty}\|\bb^{-1+\sig}\sn\Omega^2\varrho, \sn_\Omega(\bb^{-1+\sig}\sn\Omega\varrho), \sn_\Omega^2(\bb^{-1+\sig}\sn\varrho)\|_{L_u^2 L_\omega^2}\\
&\les \l t\r^{-\frac{5}{2}+3\delta}\Delta_0^\frac{3}{2}. 
\end{align*}
It follows by using Lemma \ref{5.13.11.21} (5), (\ref{8.29.9.21}), (\ref{zeh}), (\ref{3.6.2.21}) and (\ref{3.11.3.21}) that
\begin{align*}
&\sum_{i=0}^2\|\Omega^i\big(\Omega(\bb^{-1+\sig}) L\Omega^{2-i}\varrho\big)\|_{L_u^2 L_\omega^2}+\sum_{i=0}^1\|\Omega^i\big(\Omega(\bb^{-1+\sig}) \Lb\Omega^{2-i}\varrho\big)\|_{L_u^2 L_\omega^2}\\
\displaybreak[0]
&\les \|\Omega^3(\bb^{-1+\sig})\|_{L_u^2 L_\omega^2}\|L\varrho\|_{L_\omega^\infty}+\|\Omega^2(\bb^{-1+\sig})\|_{L_u^2 L_\omega^4}(\sum_{Y=L, \Lb}\|Y\Omega\varrho\|_{L_\omega^4}+\|\Omega L\varrho\|_{L_\omega^4})\\
&+\|\Omega(\bb^{-1+\sig})(\Lb\Omega^2 \varrho, \Omega\Lb \Omega \varrho)\|_{L_u^2 L_\omega^2}+\l t\r^\delta\log \l t\r\Delta_0\|\sta{\Omega^2, L}\varrho\|_{L_u^2 L_\omega^2}\\
&\les \|\Omega^2(\bb^{-1+\sig})\|_{L_u^2 L_\omega^4}\Delta_0\l t\r^{-1}+\log \l t\r\Delta_0\|\bb^{-2+\sig}(\Lb\Omega^2\varrho, \Omega \Lb \Omega\varrho)\|_{L_u^2 L_\omega^4}+\l t\r^{-2+3\delta}\Delta_0^2\\
&\les \|\Omega^3(\bb^{-1+\sig})\|_{L_u^2 L_\omega^2}^\f12\|\Omega^2(\bb^{-1+\sig})\|_{L_u^2 L_\omega^2}^\f12\l t\r^{-1}\Delta_0\\
&+\l t\r^{-1}\log \l t\r^{\frac{9}{4}}\Delta_0^\frac{3}{2}  \|\tir\Lb\Omega^3\varrho\|^\f12_{L_u^2 L_\omega^2}+\l t\r^{-1}\log \l t\r^5 \Delta_0^2
\end{align*}
where we applied (\ref{12.25.1.24}), (\ref{2.24.7.24}) and Sobolev embedding on spheres to derive the last inequality, and $\sta{\Omega^2, L}\varrho$ means $Y^3\varrho$ with exactly one of the vectorfields $Y$ being $L$ and the other two being $\Omega$. 

At last we apply Lemma \ref{5.13.11.21} (5), (\ref{2.20.2.24}), Sobolev embedding and Cauchy-Schwarz to obtain 
\begin{align*}
&\sum_{i=0}^2\|\Omega^i\big(\Omega(\bb^{-1+\sig}) L\Omega^{2-i}\varrho\big)\|_{L_u^2 L_\omega^2}+\sum_{i=0}^1\|\Omega^i\big(\Omega(\bb^{-1+\sig}) \Lb\Omega^{2-i}\varrho\big)\|_{L_u^2 L_\omega^2}\\
&+\|\Omega(\bb^{-1+\sig})\Omega^2\Lb\varrho\|_{L_u^2 L_\omega^2}+\|\Omega^2(\bb^{-1+\sig})\Omega\Lb\varrho\|_{L_u^2 L_\omega^2}\\
&\les (\|\Omega^3(\bb^{-1+\sig})\|_{L_u^2 L_\omega^2}+\|\tir\bb^{-1+\sig} \Lb \Omega^3\varrho\|_{L_u^2 L_\omega^2})\Delta_0\log \l t\r^{-\frac{3}{2}}\l t\r^{-1}+\log \l t\r^7\Delta_0^2 \l t\r^{-1}
\end{align*}
where we  treated the last term on the left-hand side similar to (\ref{12.24.1.24}), and  to treat  $\|\Omega(\bb^{-1+\sig})\Omega^2\Lb\varrho\|_{L_u^2 L_\omega^2}$, we also used (\ref{2.25.1.24}).

Hence, by using Proposition \ref{12.21.1.21} and  Gronwall's inequality, we conclude
\begin{align*}
\|\Omega^3(\bb^{-1+\sig})\|_{L_u^2 L_\omega^2}&\les \La_0+\Delta_0^\frac{5}{4}+\int_0^t \{\|\bb^{-1+\sig}\Lb\Omega^3\varrho\|_{L_u^2 L_\omega^2}+\|\Omega^3(\bb^{-1+\sig}[L\Phi])\|_{L_u^2 L_\omega^2}\}\\
\displaybreak[0]
&+\int_0^t(\|\Omega^3(\bb^{-1+\sig})\|_{L_u^2 L_\omega^2}+\|\tir\bb^{-1+\sig} \Lb \Omega^3\varrho\|_{L_u^2 L_\omega^2})\log \l t'\r^{-\frac{3}{2}}\l t'\r^{-1}\Delta_0\\
&+\log\l t'\r^7\Delta_0^2\l t'\r^{-1}\\
&\les\int_0^t \{\|\bb^{-1+\sig}\Lb\Omega^3\varrho\|_{L_u^2 L_\omega^2}+\|\Omega^3(\bb^{-1+\sig}[L\Phi])\|_{L_u^2 L_\omega^2}\}\\
&+\log\l t\r^8\Delta_0^2+\La_0+\Delta_0^\frac{5}{4}.
\end{align*}
 This gives (\ref{2.26.2.24}).  
Substituting the above estimate to (\ref{2.25.1.24}) together with using Cauchy-Schwarz and (\ref{6.5.1.21}) gives
\begin{align*}
&\|\tir \bb^{-1+\sig}\Omega^3\Lb \varrho(t)\|_{L_u^2 L_\omega^2}\\
&\les \|\tir\bb^{-1+\sig}\Lb\Omega^3\varrho\|_{L_u^2 L_\omega^2}+\|\Omega^3(\bb^{-1+\sig})\|_{L_u^2 L_\omega^2}(C\log \l t\r^{-1}\M_0+\log \l t\r^{-\frac{3}{2}}\Delta_0^\f12)+\Delta_0^\frac{3}{2}\log \l t\r^7\\
\displaybreak[0]
&\les \|\tir\bb^{-1+\sig}\Lb\Omega^3\varrho\|_{L_u^2 L_\omega^2}+\Delta_0^\frac{3}{2}\log \l t\r^7\\
\displaybreak[0]
&+(C\log \l t\r^{-1}\M_0+\log \l t\r^{-\frac{3}{2}}\Delta_0^\f12)\big(\int_0^t \{\|\bb^{-1+\sig}\Lb\Omega^3\varrho\|_{L_u^2 L_\omega^2}+\|\Omega^3(\bb^{-1+\sig}[L\Phi])\|_{L_u^2 L_\omega^2}\}\\
&+\La_0+\log\l t\r^8\Delta_0^2+\Delta_0^\frac{5}{4}\big).
\end{align*}
This gives (\ref{2.22.1.24}).

Next we show (\ref{3.7.11.24}). Using (\ref{5.13.10.21}), (\ref{5.21.1.21}) and (\ref{1.27.5.24}), we write
\begin{align*}
\Omega \bN\Omega\bT\varrho&=[\Omega, \bN]\Omega\bT\varrho+\bN\Omega^2 \bT\varrho\\ 
&=\Omega\log \bb \bN \Omega \bT\varrho+\pioh_{\bN A}\sn \Omega \bT \varrho+\bN\Omega^2 \bT\varrho\\
&=O(\log \l t\r\Delta_0)_{L_\omega^4}\bb^{-1}\bN \Omega \bT\varrho+O(\l t\r^{-\frac{3}{4}+\delta}\Delta_0)\sn\Omega\bT\varrho+\bN\Omega^2 \bT\varrho.
\end{align*}
Therefore using (\ref{2.13.3.24}), (\ref{LbBA2}), Sobolev embedding and Cauchy-Schwarz leads to  
\begin{align*}
\|\Omega \bN\Omega\bT\varrho\|_{L^2_\Sigma}&\les \log \l t\r\Delta_0 \|\tir\Omega\bN\Omega\bT \varrho\|^\f12_{L_u^2 L_\omega^2}\|\tir\bN\Omega\bT\varrho\|^\f12_{L^2_u L_\omega^2}\nn\\
&+\l t\r^{-\frac{7}{4}+\delta}\Delta_0\|\Omega^2 \bT\varrho\|_{L^2_\Sigma}+\|\bN\Omega^2\bT \varrho\|_{L^2_\Sigma}\nn\\
&\les \|\bN\Omega^2\bT\varrho\|_{L^2_\Sigma}+\Delta_0^2\log \l t\r^{\f12(\M+11)}(\La_0+\Delta_0^\frac{5}{4})
\end{align*}
as stated in (\ref{3.7.11.24}).

Now, we use (\ref{5.13.10.21}), (\ref{5.21.1.21}), (\ref{1.29.2.22}) and (\ref{1.27.5.24}) to derive
\begin{align*}
[\Omega^2, \bN]\bT\varrho&=\Big((\Omega\log \bb)^2+\bb\Omega^2(\bb^{-1})\Big)\bN \bT \varrho+\Omega\log \bb\c \bN\Omega \bT\varrho+\sum_{a=0}^1\sn_\Omega^a \pioh_{A\bN}\sn \Omega^{1-a}\bT\varrho\\
&=\bb^{-2}\Big(\Omega^2(\bb^{-1})+\bb^{-1}(\Omega \log \bb)^2\Big)O(\l t\r^{-1})_{L_\omega^4}+O(\log \l t\r\Delta_0)_{L_\omega^4}\bN\Omega \bT\varrho\\
&+O(\l t\r^{-\frac{3}{4}+\delta}\Delta_0)_{L_\omega^4}\sn \bT\varrho+O(\l t\r^{-\frac{3}{4}+\delta}\Delta_0)\sn\Omega \bT\varrho.
\end{align*}
Therefore, by Sobolev embedding and (\ref{1.27.5.24})
\begin{align*}
\|[\Omega^2,\bN]\bT\varrho\|_{L^2_\Sigma}&\les \|\Omega^2(\bb^{-1})\|_{L_u^2L_\omega^4}(1+\Delta_0\log \l t\r)+\log \l t\r\Delta_0\|\tir\Omega\bN\Omega\bT\varrho\|^\f12_{L_u^2 L_\omega^2}\|\tir \bN \Omega\bT\varrho\|^\f12_{L^2_u L_\omega^2}\\
&+\l t\r^{-\frac{3}{4}+\delta}\Delta_0(\|\bb^\f12\tir\sn\bT\varrho\|_{L_u^2 L_\omega^4}+\|\sn\Omega\bT \varrho\|_{L^2_\Sigma}).
\end{align*}
Using (\ref{3.7.11.24}), (\ref{2.13.3.24}) and (\ref{LbBA2}), we conclude (\ref{3.7.10.24}) with the help of Cauchy-Schwarz.
\end{proof}
\begin{corollary}
With $0\le \sig\le \frac{3}{2}$ and $\sig\neq 1$, there hold
\begin{align}
&\|\Omega^3(\bb^{-1+\sig}), \tir \bb^{-1+\sig}\Omega^3\Lb \varrho, \tir \Omega^3(\bb^{-1+\sig}\Lb \varrho), \tir\sn_\Omega^2(\bb^{-1+\sig}\sn\fB)\|_{L_u^2 L_\omega^2}\les \l t\r^\delta\Delta_0\label{2.26.1.24}\\
\displaybreak[0]
&\|\Omega^2\log \bb, \Omega^2(\bb^{-1+\sig})\|_{L_u^2 L_\omega^4}\les \l t\r^\frac{\delta}{2}\log\l t\r^2 \Delta_0\label{3.10.3.24}\\
&\|\Omega\bN\Omega \bT\varrho, \Omega^2 \bN\bT\varrho\|_{L^2_\Sigma}\les \l t\r^\delta\Delta_0\label{3.10.2.24}
\end{align}
\end{corollary}
\begin{proof}
Using (\ref{3.12.1.21}), the first three estimates in (\ref{2.26.1.24}) follow by using  (\ref{2.22.1.24}), (\ref{2.26.2.24}) and Proposition \ref{7.15.5.22}. 
  (\ref{3.10.3.24}) is obtained by interpolating the first estimate of (\ref{2.26.1.24}) or $\|\Omega^3\log \bb\|_{L_u^2 L_\omega^2}\les \l t\r^{\delta}\Delta_0$ (due to Proposition \ref{7.15.5.22}) with (\ref{2.20.2.24}). Using (\ref{LbBA2}), (\ref{9.14.3.22}), (\ref{LbBA2'}), (\ref{1.27.5.24}) and (\ref{3.10.3.24}), we can obtain the last estimate in (\ref{2.26.1.24}). 
  (\ref{3.10.2.24}) is obtained by using (\ref{3.7.10.24}), (\ref{3.7.11.24}), (\ref{3.10.3.24}) and (\ref{3.12.1.21}).
\end{proof}
Next we give more important commutator estimates.
\begin{lemma}\label{8.1_com}

Let $X^n=X_n X_{n-1}\cdots X_1$ be given with $X_i\in \{S, \Omega\}$. We set  $$\sta{X^n,Y}f:=X^n Yf, YX^n f, X^n\cdots X_{l+1}Y X_{l}\cdots X_1 f, \, \, l=2, \cdots n-1.$$

(1) There hold for $f=\varrho, v$
\begin{align}
\sta{\Omega^2, \bN}f-(\Omega^2 \bN f, \bN\Omega^2 f)&=\left\{\begin{array}{lll}
O((\log \l t\r)^4(\Delta_0^\frac{5}{4}+\La_0))_{L^2_\Sigma}\\
O(\l t\r^{-1+\f12\delta}(\log \l t\r)^2\Delta_0)_{L_u^2 L_\omega^4}
 \end{array}\right.
 \label{3.4.3.24}\\
\sta{\Omega^3, \bN}f-(\bN\Omega^3 f, \Omega^3 \bN f)&=O(\l t\r^{\delta}\Delta_0)_{L^2_\Sigma}.\label{3.4.1.24}
\end{align}

(2) 
\begin{align}\label{3.10.7.24}
\begin{split}
\sta{\Omega^2, \bN}{\Phi^\mu}&-O(1)\min(\mu,1)\fB=\left\{\begin{array}{lll}
O((\log \l t\r)^4(\Delta_0^\frac{5}{4}+\La_0))_{L^2_\Sigma}\\
O(\l t\r^{-1+\f12\delta}(\log \l t\r)^2\Delta_0)_{L_u^2 L_\omega^4}
 \end{array}\right.\\
\sta{X^2, \Lb}{\Phi^\mu}&-O(1)\min(\mu,1)\fB\\
&=O(\l t\r^{-1}\log \l t\r\Delta_0)_{L_u^2 L_\omega^4},O\Big(\log \l t\r(\La_0+\Delta_0^\frac{5}{4})\Big)_{L^2_\Sigma}, \mbox{ if } X^2=S\Omega, \Omega S\\
\sta{S^2, \Lb}\Phi&-O(1)\fB=\left\{\begin{array}{lll}
O(\l t\r^{-\frac{7}{4}+\delta}\Delta_0)_{L^2_u L_\omega^4}\\
O(\l t\r^{-1}\log \l t\r^{\f12\M+1}(\La_0+\l t\r^\delta\Delta_0^\frac{5}{4}))_{L^2_\Sigma}.
\end{array}\right.
\end{split}
\end{align}
 
With $n=2,3$
\begin{equation}\label{3.4.2.24}
[X^n,L]\Phi=\left\{\begin{array}{lll}
O(\l t\r^{-2}\log \l t\r^{\f12\M}(\La_0+\Delta_0^\frac{5}{4}))_{L_u^2 L_\omega^4}, \, n\le 2\\
O(\l t\r^{-1}\log \l t\r^{\f12\M}(\La_0+\Delta_0^\frac{5}{4}))_{L^2_\Sigma}, \, n\le 3.
\end{array}\right.
\end{equation}

\end{lemma}
\begin{proof}
 We write by using (\ref{5.13.10.21}) the collection of all the commutators 
\begin{align}\label{9.12.1.22}
\begin{split}
\sta{\Omega^n, \bN}f-(\bN \Omega^n f,\Omega^n \bN f)&=\sum_{\ell=0}^{n-1}\Omega^{n-\ell-1}[\Omega, \bN]\Omega^{\ell} f\\
&=\sum_{\ell=0}^{n-1}\Omega^{n-\ell-1}(\Omega\log \bb \bN \Omega^{\ell} f+\pio_\bN^A \sn_A \Omega^{\ell} f).
\end{split}
\end{align}
Recall from (\ref{2.20.2.24}), Proposition \ref{7.15.5.22} and Lemma \ref{5.13.11.21} (5)  that
\begin{equation}\label{9.12.2.22}
\begin{split}
&\|(\tir\sn)^n \ud\bA\|_{L^2_u L_\omega^2}\les\l t\r^{-1} \log \l t\r(\La_0+\log \l t\r^{3n}\Delta_0^\frac{5}{4}), n=0,1;\\
&\|\ud\bA\|_{L^\infty_x}+\|(\tir\sn)^2 \ud \bA\|_{L_u^2 L_\omega^2}\les \l t\r^{-1+\delta}\Delta_0, \|\bb \ud \bA\|_{L_\omega^4}\les \l t\r^{-1}\log \l t\r\Delta_0\\
&\|\bb\sn_X \ud\bA\|_{L_u^2 L_\omega^4}+\log \l t\r \|\bb\ud\bA\|_{L_u^2 L_\omega^\infty}\les \l t\r^{-1+\f12\delta}\log\l t\r^2\Delta_0
\end{split}
\end{equation}
where $X\in \{S, \Omega\}$.

Moreover we have from (\ref{8.25.2.21}), (\ref{LbBA2'}), (\ref{LbBA2}) and (\ref{3.12.1.21}) that
\begin{equation}\label{5.16.1.24}
\begin{split}
&\|\bN \Omega^2\varrho\|_{L_u^2 L_\omega^4}\les\l t\r^{-1} \l t\r^{\f12\delta}\Delta_0^\f12 (\Delta_0^\frac{5}{4}+\La_0)^\f12\\
&\|\Omega^2 \fB\|_{L_u^2 L_\omega^4}\les \l t\r^{-1}\l t\r^{\f12\delta}\log \l t\r^\frac{5}{2}\Delta_0^\f12(\La_0+\Delta_0^\frac{5}{4})^\f12\\
&\|\Omega\Lb \varrho\|_{L_u^2 L_\omega^4}\les \log \l t\r^3(\La_0+\Delta_0^\frac{5}{4})
\end{split}
\end{equation}
where the last estimate is obtained by using (\ref{7.13.5.22}) and (\ref{2.14.1.24}).

Using the last estimate in the above, (\ref{3.12.4.24}), Lemma \ref{5.13.11.21} (5),  (\ref{9.12.2.22}) and (\ref{12.19.1.23}), we derive
\begin{equation}\label{8.11.1.24}
\bN\Omega\Phi, \Omega \bN \Phi=O(\fB)+O(\l t\r^{-1}\log \l t\r\Delta_0)_{L_\omega^4})
\end{equation}
and thus
\begin{align*}
\bN\Omega\Phi\c\ud\bA&=\ud \bA(
(O(\fB)+O(\l t\r^{-1}\log \l t\r\Delta_0)_{L_\omega^4}))\\
&=\left\{\begin{array}{lll}
O(\l t\r^{-1}(\log \l t\r)^2(\Delta_0^\frac{5}{4}+\La_0))_{L^2_\Sigma}\\
O(\l t\r^{-2+\f12\delta}(\log \l t\r)^2(\Delta_0^\frac{5}{4}+\La_0))_{L_u^2 L_\omega^4}
\end{array}\right.\\
\sn_\Omega \ud \bA\c \bN \Phi&=O(\l t\r^{-1} \log \l t\r(\La_0+\log \l t\r^3\Delta_0^\frac{5}{4}))_{L^2_\Sigma}, O(\l t\r^{-2+\f12\delta}\log \l t\r^2\Delta_0)_{L_u^2 L_\omega^4}.
\end{align*}

From (\ref{9.12.1.22}), we bound for $f=\varrho, v$
\begin{align*}
|\sta{\Omega^2, \bN}f-(\bN \Omega^2 f, \Omega^2 \bN f)|&\les \sum_{\ell=0}^1|\sn_\Omega^{1-\ell}(\ud\bA(\Omega)\bN\Omega^\ell f),\sn_\Omega^{1-\ell}(\pio_\bN^A\sn\Omega^\ell f)|\nn\\
&\les |\tir (\sn^{\le 1}_\Omega \ud\bA\bN f, \ud\bA \bN\Omega f), \sn_\Omega\pio_{\bN A}\c \sn f, \pio_{\bN A}\sn\Omega f|.
\end{align*}
It follows  by using  (\ref{5.21.1.21}), (\ref{3.25.1.22}), (\ref{3.6.2.21}) and (\ref{3.11.3.21}),  
\begin{align*}
\sn_\Omega\pio_{\bN A}\c \sn \Phi, \pio_{\bN A}\sn\Omega \Phi=O((\l t\r^{-3+3\delta})\Delta_0^\frac{3}{2})_{L_u^2 L_\omega^2}, O(\l t\r^{-\frac{11}{4}+2\delta} \Delta_0^\frac{3}{2})_{L_u^2 L_\omega^4}.
\end{align*}
Hence we conclude (\ref{3.4.3.24}).

Similar to  treating the $n=2$ case in (\ref{3.4.3.24}), using (\ref{9.12.2.22}) and (\ref{zeh}) we compute
\begin{align*}
\bN \Omega \Phi \sn_\Omega \ud \bA&=(O(\fB)+O(\l t\r^{-1}\log \l t\r\Delta_0)_{L_\omega^4})\sn_\Omega\ud \bA\\
&=O\big(\l t\r^{-2}\log \l t\r(\La_0+\l t\r^{\f12\delta}\log \l t\r^2\Delta_0^\frac{5}{4})\big)_{L^2_u L_\omega^2}\\
\sn_\Omega^2 \ud \bA\bN \Phi&=O(\l t\r^{-1+\delta}\Delta_0)_{L^2_\Sigma}.
\end{align*}
Next we use (\ref{3.3.5.24}), (\ref{10.10.2.23}), (\ref{8.29.9.21}), Proposition \ref{7.15.5.22} and (\ref{12.19.1.23}) and to derive
\begin{align}\label{5.17.1.24}
\Omega^2 \bN v=\Omega^{\le 2}\fB+\left\{\begin{array}{lll}
O(\l t\r^{-\frac{7}{4}+\delta}\Delta_0^\frac{5}{4})_{L_u^2 L_\omega^2}\\
O(\l t\r^{-\frac{7}{4}+\delta}\Delta_0)_{L_\omega^4}.
\end{array}\right.
\end{align}
We further use (\ref{3.4.3.24}) and  to deduce
\begin{align}
\bN \Omega^2 \Phi&=\bN \Omega^2\varrho+\Omega^2\bN v+[\bN,\Omega^2]v\nn\\
&=\bN\Omega^2 \varrho+\Omega^2 \bN v+O(\l t\r^{-1+\f12\delta}(\log \l t\r)^2\Delta_0)_{L_u^2 L_\omega^4},\label{3.4.4.24}
\end{align}
and also obtain by using  (\ref{5.16.1.24}) and the above two estimates that
\begin{align*}
\bN \Omega^2 \Phi\ud \bA=O(\l t\r^{-2+\f12\delta}(\log \l t\r)^\frac{7}{2}(\Delta_0^\frac{5}{4}+\La_0))_{L^2_u L_\omega^2}.
\end{align*}
Using (\ref{5.21.1.21}), (\ref{3.25.1.22}) and Proposition \ref{7.15.5.22}, we derive
\begin{align*}
\sum_{\ell=0}^2\sn_\Omega^{2-\ell}&(\pio_\bN^A \sn_A \Omega^{\ell} \Phi)=O(\l t\r^{-\frac{5}{2}+2\delta}\Delta_0^\frac{3}{2})_{L_u^2 L_\omega^2}.
\end{align*}
In view of (\ref{9.12.1.22}), summarizing the above estimates gives (\ref{3.4.1.24}). 

Next we consider (2). The estimate of $\sta{\Omega^2, \bN}\Phi$ follows by combining (1), (\ref{8.25.2.21}), (\ref{2.14.1.24}), (\ref{5.17.1.24}) and (\ref{3.4.4.24}).

 Using (\ref{6.30.2.19}) we have for a function $f$
\begin{align}\label{7.19.2.21}
X_2 S\Lb f&=X_2\big(\tir (\sD f-\Box_\bg f-(h-k_{\bN\bN})\Lb f-\hb Lf+2\zb^A \sn_A f)\big).
\end{align}
By using Lemma \ref{5.13.11.21}, we have
\begin{align*}
&X_2(h, k_{\bN\bN})=(1-\vs(X_2))O(\l t\r^{-1}\log \l t\r\Delta_0)_{L_\omega^4}+\vs(X_2)O(\l t\r^{-1}),\\
&\|\sn_{X_2}(\tir \zb)\|_{L_\omega^4}\les \l t\r^{-1+\delta}\Delta_0.
\end{align*}
Hence we derive
\begin{align}\label{3.10.6.24}
\begin{split}
X_2 S\Lb f&=X_2\big(\tir (\sD f-\Box_\bg f)\big)+((1-\vs(X_2))O(\log \l t\r\Delta_0)_{L^4_\omega}+\vs(X_2))(\Lb f, Lf)\\
&+O(1)X_2(\Lb f, Lf)+O(\l t\r^{-1+\delta}\Delta_0)_{L_\omega^4}\sn f+O(\l t\r^{-1+\delta}\Delta_0) \sn_X \sn f.
\end{split}
\end{align}
Next we  derive by using Proposition \ref{8.29.8.21} and Proposition \ref{1steng}
\begin{align*}
\|X_2\big(\tir (\sD \Phi, \Box_\bg \Phi)\big)\|_{L^2_\Sigma}\les \l t\r^{-1}(\log \l t\r)^{\frac{\M}{2}+1} (\La_0+\l t\r^{\delta(1-\vs(X_2))}\Delta_0^\frac{5}{4}),
\end{align*}
and by using Proposition \ref{7.15.5.22} and (\ref{3.29.1.23}) that
\begin{align}
\sum_{a=0}^1\|X^{a}(\tir\sD X^{1-a}\Phi)\|_{L_u^2 L_\omega^4}+\|X_2(\tir\Box_\bg \Phi)\|_{L_u^2 L_\omega^4}\les \l t\r^{-\frac{7}{4}+\delta}\Delta_0.\label{5.22.1.24}
\end{align}
For the remaining terms, if $X_2=\Omega$, applying Proposition \ref{1steng}, (\ref{3.6.2.21}) and (\ref{9.19.5.23}) implies
\begin{align*}
\Omega S\Lb \Phi^\mu=O(1)\Omega \Lb \Phi^\mu+O(\l t\r^{-1}\log \l t\r\Delta_0)_{L_u^2 L_\omega^4}.
\end{align*}
Using Proposition \ref{1steng} and (\ref{12.19.1.23}) instead, we have
\begin{align*}
\Omega S\Lb \Phi^\mu=O(1)\Omega \Lb \Phi^\mu+O(\l t\r^{-1}\log \l t\r^{\f12\M+1}\La_0+\log \l t\r\Delta_0^\frac{5}{4}))_{L^2_\Sigma}.
\end{align*}
Substituting 
 (\ref{3.12.4.24}), (\ref{8.23.2.23}) and (\ref{8.2.1.22}) to treat $\Omega \Lb \Phi^\mu$, we have
\begin{align}\label{3.10.8.24}
\Omega S\Lb\Phi^\mu-\min(\mu,1) O(1)\fB=O(\l t\r^{-1}\log \l t\r\Delta_0)_{L_u^2 L_\omega^4}, O\Big(\log \l t\r(\La_0+\Delta_0^\frac{5}{4})\Big)_{L^2_\Sigma}. 
\end{align}

Due to (\ref{5.13.10.21}) and using (\ref{8.23.1.23}), (\ref{3.25.1.22}) and (\ref{7.26.2.22}), we derive
\begin{align*}
[S,\Omega]\Lb \Phi^\mu=\tir \pioh_{LA}\sn \Lb \Phi^\mu&= \pioh_{LA}O(\fB+\Delta_0 \l t\r^{-1+\delta})\\
 &=\left\{\begin{array}{lll}
 O( \l t\r^{-1+\delta}(\log \l t\r)^{\frac{\M+3}{2}} (\La_0+\Delta_0^\frac{5}{4}))_{L^2_\Sigma}\\
 O(\l t\r^{-2+2\delta}\Delta_0)_{L_u^2 L_\omega^4}.
 \end{array}\right.
\end{align*}
The commutator is a much better term, which is negligible. Thus the bound in (\ref{3.10.8.24}) controls $S\Omega\Lb\Phi^\mu-\min(\mu,1)O(1)\fB$ as well. 

If $X_2=S$ instead, with the help of (\ref{8.23.1.23}), Proposition \ref{1steng}, (\ref{5.22.1.24}), (\ref{8.29.9.21}) and (\ref{3.6.2.21})  we infer from (\ref{3.10.6.24}) that
\begin{align*}
S^2\Lb\Phi^\mu-O(1)\fB=O(\l t\r^{-\frac{7}{4}+\delta}\Delta_0)_{L_u^2 L_\omega^4}, O(\l t\r^{-1}(\log \l t\r)^{\frac{\M}{2}+1} (\La_0+\Delta_0^\frac{5}{4}))_{L^2_\Sigma}.
\end{align*}

Using Proposition \ref{1steng}, Proposition \ref{8.29.8.21} and (\ref{9.19.5.23}), we can improve the estimate of (\ref{3.10.4.24}) to 
\begin{equation*}
S\Lb\Omega \Phi-O(1)\Lb \Omega \Phi, \Omega\Lb S\Phi-O(1) \Omega\Lb \Phi=\left\{\begin{array}{lll}
O(\l t\r^{-\frac{7}{4}+\delta}\Delta_0)_{L_u^2 L_\omega^4}\\
O\Big(\l t\r^{-1}(\log\l t\r^{\f12\M+1}\La_0+\l t\r^{\frac{1}{4}+\delta}\Delta_0^\frac{5}{4})\Big)_{L^2_\Sigma}.
\end{array}\right.
\end{equation*}
Similar to (\ref{3.10.8.24}), using (\ref{3.12.4.24}),  (\ref{7.13.5.22}) and (\ref{8.23.2.23}), we infer from the above estimate
\begin{equation}\label{3.12.7.24}
(S\Lb\Omega\Phi^\mu, \Omega\Lb S\Phi^\mu)-\min(\mu,1)\fB=\left\{\begin{array}{lll}
O(\l t\r^{-1}\log \l t\r\Delta_0)_{L_u^2 L_\omega^4}\\
 O\big(\log \l t\r(\La_0+\Delta_0^\frac{5}{4})\big)_{L^2_\Sigma}.
 \end{array}\right.
\end{equation}
Note due to (\ref{1.27.5.24}), (\ref{10.10.2.23}) and (\ref{8.29.9.21}),
\begin{align*}
\|\tir \ud\bA\sn\Omega\Phi\|_{L^2_\Sigma}\les \l t\r^{-1}\log \l t\r^{\f12\M+1}(\La_0+\l t\r^\frac{\delta}{2}\Delta_0^\frac{5}{4})\Delta_0, 
\end{align*}
and by using (\ref{3.6.2.21}), and interpolating (\ref{10.10.2.23}) and (\ref{8.29.9.21}), we have
\begin{align*}
\|\tir \ud \bA\sn \Omega\Phi\|_{L^2_u L_\omega^4}\les \|\ud\bA\|_{L_\omega^\infty}\|\Omega^2\Phi\|_{L_u^2 L_\omega^4} \les \l t\r^{-2+2\delta}\Delta_0^2.
\end{align*}
Moreover, using (\ref{10.10.2.23}), (\ref{8.29.9.21}), (\ref{8.23.2.23}), (\ref{3.12.4.24}), (\ref{7.13.5.22}) and (\ref{L4BA1}), 
\begin{equation*}
\bN\Omega\Phi^\mu=\min(\mu, 1)\fB+\left\{ \begin{array}{lll}
O(\l t\r^{-1}\log \l t\r\Delta_0)_{L_\omega^4}\\
O\big(\log \l t\r(\La_0+\Delta_0^\frac{5}{4})\big)_{L^2_\Sigma}.
\end{array}\right.
\end{equation*}
Applying the above three estimates in view of (\ref{3.19.2}) and (\ref{3.12.7.24}) gives
\begin{align*}
\Lb S\Omega\Phi^\mu&=S\Lb\Omega\Phi^\mu+\tir \ud\bA\sn\Omega\Phi+O(1) \bN\Omega\Phi^\mu\\
&=\min(\mu, 1)\fB+\left\{ \begin{array}{lll}
O(\l t\r^{-1}\Delta_0\log \l t\r)_{L_u^2 L_\omega^4}\\
O\big(\log \l t\r(\La_0+\Delta_0^\frac{5}{4})\big)_{L^2_\Sigma}.
\end{array}\right.
\end{align*}

Using (\ref{5.13.10.21}), (\ref{9.12.2.22}), (\ref{5.21.1.21}) and Proposition \ref{7.15.5.22}, we have
\begin{align*}
[\Lb, \Omega] S\Phi^\mu&=\Omega\log \bb \bN S\Phi^\mu+\pioh_{\Lb A}\sn S\Phi\\
&=\Omega\log \bb O(\bb^{-1}\l t\r^{-1})+O(\l t\r^{-\frac{3}{4}+\delta}\Delta_0^\f12)\sn S\Phi\\
&=O(\l t\r^{-1}\log \l t\r\Delta_0)_{L_\omega^4}, O(\log \l t\r(\La_0+\Delta_0^\frac{5}{4}))_{L^2_\Sigma}
\end{align*}
where we also used (\ref{1.27.5.24})  to obtain the estimate in the last line. We hence conclude in view of (\ref{3.12.7.24}) that
\begin{align*}
\Lb \Omega S\Phi^\mu-\min(\mu,1)\fB=O(\l t\r^{-1}\log \l t\r\Delta_0)_{L_u^2 L_\omega^4}, O\big(\log \l t\r(\La_0+\Delta_0^\frac{5}{4})\big)_{L^2_\Sigma}. 
\end{align*}

It only remains to consider $\Lb SS\Phi$ and $S\Lb S\Phi$. By using (\ref{6.30.1.19}), we have
\begin{align}\label{8.25.1.21}
\Lb Sf=\tir(\sD f-\Box_\bg f+(\mho+k_{\bN\bN}) L f-h\Lb f+2\zeta^A\sn_A f).
\end{align}
Applying the above identity to $f=S\Phi$, also using Proposition \ref{8.29.8.21}, (\ref{zeh}),  (\ref{8.23.1.21}), (\ref{8.23.1.23}), the first estimate in (\ref{5.22.1.24}), (\ref{10.10.2.23}), and $\tir(\mho, k_{\bN\bN})=O(1)$, we deduce
\begin{align*}
\Lb SS\Phi&=\tir (\sD S\Phi-\Box_\bg S\Phi+(\mho+k_{\bN\bN}) LS\Phi+\ud\bA\sn S\Phi)+O(1)\Lb S\Phi\\
&=O(1)\fB+\left\{\begin{array}{lll}
O(\l t\r^{-\frac{7}{4}+\delta}\Delta_0)_{L^2_u L_\omega^4}\\
O(\l t\r^{-1}\log \l t\r^{\f12\M+1}(\La_0+\l t\r^{\delta}\Delta_0^\frac{5}{4}))_{L^2_\Sigma}.
\end{array}\right.
\end{align*} 
The estimate of $S\Lb S\Phi$ can be obtained in the same way by using (\ref{6.30.2.19}). We summarize the above estimates into (\ref{3.10.7.24}).

For the commutation with $L$, we first write with $n=2,3$,
\begin{align}
[X^n, L]f&=\sum_{a=0}^{n-1}X^{n-1-a}[X, L]X^a f\nn\\
&=\sum_{a=0}^{n-1}X^{n-1-a}\left((1-\vs(X))\pioh_{A L}\sn X^a f+\vs(X) LX^a f\right)\label{3.11.1.24}.
\end{align}
With $a=0,\cdots n-1$ and $n\le 3$ using (\ref{10.10.2.23}) and (\ref{8.29.9.21})
\begin{equation*}
X^{n-1-a}L X^a \Phi=\left\{\begin{array}{lll}
O(\l t\r^{-2}\log \l t\r^{\f12\M}(\La_0+\Delta_0^\frac{5}{4}))_{L_u^2 L_\omega^4}, n\le 2\\
O(\l t\r^{-1}\log \l t\r^{\f12\M}(\La_0+\Delta_0^\frac{5}{4}))_{L^2_\Sigma}, n\le 3.
\end{array}
\right.
\end{equation*}
Using (\ref{5.21.1.21}), (\ref{3.25.1.22}) and Proposition \ref{7.15.5.22}, we deduce
\begin{align*}
X^{n-1-a}(\pioh_{A L}\sn X^a \Phi)=\left\{\begin{array}{lll}
O(\l t\r^{-\frac{11}{4}+2\delta}\Delta_0^\frac{3}{2})_{L_u^2 L_\omega^4}, n\le 2\\
O(\l t\r^{-3+\frac{\max(n-2,0)}{2}+2\delta}\Delta_0^\frac{3}{2})_{L^2_u L_\omega^2}, n\le 3.
\end{array}\right.
\end{align*}
Hence the proof of Lemma \ref{8.1_com} is complete.
\end{proof}
\subsection{Top order comparison}
\begin{lemma}
For $X^2=X_2 X_1$ and $f=\Phi$, there hold with $Y=\bN, L$
\begin{align}
\|\tir\sn_{X_2}\sn_{X_1} \sn^2_{AB}f\|_{L_\Sigma^2} &\les \|\sn_X^{\le 2}\sn \Omega f\|_{L^2_\Sigma}+\|\sn_X^{\le 2}\sn f\|_{L^2_\Sigma}+\l t\r^{-\frac{3}{2}+2\delta}\Delta_0^\frac{3}{2}(\log \l t\r)^\f12\label{7.15.1.21}\\
\|\tir (\sn_Y \sn_A f)_{;X_1X_2}\|_{L^2_\Sigma}&\les\|X^2 Y(\Omega f), \sn_X^{\le 2}\sn f, X^{\le 1}Y\Omega f\|_{L^2_\Sigma}\nn\\
&+\Delta_0^\f12\l t\r^{-\frac{3}{4}+\delta}\|\bb^\f12\tir Y\Omega f\|_{L_u^2 L_\omega^4}\nn\\
&+(\log \l t\r)^\f12\l t\r^{-\frac{3}{2}+2\delta+\f12\max(\vs(Y),0)}\Delta_0^\frac{3}{2}.\label{7.15.2.21}
\end{align}
\end{lemma}
\begin{proof}
We start with computing for $X_1, X_2\in \{\Omega, S\}$, with $Z=\bN, L, e_B$ i.e. in null tetrad, that
\begin{align}
\Omega^A (\sn_Z \sn_A f_{;X_1 X_2})
&=\sn_{X_2}\big(\sn_{X_1}(\Omega^A \sn_Z \sn_A f)-\sn_{X_1} \Omega^A \sn_Z \sn_A f\big)-\sn_{X_2} \Omega^A \sn_{X_1} \sn_Z \sn_A f\nn\\
&=(\Omega^A \sn_Z \sn_A f)_{;X_1X_2}-\sn_{X_2}(\sn_{X_1} \Omega^A \sn_Z \sn_A f)-\sn_{X_2} \Omega^A \sn_{X_1} \sn_Z\sn_A f\nn\\
&=(\sn_Z(\Omega f))_{;X_1 X_2}-(\sn_Z \Omega^A \sn_A f)_{;X_1 X_2}-\sn_{X_2}(\sn_{X_1} \Omega^A\sn_Z\sn_A f)\nn\\
&-\sn_{X_2} \Omega^A \sn_{X_1} \sn_Z\sn_A f.\label{2.25.1.23}
\end{align}
Hence, by using Lemma \ref{3.17.2.22}, we write
\begin{align*}
&\Omega^A (\sn_Z \sn_A f_{;X_1 X_2})\\
&=\sn_X^2 \sn_Z\Omega f+O(1)\sn_X^2\sn f+\sn_X\sn_Z\Omega f+\sn_X \sn f(1+O(\l t\r^{\delta}\Delta_0^\f12)_{L_\omega^4})\\
&+Z\Omega f( 1+O(\Delta_0^\f12\l t\r^{\delta-\frac{3}{4}})_{L^4_\omega})+(O(1)+O(\l t\r^{-\f12+\delta+\f12\max(-\vs(Z), 0)}\Delta_0)_{L^2_u L_\omega^2}\\
&+O(\l t\r^{\delta-\frac{3}{4}+\frac{3}{4}(\max(-\vs(Z),0))}\Delta_0^\f12)_{L_\omega^4})\sn f.
\end{align*}
With $f=\Phi$, it follows by using Proposition \ref{7.15.5.22} that
\begin{align*}
\|\tir (\sn_Z \sn_A f_{;X_1 X_2})\|_{L^2_\Sigma}&\les \|\sn_X^{1+\le 1}\sn_Z\Omega f+\sn_X^2\sn f+\sn_X\sn f+Z\Omega f+\sn f\|_{L^2_\Sigma}\\
&+\log \l t\r^\f12\l t\r^{-\frac{3}{2}+2\delta+\f12\max(-\vs(Z),0)}\Delta_0^\frac{3}{2}+\|\bb^\f12 Z\Omega f\|_{L_u^2 L_\omega^4} \Delta_0^\f12 \l t\r^{\delta+\frac{1}{4}}\\
&+\l t\r^{2\delta-\frac{7}{4}+\frac{3}{4}(\max(-\vs(Z),0))}\log \l t\r^\f12\Delta_0^\frac{3}{2}.
\end{align*}
(\ref{7.15.1.21}) follows immediately in view of Proposition \ref{7.15.5.22}. We can similarly obtain (\ref{7.15.2.21}) with $Z=\bN, L$. 
\end{proof}

Next under the bootstrap assumptions (\ref{3.12.1.21})-(\ref{6.5.1.21}), we derive the top order derivative estimates of $\Phi$ in comparison with the top order energies. 
\begin{proposition}\label{7.16.1.21}
Let $X, X_1, X_2\in \{\Omega, S\}$ and $0<t<T_*$. We have
\begin{align}
\|\tir(\ell\bd^2_{LL}\Phi, LL\Phi)_{;X_1 X_2}\|_{L^2_\Sigma}&\les \l t\r^{-1}\log \l t\r^{\f12(\M+7\ell)}(\Delta_0^\frac{5}{4}+\La_0)\nn\\
&+\l t\r^{-\f12}W_1[X_2 X_1 S\Phi]^\f12(t), \ell=0,1.\label{7.28.2.21}\\
\|\tir(\sn^2\Phi)_{;X^2}\|_{L^2_\Sigma}&\les\l t\r^{-\f12}\sum_{\vs(Y)=\vs^-(X^2)}W_1[Y\Omega^2\Phi]^\f12(t)\nn\\
&+\l t\r^{-1}\log \l t\r^{\f12\M}(\l t\r^\delta\Delta_0^\frac{5}{4}+\La_0)\label{7.17.2.21}\\
\|\tir (\bd^2_{LA} \Phi, \sn_L \sn_A \Phi)_{;X_1X_2}\|_{L^2_\Sigma}&\les \l t\r^{-\f12} W_1[X_2 X_1 \Omega \Phi]^\f12(t)\nn\\&+\l t\r^{-1}\log \l t\r^{\f12\M}(\La_0+\l t\r^\delta\Delta_0^\frac{5}{4}).\label{8.15.2.21}
\end{align}
\begin{equation}\label{3.14.1.24}
\begin{split}
\|\tir\sn_{X_2}\sn^2 X_1\Phi\|_{L^2_\Sigma}&\les \l t\r^{-\f12}W_1[\Omega^2 X_1\Phi]^\f12(t)+\l t\r^{-1}\log \l t\r^{\f12\M}(\La_0+\l t\r^{\delta}\Delta_0^\frac{5}{4})\\
\|\tir \sn_{X_2}(\bd^2_{LA}, \sn_L \sn_A)X_1\Phi\|_{L^2_\Sigma}&\les \l t\r^{-\f12} W_1[X_2\Omega X_1\Phi]^\f12(t)+\l t\r^{-1}\log \l t\r^{\f12\M}(\La_0+\l t\r^{\delta}\Delta_0^\frac{5}{4})\\
\|\tir X_2(\bd_{LL}^2, LL)X_1\Phi\|_{L^2_\Sigma}&\les \l t\r^{-\f12} W_1[X_2S X_1\Phi]^\f12(t)\\
&+\l t\r^{-1}\log \l t\r^{\f12\M}(\La_0+(\log \l t\r)^{1-\vs(X_2)}\Delta_0^\frac{5}{4}).
\end{split}
\end{equation}

\begin{equation}\label{3.14.2.24}
\begin{split}
&\|\tir\sn^2 X^2\Phi\|_{L^2_\Sigma}\les \l t\r^{-\f12}W_1[\Omega X^2 \Phi]^\f12+\l t\r^{-1}\log \l t\r^{\f12\M}(\La_0+\l t\r^\delta\Delta_0^\frac{5}{4}) \\  
&\|\tir(\bd^2_{LL}, LL)X^2\Phi\|_{L^2_\Sigma}\les\l t\r^{-\f12} W_1[SX^2\Phi]^\f12(t)+\l t\r^{-1}\log \l t\r^{\f12\M}(\La_0+\Delta_0^\frac{5}{4}) \\
&\|\tir(\bd^2_{A L}, \sn_L \sn_A) X^2\Phi\|_{L^2_\Sigma}\les \l t\r^{-\f12}W_1[\Omega X^2\Phi]^\f12(t)+\l t\r^{-1}\log \l t\r^{\f12\M}(\La_0+\l t\r^\delta\Delta_0^\frac{5}{4}). 
\end{split}
\end{equation}
\begin{align}
&\sn_X^{2-a}(\bd^2_{LL}, \bd^2_{LA}, \sn^2)X^a\Phi=O(\l t\r^{-\frac{3}{2}+\delta}\Delta_0)_{L^2_\Sigma}, a=0, 1,2\label{3.6.5.24}\\
&\sum_{a\le 2}\|\tir\sn_X^{2-a}\sn \sn^a_X\Phi\|_{L_u^2 L_\omega^4}+\sum_{a=0}^1\|\sn_X^{1-a}\sn_S\sn^a_X\Phi\|_{L_u^2 L_\omega^4}\les \Delta_0 \l t\r^{-\frac{3}{4}+\delta}\label{3.12.2.24}\\
&\sn_X^{2-a}\bd^2_{\Lb A}\sn_X^a \Phi^\mu=O(\l t\r^{-1})\min(\mu,1)\fB+O(\l t\r^{-1+\delta}\Delta_0)_{L^2_\Sigma}, a=0,1,2.\label{3.3.4.24}
\end{align}
\begin{align}
\sta{\Omega^3, \Lb}{\Phi^\mu}&=\min(1,\mu)O(\fB)+O(\l t\r^\delta\Delta_0)_{L^2_\Sigma}\label{3.14.3.24}\\
\sta{X^2,\Lb}{\Phi^\mu}-O(\min(\vs^-(X^2)+\mu, 1))(\fB)&=\left\{\begin{array}{lll}
O(\log \l t\r^4(\La_0+\Delta_0^\frac{5}{4}))_{L^2_\Sigma}\\
O(\l t\r^{-1+\f12\delta}\log \l t\r^2\Delta_0)_{L_u^2 L_\omega^4}\end{array}\right.
\label{3.5.4.24}
\end{align}
With $a=0,1$
\begin{align}
&\|\sn_X^{1-a}\bd^2_{\Lb A}X^a\bT\varrho, \l t\r^{-1}X^{1-a}(\Lb\Omega, \Omega\Lb) X^a\bT\varrho\|_{L^2_\Sigma}\les \l t\r^{-1+\delta}\Delta_0\label{3.6.4.24}\\
&\|\sn_\Omega^{1-a}(\bd_{LA}^2, \sn^2, \bd^2_{LL}) \Omega^a\bT\varrho, \l t\r^{-2} \Omega^{1-a}X^2\Omega^a\bT\varrho\|_{L^2_\Sigma}\les \l t\r^{-2+\delta}\Delta_0\label{3.6.3.24}
\end{align}
\end{proposition}
In the sequel, we divide the proof into two parts, with the second part requiring more delicate commutator estimates.
\begin{proof}[Proof of (\ref{7.28.2.21})-(\ref{3.12.2.24})]
We apply (\ref{3.4.2.24}) and (\ref{8.29.9.21}) to derive
\begin{align}
\|\tir X_2 X_1 LL \Phi\|_{L^2_\Sigma}&\les \|\tir L X^2 Lf\|_{L_\Sigma^2}+\l t\r^{-1}\log \l t\r^{\f12\M}(\La_0+\Delta_0^\frac{5}{4}),\nn\\
&\les \|L X_2 X_1 S f\|_{L^2_\Sigma}+\l t\r^{-1}\log \l t\r^{\f12\M}(\La_0+\Delta_0^\frac{5}{4})
.\label{7.15.3.21}
\end{align}
We further note by using (\ref{7.03.3.21})
\begin{align*}
X_2 X_1(\bd^2_{LL}f)&=X_2 X_1(LL f)+X_2 X_1(k_{\bN\bN}Lf).
\end{align*}
Applying (\ref{8.9.4.22}) to $F=L\Phi$, also using (\ref{8.29.9.21}), (\ref{10.10.2.23}) and (\ref{7.25.2.22})  yields
\begin{align}
X^2(k_{\bN\bN}L\Phi)&=O(\bb^{-1} \l t\r^{-1})\sn_X^{\le 2}L\Phi+O(\Delta_0 \l t\r^{-1}\log \l t\r)_{L_\omega^4}\bb^{-1} X L\Phi\nn\\
&+O((\log \l t\r)^2\Delta_0)_{L^2_\Sigma}L\Phi\nn\\
&=O(\l t\r^{-2}(\log \l t\r)^{\f12(\M+7)}(\Delta_0^\frac{5}{4}+\La_0))_{L^2_\Sigma}\label{3.20.1.24}
\end{align} 
where we used $\M\ge 9$. 
Combining the above two estimates, we conclude (\ref{7.28.2.21}) with the help of (\ref{8.29.9.21}).

Next, we apply (\ref{7.15.1.21}), (\ref{4.22.4.22}),  (\ref{8.29.9.21}) and (\ref{10.10.2.23}) , we obtain
\begin{align*}
\|\tir\sn_{X_2}\sn_{X_1} \sn^2_{AB}\Phi\|_{L_\Sigma^2} &\les \|\sn_{X_2}\sn_{X_1}\sn \Omega \Phi\|_{L^2_\Sigma}+\l t\r^{-1}\log \l t\r^{\f12\M}(\l t\r^\delta\Delta_0^\frac{5}{4}+\La_0)\\
&\les \l t\r^{-1}\|X^2\Omega^2 \Phi\|_{L^2_\Sigma}+\l t\r^{-1}\log \l t\r^{\f12\M}(\l t\r^\delta\Delta_0^\frac{5}{4}+\La_0).
\end{align*}
For the leading term on the right-hand side,
for $X^2=X_2 X_1$, if $\vs(X_1)>\vs(X_2)$, we employ (\ref{7.17.6.21}) to commute these two vector fields, which gives
\begin{align*}
|X^2\Omega^2\Phi|\les |X_1 X_2\Omega^2\Phi|+\l t\r^{-\frac{3}{4}+\delta}\Delta_0^\f12|\Omega^{\le 1}\Omega^2\Phi|, \mbox{ if }\vs(X_1)>\vs(X_2).
\end{align*}
Otherwise, there is no need to commute the two vector fields.
We then summarize by using (\ref{8.29.9.21})
\begin{align*}
\|\tir\sn_{X_2}\sn_{X_1} \sn^2_{AB}\Phi\|_{L_\Sigma^2}&\les \l t\r^{-\f12}\sum_{\vs(Y)=\vs^-(X^2)}W_1[Y\Omega^2\Phi]^\f12(t)+\l t\r^{-1}\sum_{X=S, \Omega}\|X^{\le 1}\Omega^2\Phi\|_{L^2_\Sigma}\\
&+\l t\r^{-1}\log \l t\r^{\f12\M}(\l t\r^\delta\Delta_0^\frac{5}{4}+\La_0)\\
& \les\l t\r^{-\f12}\sum_{\vs(Y)=\vs^-(X^2)}W_1[Y\Omega^2\Phi]^\f12(t)+\l t\r^{-1}\log \l t\r^{\f12\M}(\l t\r^\delta\Delta_0^\frac{5}{4}+\La_0).
\end{align*}

In view of (\ref{7.03.3.21}) and Proposition \ref{7.15.5.22}, we derive
\begin{align*}
\|\sn_{X_2}\sn_{X_1}(\bd^2_{LA}\Phi)\|_{L^2_\Sigma}&\les \|(\sn_L\sn_A \Phi)_{;X_1 X_2}\|_{L^2_\Sigma}+\sum_{a=0}^2\|\sn_X^{2-a}\bA_{g,1}X^a L\Phi\|_{L^2_\Sigma}\\
&\les  \|(\sn_L\sn_A \Phi)_{;X_1 X_2}\|_{L^2_\Sigma}+\l t\r^{-3+2\delta}\Delta_0^\frac{3}{2}.
\end{align*}
Applying (\ref{7.15.2.21}) to the first term on the right-hand side, and using (\ref{8.29.9.21}) gives
\begin{align*}
\|\sn_{X_2}\sn_{X_1}(\bd^2_{LA}\Phi, \sn_L \sn_A \Phi)\|_{L^2_\Sigma}\les \l t\r^{-1} \|X_2 X_1 L\Omega\Phi\|_{L^2_\Sigma}+\l t\r^{-2}\log \l t\r^{\f12\M}(\La_0+\l t\r^\delta\Delta_0^\frac{5}{4}).
\end{align*}
We have obtained by the proof of (\ref{3.4.2.24}) that 
\begin{equation*}
\|[L, X^2]\Omega\Phi\|_{L^2_\Sigma}\les \l t\r^{-1}\log \l t\r^{\f12\M}(\La_0+\Delta_0^\frac{5}{4}).
\end{equation*}
Combining the above two estimates yields 
\begin{align*}
\|\sn_{X_2}\sn_{X_1}(\bd^2_{LA}\Phi, \sn_L \sn_A \Phi)\|_{L^2_\Sigma}\les \l t\r^{-1} \|L X_2 X_1 \Omega\Phi\|_{L^2_\Sigma}+\l t\r^{-2}\log \l t\r^{\f12\M}(\La_0+\l t\r^\delta\Delta_0^\frac{5}{4}).
\end{align*}
This implies (\ref{8.15.2.21}) by using (\ref{8.29.9.21}) again. As a direct consequence of (\ref{7.28.2.21})-(\ref{8.15.2.21}), we obtained the case $a=0$ in (\ref{3.6.5.24}) by using (\ref{3.12.1.21}).

Next we apply (\ref{7.03.5.21})-(\ref{7.5.5.21}) to $(X, f)=(X_2, X_1\Phi)$, using (\ref{7.03.3.21}),
\begin{align*}
\tir^2\sn_{X_2}\sn^2 X_1\Phi&=O(1) X_2^{\le 1}\Omega^{1+\le 1}X_1\Phi+O(\l t\r^{-\frac{3}{4}+\delta}\Delta_0^\f12)_{L_\omega^4}\Omega X_1\Phi\\
\tir^2 \sn_{X_2}(\bd^2_{LA}, \sn_L \sn_A)X_1\Phi&=O(1) (SX_2^{\le 1}\Omega X_1\Phi+X_2^{\le 1}\Omega X_1\Phi)+O(\l t\r^{-\frac{3}{4}+\delta}\Delta_0^\f12)_{L_\omega^4}\Omega X_1\Phi\\
&+O(1)\tir^2 \sn_{X_2}^l \bA_{g,1}\sn_{X_2}^{1-l}LX_1\Phi, \, l=0,1\\
\tir^2 X_2(\bd_{LL}^2 X_1\Phi, LLX_1\Phi)&=O(\tir)(L+h)X_2S^{\le 1}X_1\Phi+O(1)X_2X_1\Phi+O(1)X_2^{\le 1} S X_1\Phi\\
&+(1-\vs(X_2))O(\log \l t\r\Delta_0)_{L_\omega^4}\bb^{-1}S X_1\Phi
\end{align*}
where we also used (\ref{8.23.2.23}) to derive the last estimate. 

Using Proposition \ref{7.15.5.22}, (\ref{8.29.9.21}) and (\ref{10.10.2.23}), we derive
\begin{align*}
\|\tir^2\sn_{X_2}\sn^2 X_1\Phi\|_{L^2_\Sigma}&\les\|X_2^{\le 1}\Omega^{1+\le 1}X_1\Phi\|_{L^2_\Sigma}+\l t\r^{-\frac{3}{4}+2\delta}(\log \l t\r)^\f12\Delta_0^\frac{3}{2}\\
\|\tir^2 \sn_{X_2}(\bd^2_{LA}, \sn_L \sn_A)X_1\Phi\|_{L^2_\Sigma}&\les \|S^{\le 1} X_2^{\le 1}\Omega X_1\Phi\|_{L^2_\Sigma}+\l t\r^{-\frac{3}{4}+2\delta}(\log \l t\r)^\f12\Delta_0^\frac{3}{2}\\
\|\tir^2 X_2(\bd_{LL}^2, LL)X_1\Phi\|_{L^2_\Sigma}&\les\|\tir (L+h)X_2 S X_1\Phi\|_{L^2_\Sigma}\\
&+(\log \l t\r)^{\f12\M}\Big((\log \l t\r)^{1-\vs(X_2)}\Delta_0^\frac{5}{4}+\La_0\Big).
\end{align*}
Hence we can obtain (\ref{3.14.1.24}) with the help of (\ref{8.29.9.21}); and the case $a=1$ in (\ref{3.6.5.24}) is proved. 

Next we apply (\ref{7.04.7.21}) and (\ref{7.03.4.21}) to $f=X^2\Phi$ to deduce
\begin{align*}
&\tir^2 \sn^2 X^2\Phi =O(1)(\tir\sn)^{\le 1}\Omega X^2 \Phi\\  
\displaybreak[0]
&|\tir\bd^2_{LL}X^2\Phi|\les |\tir L(LX^2\Phi)|+|LX^2\Phi| \\
\displaybreak[0]
&\tir |\bd^2_{A L} X^2\Phi|\les |L \Omega X^2\Phi|+|\sn X^2\Phi|+\l t\r^{-1+\delta}\Delta_0|LX^2\Phi|
\end{align*}
from which, also applying (\ref{8.29.9.21}), we conclude (\ref{3.14.2.24}) and the $a=2$ case in (\ref{3.6.5.24}).

(\ref{3.12.2.24}) can be obtained by using (\ref{4.22.4.22}), (\ref{9.1.1.21}) and (\ref{3.12.1.21}).
\end{proof}
\begin{proof}[Proof of (\ref{3.3.4.24})-(\ref{3.6.3.24})]
$\bf\bullet$ Control of $\sn_X^{2-a} \bd^2_{A\Lb}X^a\Phi$. In this part, we prove (\ref{3.3.4.24}). Using (\ref{7.03.3.21}), we write 
\begin{align*}
\sn_X^{2-a} \bd^2_{A\Lb} X^a f&=\sn_X^{2-a}(\sn_A\sn_\Lb-\chib\sn+\bA_{g,1}\Lb)X^a f\\
 &=\sn_X^{2-a}(\sn_\Lb\sn+\bA_{g,1}L+\ze\c\bN)X^a f, \, a=0,1,2,
\end{align*}
with the lower order terms satisfying
\begin{align}\label{7.28.2.24}
\begin{split}
\sum_{a=0}^2\|\sn_X^{2-a}(\chib\sn+\bA_{g,1}(\Lb+ L)+\ud \bA L)X^a \Phi\|_{L_u^2 L_\omega^2}\les \l t\r^{-3+2\delta}\Delta_0\\
\sum_{a=0}^2\|\sn_X^{2-a}(\ud\bA\Lb X^a \Phi)\|_{L^2_\Sigma}\les \l t\r^{-1+\delta}\Delta_0.
\end{split}
\end{align}
where the first estimate is obtained by using Proposition \ref{7.15.5.22}, Lemma \ref{8.1_com} (2) and (\ref{8.23.1.23}),
which are negligible terms compared with the other terms to be treated shortly. The second estimate can be derived by using (\ref{3.10.7.24}), (\ref{3.12.4.24}) together with (\ref{8.23.2.23}), (\ref{2.9.2.24}), (\ref{3.16.1.22}) and $\bb\Lb\Phi=O(\tir^{-1})$.

Moreover using (\ref{9.8.2.22}), we have
\begin{align*}
\Omega^A\sn_X^{2-a}\sn_A \Lb X^a f&=\Omega^A(\sn_X^2\sn \Lb f, \sn_X\sn_A\Lb Xf, \sn_A \Lb X^2f)\\
&\approx \Omega^A \sn_X^2 \sn\Lb f,  X^{\le 1}\Omega \Lb Xf, \Omega\Lb X^2 f.
\end{align*}
Applying (\ref{4.22.4.22}) and (\ref{9.8.2.22}) to the first term on the right, we obtain 
\begin{align*}
\Omega^A \sn_X^2 \sn_A\Lb f=X^2\Omega \Lb f+O(1) X^{\le 1}\Omega \Lb f+O(\l t\r^{-\frac{3}{4}+\delta}\Delta_0^{\f12+\frac{1}{p}})_{L^p_u L_\omega^4} \Omega \Lb f, \, p=2,\infty.
\end{align*}
Noting that with $f=\Phi$ the lower order terms have been treated in (\ref{3.10.7.24}) with the help of (\ref{3.12.4.24}), (\ref{8.23.2.23}), (\ref{8.2.1.22}) and (\ref{8.29.9.21}) 
\begin{align}\label{3.12.9.24}
X^{\le 1}\Omega \Lb \Phi^\mu+O(\l t\r^{-\frac{3}{4}+\delta}&\Delta_0^{\f12+\frac{1}{p}})_{L_u^p L_\omega^4} \Omega \Lb \Phi=O(1)\min(\mu,1)\fB\nn\\
&+O(\log \l t\r^4(\La_0+\Delta_0^\frac{5}{4})+\l t\r^{-\frac{3}{4}+\delta}\Delta_0)_{L^2_\Sigma}.
\end{align}
Combining the above estimates, we have for $f=\Phi^\mu$
\begin{align}
&\sum_{a=0}^2\Omega^A\sn_X^{2-a}\sn_A \Lb X^a f\nn\\
&\approx \sum_{a=0}^2X^{2-a}\Omega \Lb X^a f+O(1) X^{\le 1}\Omega \Lb f+O(\l t\r^{-\frac{3}{4}+\delta}\Delta_0^{\f12+\frac{1}{p}})_{L_u^p L_\omega^4} \Omega \Lb f\nn\\
&=\sum_{a=0}^2X^{2-a}\Omega \Lb X^a \Phi^\mu+O(1)\min(\mu,1)\fB+O(\log \l t\r^4(\La_0+\Delta_0^\frac{5}{4})+\l t\r^{-\frac{3}{4}+\delta}\Delta_0)_{L^2_\Sigma}.\label{3.13.1.24}
\end{align}
For the leading terms, we list all the cases of $X^{2-a}\Omega\Lb X^a\Phi$
\begin{equation} \label{3.12.5.24}
\begin{array}{|l|l|l|}
\hline
  a=2 & a=1 & a=0 \\
 \hline
 \Omega \Lb \Omega^2\Phi &\Omega^2\Lb \Omega\Phi  &\Omega^3\Lb\Phi \\
\hline
\Omega\Lb\Omega S\Phi & \Omega^2\Lb S\Phi &  \Omega S \Omega \Lb\Phi\\
\hline
\Omega\Lb S\Omega\Phi  & S\Omega\Lb \Omega\Phi & S\Omega^2\Lb\Phi \\
\hline
\Omega\Lb S^2\Phi &S\Omega\Lb S\Phi&S^2\Omega\Lb\Phi\\
\hline
\end{array}
\end{equation}
Let us first consider the terms in (\ref{3.12.5.24}) except the first row. 
We will proceed the following commutations, keeping the rest of the terms as they were
\begin{align}\label{3.12.6.24}
\begin{split}
&\Omega\Lb \Omega S\Phi=\Omega[\Lb, \Omega]S\Phi+\Omega^2\Lb S\Phi, \Omega S\Omega\Lb\Phi=\Omega S\Lb \Omega\Phi+\Omega S[\Omega, \Lb]\Phi,\\
&S\Omega \Lb \Omega\Phi=[S, \Omega]\Lb \Omega\Phi+\Omega S\Lb \Omega\Phi, S\Omega^2\Lb\Phi=[S, \Omega^2]\Lb\Phi+\Omega^2 S\Lb \Phi, \\
&S^2\Omega\Lb\Phi=S^2\Lb \Omega\Phi+S^2[\Omega, \Lb]\Phi.
\end{split}
\end{align}
Hence, to treat the terms in (\ref{3.12.5.24}) except the first row, based on the two lists in the above, apart from treating the commutators, it suffices to prove for $X^2=X^{2-a}X^a$ with $X^{2-a}$ and $X^a$ appeared in the identity below, and $\vs^-(X^2)=0$, 
\begin{equation}\label{3.12.3.24}
\begin{split}
X^{2-a}(\Lb S, S\Lb)& X^a \Phi^\mu-X^{2-a}(\tir \sD X^a \Phi^\mu)\\
&=O(1)\min(\mu, 1) \fB+O((\log \l t\r)^5(\La_0+\Delta_0^\frac{5}{4}))_{L^2_\Sigma},\, a=0,1
\end{split}
\end{equation}
where 
\begin{align*}
\|X^{2-a}(\tir\sD X^a\Phi)\|_{L^2_\Sigma}&\les \l t\r^{-\f12+\delta}\Delta_0, a=0,1.
\end{align*}
Here the estimate in the last line is derived by using  (\ref{7.17.2.21}), (\ref{3.14.1.24}), (\ref{3.12.1.21}) and Proposition \ref{8.29.8.21}.

To prove (\ref{3.12.3.24}),
 we first consider $X^{2-a}\Lb S X^a\Phi$ with $a=0,1$, $X^2=\Omega^2, \Omega S, S\Omega$ in the sequel. It holds for scalar functions $f$, due to (\ref{8.25.1.21}), that
\begin{align}
X^{2-a}\Lb S X^a f&=X^{2-a}(\tir h\Lb X^a f)+X^{2-a}(\tir \sD X^a f+\tir \Box_\bg X^a f)\label{3.12.1.24}\\
&+X^{2-a}((\mho+k_{\bN\bN})SX^a f+\tir\sn X^a f\c \ze)\nn.
\end{align}

Recall from (\ref{8.30.3.21}) and (\ref{8.30.3.21+}) that with $a=0,1$ 
\begin{align*}
\|X^{2-a}(\tir\Box_\bg X^a \Phi)\|_{L^2_\Sigma}\les\left\{\begin{array}{lll}
\l t\r^{-\frac{1}{2}}(\La_0+\Delta_0^\frac{5}{4}\l t\r^\delta), X^2=\Omega^2\\
\l t\r^{-1+\delta}(\log \l t\r^{\f12\M}\La_0+\l t\r^{\frac{1}{4}+\delta}\Delta_0^\frac{5}{4}), X^2=S\Omega, \Omega S.
\end{array}\right.
\end{align*}
With $a=0,1$, by using (\ref{3.12.2.24}), Corollary \ref{9.2.5.23}, estimates of $\ze$ in Proposition \ref{7.15.5.22}, we bound
\begin{align*}
\sn_X^{2-a}(\sn X^a\Phi\c \ze)&=\sn_X^{2-a}(\sn X^a\Phi)\ze+\sn X^a\Phi \sn_X^{2-a}\ze+\sn_X^{1-a}\sn X^a \Phi\c \sn_X\ze\\
&=O(\l t\r^{-3+2\delta}\Delta_0^\frac{3}{2})_{L^2_u L_\omega^2}.
\end{align*}
With $a=0,1$, using Proposition \ref{7.22.2.22}, (\ref{8.23.1.23}) and Proposition \ref{8.29.8.21}, we have
\begin{align*}
X^{2-a}\big((\mho+k_{\bN\bN})S X^a\Phi\big)&=X^{2-a}(\mho+k_{\bN\bN})SX^a\Phi+(\mho+k_{\bN\bN})X^{2-a}S X^a\Phi\\
&+X^{1-a}(\mho+k_{\bN\bN})XSX^a\Phi\\
&=O(\l t\r^{-1}\log \l t\r^{\f12\M}\La_0+\l t\r^{-1+2\delta}\log \l t\r^\f12\Delta_0^\frac{5}{4})_{L^2_\Sigma}.
\end{align*}
It only remains to consider the first term on the right-hand side of (\ref{3.12.1.24}). We first expand it as
\begin{align*}
X^{2-a}(\tir h\Lb X^a \Phi)=X^2(\tir h) \Lb\Phi+X(\tir h)(\Lb X\Phi+X\Lb\Phi)+\tir h \sta{X^2, \Lb}\Phi.
\end{align*}
Using (\ref{2.17.2.24}) and (\ref{L2conndrv'}), we derive
\begin{align*}
X^2(\tir h)\Lb\Phi&=\left\{\begin{array}{lll}
O(\l t\r^{-\f12+\delta}\Delta_0^\frac{5}{4})_{L^2_\Sigma}, X^2=\Omega^2\\
O(\l t\r^{-1}\log \l t\r^{\f12\M+1}(\Delta_0^\frac{5}{4}\l t\r^\delta+\La_0))_{L^2_\Sigma}, X^2=S\Omega, \Omega S.
\end{array}\right.
\end{align*}
It follows by using (\ref{8.23.1.23}), (\ref{8.8.6.22}) and (\ref{2.18.2.24}) that
\begin{align*}
X(\tir h) (X\Lb \Phi+\Lb X\Phi)=O(\l t\r^{-1}\log \l t\r^{\f12\M+1}(\La_0+\log\l t\r^\f12\l t\r^{2\delta}\Delta_0^\frac{5}{4}))_{L^2_\Sigma}
\end{align*}
Finally using (\ref{3.10.7.24}), we have
\begin{equation*}
\tir h \sta{X^2, \Lb}{\Phi^\mu}= O(1)\min(\mu,1) \fB+O\Big((\log \l t\r)^4(\La_0+\Delta_0^\frac{5}{4})\Big)_{L^2_\Sigma}.
\end{equation*} 
Combining the above estimates, noting that the last one dominates, we conclude the first estimate in (\ref{3.12.3.24}).

For the other case in (\ref{3.12.3.24}), using (\ref{6.30.2.19}) instead, instead of (\ref{3.12.1.24}) we write for $f=\Phi$,
\begin{align}\label{3.22.5.24}
X^{2-a} S\Lb X^a f&=X^{2-a}(\tir (h-k_{\bN\bN})\Lb X^a f)+X^{2-a}(\tir \sD X^a f+\tir \Box_\bg X^a f)\\
&+X^{2-a}(\hb SX^a f+\tir\sn X^a f\c \zb)\nn.
\end{align}
We decompose $\hb=-h+\fB$. Note the term $X^{2-a}(\fB SX^a \Phi)$ can be included in the last term of the first line; while the term $X^{2-a}(\tir h LX^a \Phi)$ verifies better estimates than $X^{2-a}(\tir h \Lb X^a\Phi)$. Moreover,
noting that $\zb$ verifies better estimates than $\ud \bA$, the only additional term we need to take care of is $X^{2-a}(\tir k_{\bN\bN}\Lb X^a\Phi^\mu), a=0,1$. Using (\ref{3.10.7.24}), (\ref{7.29.1.22}), (\ref{8.23.2.23}) and (\ref{3.12.4.24}), we infer
\begin{align*}
X^{2-a}(\tir k_{\bN\bN}\Lb X^a\Phi^\mu)&=\tir k_{\bN\bN}\sta{X^2, \Lb}{\Phi^\mu}+X(\tir k_{\bN\bN})(\Lb X\Phi^\mu+X\Lb\Phi^\mu)+X^2(\tir k_{\bN\bN})\Lb\Phi^\mu\\
&= O(\log \l t\r^4(\La_0+\Delta_0^\frac{5}{4}))_{L^2_\Sigma}+X_2(\tir k_{\bN\bN}) (\min(\mu+\vs(X_1),1)O(\fB)\\
&+(1-\vs(X_1))O(\l t\r^{-1}\log \l t\r(\log \l t\r^3\Delta_0^\frac{5}{4}+\La_0))_{L_u^2 L_\omega^4})\\
&=O(\log \l t\r^5(\La_0+\Delta_0^\frac{5}{4}))_{L^2_\Sigma}+\min(\mu,\vs^+(X^2))O(\fB)
\end{align*}
where $X^2=X_2 X_1$ or $X_1 X_2$, and we used  $\vs^-(X^2)=0$ in the above.  We proved the second case in (\ref{3.12.3.24}).

Next we treat the commutators in (\ref{3.12.6.24}), which are 
\begin{align}\label{3.12.8.24}
&X^{2-a}[\Lb, \Omega]X^a\Phi, \mbox{ with } X^2=\Omega S, S^2; \mbox{ and } [S,\Omega^{2-a}]\Lb \Omega^a\Phi,\, a=0,1.
\end{align}
For the first term, we apply (\ref{5.13.10.21}), with $X^2=\Omega S, SS$ and $a=0,1$, to write 
\begin{align*}
X^{2-a}[\Lb, \Omega]X^a\Phi&=X^{2-a}(\Omega\log \bb\bN X^a\Phi+\pioh_{A\Lb}\sn X^a\Phi)\\
&=X^2(\Omega \log \bb \bN \Phi)+X(\Omega\log \bb \bN X\Phi)+\sum_{a=0}^1X^{2-a}(\pioh_{A\Lb}\sn X^a\Phi)\\
&=O(\fB)X^{1+\le 1}\Omega \log \bb+\Omega\log \bb(O(\fB)+O(\l t\r^{-1}\log \l t\r\Delta_0)_{L_u^2 L_\omega^4})\\
&+S\Omega\log \bb (\bN \Omega\Phi+\Omega\bN\Phi)+\sum_{a=0}^1X^{2-a}(\pioh_{A\Lb}\sn X^a\Phi)
\end{align*}
where we also applied (\ref{3.10.7.24}) and (\ref{3.12.4.24}) to obtain the last estimate and the vector-fields $X$ appearing in the last term of the formula keep the original ordering. Using (\ref{8.11.1.24}), (\ref{2.20.2.24}) and Lemma \ref{5.13.11.21} (5), we obtain
\begin{equation*}
S\Omega\log \bb (\bN \Omega\Phi+\Omega\bN\Phi)=O(\log \l t\r(\La_0+\log \l t\r^3\Delta_0^\frac{5}{4})_{L^2_\Sigma}.
\end{equation*}
 Using (\ref{5.21.1.21}) and (\ref{3.25.1.22}), Corollary \ref{9.2.5.23}, (\ref{3.6.2.21}) and (\ref{3.12.2.24}) we infer 
\begin{align*}
\sum_{a=0}^1\|X^{2-a}(\pioh_{A\Lb}\sn X^a\Phi)\|_{L^2_u L_\omega^2}\les \l t\r^{-\frac{11}{4}+2\delta}\Delta_0^\frac{3}{2}.
\end{align*}
Combining the above two sets of estimates, also using (\ref{2.20.2.24}) and (\ref{1.27.5.24}), we conclude
\begin{align}\label{3.12.10.24}
X^{2-a}[\Lb, \Omega]X^a\Phi=O(\log \l t\r(\La_0+\log \l t\r^3\Delta_0^\frac{5}{4}))_{L^2_\Sigma}, a=0,1, X^2=\Omega S, S^2.
\end{align}

Next we consider the other term in  (\ref{3.12.8.24}). With $a=0,1$, applying (\ref{5.13.10.21}), (\ref{3.5.8.24}) and (\ref{3.25.1.22}) leads to  
\begin{align}\label{3.13.10.24}
&[S, \Omega^{2-a}]\Lb \Omega^a\Phi\nn\\
\displaybreak[0]
&=\tir(\pioh_{AL}(\sn\Omega\Lb\Phi+\sn\Lb\Omega\Phi)+\sn_\Omega(\pioh_{AL}\sn \Lb\Phi))\nn\\
&=O(1)\pioh_{AL}(\Omega^2\Lb\Phi+\Omega\Lb \Omega\Phi)+O( \l t\r^{-1}\log \l t\r^{\f12\M+1}(\La_0+\Delta_0^\frac{5}{4}\l t\r^\delta))_{L_u^2 L_\omega^2}\Omega^{\le 1}\Lb \Phi\nn\\
&=O((\log \l t\r)^{\f12\M+1}\l t\r^{-2+2\delta}(\La_0+\Delta_0^\frac{5}{4}))_{L^2_u L_\omega^2}, 
\end{align}
where we applied (\ref{8.23.1.23}) and (\ref{3.10.7.24}). This is a better error term compared with the estimate (\ref{3.12.10.24}) for the other type of commutators.

Therefore, combining (\ref{3.12.10.24}), (\ref{3.13.10.24}) with (\ref{3.12.3.24}), we conclude for the terms in (\ref{3.12.5.24}) excluding the first row, with $\vs^+(X^2)=1$
\begin{equation}\label{3.12.11.24}
\sum_{a=0}^2X^{2-a}\Omega \Lb X^a \Phi^\mu=O(1)\min(\mu, 1)\fB+O\Big(\log \l t\r^5(\La_0+\Delta_0^\frac{5}{4})+\l t\r^{-\f12+\delta}\Delta_0\Big)_{L^2_\Sigma}.
\end{equation}
Now we consider the first row in (\ref{3.12.5.24}).
Using Lemma \ref{6.30.4.23}, we derive  
\begin{align*}
[\tir \sn\Omega\Omega \bN v]&=\tir \sn\Omega\Omega \fB+O(1)\sn_\Omega^{\le 2}[\sn v]+O(\l t\r^{-\frac{3}{4}+\delta}\Delta_0)X^{1+\le 1}\fB\\
&+O(\l t\r^{-\frac{3}{4}+\delta}\Delta_0)_{L_u^2 L_\omega^4}(\sn_\Omega[\sn v]+\fB)+O(\l t\r^{-\f12+\delta}\Delta_0)_{L_u^2 L_\omega^2}[\sn v],\\
\tir \sn\Omega\Omega \bN v^\|&=\tir \sn\sn_\Omega^2[\sn v]+O(1)\sn_\Omega^{\le 2}\fB+O(\l t\r^{-\frac{3}{4}+\delta}\Delta_0)\Omega^{1+\le 1}[\sn v]\\
&+O(\l t\r^{-\frac{3}{4}+\delta}\Delta_0)_{L_u^2 L_\omega^4}(\sn_\Omega\fB+[\sn v])+O(\l t\r^{-\f12+\delta}\Delta_0)_{L_u^2 L_\omega^2}\fB.
\end{align*}
Using Proposition \ref{7.15.5.22}, we infer from the above estimates that 
\begin{equation*}
\tir \sn\Omega^2 \bN v=\tir \sn \Omega^2\fB+O(1)\Omega^{\le 2}\fB+O(\l t\r^{-\f12+\delta}\Delta_0)_{L^2_\Sigma}.
\end{equation*}
Applying (\ref{3.4.1.24}) and (\ref{LbBA2}) to treat the terms on the right-hand side leads to
\begin{align*} 
\sta{\Omega^3, \bN}\Phi&=\min(\mu, 1)O(\fB)+\tir \sn \Omega^2\fB+O(1)(\Omega^{1+\le 1}\fB+\mu \fB)+O(\l t\r^{-\f12+\delta}\Delta_0)_{L^2_\Sigma}+O(\l t\r^\delta\Delta_0)_{L^2_\Sigma}\\
&=\min(\mu, 1)O(\fB)+O(\l t\r^\delta\Delta_0)_{L^2_\Sigma}.
\end{align*}
Using $\|\sta{\Omega^3, L}\Phi\|_{L^2_\Sigma}\les \l t\r^{-\f12+\delta} \Delta_0$ due to (\ref{3.12.1.21}), we have
\begin{equation*}
\sta{\Omega^3, \Lb}\Phi=\min(1,\mu)O(\fB)+O(\l t\r^\delta\Delta_0)_{L^2_\Sigma}, 
\end{equation*}
as stated in (\ref{3.14.3.24}), which controls the term in the first row in (\ref{3.12.5.24}). Combing this estimate with (\ref{3.12.11.24}) leads to
\begin{align*}
X^{2-a}\Omega\Lb X^a\Phi=O(1)\min(1, \mu)\fB+O(\l t\r^\delta\Delta_0)_{L^2_\Sigma}
\end{align*}
Substituting the above estimate to (\ref{3.13.1.24}), in view of (\ref{7.28.2.24}) we then conclude (\ref{3.3.4.24}).

(\ref{3.5.4.24}) is derived from  Lemma \ref{8.1_com} (2).

Next we prove (\ref{3.6.4.24}). We first apply (\ref{8.22.2.21}) and (\ref{7.11.7.21}) to $\bT\varrho$, which yields, 
\begin{equation}\label{3.13.2.24}
\begin{split}
\Omega\Lb S\bT \varrho&=O(1)\tir \Omega\sD\bT\varrho+\tir \Omega\Box_\bg\bT\varrho+\sn_\Omega(\ze, \mho, k_{\bN\bN}) \sum_{X=S, \Omega}X\bT\varrho\\
&+O(1)\Omega \Lb \bT \varrho+\tir \Omega h\c \Lb \bT \varrho+\Omega^2\bT \varrho \c \ze,
\end{split}
\end{equation}
\begin{equation}\label{3.13.3.24}
\begin{split}
S\Lb \Omega\bT \varrho&=O(1)(\tir \sD \Omega\bT\varrho+\tir \Box_\bg\Omega \bT \varrho+\Lb \Omega \bT\varrho+L\Omega \bT\varrho)\\
&+O(\l t\r^{-1+\delta}\Delta_0)\sn\Omega\bT\varrho+(O(1)+O(\l t\r^{-1+2\delta}\Delta_0^\f12)_{L_\omega^4})\sn\bT\varrho.
\end{split}
\end{equation}
Using (\ref{LbBA2}), (\ref{2.13.3.24}) and (\ref{3.9.8.24}), we have 
\begin{align*}
&\tir \sD\Omega\bT\varrho, \tir\Omega\sD\bT\varrho=O(\l t\r^{-1+\delta})_{L^2_\Sigma}\\
\displaybreak[0]
&\Box_\bg\Omega \bT \varrho,\Box_\bg \Omega\bT\varrho=O(\l t\r^{-1}\log \l t\r^{\f12(\M+7)}(\La_0+\Delta_0^\frac{5}{4}))_{L^2_\Sigma} \\ 
 &(\Lb \Omega \bT\varrho, \Omega\Lb \bT\varrho)=O(\log \l t\r^{\f12(\M+7)}(\La_0+\Delta_0^\frac{5}{4}))_{L^2_\Sigma}.
\end{align*}
Using Proposition \ref{7.22.2.22}, Proposition \ref{7.15.5.22} and $X\bT\varrho=O(\l t\r^{-1+\delta}\Delta_0)+O(\fB)$, we have
\begin{align*}
\sn_\Omega(\ze, \mho, k_{\bN\bN})\sum_{X=S, \Omega}X\bT\varrho=O(\l t\r^{-2+2\delta}\Delta_0)_{L_u^2 L_\omega^2}\\
\Omega^2\bT\varrho\c \ze=O(\l t\r^{-1+2\delta}\Delta_0^2)_{L^2_\Sigma}.
\end{align*}
It follows by using (\ref{L4conn}) and (\ref{1.29.2.22}) that
\begin{align*}
\|\tir\Omega h\Lb \bT\varrho\|_{L^2_u L_\omega^2}\les \|\Omega h\|_{L_u^2 L_\omega^4}\les \l t\r^{-\frac{7}{4}+\delta}\Delta_0.
\end{align*}
Using the above estimates and applying (\ref{LbBA2}), we conclude
\begin{align}\label{3.13.7.24}
\Omega \Lb S\bT\varrho, S\Lb \Omega \bT\varrho=O(\log \l t\r^{\f12(\M+7)}\Delta_0)_{L^2_\Sigma}.
\end{align}
This gives two estimates stated  in the last estimate in (\ref{3.6.4.24}) with $X=S$. Applying (\ref{7.03.3.21}) to $f=\bT\varrho$,  we write
\begin{align*}
\sn_X^{1-a}\bd^2_{\Lb A} X^a \bT\varrho&=\sn_X^{1-a}(\sn_A\Lb X^a \bT\varrho+\chib\c \sn X^a\bT\varrho+\zb\Lb X^a \bT\varrho)\\
&=\sn_X^{1-a}(\sn_\Lb\sn_A X^a \bT\varrho+k_{A\bN} L X^a\bT\varrho+\ze \bN X^a\bT\varrho).
\end{align*}
If $a=1$, we take the first identity in the above, if $a=0$, we will employ the second identity.  Note that using (\ref{6.30.1.19}), (\ref{6.30.2.19}), (\ref{8.23.1.23}), (\ref{1.29.2.22}), (\ref{3.9.6.24}) and Proposition \ref{7.15.5.22} we have
\begin{align*}
\Lb S\bT\varrho, S\bT^2 \varrho=O(\l t\r^{-1})_{L_\omega^4}.
\end{align*}
Thus using the above estimates, (\ref{7.29.1.22}), (\ref{2.13.3.24}), Proposition \ref{6.27.1.24} and Proposition \ref{7.15.5.22}, we obtain the rough bound for errors 
\begin{align}
&\chib\c \sn X\bT\varrho+\zb\Lb X\bT\varrho, \sn_\Omega(\chib\c \sn \bT\varrho+\zb\Lb \bT\varrho)=O(\l t\r^{-3+2\delta}\Delta_0)_{L^2_u L_\omega^2},\label{8.12.1.24}\\
&\sn_S(k_{A\bN}L\bT\varrho+\ze\bN\bT\varrho)=O(\l t \r^{-1}\log \l t\r\Delta_0)_{L^2_\Sigma},\nn
\end{align}
where the last line follows by also using (\ref{8.23.1.23}) and Lemma \ref{5.13.11.21} (5).

Moreover due to (\ref{4.22.4.22}) and Lemma \ref{5.13.11.21} (5)
\begin{align*}
\tir \sn_S\sn_\Lb \sn \bT\varrho=S^{\le 1}\Lb \Omega \bT\varrho+O(\l t\r^{-3+2\delta}\Delta_0^\frac{3}{2}\log \l t\r)_{L_u^2 L_\omega^2}.
\end{align*}
Combining the above estimates with (\ref{3.13.7.24}) implies, with $a=0,1$, 
\begin{align*}
\tir \sn_S^{1-a}\bd^2_{\Lb A} S^a \bT\varrho&=\Omega \Lb S\bT\varrho, S^{\le 1}\Lb \Omega \bT\varrho+O(\log \l t\r\Delta_0)_{L^2_\Sigma}\nn\\
&=O(\log \l t\r^{\f12(\M+7)}\Delta_0)_{L^2_\Sigma}.
\end{align*}

Next we check  that 
\begin{align}\label{3.13.8.24}
\|\Omega\Lb\Omega \bT\varrho, \Omega^2 \Lb\bT\varrho\|_{L^2_\Sigma}\les \l t\r^\delta\Delta_0
\end{align}
which is the second set of estimates in (\ref{3.6.4.24}) with $a=0$ and $X=\Omega$. 

In view of (\ref{3.10.2.24}), we only need to prove
\begin{align}\label{3.13.5.24}
\|\Omega L\Omega\bT \varrho, \Omega^2 L\bT\varrho\|_{L^2_\Sigma}\les\l t\r^{-1}(\log \l t\r)^5\Delta_0.
\end{align}
Indeed, 
\begin{align*}
&\Omega L\Omega \bT\varrho=\Omega L\Omega L\varrho+\Omega L[\Omega,\Lb]\varrho+\Omega L\Lb \Omega \varrho\\
&\Omega^2 L\bT\varrho=\Omega^2 LL \varrho+\Omega^2 L\Lb \varrho.
\end{align*}
The last terms in the above two lines have been estimated in (\ref{3.12.3.24}), the commutator in the first line has been treated in (\ref{3.12.10.24}). Also using (\ref{3.12.1.21}), we conclude (\ref{3.13.5.24}).

Therefore, due to (\ref{3.13.8.24}) and (\ref{8.12.1.24}),
\begin{align*}
\tir \sn_\Omega^{1-a}\bd^2_{\Lb A} \Omega^a \bT\varrho&=(\Omega \Lb \Omega\bT\varrho, \Omega^2 \Lb\bT\varrho)+O(\l t\r^{-2+2\delta}\Delta_0)_{L^2_u L_\omega^2}\\
&=O(\l t\r^\delta\Delta_0)_{L^2_\Sigma}.
\end{align*}
Thus the first estimate in (\ref{3.6.4.24}) is proved.

It remains to consider the $L^2_\Sigma$ norm of the following terms to complete the second set of estimates in (\ref{3.6.4.24})
\begin{equation*}
\Lb \Omega S \bT \varrho, S\Omega \Lb \bT \varrho, \Lb \Omega^2 \bT \varrho.
\end{equation*}
To bound them by using the proved four estimates in the second set of estimates of (\ref{3.6.4.24}), it suffices to show
\begin{equation}\label{3.13.6.24}
X^{1-a}[\Lb, \Omega]X^a\bT\varrho=O(\l t\r^\delta\Delta_0)_{L^2_\Sigma}, a=0, 1.
\end{equation}
Indeed, we first note the following results, 
\begin{equation}\label{3.13.9.24}
\begin{split}
&\|\bN\Omega \bT\varrho, \Omega\bN \bT\varrho\|_{L_u^2 L_\omega^4}\les\l t\r^{\f12\delta-1}\log \l t\r^{\frac{\M}{4}+2}\Delta_0\\
&\bN S\bT\varrho, S\bN\bT\varrho=\Lb\bT\varrho+O(\l t\r^{-1})\fB+O(\l t\r^{-2+\delta}\Delta_0)_{L^2_u L_\omega^2}.
\end{split}
\end{equation}
The first line in the above is obtained by using (\ref{2.13.3.24}), (\ref{2.19.1.24}) and (\ref{3.10.2.24}). 

Due to (\ref{2.14.1.24}), (\ref{7.13.5.22}) and Lemma \ref{5.13.11.21} (5), $$\ud \bA\sn \bT \varrho=O(\l t\r^{-2+\delta}\Delta_0)_{L_u^2 L_\omega^2}.$$ 
 Hence the second line of (\ref{3.13.9.24}) is a result of the above estimate, (\ref{3.9.6.24}), (\ref{7.29.1.22}), (\ref{2.14.1.24}), (\ref{6.30.1.19}), (\ref{6.30.2.19}), (\ref{6.22.1.21}) and the following straightforward calculation,
\begin{align*}
\bN S\bT\varrho, S\bN\bT\varrho&=O(1)\tir(\tir^{-1}\Lb\bT\varrho+\sD\bT\varrho+\ud \bA \sn \bT\varrho+\Box_\bg \bT\varrho+LS \bT\varrho+S^{\le 1}L\bT\varrho).
\end{align*}

To see (\ref{3.13.6.24}), we derive, with $a=0,1$, that
\begin{align*}
X^{1-a}[\Lb, \Omega]X^a\bT\varrho&=X^{1-a}(\Omega\log \bb \bN X^a\bT \varrho+\pioh_{A\Lb}\sn X^a\bT\varrho).
\end{align*}
Hence
\begin{align*}
\sum_{a=0}^1 X^{1-a}[\Lb, \Omega]X^a\bT\varrho&=\Omega\log \bb (\bN X\bT\varrho+X\bN \bT\varrho)+X\Omega\log \bb \bN \bT\varrho\\
&+\sum_{a=0}^1 \sn_X^{1-a}(\pioh_{A\Lb} \sn X^a\bT \varrho).
\end{align*}
Next we treat the first two terms by using (\ref{3.13.9.24}), (\ref{1.27.5.24}), (\ref{2.20.2.24}), (\ref{1.29.2.22}) and the $\zeta$ estimates in Proposition \ref{7.15.5.22} that
\begin{align*}
\|\Omega\log \bb (\bN \Omega\bT\varrho+\Omega\bN \bT\varrho)\|_{L^2_\Sigma}&\les \l t\r^{\f12\delta}\log \l t\r^{\frac{\M}{4}+3}\Delta_0^2\\
\|\Omega\log\bb (\bN S\bT\varrho, S\bN\bT\varrho)\|_{L^2_\Sigma}&\les \|\Omega\log \bb\|_{L^2_u L_\omega^4}+\l t\r^{-1+2\delta}\Delta_0^2\log \l t\r^\f12\\
&\les(\log \l t\r)^3\Delta_0\\
\|X\Omega\log \bb \bN \bT\varrho\|_{L^2_\Sigma}&\les \|X\Omega\log \bb\|_{L_u^2 L_\omega^4}\les \l t\r^{\delta}\Delta_0.
\end{align*}
Using (\ref{5.21.1.21}), (\ref{3.25.1.22}), (\ref{LbBA2}) and (\ref{3.6.2.21}), it is straightforward to compute
\begin{align*}
\sum_{a=0}^1 \sn_X^{1-a}(\pioh_{A\Lb} \sn X^a\bT \varrho)&=\pioh_{\Lb A}(\sn X \fB+\sn_X \sn \fB)+\sn_X\pioh_{A\Lb}\sn\fB\\
&=O(\l t\r^{3\delta-1}\Delta_0\log \l t\r^\f12)_{L^2_\Sigma}.
\end{align*}
Summarizing the above estimates gives (\ref{3.13.6.24}). Using (\ref{3.13.7.24}), (\ref{3.13.8.24}) and (\ref{3.13.6.24}), we can obtain
\begin{equation*}
\sum_{a=0}^1\|X^{1-a}(\Lb\Omega, \Omega\Lb) X^a\bT\varrho\|_{L^2_\Sigma}\les \l t\r^\delta\Delta_0
\end{equation*}
as stated in (\ref{3.6.4.24}).

Next we prove the second estimate in (\ref{3.6.3.24}). Noting that due to (\ref{3.12.1.21}), $\Omega^{1-a}X^2\Omega^aL\varrho=O(\l t\r^{-\frac{1}{2}+\delta}\Delta_0)_{L^2_\Sigma}$, it then suffices to consider
\begin{align*}
\Omega^{1-a}X^2\Omega^a \Lb \varrho=\Omega^3\Lb\varrho,\sum_{\vs^+(X^2)=1} \Omega^{1-a}X^2 \Omega^a \Lb\varrho,
\end{align*}
where the first term on the right hand side is $O(\l t \r^\delta\Delta_0)_{L^2_\Sigma}$ due to (\ref{3.14.3.24}). It only remains to consider the term $\sum_{\vs^+(X^2)=1} \Omega^{1-a}X^2 \Omega^a \Lb\varrho.$

Due to (\ref{3.12.3.24})
\begin{align*}
\Omega^2 S\Lb \varrho, \Omega SS\Lb \varrho=O((\log \l t\r)^5(\La_0+\Delta_0^\frac{5}{4}))_{L^2_\Sigma}.
\end{align*}
To exhaust all the cases in $\sum_{\vs^+(X^2)=1} \Omega^{1-a}X^2 \Omega^a \Lb\varrho$, apart from those bounded in the above, it suffices to consider the following three terms, which can be reduced to the known estimates modulo commutators. 
\begin{align*}
&\Omega S\Omega\Lb\varrho=\Omega^2 S\Lb \varrho+\Omega[S, \Omega]\Lb \varrho, S\Omega^2 \Lb \varrho=\Omega^2 S\Lb \varrho+[S, \Omega^2]\Lb \varrho \\
 &S^2\Omega \Lb \varrho=[S^2, \Omega]\Lb \varrho+\Omega S^2\Lb \varrho
\end{align*}
The commutator $\Omega[S,\Omega]\Lb \varrho$ is one of terms treated while controlling $[S, \Omega^2]\Lb\varrho$, with the latter given in (\ref{3.13.10.24}). Hence the only commutator we need to consider is the following one,  written in view of (\ref{5.13.10.21}), 
\begin{align}\label{3.16.10.24}
[S^2, \Omega]\Lb \varrho&=\sum_{a=0}^1\sn_S^{1-a}(\tir\pioh_{LA}\sn S^a \Lb\varrho)\nn\\
&=O(\l t\r^{-2+2\delta}(\La_0+\Delta_0^\frac{5}{4}))_{L^2_u L_\omega^2}
\end{align}
which is treated similar to (\ref{3.13.10.24}) with the help of (\ref{3.10.7.24}).

Note that the above commutators all verify better estimates than the leading terms. We thus conclude 
\begin{align}\label{8.13.1.24}
\sum_{\vs^+(X^2)=1}\Omega^{1-a}X^2\Omega^a\Lb\varrho=O(\log \l t\r^5(\La_0+\Delta_0^\frac{5}{4}))_{L^2_\Sigma}.
\end{align}
Therefore 
\begin{equation*}
\Omega^{1-a}X^2\Omega^a \Lb \varrho=O(\l t\r^\delta\Delta_0)_{L^2_\Sigma},\, a=0, 1.
\end{equation*}
Hence we completed the second estimate in (\ref{3.6.3.24}). To obtain the first estimate in (\ref{3.6.3.24}), using (\ref{7.03.3.21}), (\ref{7.29.1.22}), (\ref{6.22.1.21}), (\ref{2.19.1.24}), (\ref{8.29.9.21}) and Proposition \ref{1steng}, we control the error terms
\begin{align*}
\Omega^{1-a}(k_{\bN A}L \Omega^a \bT\varrho), \Omega^{1-a}(k_{\bN\bN}L\Omega^a \bT\varrho)=O(\l t\r^{-2}(\Delta_0^\frac{5}{4}+\La_0)\log \l t\r)_{L^2_\Sigma}
\end{align*}
which are negligible  terms.  
Finally, using (\ref{4.22.4.22}), we can obtain the first estimates in (\ref{3.6.3.24}) from the second one.  

\end{proof}

\subsection{Top order geometric estimates}
Under the bootstrap assumptions (\ref{3.12.1.21})-(\ref{6.5.1.21}), we first derive a useful preliminary result. 
\begin{lemma}\label{2.26.3.24}
There hold for $n=0,1$ that
\begin{align*}
 &\|(\tir \sn)^\ell(\bb^{-n}\sn\chih)\|_{L^2_u L_\omega^2}\les\log \l t\r\|(\tir \sn)^\ell\sF\|_{L_u^2 L_\omega^2}+\l t\r^{-2+\delta}\Delta_0, \ell=1,2,\\
 &\|(\tir \sn)^\ell(\bb^{-n} \sn\chih)\|_{L_\omega^4}\les \log \l t\r\|(\tir \sn)^\ell\sF\|_{L_\omega^4}+\Delta_0\l t\r^{-2+\delta}, \ell=0,1.
\end{align*}
\end{lemma}
\begin{proof}
Recalling from (\ref{3.31.5.22}) the definition of $\sF$,  using (\ref{1.27.5.24}), (\ref{L4BA1}) and Sobolev embedding on spheres, we directly obtain
\begin{align}\label{5.31.1.24}
\|(\tir\sn)^\ell(\sn \tr\chi)\|_{L_\omega^4}&\les \|(\tir\sn)^\ell(\bb \sF-\sn\Xi_4)\|_{L_\omega^4}\\
&\les \log \l t\r\|(\tir \sn)^{\le 1} \sF\|_{L_\omega^4}+\l t\r^{-2+\delta}\Delta_0, \ell=0,1.\nn
\end{align}
Similarly, by also using (\ref{LbBA2}), we have for $\ell=1,2$ that
\begin{align*}
\|(\tir \sn)^\ell(\sn \tr\chi)\|_{L_u^2 L_\omega^2}&\les \|(\tir \sn)^\ell (\bb \sF-\sn \Xi_4)\|_{L_u^2 L_\omega^2}\\
&\les\|(\tir \sn)^\ell(\bb \sF)\|_{L_u^2 L_\omega^2}+\|(\tir \sn)^\ell\sn \Xi_4\|_{L_u^2 L_\omega^2}\\
&\les\|\bb (\tir \sn)^\ell \sF\|_{L_u^2 L_\omega^2}+\max(\ell-1, 0)\|(\tir \sn)^{\ell-1}\bb (\tir \sn)\sF\|_{L_u^2 L_\omega^2}\\
&+\|\sF\c (\tir \sn)^\ell \bb \|_{L_u^2 L_\omega^2}+\l t\r^{-2+\delta}\Delta_0 
\end{align*}
It follows by using Proposition \ref{7.15.5.22}, Lemma \ref{5.13.11.21} (5) and (\ref{9.12.2.22}) that
\begin{align*}
\|\sF \c (\tir \sn)^\ell \bb \|_{L_u^2 L_\omega^2}&\les \l t\r^{-2+\delta}\Delta_0^2,\, \ell=1,2
\end{align*}
and with $\ell=2$, 
\begin{align*}
\|(\tir \sn)^{\ell-1}\bb (\tir \sn)\sF\|_{L_u^2 L_\omega^2}&\les \|(\tir \sn)\sF\|_{L_u^2 L_\omega^4}\|(\tir \sn)\bb\|_{L_u^\infty L_\omega^4}\\
&\les \log \l t\r\Delta_0\|(\tir \sn)^{1+\le 1}\sF\|_{L_u^2 L_\omega^2}\\
&\les\log \l t\r\Delta_0|(\tir \sn)^2\sF\|_{L_u^2 L_\omega^2}+\l t\r^{-2+\delta}\Delta_0^2.
\end{align*}
Hence, 
\begin{align}\label{5.31.2.24}
\|(\tir \sn)^\ell(\sn \tr\chi)\|_{L_u^2 L_\omega^2}\les \log \l t\r\|(\tir \sn)^\ell \sF\|_{L_u^2 L_\omega^2}+\l t\r^{-2+\delta}\Delta_0.
\end{align}
  Applying Proposition \ref{11.4.1.22} (1), Proposition \ref{7.15.5.22},  Lemma \ref{5.13.11.21} (3), and using (\ref{5.31.1.24}) and (\ref{5.31.2.24}), we can conclude Lemma \ref{2.26.3.24}. 
\end{proof}

Next we prove the following key estimates for derivatives of $\sF$. 
\begin{proposition}
Under the bootstrap assumptions (\ref{3.12.1.21})-(\ref{6.5.1.21}), we have
\begin{align}\label{9.6.3.22}
\begin{split}
\|\tir^2 \sF\|_{L_\omega^p}&\les \log \l t\r^{\f12\M+1}(\Delta_0^\frac{5}{4}
+\La_0), \, 2\le p\le 4\\
\|\tir^3(\tir\sn)^\ell\sF\|_{L_\omega^4} &\les\|\tir^3 (\tir \sn)^\ell \sF(0)\|_{L_\omega^4}+\int_0^t \tir^2 \{\|(\tir \sn)^\ell (\bb^{-1}\sn \fB)\|_{L_\omega^4}\\
&+\ell \l t\r^{-1}\Delta_0^\f12(\|\bb^{-2}(\tir \sn)^\ell \Xi_4\|_{L_\omega^\infty}+\|\tir^2\sn\sF\|_{L_\omega^4})\}\\
&+\l t\r^{\frac{1}{4}+2\delta}\Delta_0^\frac{3}{2};\ell=1
\end{split}
\end{align}

\begin{align}\label{10.26.3.23}
\begin{split}
\|\tir^3 (\tir\sn)^\ell \sF\|_{L_u^2 L_\omega^2}&\les\int_0^t \l t'\r\log \l t'\r\Delta_0\|(\tir \sn)^{\ell+1}(\bb^{-1}\Xi_4)\|^\f12_{L_u^2 L_\omega^2}\|(\tir \sn)^\ell (\bb^{-1}\Xi_4)\|^\f12_{L_u^2 L_\omega^2}\\
&+\|\tir^3 \sn^\ell \sF(0)\|_{L_u^2 L_\omega^2}+\l t\r^{\f12+2\delta}\Delta_0^\frac{3}{2}+\int_0^t \|\tir^2\bb^{-1}(\tir \sn)^\ell \sn \Xi_4\|_{L_u^2 L_\omega^2}\\
&+\int_0^t \l t'\r^2 \log \l t'\r\|(\tir\sn)^{\le \ell}\sn[L\Phi]\|_{L_u^2 L_\omega^2}\\
&+\int_0^t \Delta_0\l t'\r^{-1}(\log \l t'\r)^4\|\tir^3(\tir \sn)^{\ell-1}\sF\|_{L_u^2 L_\omega^2},\quad\ell=1,2,
\end{split}
\end{align}
and
\begin{equation}\label{1.29.1.24}
\|\tir^2\sF\|_{L_u^2 L_\omega^2}\les \log \l t\r(\La_0+\Delta_0^\frac{5}{4}).
\end{equation}

And there also holds the rough estimate, 
\begin{align}\label{3.15.6.24}
\|\tir^2 (\tir\sn)^2 \sF\|_{L^2_u L_\omega^2}+\log \l t\r^{-1}\|\tir^2(\tir \sn)^{\le 2}(\bb^{-1}\sn\chih)\|_{L^2_u L_\omega^2}\les \l t\r^\delta\Delta_0.
\end{align}

\end{proposition}

\begin{proof}
For the term $(\tir\sn)^{\ell}\sF$, we will provide the comparison estimate by using the transport equations  (\ref{12.5.1.21}) and  (\ref{12.5.2.21}).
Using the transport lemma and (\ref{3.6.2.21}), with $\ell=0,1,2$, we bound
\begin{align}
\|\tir^3(\tir\sn)^\ell \sF\|_{L_\omega^p}&\les \|\tir^3(\tir\sn)^\ell \sF\|_{L_\omega^p(S_{0,u})}+\int_0^t \tir^3\|(\tir\sn)^\ell G_2\|_{L_\omega^p}\nn\\
&+\ell\sum_{a=0}^{\ell-1}\|\tir^3(\tir \sn)^{\ell-1-a}(\tir(\bR_{BCLA}, \chi\zb)(\tir \sn)^a \sF)\|_{L_\omega^p}\nn\\
&+\max(\ell-1, 0)\|\tir^3\tir\sn\tr\chi (\tir\sn) \sF\|_{L_\omega^p},\label{12.13.1.21}
\end{align}
where 
\begin{align*}
(\tir\sn)^{\ell} G_2&=(\tir\sn)^{\ell} {\ti G}_1+(\tir\sn)^\ell\{\bb^{-1}\sn[L\Phi]\sX+\bb^{-1}(-2\sn \chih\c \chih\\
&+\sn \N(\Phi, \p\Phi))-\chih\c \sF\}, 
\end{align*}
with ${\ti G}_1=\bb^{-1}\sn\Xi_4(\Xi_4+\f12\sX)$, and we used $k_{\bN\bN}-\f12\Xi_4=[L\Phi]$. 

Using $\tir |\Xi_4, \sX|\les 1$, it is straightforward to bound
\begin{align}\label{2.1.3.24}
\begin{split}
\|(\tir\sn)^\ell G_2\|_{L_\omega^p}&\les\|(\tir\sn)^\ell(\sn[L\Phi]\bb^{-1}\sX)\|_{L_\omega^p}\\
&+\|(\tir\sn)^\ell(\bb^{-1}\sn \chih\c \chih)\|_{L_\omega^p}+\|(\tir\sn)^\ell(\chih\c\sF)\|_{L_\omega^p}\\
&+\|(\tir\sn)^\ell\big(\bb^{-1}\sn\N(\Phi, \p \Phi)\big)\|_{L_\omega^p}+\|(\tir\sn)^\ell{\ti G}_1\|_{L_\omega^p}.
\end{split}
\end{align}

Note symbolically
$
\Xi_4, \fB=[L\Phi]+\Lb \varrho.
$
We derive by using (\ref{3.6.2.21}) and (\ref{3.11.3.21}) that
\begin{equation}\label{10.26.4.23}
\|\sn(\Xi_4), \sn \fB, \sn \sX\|_{L^4_\omega}\les \l t\r^{-2+\delta}\Delta_0, \quad\|\tir\sn \chih, \tir\sn \tr\chi\|_{L_\omega^4}\les\Delta_0 \l t\r^{-\frac{7}{4}+\delta}.
\end{equation}
Using Proposition \ref{7.15.5.22}, we have 
\begin{equation}\label{1.28.1.24}
|\sF, \sn \chih|\les \l t\r^{-2+\delta}\Delta_0, \|(\tir \sn)^m \sF\|_{L_u^2 L_\omega^2}\les \l t\r^{-2+\delta}\Delta_0, m\le 1.
\end{equation}
Due to (\ref{10.26.4.23}), (\ref{1.28.1.24}) and Sobolev embedding on spheres, 
\begin{align*}
\|\tir^2 \sn \tr\chi\c \sn \sF\|_{L_u^2 L_\omega^2}&\les \l t\r^{-\frac{7}{4}+\delta}\Delta_0\|(\tir\sn)^2\sF\|_{L_u^2 L_\omega^2}^\f12 \|\tir \sn \sF\|^\f12_{L_u^2 L_\omega^2}\\ 
&\les \|(\tir \sn)^2 \sF\|_{L_u^2 L_\omega^2}\c\l t\r^{-\frac{7}{4}+\delta}\Delta_0+ \l t\r^{-3-\frac{3}{4}+2\delta}\Delta_0^2.
\end{align*}
Using (\ref{3.6.2.21}), (\ref{L2BA2}), (\ref{1.28.1.24}), (\ref{1.29.4.22}) and (\ref{8.24.4.23}), we derive
\begin{align*}
\sum_{a=0}^{\ell-1}&\|(\tir \sn)^{\ell-1-a}(\tir(\bR_{BCLA}, \chi\c \zb)(\tir \sn)^a\sF)\|_{L^2_u L_\omega^2}\\
 &\les \|(\tir \sn)^{\ell-1}(\tir(\bR_{BCLA}, \chi\c \zb))\|_{L_u^2 L_\omega^2}\|\sF\|_{L_\omega^\infty}+\|\tir(\bR_{BCLA}, \chi\c \zb)\|_{L_\omega^\infty}\|(\tir \sn)^{\ell-1}\sF\|_{L_u^2 L_\omega^2}\\
 &\les \l t\r^{-4+2\delta}\Delta_0^2+\l t\r^{-\frac{7}{4}+\delta}\Delta_0\|(\tir \sn)^{\ell-1}\sF\|_{L_u^2 L_\omega^2}\les \l t\r^{-3-\frac{3}{4}+2\delta}\Delta_0^2. 
\end{align*}
The above two estimates control the errors except the term of $G_2$ in (\ref{12.13.1.21}). Next we consider the bound on $\|(\tir\sn)^\ell G_2\|_{L_u^2 L_\omega^2}$ in view of (\ref{2.1.3.24}). 
Using Proposition \ref{7.15.5.22}, Sobolev embedding on spheres and Cauchy-Schwarz inequality,  we deduce with $\ell=0,1,2$
\begin{align*}
&\|(\tir \sn)^\ell(\bb^{-1}\sn\chih\c \chih)\|_{L_u^2 L_\omega^2}\\
&\les \|(\tir \sn)^\ell(\bb^{-1}\sn\chih)\|_{L_u^2 L_\omega^2}\l t\r^{-\frac{7}{4}+\delta}\Delta_0+\min(\ell-1, 0)\|(\tir \sn)^{\ell-1}(\bb^{-1}\sn\chih)\|_{L_u^2 L_\omega^4}\|\tir\sn \chih\|_{L_\omega^4}\\
&+\|\bb^{-1}\sn\chih \c (\tir \sn)^\ell \chih\|_{L_u^2 L_\omega^2}\\
&\les\l t\r^{-\frac{7}{4}+\delta}\Delta_0\sum_{n=0}^1\|(\tir \sn)^\ell(\bb^{-n}\sn\chih)\|_{L^2_u L_\omega^2}+\Delta_0^2 \l t\r^{-\frac{5}{4}-3+2\delta}.
\end{align*}
It follows by using Lemma \ref{2.26.3.24} and  the above estimate that
\begin{align*}
\|(\tir \sn)^\ell&(\bb^{-1}\sn\chih\c \chih)\|_{L_u^2 L_\omega^2}\les \|(\tir \sn)^\ell\sF\|_{L^2_u L_\omega^2}\l t\r^{-\frac{7}{4}+\delta+}\Delta_0+\Delta_0^2 \l t\r^{-\frac{3}{4}-3+2\delta}.
\end{align*}
Moreover, by using (\ref{3.6.2.21}), (\ref{10.26.4.23}) and (\ref{L2conndrv}) for $\chih$ and its derivatives, (\ref{1.28.1.24}) and Sobolev embedding, we derive that 
\begin{align*}
\|(\tir \sn)^\ell (\chih\c \sF)\|_{L_u^2 L^2_\omega}
&\les \|(\tir \sn)^\ell \sF\c\chih\|_{L_u^2 L_\omega^2}+\sum_{a=1}^\ell\|(\tir \sn)^{\ell-a}\sF(\tir \sn)^a \chih\|_{L_u^2 L_\omega^2}\\
&\les \|(\tir \sn)^\ell \sF\|_{L_u^2 L_\omega^2}\c\l t\r^{-\frac{7}{4}+\delta}\Delta_0+ \l t\r^{-\frac{7}{2}+2\delta}\Delta_0^2\\
&+\l t\r^{-\frac{7}{4}+\delta}\Delta_0 \|(\tir \sn)^{\ell -1}\sF\|_{L_u^2 L_\omega^4}\\
&\les  \|(\tir \sn)^\ell \sF\|_{L_u^2 L_\omega^2}\c\l t\r^{-\frac{7}{4}+\delta}\Delta_0+ \l t\r^{-\frac{7}{2}+2\delta}\Delta_0^2. 
\end{align*}

 Using (\ref{10.26.4.23}), (\ref{L2BA2}) and (\ref{3.6.2.21}) for $\sn\log \bb$, we bound
\begin{align*}
&\|(\tir \sn)^\ell(\sn[L\Phi]\bb^{-1}\sX)\|_{L_u^2 L_\omega^2}\\
&\les \l t\r^{-1}\|\bb^{-1}(\tir\sn)^{\ell} \sn[L\Phi]\|_{L_u^2 L_\omega^2}+\|(\tir \sn)^{\ell-1}\sn[L\Phi]\|_{L_u^2 L_\omega^4}\|(\tir \sn)(\bb^{-1}\sX)\|_{L_\omega^4}\\
&+\|\sn[L\Phi](\tir \sn)^\ell(\bb^{-1}\sX)\|_{L_u^2 L_\omega^2}\\
&\les \l t\r^{-1}\|\bb^{-1}(\tir\sn)^{\ell} \sn[L\Phi]\|_{L_u^2 L_\omega^2}+ \l t\r^{-4+\frac{1}{4}+2\delta}\Delta_0^2+\|\sn[L\Phi](\tir \sn)^\ell(\bb^{-1}\sX)\|_{L_u^2 L_\omega^2}. 
\end{align*}
For the last term, it is direct to obtain
\begin{align}\label{6.1.1.24}
(\tir \sn)^\ell(\bb^{-1}\sX)=(\tir \sn)^{\ell-1}\Big((\tir \sn)(\bb^{-1})\sX\Big)+\tir(\tir \sn)^{\ell-1}\sF.
\end{align}
In view of the above identity, noting that by using (\ref{zeh}), (\ref{3.6.2.21}), (\ref{10.26.4.23}) and (\ref{1.28.1.24}), we have
\begin{align*}
&\|\sn L[\Phi](\tir \sn)^{\ell}(\bb^{-1}\sX)\|_{L^2_u L_\omega^2}\\
&\les\|\sn[L\Phi]\|_{L_\omega^4}\|\tir (\tir \sn)^{\ell-1}\sF\|_{L_u^2 L_\omega^4}+\|\sn[L\Phi]\|_{L^\infty_\omega}\|(\tir \sn)^{\ell -1}(\tir \sn \log \bb \bb^{-1}\sX)\|_{L_u^2 L_\omega^2}\\
&\les\l t\r^{-2+\delta}\Delta_0\|(\tir \sn)^\ell \sF\|^\f12_{L_u^2 L_\omega^2}\|(\tir \sn)^{\ell -1}\sF\|^\f12_{L_u^2 L_\omega^2}+\Delta_0 \l t\r^{-\frac{7}{4}+\delta}\|(\tir \sn)^{\ell-1}( \sn\log \bb \bb^{-1}\sX)\|_{L_u^2 L_\omega^2}\\
&\les \l t\r^{-2+\delta}\Delta_0\|(\tir \sn)^\ell \sF\|_{L_u^2 L_\omega^2}+\l t\r^{-\frac{15}{4}+3\delta}\Delta_0^2. 
\end{align*}
Combining both yields
\begin{align*}
&\|(\tir \sn)^\ell(\sn[L\Phi]\bb^{-1}\sX)\|_{L_u^2 L_\omega^2}\\
&\les  \l t\r^{-1}\|\bb^{-1}(\tir\sn)^{\ell} \sn[L\Phi]\|_{L_u^2 L_\omega^2}+ \l t\r^{-2+\delta}\Delta_0\|(\tir \sn)^\ell \sF\|_{L_u^2 L_\omega^2}+\l t\r^{-\frac{15}{4}+3\delta}\Delta_0^2 .
\end{align*}
This proves
\begin{align}
\|(\tir\sn)^\ell G_2\|_{L_u^2 L_\omega^2}&\les \l t\r^{-\frac{7}{4}+\delta+}\Delta_0\|(\tir \sn)^\ell \sF\|_{L_u^2 L_\omega^2}+ \|\bb^{-1}(\tir\sn)^\ell(\sn[L\Phi])\tir^{-1}\|_{L^2_u L_\omega^2}\nn\\
&+\|(\tir\sn)^ \ell\big(\bb^{-1}\sn\N(\Phi, \p \Phi)\big)\|_{L_u^2 L_\omega^2}+\|(\tir\sn)^\ell{\ti G}_1\|_{L_u^2 L_\omega^2}+ \Delta_0^2\l t\r^{-\frac{7}{2}+2\delta}.\label{9.6.1.22}
\end{align}
At last we consider the crucial term $(\tir \sn)^\ell {\ti G}_1$. 
\begin{align*}
\|(\tir \sn)^\ell&(\sn \Xi_4\c \bb^{-1}(\Xi_4+\f12\sX))\|_{L_u^2 L_\omega^2}\\
&\les \|\tir^{-1}\bb^{-1}(\tir \sn)^\ell \sn \Xi_4\|_{L_u^2 L_\omega^2}+\|(\tir \sn)^{\ell-1}\sn \Xi_4\c \tir \sn\Big(\bb^{-1}(\frac{3}{2}\Xi_4+\f12\tr\chi)\Big)\|_{L_u^2 L_\omega^2}\\
&+\|\bb\sn\Xi_4\|_{L_u^\infty L_\omega^4}\|\bb^{-1}(\tir \sn)^\ell(\bb^{-1}(\Xi_4+\f12 \sX))\|_{L_u^2 L_\omega^4}. 
\end{align*}
Using (\ref{zeh}), Lemma \ref{5.13.11.21} (5), (\ref{10.26.4.23}), (\ref{6.1.1.24}) and Sobolev embedding on spheres, schematically we obtain for $\ell\ge 1$ that
\begin{align*}
\|(\tir \sn)^\ell&(\sn \Xi_4\c \bb^{-1}(\Xi_4+\f12\sX))\|_{L_u^2 L_\omega^2}\\
&\les \|\tir^{-1}\bb^{-1}(\tir \sn)^\ell \sn \Xi_4\|_{L_u^2 L_\omega^2}+\l t\r^{-2}\log \l t\r\Delta_0(\|\tir(\tir \sn)^{\ell-1}\sF\|_{L_u^2 L_\omega^4}+\l t\r^{-1+2\delta}\Delta_0^2\\
&+\|(\tir \sn)^\ell(\bb^{-1}\Xi_4)\|_{L_u^2 L_\omega^4}) \\
&\les \l t\r^{-2}\log \l t\r^{-\frac{5}{4}}\Delta_0 \|\tir (\tir\sn)^\ell \sF\|_{L_u^2 L_\omega^2}+\Delta_0\l t\r^{-2}\log \l t\r^\frac{13}{4}\|\tir(\tir \sn)^{\ell-1}\sF\|_{L_u^2 L_\omega^2}\\
&+\l t\r^{-2}\log \l t\r\Delta_0\|(\tir \sn)^{\ell+1}(\bb^{-1}\Xi_4)\|^\f12_{L_u^2 L_\omega^2}\|(\tir \sn)^\ell (\bb^{-1}\Xi_4)\|^\f12_{L_u^2 L_\omega^2}\\
&+\|\tir^{-1}\bb^{-1}(\tir \sn)^\ell \sn \Xi_4\|_{L_u^2 L_\omega^2}+\l t\r^{-3+2\delta}\log \l t\r\Delta_0^3.
\end{align*}
In view of transport lemma, we thus summarize the above estimates as
\begin{align*}
\|\tir^3 (\tir\sn)^\ell \sF\|_{L_u^2 L_\omega^2}&\les\int_0^t \l t'\r\log \l t'\r\Delta_0\|(\tir \sn)^{\ell+1}(\bb^{-1}\Xi_4)\|^\f12_{L_u^2 L_\omega^2}\|(\tir \sn)^\ell (\bb^{-1}\Xi_4)\|^\f12_{L_u^2 L_\omega^2}\\
&+\int_0^t \{\|\tir^2\bb^{-1}(\tir \sn)^\ell \sn (\Xi_4, [L\Phi])\|_{L_u^2 L_\omega^2}+ \|\tir^3(\tir\sn)^\ell\big(\bb^{-1}\sn\N(\Phi, \p \Phi)\big)\|_{L_u^2 L_\omega^2}\}\\
&+\int_0^t \Delta_0\l t'\r^{-1}\log \l t'\r^\frac{13}{4}\|\tir^3(\tir \sn)^{\ell-1}\sF\|_{L_u^2 L_\omega^2}+ \|\tir^3 \sn^\ell \sF(0)\|_{L_u^2 L_\omega^2}+\l t\r^{\f12+2\delta}\Delta_0^2.
\end{align*}
Applying Proposition \ref{6.24.10.23}, (\ref{3.6.2.21}) and Lemma \ref{5.13.11.21} (2) and (5), we derive
\begin{align*}
&\|(\tir\sn)^\ell\big(\bb^{-1}\sn\N(\Phi, \p \Phi)\big)\|_{L_u^2 L_\omega^2}\\
 &\les\|(\tir \sn)^\ell(\bb^{-1})\|_{L_u^2 L_\omega^2}\|\sn\N(\Phi, \bp\Phi)\|_{L^\infty_\omega}+\|(\tir \sn)^{\ell-1+\le 1}\sn \N(\Phi, \bp\Phi)\|_{L^2_u L_\omega^2}\log \l t\r\Delta_0\\
 &+\|\bb^{-1}(\tir \sn)^\ell \sn\N(\Phi, \bp\Phi)\|_{L_u^2 L_\omega^2}\\
 &\les \l t\r^{-1-\frac{11}{4}+2\delta}\Delta_0^\frac{3}{2}+\log \l t\r\l t\r^{-1}\Delta_0 \|(\tir\sn)^{\ell-1+\le 1}\sn[L\Phi]\|_{L_u^2L_\omega^2}\\
 &+\l t\r^{-1}\|\bb^{-1}(\tir \sn)^{\ell}\sn[L\Phi]\|_{L_u^2 L_\omega^2}. 
\end{align*}
Integrating in $t$ gives
\begin{align*}
 \int_0^t \|\tir^3&(\tir\sn)^\ell\big(\bb^{-1}\sn\N(\Phi, \p \Phi)\big)\|_{L_u^2 L_\omega^2}\\
 &\les \int_0^t \l t'\r^2 \log \l t'\r \|(\tir \sn)^{\ell-1+\le 1}\sn[L\Phi]\|_{L_u^2 L_\omega^2} +\Delta_0^\frac{3}{2}\l t\r^{\frac{1}{4}+2\delta}.
 \end{align*}
 Hence we conclude (\ref{10.26.3.23}) for $\ell\ge 1$. 
 
With $\ell=0$ in (\ref{12.13.1.21}), we have the straightforward estimate by using the transport lemma, Proposition \ref{12.21.1.21} , Proposition \ref{7.15.5.22}, Proposition \ref{1steng} and Proposition \ref{6.24.10.23} that
\begin{align*}
\|\tir^3 \sF\|_{L_u^2 L_\omega^2}&\les \|\tir^3\sF(0)\|_{L_u^2 L_\omega^2}+\int_0^t \l t\r\|\bb^{-1}\tir\sn \fB\|_{L_u^2 L_\omega^2}\\
&+\int_0^t\tir^3\|\bb^{-1}\sn\N(\Phi, \p \Phi), \bb^{-1}\sn\chih\c \chih\|_{L_u^2 L_\omega^2}\\
&\les  \l t\r\log \l t\r(\La_0+\Delta_0^\frac{5}{4}).
\end{align*}
This gives (\ref{1.29.1.24}).

Next we prove (\ref{9.6.3.22}) by using (\ref{12.13.1.21}). For $\ell=0$ and $2\le p\le 4$ using the transport lemma,
\begin{align*}
\|\tir^3 \sF\|_{L_\omega^p}&\les \|\tir^3 \sF(0)\|_{L_\omega^p}+\int_0^t \tir^3\|G_2\|_{L_\omega^p}\\
&\les \|\tir^3 \sF(0)\|_{L_\omega^p}+\int_0^t \tir^3\|{\ti G}_1\|_{L_\omega^p}\\
&+\int_0^t \{\tir^2\|\bb^{-1}\sn[L\Phi]\|_{L_\omega^p}+\tir^3 \|\bb^{-1}\sn\chih\c \chih\|_{L_\omega^p}+\tir^3\|\sn\N(\Phi,\bp\Phi)\|_{L_\omega^p}\}.
\end{align*}
Using (\ref{3.6.2.21}), (\ref{L4conn'}), (\ref{10.26.4.23}), (\ref{9.14.3.22}) and (\ref{10.22.5.22}), we derive
\begin{align}
\|\tir^3 \sF\|_{L_\omega^p}
&\les \|\tir^3 \sF(0)\|_{L_\omega^p}+\int_0^t \tir^2\|\sn\fB\|_{L_\omega^p}+\l t\r^{2\delta}\Delta_0^\frac{3}{2}\nn\\
&\les (\log \l t\r)^{\f12\M+1}\l t\r(\Delta_0^\frac{5}{4}
+\La_0)\label{10.26.2.23}
\end{align}
This gives the $\ell=0$ case in (\ref{9.6.3.22}). 

Next we consider $\ell=1$. Using (\ref{12.13.1.21}), the main term on the right-hand side is given in (\ref{2.1.3.24}). With $\Delta_0^\f12 \M_0\les 1$, using (\ref{6.5.1.21}) and Lemma \ref{5.13.11.21} (5) we derive
\begin{align*}
&\|(\tir \sn)^\ell{\ti G}_1\|_{L_\omega^4}\\
&\les\|(\tir \sn)^\ell (\bb^{-1}\sn \Xi_4)\|_{L_\omega^4}\l t\r^{-1}+\ell \|\bb^{-1}\sn \Xi_4 (\tir \sn)^\ell(\Xi_4+\f12\sX)\|_{L_\omega^4}\\
&\les \|(\tir \sn)^\ell (\bb^{-1}\sn\Xi_4)\|_{L_\omega^4}\l t\r^{-1}+\ell\|(\bb^{-2}\sn(\bb \Xi_4)-\bb^{-2}\sn \bb \Xi_4) (\tir \sn)^\ell(\Xi_4+\f12\sX)\|_{L_\omega^4}\\
&\les \|(\tir \sn)^\ell (\bb^{-1}\sn \Xi_4)\|_{L_\omega^4}\l t\r^{-1}+\ell(\l t\r^{-2}\Delta_0 (\|\bb^{-2}(\tir \sn)^\ell \Xi_4\|_{L_\omega^\infty}+\|\bb^{-1}\tir \sF\|_{L_\omega^\infty})\\
&+\l t\r^{-2}\Delta_0^\f12|\bb^{-2}(\tir \sn)^\ell (\Xi_4+\f12\sX)\|_{L_\omega^\infty})\\
&\les \|(\tir \sn)^\ell (\bb^{-1}\sn \Xi_4)\|_{L_\omega^4}\l t\r^{-1}+\ell \l t\r^{-2}\Delta_0^\f12(\|\bb^{-2}(\tir \sn)^\ell \Xi_4\|_{L_\omega^\infty}+\|\bb^{-1}\tir \sF\|_{L_\omega^\infty}).
\end{align*}
Using Lemma \ref{2.26.3.24}, we deduce by using (\ref{10.26.4.23}), (\ref{1.28.1.24}) and (\ref{3.6.2.21}) that
\begin{align*}
\|(\tir \sn)^\ell(\bb^{-1}\sn\chih\c \chih)\|_{L_\omega^4}
&\les \|(\tir \sn)^\ell(\bb^{-1}\sn \chih)\|_{L_\omega^4}\|\chih\|_{L_\omega^\infty}+\|(\tir \sn)^\ell \chih \bb^{-1}\sn\chih\|_{L_\omega^4}\\
&\les \|(\tir \sn)^\ell(\bb^{-1}\sn \sX)\|_{L_\omega^4}\l t\r^{-\frac{7}{4}+\delta+}\Delta_0+\l t\r^{-2-\frac{7}{4}+2\delta}\Delta_0^2.
\end{align*}
Next we derive by using (\ref{3.6.2.21}) and (\ref{3.11.3.21})
\begin{align*}
\|(\tir \sn)^\ell(\bb^{-1}\chih\c \sF)\|_{L_\omega^4}&\les \|(\tir \sn)^{\ell}\sF\|_{L_\omega^4}\l t\r^{-\frac{7}{4}+\delta}\Delta_0+\|(\tir \sn)^\ell(\bb^{-1}\chih)\|_{L_\omega^4}\|\sF\|_{L^\infty_\omega}\\
&\les \|(\tir \sn)^{\ell}\sF\|_{L_\omega^4}\l t\r^{-\frac{7}{4}+\delta}\Delta_0+\l t\r^{-\frac{7}{4}+\delta}\Delta_0\|\sF\|_{L_\omega^\infty}\\
&\les \l t\r^{-\frac{7}{4}+\delta}\Delta_0 \|(\tir \sn)^{\le\ell} \sF\|_{L_\omega^4};
\end{align*}
and we have, also due to (\ref{10.26.4.23})
\begin{align*}
\|(\tir\sn)^\ell(\sn[L\Phi] \bb^{-1}\sX)\|_{L_\omega^4}&\les\l t\r^{-1} \|(\tir \sn)^\ell(\bb^{-1}\sn[L\Phi])\|_{L_\omega^4}+\ell\l t\r^{-3+\frac{1}{4}+\delta}\Delta_0\|(\tir \sn)^\ell\sX\|_{L_\omega^4}\\
&\les\l t\r^{-1} \|(\tir \sn)^\ell(\bb^{-1}\sn[L\Phi])\|_{L_\omega^4}+\l t\r^{-4+\frac{1}{4}+2\delta}\Delta_0^2.
\end{align*}
Finally, in view of (\ref{2.1.3.24}), by using (\ref{10.22.5.22}) and combining the above estimates, we conclude
\begin{align*}
\|(\tir \sn)^\ell{\ti G}_2\|_{L_\omega^4}&\les  \l t\r^{-\frac{7}{4}+\delta+}\Delta_0 \|(\tir \sn)^{\ell} \sF\|_{L_\omega^4}+\l t\r^{-2-\frac{7}{4}+2\delta}\Delta_0^\frac{3}{2}\\
&+ \|(\tir \sn)^\ell (\bb^{-1}\sn\fB)\|_{L_\omega^4}\l t\r^{-1}+\ell \l t\r^{-2}\Delta_0^\f12(\|\bb^{-2}(\tir \sn)^\ell \Xi_4\|_{L_\omega^\infty}+\|\tir^2\sn\sF\|_{L_\omega^4}).
\end{align*}
Using (\ref{1.29.4.22}), we derive
\begin{align*}
\|\tir (\bR_{BCLA}, \chi\c \zb)\sF\|_{L_\omega^4}\les \l t\r^{-2+\delta}\Delta_0\|(\tir\sn)^{\le 1}\sF\|_{L_\omega^4}.
\end{align*}
Combining the above estimates by using (\ref{12.13.1.21}), using the transport lemma we derive for $\ell=1$
\begin{align*}
\|\tir^3 (\tir \sn)^\ell \sF\|_{L_\omega^4}&\les \|\tir^3 (\tir \sn)^\ell \sF(0)\|_{L_\omega^4}+\int_0^t \tir^2 \{\|(\tir \sn)^\ell (\bb^{-1}\sn\fB)\|_{L_\omega^4}\\
&+\ell \l t'\r^{-1}\Delta_0^\f12(\|\bb^{-2}(\tir \sn)^\ell \Xi_4\|_{L_\omega^\infty}+\|\tir^2\sn\sF\|_{L_\omega^4})\}+\l t\r^{\frac{1}{4}+2\delta}\Delta_0^\frac{3}{2}\\
&+\int_0^t \l t'\r^{\frac{5}{4}+\delta}\Delta_0\|\sF\|_{L_\omega^4}. 
\end{align*}
We then substitute (\ref{10.26.2.23}) into the last line in the above to obtain (\ref{9.6.3.22}). 
\end{proof}
By using (\ref{10.26.3.23}) and (\ref{1.29.1.24}), we further obtain the following result.
\begin{proposition}
\begin{enumerate}
\item There hold
\begin{align}\label{3.31.2.22}
\begin{split}
\|\tir^3(\tir\sn)^2\sF\|_{L_u^2 L_\omega^2}&\les \int_0^t\|\tir\{(\tir \sn)^3(\bb^{-1}\Lb\varrho), \bb^{-1}(\tir\sn)^3\Lb\varrho\}\|_{L_u^2 L_\omega^2}\\
&+\int_0^t \l t'\r^2 \log \l t'\r\|(\tir \sn)^2\sn[L\Phi]\|_{L_u^2 L_\omega^2}+(\log \l t\r)^\frac{15}{2}(\La_0+\Delta_0^\frac{5}{4})\l t\r
\end{split}
\end{align}
and 
\begin{align}\label{6.2.3.24}
\|\bb^{-1}\tir (\tir\sn)^{1+\le 2}&\Lb \varrho,\,\tir(\tir \sn)^3(\bb^{-1}\Lb\varrho)\|_{L_u^2 L_\omega^2}\nn\\
&\les \|\bb^{-1}\tir\Lb \Omega^3\varrho\|_{L_u^2 L_\omega^2}+\Delta_0^\frac{5}{4}\log \l t\r^7+\log \l t\r\La_0\\
&+(C\log \l t\r^{-1}\M_0+\log \l t\r^{-\frac{3}{2}}\Delta_0^\f12)\big(\La_0+\Delta_0^\frac{5}{4}\nn\\
&+\int_0^t \|\bb^{-1}\Lb\Omega^3\varrho, \bb^{-1}\Omega^3[L\Phi]\|_{L_u^2 L_\omega^2}\big).\nn
\end{align}
 We also have
\begin{align}\label{2.25.2.24'}
\|\tir^2(\tir \sn)\sF\|_{L^2_u L_\omega^2}\les(\log \l t\r)^\frac{17}{4}(\Delta_0^\frac{5}{4}+\La_0).
\end{align}
\item

 For $\vs^+(X^2)=1$, there holds with $\vs(V)=\vs^-(X^2)$
\begin{align}
\|\sn_X^2 \Omega \tr\chi-(1-\vs^-(X^2))\Omega^2\tr\chi\|_{L^2_\Sigma}&\les\l t\r^{-1}(\l t\r^{2\delta}\Delta_0^\frac{5}{4}+\log \l t\r^{\f12\M+1}\La_0)\nn\\
&+\l t\r^{-\frac{1}{2}}W_1[V\Omega S\Phi, \Omega^3\varrho]^\f12(t).\label{9.15.4.22}
\end{align}
\end{enumerate}
\end{proposition}
\begin{proof}
Substituting  Lemma \ref{5.13.11.21} (5) and (\ref{9.20.4.22}), (\ref{2.14.1.24}), (\ref{L2BA2'}) to (\ref{10.26.3.23}) gives (\ref{2.25.2.24'}).

Due to $\Xi_4=\Lb \varrho+[L\Phi]$, also using Lemma \ref{5.13.11.21} (5) and Proposition \ref{7.15.5.22}, we derive
\begin{align}
\|\tir(\tir \sn)^3(\bb^{-1}\Xi_4)\|_{L_u^2 L_\omega^2}&\les\|\tir(\tir \sn)^3(\bb^{-1}\Lb\varrho)\|_{L_u^2 L_\omega^2}+\|\tir(\tir \sn)^3(\bb^{-1}[L\Phi])\|_{L_u^2 L_\omega^2}\nn\\
&\les \|\bb^{-1}\tir(\tir \sn)^3[L\Phi], \tir(\tir \sn(\bb^{-1}) (\tir \sn)^2 [L\Phi])\|_{L_u^2 L_\omega^2}\nn\\
&+\|\tir(\tir \sn)^3(\bb^{-1}\Lb\varrho)\|_{L_u^2 L_\omega^2}+\l t\r^{-1+2\delta}\Delta_0^2\label{6.2.1.24}\\
&\les \|\tir(\tir \sn)^3(\bb^{-1}\Lb\varrho), \bb^{-1}\tir(\tir \sn)^3[L\Phi]\|_{L_u^2 L_\omega^2}+\l t\r^{-1+2\delta}\Delta_0^2.\nn 
\end{align} 
For the lower order estimate, we similarly bound by using (\ref{2.14.1.24}) that
\begin{align}\label{6.2.2.24}
\|(\tir \sn)^2(\bb^{-1}\Xi_4)\|_{L_u^2 L_\omega^2}\les \l t\r^{-1}\log \l t\r(\La_0+\log \l t\r^\frac{5}{2}\Delta_0^\frac{5}{4}).
\end{align}
It remains to treat the most crucial terms $(\tir \sn)^3(\bb^{-1}\Lb \varrho)$ and $\bb^{-1}(\tir \sn)^3 \Lb \varrho$. We use (\ref{7.03.5.21}), (\ref{7.13.5.22}), (\ref{2.14.1.24}) and Sobolev embedding on spheres to derive
\begin{align*}
(\tir\sn)^3\Lb \varrho&=\Omega^{1+\le 2}\Lb \varrho+O(\l t\r^{-\frac{3}{4}+\delta}\Delta_0^\f12)_{L_\omega^4}\Omega \Lb \varrho\\
&=\Omega^3\Lb \varrho+O(\log \l t\r\l t\r^{-1}(\La_0+\log \l t\r^\frac{5}{2}\Delta_0^\frac{5}{4}))_{L_u^2 L_\omega^2}.
\end{align*}
Similarly, with $m=2, 3$,
\begin{align*}
(\tir\sn)^m(\bb^{-1}\Lb \varrho)&=\Omega^{1+\le m-1}(\bb^{-1}\Lb\varrho)+O(\l t\r^{-\frac{3}{4}+\delta}\Delta_0^\f12)_{L_\omega^4}\Omega(\bb^{-1}\Lb\varrho)\\
&=\Omega^m(\bb^{-1}\Lb\varrho)+O(\log \l t\r\l t\r^{-1}(\La_0+\log \l t\r^{\frac{5}{2}(m-2)}\Delta_0^\frac{5}{4}))_{L_u^2 L_\omega^2}.
\end{align*} 
Hence, 
 applying (\ref{2.22.1.24}) we bound 
\begin{align*}
\|\bb^{-1}\tir (\tir\sn)^{1+\le 2}\Lb \varrho,\,& \tir(\tir \sn)^3(\bb^{-1}\Lb\varrho)\|_{L_u^2 L_\omega^2}\\&\les \|\bb^{-1}\tir\Omega^3\Lb \varrho, \tir \Omega^3(\bb^{-1}\Lb\varrho)\|_{L_u^2 L_\omega^2}+\log \l t\r(\La_0+\log \l t\r^\frac{5}{2}\Delta_0^\frac{5}{4})\\
&\les \|\bb^{-1}\tir\Lb \Omega^3\varrho\|_{L_u^2 L_\omega^2}+\Delta_0^\frac{5}{4}\log \l t\r^7+\log \l t\r\La_0\\
&+(C\log \l t\r^{-1}\M_0+\log \l t\r^{-\frac{3}{2}}\Delta_0^\f12)\big(\La_0+\Delta_0^\frac{5}{4}\\
&+\int_0^t  \{\|\bb^{-1}\Lb\Omega^3\varrho\|_{L_u^2 L_\omega^2}+\|\Omega^3(\bb^{-1}[L\Phi])\|_{L_u^2 L_\omega^2}\}\big).
\end{align*}
Here we note that the last term in the above can be similarly treated as in (\ref{6.2.1.24}).
Substituting the above estimate and (\ref{2.25.2.24'}) into (\ref{10.26.3.23}), using (\ref{8.29.9.21}), (\ref{6.2.1.24}), Cauchy-Schwarz and (\ref{6.2.2.24}) we obtain
\begin{align*}
&\|\tir^3(\tir\sn)^2\sF\|_{L_u^2 L_\omega^2}\\
&\les \La_0+\int_0^t\|\tir\{(\tir \sn)^3(\bb^{-1}\Xi_4), \bb^{-1}(\tir\sn)^3\Xi_4\}\|_{L_u^2 L_\omega^2}+(\log \l t\r)^5(\La_0+\Delta_0^\frac{5}{4})\l t\r\\
&+\int_0^t \l t'\r^2 \log \l t'\r\|(\tir \sn)^{\le 2}\sn[L\Phi]\|_{L_u^2 L_\omega^2}+\int_0^t \Delta_0 \l t'\r^{-1}(\log \l t'\r)^\frac{13}{4}\|\tir^3(\tir \sn)\sF\|_{L_u^2 L_\omega^2}\\
\displaybreak[0]
&\les \int_0^t\|\tir\{(\tir \sn)^3(\bb^{-1}\Lb\varrho), \bb^{-1}(\tir\sn)^3\Lb\varrho\}\|_{L_u^2 L_\omega^2}+\int_0^t \l t'\r^2 \log \l t'\r\|(\tir \sn)^2\sn[L\Phi]\|_{L_u^2 L_\omega^2}\\
&+\l t\r(\log \l t\r)^\frac{15}{2}(\La_0+\Delta_0^\frac{5}{4}).
\end{align*}
 Hence we conclude (\ref{3.31.2.22}) and (\ref{6.2.3.24}).


Next we consider (2).
To prove (\ref{9.15.4.22}), we rewrite (\ref{9.15.3.22}) below
\begin{align*}
\sn_X \sn_L (\tir^3\sn \tr\chi)&=\sn_X\{\tir^3(\bA_b+\bA_{g,2})\c \sn \tr\chi+\tir^3 (\chih\sn \chih+\sn(\tr\chi[L\Phi])+\sn \wt{L \Xi_4})\\
&+\tir^3 \sn \N(\Phi, \bp\Phi)\}.
\end{align*}
Using (\ref{9.8.2.22}) and (\ref{3.20.1.22}), we first directly bound
\begin{align*}
|\sn_X\big(\tir^3 \sn\wt{L \Xi_4}\big)|&\les|\tir^2 X^{\le 1}\big(\Omega(\sD \varrho-\hb L \varrho+\zb\c \sn \varrho-\Box_\bg \varrho)\big)|\\
&+|\tir^2 X^{\le 1}\Omega((L+(h-k_{\bN\bN}))L \varrho)|.
\end{align*}
Taking $L^2_\Sigma$, by using (\ref{3.16.1.22}), (\ref{2.27.1.24}), (\ref{12.22.4.23}), Proposition \ref{1steng} and Proposition \ref{8.29.8.21}, 
\begin{equation*}
\|X^{\le 1}\{\Omega(h L \varrho), \Omega(\zb\sn\varrho)\}\|_{L^2_\Sigma}\les  \l t\r^{-2}(\log \l t\r)^{\f12\M}(\La_0+\Delta_0^\frac{5}{4}).
\end{equation*}
Taking $L^2_\Sigma$, by using  (\ref{7.28.2.21}), (\ref{7.17.2.21}) and (\ref{8.29.9.21}), we derive
\begin{align*}
\|X^{\le 1}\Omega L^2\varrho\|_{L^2_\Sigma}&\les \l t\r^{-\frac{3}{2}}W_1[X\Omega S\varrho]^\f12(t)+\l t\r^{-2}\log \l t\r^{\f12\M}(\La_0+\Delta_0^\frac{5}{4})\\
\|X^{\le 1}\Omega\sD \varrho\|_{L^2_\Sigma}&\les \l t\r^{-\frac{3}{2}}W_1[\Omega^3\varrho]^\f12(t)+\l t\r^{-2}\log \l t\r^{\f12\M}(\La_0+\l t\r^\delta \Delta_0^\frac{5}{4}).
\end{align*}
Recall from (\ref{8.30.3.21}) and (\ref{12.21.3.22}) that
\begin{align*}
\|X^{\le 1}\Omega\left(\Box_\bg \varrho, \fB L \varrho, \N(\Phi, \bp\Phi)\right)\|_{L^2_\Sigma}\les   \l t\r^{-2}((\log \l t\r)^{\f12\M}\La_0+\l t\r^{2\delta}\Delta_0^\frac{5}{4}).
\end{align*}
where we used the fact that (\ref{12.21.3.22}) can be extended with $\Box_\bg \Phi$ replaced by $\N(\Phi, \bp\Phi)$ without changing the bound on the right-hand side. 

We directly derive, by using (\ref{3.16.1.22}), (\ref{2.27.1.24}), Proposition \ref{7.15.5.22} and Proposition \ref{8.29.8.21}, that
\begin{align*}
&\|\sn_X((\tir^3(\bA_b+\bA_{g,2})+\tir^3[L\Phi])\c \sn \tr\chi+\tir^3(\tr\chi \sn[L\Phi], \chih\sn\chih)\|_{L^2_\Sigma}\\
&\les \|\tir\sn_X^{\le 1}(\tir\sn_\Omega\bA\c \bA)\|_{L^2_\Sigma}+\|\sn_X(\tir\tr\chi)\c\tir \Omega[L\Phi]\|_{L^2_\Sigma}+\|\tir\sn_X\Omega[L\Phi]\|_{L^2_\Sigma}\\
&\les \log \l t\r^{\f12\M}(\La_0+\Delta_0^\frac{5}{4}).
\end{align*}
Therefore we conclude
\begin{align}\label{6.2.5.24}
\begin{split}
\|\sn_X \sn_L(\tir^3 \sn\tr\chi)\|_{L_\Sigma^2}&\les \log \l t\r^{\f12\M}\La_0+\l t\r^{2\delta}\Delta_0^\frac{5}{4}\\
&+\l t\r^\f12(W_1[X\Omega S\varrho]^\f12(t)+W_1[\Omega^3\varrho]^\f12(t)).
\end{split}
\end{align}
Due to (\ref{4.22.4.22}) and (\ref{1.25.2.22}), we write
\begin{align*}
XS\Omega\tr\chi&=X^{\le 1}\Omega \tr\chi+\tir^{-1}XL(\tir^2\Omega\tr\chi)\\
&=X^{\le 1}\Omega\tr\chi+\tir^{-1}\sn_X \sn_L(\tir^3 \sn\tr\chi)+O(\l t\r^{-\frac{3}{4}+\delta}\Delta_0^\f12)_{L_\omega^4}\Omega \tr\chi.
\end{align*}
Hence, we conclude by using (\ref{3.11.3.21}), (\ref{2.18.2.24}) and (\ref{L2conndrv'}) that
\begin{align*}
\|XS\Omega \tr\chi&-O(1)(1-\vs(X)) X\Omega \tr\chi\|_{L^2_\Sigma}\\
&\les \l t\r^{-1}\|\sn_X \sn_L(\tir^3\sn\tr\chi)\|_{L^2_\Sigma}+\|S^{\le 1}\Omega\tr\chi\|_{L^2_\Sigma}+\l t\r^{-\frac{3}{2}+2\delta}\Delta_0^\frac{3}{2}\log \l t\r\\
&\les \l t\r^{-1}\|\sn_X \sn_L(\tir^3\sn\tr\chi)\|_{L^2_\Sigma}+\l t\r^{-1}\log \l t\r^{\f12\M+1}(\l t\r^\delta\Delta_0^\frac{5}{4}+\La_0).
\end{align*}
Thus we completed the proof of (\ref{9.15.4.22}) for $X^2=\Omega S, SS$ by substituting (\ref{6.2.5.24}) to the last line in the above. For $X^2=S\Omega$, we take direct commutation in (\ref{3.21.1.23}) by using the result of $X^2=\Omega S$. We will skip the details here.
\end{proof}
Next we provide more top order derivative estimates of $\tr\chi-\frac{2}{\tir}$ and $\chih$. 
\begin{proposition}\label{2.28.2.24}
Let $X\in \{\Omega, S\}$ and $0<t<T_*$. Under the assumptions of (\ref{3.12.1.21})-(\ref{6.5.1.21}), we have
\begin{align}
\|\bb^{-\f12}\sn_X^2\sn_L\chih\|_{L^2_\Sigma}&\les\|\sn_X^2\chih, \tir\sn_X^2(\sn, \sn_L)\bA_{g,1}\|_{L^2_u L_\omega^2}\nn\\
&+ \l t\r^{-2+\delta}(\log \l t\r\La_0+\l t\r^{\delta+\frac{1}{4}}\Delta_0^\frac{5}{4}),\label{9.17.10.22}\\
\|\bb^{-\f12}\sn_X^2 \sn_L(\tr\chi-\frac{2}{\tir})\|_{L^2_\Sigma}&\les \|\bb^{-\f12}(X^2\sD\varrho, X^2L^2\varrho)\|_{L^2_\Sigma}\nn\\
&+\l t\r^{-2}(\log \l t\r^{\f12\M}\La_0+\l t\r^{2\delta}\Delta_0^\frac{5}{4}),\label{9.18.1.22}\\
\|\bb^{-\f12}\sn_X^2\sn_\Lb\chih\|_{L^2_\Sigma}&\les \|\bb^{-\f12}(\sn_X^2\bar\bp\bA_{g,1}, \sn_X^2\sn\ze, \tir^{-1}\sn_X^2\bA_g)\|_{L^2_\Sigma}+\l t\r^{-1+2\delta}\Delta_0^\frac{3}{2}.\label{9.18.2.22}
\end{align}
\end{proposition}
\begin{proof}
We first consider (\ref{9.17.10.22}) by using (\ref{s2}) and (\ref{1.30.2.22})
\begin{align*}
\sn_X^2\sn_L\chih&=\sn_X^2\big((\tir^{-1}+\bA+\fB)\chih\big)+\sn_X^2(\widehat{\bR_{4A4B}})\\
&=\sn_X^2\{\sn[\sn\Phi]+\sn_L \hk_{AB}+[\sn\Phi]^2+\hk_{AB}(\fB+\tr\chi)+(\chih, \chibh)(\fB+\hk_{AB})\}\\
&+\tir^{-1}\sn_X^{\le 2}\chih+\sn_X^2((\bA+\fB)\c \chih).
\end{align*}
Rewriting the lower order terms in the above, we further deduce, with the help of (\ref{8.9.4.22}), Proposition \ref{7.15.5.22},  (\ref{2.18.2.24}),  (\ref{8.8.6.22'}), (\ref{8.28.2.23}) and (\ref{L2BA2'}), that
\begin{align*}
\|\sn_X^2((\bA+\fB)\bA_g+\bA_{g,1}\tir^{-1})\|_{L^2_u L_\omega^2}&\les \l t\r^{-1}\|\sn_X^2 \bA_{g,2}\|_{L_u^2 L_\omega^2}+\l t\r^{-3+\delta}(\log \l t\r\La_0+\l t\r^{\frac{1}{4}+\delta}\Delta_0^\frac{5}{4}).
\end{align*}
Combining the above two calculations gives (\ref{9.17.10.22}).

Next, we rewrite (\ref{6.3.1.23}) as
\begin{align*}
(L+\frac{2}{\tir}+\f12\bA_b)(\tr\chi-\frac{2}{\tir})&=\bA_{g,2}^2+\wt{L\Xi_4} +\tr\chi[L \Phi]+\N(\Phi, \bp\Phi).
\end{align*}
Using (\ref{3.20.1.22}), applying $X^2$ to the above identity, we have
\begin{align*}
X^2 L\bA_b&=X^2\big((\tir^{-1}, \bA_b)\bA_b\big)+X^2 (\bA^2_{g,2}+\tr\chi [L \Phi]+\N(\Phi, \bp\Phi))\\
&+X^2\{\sD \varrho-\hb L \varrho+2\zb \sn\varrho+(L+(h-k_{\bN\bN}))L \varrho\}\\
&=X^2\big(\tir^{-1}\bA_b, \bA^2, \N(\Phi, \bp\Phi)\big)+X^2\big([L\Phi](\tir^{-1}+\bA+\fB)\big)+X^2(\sD\varrho, L^2\varrho).
\end{align*}
Using (\ref{8.29.9.21}) and (\ref{2.27.1.24}), we derive
\begin{align*}
\|\bb^{-\f12} X^2\{(\bA, \tir^{-1})[L\Phi], \bA^2\}\|_{L^2_\Sigma}&\les\l t\r^{-2}\log \l t\r^{\f12\M}(\Delta_0^\frac{5}{4}+\La_0).
\end{align*}
Note $\fB[L\Phi]$ can be included in $\N(\Phi, \bp\Phi)$. Recall from (\ref{8.30.3.21}) that
\begin{equation}\label{2.29.4.24}
\|X^2(\N(\Phi, \bp \Phi), \fB [L\Phi])\|_{L^2_\Sigma}\les \l t\r^{-2} \log \l t\r^{\f12\M}(\La_0+\Delta_0^\frac{5}{4})+(1-\vs^-(X^2))\l t\r^{-2+2\delta}\Delta_0^\frac{3}{2}.
\end{equation}
 Hence we conclude
\begin{align*}
X^2 L(\tr\chi-\frac{2}{\tir})&=\tir^{-1}X^{\le 2} \bA_b+X^2\sD\varrho+X^2L^2\varrho\\
&+O(\l t\r^{-3}(\log \l t\r^{\f12\M}\La_0+\l t\r^{2\delta}\Delta_0^\frac{5}{4}))_{L^2_u L_\omega^2}.
\end{align*}
This implies (\ref{9.18.1.22}).

To see (\ref{9.18.2.22}), using (\ref{3chi}), we write
\begin{equation*}
\sn_X^2 \sn_\Lb \chih=\sn_X^2\left((\bA_b+\fB+\tir^{-1})\chih+\sn \ud\bA+\ud\bA^2+\widehat{\bR_{A43B}}\right).
\end{equation*}
It follows by using (\ref{8.9.4.22}), Proposition \ref{7.15.5.22} for the estimates of $\bA_g$ and (\ref{8.28.2.23}) that 
\begin{align*}
\sn_X^2\left((\bA+\fB)\bA_g\right)= O(\bb^{-1}\l t\r^{-1})\sn_X^{\le 2}\bA_g+O(\l t\r^{-\frac{11}{4}+2\delta}\Delta_0^2)_{L_u^2 L_\omega^2}.
\end{align*}
Recall from (\ref{1.30.1.24}) that $$
\|\sn_X^2(\ud \bA^2)\|_{L_u^2 L_\omega^2}\les \l t\r^{-2+2\delta}\Delta_0^2.$$ 
Hence 
\begin{align*}
\|\sn_X^2&\left((\bA_b+\fB+\tir^{-1})\bA_g+\ud\bA^2\right)\|_{L_u^2 L_\omega^2}\les\|\tir^{-1}\sn_X^2 \bA_g\|_{L_u^2 L_\omega^2} + \l t\r^{-2+2\delta}\Delta_0^2.
\end{align*}

Finally, using (\ref{5.21.3.23}), the above error estimate and (\ref{8.28.2.23}), we bound
\begin{align*}
\sn_X^2(\widehat{\bR_{A43B}})&=\sn_X^2  \{\sn[\sn\Phi]+\sn_L \hk_{AB}+[\sn\Phi]^2+\bA_g(\fB+\bA+\tir^{-1})\}\nn\\
&=\sn_X^2\bar\bp\bA_{g,1}+O(\bb^{-1}\l t\r^{-1})\sn_X^{\le 2}\bA_g+ O(\l t\r^{-\frac{11}{4}+2\delta}\Delta_0^\frac{3}{2})_{L^2_u L_\omega^2}
\end{align*}
Thus we conclude (\ref{9.18.2.22}).
\end{proof}
\begin{corollary}\label{2.29.1.24}
Let $X\in \{\Omega, S\}$ and $0<t<T_*$. Under the assumptions of (\ref{3.12.1.21})-(\ref{6.5.1.21}), we have
\begin{align*}
\|X^2 \Omega \tr\chi\|_{L_u^2 L_\omega^2}\les \l t\r^{-1+\delta-\f12\vs^+(X^2)}\Delta_0
\end{align*}
\end{corollary}
\begin{proof}
The case $\vs^+(X^2)=1$ follows by using (\ref{9.18.1.22}), (\ref{7.17.7.21}) and (\ref{L2BA2}). If $\vs^+(X^2)=0$, the estimate can be derived by using Proposition \ref{7.15.5.22} and Lemma \ref{3.17.2.22}. 
\end{proof}

In the sequel, we prove the top order geometric estimates.
\begin{proposition}\label{8.10.2.21}
Let $X\in \{\Omega, S\}$ and $u_0\le u\le u_*$ and $0<t<T_*$. Under the assumptions of (\ref{3.12.1.21})-(\ref{6.5.1.21}),

(1)Symbolically, there hold
\begin{equation}\label{8.8.6.21}
\begin{split}
\sn_X^2\bJ[{}\rp{a}\Omega]_L&=\sn^2_X(\ckr^{-1}{}\rp{a}v^\sharp)+ X^2 (\Omega\tr\thetac+\tir^{-2}\la)\\
&+[LX^2\Omega\Phi]+O(\l t\r^{-\frac{3}{2}+2\delta}\Delta_0^\frac{3}{2}+\l t\r^{-2}\La_0\log \l t\r^{\f12\M})_{L_u^2 L_\omega^2}\\
\sn_X^2 \bJ[{}\rp{a}\Omega]_\Lb&= O(\l t\r^{3\delta-1}\Delta_0)_{L_u^2 L_\omega^2}+[LX^2\Omega\Phi]\\
\sn_X^2 {}\rp{a}\bJ_B&= O(\fB)+O(\Delta_0\l t\r^{-1+2\delta})_{L^2_u L_\omega^2}.
\end{split}
\end{equation}
(2)
\begin{align}
\|\sn_X^2 \bJ[S]_B\|_{L^2_u L_\omega^2}&\les \l t\r^{-1+\delta}(\log \l t\r)^2\Delta_0\label{8.13.1.22}\\
X^2\bJ[S]_\Lb&=O(\l t\r^{-1})+O(\l t\r^{-1+2\delta}\Delta_0)_{L^2_u L_\omega^2}\label{8.13.2.22}\\
 \|X^2L{}\rp{S}\ss\|_{L^2_\Sigma}&\les\l t\r^{-\f12}\{W_1[X_2 X_1 S\varrho]^\f12(t)+ \sum_{\vs(V)=\vs^-(X^2)}W_1[V\Omega^2\varrho]^\f12(t)\}\nn\\
&+\l t\r^{-1}(\log \l t\r)^\frac{\M}{2}(\La_0+\l t\r^\delta\Delta_0^\frac{5}{4}).\label{8.14.1.22}
\end{align}
\end{proposition}
We will rely on the following geometric estimates to prove the above results.
\begin{proposition}\label{gmtrc_high_od}
Let $X\in \{\Omega, S\}$ and $0<t<T_*$. Under the assumptions of (\ref{3.12.1.21})-(\ref{6.5.1.21}), we have
\begin{align}
&\|\bb^{-\f12}\sn_X^2(L\bA_b, \sn_L\bA_{g,2}), \bb^{-\f12}\l t\r^{-\f12}\sn_X^2 \sn_L \ze\|_{L^2_\Sigma}\les \l t\r^{-\frac{3}{2}+\delta}\Delta_0 \label{8.7.2.22}\\
&\|\bb^{-\f12}\sn_X^2\sn (\bA_{g,2}, \bA_b)\|_{L^2_\Sigma}\les \l t\r^{-1-\f12\vs^+(X^2)+\delta}\Delta_0, \label{8.7.1.22}\\
\displaybreak[0]
&\sn_X^2(\Lb[L\Phi], L[\Lb\Phi])=\sn_X^2(\sD\varrho+L[L\Phi])+O(\l t\r^{-1}\fB)\vs^-(X^2)\nn\\
&\qquad\qquad\qquad\qquad+O( \l t\r^{-2}(\log \l t\r)^4\Delta_0)_{L_u^2 L_\omega^2}\label{2.29.5.24}\\
&\|\sn_X^2(\mu-\tr\chi k_{\bN\bN}-\varpi)\|_{L^2_u L_\omega^2}\les \l t\r^{-1-\vs^+(X^2)+2\delta}\Delta_0\label{6.5.3.24}\\
&\sn_X^2\Lb\bA_b=O(\l t\r^{-2})\vs^-(X^2)+O(\l t\r^{-1+2\delta-\vs^+(X^2)}\Delta_0)_{L_u^2 L_\omega^2}\label{8.12.2.22}\\
&\|\sn\sn_\Omega^2\ze, \sn_\Omega^2 \sn_\Lb\bA_{g,2}\|_{L_u^2 L_\omega^2}\les \l t\r^{-1+2\delta}\Delta_0, \label{3.15.11.24}\\
&\|\bb^{-\f12}\sn_X^2\sn_L \sn \la\|_{L^2_\Sigma}\les \l t\r^{\delta-\f12}\Delta_0.\label{8.9.2.22}\\
&\|\Omega^2\Lb\N(\Phi, \bp\Phi)\|_{L^2_\Sigma}\les \|\bb^{-2}\Omega^2[L\Phi]\|_{L_u^2 L_\omega^4}+\|\fB\Omega(\sD\varrho, L[L\Phi])\|_{L^2_\Sigma}\nn\\
&\qquad\qquad\qquad\qquad\,\,+\l t\r^{-2} (\log \l t\r\La_0+\l t\r^{2\delta}\Delta_0^\frac{5}{4}).\label{3.8.1.24}
\end{align}
\end{proposition}

\begin{remark}
The decay estimate (\ref{3.15.11.24}) is nearly $\l t\r^{-1}$ weaker than the corresponding lower order behaviors. Fortunately, they do not appear as the highest order linear term in the top order commutators, due to the delicate cancelation structure in $\bJ[\Omega]_A$. 
\end{remark}
\begin{proof} 
 (\ref{8.7.2.22}) follows by applying Proposition \ref{7.15.5.22} to (\ref{9.17.10.22}) and (\ref{9.18.1.22}) and  the derivative estimate of $\sn_L \ze$ is from Proposition \ref{7.15.5.22}.  The estimates (\ref{8.7.1.22}) with $\vs^+(X^2)=1$ can be obtained by  employing the first set of estimates in (\ref{8.7.2.22}) followed with applying (\ref{3.20.2.23}) to $G=\bA_{g,2}, \bA_b$.  If $\vs^+(X^2)=0$ in (\ref{8.7.1.22}), they follow from Proposition \ref{7.15.5.22}.

Noting that $|X^l(\tr\chi, \fB)|\les \l t\r^{-1} (\l t\r^\delta\Delta_0)^{l(1-\vs(X^l))}$, for $l=0,1$,  and using (\ref{7.29.1.22}) and  Proposition \ref{7.15.5.22}, we have
\begin{align*}
X^2((\tr\chi, \fB)\bA)&=O(\l t\r^{-1})X^{1+\le 1}\bA+X^2(\tr\chi, k_{\bN\bN})\bA+O(\l t\r^{-3+2\delta}\Delta_0^2)_{L_u^2 L_\omega^2}\\
&=O(\l t\r^{-1})X^{\le 2}\bA+O(\l t\r^{-3+2\delta}\Delta_0^\frac{3}{2})_{L^2_u L_\omega^2}.
\end{align*}
Hence using Proposition \ref{7.15.5.22} again, we deduce
\begin{align}
\|X^2((\tr\chi, \fB)\bA)\|_{L^2_u L_\omega^2}\les \l t\r^{-\frac{5}{2}+\delta}\Delta_0.\label{9.19.2.22}
\end{align}
 Next we consider $X^2\big((\Lb [L\Phi], L[\Lb\Phi])-\sD\varrho-L[L\Phi]\big)$ by using (\ref{8.26.1.23}), (\ref{1.30.1.24}),  (\ref{7.29.1.22}) and  (\ref{2.29.4.24})  
\begin{align}\label{3.16.4.24}
\begin{split}
X^2\Big(&(\Lb [L\Phi], L[\Lb\Phi])-\sD\varrho-L[L\Phi]\Big)\\
&=X^2\Big((\tir^{-1}+\fB)\fB\Big)+X^2(\bA_b\fB, \ud \bA\bA_{g,1}, \N(\Phi, \bp\Phi))\\
&=\vs^-(X^2)O(\l t\r^{-1}\fB)+X^2(\bA_b \fB)+O(\l t\r^{-2}\log \l t\r(\La_0+\log \l t\r^3\Delta_0^\frac{5}{4}))_{L_u^2 L_\omega^2}.
\end{split}
\end{align}
 Applying (\ref{9.19.2.22}) to the term $X^2(\bA_b\fB)$, we then conclude (\ref{2.29.5.24}).

Next we prove (\ref{8.12.2.22}) and (\ref{3.15.11.24}). In view of (\ref{8.11.9.22}), we derive
\begin{align}\label{3.16.2.24}
\begin{split}
X^2 \Lb (\tr\chi-\frac{2}{\tir})&=X^2(\mu-\tr\chi k_{\bN\bN}-\varpi)+X^2(\tir^{-1} k_{\bN\bN}+\varpi)\\
&+\f12 X^2((\tr\chi, k_{\bN\bN}) \bA_b)+X^2(\tir^{-1}\mho).
\end{split}
\end{align}

Using (\ref{dze}) and (\ref{4.17.1.24}), we write
\begin{align}\label{6.5.2.24}
\sn_X^2(\mu-\tr\chi k_{\bN\bN}-\varpi)=\sn_X^2(\sdiv\ze+\sn\bA_{g,1}+\ud\bA^2+\bA_g\c \bA_g+\fB\bA_{g,2}).
\end{align}
Applying Proposition \ref{7.15.5.22}, (\ref{8.28.2.23}), (\ref{1.30.1.24}) and (\ref{9.19.2.22})
we obtain the case $\vs^+(X^2)=1$ in (\ref{6.5.3.24}).
Moreover due to (\ref{9.19.2.22}), for $\ell=1,2$
\begin{align}\label{6.5.1.24}
&X^\ell\varpi=X^{\ell}((\bA+\frac{1}{\tir}+\fB)\fB)=X^{\ell}\Big((\tir^{-1}+\fB)\fB\Big)+O(\l t\r^{-\frac{5}{2}+\delta})_{L_u^2 L_\omega^2}.
\end{align}
Then in the case of $\vs^+(X^2)=1$, (\ref{8.12.2.22}) follows by using (\ref{9.19.2.22}), Proposition \ref{7.22.2.22}, (\ref{7.29.1.22}) and (\ref{8.23.1.23}). 
For the case $\vs^+(X^2)=0$, i.e. $X^2=\Omega^2$, we similarly derive in view of (\ref{3.16.2.24}) and (\ref{8.23.2.23}) that
\begin{align*}
\Omega^2\Lb\bA_b&=\Omega^2(\mu-\tr\chi k_{\bN\bN}-\varpi)+O(\l t\r^{\delta-2}\Delta_0)_{L^2_u L_\omega^2}.
\end{align*}
It remains to prove (\ref{6.5.3.24}) for $X^2=\Omega^2$. For this purpose, we consider $\|\Omega^2 \Lb \sX\|_{L^2_u L_\omega^2}$ by using the transport equation (\ref{8.31.4.19}). It is direct to write
\begin{align*}
L(\tir^2 \Lb \sX)=\tir^2\{(L+\tr\chi) \Lb \sX-\bA_b \Lb \sX\}.
\end{align*}
Using (\ref{5.13.10.21}) we derive
\begin{align*}
L\Omega^2(\tir^2 \Lb \sX)&=[L,\Omega^2](\tir^2 \Lb \sX)+\tir^2 \Omega^2((L+\tr\chi) \Lb \sX-\bA_b \Lb \sX)\\
&=\tir^2\{ \pioh_{A L}\sn(\Omega \Lb \sX)+\Omega(\pioh_{AL}\sn(\Lb \sX))+\Omega^2((L+\tr\chi) \Lb \sX-\bA_b \Lb \sX)\}.
\end{align*}
Using (\ref{5.21.1.21}) and taking $L^2_u L_\omega^2$ norm, we infer that
\begin{align*}
\|L(\tir^2\Omega^2 \Lb \sX)\|_{L_u^2 L_\omega^2}&\les \l t\r^{-\frac{7}{4}+\delta}\Delta_0^\f12 \|\tir^2 \Omega^{1+\le 1} \Lb \sX\|_{L_u^2 L_\omega^2}\\
&+\|\tir^2\Omega^2((L+\tr\chi) \Lb \sX-\bA_b \Lb \sX)\|_{L^2_u L_\omega^2}
\end{align*}
where we applied Sobolev embedding in the above. Using the transport lemma and Proposition \ref{12.21.1.21}, we deduce
\begin{align}\label{3.16.1.24}
\|\tir^2 \Omega^2 \Lb \sX\|_{L^2_u L_\omega^2}&\les\La_0+ \int_0^t \{\l t'\r^{-\frac{7}{4}+\delta}\Delta_0^\f12 \|\tir^2 \Omega\Lb\sX\|_{L^2_u L_\omega^2}\\
&+\|\tir^2\Omega^2((L+\tr\chi) \Lb \sX-\bA_b \Lb \sX)\|_{L^2_u L_\omega^2}\} dt.\nn
\end{align}
For the last term in the above, we compute
\begin{align*}
\|\Omega^2(\bA_b \Lb\sX)\|_{L^2_u L_\omega^2}&\les \|\bA_b\|_{L^\infty_\omega}\|\Omega^2\Lb \sX\|_{L_u^2 L_\omega^2}+\|\Omega\bA_b\|_{L_\omega^4}\|\Omega^{1+\le 1}\Lb\sX\|_{L_u^2 L_\omega^2}\\
&+\|\Omega^2\bA_b\|_{L_u^2 L_\omega^4}\|\Lb \sX\|_{L_\omega^4}
\end{align*}
Due to Proposition \ref{7.22.2.22} and (\ref{1.29.2.22}), symbolically,
\begin{align}
\Lb\sX=\Lb(\Xi_4)+\Lb \bA_b+\tir^{-1}\tr\chib+\tir^{-1}\mho=O(\l t\r^{-1})_{L^4_\omega}.\label{6.4.2.24}
\end{align}
Combining the above two results and using (\ref{L4conn}), (\ref{L2conndrv}) and (\ref{ConnH}), we have
\begin{align*}
\|\Omega^2(\bA_b \Lb\sX)\|_{L^2_u L_\omega^2}&\les \l t\r^{-\frac{7}{4}+\delta}\Delta_0^\f12\|\Omega^{1+\le 1}\Lb\sX\|_{L_u^2 L_\omega^2}+\l t\r^{-\frac{9}{4}+\delta}\Delta_0.
\end{align*}
Using Proposition \ref{7.22.2.22}, (\ref{LbBA2}), (\ref{2.19.1.24}) and (\ref{2.13.3.24}), in view of (\ref{6.4.2.24}), it is direct to obtain
\begin{equation*}
\|\Omega\Lb \sX\|_{L^2_u L_\omega^2}\les \l t\r^{-1}(\log \l t\r)^{\f12(\M+7)}\Delta_0.
\end{equation*}
It follows by substituting the above two estimates into (\ref{3.16.1.24}) and using Gronwall's inequality 
\begin{align*}
\|\tir^2 \Omega^2 \Lb \sX\|_{L^2_u L_\omega^2}&\les\La_0+ \int_0^t\{\|\tir^2\Omega^2((L+\tr\chi) \Lb \sX)\|_{L_u^2 L_\omega^2}+\l t'\r^{-\frac{1}{4}+\delta}\Delta_0\}.
\end{align*}

Next, using (\ref{8.31.4.19}), we write
\begin{align}
\Omega^2((L+\tr\chi) \Lb \sX)&=\Omega^2\{(\f12\sX+\Xi_4+[L\Phi])\Lb \Xi_4+\ud \bA\sn\sX+\fB L\sX+\Lb[L\Phi]\tr\chi\nn\\
&+\Lb(|\chih|^2+\N(\Phi, \bp\Phi))\}.\label{3.15.5.24}
\end{align}
For the first term on the right-hand side, it is direct to derive
\begin{align*}
\|\Omega^2(\Lb \Xi_4 (\tr\chi+\fB))\|_{L^2_u L_\omega^2}&\les \l t\r^{-1}\|\Omega^2\Lb\bT\varrho, \Omega^2\Lb[L\Phi]\|_{L^2_u L_\omega^2}+\|\Omega\Lb\Xi_4\c \Omega(\bA_b+\fB)\|_{L_u^2 L_\omega^2}\\
&+\|\Lb \Xi_4\Omega^2(\bA_b+\fB)\|_{L^2_u L_\omega^2}, 
\end{align*}
where the first term on the right-hand side can be bounded by $\l t\r^{-2+\delta}\Delta_0$ by using (\ref{3.10.2.24}) and (\ref{2.29.5.24}). 
Using (\ref{L2conndrv}), (\ref{ConnH}), (\ref{3.11.3.21}), (\ref{LbBA2}) and (\ref{2.13.1.24}), we have for the two error terms 
\begin{align*}
\|\Omega\Lb\Xi_4\c \Omega(\bA_b+\fB)\|_{L_u^2 L_\omega^2}&\les \l t\r^{-2+2\delta}\Delta_0^2\\
\|\Lb \Xi_4 \Omega^2(\bA_b+\fB)\|_{L_u^2 L_\omega^2}&\les \l t\r^{-1}(\|\Omega^2\fB\|_{L_u^2 L_\omega^2}+\|\Omega^2\bA_b\|_{L^2_u L_\omega^2})+\log \l t\r\l t\r^{-1+\f12\delta}\Delta_0^2\l t\r^{-1+\delta}\\
&\les\l t\r^{-2+\delta}\Delta_0+\l t\r^{-2+2\delta}\Delta_0^2.  
\end{align*}
Hence we conclude from the above two sets of calculations that
\begin{align*}
\|\Omega^2(\Lb \Xi_4 (\tr\chi+\fB))\|_{L^2_u L_\omega^2}\les \l t\r^{-2+2\delta}\Delta_0. 
\end{align*}
Using (\ref{ConnH}), (\ref{3.6.2.21}), (\ref{zeh}) and (\ref{LbBA2}), we obtain
\begin{equation*}
\|\Omega^2(\ud\bA\sn\sX)\|_{L_u^2 L_\omega^2}\les \l t\r^{-3+2\delta}\Delta_0^2. 
\end{equation*}
Due to (\ref{6.22.1.21}) and (\ref{3.16.1.22}), $L\sX=O(\l t\r^{-2})$. Using  (\ref{8.7.2.22}), (\ref{2.29.5.24}) and Proposition \ref{7.15.5.22}, we bound
\begin{align*}
\|\Omega^2(\fB L\sX)\|_{L^2_u L_\omega^2}&\les\l t\r^{-3+2\delta}\Delta_0.  
\end{align*}
Similarly, using in addition (\ref{2.19.1.24}), we have
\begin{equation*}
\|\Omega^2(\Lb[L\Phi]\tr\chi)\|_{L^2_u L_\omega^2}\les \l t\r^{-3+2\delta}\Delta_0.
\end{equation*}
Assuming (\ref{3.8.1.24}) and using (\ref{L2BA2}) we have
\begin{align*}
\|\Omega^2\Lb\N(\Phi, \bp\Phi)\|_{L^2_\Sigma}\les \l t\r^{-\frac{3}{2}+\delta}\Delta_0.
\end{align*}
Assuming
\begin{align}\label{3.15.9.24}
\|\Omega^2(\sn_\Lb \chih \c \chih)\|_{L^2_u L_\omega^2}&\les \l t\r^{-\frac{13}{4}+4\delta}\Delta_0^2+\l t\r^{-\frac{7}{4}+\delta}\Delta_0\|\Omega^2(\mu-k_{\bN\bN}\tr\chi-\varpi)\|_{L^2_u L_\omega^2},
\end{align}
summarizing the above estimates leads to 
\begin{align*}
\|\tir^2 \Omega^2\Lb \sX\|_{L_u^2 L_\omega^2}&\les \La_0+\l t\r^{1+2\delta}\Delta_0+\int_0^t\l t'\r^{-\frac{7}{4}+\delta}\|\tir^2 \Omega^2(\mu-k_{\bN\bN}\tr\chi-\varpi)\|_{L^2_u L_\omega^2}.
\end{align*}
By the definition of $\mu$, we have by using (\ref{L4conn}), (\ref{L2BA2}) and (\ref{LbBA2})
\begin{align*}
\|\Omega^2(\mu-k_{\bN\bN}\tr\chi-\varpi)\|_{L_u^2 L_\omega^2}&\les \|\Omega^2(\Lb \sX, \Lb \Xi_4, \tr\chi^2, \tr\chi\fB, \varpi)\|_{L^2_u L_\omega^2}\\
&\les \|\Omega^2\Lb\sX\|_{L^2_u L_\omega^2}+\|\Omega^2\Lb^2\varrho, \Omega^2\Lb[L\Phi]\|_{L^2_u L_\omega^2}+\l t\r^{-2+2\delta}\Delta_0\\
&\les\|\Omega^2\Lb \sX\|_{L_u^2 L_\omega^2}+\l t\r^{-1+\delta}\Delta_0
\end{align*}
where we employed (\ref{2.29.5.24}), (\ref{3.10.2.24}), (\ref{6.5.1.24}) and (\ref{8.23.1.23}). 

Combining the above two estimates we conclude
\begin{equation}\label{3.15.10.24}
\|\Omega^2\Lb \sX, \Omega^2(\mu-k_{\bN\bN}\tr\chi-\varpi)\|_{L^2_u L_\omega^2}\les \l t\r^{-1+2\delta}\Delta_0.
\end{equation}
This gives the last case in (\ref{6.5.3.24}) and (\ref{8.12.2.22}) . It remains to prove (\ref{3.8.1.24}) and (\ref{3.15.9.24}).

Using (\ref{L4conn}), (\ref{3.6.2.21}), (\ref{7.22.1.22}) and (\ref{3.15.6.24}), we bound
\begin{align}\label{3.15.8.24}
\begin{split}
\|\Omega^2(\sn_\Lb \chih \c \chih)\|_{L^2_u L_\omega^2}&\les \|\sn_\Omega^2\sn_\Lb \chih\|_{L_u^2 L_\omega^2}\|\chih\|_{L_x^\infty}+\|\sn_\Omega \sn_\Lb \chih\c\sn_\Omega\chih\|_{L^2_u L_\omega^2}\\
&+\|\sn_\Lb \chih\c \sn_\Omega^2\chih\|_{L^2_u L_\omega^2}\\
&\les \l t\r^{-\frac{7}{4}+\delta}\Delta_0\|\sn_\Omega^{1+\le 1}\sn_\Lb \chih\|_{L_u^2 L_\omega^2}+\|\sn_\Lb \chih\|_{L_\omega^4}\|\sn_\Omega^2 \chih\|_{L_u^2 L_\omega^4}\\
&\les \l t\r^{-\frac{7}{4}+\delta}\Delta_0\|\sn_\Omega^2\sn_\Lb \chih\|_{L_u^2 L_\omega^2}+\l t\r^{-\frac{13}{4}+4\delta}\Delta_0^2.
\end{split}
\end{align}
Applying Proposition \ref{7.15.5.22}, we derive in view of (\ref{3chi}) and (\ref{5.21.3.23})
\begin{align}
\|\sn_\Omega^2 \sn_\Lb \chih\|_{L_u^2 L_\omega^2}&\les \l t\r^{-2+2\delta}\Delta_0^2+\|\sn_\Omega^2\sn \ze, \l t\r^{-1}\sn_\Omega^2\bA_g, \sn_\Omega^2\bar\bp\bA_{g,1}\|_{L^2_u L_\omega^2}\nn\\
&\les \l t\r^{-2+2\delta}\Delta_0^2+\l t\r^{-\frac{5}{2}+\delta}\Delta_0+\|\sn_\Omega^2\sn \ze\|_{L^2_u L_\omega^2}.\label{3.15.7.24}
\end{align}
Applying (\ref{9.14.10.22}) to $F=\ze$ and using (\ref{zeh}) imply
\begin{align*}
\|\tir\sn\sn_\Omega^2\ze\|_{L^2_u L_\omega^2}+\|\sn_\Omega^2 \ze\|_{L^2_u L_\omega^2}&\les \|\bb^{-\f12}\{\sn_\Omega^2 \D_1 \ze, \l t\r^{-1}\sn_\Omega^{1+\le 1}\ze,  \sn_\Omega^{\le 1}(\D_1\ze, \tir^{-1}\ze)\}\|_{L^2_\Sigma}\nn\\
&+\l t\r^{-\frac{7}{4}+\delta}\Delta_0^\f12 \|\tir \Big(\tir\sn \ze,  \ze\Big)\|_{L^2_u L_\omega^4}\\
&\les \|\tir \Omega^2\D_1\ze\|_{L^2_u L_\omega^2}+\l t\r^{-1+\delta}\Delta_0.
\end{align*}
Due to (\ref{6.5.2.24}), (\ref{8.28.2.23}), (\ref{1.30.1.24}), (\ref{9.19.2.22}) and (\ref{L2BA2}), we have
\begin{align*}
\|\Omega^2 \D_1 \ze\|_{L_u^2 L_\omega^2}&\les \|\Omega^2(\mu-k_{\bN\bN}\tr\chi-\varpi, \bA_{g}^2, \ud \bA^2)\|_{L^2_u L_\omega^2}+\|\Omega^2(\sn\bA_{g,1}, \fB \bA_{g,2})\|_{L^2_u L_\omega^2}\\
&\les\|\Omega^2(\mu-k_{\bN\bN}\tr\chi-\varpi)\|_{L^2_u L_\omega^2}+\l t\r^{-2+2\delta}\Delta_0. 
\end{align*}
Combining the above two estimates, we conclude
\begin{align*}
\|\tir\sn\sn_\Omega^2\ze\|_{L^2_u L_\omega^2}+\|\sn_\Omega^2 \ze\|_{L^2_u L_\omega^2}\les \|\tir\Omega^2(\mu-k_{\bN\bN}\tr\chi-\varpi)\|_{L^2_u L_\omega^2}+\l t\r^{-1+2\delta}\Delta_0. 
\end{align*}
Substituting the above estimate to (\ref{3.15.7.24}) with the help of (\ref{4.22.4.22}) yields
\begin{align*}
\|\sn_\Omega^2 \sn_\Lb \chih\|_{L_u^2 L_\omega^2}\les \l t\r^{-2+2\delta}\Delta_0+\|\Omega^2(\mu-k_{\bN\bN}\tr\chi-\varpi)\|_{L^2_u L_\omega^2}.
\end{align*}
Substituting the above estimate to (\ref{3.15.8.24}), we conclude (\ref{3.15.9.24}).

It remains to check (\ref{3.8.1.24}) by using (\ref{2.10.7.24}).
\begin{align*}
\Omega^2\big(\Lb\N(\Phi, \bp\Phi)\big)&=\Omega^2(\Lb^2\varrho[L\Phi]+\fB\Lb[L\Phi]+[\bar\bp\Phi]\sn_\Lb[\bar\bp\Phi]+\sn_\Lb \eh \c \eh)
\end{align*}
The last two terms on the right-hand side have better decay behavior than the first two.  
Using (\ref{1.29.2.22}), (\ref{3.10.2.24}), (\ref{2.13.3.24}), (\ref{2.29.5.24}) and Proposition \ref{7.15.5.22}, we derive  for the first term that
\begin{align*}
\|\Omega^2(\Lb^2\varrho [L\Phi])\|_{L^2_\Sigma}&\les \|\bb^{-2}\Omega^2[L\Phi]\|_{L_u^2 L_\omega^4}+\l t\r^{-2+2\delta}\Delta_0^2
\end{align*}
Using (\ref{3.6.2.21}), (\ref{6.22.1.21}), (\ref{2.29.5.24}), (\ref{2.19.1.24}) and (\ref{2.14.1.24}), we bound the second term on the right-hand side tht
\begin{align*}
\|\Omega^2(\fB\Lb[L\Phi])\|_{L^2_\Sigma}&\les \|\fB\Omega^2(\sD\varrho, L[L\Phi])\|_{L^2_\Sigma}+\l t\r^{-1}\|\Omega^2\fB\|_{L_u^2 L_\omega^2}+\|\Omega \fB\|_{L_\omega^\infty}\|\Omega\Lb[L\Phi]\|_{L^2_\Sigma}\\
&\les\l t\r^{-2} \log \l t\r(\La_0+\l t\r^\delta\Delta_0^\frac{5}{4})+\|\fB\Omega^2(\sD\varrho, L[L\Phi])\|_{L^2_\Sigma}.
\end{align*}
We omit the checking details for the other two terms, which can be done by using Proposition \ref{7.15.5.22}. Thus (\ref{3.8.1.24}) is proved.

Next we prove (\ref{8.9.2.22}).
Using (\ref{3.22.5.21}) and (\ref{cmu_2}), symbolically,
\begin{equation*}
\sn_L\sn\la=\chi\sn \la+(\bA_{g,1}\c \Omega+\sn\Omega)(\log c).
\end{equation*}
Applying (\ref{12.22.4.23}) and (\ref{4.22.4.22}) yields
\begin{equation*}
\|\sn_X^2\left(\bA_{g,1}\c\Omega\log c\right)\|_{L^2_\Sigma}\les \l t\r^{-2+2\delta}\Delta_0^\frac{3}{2}.
\end{equation*}
In view Proposition \ref{7.15.5.22} and Lemma \ref{3.17.2.22}, we deduce
\begin{align*}
\sn_X^2 \sn\Omega\log c&=O(\l t\r^{-\f12+\delta}\Delta_0)_{L^2_\Sigma}.
\end{align*}
Combining the above two estimates gives
\begin{equation}\label{8.9.5.22}
\|\sn_X^2(\bA_{g,1}\c \Omega+\sn\Omega)\log c\|_{L^2_\Sigma}\les \l t\r^{\delta-\f12}\Delta_0.
\end{equation}
Next applying (\ref{8.9.3.22+})  to $F=\sn \la$, using Proposition \ref{10.16.1.22}, we infer
\begin{align}
\|\bb^{-\f12}\sn_X^2(\bA \c \sn \la)\|_{L^2_\Sigma}\les \l t\r^{-\frac{5}{4}+2\delta}\Delta_0^\frac{3}{2}\label{3.1.1.24}
\end{align}
Finally, using Proposition \ref{10.16.1.22} again, we have
\begin{equation}\label{8.10.1.22}
\|\bb^{-\f12}\sn_X^2(\tir^{-1}\sn\la)\|_{L^2_\Sigma}\les\l t\r^{\delta-1+\f12(1-\vs^+(X^2))}\Delta_0.
\end{equation}
Combining the above three estimates,  in view of $\chi=\bA+\tir^{-1}$, we conclude
\begin{equation*}
\|\bb^{-\f12}\sn_X^2\left(\chi\sn \la\right)\|_{L^2_\Sigma}\les\l t\r^{\delta-1+\f12(1-\vs^+(X^2))}\Delta_0.
\end{equation*}
 Using the above estimate together with (\ref{8.9.5.22}) we conclude (\ref{8.9.2.22}).
\end{proof}

\begin{proof}[Proof of Proposition \ref{8.10.2.21}] We will rely on the $l=2$ case in Lemma \ref{3.23.1.23} for the proof.

\noindent{\bf$\bullet$ Estimates of $\sn_X^2 {}\rp{a}\bJ_L$ and $\sn_X^2{}\rp{a}\bJ_\Lb$.}
We first consider
 $\sn_X^2\sn^A {}\rp{a}\pih_{A L}$ in view of (\ref{xdjo}). Recall that
 \begin{align*}
 \sn_X^2 {}\rp{a}\ckk \J_L
 &=\sn_X^2(\ud \bA{}\rp{a}\pih_{AL})+\sn_X^2(\sdiv{}\rp{a}\pih_b)\\
& +c^{-1}(X\varrho+\sn_X)^2\circ(\sn \varrho+\sdiv)(c{}\rp{a}\pih^+_{AL})
 \end{align*}
 where $c\pioh^+_{AL}=\sn \la^a+\la \bA_{g,1}+\eta({}\rp{a}\Omega)$.
 
 Using (\ref{5.21.1.21}), (\ref{3.25.1.22}) and the estimate of $\ud\bA$ in Proposition \ref{7.15.5.22}, we bound the first term on the right-hand side by
\begin{equation}\label{2.28.4.24}
X^2(\ud \bA\c \pioh_{A L})=O(\l t\r^{-\frac{3}{2}+2\delta}\Delta_0^\frac{3}{2})_{L_u^2 L_\omega^2}.
\end{equation}
For simplicity, for the last line, we only present the treatment of the higher order terms with details, since $X\varrho$ has much better decay. Thus the last line is shortened as 
\begin{align*}
\sn_X^2\{(\sn \varrho+\sdiv)(c{}\rp{a}\pih^+_{AL})\}&=c^{-1}(\sn_X^2(\sn\varrho\c c\pioh^+_{A L})+\sn_X^2(\sdiv {}\rp{a}(c\pih^+_{A L})))\\
&=O(\l t\r^{-\frac{5}{2}+2\delta}\Delta_0^2)_{L_u^2 L_\omega^2}+c^{-1}\sn_X^2\sdiv(c {}\rp{a}\pih^+_{A L}),
\end{align*}   
where the estimate is obtained by using (\ref{3.6.2.21}), (\ref{3.11.3.21}), (\ref{5.21.1.21}) and (\ref{3.25.1.22}). 
It remains to treat the last term in the above. In view of (\ref{12.17.1.23}), we write schematically
\begin{align*}
X^2\sD \la&=c^2(X\varrho+X)^2\{(\Omega\tr\thetac+\bA\la \thetac)\}.
\end{align*}

It follows by using Proposition \ref{10.16.1.22} and Proposition \ref{7.15.5.22}
\begin{align*}
\|\sn_X^2(\bA \la\thetac)\|_{L^2_u L_\omega^2}\les \l t\r^{-\frac{5}{2}+2\delta}\Delta_0^\frac{3}{2}. 
\end{align*}
Hence we summarize the above calculations below
\begin{align*}
\sn_X^2{}\rp{a}\ckk\J_L&=\sn_X^2(\sdiv{}\rp{a}\pih_b)+c\sn_X^2\big(\Omega\tr\thetac, \sdiv (\eta({}\rp{a}\Omega))\big)+O(\l t\r^{-\frac{3}{2}+2\delta}\Delta_0^\frac{3}{2})_{L_u^2 L_\omega^2}.
\end{align*}
Moreover, recall from (\ref{8.2.2.23}) that
\begin{align*}
\sn_X^2\big(\sdiv\eta({}\rp{a}\Omega)\big)
&=\tir^{-1}\sn_X^{\le 2}[\Omega v]+X^2(\Omega[L\Phi])+O(\l t\r^{-\frac{5}{2}+2\delta}\Delta_0^\frac{3}{2})_{L_u^2 L_\omega^2}.
\end{align*}

Using the second identity in Lemma \ref{3.23.1.23} and (\ref{8.29.9.21}), we summarize the above calculations as
\begin{align}\label{2.28.3.24}
\begin{split}
\sn_X^2{}\rp{a}\ckk\J_L&=X^2(\ckr^{-1} {}\rp{a}v^\sharp+\bA_b {}\rp{a}v^\sharp+\sn\varrho\c {}\rp{a} v^\star)+c(X^2 \Omega\tr\thetac+X^2(\Omega[L\Phi]))\\
&+O(\l t\r^{-\frac{3}{2}+2\delta}\Delta_0^\frac{3}{2}+\l t\r^{-2}\log \l t\r^{\f12\M}\La_0)_{L_u^2 L_\omega^2}.
\end{split}
\end{align}
For the first term on the right-hand side,  we claim
\begin{align}
&\|\bb^{-\f12}\sn_X^2(\bA {}\rp{a}v^\sharp)\|_{L^2_\Sigma}\les \l t\r^{-\frac{7}{4}+2\delta}\Delta_0^\frac{3}{2}\label{8.10.7.22}\\
&\|\bb^{-\f12}\sn_X^2(\bA_{g,1} {}\rp{a}v^*)\|_{L^2_\Sigma}\les \l t\r^{-2+2\delta}\Delta_0^\frac{3}{2}.\label{8.10.8.22}
\end{align}
Using (\ref{3.28.3.24}), Proposition \ref{7.15.5.22} and (\ref{6.24.1.21}) we derive
\begin{align*}
&\|\bb^{-\f12}\sn_X^2(\bA {}\rp{a}v^\sharp)\|_{L^2_\Sigma}\\
&\les \|\bb^{-\f12}\sn_X^2\bA {}\rp{a}v^\sharp\|_{L^2_\Sigma}+\|\bb^{-\f12}\bA\c \sn_X^2 {}\rp{a}v^\sharp\|_{L^2_\Sigma}+\|\bb^{-\f12}\sn_X \bA \sn_X{}\rp{a}v^\sharp\|_{L^2_\Sigma}\\
&\les \l t\r^{-\frac{7}{4}+2\delta}\Delta_0^\frac{3}{2}
\end{align*}
 Thus (\ref{8.10.7.22}) is proved. 
(\ref{8.10.8.22}) can be obtained similarly. Hence the first term on the right-hand side of (\ref{2.28.3.24}) is bounded by 
\begin{align*}
X^2(\ckr^{-1}{}\rp{a}v^\sharp)+O(\l t\r^{-2+2\delta}\Delta_0^\frac{3}{2})_{L_u^2 L_\omega^2}.
\end{align*}

Note that by using Proposition \ref{2.28.2.24},  Proposition \ref{10.16.1.22} and Proposition \ref{7.15.5.22}, we have
\begin{align}
&\|\sn_X^2\sn_L(\la \c \bA)\|_{L^2_u L_\omega^2}\les \l t\r^{-\frac{5}{2}+2\delta}\Delta_0^\frac{3}{2}.\label{8.11.2.22}
\end{align}
In view of (\ref{xdjo}), we further derive by using (\ref{8.11.2.22}), Proposition \ref{1steng}, Proposition \ref{8.29.8.21} and Lemma \ref{10.10.3.23},
\begin{align*}
X^2L{}\rp{a}\ss&=X^2L\Omega \varrho+c^{-1}X^2 L(\tir^{-1}\la+\bA\la)\\
&=X^2 L\Omega\varrho+c^{-1}X^2(\tir^{-2}\la+\tir^{-1}\Omega c)+O(\l t\r^{-\frac{5}{2}+2\delta}\Delta_0^\frac{3}{2})_{L^2_u L_\omega^2}\\
&=LX^2\Omega\varrho+c^{-1}X^2(\tir^{-2}\la)+ O(\l t\r^{-1}(\Delta_0^\frac{5}{4}\l t\r^{\delta}+\La_0)\log \l t\r^{\f12\M})_{L^2_\Sigma}
\end{align*}
where in the first identity, we dropped lower order quadratic error terms. 
Combining the above two calculations with (\ref{2.28.3.24}) yields
\begin{align*}
\sn_{X_2}\sn_{X_1}{}\rp{a}\bJ_L&=\sn^2_X(\ckr^{-1}{}\rp{a}v^\sharp)+ X^2 \Omega\tr\thetac+X^2(\tir^{-2}\la)\\
&+[LX^2\Omega\Phi]+O(\l t\r^{-\frac{3}{2}+2\delta}\Delta_0^\frac{3}{2}+\l t\r^{-2}\La_0\log \l t\r^{\f12\M})_{L_u^2 L_\omega^2}.
\end{align*}
where we dropped the coefficients $c^{-1}$ or $c$ without affecting the subsequent result, as stated in (\ref{8.8.6.21}).

Next we consider the estimate of $\sn_X^2{}\rp{a}\bJ_\Lb$. In view of (\ref{xdjo}), regarding $\pioh_{A\bN}=\pioh_{A L}+\ud\bA\la$,  in comparison with the terms in $\sn_X^2{}\rp{a}\ckk \J_L$, we need to bound the additional term 
\begin{align*}
\sn_X^2\{(\ud\bA+\sdiv)(\la \ud \bA)\}&=\sn_X^2 \sdiv(\la \ud \bA)+\sn^2_X(\la \ud\bA^2).
\end{align*} 
Again for simplicity, we present the control of the higher order term here. 
Using Proposition \ref{10.16.1.22} and Proposition \ref{7.15.5.22}, we derive
\begin{align*}
\|\sn_X^2(\la \sn\ud \bA)\|_{L^2_u L_\omega^2}&\les \l  t\r^{-1-\vs^+(X^2)+2\delta}\Delta_0^2\\
\displaybreak[0]
\|\sn_X^2(\sn\la \c \ud \bA)\|_{L^2_u L_\omega^2}&\les \l t\r^{-\frac{7}{4}+2\delta}\Delta_0^2\\
\|\sn_X^2( \la \ud \bA^2)\|_{L^2_u L_\omega^2}&\les \l t\r^{-2+3\delta}\Delta_0^3. 
\end{align*}
Hence using the bound of the above terms, (\ref{2.28.4.24}) and (\ref{3.28.3.24}), symbolically, we infer 
\begin{align*}
&\sn_X^l((\bA_{g,1}+\sn_A){}\rp{a}\pih_{AL})-\sn_X^l((\ud\bA+2\sn_A){}\rp{a}\pih_{A\bN})\\
&=\sn^2_X(\ckr^{-1}{}\rp{a}v^\sharp)+ X^2(\Omega\tr\thetac)+[LX^2\Omega\Phi]+O(\l t\r^{-1-\f12\vs^+(X^2)+2\delta}\Delta_0^\frac{3}{2}+\l t\r^{-2}\La_0\log \l t\r^{\f12\M+1})_{L_u^2 L_\omega^2}\\
&= O(\l t\r^{2\delta-1-\f12\vs^+(X^2)}\Delta_0)_{L_u^2 L_\omega^2}+[LX^2\Omega\Phi].
\end{align*}

Next we show
\begin{align}\label{8.10.9.22}
\|\bb^{-\f12}\sn_X^2((L+k_{\bN\bN}+\tr\chi){}\rp{a}\pih_{\Lb\Lb})\|_{L^2_\Sigma}\les \l t\r^{\delta}\Delta_0(\log \l t\r)^2.
\end{align}
Recall that ${}\rp{a}\pih_{\Lb \Lb}=\ud \bA(\Omega)$.
 Using (\ref{2.29.6.24}), (\ref{8.5.1.22+}), (\ref{8.9.4.22}), Lemma \ref{3.17.2.22} and Proposition \ref{7.15.5.22}, we derive
\begin{align*}
&\|\bb^{-\f12}\sn_X^2 (L+\tir^{-1})\pioh_{\Lb \Lb}\|_{L^2_\Sigma}\les \l t\r^{\delta}\Delta_0;\\
&\|\bb^{-\f12}\sn_X^2(k_{\bN\bN} \pioh_{\Lb \Lb})\|_{L^2_\Sigma}\les \l t\r^\delta(\log\l t\r)^2\Delta_0;\\
&\|\bb^{-\f12}\sn_X^2\big(\bA\pioh_{\Lb\Lb}\big)\|_{L^2_\Sigma}\les  \l t\r^{-\f12+2\delta}\Delta_0^\frac{3}{2}.
\end{align*}
Combining the above estimates implies (\ref{8.10.9.22}).

Hence, 
\begin{align}
\sn_X^2{}\rp{a}\ckk\J_\Lb&= O(\l t\r^{2\delta-1}\Delta_0)_{L_u^2 L_\omega^2}+[LX^2\Omega\Phi].\label{8.11.5.22}
\end{align}

Finally, we consider the term $X^2\Lb {}\rp{a}\ss$. In view of (\ref{xdjo}), we prove the following estimates
\begin{align}
&\|X^2 \big(\Lb\varrho ([L\Phi]+\tir^{-1}+\bA_b)\la\big)\|_{L^2_u L_\omega^2}\les \l t\r^{-2+3\delta}\Delta_0^\frac{3}{2}\label{8.11.1.22}\\
&\|\sn_X^2\Lb \big(([L\Phi]+\tir^{-1}+\bA_b)\la\big)\|_{L^2_u L_\omega^2}\les \Delta_0 \l t\r^{-1+3\delta}.\label{8.11.6.22}
\end{align}
Combining the above two estimates in view of (\ref{xdjo}) and (\ref{LbBA2}) yields
\begin{align*}
X^2\Lb{}\rp{a}\ss=O(\l t\r^{3\delta-1}\Delta_0)_{L_u^2 L_\omega^2}.
\end{align*}
In view of (\ref{8.11.5.22}), this implies
\begin{align*}
\sn_{X_2}\sn_{X_1}{}\rp{a}\bJ_\Lb&= O(\l t\r^{3\delta-1}\Delta_0)_{L_u^2 L_\omega^2}+[LX^2\Omega\Phi].
\end{align*}
Thus (\ref{8.8.6.21}) is proved.

Now we prove (\ref{8.11.1.22}) and (\ref{8.11.6.22}). Note
\begin{align}
X^l \Lb \bA_b&=O(\l t\r^{-2})+O(\l t\r^{-1-\vs^+(X^l)+2\delta}\Delta_0)_{L^2_u L_\omega^2}\nn\\
 X^l\Lb [L\Phi]&=O(\l t\r^{-1}\fB)+O(\l t\r^{-2+\delta}\Delta_0)_{L_u^2 L_\omega^2}, \quad l=1,2\label{3.1.4.24}\\
  \Lb (\bA_b,[L\Phi])&=O(\l t\r^{-2})+O(\l t\r^{-2+\delta}\log \l t\r\Delta_0)_{L_\omega^4}.\nn
 \end{align}
Indeed, the first estimate can be seen from Proposition \ref{7.22.2.22} and (\ref{8.12.2.22}). The second estimate can be derived by using (\ref{2.29.5.24}) and Proposition \ref{7.15.5.22}. The last estimates are from Proposition \ref{7.22.2.22} and Lemma \ref{5.13.11.21}.
Using the above estimates, for $F=\bA_b, [L\Phi]$, we derive by using Proposition \ref{7.22.2.22} , Proposition \ref{7.15.5.22} and Proposition \ref{10.16.1.22} that
\begin{align*}
\|X^2\Lb(F\c \la)\|_{L^2_u L_\omega^2}&\les\|X^2\Lb F\|_{L_u^2 L_\omega^2}\|\la\|_{L^\infty_x}+\|\Lb F\c X^2\la\|_{L_u^2 L_\omega^2}+\|X^2 F\|_{L^2_u L_\omega^2}\|\Lb\la\|_{L^\infty_x}\\
&+\|X F\c X\Lb\la, X\Lb F\c X\la\|_{L^2_u L_\omega^2}+\|X^2\Lb \la\|_{L^2_u L_\omega^2}\|F\|_{L^\infty_x}\\
&\les \l t\r^{-1+3\delta}\Delta_0.
\end{align*}
Applying (\ref{8.9.4.22}) to $F=\bA\la$ implies
\begin{align*}
\|X^2(\Lb \varrho\bA\la)\|_{L^2_u L_\omega^2}&\les \l t\r^{-1}\|\sn_X^{\le 2}\big( \bA\la\big)\|_{L^2_u L_\omega^2}+\Delta_0\l t\r^{-1}(\log \l t\r)^2\|\bA\la\|_{L^\infty_x}\\
&+\Delta_0\l t\r^{-1}\log \l t\r\|X(\bA\la)\|_{L_u^2L_\omega^4}\les \l t\r^{-\frac{5}{2}+3\delta}\Delta_0^\frac{3}{2},
\end{align*}
where, to obtain the last estimate, we applied (\ref{dpio}), Proposition \ref{7.15.5.22} and Proposition \ref{10.16.1.22}. 
Moreover, we can check by using (\ref{8.9.4.22}), (\ref{dpio}) and Proposition \ref{10.16.1.22},
\begin{align*}
\|X^2(\la \fB)\|_{L_u^2 L_\omega^2}\les \l t\r^{-1+2\delta}\Delta_0, \|X^2\Lb (\la \tir^{-1})\|_{L_u^2 L_\omega^2}\les\l t\r^{-1+\delta}\Delta_0 \end{align*}
Hence we obtained (\ref{8.11.1.22}) and (\ref{8.11.6.22}).

\noindent{\bf$\bullet$ Estimate of $\sn_X^2 {}\rp{a}\bJ_B$.}
Let us first consider the estimate of $\sn_X^2 {}\rp{a}\bJ_B$ in (\ref{8.8.6.21}).

 For this purpose, we will prove
\begin{align}
&\sn_X^2 {}\rp{a}\wt\eth_B=\sn_X^2 \sn_\bT {}\rp{a}\pih_b+O(\fB)+O(\l t\r^{-1+2\delta}\Delta_0)_{L_u^2 L_\omega^2} 
\label{8.8.2.22}\\ 
&\|\sn_X^2((k_{\bN\bN}+\tir^{-1}+\bA)(\pioh_{AB}, \pioh_{A\bN},\pioh^+_{A\bT})),\nn\\
&\qquad\qquad\sn_X^2\big((\bA_g+\fB){}\rp{a}\pih_b\big)\|_{L^2_u L_\omega^2}\les \l t\r^{2\delta-1}\Delta_0.\label{8.10.3.22}
\end{align}
To see (\ref{8.8.2.22}), we rewrite the formula of ${}\rp{a}\wt\eth$ in Proposition \ref{error_terms} with the help of (\ref{dchi}) as
\begin{align}\label{3.1.2.24}
\begin{split}
{}\rp{a}\wt\eth&=\sn{}\rp{a}\Omega \log c+c^{-1}\la_a\Big(\sn\tr\chi+\bR_{ABAL}+\sdiv\eh+\bA_{g,1}(\bA+\tir^{-1})\Big)\\
&+c^{-1}\sn\la\bA+c^{-1}\tir^{-1}\sn\la+\ud\bA^2 \Omega\\
&+\sn_L(c^{-1}\la\ud\bA)+\sn_\bT {}\rp{a}\pih_b+\sn_\bT(c^{-1}\la \bA_{g,1}+c^{-1}\eta(\Omega)).
\end{split}
\end{align} 
 We will drop cubic error terms in the sequel since they decay better. 

Using (\ref{1.30.1.24}), Lemma \ref{3.17.2.22}, Proposition \ref{7.15.5.22}, (\ref{8.7.1.22}), Proposition \ref{10.16.1.22} and Lemma \ref{5.13.11.21} (3), we derive
\begin{align*}
&\|\sn_X^2(\sn(\bA_{g,1}\c\Omega))\|_{L_u^2 L_\omega^2}\les \l t\r^{-\frac{3}{2}+\delta}\Delta_0,\quad \|\sn_X^2(c^{-1}\la \sn\tr\chi)\|_{L_u^2 L_\omega^2}\les\l t\r^{-2+2\delta}\Delta_0^2\\
&\|\sn_X^2\big(c^{-1}\la (\bR_{ABAL}, \sn\bA_{g,1}+\bA_{g,1}(\bA+\tir^{-1}))\big)\|_{L_u^2 L_\omega^2}\les \l t\r^{-\frac{5}{2}+2\delta}\Delta_0^2\\
&\|\sn_X^2(\ud \bA^2\Omega)\|_{L^2_u L_\omega^2}\les \l t\r^{-1+2\delta}\Delta_0^2. 
\end{align*}
Using (\ref{3.1.1.24}) and (\ref{8.10.1.22}), combining the above estimates, we can bound the terms of $\sn_X^2{}\rp{a}\wt\eth$ contributed from the first two lines in (\ref{3.1.2.24}) by 
$
O(\l t\r^{-1+2\delta}\Delta_0)_{L_u^2 L_\omega^2}.
$
 
To bound the remaining terms, we need to show with $Y=L,\Lb$
\begin{align}
&\|\sn_X^2\sn_L(c^{-1}\la\c \ud \bA)\|_{L^2_u L_\omega^2}\les \l t\r^{-2+2\delta}\Delta_0^2,\label{3.1.3.24}\\
&\|\bb^{-\f12}\sn_X^2 \sn_Y(\la\c \bA_{g,1})\|_{L^2_\Sigma}\les \l t\r^{-1-\f12\max(\vs(Y),0)+2\delta}\Delta_0^2,\label{8.9.7.22}\\
&\sn_X^2\sn_Y(c^{-1}\eta(\Omega))=O(\fB)+O(\l t\r^{-1+\delta}\Delta_0)_{L_u^2 L_\omega^2}\label{3.1.5.24}.
\end{align}
(\ref{3.1.3.24}) can be similarly obtained as (\ref{8.11.2.22}) by using (\ref{8.7.2.22}), (\ref{3.16.1.22}), Proposition \ref{10.16.1.22} and Proposition \ref{7.15.5.22}. (\ref{8.9.7.22}) can be obtained by using Proposition \ref{7.15.5.22} and Proposition \ref{10.16.1.22}. Noting that $\tr\eta=[L\Phi]$, (\ref{3.1.5.24}) follows by using Proposition \ref{7.15.5.22}, (\ref{3.1.4.24}) and Lemma \ref{3.17.2.22}. Summarizing the above estimates, we controlled the terms contributed by the last line of (\ref{3.1.2.24}) by $$\sn_X^2\sn_\bT \pioh_b+O(\fB)+O(\l t\r^{-1+2\delta}\Delta_0)_{L_u^2 L_\omega^2}.$$ 
Hence we conclude (\ref{8.8.2.22}).

Next, applying (\ref{8.9.3.22+}) to $F=\pioh$, also using (\ref{5.21.1.21}) and (\ref{3.25.1.22}), we deduce
\begin{align*}
&\|\sn_X^2(\bA \c \pioh)\|_{L^2_u L_\omega^2}\les \l t\r^{-\frac{3}{2}+2\delta}\Delta_0^\frac{3}{2}.
\end{align*}
Using (\ref{8.9.4.22}), (\ref{5.21.1.21}), (\ref{3.25.1.22}) and the estimates of $\ud \bA=\bA_{g,1}+\ze$ in Proposition \ref{7.15.5.22}, we have
\begin{align}\label{3.2.1.24}
\|\sn_X^2\big((\fB, \tir^{-1})\pioh\big)\|_{L_u^2 L_\omega^2}\les \l t\r^{-1+2\delta}\Delta_0; \|\sn_X^2\big((\fB, \tir^{-1})\ud\bA\big)\|_{L_u^2 L_\omega^2}\les \l t\r^{-2+\delta}(\log \l t\r)^2\Delta_0
\end{align}
Hence (\ref{8.10.3.22}) is proved. 

Next, we use (\ref{6.28.6.21}) to derive
\begin{align*}
\sn_X^2 \sn_\bT {}\rp{a}\pih_b&=\sn_X^2(\bT v\sta{a}\wedge\Pi)+\sn_X^2(\ud\bA v^\sharp+[\bT\Phi]v^*).
\end{align*}
For the two error terms in the above, we derive by using Proposition \ref{1steng}, (\ref{7.25.2.22}), (\ref{8.9.4.22}) and (\ref{3.28.3.24}) that
\begin{align*}
\sn_X^2([\bT\Phi]v^*)=O(\l t\r^{-2+\delta}(\Delta_0^\frac{5}{4}+\La_0))_{L_u^2 L_\omega^2} 
\end{align*}
and by using (\ref{3.28.3.24}) and Proposition \ref{7.15.5.22} the following estimate
 \begin{align*}
 \sn_X^2(\ud\bA v^\sharp)=O(\l t\r^{-2+2\delta}\Delta_0^\frac{3}{2})_{L^2_u L_\omega^2}.
 \end{align*}
To treat the leading term, we apply (\ref{10.17.1.23}) to $\bT v^i\tensor{\ud\ep}{^a_i_j}$  and then to $X_1 \bT v^i \tensor{\ud\ep}{^a_i_j}$ to obtain
\begin{align*}
\sn_X^2 (\bT v\sta{a}\wedge\Pi)
&=X^2\bT v\sta{a}\wedge \Pi+[X_2\Phi]X_1\bT v\sta{a}\wedge \Pi+\big(\tir\ud \bA+\thetac(X_2)\big)X_1 \bT v\sta{a}\wedge \bN\\
&+\sn_{X_2}\Big([X_1\Phi]\bT v\sta{a}\wedge\Pi+(\tir\ud \bA+\thetac(X_1))\c \bT v\sta{a}\wedge \bN\Big).
\end{align*} 
Using (\ref{3.6.2.21}), (\ref{8.23.1.23'}) and (\ref{8.29.9.21}), the first line in the above is bounded by 
$$O(\fB)+O(\l t\r^{-1+2\delta}\Delta_0)_{L_u^2 L_\omega^2}.$$

Replacing ${}\rp{a} v^*$ in (\ref{6.24.2.21}) and (\ref{6.28.6.21}) by ${}\rp{a}(\bT v)^*=\bT v\sta{a}\wedge \Pi$, and replacing ${}\rp{a}v^\sharp$ in (\ref{7.25.1.22}) by ${}\rp{a}(\bT v)^\sharp=\bT v\sta{a}\wedge\bN$, also using Proposition \ref{7.15.5.22} and (\ref{8.23.1.23}), we obtain
\begin{equation*}
\sn_X ((\bT v)^*), \sn_X( (\bT v)^\sharp)=O(\fB)+O(\l t\r^{-1+\delta}\Delta_0)_{L_u^2 L_\omega^2}.
\end{equation*}
Thus, using Proposition \ref{7.15.5.22}, we have
\begin{align*}
\sn_{X_2}\big([X_1\Phi]\bT v \sta{a}\wedge \Pi+(\tir\ud \bA+\thetac(X_1))\bT v\sta{a}\wedge \bN\big)=O(\fB)+O(\l t\r^{-1+2\delta}\Delta_0)_{L_u^2 L_\omega^2}
\end{align*}
Therefore
\begin{align*}
\sn_X^2 (\bT v\sta{a}\wedge\Pi)=O(\fB)+O(\l t\r^{-1+2\delta}\Delta_0)_{L_u^2 L_\omega^2}.
\end{align*}
We conclude 
\begin{equation*}
\sn_X^2 \sn_\bT {}\rp{a}\pih_b=O(\fB)+O(\l t\r^{-1+2\delta}\Delta_0)_{L_u^2 L_\omega^2}.
\end{equation*}
Substituting the above estimate to (\ref{8.8.2.22}) and also using (\ref{8.10.3.22}), in view of (\ref{xdjo}), we have
\begin{equation*}
\sn_X^2 {}\rp{a}\bJ_B=O(\fB)+O(\Delta_0\l t\r^{-1+2\delta})_{L^2_u L_\omega^2}
\end{equation*}
as desired in (\ref{8.8.6.21}).

\noindent{\bf $\bullet$ Estimates of $\sn_X^2\bJ[S]$.} We first prove (\ref{8.13.2.22}).
Recall from (\ref{xdjs}) with $l=2$ that
\begin{align*}
\sn_X^2 \bJ[S]_\Lb&=\sn_X^2\left((L+\tr\chi-2k_{\bN\bN})(\tir k_{\bN\bN}+\Lb \tir)\right)\\
&+\sn_X^2\Big(\tir\sn \bA_{g,1}+\tir\varpi+\tir k_{\bN\bN}\tr\chi+\tir\tr\chi\mho+\tir (\ud\bA^2+\bA\c \bA)\Big).
\end{align*}
Using Proposition \ref{7.15.5.22}, (\ref{2.27.1.24}) and (\ref{1.30.1.24}), we infer
\begin{align*}
\|\sn_X^2\Big(\tir(\sn \bA_{g,1}+ \ud\bA^2+\bA\c\bA)\Big)\|_{L^2_u L_\omega^2}\les \l t\r^{-1+2\delta}\Delta_0.
\end{align*}
For the first term in $\sn_X^2 \bJ[S]_\Lb$,  it is straightforward to derive
\begin{align*}
\sn_X^2&\left((L+\tr\chi-2k_{\bN\bN})(\tir k_{\bN\bN}+\Lb \tir)\right)\\
&=\sn_X^2\left(L(\tir k_{\bN\bN}+\Lb \tir)\right)+\sn_X^2\left((\tr\chi-2k_{\bN\bN})(\tir k_{\bN\bN}+\Lb \tir)\right).
\end{align*}
Using (\ref{6.30.2.19}) and $k_{\bN\bN}=\Lb\varrho+[L\Phi]$, we derive
\begin{align*}
\sn_X^2 L(\tir k_{\bN\bN})&=X^2\big(\tir\{-\Box_\bg \varrho+\sD \varrho-\hb L \varrho+2\zb\c \sn \varrho+\bA_b \Lb\varrho+\fB^2\}\big)+X^2L[S\Phi]\\
&=X^2S^{\le 1}[L\Phi]+X^2(\tir\sD\varrho+\tir \fB^2)+O(\l t\r^{-\frac{3}{2}+2\delta}\Delta_0)_{L_u^2 L_\omega^2}
\end{align*}
where for the last line, we employed (\ref{9.19.2.22}), (\ref{8.26.4.21}), (\ref{8.30.3.21}) and (\ref{12.22.4.23}).

In view of (\ref{7.16.3.22}) and using $\mho=\Lb \log \tir-\f12 \tr\chib=\Lb\log\tir+\tir^{-1}+\fB+\bA$, we write
\begin{align*}
X^2L \Lb \tir&=X^2\left(k_{\bN\bN}(\Lb \tir -1)\right)=X^2(\tir \fB\c (\mho, \bA))+X^2(\fB, \tir \fB^2).
\end{align*}
Moreover, by our convention, we directly write
\begin{align*}
X^2\left((\tr\chi-2k_{\bN\bN})(\tir k_{\bN\bN}+\Lb \tir)\right)=X^2(\delta, \tir \delta^2,\Lb \log \tir, \delta \Lb \tir, \tir\delta\bA_b,\bA_b \Lb \tir).
\end{align*}
Summarizing the above calculations, using the definition of $\mho$ again,  we deduce
\begin{align*}
\sn_X^2&\left((L+\tr\chi-2k_{\bN\bN})(\tir k_{\bN\bN}+\Lb \tir)\right)\\
&=X^2(\fB, \tir \fB^2, \tir^{-1}, \tir \fB(\bA+\mho))+X^2(\mho+\bA)+X^2\left(\bA_b(\tir \fB, \tir \Lb\log \tir)\right)\\
&+X^2\left(S^{\le 1}[L\Phi]+\tir \sD \varrho)\right)+O(\l t\r^{-\frac{3}{2}+2\delta}\Delta_0)_{L_u^2 L_\omega^2}.
\end{align*}
Due to (\ref{7.29.1.22}) and (\ref{8.23.2.23})
\begin{equation*}
X^2 (\fB, \tir \fB^2)= O(\fB)+O(\l t\r^{-1+\delta}\Delta_0)_{L_u^2 L_\omega^2}.
\end{equation*}
Applying (\ref{8.9.4.22}) to $F=\tir(\bA+\mho)$, also using Proposition \ref{7.22.2.22}, Proposition \ref{7.15.5.22} and (\ref{2.27.1.24}), we derive
\begin{align}\label{3.27.1.23}
X^2\big( \tir(\bA+\mho)(\fB, \tir^{-1}, \bA)\big)=O(\l t\r^{-1})+ O(\l t\r^{-1+2\delta}\Delta_0)_{L_u^2 L_\omega^2}.
\end{align}
Due to $\varpi=(\bA+\frac{1}{\tir}+\fB)\fB$, $X^2(\tir \varpi)$ has been estimated in the above. Finally, using (\ref{L2BA2}) we obtain
\begin{align*}
\sn_X^2&\left((L+\tr\chi-2k_{\bN\bN})(\tir k_{\bN\bN}+\Lb \tir)\right)=O(\l t\r^{-1})+O\left(\l t\r^{-1+2\delta}\Delta_0\right)_{L_u^2 L_\omega^2}.
\end{align*}
Hence we proved (\ref{8.13.2.22}).

We recall the definition of $\sn_X^2\bJ[S]_B$ from (\ref{xdjs}). Using Proposition \ref{7.15.5.22}, (\ref{2.29.6.24}), (\ref{8.7.2.22}) and (\ref{3.2.1.24}), we bound 
\begin{align*}
\sn_X^2\bJ[S]_B&=\sn_X^l\eth[S]_B+\sn_X^2\big(\tir(\fB+\bA)\ud\bA\big)+\sn_X^2(\sn_S^{\le 1} \ud\bA)\\
&=\sn_X^l\eth[S]_B+O(\l t\r^{-1+\delta}\Delta_0(\log \l t\r)^2)_{L_u^2 L_\omega^2}.
\end{align*}
Due to (\ref{8.28.2.23}), Proposition \ref{1steng} and Proposition \ref{8.29.8.21}, it holds that 
\begin{equation*}
\|\sn_X^2(\chi\c \zb)\|_{L^2_\Sigma}\les \l t\r^{-2}(\log \l t\r)^{\f12\M}(\Delta_0^\frac{5}{4}+\La_0).
\end{equation*}
Using the above estimate, Corollary \ref{2.29.1.24}, Lemma \ref{5.13.11.21} (4), (\ref{L2BA2}) and (\ref{LbBA2}) 
\begin{align*}
\sn_X^2 \eth[S]_B&=\sn_X^2\big(\tir(\bR_{B4BA}+k_{A\bN}\c \chi+\sn \bA_b+\sn k_{\bN\bN})\big)\\
&=O(\l t\r^{-1+\delta})_{L_u^2 L_\omega^2}
\end{align*}
Hence (\ref{8.13.1.22}) is proved. 

Recall the definition of $X^2L{}\rp{S}\ss$ from (\ref{xdjs}). We derive by using (\ref{2.27.1.24}) and Proposition \ref{7.15.5.22} that
\begin{align*}
X^2L{}\rp{S}\ss&=X^2\{ \tir(\bA^2+(\tir^{-1}, \fB, \bA) \bAn+\sD\varrho+LL\varrho+\tir^{-1}\bA_b)\}\\
&=X^2(\tir(\sD\varrho+LL\varrho))+X^2\bA+O\Big((\l t\r^{-1}\log \l t\r)^{\frac{\M}{2}} (\La_0+\l t\r^\delta\Delta_0^\frac{5}{4})\Big)_{L^2_\Sigma}
 \end{align*}
 where we applied the following estimate
 \begin{align*}
 X^2(\fB\bAn)=O( (\l t\r^{-2}\log \l t\r)^{\frac{\M}{2}} (\La_0+\l t\r^\delta\Delta_0^\frac{5}{4}))_{L^2_\Sigma}
 \end{align*}
 which is obtained by using (\ref{8.9.4.22}) and Proposition \ref{9.8.6.22}.
Also by using Proposition \ref{8.29.8.21} and Proposition \ref{7.16.1.21}, we obtain
\begin{align*}
\|X^2L{}\rp{S}\ss\|_{L^2_\Sigma}\les &\l t\r^{-\f12}\{W_1[X_2 X_1 S\varrho]^\f12(t)+ \sum_{\vs(V)=\vs^-(X^2)}W_1[V\Omega^2\varrho]^\f12(t)\}\\
&+\l t\r^{-1}(\log \l t\r)^\frac{\M}{2}(\La_0+\l t\r^\delta\Delta_0^\frac{5}{4}).
\end{align*}
This is (\ref{8.14.1.22}).
\end{proof}
\subsection{Top order commutator estimates}

\begin{proposition}
With $X^2=X_3 X_2$ and $X_2, X_3\in \{\Omega, S\}$, there hold with $\ell=0,1$ 
\begin{align}
&\|\sP[\Omega, X^2\Phi]\|_{L^2_\Sigma}+\ell\| X_3(\sP[\Omega, X_2 \Phi])\|_{L^2_\Sigma}\nn \\
&\qquad\quad\les\l t\r^{-\frac{3}{2}+2\delta}\log \l t\r^{\f12\ell}\Delta_0^\frac{5}{4}+\l t\r^{-2}\log \l t\r^{\f12\M+1}\La_0,\label{10.1.1.21}\\
&\ell\|\sP[S, X^2\Phi]\|_{L^2_\Sigma}+(1-\ell)\|X_3(\sP[S, X_2\Phi])\|_{L^2_\Sigma}\nn\\
&\qquad\quad\les  \l t\r^{-\frac{3}{2}}\{\ell W_1[SX^2\Phi]^\f12(t)+(1-\ell)(W_1[X_3 S X_2\Phi]^\f12(t)+W_1[\Omega^2 X_2\Phi]^\f12(t))\nn\\
&\qquad\quad+\l t\r^{\delta}(\ell W_1[\Omega X^2\Phi]^\f12(t)+(1-\ell)W_1[X_3\Omega X_2\Phi]^\f12(t))\}\nn\\
 &\qquad\quad+\l t\r^{-2}\log \l t\r^{\f12\M+1}\La_0+\Delta_0^\frac{5}{4}\l t\r^{-\frac{7}{4}+3\delta}(\log \l t\r)^\f12,\label{3.5.1.24}\\
 &\begin{array}{lll}
\|\sP[\Omega, \Omega \bT\varrho]\|_{L^2_\Sigma}\les \Delta_0^\frac{5}{4}\l t\r^{-\frac{7}{4}+2\delta}\log \l t\r^\f12+\l t\r^{-2}\log \l t\r^\frac{3}{2}\La_0\\
\|\Omega(\sP[\Omega, \bT\varrho])\|_{L^2_\Sigma}\les \l t\r^{-1}(\log \l t\r)^5(\La_0+\Delta_0^\frac{5}{4}).
\end{array}\label{3.6.2.24}
\end{align}
\end{proposition}
\begin{proof}
We first consider (\ref{10.1.1.21}). Applying (\ref{2.14.3.24}) to $f=X^2\Phi$ and using Proposition \ref{7.16.1.21}, (\ref{3.12.1.21}) and Proposition \ref{7.15.5.22}, we have
\begin{align*}
\|\sP[\Omega, f]&+\f12 {}\rp{a}\pih_{LA}\bd^2_{\Lb A} f+\f12\Lb f \bJ[\Omega]_L\|_{L^2_\Sigma}\\
&\les\l t\r^{\delta}\Delta_0\|\sn^2 f+\bd^2_{LL}f\|_{L^2_\Sigma} +\l t\r^{-\frac{3}{4}+\delta}\Delta_0^\f12\|\bd^2_{LA} f\|_{L^2_\Sigma}\nn\\
&+\l t\r^{\delta}\Delta_0\|\bb^\f12(|L f|+\l t\r^{-1+\delta}|\sn f|)\|_{L^2_u L_\omega^4}\nn\\
&+\l t\r^{-1}\|\sn f\|_{L^2_\Sigma}+\l t\r^{-1+2\delta}\Delta_0\|\sn f\|_{L^2_\Sigma}\\
&\les \l t\r^{-\frac{3}{2}+2\delta}\Delta_0^\frac{3}{2}.
\end{align*}
Using Lemma \ref{8.1_com}, (\ref{5.21.1.21}), (\ref{5.25.1.21}), (\ref{12.19.1.23}),  (\ref{7.26.2.22}) and Proposition \ref{7.16.1.21} we derive
\begin{align}
\|\bJ[\Omega]_L\c \sta{X^2,\Lb}{\Phi}\|_{L^2_\Sigma}&\les \l t\r^{\f12\delta}\log \l t\r^\frac{5}{2}\Delta_0\|\bJ[\Omega]_L\|_{L_u^\infty L_\omega^4}+\|\fB\bJ[\Omega]_L\|_{L^2_\Sigma}\nn\\
&\les \l t\r^{-\frac{7}{4}+2\delta}\Delta_0^\frac{5}{4},\label{3.5.13.24}\\
\| \pioh_{LA}\bd^2_{\Lb A}X^2\Phi\|_{L^2_\Sigma}&\les \l t\r^{-1}\|\pioh_{LA}\fB\|_{L^2_\Sigma}+\l t\r^{-\frac{7}{4}+2\delta}\Delta_0^\frac{3}{2}\nn\\
&\les \l t\r^{-2}\log \l t\r^{\f12\M+1}\La_0+\l t\r^{-\frac{7}{4}+2\delta}\Delta_0^\frac{5}{4}.\nn
\end{align}
Combining the above estimates yields
\begin{align*}
&\|\sP[\Omega, X^2\Phi]\|_{L^2_\Sigma}\les\l t\r^{-2}\log \l t\r^{\f12\M+1}\La_0+ \l t\r^{-\frac{3}{2}+2\delta}\Delta_0^\frac{5}{4}.
\end{align*}

Applying (\ref{3.4.6.24}) to $f=X_2\Phi$ and using  Proposition \ref{7.16.1.21}, Corollary \ref{9.2.5.23} for the lower order estimates, and (\ref{3.12.1.21}) imply
 \begin{align*}
\|X_3&(\sP[\Omega, X_2\Phi]-\bJ[\Omega]^\a \p_\a X_2\Phi)+\f12 X_3(\pioh_{LA}\bd^2_{\Lb A}X_2\Phi)\|_{L^2_\Sigma} \\
&\les \l t\r^\delta\Delta_0\|\sn_{X_3}(\bd^2_{LL}+\sn^2 )X_2\Phi\|_{L^2_\Sigma}+\Delta_0^\frac{1}{2} \l t\r^{-\frac{3}{4}+\delta} \|\sn_{X_3}\bd^2_{LA}X_2\Phi\|_{L^2_\Sigma}\\
&+\Delta_0\l t\r^\delta\|\tir\bb^\f12(\bd^2_{LL}X_2\Phi+\sn^2 X_2\Phi)\|_{L_u^2 L_\omega^4}+\l t\r^{\delta-\frac{3}{4}}\Delta_0^\f12\|\bb^\f12 \tir \bd^2_{LA}X_2\Phi\|_{L_u^2 L_\omega^4}\\
&\les \l t\r^{-\frac{3}{2}+2\delta}\Delta_0^\frac{3}{2}.
\end{align*}
It follows by using (\ref{5.21.1.21}), (\ref{3.5.8.24}), (\ref{9.5.5.23}), (\ref{7.26.2.22}) and (\ref{3.3.4.24}) that
\begin{align*}
X_3(\pioh_{LA}\bd^2_{\Lb A} X_2\Phi)&=O(\l t\r^{-1}\fB)\sn_{X_3}^{\le 1}\pioh_{LA}+\sn_{X_3}\pioh_{LA} O(\l t\r^{-2+\delta}\Delta_0)_{L_u^2 L_\omega^4}\\
&+\pioh_{LA}O(\l t\r^{-1+\delta}\Delta_0)_{L^2_\Sigma}\\
&=O(\l t\r)^{-1}\fB\sn_\Omega\pioh_{LA} +O(\l t\r^{-2}\log \l t\r^{\f12\M+1}(\La_0+\Delta_0^\frac{5}{4})\\
&+\l t\r^{-\frac{7}{4}+2\delta}\Delta_0^\frac{3}{2})_{L^2_\Sigma}\\
&=O(\l t\r^{-2}\log \l t\r^{\f12\M+1}\La_0+\log \l t\r^\f12\l t\r^{-\frac{7}{4}+2\delta}\Delta_0^\frac{5}{4})_{L^2_\Sigma}.
\end{align*}
Using (\ref{5.25.1.21}), Proposition \ref{7.16.1.21} and (\ref{8.29.9.21}) we obtain
\begin{align}\label{3.5.15.24}
\begin{split}
\bJ[\Omega]^\a \sn_{X_3}\p_\a X_2\Phi&=O(\l t\r^{-\frac{7}{4}+\delta}\Delta_0)_{L_\omega^4}X_3 \Lb X_2\Phi+O(\l t\r^{-1+\delta}\Delta_0)_{L_\omega^4}X_3 L X_2\Phi\\
&+(O(\fB)+O(\l t\r^{-1+2\delta}\Delta_0)_{L_\omega^4})\sn_{X_3} \sn X_2\Phi\\
&=O(\l t\r^{-2}\log \l t\r^{\frac{\M}{2}}\La_0+\l t\r^{-\frac{7}{4}+3\delta}\Delta_0^\frac{5}{4}\log \l t\r^\f12)_{L^2_\Sigma}.
\end{split}
\end{align}
Using Proposition \ref{djoest}, Proposition \ref{7.15.5.22} and (\ref{10.10.2.23}), we bound
\begin{align}\label{3.5.16.24}
\begin{split}
\sn_{X_3} \bJ[\Omega]^\a \p_\a X_2\Phi&=\sn_{X_3}\bJ[\Omega]_L\Lb X_2\Phi+(O(\l t\r^{-1}|v|)+O(\l t\r^{-1+\delta}\Delta_0)_{L^2_u L_\omega^2})LX_2\Phi\\
&+(O(\fB)+O(\l t\r^{-1+\delta}\Delta_0)_{L_u^2 L_\omega^2})\sn X_2\Phi\\
&=\sn_X\bJ[\Omega]_L\Lb X_2\Phi+O(\l t\r^{-2}\log \l t\r^{\frac{\M}{2}}\La_0+\l t\r^{-\frac{7}{4}+2\delta}\Delta_0^\frac{5}{4}\log \l t\r^\f12)_{L^2_\Sigma}.
\end{split}
\end{align}
Note it follows from using (\ref{1.26.2.23}), (\ref{7.11.5.21}),  (\ref{8.2.2.23}), (\ref{8.29.9.21}) and (\ref{L2conndrv'}) that
\begin{align*}
\sn_{X_3}\bJ[\Omega]_L&=(1-\vs(X_3))X_3\Omega \wt{\tr\chi}+O(\l t\r^{-2}(\log \l t\r^{\f12\M+1}\La_0+\Delta_0^\frac{5}{4}\l t\r^{2\delta}))_{L^2_u L_\omega^2}\\
&=(1-\vs(X_3))O(\l t\r^{-\frac{3}{2}+\delta}\Delta_0)_{L_\omega^2}+O(\l t\r^{-2}(\log \l t\r^{\f12\M+1}\La_0+\Delta_0^\frac{5}{4}\l t\r^{2\delta}))_{L^2_u L_\omega^2}
\end{align*}
  
Using (\ref{8.23.1.23}), (\ref{12.19.1.23}), (\ref{10.10.2.23}), (\ref{8.29.9.21}) and the above estimate,  we obtain 
\begin{align}\label{3.5.14.24}
\begin{split}
&\sn_{X_3}\bJ[\Omega]_L\Lb X_2\Phi\\
&=O((1-\vs(X_2))\Delta_0\l t\r^{-1+\delta}+[\Lb\Phi])\\
&\times\{ (1-\vs(X_3))O(\l t\r^{-\frac{3}{2}+\delta}\Delta_0)_{L_\omega^2}+O(\l t\r^{-2}(\log \l t\r^{\f12\M+1}\La_0+\Delta_0^\frac{5}{4}\l t\r^{2\delta}))_{L^2_u L_\omega^2}\}\\
&=O(\l t\r^{-\frac{5}{2}+2\delta}\Delta_0^\frac{5}{4}+\bb^{-1}\l t\r^{-3}\log \l t\r^{\f12\M+1}\La_0)_{L^2_u L_\omega^2}.
\end{split}
\end{align}
Therefore
\begin{align*}
\sn_{X_3} \bJ[\Omega]^\a \p_\a X_2\Phi&=O(\l t\r^{-\frac{5}{2}+2\delta}\Delta_0^\frac{5}{4}+\bb^{-1}\l t\r^{-3}\log \l t\r^{\f12\M+1}\La_0)_{L^2_u L_\omega^2}.
\end{align*}
Thus we obtained (\ref{10.1.1.21}). 

Next we prove (\ref{3.5.1.24}). Applying (\ref{3.5.2.24}) to $f=X^2\Phi$, (\ref{7.03.4.21}), Proposition \ref{8.29.8.21} and Proposition \ref{7.16.1.21}, we derive
\begin{align*}
\|\sP[S, X^2\Phi]-&\frac{1}{4}L{}\rp{S}\ss \Lb X^2\Phi\|_{L^2_\Sigma}\les \|\sn^2 f+\bd^2_{LL} f+\l t\r^{-1}Lf+O(\Delta_0 \l t\r^\delta)\bd^2_{LA}f\|_{L^2_\Sigma}\\
 &+\l t\r^{-1+\delta}\Delta_0^\f12(\Delta_0^\f12\|\bb^\f12\tir\sn f\|_{L_u^2 L_\omega^4}+\l t\r^\delta\|\bb^\f12 \tir Lf\|_{L_u^2 L_\omega^4})\\
 &\les \l t\r^{-\frac{3}{2}}(W_1[SX^2\Phi]^\f12(t)+\l t\r^{\delta}W_1[\Omega X^2\Phi]^\f12(t))\\
 &+\Delta_0^\frac{5}{4}\l t\r^{-\frac{7}{4}+3\delta}(\log \l t\r)^\f12+\l t\r^{-2}\log \l t\r^{\f12\M}\La_0.
 \end{align*}
Similarly, it follows by applying (\ref{3.5.3.24}) to $f=X_2\Phi$, Proposition \ref{7.16.1.21}, (\ref{9.2.1.23}) and Proposition \ref{8.29.8.21} that
 \begin{align*}
\|X_3(\sP[S, X_2 \Phi]&-\bJ[S]^\a\p_\a X_2\Phi)\|_{L^2_\Sigma}\\
&\les \|\sn_{X_3}^{\le 1}(\sn^2 X_2\Phi, \bd^2_{LL}X_2\Phi)\|_{L^2_\Sigma}+\l t\r^{\delta}\Delta_0\|\bb^\f12\tir\big(\sn^2 X_2\Phi, \bd^2_{LL}X_2\Phi\big)\|_{L_u^2 L_\omega^4}\\
&+\Delta_0\l t\r^{\delta}\|\sn_{X_3}\bd^2_{LA}X_2\Phi\|_{L^2_\Sigma}+\l t\r^{\delta}\Delta_0\|\bb^\f12\tir\bd^2_{LA}X_2\Phi\|_{L_u^2 L_\omega^4}\\
&\les \l t\r^{-\frac{3}{2}}(W_1[\Omega^2 X_2\Phi]^\f12(t)+W_1[X_3SX_2\Phi]^\f12(t)+\l t\r^{\delta}W_1[X_3\Omega X_2\Phi]^\f12(t))\\
&+\Delta_0^\frac{5}{4}\l t\r^{-\frac{7}{4}+2\delta}(\log \l t\r)^\f12+\l t\r^{-2}\log \l t\r^{\f12\M}\La_0.
 \end{align*}
  Recall from (\ref{7.3.1.22}) that $L{}\rp{S}\ss=O(\l t\r^{-\frac{7}{4}+\delta}\Delta_0^\f12)$. Using this estimate, (\ref{8.21.1.22}) and (\ref{3.5.4.24}), we derive
 \begin{align}\label{3.6.1.24}
 \|L{}\rp{S}\ss (\Lb X^2\Phi, X_3 \Lb X_2\Phi)\|_{L^2_\Sigma}\les  \l t\r^{-2}\log \l t\r^{\f12\M+1}\La_0+\Delta^\frac{5}{4}\l t\r^{-\frac{7}{4}+2\delta}.
 \end{align}
Using (\ref{3.5.5.24}) and (\ref{8.23.1.23}), we have
 \begin{align*}
 \|X_3L{}\rp{S}\ss \Lb X_2\Phi\|_{L^2_\Sigma}\les \l t\r^{-2}\log \l t\r^{\f12\M+1}\La_0+\l t\r^{-2+2\delta}\log \l t\r^{\f12\M+1}\Delta_0^\frac{5}{4}.
 \end{align*}
 Similar to the proof for (\ref{8.7.1.21}), using (\ref{6.29.3.21}), (\ref{6.29.4.21}), Proposition \ref{djoest}, Proposition \ref{8.29.8.21} and Proposition \ref{7.16.1.21}, we bound
\begin{align*}
&\|\sn_{X_3}\bJ[S]^A \sn_A X_2\Phi+\sn_{X_3}\bJ[S]_\Lb LX_2\Phi)\|_{L^2_\Sigma}\les \l t\r^{-2}\log \l t\r^{\f12\M}\La_0+\Delta^\frac{5}{4}\l t\r^{-\frac{7}{4}+2\delta}(\log \l t\r)^\f12\\
&\|\bJ[S]^A\sn_{X_3}\sn X_2\Phi+\bJ[S]_\Lb X_3 L X_2\Phi\|_{L^2_\Sigma}\les\l t\r^{-1+\delta}\Delta_0^\f12 \|\bb^\f12\tir(\sn_{X_3}\sn X_2\Phi, \log \l t\r \bb^{-1}X_3 LX_2\Phi)\|_{L_u^2 L_\omega^4}\\
&\qquad\qquad\qquad\qquad\qquad\qquad\qquad\qquad\qquad+ \l t\r^{-1}\|X_3L X_2\Phi\|_{L^2_\Sigma}\\
&\qquad\qquad\qquad\qquad\qquad\qquad\qquad\qquad\quad\les  \l t\r^{-2}\log \l t\r^{\f12\M}\La_0+\Delta^\frac{5}{4}\l t\r^{-\frac{7}{4}+2\delta}\log \l t\r.
\end{align*} 
Summarizing the above estimates, we conclude (\ref{3.5.1.24}).

Next we prove (\ref{3.6.2.24}). Applying (\ref{2.14.3.24}) to $f=\Omega \bT\varrho$ with the help of (\ref{3.6.3.24}), (\ref{8.29.9.21}) and (\ref{2.14.1.24}) and (\ref{2.19.1.24})
\begin{align*}
\|&\sP[\Omega, \Omega\bT\varrho]+\f12 {}\rp{a}\pih_{LA}\bd^2_{\Lb A} \Omega\bT\varrho+\f12\Lb \Omega\bT \varrho \bJ[\Omega]_L\|_{L^2_\Sigma}\\
&\les \l t\r^{2\delta-2}\Delta_0^\frac{3}{2}+\l t\r^{-2+\delta}\Delta_0\sum_{X=\Omega, S}\|\bb^\f12\tir X\Omega\bT \varrho\|_{L_u^2 L_\omega^4}+(\l t\r^{-1}+\l t\r^{-1+2\delta}\Delta_0)\|\sn \Omega\bT\varrho\|_{L^2_\Sigma}\\
&\les \l t\r^{-2}\log \l t\r^\frac{3}{2}(\La_0+ \l t\r^{2\delta}(\log \l t\r)^3\Delta_0^\frac{5}{4}).
\end{align*}
Using (\ref{5.21.1.21}) and (\ref{3.6.4.24}), we obtain
\begin{align*}
\|{}\rp{a}\pih_{LA}\bd^2_{\Lb A} \Omega\bT\varrho\|_{L^2_\Sigma}\les \l t\r^{-\frac{7}{4}+2\delta}\Delta_0^\frac{3}{2}.
\end{align*}
Using (\ref{5.25.1.21}), (\ref{2.13.3.24}) and (\ref{3.6.4.24}), we infer
\begin{align*}
\|\Lb \Omega\bT \varrho \bJ[\Omega]_L\|_{L^2_\Sigma}&\les\l t\r^\delta\Delta_0\log \l t\r^\f12\|\bJ[\Omega]_L\|_{L^4_\omega}\\
&\les \Delta_0^2\l t\r^{-\frac{7}{4}+2\delta}\log \l t\r^\f12.
\end{align*}
Hence, we conclude 
\begin{equation*}
\|\sP[\Omega, \Omega\bT\varrho]\|_{L^2_\Sigma}\les \Delta_0^\frac{5}{4}\l t\r^{-\frac{7}{4}+2\delta}\log \l t\r^\f12+\l t\r^{-2}\log \l t\r^\frac{3}{2}\La_0
\end{equation*}
as stated in (\ref{3.6.2.24}).

Applying (\ref{3.4.6.24}) to $f=\bT\varrho$, 
\begin{align*}
&\Omega(\sP[{}\rp{a}\Omega, f]-\bJ[{}\rp{a}\Omega]^\a \p_\a f)+\f12 \Omega({}\rp{a}\pih_{LA}\bd^2_{\Lb A}f)\\
&=O(\l t\r^\delta)(\Delta_0\sn_\Omega(\bd^2_{LL} f+\sn^2 f)+\Delta_0^\f12 \l t\r^{-\frac{3}{4}} \sn_\Omega\bd^2_{LA}f)+O(\Delta_0\l t\r^\delta)_{L_\omega^4}\sn^2 f\\
&+O(\l t\r^{\delta-\frac{3}{4}}\Delta_0^\f12)_{L_\omega^4}\bd^2_{LA} f+\sn_\Omega{}\rp{\Omega}\ud\pih\bd^2_{LL}f.
\end{align*} 
For the last term in the above, using (\ref{7.03.4.21}), (\ref{8.23.1.23}), (\ref{2.20.2.24}) and (\ref{zeh}),  we infer
\begin{align*}
\sn_\Omega{}\rp{\Omega}\ud\pih\bd^2_{LL}\bT\varrho&=\sn_\Omega{}\rp{\Omega}\ud \pih (\l t\r^{-2}(O(1)\fB+S[L\Phi]+O(\l t\r^{-\frac{7}{4}+\delta}\Delta_0)_{L_\omega^4})+\bd^2_{LL} L\varrho)\\
&=O(\l t\r^{-2}\log \l t\r(\La_0+\log \l t\r^{3}\Delta_0^\frac{5}{4}))_{L^2_\Sigma}+O(\l t\r^\delta\Delta_0)_{L_\omega^4} (\l t\r^{-1}L[L\Phi]\\
&+\bd^2_{LL} L\varrho+O(\l t\r^{-\frac{7}{4}-2+\delta}\Delta_0)_{L_\omega^4}))
\end{align*}
Hence using Corollary \ref{9.2.5.23}, (\ref{7.03.4.21}), (\ref{3.14.2.24}) and Sobolev embedding we obtain 
\begin{align*}
\|\sn_\Omega{}\rp{\Omega}\ud\pih\bd^2_{LL}\bT\varrho\|_{L^2_\Sigma}\les \l t\r^{-2}\log \l t\r(\La_0+\log \l t\r^{3}\Delta_0^\frac{5}{4})
\end{align*}
and consequently by using Proposition \ref{7.16.1.21}  we conclude
\begin{align*}
\|\Omega(\sP&[{}\rp{a}\Omega, \bT\varrho]-\bJ[{}\rp{a}\Omega]^\a \p_\a \bT\varrho)+\f12 \Omega({}\rp{a}\pih_{LA}\bd^2_{\Lb A}\bT\varrho)\|_{L^2_\Sigma}\\
&\les\l t\r^{-2+2\delta}\Delta_0^\frac{3}{2}+\l t\r^\delta\Delta_0\|\bb^\f12\tir \sn^2\bT\varrho\|_{L_u^2 L_\omega^4}+\l t\r^{\delta-\frac{3}{4}}\Delta_0^\f12\|\bb^\f12\tir \bd^2_{LA}\bT \varrho\|_{L_u^2 L_\omega^4}\\
&+\|\sn_\Omega{}\rp{\Omega}\ud\pih\bd^2_{LL}\bT\varrho\|_{L^2_\Sigma}\\
&\les \l t\r^{-2+2\delta}\Delta_0^\frac{3}{2}\log \l t\r^\f12+\l t\r^{-2}\log \l t\r\La_0.
\end{align*}
Note that, due to (\ref{9.2.1.23}), (\ref{2.15.3.24}), (\ref{2.13.3.24}), (\ref{1.29.2.22}), (\ref{3.11.3.21}) and (\ref{LbBA2})
\begin{align*}
\|\bd^2_{\Lb A}\bT\varrho\|_{L^2_\Sigma}&\les \|\sn_A \Lb \bT\varrho+\sn \bT \varrho+\bA_{g,1}\Lb \bT\varrho\|_{L^2_\Sigma}\\
&\les \l t\r^{-1}\|\Lb\Omega\bT\varrho\|_{L^2_\Sigma}+\l t\r^{-1+\delta}\Delta_0\les \l t\r^{-1+\delta}\Delta_0.
\end{align*}
Using the above estimate, (\ref{5.21.1.21}), Proposition \ref{7.16.1.21} and Sobolev embedding, we derive
\begin{align*}
&\|\Omega({}\rp{a}\pih_{LA}\bd^2_{\Lb A}\bT\varrho)\|_{L^2_\Sigma}\\
&\les\|\sn_\Omega\pioh_{LA}\|_{L_\omega^4}\|\bb^\f12\tir\bd^2_{\Lb A}\bT\varrho\|_{L_u^2 L_\omega^4}+\|\pioh_{LA}\|_{L^\infty_\omega}\|\sn_\Omega\bd^2_{\Lb A}\bT\varrho\|_{L^2_\Sigma}\\
&\les\l t\r^{-\frac{7}{4}+2\delta}\log \l t\r^\f12\Delta_0^\frac{3}{2},
\end{align*}
It remains to consider the term $\Omega(\bJ[{}\rp{a}\Omega]^\a \p_\a \bT\varrho)$, which is expanded as
\begin{align*}
\Omega(\bJ[{}\rp{a}\Omega]^\a \p_\a \bT\varrho)&=\Omega(\bJ[\Omega]_L)\Lb \bT\varrho+\bJ[\Omega]_L \Omega \Lb \bT\varrho+\Omega(\bJ[\Omega]_A\sn_A\bT \varrho+\bJ[\Omega]_\Lb L \bT\varrho).
\end{align*}
Using (\ref{2.15.3.24}), (\ref{2.13.3.24}), (\ref{3.6.4.24}) and  Sobolev embedding, we have
\begin{equation*}
\|\Omega\Lb \bT\varrho\|_{L_u^2 L_\omega^4}\les \l t\r^{-1+\delta}\Delta_0
\end{equation*}
Using the above estimate and (\ref{5.25.1.21}) we derive
\begin{align*}
\|\bJ[\Omega]_L \Omega \Lb \bT\varrho\|_{L^2_\Sigma}\les \l t\r^{-\frac{7}{4}+2\delta}\log \l t\r^\f12\Delta_0^2.
\end{align*}

Using (\ref{2.20.2.24}), Proposition \ref{7.22.2.22} and (\ref{2.4.4.22}) , we obtain 
\begin{equation*}
\|X\Lb(\tir^{-1}\la)\|_{L_u^2 L_\omega^2}\les \l t\r^{-1}\log \l t\r(\La_0+\log \l t\r^3\Delta_0^\frac{5}{4}).
\end{equation*}
Hence, in view of (\ref{xdjo}), also using (\ref{2.22.3.24}), (\ref{2.19.1.24}) and (\ref{3.10.7.24}), (\ref{7.10.11.22}) is refined to 
\begin{align*}
\|X\Lb{}\rp{a}\ss\|_{L_u^2 L_\omega^2}\les \l t\r^{-1}\log \l t\r(\La_0+\log \l t\r^3\Delta_0^\frac{5}{4}).
\end{align*}
Using the above estimate, it follows by refining (\ref{7.22.8.22}) using (\ref{10.10.2.23}), (\ref{3.28.3.24}),  (\ref{3.7.1.24}) and (\ref{3.7.2.24}) that
\begin{align*}
\sn_\Omega \bJ[\Omega]_\Lb&=\sn_\Omega \ckk \J[\Omega]_\Lb +O(\l t\r^{-1}\log \l t\r(\La_0+(\log \l t\r)^3\Delta_0^\frac{5}{4}))_{L^2_u L_\omega^2}\\
&=\Omega^2\tr\chi+\l t\r^{-1}O(\Omega^{\le 1}v)+\sn_\Omega \sdiv\eta(\Omega)+O(\l t\r^{-1}\log \l t\r(\La_0+(\log \l t\r)^3\Delta_0^\frac{5}{4}))_{L^2_u L_\omega^2}\\
&=\Omega^2\tr\chi+O(\l t\r^{-1}\log \l t\r(\La_0+\log \l t\r^3\Delta_0^\frac{5}{4}))_{L^2_u L_\omega^2}.
\end{align*} 
Using (\ref{6.22.1.21}) and (\ref{3.11.3.21}), we infer
\begin{align*}
\sn_\Omega \bJ[\Omega]_\Lb L\bT \varrho&=O(\l t\r^{-1}\fB)\Omega^2\tr\chi+O(\l t\r^{-2}\log \l t\r(\La_0+(\log \l t\r)^3\Delta_0^\frac{5}{4}))_{L^2_\Sigma}
\end{align*}
Using  (\ref{12.19.1.23}) and (\ref{1.26.2.23}), for the first term on the right-hand side, we derive
\begin{align}\label{3.21.1.24}
\|\fB \Omega^2\tr\chi\|_{L^2_\Sigma}\les \l t\r^{-\frac{3}{2}+\delta}\Delta_0^\frac{5}{4}.
\end{align}
Hence by also using (\ref{8.23.1.23}) and (\ref{5.25.1.21}) we conclude 
\begin{align*}
\Omega( \bJ[\Omega]_\Lb L\bT \varrho)=O(\l t\r^{-2}\log \l t\r(\La_0+\l t\r^{2\delta}\Delta_0^\frac{5}{4}))_{L^2_\Sigma}.
\end{align*}
Using (\ref{5.25.1.21}) and (\ref{2.24.1.22}), (\ref{3.6.2.21}), (\ref{3.11.3.21}), Proposition \ref{1steng}, Proposition \ref{8.29.8.21} and Proposition \ref{7.16.1.21}, we obtain
 \begin{align*}
\Omega(\bJ[\Omega]_A \sn\bT\varrho)&=(\fB\c O(1)+O(\l t\r^{-1+2\delta}\Delta_0)_{L^2_u L_\omega^2})\sn\bT\varrho+O(\Delta_0 \l t\r^{-2+2\delta})_{L_\omega^4}\sn_\Omega \sn \bT\varrho\\
&+(O(\fB)+O(\l t\r^{-1+2\delta}\Delta_0))\sn_\Omega \sn \bT\varrho\\
&=O(\l t\r^{-2}\log \l t\r(\La_0+\l t\r^{3\delta}\Delta_0^\frac{5}{4}))_{L^2_\Sigma}.
\end{align*}
For the most crucial term, we use (\ref{3.28.3.24}), (\ref{7.20.7.22}), (\ref{7.22.6.22}) and  (\ref{7.22.7.22}) to write
\begin{align*}
X(\bJ[\Omega]_L)&-\{c^2 X(\Omega \tr\thetac)+X\sdiv\eta({}\rp{a}\Omega)+X(L+\tir^{-1}) \Omega \varrho+\tir^{-2}X^{\le 1}\la+\l t\r^{-1}X^{\le 1}v\}\\
&=O(\l t\r^{\frac{-7-\vs(X)}{4}+2\delta}\Delta_0)_{L^\infty_u L^2_\omega}, O(\l t\r^{-2+2\delta}\Delta_0^\frac{3}{2})_{L^2_u L_\omega^2}.
\end{align*}
From Proposition \ref{7.15.5.22}, (\ref{8.2.2.23}) and Proposition \ref{10.16.1.22}, we infer
\begin{align*}
X\sdiv \eta(\Omega), XL\Omega\varrho=O(\l t\r^{-\frac{7}{4}+\delta}\Delta_0)_{L_\omega^4}, O(\l t\r^{-1}\log\l t\r^{\f12\M}(\Delta_0^\frac{5}{4}+\La_0))_{L^2_\Sigma}
\end{align*}
where the second bound on the right-hand side is obtained in view of (\ref{8.29.9.21}) and (\ref{3.7.2.24}). 

Using Proposition \ref{9.8.6.22} and (\ref{11.26.2.23}), we have
\begin{align*}
\l t\r X\Omega\varrho, X^{\le 1}\la&
=\left\{\begin{array}{lll}
O((\Delta_0^\frac{5}{4} \l t\r^{\f12\delta}+\La_0)\log \l t\r^{\f12\M+1})_{L_\omega^4},\\
O((\Delta_0^\frac{5}{4}+\La_0)\log \l t\r^{\f12\M+1})_{L_u^2 L_\omega^2}\\
\end{array}\right.\\
 X^{\le 1}v&=O((\Delta_0^\frac{5}{4}+\La_0^\f12)\l t\r^{-1}\log \l t\r^{\f12\M})_{L_\omega^4}, O(\l t\r^{-1}(\Delta_0^\frac{5}{4}+\La_0)\log \l t\r^{\f12\M})_{L_u^2 L_\omega^2}.
\end{align*}
Combining the above two lines of estimates yields
\begin{align*}
\Omega\bJ[\Omega]_L-c^2\Omega^2\tr\thetac&=\left\{\begin{array}{lll}
 O(\l t\r^{-2}(\log \l t\r^{\f12\M+1}\La_0+ \l t\r^{2\delta}\Delta_0^\frac{5}{4}))_{L^2_u L_\omega^2}\\
O(\l t\r^{-\frac{7}{4}+2\delta}\Delta_0+\l t\r^{-2}\La_0^\f12\log \l t\r^{\f12\M})_{L_\omega^2}.
\end{array}\right.
\end{align*}
Hence using (\ref{2.13.1.24}), (\ref{1.26.2.23}) and the above estimate, we deduce
\begin{align}\label{3.7.3.24}
\begin{split}
\|\Omega(\bJ[\Omega]_L)\Lb \bT\varrho\|_{L^2_\Sigma}&\les \|\bb^{-1} \tir^{-1}\Omega\bJ[\Omega]_L\|_{L^2_\Sigma}+\|\Omega\bJ[\Omega]_L\|_{L_\omega^2}\log \l t\r\l t\r^{\f12\delta}\Delta_0\\
&\les\|\bb^{-\f12}\Omega^2\tr\thetac\|_{L^2_u L_\omega^2}+\l t\r^{-2}(\log \l t\r^{\f12\M+1}\La_0+ \l t\r^{2\delta}\Delta_0^\frac{5}{4})\\
&+\l t\r^{-\frac{3}{2}+\frac{3}{2}\delta} \Delta_0^2 \log \l t\r.
\end{split}
\end{align}

It remains to bound the first term on the right-hand side of (\ref{3.7.3.24}).  
Using (\ref{2.25.2.24'}), (\ref{1.29.1.24}) and (\ref{1.27.5.24}), we have
\begin{align*}
\|\Omega\log \bb \bb^{-1}\sn \sX\|_{L^2_u L_\omega^2}&\les \log \l t\r\Delta_0\|\sF\|_{L_u^2 L_\omega^4}\les\l t\r^{-2}\log \l t\r^4(\Delta_0^\frac{5}{4}+\La_0).
\end{align*}
Using the above estimate, (\ref{2.25.2.24'}) and (\ref{1.29.1.24}), we bound that
\begin{align*}
\|\bb^{-1}\Omega^2\sX\|_{L^2_u L_\omega^2}&\les \l t\r^{-1}\log \l t\r^\frac{17}{4}(\Delta_0^\frac{5}{4}+\La_0).
\end{align*}
Using (\ref{2.14.1.24}) and (\ref{8.29.9.21}), we then obtain 
\begin{align}\label{3.7.9.24}
\|\bb^{-1}\Omega^2\tr\thetac\|_{L^2_u L_\omega^2}\les\l t\r^{-1}\log \l t\r^\frac{17}{4}(\La_0+\Delta_0^\frac{5}{4}).
\end{align}
Substituting the above estimate to (\ref{3.7.3.24}) implies
\begin{align*}
\|\Omega(\bJ[\Omega]_L)\Lb \bT\varrho\|_{L^2_\Sigma}\les \l t\r^{-1}(\log \l t\r)^5(\La_0+\Delta_0^\frac{5}{4}).
\end{align*}
Hence we conclude (\ref{3.6.2.24}).
\end{proof}

\begin{proposition}\label{7.16.2.21}
Let $X\in \{\Omega, S\}$ and $0<t<T_*$. Under the assumptions of (\ref{3.12.1.21})-(\ref{6.5.1.21}), we have
\begin{align}
&\|X^2(\sP[{}\rp{a}\Omega, \Phi])-(1-\vs^+(X^2))(X^2{}\rp{a}\Omega\tr\chi\Lb \Phi)\|_{L^2_\Sigma}\nn\\
&\les  \l t\r^{-\frac{3}{2}+2\delta}\Delta_0^\frac{5}{4}(\log \l t\r)^\f12+\l t\r^{-2}\log \l t\r^{\f12\M+1}\La_0\nn\\
&+\l t\r^{-\frac{3}{2}}(W_1[X^2\Omega \Phi]^\f12(t)+\vs^+(X^2)\sum_{\vs(V)=\vs^-(X^2)}W_1[V\Omega S\Phi, \Omega^3\varrho]^\f12(t))\label{8.14.4.21}\\
&\|X^2(\sP[S, \Phi])\|_{L^2_\Sigma}\nn\\
&\les \l t\r^{-\frac{3}{2}}\{W_1[X^2 S\Phi, X^2 S\varrho]^\f12(t)+W_1[X^2 \Omega\Phi]^\f12(t)+ \sum_{\vs(V)=\vs^{-}(X^2)}W_1[V\Omega^2\Phi, V\Omega^2 \varrho]^\f12(t)\}\nn\\
&+\l t\r^{-2}\log \l t\r^{\f12(\M+7)}\La_0+\Delta^\frac{5}{4}\l t\r^{-\frac{7}{4}+2\delta}\log \l t\r.\label{8.17.5.21}
\end{align}
\end{proposition}
\begin{proof}
In what follows, we will frequently use Corollary \ref{9.2.5.23} and Lemma \ref{3.17.2.22} without mentioning.
We first prove (\ref{8.14.4.21}).
In view of (\ref{6.23.1.23}), using (\ref{5.21.1.21}) and (\ref{3.25.1.22}), we bound
\begin{align*}
X^2(\sP[{}\rp{a}\Omega, \Phi]&-{}\rp{a}\bJ^\a\p_\a \Phi+\f12 {}\rp{a}\pih_{LA}\bd^2_{\Lb A} \Phi)\\
&=X^2\left({}\rp{\Omega}\ud\pih\c(\sn^2\Phi, \bd^2_{LL}\Phi)+\pioh_{\Lb A} \bd^2_{LA}\Phi\right)\\
&=O(\l t\r^\delta\Delta_0)_{L_u^2 L_\omega^2}(\sn^2\Phi+\bd_{LL}^2\Phi)+O(\l t\r^{\delta-\f12}\Delta_0)_{L_u^2 L_\omega^2}\bd_{LA}^2\Phi\\
&+O(\l t\r^\delta\Delta_0)_{L_\omega^4}\sn_X(\sn^2\Phi+\bd_{LL}^2\Phi)+O(\l t\r^{-\frac{3}{4}+\delta}\Delta_0^\f12)_{L_\omega^4}\sn_X\bd^2_{LA}\Phi\\
&+O(\l t\r^\delta \Delta_0)\sn_X^2(\sn^2\Phi+\bd_{LL}^2\Phi)+O(\l t\r^{-\frac{3}{4}+\delta}\Delta_0^\f12)\sn_X^2\bd^2_{LA}\Phi\\
&=O(\l t\r^{-\frac{3}{2}+2\delta}\Delta_0^\frac{3}{2})_{L^2_\Sigma}
\end{align*}
where we also used (\ref{7.31.4.22}), Proposition \ref{7.16.1.21} and the decay estimates 
\begin{equation}\label{6.10.2.24}
|\sn^2\Phi, \bd^2_{LL}\Phi, \bd^2_{LA}\Phi|\les \l t\r^{-\frac{11}{4}+\delta}\Delta_0^\f12
\end{equation}
 in the above. 

Using (\ref{12.25.3.23}), (\ref{9.5.5.23}) and Proposition \ref{7.16.1.21}, we derive
\begin{align}\label{3.5.9.24}
\begin{split}
\sn_X^2({}\rp{a}\pih_{LA}\bd^2_{\Lb A} \Phi)&=\sn_X^{\le 2}\pioh_{LA}O(\l t\r^{-1}\fB)+\sn_X^2\pioh_{LA}O(\l t\r^{-2+\delta}\Delta_0)\\
&+\sn_X\pioh_{LA}O(\l t\r^{-2+\delta}\Delta_0)_{L_u^2 L_\omega^4}+\pioh_{LA} O(\l t\r^{-1+\delta}\Delta_0)_{L^2_\Sigma}.
\end{split}
\end{align}
Using  (\ref{10.8.2.23}), (\ref{3.5.8.24}), (\ref{3.28.3.24}), (\ref{L2BA2'}) and (\ref{10.10.2.23}), we derive
\begin{align*}
\sn_X^2\pioh_{LA}&=O(1)\sn_X^2(\eta(\Omega), {}\rp{a}v^*, \sn\la)+O(\Delta_0^2\l t\r^{-\frac{7}{4}+2\delta})_{L^4_\omega}\\
&=(1-\vs^+(X^2))O(1)\sn_X^2 \sn\la +O(\l t\r^{-1}\log \l t\r^{\f12\M+1}(\l t\r^\delta\Delta_0^\frac{5}{4}+\La_0))_{L^2_u L_\omega^2}.
\end{align*}
Hence we obtain from the above estimate that
\begin{equation*}
\l t\r^{-1}\|(\sn_X^2\pioh_{LA}-(1-\vs^+(X^2))\sn_X^2 \sn\la) \fB\|_{L^2_\Sigma}\les\l t\r^{-2}\log \l t\r^{\f12\M+1}(\l t\r^\delta\Delta_0^\frac{5}{4}+\La_0).
\end{equation*}
Using (\ref{2.24.5.24}) and (\ref{12.19.1.23}), we infer
\begin{align*}
\l t\r^{-1}\|\sn_\Omega^2 \sn\la \fB\|_{L^2_\Sigma}&\les\| \|\tir \bb^\f12\fB\|_{L_u^2 L_\omega^\infty}\l t\r^{-\frac{3}{2}+\delta}\Delta_0\les \l t\r^{-\frac{3}{2}+\delta}\Delta_0^\frac{5}{4}.
\end{align*}

Substituting the above two estimates into (\ref{3.5.9.24}), also applying (\ref{7.26.2.22}), (\ref{3.5.8.24}), (\ref{5.21.1.21}) and (\ref{3.25.1.22}), yields 
\begin{align}
\|\sn_X^2({}\rp{a}\pih_{LA}\bd^2_{\Lb A} \Phi)\|_{L^2_\Sigma}&\les  \l t\r^{-\frac{3}{2}+2\delta}(\log \l t\r)^\f12\Delta_0^\frac{5}{4}+\l t\r^{-2}\log \l t\r^{\f12\M+1}\La_0.\label{3.5.10.24}
\end{align}

Now we write
\begin{align}
X^2({}\rp{a}\bJ^\a \bd_\a \Phi)&=X^2{}\rp{a}\bJ_L\c \bd_\Lb \Phi+ X{}\rp{a}\bJ_L X \Lb \Phi+{}\rp{a}\bJ_L X^2 \Lb \Phi\nn\\
&+X^2({}\rp{a}\bJ_\Lb\c L \Phi+{}\rp{a}\bJ_A \sn_A \Phi)\label{9.30.3.21}
\end{align}
We first prove 
\begin{equation}\label{3.5.6.24}
\|X^2({}\rp{a}\bJ_\Lb\c L \Phi+{}\rp{a}\bJ_A \sn_A \Phi)\|_{L^2_\Sigma}\les\log \l t\r^{\f12\M}( \l t\r^{-2}\La_0+\l t\r^{-\frac{3}{2}+\delta}\Delta_0^\frac{5}{4}).
\end{equation}
Indeed, using (\ref{8.8.6.21}), Proposition \ref{1steng}, (\ref{8.29.9.21}), (\ref{3.12.1.21}) and (\ref{7.25.2.22}), we derive
\begin{align*}
&X^2{}\rp{a}\bJ_\Lb\c L \Phi+\sn_X^2{}\rp{a}\bJ_A \sn_A \Phi\\
&=O(\l t\r^{3\delta-1}\Delta_0)_{L^2_u L_\omega^2}+[LX^2\Omega\Phi])L\Phi+(O(\fB)+(\Delta_0\l t\r^{-1+2\delta})_{L^2_u L_\omega^2})\sn\Phi\\
&=\l t\r^{-2}O(\log \l t\r^{\f12\M}\La_0+\l t\r^{\f12+\delta}\log \l t\r^{\frac{\M}{2}}\Delta_0^\frac{5}{4})_{L^2_\Sigma}.
\end{align*}
Similarly to (\ref{3.5.16.24}), we bound
\begin{align*}
X{}\rp{a}\bJ_\Lb\c X L \Phi+\sn_X{}\rp{a}\bJ_A \sn_X\sn_A \Phi=O(\l t\r^{-2}\log \l t\r^{\f12\M}\La_0+\l t\r^{-\frac{7}{4}+2\delta}\log \l t\r^{\f12})_{L^2_\Sigma}
\end{align*}
Similar to (\ref{3.5.15.24}), we have
\begin{align*}
{}\rp{a}\bJ_\Lb\c X^2 L \Phi&+{}\rp{a}\bJ_A \sn_X^2\sn_A \Phi
=O(\l t\r^{-2}\log \l t\r^{\f12\M}\La_0+\l t\r^{-\frac{7}{4}+3\delta}\log \l t\r^{\f12})_{L^2_\Sigma}.
\end{align*}
Hence (\ref{3.5.6.24}) is proved.

In view of (\ref{8.8.6.21}) and (\ref{7.1.1.23}) we decompose the first term on the right-hand side of (\ref{9.30.3.21}) into
\begin{align*}
X^2{}\rp{a}\bJ_L\c\Lb \Phi&= \{\sn_X^2(\ckr^{-1}{}\rp{a}v^\sharp)+X^2(\tir^{-2}\la)+O(1)X^2\Omega\tr\chi\}\Lb\Phi\\
&+ \bb^{-1}O(\l t\r^{-\frac{3}{2}}W_1[X^2\Omega \Phi]^\f12(t))_{L^2_\Sigma}
+O(\l t\r^{-\frac{3}{2}+2\delta}\Delta_0^\frac{3}{2}+\l t\r^{-2}\La_0\log \l t\r^{\f12\M})_{L^2_\Sigma}.
\end{align*}

For the first two terms, using (\ref{3.5.7.24})  we bound
\begin{align*}
\|X^2(\tir^{-2}\la)\fB\|_{L^2_\Sigma}\les \l t\r^{-2}\log \l t\r^{\f12\M+1}(\La_0+\l t\r^\delta\Delta_0^\frac{5}{4}),
\end{align*}
and we obtain by  using (\ref{3.28.3.24}) and (\ref{10.10.2.23}) that
\begin{align*}
\|X^2(\tir^{-1}v^\sharp)\fB\|_{L^2_\Sigma}&\les \l t\r^{-1}\||X^{\le 2} v| \fB\|_{L^2_\Sigma}+\l t\r^{-3+2\delta}\Delta_0^\frac{3}{2}\\
&\les \l t\r^{-2}(\log \l t\r)^{\f12\M}(\La_0+\Delta_0^\frac{5}{4}).
\end{align*}
If $\vs^+(X^2)=1$, we apply (\ref{9.15.4.22});
also by using (\ref{3.21.1.24}) and (\ref{7.1.1.23}), we derive
\begin{align*}
\|X^2{}\rp{a}\bJ_L\c\Lb \Phi&-O(1)(1-\vs^+(X^2)) X^2\Omega\wt{\tr\chi}\Lb \Phi\|_{L^2_\Sigma}\\
&\les \l t\r^{-\frac{3}{2}}(W_1[X^2\Omega \Phi]^\f12(t)+\vs^+(X^2)\sum_{\vs(V)=\vs^-(X^2)}W_1[V\Omega S\Phi, \Omega^3\varrho]^\f12(t))\\
&+\l t\r^{-2}(\log \l t\r^{\f12 \M+1}\La_0+\l t\r^{\f12+2\delta} \Delta_0^\frac{5}{4}).
\end{align*} 

Similar to (\ref{3.5.14.24})
\begin{equation*}
\sn_X\bJ[\Omega]_L X_2\Lb\Phi=O(\l t\r^{-\frac{5}{2}+2\delta}\Delta_0^\frac{5}{4}+\bb^{-1}\l t\r^{-3}\log \l t\r^{\f12\M}\La_0)_{L^2_u L_\omega^2}.
\end{equation*}
Recall from (\ref{3.5.13.24}) that
\begin{align*}
\|{}\rp{a}\bJ_L \sn_X^2 \Lb \Phi\|_{L^2_\Sigma}&\les \l t\r^{-\frac{7}{4}+2\delta}\Delta_0^\frac{5}{4}.
\end{align*}
We therefore conclude (\ref{8.14.4.21}).

Next we consider the estimate for $X_2 X_1(\sP[S, \Phi])$. Recalling from (\ref{6.24.2.23})  and applying (\ref{5.25.2.21})-(\ref{7.5.2.21}) and Proposition \ref{7.22.2.22}, we deduce
\begin{align*}
X^2(\sP[S, \Phi]-\bJ[S]^\mu \p_\mu \Phi)&=X^2(\piSh^{AB}\sn_A\sn_B\Phi+\piSh_{\Lb \Lb}\bd^2_{LL}\Phi+{}\rp{S}\ud\pih\bd^2_{L A}\Phi)\\
&=O(\l t\r^\delta\Delta_0)\sn_X^2\bd^2_{LA}\Phi+\sn_X^{\le 2}(\sn^2\Phi, \bd^2_{LL}\Phi)\\
&+O(\l t\r^\delta\Delta_0^\f12)_{L_\omega^4} \sn_X (\bd^2_{LA}\Phi+\sn^2\Phi+\bd^2_{LL}\Phi)\\
&+O(\l t\r^\delta\Delta_0)_{L^2_u L_\omega^2}(\bd^2_{LA}\Phi+\sn^2\Phi+\bd^2_{LL}\Phi).
\end{align*}

Using Proposition \ref{7.16.1.21}, (\ref{7.31.4.22}), Proposition \ref{8.29.8.21} and Proposition \ref{9.8.6.22}  we can obtain
\begin{align*}
&\sn_X (\bd^2_{LA}\Phi, \sn^2\Phi, \bd^2_{LL}\Phi)=O(\l t\r^{-\frac{11}{4}+\delta}\Delta_0)_{L^2_u L_\omega^4},
\end{align*}
and then 
\begin{align*}
&\|X^2(\sP[S, \Phi]-\bJ[S]^\mu \p_\mu \Phi)\|_{L^2_\Sigma}\\
&\les  \l t\r^{-\frac{3}{2}}\big(\sum_{\vs(V)=\vs^{-}(X^2)}W_1[V\Omega^2\Phi]^\f12(t)+W_1[X^2S\Phi]^\f12(t)+\l t\r^\delta\Delta_0 W_1[X^2\Omega\Phi]^\f12(t)\big)\\
&+\l t\r^{-2}(\log \l t\r^{\f12(\M+7)}\La_0+\log \l t\r^\f12 \l t\r^{\frac{1}{4}+2\delta}\Delta_0^\frac{5}{4}).
\end{align*}
where we also employed (\ref{6.10.2.24}).

Using Proposition \ref{8.10.2.21},  Proposition \ref{1steng} and Proposition \ref{9.8.6.22}, we derive
\begin{align*}
\|\sn_X^2\bJ[S]^\a \p_\a \Phi\|_{L^2_\Sigma}&\les\l t\r^{-\frac{3}{2}}\{W_1[X^2S\varrho]^\f12(t)+\sum_{\vs(V)=\vs^-(X^2)}W_1[V\Omega^2\varrho]^\f12(t)\}\nn\\
&+\l t\r^{-2}(\l t\r^{3\delta}\Delta_0^\frac{5}{4}+\log \l t\r^{\f12\M}\La_0).
\end{align*}
The remaining estimate follows by repeating the proof of (\ref{3.5.1.24}) from (\ref{3.6.1.24}) to the end, with the help of Proposition \ref{7.16.1.21}, 
\begin{align*}
\sum_{a=0}^1\|\sn_X^a\bJ[S]^\a \sn_X^{2-a}\p_\a \Phi\|_{L^2_\Sigma}\les \l t\r^{-2}\log \l t\r^{\f12\M+1}\La_0+\Delta^\frac{5}{4}\l t\r^{-\frac{7}{4}+2\delta}\log \l t\r.
\end{align*}
Thus by combining the above estimates, (\ref{8.17.5.21}) is proved. 
\end{proof} 

Next we give an important comparison result between $\bb^{-1}\Omega^3\tr\chi\Lb \Phi$ and $\tir (\tir\sn)^2 \sF\c \Lb \varrho$ in $L^2_\Sigma$.
\begin{lemma}\label{3.8.3.24}
With $\M\ge 15$, we have
\begin{align*}
&\| \aaa^{-\frac{\M}{2}}\Omega^3\tr\chi\c \Lb\Phi\|_{L^2_\Sigma}\\
&\les (C\M_0+1)^2 \sup_{t'\in(0, t]}(E_{-\M}[\Omega^3\varrho]^\f12(t')+W_{1,-\M-1}[\Omega^3\Phi]^\f12(t'))\log \l t\r^{-1}\l t\r^{-1}\\
&+\l t\r^{-1}\log \l t\r^{-1}(\La_0+\Delta_0^\frac{5}{4})+\l t\r^{-\frac{3}{2}}\log \l t\r^{\f12}\sup_{ 0<t'\le t} W_{1,-\M-1}[\Omega^3\Phi]^\f12(t').
\end{align*}
\end{lemma}
\begin{proof}
 Using (\ref{8.7.1.22}), Lemma \ref{3.17.2.22} and Proposition \ref{7.15.5.22}, we have
\begin{equation*}
\|\Omega^3\tr\chi\|_{L_u^2 L_\omega^2}\les \l t\r^{-1+\delta} \Delta_0.
\end{equation*}
Then due to $\Lb\Phi=[\bar\bp\Phi]+\Lb\varrho$, and in view of (\ref{3.6.2.21})
\begin{align*}
\bb^\f12 \aaa^{-\frac{\M}{2}}\Omega^3 \tr\chi\c \Lb \Phi=\bb^\f12 \aaa^{-\frac{\M}{2}}\Omega^3 \tr\chi \Lb\varrho+ O(\l t\r^{-3+2\delta}\log \l t\r^{-\frac{\M-1}{2}}\Delta_0^\frac{3}{2})_{L_u^2 L_\omega^2}
\end{align*}
where the error is negligible. We then focus on treating the first term on the right-hand side. It is direct to expand 
\begin{align*}
\bb^{\f12}\Omega^3\tr\chi=\bb^\frac{3}{2}\{\Omega^2(\bb^{-1}\Omega\tr\chi)+\Omega(\bb^{-1})\Omega^2\tr\chi+\Omega^2(\bb^{-1})\Omega\tr\chi\}.
\end{align*}
Using (\ref{1.26.2.23}) and (\ref{9.12.2.22}), we derive
\begin{align*}
\|\bb^\frac{3}{2}\Omega(\bb^{-1})\Omega^2 \tr\chi\|_{L_u^2 L_\omega^2}&\les \l t\r^{-\frac{3}{2}+\delta}\Delta_0 \|\bb^{-\f12}\Omega(\bb)\|_{L_u^2 L_\omega^\infty}\les \l t\r^{-\frac{3}{2}+\frac{3}{2}\delta}\Delta_0^2 \log\l t\r
\end{align*}
and due to (\ref{3.10.3.24})
\begin{align*}
\|\bb^\frac{3}{2}\Omega^2(\bb^{-1})\Omega\tr\chi\|_{L_u^2 L_\omega^2}\les \l t\r^{-\frac{7}{4}+\delta}\Delta_0\|\bb^\frac{3}{2}\Omega^2(\bb^{-1})\|_{L_u^2 L_\omega^4}\les \l t\r^{-\frac{7}{4}+\frac{3}{2}\delta}\Delta_0^2\log \l t\r^4.
\end{align*}
Hence 
\begin{align*}
\bb^{\f12}\Omega^3\tr\chi=\bb^{\frac{3}{2}}\Omega^2(\bb^{-1}\Omega\tr\chi)+O( \l t\r^{-\frac{7}{4}+\frac{3}{2}\delta}\Delta_0^2\log \l t\r^4)_{L^2_u L_\omega^2}.
\end{align*}
To relate the leading term on the right-hand side with $\tir(\tir\sn)^2\sF$, we further compute
\begin{align}\label{12.22.1.24}
\Omega^2(\bb^{-1}\Omega\tr\chi)=\Omega^2(\sF(\Omega))+\Omega^2(\bb^{-1}\Omega\Xi_4).
\end{align}
Applying (\ref{4.22.4.22}), (\ref{2.25.2.24'}), (\ref{1.29.1.24}) together with Sobolev embedding, we bound 
\begin{align*}
\|\bb^\frac{3}{2}\Omega^2(\sF(\Omega)) \Lb \varrho\tir\|_{L^2_u L_\omega^2}&\les\|\bb^\frac{3}{2}\tir^2 (\tir\sn)^{1+\le 1}\sF \Lb \varrho\|_{L_u^2 L_\omega^2}+\l t\r^{\delta+\frac{1}{4}}\log \l t\r^\f12\Delta_0^\f12 \| \sF\|_{L_u^2 L_\omega^4}\\
&\les \|\bb^\frac{3}{2} \tir^2(\tir \sn)^2\sF\c \Lb \varrho\|_{L_u^2 L_\omega^2}+\l t\r^{-1}(\log \l t\r)^5(\Delta_0^\frac{5}{4}+\La_0).
\end{align*}
We will bound the first term on the right-hand side by using transport equations. Note that for $S_{t,u}$ tangent tenor $F$, in view of (\ref{lb}), there holds 
\begin{align*}
\sn_L(\bb \aaa^{-\M} |F|^2)=(- k_{\bN\bN}-\M\frac{1}{(\tir+3)\aaa})\bb\aaa^{-\M}|F|^2+2\bb \aaa^{-\M}\sn_L F\c F.
\end{align*}
Due to (\ref{11.11.2.23}) and (\ref{4.12.1.24}), with $\M=C\M_0$ sufficiently large, the first term on the right-hand side is negative. Hence
\begin{align*}
\|F\bb^\f12\aaa^{-\f12\M}(t)\|_{L_{[u_1, u_*]}^2 L_\omega^2}\les \|(F \bb^\f12 \aaa^{-\f12 \M})(0)\|_{L^2_{[u_1, u_*]} L_\omega^2} +\int_0^t \|\sn_L F \bb^\f12 \aaa^{-\f12\M}\|_{L^2_{[u_1, u_*]}  L_\omega^2}d t'. 
\end{align*}
Applying the above inequality to $F=\tir^4(\tir \sn)^2\sF\c \bb\Lb \varrho$, we obtain 
\begin{align*}
\|\bb^\frac{3}{2} \aaa^{-\f12\M}\tir^4(\tir \sn)^2\sF\c \Lb \varrho\|_{L_u^2 L_\omega^2}&\les \|( \bb^\frac{3}{2} \aaa^{-\f12\M}\tir^4(\tir \sn)^2\sF\c \Lb \varrho)(0)\|_{L^2_{[u_1, u_*]} L_\omega^2}\\
&+\int_0^t \|\sn_L F \bb^\f12 \aaa^{-\f12\M}\|_{L^2_{[u_1, u_*]}  L_\omega^2}d t'. 
\end{align*}
Similar to (\ref{12.13.1.21}), differentiating (\ref{12.5.2.21}), we derive
\begin{align*}
||\sn_L (\tir^3(\tir \sn)^2 \sF)|&\les |(\chih, \bA_b)\c\tir^3 (\tir \sn)^2 \sF+
+\tir^3\tir\sn\tr\chi (\tir\sn) \sF\\
&+\sum_{a=0}^1\tir^3(\tir \sn)^{1-a}(\tir(\bR_{BCLA}, \chi\zb)(\tir \sn)^a \sF)|+|\tir^3 (\tir \sn)^2 G_2|.
\end{align*}
Hence the most critical term to bound $|\sn_L F|$ is 
\begin{align*}
|\tir^3(\tir \sn)^2 G_2 \c \bb \tir\Lb \varrho|\le |I_1|+|I_2|
\end{align*}
where 
\begin{align*}
I_1&=\tir^3 (\tir\sn)^2 {\ti G}_1 \bb\tir \Lb \varrho,\quad {\ti G}_1=\bb^{-1}\sn\Xi_4(\Xi_4+\f12\sX)\\
I_2&=\tir^3
(\tir\sn)^2\{\bb^{-1}\sn[L\Phi]\sX+\bb^{-1}(-2\sn \chih\c \chih+\sn \N(\Phi, \p\Phi))-\chih\c \sF\}\c \bb \tir \Lb \varrho.
\end{align*}
With 
\begin{align*}
I_3&= \tir((\chih, \bA_b)\c\tir^3 (\tir \sn)^2 \sF+\tir^3\tir\sn\tr\chi (\tir\sn) \sF+\sum_{a=0}^1\tir^3(\tir \sn)^{1-a}(\tir(\bR_{BCLA}, \chi\zb)(\tir \sn)^a \sF))\bb \Lb \varrho\\
I_4&=\tir^3 (\tir \sn)^2 \sF L(\bb\tir\Lb \varrho),   
\end{align*}
we have
\begin{equation}\label{12.20.1.24}
|\sn_L F|\les |I_1|+|I_2|+|I_3|+|I_4|. 
\end{equation}
The estimates of $I_2, I_3$ can be similarly obtained as in the proof of (\ref{10.26.3.23}).
It follows by using Proposition \ref{7.15.5.22}, (\ref{1.29.4.22}) and (\ref{11.28.5.23})
\begin{align*}
\int_0^t \|\bb^\f12 \aaa^{-\frac{\M}{2}}(I_3, I_4)\|_{L_u^2 L_\omega^2}&\les \Delta_0\int_0^t \tir^3 \l t\r^{-\frac{7}{4}+2\delta}\|(\tir \sn)^{\le 2} \sF\|_{L_u^2 L_\omega^2} dt'\\
&\les \Delta_0^2 \int_0^t \l t'\r^{-\frac{3}{4}+3\delta}\les \l t\r^{\frac{1}{4}+3\delta}\Delta_0^2. 
\end{align*}
\begin{align*}
\int_0^t \|\bb^\f12 \aaa^{-\frac{\M}{2}}I_2\|_{L_u^2 L_\omega^2}\les \l t\r^{\frac{1}{2}+3\delta}\Delta_0^2+\int_0^t \tir^2\| (\tir \sn)^2(\sn[L\Phi])\bb^{-\f12} \aaa^{-\f12\M}\|_{L_u^2 L_\omega^2} dt'.
\end{align*}
Using $|\bb\tir \Lb \varrho|\les 1$ and symbolically $\Xi_4=\Lb\varrho+[L\Phi]$, we derive
\begin{align*}
&\int_0^t \|\bb^\f12 \aaa^{-\frac{\M}{2}} I_1 \|_{L_u^2 L_\omega^2}\\
\displaybreak[0]
&\les \int_0^t \|\tir^2 (\tir \sn)^2 \sn\Lb\varrho \bb^\f12\aaa^{-\frac{\M}{2}}\|_{L_u^2 L_\omega^2}\|\tir \Lb\varrho\|_{L^\infty_\Sigma}\\
&+\int_0^t  \|\tir^2 (\tir \sn)^2 \sn [L\Phi]\bb^{-\f12}\aaa^{-\frac{\M}{2}}, \tir^3 \sum_{i=0}^1(\tir\sn)^i\sn\Xi_4(\tir\sn)^{2-i}(\bb^{-1}(\Xi_4+\tr\chi))\bb^\frac{1}{2}\aaa^{-\f12\M}\|_{L_u^2 L_\omega^2}.
\end{align*}
We derive in view of (\ref{2.22.1.24}), Lemma \ref{3.17.2.22}, Lemma \ref{6.30.4.23}, Lemma \ref{10.10.3.23}, the lower order estimates in (\ref{8.29.9.21}), Proposition \ref{1steng} and (\ref{3.6.2.21}) that
\begin{align*}
\|\tir^2 (\tir\sn)^2 \sn \Lb \varrho\bb^\f12 \aaa^{-\frac{\M}{2}}\|_{L_u^2 L_\omega^2}&\les (C\M_0+1)\sup_{t'\in(0, t]}(E_{-\M}^\f12[\Omega^3\varrho](t')+W_{1,-\M-1}^\f12[\Omega^3\Phi](t'))\\
&+\log \l t\r^{7-\frac{\M}{2}}(\La_0+\Delta_0^\frac{5}{4})
\end{align*}
and 
\begin{align*}
\int_0^t  \|\tir^2 (\tir \sn)^2 \sn [L\Phi]\bb^{-\f12}\aaa^{-\frac{\M}{2}}\|_{L_u^2 L_\omega^2}\les \l t\r^{\f12}\log \l t\r^\f12\sup_{ 0<t'\le t} W^\f12_{1,-\M-1}[\Omega^3\Phi](t)+(\La_0+\Delta_0^\frac{5}{4})\l t\r^\delta
\end{align*}
The remaining terms are nonlinear error terms.   Due to $\Xi_4=\Lb \varrho+[L\Phi]$, similar to (\ref{6.2.1.24}) and (\ref{6.2.2.24}), 
we derive by using Proposition \ref{7.15.5.22}, Lemma \ref{5.13.11.21} (5), (\ref{2.20.2.24}), (\ref{2.14.1.24}) and Sobolev embedding that
\begin{align*}
\sum_{i=0}^1&\|\tir^3 (\tir\sn)^i\sn\Xi_4(\tir\sn)^{2-i}(\bb^{-1}(\Xi_4+\tr\chi))\bb^\frac{1}{2}\aaa^{-\f12\M}\|_{L_u^2 L_\omega^2}\\
&\les\Delta_0^\f12((1+\log \l t\r\Delta_0)\|\Omega^{1+\le 2} \Lb\varrho \aaa^{-\frac{\M}{2}}\tir\|_{L_u^2 L_\omega^2}^\f12 \|\Omega^{1+\le 1}\Lb\varrho \aaa^\frac{-\M}{2}\tir\|_{L_u^2 L_\omega^2}^\f12\\
&+\|(\tir\sn)^2 (\bb^{-1}) \bb^{-\f12} \aaa^{-\frac{\M}{2}}\|_{L_u^2 L_\omega^4}+\l t\r^{-\frac{1}{2}+2\delta}\Delta_0)\\
&\les \Delta_0^\f12 \log \l t\r^{-\frac{\M}{2}+1}(\log \l t\r^\frac{7}{4}(\La_0+\Delta_0^\frac{5}{4})^\f12\|\Omega^3 \Lb\varrho \aaa^{-\frac{\M}{2}}\tir\|_{L_u^2 L_\omega^2}^\f12+\log \l t\r^{\frac{7-\M}{2}}\Delta_0)\\
&+\Delta_0^\f12 \|\Omega^3(\bb^{-1})\aaa^{-\frac{\M}{2}}\|_{L_u^2 L_\omega^2}^\f12\log \l t\r^{\frac{7}{4}-\frac{\M}{2}}(\La_0+\Delta_0^\frac{5}{4})^\f12.
\end{align*} 
Using Lemma \ref{3.24.3.24}, we have
\begin{align*}
\sum_{i=0}^1&\|\tir^3 (\tir\sn)^i\sn\Xi_4(\tir\sn)^{2-i}(\bb^{-1}(\Xi_4+\tr\chi))\bb^\frac{1}{2}\aaa^{-\f12\M}\|_{L_u^2 L_\omega^2}\\
&\les \sup_{t'\in (0, t]} E_{-\M}[\Omega^3\varrho]^\frac{1}{4}(t')\Delta_0 \log \l t\r^{-\frac{\M}{2}+\frac{13}{4}}+\Delta_0^\frac{3}{2}\log \l t\r^{\frac{13}{2}-\M}.
\end{align*}
Therefore
\begin{align*}
\int_0^{t_1} \|\bb^\f12 \aaa^{-\frac{\M}{2}} I_1 \|_{L_u^2 L_\omega^2}&\les \int_0^{t_1}\log \l t\r^{-1}\Big( (C\M_0+1)^2 \sup_{t'\in(0, t]}(E^\f12_{-\M}[\Omega^3\varrho](t')\\
&+W_{1,-\M-1}^\f12[\Omega^3\Phi](t'))+\log \l t\r^{7-\frac{\M}{2}}(\La_0+\Delta_0^\frac{5}{4})\Big)  dt \\
&+\int_0^{t_1}\{\sup_{t'\in (0, t]} E_{-\M}[\Omega^3\varrho]^\frac{1}{4}(t')\Delta_0 (\log \l t\r)^{-\frac{\M}{2}+\frac{13}{4}}+\Delta_0^\frac{3}{2}\log \l t\r^{\frac{13}{2}-\M}\} dt\\
&+\l t_1\r^{\f12}\log \l t_1\r^\f12\sup_{ 0<t'\le t_1} W^\f12_{1,-\M-1}[\Omega^3\Phi](t').
\end{align*}
Summarizing the estimates of $I_1, I_2, I_3$ and $I_4$, due to $\M\ge 15$, we arrive at
\begin{align}\label{12.22.2.24}
\begin{split}
&\|\bb^\frac{3}{2}\aaa^{-\frac{\M}{2}} \tir^2(\tir \sn)^2\sF\c \Lb \varrho\|_{L_u^2 L_\omega^2}\\
&\les (C\M_0+1)^2 \sup_{t'\in(0, t]}(E_{-\M}[\Omega^3\varrho]^\f12(t')+W_{1,-\M-1}[\Omega^3\Phi]^\f12(t'))\log \l t\r^{-1}\l t\r^{-1}\\
&+\l t\r^{-1}\log \l t\r^{-1}(\La_0+\Delta_0^\frac{5}{4})+\l t\r^{-\frac{3}{2}}\log \l t\r^\f12\sup_{ 0<t'\le t} W_{1,-\M-1}[\Omega^3\Phi]^\f12(t').
\end{split}
\end{align}
Finally we consider the remaining term in (\ref{12.22.1.24}). It is direct to calculate
\begin{align*}
\bb^\frac{3}{2}\Omega^2(\bb^{-1}\Omega \Xi_4)=\bb^{-\frac{1}{2}}\Omega\bb \Omega^2\Xi_4+\bb^\frac{3}{2}\Omega^2(\bb^{-1})\Omega\Xi_4+\bb^{\f12}\Omega^3\Xi_4.
\end{align*}
We bound the first term on its right-hand side by using (\ref{1.27.5.24})
\begin{align*}
\|\Omega(\bb) \bb^{-\frac{1}{2}}\Omega^2\Xi_4\tir\|_{L_u^2 L_\omega^2}\les \log \l t\r\Delta_0\|\bb^{-\frac{1}{2}}\Omega^2\Xi_4\tir\|_{L_u^2 L_\omega^4}.
\end{align*}
For the second term, noting that symbolically, 
\begin{equation*}
\bb^\frac{3}{2}\Omega^2(\bb^{-1})=\Omega^2(\bb^\f12)+\bb^\f12 (\Omega\log \bb)^2,
\end{equation*}
we derive by using Lemma \ref{5.13.11.21} (5) that
\begin{align*}
\|\bb^\frac{3}{2}\Omega^2(\bb^{-1})\Omega \Xi_4 \tir\|_{L_u^2 L_\omega^2}\les \Delta_0^\f12\|\Omega^2(\bb^\f12)\|_{L_u^2 L_\omega^4}+(\log \l t\r\Delta_0)^2\|\Omega^{1+\le 1} \Xi_4 \tir\|_{L_u^2 L_\omega^4}.
\end{align*}
Hence
\begin{align*}
\|\bb^\frac{3}{2}\Omega^2(\bb^{-1}\Omega \Xi_4)\tir\|_{L_u^2 L_\omega^2}&\les \|\bb^\f12 \Omega^3\Xi_4 \tir\|_{L_u^2 L_\omega^2}+\Delta_0^\f12 \|\Omega^2(\bb^\f12)\|_{L_u^2 L_\omega^4}\\
&+(\log \l t\r)^2\Delta_0\|\Omega^{1+\le 1} \Xi_4 \tir\|_{L_u^2 L_\omega^4}.
\end{align*}
Using Lemma \ref{3.24.3.24} and Lemma \ref{10.10.3.23}, we derive
\begin{align*}
\|\bb^\frac{3}{2}&\Omega^2(\bb^{-1}\Omega \Xi_4)\tir \aaa^{-\frac{\M}{2}}\|_{L_u^2 L_\omega^2}\\
&\les \sup_{t'\in(0, t]}(C \M_0+1)E_{-\M}[\Omega^3\varrho]^\f12(t')+\int_0^t \l t'\r^{-\frac{3}{2}+\delta} W_{1,-\M-1}[\Omega^3\Phi]^\f12(t')\\
&+\l t\r^{-\f12}\log \l t\r^\f12 W_{1,-\M-1}[\Omega^3\Phi]^\f12(t)+\Delta_0^2\log \l t\r^\frac{15-\M}{2}+(\log \l t\r)^{-\f12\M}(\La_0+\Delta_0^\frac{5}{4}).
\end{align*}
In view of (\ref{11.11.2.23}), combining the above estimate with (\ref{12.22.2.24}), we conclude Lemma \ref{3.8.3.24}. 
\end{proof}
\begin{proposition}
\begin{align}\label{3.7.5.24}
\begin{split}
&\|\Omega^2\Box_\bg \bT\varrho\|_{L^2_\Sigma}\\
&\les (\|\fB\|_{L_x^\infty}+\Delta_0\l t\r^{-1}\log \l t\r^{-\frac{3}{2}})\{\|\Omega^2(\bb^{-1})\|_{L_u^2L_\omega^4}(1+\Delta_0\log \l t\r)\\
&+\|\bN\Omega^2\bT\varrho, \sn_\Omega^2\sn \bA_{g,1}\|_{L^2_\Sigma}\}+\|\bb^{-2}\Omega^2\fB\|_{L_u^2 L_\omega^4}+\l t\r^{-1}\log \l t\r^{\frac{\M}{2}+7}(\La_0+\Delta_0^\frac{5}{4}).
\end{split}
\end{align}
\end{proposition}
\begin{proof}
Differentiating (\ref{3.7.4.24}) and (\ref{2.19.3.24}) yields
\begin{align}\label{3.7.8.24}
\begin{split}
\Omega^2\Box_\bg \bT\varrho&=\Omega^2 \bT\Box_\bg \varrho+\Omega^2(\fB(\div_g k)_\bN)+\Omega^2(\sn\varrho (\div_g k)_A)+\Omega^2(\fB\ud \bA\c \sn \varrho)\\
&+\Omega^2(\fB(\bN\fB+\Box_\bg \varrho+\sn^2 \varrho+\theta \fB))+\Omega^2(\bA_{g,1}(\sn \fB+\theta\sn\varrho))
\end{split}
\end{align}
and
\begin{align}\label{3.7.7.24}
\begin{split}
\Omega^2\Big((\div_g k)_\bN\Big)&=\sn_\Omega^2(\fB^2+[\sn\Phi]^2+\bN \fB+\sn \bA_{g,1}+\tr\theta\fB+\ud \bA \bA_{g,1})\\
&=O(\l t\r^{-2}\log \l t\r(\La_0+\log \l t\r^\frac{5}{2}\Delta_0^\frac{5}{4}))_{L^2_u L_\omega^2}+\sn_\Omega^2(\sn\bA_{g,1}+\bN \fB)\\
\sn_\Omega^2(\div_g k)_A&=\Omega^2(\fB \bA_{g,1})+\sn_\Omega^2\sn \fB+\sn_\Omega^2 \sn \bAn+\sn_\Omega^2\Big((\tr\theta+\bA_g)\bA_{g,1}\Big)\\
&+\sn_\Omega^2 \sn_\bN \bA_{g,1}+\sn_\Omega^2((\fB+\bAn)\ud \bA)\\
&=O(\l t\r^{-2+\delta}\Delta_0)_{L^2_u L_\omega^2}
\end{split}
\end{align}
where we employed (\ref{9.12.2.22}), Lemma \ref{5.13.11.21} (5), Proposition \ref{8.29.8.21} and Proposition \ref{7.15.5.22} to derive the above two estimates. 

Recall from (\ref{2.19.2.24}) and (\ref{2.19.3.24}) for estimates of $\sn_\Omega^{\le 1}(\div_g k)_A$ and $\Omega^{\le 1}((\div_g k)_\bN)$.
\begin{align}\label{3.7.6.24}
\begin{split}
(\div_g k)_A&=O(\l t\r^{-2+\delta}\Delta_0), (\div_g k)_\bN=\bN \fB+O(\tir^{-1})\fB\\\
\sn_\Omega(\div_g k)_A&=O(\l t\r^{-1+\delta}\Delta_0)_{L^2_\Sigma}, \Omega \bN \fB+\bb^{-1}(\l t\r^{-2}\log \l t\r\Delta_0)_{L_\omega^4}\\
 \sn_\Omega(\div_g k)_\bN&=O(\log \l t\r^{\f12(\M+7)}(\La_0+\Delta_0^\frac{5}{4}))_{L^2_\Sigma}
\end{split}
\end{align}
where we also applied (\ref{2.13.3.24}) to derive the last estimates, and employed (\ref{8.23.2.23}) and Proposition \ref{7.15.5.22} to derive the $L^4_\omega$ estimate.

Using Proposition \ref{7.15.5.22},  (\ref{3.7.7.24}), (\ref{3.7.6.24})  and Sobolev embedding,
\begin{align*}
\sum_{a=0}^2(\sn_\Omega^a \sn \varrho \sn_\Omega^{2-a}&((\div_g k)_A))=O(\l t\r^{-4+2\delta}\Delta_0^2)_{L^2_u L_\omega^2}\\
\Omega^2(\fB\ud\bA\c \bA_{g,1})&=O(\l t\r^{-4+3\delta}\Delta_0^2)_{L^2_u L_\omega^2}, \Omega^2(\bA_{g,1}\sn \fB)=O(\l t\r^{-4+2\delta}\Delta_0^2)_{L^2_u L_\omega^2}\\
\Omega^2(\fB\sn^2\varrho)&=\fB \sn_\Omega^2 \sn^2\varrho+O(\l t\r^{-4+2\delta}\Delta_0^2)_{L_u^2 L_\omega^2}\\
\Omega^2(\fB^2\theta)&=\fB \Omega^2\fB \theta+O(\tir^{-1})(\Omega\fB)^2+\fB^2\Omega^2(\thetac+\fB)+\fB \Omega\fB \Omega\bA_b\\
&=O(\l t\r^{-1})\fB \Omega^2\fB+O(\l t\r^{-3}\log \l t\r^2\Delta_0^2)_{L_\omega^2}\\
&+O(\fB)O(\l t\r^{-2}\log \l t\r^\frac{7}{2}(\La_0+\Delta_0^\frac{5}{4}))_{L^2_u L_\omega^2}\\
&=\bb^{-1}O(\l t\r^{-3}) \log \l t\r^\frac{7}{2}(\La_0+\Delta_0^\frac{5}{4})_{L^2_u L_\omega^2}
\end{align*}
where we used (\ref{3.7.9.24}), (\ref{1.27.2.24}) and (\ref{2.14.1.24}). 

Substituting the above estimates to (\ref{3.7.8.24}) gives 
\begin{align}
\Omega^2\Box_\bg \bT\varrho&=\Omega^2 \bT\Box_\bg \varrho+\fB\sn_\Omega^2(\bN \fB+\sn\bA_{g,1})+O(\l t\r^{-1})_{L_\omega^4}\bb^{-3}\Omega^2 \fB\nn\\
&+\sum_{a=1}^2\Omega^a \fB \Omega^{2-a}((\div_g k)_\bN)+O(\l t\r^{-2}\log \l t\r^\frac{7}{2}(\La_0+\Delta_0^\frac{5}{4}))_{L^2_\Sigma}\label{3.8.2.24}
\end{align}
 where other error terms including $\Omega^2(\fB \Box_\bg\varrho)$ in (\ref{3.7.8.24})  are included in the last term in the above. 
  
 Using (\ref{1.29.2.22}), (\ref{L2BA2}), (\ref{3.7.7.24}), (\ref{3.7.6.24}) and Sobolev embedding, we deduce
\begin{align*}
\sum_{a=1}^2\Omega^a \fB \Omega^{2-a}((\div_g k)_\bN)&=\bb^{-3}\Omega^2\fB O(\l t\r^{-1})_{L_\omega^4}+\bb^{-1}O(\l t \r^{-3}\log \l t\r^2\Delta_0^2)_{L_\omega^2}\\
&+\Omega\fB\c O\Big(\l t\r^{-\f12}\log \l t\r^{\frac{\M+7}{4}}(\Delta_0^\frac{5}{4}+\La_0)^\f12\|\Omega^2 \bN \fB\|^\f12_{L_u^2 L_\omega^2}\Big)_{L^2_u L_\omega^4}.
\end{align*}
Hence, due to Lemma \ref{5.13.11.21} (5), and Cauchy-Schwarz
\begin{align*}
\|\bb^\f12\sum_{a=1}^2&\Omega^a \fB \Omega^{2-a}((\div_g k)_\bN)\|_{L_u^2 L_\omega^2}\\
\displaybreak[0]
&\les\l t\r^{-2}\log \l t\r^{\frac{\M}{2}+7}(\Delta_0^\frac{5}{2}+\La_0)\Delta_0+\l t\r^{-1}\|\bb^{-\frac{5}{2}}\Omega^2 \fB\|_{L_u^2 L_\omega^4}\\
&+\Delta_0\l t\r^{-1}\log \l t\r^{-\frac{3}{2}}\|\Omega^2\bN \fB\|_{L_u^2 L_\omega^2}.
\end{align*}

Combining (\ref{3.8.2.24}), the above two estimates, (\ref{3.7.10.24}), (\ref{3.8.1.24}), (\ref{10.10.2.23}) and (\ref{8.29.9.21}) gives
\begin{align*}
&\|\Omega^2\Box_\bg \bT\varrho\|_{L^2_\Sigma}\\
&\les (\|\fB\|_{L_x^\infty}+\Delta_0\l t\r^{-1}\log \l t\r^{-\frac{3}{2}})\{\|\Omega^2(\bb^{-1})\|_{L_u^2L_\omega^4}(1+\log \l t\r\Delta_0)+\|\bN\Omega^2\bT\varrho, \sn_\Omega^2\sn \bA_{g,1}\|_{L^2_\Sigma}\}\\
&+\|\bb^{-2}\Omega^2\fB\|_{L_u^2 L_\omega^4}+\l t\r^{-1}\log \l t\r^{\frac{\M}{2}+7}(\La_0+\Delta_0^\frac{5}{4})
\end{align*}
as stated in (\ref{3.7.5.24}).
\end{proof}
\subsection{The control of the energies up to the top order}\label{4.2.1.23}
\begin{proposition}[Top order energy estimates]\label{10.30.4.21}
Let $X\in \{\Omega, S\}$ and $u_0\le u\le u_*$ and $0<t<T_*$. Under the assumptions of (\ref{3.12.1.21})-(\ref{6.5.1.21}), there hold

(1) \begin{align*}
&E_{-\M}[\Omega^3 \varrho](t)+F_{0,-\M}[\Omega^3\varrho](\H_u^t)+W_{1,-\M-1}[\Omega^3\Phi](t)+F_{1,-\M-1}[\Omega^3 \Phi](\H_u^t)\\
&\qquad\qquad\qquad\qquad\qquad\qquad\les (\log \l t\r)^{2\M^2-2\M-1}(\La_0^2+\Delta_0^\frac{5}{2})\\
&W_1[X^3\Phi](t)+F_1[X^3\Phi](\H_u^t)\les \log \l t\r^{\M}(\La_0^2+\Delta_0^\frac{5}{2}), \mbox{ if }\vs^+(X^3)=1\\
&WFIL_2[\Omega^2 \Phi](\D_u^t)\les (\log \l t\r)^\M(\La_0^2+(\log \l t\r)^{2\M^2-\M-1}\Delta_0^\frac{5}{2}).
\end{align*}
\begin{align}\label{3.10.1.24}
\begin{split}
\|\Sc(\sta{X^3, S}\Phi), \Ac(\sta{X^3,S} \Phi)\|^2_{L^2_\Sigma}&\les \l t\r\log \l t\r^\M(\La_0^2+\Delta_0^\frac{5}{2}),\\
\|\Sc(\Omega^4\Phi), \Ac(\Omega^4\Phi)\|^2_{L^2_\Sigma}&\les \l t\r (\log \l t\r)^{2\M^2}(\La_0^2+\Delta_0^\frac{5}{2}).
\end{split}
\end{align}
(2) There hold, with $0\le \sig\le \frac{3}{2}$ and $\sig\neq 1$,  that
\begin{align}\label{3.9.7.24}
\begin{split}
&\|\Omega^3(\bb^{-1+\sig}), \Omega^3\log \bb\|_{L^2_u L_\omega^2}\les \log \l t\r^{\M^2-\f12\M+\f12}(\La_0+\Delta_0^\frac{5}{4})\\
&\|\Omega^2(\bb^{-1+\sig}), \log \l t\r^{-1}\Omega^2\log \bb\|_{L^2_u L_\omega^4}\les\log \l t\r^{\f12\M^2-3}(\La_0+\Delta_0^\frac{5}{4})
\end{split}
\end{align}
\begin{align}\label{3.9.9.24}
\begin{split}
&\|\tir(\bb^{-1+\sig}\Omega^3\fB, \Omega^3(\bb^{-1+\sig}\fB))\|_{L^2_u L_\omega^2}\les \log \l t\r^{\M^2-\f12\M}(\La_0+\Delta_0^\frac{5}{4})\\
&\|\Omega^2 \fB \tir\|_{L_u^2 L_\omega^4}\les \log \l t\r^{\f12\M^2-2}(\La_0+\Delta_0^\frac{5}{4});
\end{split}
\end{align}

\begin{align}\label{3.9.10.24}
E[\Omega^2 \bT\varrho]^\f12(t)+F_0[\Omega^2 \bT\varrho]^\f12(\H_u^t)\les \log \l t\r^{\f12\M^2-1}\Big(\La_0+\Delta_0^\frac{5}{4}\Big).
\end{align}
\end{proposition}
In the sequel, we will frequently use Lemma \ref{3.17.2.22} without explicit mentioning.
\begin{proof}
{\bf $\bullet$ Proof of (1)}
Recall from (\ref{3.29.2.23'}), we write
\begin{align*}
\Er_3(\Phi, X_3 X_2 X_1)&=X_3 X_2(\fm{X_1}\Box_\bg \Phi+\sP[X_1, \Phi])\\
&+X_3(\fm{X_2}\Box_\bg X_1 \Phi+\sP[X_2, X_1 \Phi])\\
&+\fm{X_3}\Box_\bg X_2 X_1 \Phi+\sP[X_3, X_2 X_1\Phi].
\end{align*}

Recall that $\fm{\Omega}=\Omega\log \bb$ and $\fm{S}= 1+\tir k_{\bN\bN}$. Using Proposition \ref{8.29.8.21}, Proposition \ref{7.15.5.22}, Lemma \ref{5.13.11.21} and (\ref{10.30.1.21}) we have
\begin{equation}\label{11.13.1.21}
\begin{split}
& |\fm{S}|\les 1, |\fm{\Omega}|+\|\Omega \fm{\Omega}\|_{L_\omega^4}+\|\Omega^2\fm{\Omega}\|_{L_u^2 L_\omega^2}\les \l t\r^\delta\Delta_0\\
&\|Y^2\fm{X}\|_{L^2_u L_\omega^2}\les \log\l t\r(\La_0+\log \l t\r^3\Delta_0^\frac{5}{4}), \mbox{ if } \vs^-(Y^2X)=0 \mbox{ and }\vs^+(Y^2X)=1\\
&S^2 \fm{S}=O(1)\tir k_{\bN\bN}+O(\l t\r^{\delta}\Delta_0)_{L^2_\Sigma}\\
&|Y\fm{X}|\les \l t\r^{\delta}\Delta_0, \mbox{ if }\vs^-(YX)=0 \mbox{ and } \vs^+(YX)=1\\
& Y\fm{X}=O(\tir \fB) \mbox{ if }\vs^-(YX)=1.
\end{split}
\end{equation}
We divide our analysis into four cases according to various combinations of $X_2, X_3$ in (\ref{3.29.2.23'}).

{\bf Step 1: Let $X_2, X_3=\Omega$.}  We first consider the terms
\begin{align*}
I[\Omega^2 X_1\Phi]=\sum_{\ell=0}^2\Omega^\ell\fm{X_1}\Omega^{2-\ell}\Box_\bg \Phi+\Omega(\fm{\Omega}\Box_\bg X_1\Phi)+\fm{\Omega}\Box_\bg \Omega X_1\Phi.
\end{align*}
Using (\ref{11.13.1.21}), (\ref{4.3.3.21}), (\ref{8.30.3.21}), (\ref{8.30.3.21+}) and Proposition \ref{1steng} (1), we have
\begin{equation}\label{8.20.1.22}
\|I[\Omega^2 X_1\Phi]\|_{L^2_\Sigma}\les \l t\r^{\f12(1-\vs(X_1))+\delta}
\Delta_0\cdot\left\{\begin{array}{lll}
\l t\r^{-2}\La_0+\l t\r^{-2+\delta}\Delta_0^\frac{5}{4}, X_1=\Omega\\
(\l t\r^{-2}\log \l t\r^{\f12\M+1}\La_0+\l t\r^{-\frac{7}{4}+2\delta}\Delta_0^\frac{5}{4}), X_1=S.
\end{array}\right.
\end{equation}

For the remaining terms, we will divide the analysis into two cases.

\noindent{\bf Case 1  $X_1=\Omega$.} As seen in Proposition \ref{7.16.2.21}, we  decompose
\begin{equation*}
\Omega^2(\sP[\Omega, \Phi])=\Omega^3\tr\chi\Lb \Phi+\big(\Omega^2(\sP[\Omega, \Phi])\big)_g.
\end{equation*}
In view of the above decomposition and (\ref{3.29.2.23'}), symbolically we can write
\begin{equation*}
\Er_3(\Phi, \Omega^3)=\Er_3(\Phi, \Omega^3)_g+\Omega^3\tr\chi\Lb \Phi
\end{equation*}
with
\begin{align*}
\Er_3(\Phi, \Omega^3)_g&=(\Omega^2\sP[{}\rp{a}\Omega, \Phi])_g+\sP[\Omega, \Omega^2 \Phi]+\Omega\sP[\Omega, \Omega \Phi]+I[\Omega^3\Phi].
\end{align*}
The commutators have been treated in (\ref{10.1.1.21}) and Proposition \ref{7.16.2.21}. Hence, also using (\ref{8.20.1.22}),  we conclude
\begin{equation}\label{10.7.2.21}
\|\Er_3(\Phi, \Omega^3)_g\|_{L^2_\Sigma}\les \l t\r^{-\frac{3}{2}}W_1[\Omega^3\Phi]^\f12(t)+\l t\r^{-\frac{3}{2}+2\delta}\Delta_0^\frac{5}{4}\log \l t\r^\f12+\l t\r^{-2}\log \l t\r^{\f12\M+1}\La_0.
\end{equation}
For the remaining  terms which are the most crucial ones, we treat the following two terms separately
\begin{align*}
\fY_1=\int_{\D_{u_1}^{t_1}}&|\Omega^3\tr\chi \Lb\varrho \bT\Omega^3\varrho| \aaa^{-\M} d\mu_g dt\\
\fY_2=\int_{\D_{u_1}^{t_1}}&|\Omega^3\tr\chi \Lb\Phi (L+h)\Omega^3\Phi| \tir\aaa^{-\M-1} d\mu_g dt.
\end{align*}
In view of (\ref{10.7.2.21}), by using Cauchy-Schwarz inequality, we derive 
 \begin{align}\label{3.8.6.24}
 \begin{split}
 \int_{\D_{u_1}^{t_1}}|\Er_3(\varrho, \Omega^3)\bT\Omega^3\varrho| \aaa^{-\M}d\mu_g dt &\les\La_0^2+\Delta_0^\frac{5}{2}+ \fY_1+\int_0^t \l t'\r^{-1-\delta}E_{-\M}[\Omega^3\varrho](t') dt'\\
 &+\int_0^t \l t'\r^{-2+\delta}\log \l t'\r^{\M+1} W_{1,-\M-1}[\Omega^3\Phi](t')dt',
 \end{split}
 \end{align}
  \begin{align}\label{3.8.5.24}
  \begin{split}
\int_{\D_{u_1}^{t_1}}|\Er_3(\Phi, \Omega^3) &(L+h)\Omega^3\Phi |\tir \aaa^{-\M-1} d\mu_g dt\\
&\les \fY_2+\int_0^t \l t'\r^{-2}\log \l t'\r^{\M+1} W_{1,-\M-1}[\Omega^3\Phi](t') dt'\\
&+\Delta_0^\frac{5}{2}+\La_0^2+\int_{u_1}^{u_*} F_{1,-\M-1}[\Omega^3\Phi](\H_u^t) du.
\end{split}
  \end{align}
For $\fY_1$ and $\fY_2$, we derive by using Lemma \ref{3.8.3.24} that
\begin{align*}
\fY_1&\les \int_0^{t_1}\{ (C\M_0+1)^2 \sup_{t'\in(0, t]}(E_{-\M}[\Omega^3\varrho]^\f12(t')+W_{1,-\M-1}[\Omega^3\Phi]^\f12(t'))\log \l t\r^{-1}\l t\r^{-1}\\
&+\l t\r^{-1}(\log \l t\r)^{-1}(\La_0+\Delta_0^\frac{5}{4})+\l t\r^{-\frac{3}{2}}(\log \l t\r)^\f12\sup_{ 0<t'\le t} W_{1,-\M-1}[\Omega^3\Phi]^\f12(t')\}\\
&\times E_{-\M}[\Omega^3\varrho]^\f12(t) dt,
\end{align*}
  \begin{align*}
  \fY_2&\les \int_0^{t_1} \l t'\r^{\f12}\|\bb^\f12\Omega^3\tr\chi\tir\Lb \Phi \aaa^{-\frac{\M}{2}}\|_{L_u^2 L_\omega^2}\|\tir^\frac{3}{2}\aaa^{-\f12(\M+1)}(L+h)\Omega^3\Phi\|_{L_u^2 L_\omega^2} dt'\\
  \displaybreak[0]
&\les \int_0^{t_1} \l t'\r^{-1}\log \l t'\r^{-2}(C\M_0+1)^4\sup_{t''\le t'} \Big(E_{-\M}[\Omega^3\varrho](t'')+W_{1,-\M-1}[\Omega^3\Phi](t'')\Big) dt'\\
&+\int_{u_1}^{u_*}  F_{1,-\M-1}[\Omega^3\Phi](\H_u^{t_1}) du+\La_0^2+\Delta_0^\frac{5}{2}
  \end{align*}  
  where we used $1\les\bb\les \log \l t\r$. 
Substituting (\ref{3.8.6.24}), (\ref{3.8.5.24}) to (\ref{6.21.2.21})  and Proposition \ref{MA2} leads to 
\begin{align*}
&E_{-\M}[\Omega^3\varrho](t_1)+ F_{0,-\M}[\Omega^3\varrho](\H_{u_1}^{t_1})+F_{1,-\M-1}[\Omega^3\Phi](\H_{u_1}^{t_1})+W_{1,-\M-1}[\Omega^3\Phi](t_1, [u_1,u_*])
\\
&\les\int_{u_1}^{u_*}(F_{0,-\M}[\Omega^3\varrho](\H_u^{t_1})+F_{1,-\M-1}[\Omega^3\Phi](\H_u^{t_1})) du+E_{-\M}[\Omega^3\varrho](0)+F_{0,-\M}[\Omega^3\varrho](\H_{u_*}^{t_1})\\
\displaybreak[0]
&+W_{1,-\M-1}[\Omega^3\Phi](0, [u_1,u_*])+F_{1,-\M-1}[\Omega^3\Phi](\H_{u_*}^{t_1})+\fY_1+\fY_2+\La_0^2+\Delta_0^\frac{5}{2}\\
&+\int_0^{t_1} \l t'\r^{-1-\delta}(W_{1,-\M-1}[\Omega^3\Phi](t', [u_1,u_*])+E_{-\M}[\Omega^3\varrho](t'))dt'.\end{align*}

It follows by using the estimates on $\fY_1$, $\fY_2$ and Gronwall's inequality that
\begin{align}\label{3.9.1.24}
E_{-\M}[\Omega^3\varrho](t_1)+ F_{0,-\M}[\Omega^3\varrho](\H_{u_1}^{t_1})&+F_{1,-\M-1}[\Omega^3\Phi](\H_{u_1}^{t_1})+W_{1,-\M-1}[\Omega^3\Phi](t_1, [u_1,u_*])\nn\\
&\les (\log \l t_1\r)^{2\M(\M-1)-1} (\La_0^2+\Delta_0^\frac{5}{2}),
\end{align}
where $\M \ge 15$, being a constant multiple of $\M_0$, chosen similar as in Lemma \ref{1.6.4.18}.

\noindent{\bf Case 2 $X_1=S$.}
 \begin{align*}
 \Er_3(\Phi,\Omega^2 S)=\Omega^2\sP[S, \Phi]+\Omega\sP[\Omega, S\Phi]+\sP[\Omega, \Omega S\Phi]+I[\Omega^2 S\Phi].
 \end{align*}
 Using (\ref{8.20.1.22}), (\ref{10.1.1.21}) and (\ref{8.17.5.21}), we bound
 \begin{align*}
 \|\Er_3(\Phi, \Omega^2S)\|_{L^2_\Sigma}&\les \l t\r^{-\frac{3}{2}+2\delta}\Delta_0^\frac{5}{4}\log \l t\r^\f12+\l t\r^{-2}\log \l t\r^{\f12\M+7}\La_0\\
 &+ \l t\r^{-\frac{3}{2}}\{W_1[\Omega^2 S\Phi]^\f12(t)+W_1[\Omega^3\Phi]^\f12(t)\}
 \end{align*}
Applying Proposition \ref{MA2} to $\psi=\Omega^2 S\Phi$, using (\ref{3.9.1.24}) and Gronwall's inequality, we have
\begin{align}\label{3.9.2.24}
W_1[\Omega^2 S\Phi](t)+F_1[\Omega^2S\Phi](\H_u^t)\les \log \l t\r^{\M}(\La_0^2+\Delta_0^\frac{5}{2}). 
\end{align}

{\bf Step 2: Let $X_3=\Omega, X_2=S$.} Recall from (\ref{3.29.2.23'})
\begin{align*}
I[\Omega S X_1\Phi]&=\Omega S(\fm{X_1}\Box_\bg \Phi)+\Omega(\fm{S}\Box_\bg X_1\Phi)+\fm{\Omega}\Box_\bg S X_1\Phi\\
&=\Omega S\fm{X_1}\Box_\bg \Phi+\fm{X_1}\Omega S\Box_\bg \Phi+\Omega \fm{X_1}S\Box_\bg \Phi+S\fm{X_1}\Omega \Box_\bg\phi\\
&+\Omega \fm{S}\Box_\bg X_1\Phi+\fm{S}\Omega\Box_\bg X_1\Phi+\fm{\Omega}\Box_\bg S X_1\Phi.
\end{align*}
Using (\ref{11.13.1.21}), (\ref{3.29.1.23}), Sobolev embedding, (\ref{8.30.3.21}) and (\ref{8.30.3.21+}) we infer
\begin{align*}
\|I[\Omega S X_1\Phi]\|_{L^2_\Sigma}&\les (1+\l t\r^\delta\Delta_0)\|\Omega S^{\le 1}\Box_\bg \Phi, \Omega^{\le 1}\Box_\bg X_1\Phi, \Box_\bg S X_1\Phi\|_{L^2_\Sigma}\\
&+\l t\r^{\delta}\Delta_0\|\bb^\f12 \tir S\Box_\bg \Phi\|_{L^2_u L_\omega^4}+\|\Omega S\fm{X_1}\|_{L^2_\Sigma}\|\Box_\bg \Phi\|_{L^\infty_x}\\
&\les
 \l t\r^{\f12(1-\vs(X_1))+\delta}
\cdot\left\{\begin{array}{lll}
(\l t\r^{-2}\La_0+\l t\r^{-2+\delta}\Delta_0^\frac{5}{4}), X_1=\Omega\\
(\l t\r^{-2}\log \l t\r^{\f12\M+1}\La_0+\l t\r^{-\frac{7}{4}+2\delta}\Delta_0^\frac{5}{4}), X_1=S.
\end{array} \right.
\end{align*}
Next we consider the other terms in (\ref{3.29.2.23'}). We first apply (\ref{3.5.1.24}) to $X^2=\Omega X_1$, and Proposition \ref{8.29.8.21}  to obtain
\begin{align*}
\|\Omega(\sP[S, X_1 \Phi])\|_{L^2_\Sigma}&\les \l t\r^{-\frac{3}{2}}(W_1[\Omega S X_1 \Phi]^\f12(t)+\l t\r^\delta W_1[\Omega^2 X_1\Phi]^\f12(t))\\
&+\l t\r^{-2}(\log \l t\r)^{\frac{\M}{2}+1}\La_0+\l t\r^{-\frac{7}{4}+3\delta}\Delta_0^\frac{5}{4}(\log \l t\r)^\f12. 
\end{align*}
Applying (\ref{10.1.1.21}) to $f=S X_1\Phi$ gives
 \begin{align*}
 \|\sP [\Omega,S X_1\Phi]\|_{L^2_\Sigma} &\les \l t\r^{-\frac{3}{2}+2\delta}\Delta_0^\frac{5}{4}+\l t\r^{-2}\log \l t\r^{\f12\M+1}\La_0.
\end{align*}
Now we apply Proposition \ref{7.16.2.21} to $X^2 =\Omega S$ to derive
\begin{align*}
\sum_{X_1=S, \Omega}\|\Omega S \sP[X_1, \Phi])\|_{L^2_\Sigma}&\les \l t\r^{-\frac{3}{2}+2\delta}\Delta_0^\frac{5}{4}\log \l t\r^\f12+\l t\r^{-2}\log \l t\r^{\f12(\M+7)}\La_0\\
&+\l t\r^{-\frac{3}{2}}(W_1[\Omega S\Omega\Phi, \Omega S^2\Phi, \Omega^2 S\Phi, \Omega^3\Phi]^\f12(t)).
\end{align*}
Using (\ref{3.9.1.24}) and (\ref{3.9.2.24}), combining the above estimates, we deduce
\begin{align*}
\sum_{X_1=S, \Omega}\|\Er_3(\Phi, \Omega S X_1)\|_{L^2_\Sigma}\les \l t\r^{-\frac{3}{2}}(\l t\r^{2\delta}\Delta_0^\frac{5}{4}\log \l t\r^\f12+\La_0)+\l t\r^{-\frac{3}{2}}\sum_{X_1=S, \Omega}W_1[\Omega SX_1\Phi]^\f12(t).
\end{align*}
Applying Proposition \ref{MA2} to $\psi=\Omega SX_1\Phi$ with the help of the above estimate leads to 
\begin{align}\label{3.9.4.24}
\sum_{X_1=S, \Omega}&(W_1[\Omega S X_1\Phi](t)+F_1[\Omega S X_1\Phi](\H_u^t))\les \log \l t\r^\M(\La_0^2+\Delta_0^\frac{5}{2}).
\end{align}

{\bf Step 3: Let $X_3=S,\, X_2=\Omega$.} Let us first consider the error term in (\ref{3.29.2.23'})
\begin{align*}
I[S\Omega X_1\Phi]=S\Omega(\fm{X_1}\Box_\bg\Phi)+S(\fm{\Omega}\Box_\bg X_1\Phi)+\fm{S}\Box_\bg \Omega X_1\Phi.
\end{align*}
Treating the error term similar to the previous case,  we derive
\begin{align*}
\|I[S\Omega X_1\Phi]\|_{L^2_\Sigma}&\les\|S\Omega\fm{X_1}\Box_\bg \Phi,\fm{X_1}S\Omega \Box_\bg\Phi, S\fm{X_1}\Omega \Box_\bg \Phi, \Omega\fm{X_1}S\Box_\bg\Phi\|_{L^2_\Sigma}\\
&+\|S\fm{\Omega}\Box_\bg X_1\Phi\|_{L^2_\Sigma}+\|\fm{\Omega}S\Box_\bg X_1\Phi\|_{L^2_\Sigma}+\|\fm{S}\Box_\bg \Omega X_1\Phi\|_{L^2_\Sigma}\\
&\les  \l t\r^{\f12(1-\vs(X_1))+\delta}\cdot
\left\{\begin{array}{lll}
\l t\r^{-2}\La_0+\l t\r^{-2+\delta}\Delta_0^\frac{5}{4}, X_1=\Omega\\
(\l t\r^{-2}\log \l t\r^{\f12\M+1}\La_0+\l t\r^{-\frac{7}{4}+2\delta}\Delta_0^\frac{5}{4}), X_1=S.
\end{array}\right.
\end{align*}
Recall from  (\ref{3.29.2.23'}) that
\begin{align*}
\Er_3(\Phi, S\Omega X_1)&=S\Omega(\sP[X_1, \Phi])+S(\sP[\Omega, X_1\Phi])+\sP[S, \Omega X_1\Phi]+I[S\Omega X_1\Phi].
\end{align*}
It remains to control the other three terms. Applying  (\ref{3.5.1.24}) to $f=\Omega X_1\Phi$ yields
\begin{align*} 
\|\sP[S, \Omega X_1\Phi]\|_{L^2_\Sigma}&\les \l t\r^{-\frac{3}{2}}\{ W_1[S\Omega X_1\Phi]^\f12(t)+\l t\r^{\delta}( W_1[\Omega^2 X_1\Phi]^\f12(t))\}\nn\\
 &+\l t\r^{-2}\log \l t\r^{\f12\M+1}\La_0+\Delta^\frac{5}{4}\l t\r^{-\frac{7}{4}+3\delta}(\log \l t\r)^\f12.
\end{align*}
It follows by applying (\ref{3.9.1.24}) and (\ref{3.9.2.24}) to the above that 
\begin{align*}
\|\sP[S, \Omega X_1\Phi]\|_{L^2_\Sigma}\les \l t\r^{-\frac{3}{2}}W_1[S\Omega X_1\Phi]^\f12(t)+ \l t\r^{-\frac{3}{2}+\delta}\log \l t\r^{\M^2}(\Delta_0^\frac{5}{4}+\La_0).
\end{align*}
Recall from (\ref{10.1.1.21}) that
\begin{align*} 
\|S(\sP[\Omega, X_1\Phi])\|_{L^2_\Sigma}&\les\l t\r^{-\frac{3}{2}+2\delta}\Delta_0^\frac{5}{4}\log \l t\r^\f12+\l t\r^{-2}\log \l t\r^{\f12\M+1}\La_0.
\end{align*}
Using Proposition \ref{7.16.2.21}, we derive
\begin{align*}
\sum_{X_1=S, \Omega}\|S{}\Omega\sP[X_1, \Phi]\|_{L^2_\Sigma}&\les \l t\r^{-\frac{3}{2}}\sum_{X_1=S, \Omega}(W_1[S\Omega X_1\Phi]^\f12(t)+W_1[\Omega^2 X_1\Phi]^\f12(t))\nn\\
&+\l t\r^{-\frac{3}{2}+2\delta}\Delta_0^\frac{5}{4}(\log \l t\r)^\f12+\l t\r^{-2}\log \l t\r^{\f12(\M+7)}\La_0.
\end{align*} 
Using  the above three estimates and the estimate of $I[S\Omega X_1\Phi]$, we obtain
\begin{align*}
\sum_{X_1=\Omega, S}\|\Er_3(\Phi, S\Omega X_1)\|_{L^2_\Sigma}&\les \l t\r^{-\frac{3}{2}}\sum_{X_1=S, \Omega}\{W_1[S\Omega X_1\Phi]^\f12(t)+\l t\r^\delta W_1[\Omega^2 X_1\Phi]^\f12(t)\}\nn\\
&+ \l t\r^{-\frac{3}{2}+\delta}(\l t\r^\delta\Delta_0^\frac{5}{4}\log \l t\r^\f12+(\log \l t\r)^{\M^2}\La_0).
\end{align*}
Consequently, in view of (\ref{3.9.1.24}) and (\ref{3.9.2.24}),
applying Proposition \ref{MA2} to $\psi=S\Omega X_1\Phi$, with $X_1=\Omega, S$ and the above estimate yields
\begin{align}\label{3.9.5.24}
W_1[S\Omega X_1\Phi](t)+F_1[S\Omega X_1\Phi](\H_u^t)\les \log \l t\r^\M(\La_0^2+\Delta_0^\frac{5}{2}).
\end{align}

{\bf Step 4: Let $X_3=X_2=S$.}
Recall from (\ref{3.29.2.23'})
\begin{align*}
I[S^2 X_1\Phi]=S^2(\fm{X_1}\Box_\bg \Phi)+S(\fm{S}\Box_\bg X_1\Phi)+\fm{S}\Box_\bg S X_1\Phi.
\end{align*}
Using (\ref{4.3.3.21}), (\ref{8.26.4.21}) and (\ref{11.13.1.21}) we infer
\begin{align*}
\|S^2\fm{X_1}\c \Box_\bg\Phi\|_{L^2_\Sigma}&\les \|\Box_\bg \Phi\|_{L^2_\Sigma}+\l t\r^{1+\delta}\Delta_0\|\bb^\f12\Box_\bg \Phi\|_{L^\infty_x}\\
&\les \l t\r^{-2}(\log \l t\r^{\f12\M}\La_0+\l t\r^{2\delta}\Delta_0^\frac{5}{4}).
\end{align*}
In view of (\ref{11.13.1.21}), (\ref{8.30.3.21}) and (\ref{8.30.3.21+}), we bound
\begin{align*}
\|I[S^2 X_1\Phi]\|_{L^2_\Sigma}&\les \|S^2\fm{X_1}\c \Box_\bg\Phi\|_{L^2_\Sigma}+\sum_{a=0}^1\|S^a\fm{X_1}\|_{L^\infty_\Sigma}\|S^{2-a}\Box_\bg \Phi\|_{L^2_\Sigma}\\
&+\|S^{\le 1}\Box_\bg X_1\Phi\|_{L^2_\Sigma}+\|\Box_\bg S X_1\Phi\|_{L^2_\Sigma}\\
&\les \l t\r^{-2}(\log \l t\r^{\f12\M+1}\La_0+\l t\r^{\frac{1}{4}+2\delta-\vs(X_1)\frac{1}{4}}\Delta_0^\frac{5}{4}).
\end{align*}
Recall from (\ref{3.29.2.23'}) that
\begin{align*}
\Er_3(\Phi, S^2X_1)=I[S^2 X_1\Phi]+S^2(\sP[X_1, \Phi])+S(\sP[S, X_1\Phi])+\sP[S, S X_1\Phi].
\end{align*}
Using Proposition \ref{7.16.2.21} we have
\begin{align*}
\|S^2(\sP[X_1, \Phi])\|_{L^2_\Sigma}&\les \l t\r^{-\frac{3}{2}} (\sum_{\vs(Y)\le \vs(X_1)}W_1[S^2 Y \Phi]^\f12(t)+W_1[SX_1 S\Phi]^\f12(t)+ W_1[X_1\Omega^2\Phi]^\f12(t))\\
&+ \l t\r^{-\frac{3}{2}+2\delta}(\log \l t\r)^\f12\Delta_0^\frac{5}{4}+\l t\r^{-2}\log \l t\r^{\f12(\M+7)}\La_0.
\end{align*}
Applying(\ref{3.5.1.24}) to $X_3 X_2=S X_1$, we infer
\begin{align*}
\|S(\sP[S, X_1\Phi])\|_{L_\Sigma^2}&\les \l t\r^{-\frac{3}{2}}(W_1[S^2 X_1\Phi]^\f12(t)+W_1[\Omega^2 X_1\Phi]^\f12(t)+\l t\r^\delta W_1[S\Omega X_1\Phi]^\f12(t))\\
&+\l t\r^{-2 }\log \l t\r^{\f12\M+1}\La_0+\Delta_0^\frac{5}{4}\l t\r^{-\frac{7}{4}+3\delta}(\log \l t\r)^\f12.
\end{align*}
Next applying (\ref{3.5.1.24}) to $f=S X_1\Phi$ yields
\begin{align*}
 \|\sP[S, SX_1\Phi]\|_{L^2_\Sigma}&\les  \l t\r^{-\frac{3}{2}} (W_1[\Omega S X_1\Phi]^\f12(t)+W_1[S^2X_1\Phi]^\f12(t))\\
 &+\l t\r^{-2 }\log \l t\r^{\f12\M+1}\La_0+\Delta_0^\frac{5}{4}\l t\r^{-\frac{7}{4}+3\delta}(\log \l t\r)^\f12.
\end{align*}
Combining the above estimates, (\ref{3.9.1.24})-(\ref{3.9.5.24}), we conclude
\begin{align*}
\sum_{X_1=S, \Omega}\|\Er_3(\Phi, S^2 X_1)\|_{L^2_\Sigma}&\les \l t\r^{-\frac{3}{2}}\sum_{Y=S, \Omega}W_1[S^2 Y\Phi]^\f12(t)\\
&+\l t\r^{-\frac{3}{2}+2\delta}(\La_0+(\log \l t\r)^\f12\Delta_0^\frac{5}{4}).
\end{align*}
Applying Proposition \ref{MA2} to $\psi=S^2 X_1\Phi,$ and substituting the above estimates into the bound, we can obtain 
\begin{equation*}
W_1[S^2X_1\Phi](t)+F_1[S^2 X_1\Phi](\H_u^t)\les \log \l t\r^\M(\La_0^2+\Delta_0^\frac{5}{2}).
\end{equation*}
Hence the first two sets of estimates in Proposition \ref{10.30.4.21} (1) are proved. The lower order estimate follows as its consequence by using (\ref{8.25.2.21}). The remaining two sets of estimates follow by using the first three sets of estimates, Lemma \ref{6.30.4.23}, Lemma \ref{10.10.3.23}, (\ref{8.29.9.21}), (\ref{10.10.2.23}) and Proposition \ref{7.15.5.22}. Thus we completed the proof of Proposition \ref{10.30.4.21} (1).\\

{\bf $\bullet$ Proof of (2)}  The first estimate of (\ref{3.9.7.24}) can be obtained by using (\ref{2.26.2.24}), the first estimate in Proposition \ref{10.30.4.21} (1) and (\ref{3.10.1.24}). The second one is an interpolation between the first one and (\ref{2.20.2.24}), with the help of (\ref{1.27.5.24}). The first estimate in (\ref{3.9.9.24}) is obtained by using (\ref{2.22.1.24}) and Proposition \ref{10.30.4.21} (1). The second one is an interpolation between the first one and (\ref{2.14.1.24}).

Similar to (\ref{3.29.2.23}) and (\ref{5.02.4.21}), for $f=\bT\varrho$ in (\ref{10.30.2.21}) we write
\begin{align*}
\Er_2(\bT\varrho, \Omega^2)=\Omega(\sP[\Omega, \bT \varrho])+\sP[\Omega, \Omega \bT\varrho]+I[\Omega^2\bT\varrho]
\end{align*}
with
\begin{align}\label{3.9.11.24}
I[\Omega^2 \bT\varrho]&=\Omega \fm{\Omega}\Box_\bg \bT\varrho+\fm{\Omega}(\Omega\Box_\bg \bT\varrho+[\Box_\bg,\Omega] \bT\varrho).
\end{align}
Using (\ref{1.29.2.22}), (\ref{11.11.2.23}) and (\ref{3.9.7.24}), we have
\begin{align*}
\|\Omega \fm{\Omega}\Lb^2\varrho\c \fB\|_{L^2_\Sigma}\les\|\bb^{-\frac{5}{2}}\Omega \fm{\Omega}\c \fB\|_{L_u^2 L_\omega^4}\les (\log \l t\r)^{\f12\M^2-2}\l t\r^{-1}(\La_0+\Delta_0^\frac{5}{4}). 
\end{align*}
Then using (\ref{3.9.6.24}) and (\ref{2.20.2.24}), the first term on the right-hand side of (\ref{3.9.11.24}) can be bounded by
\begin{align*}
\|\Omega \fm{\Omega}\Box_\bg \bT\varrho\|_{L^2_\Sigma}&\les \l t\r^{-1}\log \l t\r^{\f12\M^2-2}(\La_0+\Delta_0^\frac{5}{4}).
\end{align*}
We bound the last term in (\ref{3.9.11.24}) by using (\ref{1.27.5.24}), (\ref{3.9.12.24}) and (\ref{3.6.2.24})
\begin{align*}
\|\fm{\Omega}[\Box_\bg,\Omega] \bT\varrho\|_{L^2_\Sigma}\les \Delta_0\l t\r^{-1}\log \l t\r^4(\La_0+\Delta_0^\frac{5}{4}). 
\end{align*}
Next we compute by using (\ref{1.27.5.24}) and Cauchy-Schwarz inequality that 
\begin{align*}
\|\fm{\Omega}\Omega\Box_\bg \bT\varrho\|_{L^2_\Sigma}&\les \log \l t\r\Delta_0\|\tir\Omega \Box_\bg \bT\varrho\|_{L_u^2 L_\omega^4}\\
&\les \log \l t\r\Delta_0\|\tir\Omega^2\Box_\bg \bT\varrho\|_{L_u^2 L_\omega^2}^\f12\|\tir\Omega \Box_\bg \bT\varrho\|_{L_u^2 L_\omega^2}^\f12\\
&\les\Delta_0(\log \l t\r^{-1}\|\tir\Omega^2\Box_\bg \bT\varrho\|_{L_u^2 L_\omega^2}+\log \l t\r^3 \|\tir\Omega \Box_\bg \bT\varrho\|_{L_u^2 L_\omega^2}).
\end{align*}
Note that $\bA_{g,1}$ are combinations of various components of $\sn\Phi^\dagger$. Using (\ref{3.10.1.24}), we have 
$$
\|\sn_\Omega^2\sn\bA_{g,1}\|_{L^2_\Sigma}\les \l t\r^{-\frac{3}{2}} (\log \l t\r)^{\M^2}(\La_0+\Delta_0^\frac{5}{4}).
$$
Substituting the above estimate, (\ref{3.9.7.24}) and (\ref{3.9.9.24}) to (\ref{3.7.5.24}), we obtain
\begin{align}\label{1.1.1.25}
\begin{split}
\|\Omega^2\Box_g \bT\varrho\|_{L^2_\Sigma}&\les\|\fB\|_{L^\infty_x}(\|\bN\Omega^2\bT\varrho\|_{L^2_\Sigma}+\l t\r^{-\frac{3}{2}+\delta}(\La_0+\Delta_0^\frac{5}{4}))\\
&+\l t\r^{-1} \log \l t\r^{\f12\M^2-2}(\La_0+\Delta_0^\frac{5}{4}).
\end{split}
\end{align}
Combining the above two estimates and (\ref{3.9.8.24}) implies
\begin{align*}
\|\fm{\Omega}\Omega\Box_\bg \bT\varrho\|_{L^2_\Sigma}\les \Delta_0\|\fB\|_{L^\infty_x}\log \l t\r^{-1}\|\bN\Omega^2\bT\varrho\|_{L^2_\Sigma}+\l t\r^{-1}\log \l t\r^{\f12\M^2-3}\Delta_0^2.
\end{align*}
Summarizing the above estimates in view of (\ref{3.9.11.24}), we conclude
\begin{equation*}
\|I[\Omega^2 \bT\varrho]\|_{L^2_\Sigma}\les  \Delta_0\|\fB\|_{L^\infty_x}\log \l t\r^{-1}\|\bN\Omega^2\bT\varrho\|_{L^2_\Sigma}+\l t\r^{-1}\log \l t\r^{\f12\M^2-2}(\La_0+\Delta_0^\frac{5}{4}).
\end{equation*}
Using (\ref{3.6.2.24}) and the above estimate leads to 
\begin{align*}
\|\Er_2(\bT\varrho, \Omega^2)\|_{L^2_\Sigma}\les  \Delta_0\|\fB\|_{L^\infty_x}\log \l t\r^{-1}\|\bN\Omega^2\bT\varrho\|_{L^2_\Sigma}+ \l t\r^{-1}\log \l t\r^{\f12\M^2-2}(\La_0+\Delta_0^\frac{5}{4}).
\end{align*}
Applying (\ref{6.21.2.21}) to $\Omega^2\bT\varrho$ with the help of the above estimate and (\ref{1.1.1.25})  gives (\ref{3.9.10.24}).

 \end{proof}

\subsubsection{Improvement of Assumption \ref{5.13.11.21+}} \label{Improve}
Using the result of Proposition \ref{10.30.4.21}, we can improve an important set of the second order estimates in Proposition \ref{8.29.8.21} and Proposition \ref{9.8.6.22}, which proves the bootstrap assumption for all second order quantities. 
\begin{proposition}\label{3.14.4.24}
Let $X\in \{\Omega, S\}$ and $u_0\le u\le u_*$ and $0<t<T_*$. Under the assumptions of (\ref{3.12.1.21})-(\ref{6.5.1.21}), we have
\begin{align*}
&WFIL_2[\Omega^2 \Phi](\D_u^t)\les (\log \l t\r)^\M(\La_0^2+(\log \l t\r)^{2\M^2-\M-1}\Delta_0^\frac{5}{2})\\
\displaybreak[0]
&\|\Omega^3\Phi, \Sc(\Omega^3\Phi), \Ac(\Omega^3\Phi), \tir\sn^2_\Omega \bA_{g,1}\|^2_{L^2_\Sigma}\les(\log \l t\r)^\M(\La_0^2+(\log \l t\r)^{2\M^2-\M-1}\Delta_0^\frac{5}{2})\\
&\|\Omega \wt{\tr\chi}, \sn^{\le 1}_\Omega \chih\|_{L_u^2 L_\omega^2}\les  \l t\r^{-2}\log \l t\r^{\M^2+\f12}(\La_0+\Delta_0^\frac{5}{4}) \\
&\|\tir(\sn_S \sn, \sn \sn_S)(\tr\chi, \bA_{g,2})\|_{L_u^2 L_\omega^2}\les \l t\r^{-2}\log \l t\r^{\M^2-\f12}(\La_0+\Delta_0^\frac{5}{4})\\
&\|\tir \Omega^2\Phi\|_{L_\omega^4}^4\les (\La_0^2+\Delta_0^5)(\log \l t\r)^{2\M^2+10}\\
&\|\tir \Omega^2\varrho\|_{L_\omega^4}^4\les (\Delta_0^5(\log \l t\r)^{2\M^2-\M}+\La_0^4)\log \l t\r^{\M}\\
&\|\sn_X\bA_{g,1}\|_{L_\omega^4}\les\l t\r^{-2} \log \l t\r^{\frac{\M^2+5}{2}}(\Delta_0^\frac{5}{4}+\La_0)
\end{align*}
\end{proposition}
The first estimate was proved in Proposition \ref{10.30.4.21} (1). The rest of the above improved estimates can be proved in the same way as before, by also using  the results in (1) in Proposition \ref{10.30.4.21}. 

Next we establish the top order decay estimates. 
\begin{proposition}\label{imp_decay}
Let $X\in \{\Omega, S\}$ and $u_0\le u\le u_*$ and $0<t<T_*$. Under the assumptions of (\ref{3.12.1.21})-(\ref{6.5.1.21}), with  $Y=\sn, \sn_L$, we have
\begin{align}
&\|\tir SX^2\Phi\|_{L_\omega^4}^4+\sum_{\vs^+(X^2)=1}\|\tir \Omega X^2\Phi\|_{L_\omega^4}^4\les\log \l t\r^{2\M^2}\l t\r(\La_0^2+\Delta_0^5)\label{3.15.2.24}\\
&\|S^2 X\Phi, \Sc(S^2 X\Phi), \Ac(S^2X\Phi)\|_{L_\omega^4}^4\les\l t\r^{-3}\log \l t\r^{2\M}(\La_0^2+\Delta_0^5)\label{3.15.3.24}\\
&\|\tir \Omega^3\varrho\|_{L_\omega^4}\les \l t\r^{\frac{1}{4}}\log \l t\r^{(\M^2-\f12\M)}(\La_0+\Delta_0^\frac{5}{4}) \label{3.16.15.24}\\ 
&\|\tir \Omega^3\Phi\|_{L_\omega^4}\les\l t\r^{\frac{1}{4}}\log \l t\r^{\M^2}(\Delta_0^\frac{5}{4}+\La_0^\f12)\label{3.18.12.24}\\
\displaybreak[0]
&\|\bb^{-\frac{1}{2}}(\tir \sn)^{\ell+\le 1} (\bA_b, \bA_{g,2})\|_{L^2_\Sigma}\le \l t\r^{-1+\f12\ell}\log \l t\r^{\M^2+1-\ell}(\La_0+\Delta_0^\frac{5}{4}), \ell=0,1\label{L2conn'}\\
&\|\l t\r^{-\frac{l}{2}}X^{l+\le 1} Y (\bA^\natural), \l t\r^{-\frac{l}{2}-1} X^{l+\le 2}Y\Phi\|_{L^2_\Sigma}\les \l t\r^{-2} (\log \l t\r)^{\M^2}(\La_0+\Delta_0^\frac{5}{4}),\, l=0,1\label{3.16.7.24}\\
\displaybreak[0]
&\|X \Omega[L\Phi], \Omega X[L\Phi], \sn_X(\tir \sn)[L\Phi], \tir \sn \sn_X[L\Phi]\|_{L_\omega^4}\les  \l t\r^{-\frac{7}{4}}\log \l t\r^{\frac{\M+3}{2}} (\La_0+\Delta_0^\frac{5}{4})\label{3.16.3.24}\\
&\|\bb^{-\f12}\{\sn_X^2\sn[\Lb\Phi], \sn_X^2 \sn_\Lb\bA_{g,1},\sn_X \sn_\Lb \sn_X\bA_{g,1}, \sn_\Lb \sn_X^2 \bA_{g,1}\}\|_{L^2_\Sigma}\nn\\
&\qquad\qquad\qquad\les  \l t\r^{-1}\log \l t\r^{\M^2-\f12\M+\f12}(\La_0+\Delta_0^\frac{5}{4})\label{3.16.13.24}\\
&\| \sn_X^2\bA_{g,1}\|_{L_\omega^4}+|\sn_X \bA_{g,1}|\les \l t\r^{-\frac{7}{4}}\log \l t\r^{\M^2}(\Delta_0^\frac{5}{4}+\La_0)\label{3.18.2.24}\\
\displaybreak[0]
&\left\{\begin{array}{lll}
\|\tir^2 \sn(\bA_b, \bA_{g,2})\|_{L_\omega^4}\les\l t\r^{-\frac{3}{4}}\log \l t\r^{\M^2-\f12\M}(\La_0+\Delta_0^\frac{5}{4})\\
\|\bA_{g,2}\|_{L_\omega^4}\les \l t\r^{-2} \log\l t\r^{\f12(\M^2+6)}(\Delta_0^\frac{5}{4}+\La_0)\\
\|\bA_b\|_{L_\omega^4}\les  (\La_0^\f12+ \Delta_0^\frac{5}{4})\l t\r^{-2}(\log \l t\r)^{\f12\M^2}
\end{array}\right.
\label{9.11.3.22}\\
&\left\{
\begin{array}{lll}
|\bA_b|\les \l t\r^{-\frac{7}{4}}\log \l t\r^{\M^2-\f12\M}(\La_0^\f12+\Delta_0^\frac{5}{4})\\
|\chih|\les \l t\r^{-\frac{7}{4}}\log \l t\r^{\M^2-\f12\M}(\La_0+\Delta_0^\frac{5}{4}).
\end{array}\right.\label{3.17.1.24}\\
&|\Omega \fB|\les(\log \l t\r)^{\frac{3\M^2}{4}-1}\l t\r^{-1}(\La_0+\Delta_0^\frac{5}{4})\label{3.18.13.24}
\end{align}
\begin{align}
&\|\sn_X(\tir\sn)\fB, \log \l t\r^{-1}\tir\sn_X \sn_\Lb \bA_{g,1}\|_{L_\omega^4}\les (\log \l t\r)^{\frac{3\M^2}{4}-1}\l t\r^{-1}(\La_0+\Delta_0^\frac{5}{4})\label{3.18.3.24}\\
\displaybreak[0]
&\left\{
\begin{array}{lll}
\|\bb^{-1}\tir^2(\tir\sn)\sn\tr\chi\|_{L_\omega^4}&\les \log \l t\r^{\frac{3\M^2}{4}}(\La_0+\Delta_0^\frac{5}{4})\\
\|\tir \bb^{-1}\Omega^3\tr\chi\|_{L^2_u L_\omega^2}&\les (\log \l t\r)^{\M^2-\f12\M}(\La_0+\Delta_0^\frac{5}{4})
\end{array}\right.\label{3.19.2.24}\\
&\left\{\begin{array}{lll}
\|\tir(\tir \sn)\ze\|_{L_\omega^4}+|\tir \ze|\les\log \l t\r^{\frac{3\M^2}{4}}(\La_0+\Delta_0^\frac{5}{4})\\
\|\tir(\tir\sn)^2 \ze\|_{L_u^2 L_\omega^2}\les (\log \l t\r)^{\M^2-\f12\M+\f12}(\La_0+\Delta_0^\frac{5}{4})
\end{array}\right.\label{3.19.4.24}\\
& \begin{array}{lll}
X\Box_\bg \Phi=O(\l t\r^{-2}\log \l t\r^{\f12\M+1}(\Delta_0^\frac{5}{4}+\La_0))_{L^2_\Sigma}\\
\Box_\bg \Omega\Phi=O(\l t\r^{-2+\delta}\log \l t\r^{\f12\M+\frac{3}{2}}(\La_0+\Delta_0^\frac{5}{4}))_{L^2_\Sigma}
\end{array}
\label{wave_ass'}
\end{align}
\end{proposition}
\begin{proof}
We first prove  the first estimate in (\ref{3.15.2.24}) by considering the two cases (a) $SSX\Phi$, (b) $S\Omega^2\Phi$. For (a), we apply (\ref{9.19.3.23}) to $f=SX\Phi$, followed with applying (\ref{9.19.2.23}) and (\ref{10.22.2.23}) to $f=X\Phi$, which leads to
\begin{align*}
\|\tir SSX\Phi\|_{L_\omega^4}^4\les \La_0^2+\l t\r\log \l t\r^{2\M}(\La_0^2+\Delta_0^5), X=S, \Omega,
\end{align*}
 where we also employed (\ref{8.25.2.21}) and (1) in Proposition \ref{10.30.4.21}.
 
 Using the above estimate, (\ref{9.19.5.23}) and Lemma \ref{6.30.4.23}, we obtain (\ref{3.15.3.24}).

 For (b), we first apply (\ref{9.19.3.23}), followed with using (\ref{10.22.2.23}) with $X^l\Phi=\Omega^2\Phi$, (1) in Proposition \ref{10.30.4.21}, which yields
 \begin{align*}
 \|\tir S\Omega^2\Phi\|_{L_\omega^4}^4\les\La_0^2+\log \l t\r^{2\M^2}\l t\r(\La_0^2+\Delta_0^5).
 \end{align*}
 
Note due to (\ref{5.21.1.21}), (\ref{7.25.2.22})  and (\ref{9.19.5.23})
 \begin{align}
 \|X[\Omega, S]\Phi\|_{L_\omega^4}&\les \|\sn_X\pioh_{A L}\|_{L_\omega^4}\|\Omega\Phi\|_{L^\infty_x}+\l t\r^{-\frac{3}{4}+\delta}\Delta_0^\f12 \|X \Omega\Phi\|_{L_\omega^4}\nn\\ 
 &\les \l t\r^{-\frac{3}{4}+\delta}\Delta_0^\f12\l t\r^{-1+\f12\delta+}(\Delta_0^\frac{5}{4}+\La_0^\f12)\nn\\
 &\les \l t\r^{-\frac{7}{4}+\frac{3}{2}\delta+}(\La_0^\f12+\Delta_0^\frac{3}{2}).\label{3.15.4.24}
 \end{align}
Combining the above three estimates, we conclude the first estimate of (\ref{3.15.2.24}).

Next, we obtain by using (\ref{5.21.1.21}) and (\ref{9.19.5.23})
\begin{align*}
\|[\Omega, S]X\Phi\|_{L_\omega^4}&\les \l t\r^{-\frac{3}{4}+\delta}\Delta_0^\f12\|\Omega X\Phi\|_{L_\omega^4}\les \l t\r^{-\frac{7}{4}+\frac{3\delta}{2}+}(\La_0^\f12+\Delta_0^\frac{3}{2}).
\end{align*} 
The second estimate in (\ref{3.15.2.24}) follows by using the proved part of (\ref{3.15.2.24}), the above estimate and (\ref{3.15.4.24}).

Applying (\ref{9.19.4.23}) to $f=\Omega^2\varrho$, using (1) in Proposition \ref{10.30.4.21}, we can obtain (\ref{3.16.15.24}).

Using Proposition \ref{11.4.1.22}, Proposition \ref{9.8.6.22}, (\ref{1.21.2.22}) and (\ref{3.10.1.24}), we can obtain (\ref{L2conn'}).

Next we prove (\ref{3.16.7.24}) which improves the higher order estimates in (\ref{L2BA2}). The the second set, derivative estimates of $\Phi$, can be obtained by using Lemma \ref{3.17.2.22}, (\ref{10.10.2.23}), (\ref{8.29.9.21}) together with the second estimate in Proposition \ref{3.14.4.24} and (\ref{3.10.1.24}).  For $\bAn=[L\Phi], [\sn\Phi]$ the result follows from (\ref{3.10.1.24}). It remains to consider the case for $\bAn=\eh$, which can be done by using (\ref{8.9.5.23})  and the proved cases for $\bAn=[L\Phi], [\sn\Phi]$. Due to $\tr\eta=[L\Phi]$, using Lemma \ref{3.17.2.22}, (\ref{3.10.1.24}), (\ref{8.29.9.21}) together with the second estimate in Proposition \ref{3.14.4.24}, we derive
\begin{align*}
\|\sn_X^2(\tir Y)^l \eh\|_{L^2_\Sigma}\les  \l t\r^{\frac{l}{2}-1}\log \l t\r^{\M^2}(\La_0+\Delta_0^\frac{5}{4}), \, l=0,1.
\end{align*}
 We thus conclude (\ref{3.16.7.24}).

Next we prove (\ref{3.16.3.24}) by applying (\ref{9.20.1.23}) to $F=X\Omega[L\Phi]$. For this purpose, we need to bound
 $\|\Lb X\Omega[L\Phi]\|_{L^2_\Sigma}$. We first write
\begin{align*}
\Lb X\Omega[L\Phi]=[\Lb, X]\Omega[L\Phi]+X[\Lb, \Omega][L\Phi]+X\Omega\Lb[L\Phi].
\end{align*}
 In view of (\ref{3.19.2}) and (\ref{5.13.10.21}), the commutators can be written as
\begin{align*}
[\Lb, X]\Omega[L\Phi]&=\vs(X)(O(1) L\Omega[L\Phi]+\tir (\ud\bA\sn\Omega[L\Phi]+k_{\bN\bN}\bN \Omega[L\Phi]))\\
&+(1-\vs(X))(\pioh_{A\Lb}\sn\Omega [L\Phi]+\Omega\log \bb \bN \Omega[L\Phi]),\\
X[\Lb, \Omega][L\Phi]&=X(\pioh_{\Lb A}\sn[L\Phi]+\Omega\log \bb \bN [L\Phi]).
\end{align*}
Note the following improved estimates over (\ref{6.7.4.23})
\begin{align*}
X\Lb[L\Phi]&=\vs(X)(\l t\r^{-1}[\Lb\Phi]+O(\l t\r^{-\frac{11}{4}+\delta}\Delta_0^\f12)_{L_\omega^4})+(1-\vs(X))O(\l t\r^{-2}\log \l t\r^{\f12\M+1}\Delta_0)_{L_\omega^4}\\
\Lb\Omega[L\Phi]&=O(\l t\r^{-2}\log \l t\r^{\f12\M+1}\Delta_0)_{L_\omega^4}
\end{align*}
which are obtained by using Proposition \ref{7.15.5.22}, (\ref{5.13.10.21}) and (\ref{9.14.3.22}) in view of (\ref{8.26.1.23}).
Using the above estimate, (\ref{8.29.9.21}), (\ref{2.20.2.24}), (\ref{2.19.1.24}), (\ref{5.21.1.21}), (\ref{3.6.2.21}), (\ref{6.22.1.21}) and (\ref{1.27.5.24}), we bound
\begin{align}\label{3.16.5.24}
\begin{split}
\|[\Lb, X]\Omega[L\Phi]\|_{L^2_\Sigma}&\les \l t\r^{-2}\log \l t\r^{\f12\M}(\La_0+\Delta_0^\frac{5}{4})(\l t\r^{\delta}\Delta_0^\f12+1)+\l t\r^{-1}\log \l t\r^{\f12\M+2}(\Delta_0^\frac{5}{4}+\La_0)\\
\|X[\Lb, \Omega][L\Phi]\|_{L^2_\Sigma}&\les \l t\r^{-1}\log \l t\r^4(\La_0+\Delta_0^\frac{5}{4})+\|\Omega\log\bb X\Lb[L\Phi]\|_{L^2_\Sigma}\\
&\les \l t\r^{-1}\log \l t\r^{\f12\M+2}(\Delta_0^\frac{5}{4}+\La_0).
\end{split}
\end{align}
Next we will use (\ref{3.16.4.24}) to bound $\|X\Omega\Lb[L\Phi]\|_{L^2_\Sigma}$, which requires us to bound $\|X^2(\bA_b\c \fB)\|_{L^2_\Sigma}$. This can be done by applying (\ref{8.9.4.22}) to $F=\bA$, 
\begin{align*}
X^2(\fB \bA)&=O(\bb^{-1} \l t\r^{-1})\sn_X^{\le 2}\bA+O(\Delta_0\l t\r^{-1}\log \l t\r)_{L_\omega^4}X \bA+O(\Delta_0\log \l t\r^2)_{L^2_\Sigma}\bA.
\end{align*}
Substituting (\ref{3.16.7.24}), (\ref{L2conn'}), (\ref{L2conndrv'}), (\ref{8.8.6.22'}) and Proposition \ref{7.15.5.22} to the above leads to 
\begin{equation}\label{6.15.3.24}
\|X^2(\fB\bA)\|_{L^2_\Sigma}\les \l t\r^{-\frac{3}{2}}\log \l t\r^{\M^2}(\La_0+\Delta_0^\frac{5}{4}).
\end{equation}
Substituting the above estimate to (\ref{3.16.4.24}) yields
\begin{align}\label{3.16.8.24}
\begin{split}
&X^2(\Lb [L\Phi]-\sD\varrho-L[L\Phi])\\
&=\vs^-(X^2)O(\l t\r^{-2}\bb^{-1})+O(\l t\r^{-2}\log \l t\r(\La_0+\log \l t\r^3\Delta_0^\frac{5}{4}))_{L_u^2 L_\omega^2}
\end{split}
\end{align}
It follows by further applying (\ref{3.10.1.24}) to the term $X^2(\sD\varrho, L[L\Phi])$ that
\begin{align*}
X\Omega\Lb[L\Phi]=O(\l t\r^{-2}\log \l t\r(\La_0+\log \l t\r^3\Delta_0^\frac{5}{4}))_{L_u^2 L_\omega^2}.
\end{align*}
Hence, we conclude by combining the above estimate with the commutator estimates
\begin{align}\label{6.23.10.24}
\|\Lb X\Omega[L\Phi]\|_{L^2_u L_\omega^2}\les \l t\r^{-2}\log \l t\r^{\f12\M+2}(\La_0+\Delta_0^\frac{5}{4}).
\end{align}
Substituting the above estimate and (\ref{3.10.1.24}) to (\ref{9.20.1.23}) with $F=X\Omega[L\Phi]$, 
\begin{align*}
\|\tir F\|_{L_\omega^4}^2&\les \l t\r^{-\f12-1}\log \l t\r^{\M+3}(\La_0+\Delta_0^\frac{5}{4})^2. 
\end{align*}
The other estimate in (\ref{3.16.3.24}) follows by commutation (\ref{7.17.6.21}) and using the above estimate. Thus we conclude (\ref{3.16.3.24}).    

Next we prove
\begin{equation}\label{3.16.6.24}
\|\bb^{-\f12}\{\sn_X^2\sn[\Lb\Phi], \sn_X^2 \big(\sn_\Lb[\sn\Phi],\sn_\Lb\eh\big)\}\|_{L^2_\Sigma}\les \l t\r^{-1} (\log \l t\r)^{\M^2-\f12\M+\f12}(\La_0+\Delta_0^\frac{5}{4})
\end{equation}
which improves (\ref{LbBA2}), as stated in (\ref{3.16.13.24}).

Consider the first estimate in the case $X^2=\Omega^2$.  Using the first estimate of (\ref{3.9.9.24}), Lemma \ref{3.17.2.22}, Proposition \ref{1steng}, Proposition \ref{8.29.8.21} and Proposition \ref{7.15.5.22}, we obtain the first estimate in (\ref{3.16.6.24}) when $X^2=\Omega^2$. Next we consider $X_2 X_1=X S$. 
Applying (\ref{8.5.1.22+}), (\ref{LbBA2'}) and (\ref{3.16.7.24}) gives 
\begin{align}\label{9.12.1.24}
\|\sn_X(\tir\sn) L\fB, \sn_X \sn_L(\tir \sn) \fB\|_{L^2_u L_\omega^2}\les \l t\r^{-2}\log \l t\r(\La_0+\log \l t\r^3\Delta_0^\frac{5}{4}).
\end{align}
If $X=\Omega$, applying (\ref{3.21.1.23}) to $\sn \fB$, and using (\ref{LbBA2'}) and the second estimate in the above, we have 
\begin{equation*}
\|\sn_S \sn_\Omega\sn \fB\|_{L_u^2 L_\omega^2}\les\l t\r^{-2}\log \l t\r(\La_0+\log \l t\r^3\Delta_0^\frac{5}{4}).
\end{equation*} 
Thus the first estimate in (\ref{3.16.6.24}) is proved. In view of (\ref{4.22.4.22}) and  (\ref{LbBA2'}), this estimate implies
\begin{equation}\label{6.15.1.24}
\vs^+(X^2)\|\sn_X^2 \Omega \fB\|_{L_u^2 L_\omega^2}\les \l t\r^{-1}\log \l t\r(\La_0+\log \l t\r^3\Delta_0^\frac{5}{4}).
\end{equation}
To see the second estimate, we consider the cases that $[\sn\Phi]=\ep$ and $\sn \varrho$ separately. For this purpose, we first give a preliminary estimate. 

Due to (\ref{9.12.1.24}) and (\ref{2.19.1.24}) 
\begin{equation*}
\|\Omega S\fB, S\Omega \fB\|_{L^2_u L_\omega^4}\les  \l t\r^{-1}\log \l t\r(\La_0+\log \l t\r^3\Delta_0^\frac{5}{4}).
\end{equation*} 
Note by using the above estimate, (\ref{8.23.1.23}) and (\ref{3.9.9.24}) for an $S$-tangent tensor or scalar, (\ref{8.9.4.22}) can be improved to 
\begin{align*} 
X^2(\fB F)&=O(\bb^{-1} \l t\r^{-1})\sn_X^{\le 2}F+O\Big(\Delta_0(\l t\r^{-1}\log \l t\r)\Big)_{L_\omega^4}X F\nn\\
&+O(\l t\r^{-1}\log \l t\r^{\f12\M^2-2}\Delta_0)_{L^2_u L_\omega^4}\bb F
\end{align*}
Applying the above estimate to $F=\ud \bA$, also using (\ref{3.9.7.24}), (\ref{2.20.4.24}), (\ref{1.27.5.24}) and Proposition \ref{8.29.8.21}, leads to 
\begin{align*}
X^2(\fB \ud\bA)=O\Big(\l t\r^{-2}\log \l t\r^{\M^2-\f12\M+\f12}(\La_0+\Delta_0^\frac{5}{4})\Big)_{L_u^2 L_\omega^2}.
\end{align*}

In view of (\ref{10.1.8.23}), similar to the proof of (\ref{LbBA2'}), we use (\ref{8.28.2.23}), (\ref{1.30.1.24}), the second line of Proposition \ref{3.14.4.24}, the first estimate in (\ref{3.16.13.24}) and (\ref{3.16.7.24}) to derive that
\begin{align*}
\|\sn_X^2\sn_\bN \ep\|_{L^2_u L_\omega^2}&\les\|\sn_X^{\le 2}\Big((\chi+\fB)\c \ep+\sn \fB+\ud \bA (\fB+\bA_{g,1})\Big)\|_{L^2_u L_\omega^2}\\
&\les\l t\r^{-2}\log \l t\r^{\M^2-\f12\M+\f12}(\La_0+\Delta_0^\frac{5}{4}).
\end{align*}
Thus we completed the proof for the case that $[\sn\Phi]=\ep$.

Next, due to the commutator estimates given in (\ref{3.12.10.24}) and (\ref{3.13.10.24}), we have, for $\vs^+(X^2)=1$, that 
\begin{align*}
\|\bb^\f12 X^2[\Omega, \Lb]\varrho\|_{L_u^2L_\omega^2}\les\l t\r^{-1}\log \l t\r(\La_0+\log \l t\r^3\Delta_0^\frac{5}{4}).
\end{align*}
Moreover, using (\ref{3.9.7.24}), the first estimate in Proposition \ref{10.30.4.21} (1) and (\ref{6.15.2.24}), we deduce
\begin{equation*}
\|\bb^\f12\Omega^2[\Omega, \Lb]\varrho\|_{L_u^2 L_\omega^2}\les \l t\r^{-1}\log \l t\r^{\M^2-\f12\M}(\La_0+\Delta_0^\frac{5}{4}).
\end{equation*}
Combining the above estimates with (\ref{6.15.1.24}) and the first estimate in (\ref{3.9.9.24}), we obtain
\begin{align*}
\|\bb^\f12 X^2\Lb \Omega\varrho\|_{L^2_u L_\omega^2}\les \l t\r^{-1}\log \l t\r^{\M^2-\f12\M}(\La_0+\Delta_0^\frac{5}{4})
\end{align*}
from which we can obtain the second estimate for the case $[\sn\Phi]=\sn\varrho$ in (\ref{3.16.6.24}) with the help of Lemma \ref{3.17.2.22}. 
It only remains to consider the last estimate in (\ref{3.16.6.24}). By using (\ref{1.30.2.24}), (\ref{1.30.1.24}), (\ref{6.15.3.24}) and (\ref{3.16.7.24}), we derive the stronger estimate
\begin{align*}
\|\sn_X^2 \sn_\Lb\eh\|_{L^2_\Sigma}\les \l t\r^{-\frac{3}{2}} \log \l t\r^{\M^2}(\Delta_0^\frac{5}{4}+\La_0).
\end{align*}

Next we prove the last two estimates in (\ref{3.16.13.24}) with the help of the second one. 
Using (\ref{7.03.1.19}), (\ref{7.04.5.21}) and (\ref{5.13.10.21}), we write
\begin{align}
&[\sn_\Lb, \sn_\Omega]F=\pioh_{A\Lb}\sn_A F+\ud\bA(\Omega) \sn_\bN F+((\chi, \chib)\ud\bA+\bR_{AC\Lb B})F\c \Omega\label{3.18.4.24}\\
&[\sn_\Lb, \sn_S]F=\Lb\tir L F+\tir(\ud \bA\sn F+k_{\bN\bN}\sn_\bN F+(\bR_{ABL\Lb}+\bA_{g,1}\ud \bA)F) \label{3.18.5.24}\\
&\bR_{AB43}=\ud \bA[\sn\Phi], \bR_{AB C\Lb}=\sn\Lb \varrho+\Lb \sn\varrho+[\sn\Phi](\bA+\tir^{-1})+\ud \bA\fB\label{3.18.6.24}
\end{align}
where the last line are recalled from (\ref{1.27.2.22}) and (\ref{9.2.1.22}). 
We will prove with $a=0,1$ that
\begin{align}
 I_a:&=\|\sn_X^{1-a}[\sn_\Lb, \sn_\Omega]\sn_X^a \bA_{g,1}\|_{L^2_\Sigma}\les \log \l t\r^{\f12(\M^2+\M)}\l t\r^{-1} (\La_0+\Delta_0^\frac{5}{4}),\label{3.18.8.24}\\
 II_a:&=\|\sn_X^{1-a}[\sn_\Lb, \sn_S]\sn_X^a \bA_{g,1}\|_{L^2_\Sigma}\les \l t\r^{-1}(\log \l t\r)^2(\log \l t\r^3\Delta_0^\frac{5}{4}+\La_0).\label{3.18.9.24}
 \end{align}
Note that using the estimates of $I_0, II_0$ and the second estimate in (\ref{3.16.13.24}) gives the third estimate in (\ref{3.16.13.24}); then combining the estimates of $I_1$, $II_1$ with the third estimate of (\ref{3.16.13.24}) yields the last estimate of (\ref{3.16.13.24}).

We first prove (\ref{3.18.9.24}). Setting  $II=II_0+II_1$, it follows by using (\ref{3.18.5.24}) that
\begin{align*}
II&\les\sum_{a=0}^1\|\sn_X^{1-a}\Big(\tir \ud\bA\sn \sn_X^a \bA_{g,1}+\tir k_{\bN\bN}\sn_\bN\sn_X^a \bA_{g,1}+\tir\ud \bA\bA_{g,1}\sn_X^a \bA_{g,1}\Big), \\
&\sn_X^{1-a}(\Lb \tir\sn_L \sn_X^a \bA_{g,1})\|_{L^2_\Sigma}.
\end{align*}
Using (\ref{LbBA2'}), (\ref{7.13.5.22}) and (\ref{8.23.2.23}), we obtain
\begin{align*}
\sum_{a=0}^1\|\sn_X^{1-a}(\tir k_{\bN\bN}\sn_\bN \sn_X^a \bA_{g,1})\|_{L^2_\Sigma}\les\l t\r^{-1}(\log \l t\r)^2(\log \l t\r^3\Delta_0^\frac{5}{4}+\La_0).
\end{align*}
Using Proposition \ref{7.22.2.22}, (\ref{3.16.1.22}), (\ref{3.11.3.21}), (\ref{8.23.2.23}) and (\ref{L2BA2'}), we derive
\begin{align*}
\sum_{a=0}^1&\|\sn_X^{1-a}(\Lb \tir \sn_L \sn_X^a \bA_{g,1})\|_{L^2_\Sigma}\\
&\les \|\sn_X(\tir \mho+\f12\tir \tr\chib)\sn_L \bA_{g,1}\|_{L^2_\Sigma}+\|\sn_X \sn_L \bA_{g,1}, \sn_L \sn_X \bA_{g,1}\|_{L^2_\Sigma}\\
&\les \log \l t\r\Delta_0^\f12\|\bb^\f12\tir\sn_L \bA_{g,1}\|_{L_u^2 L_\omega^4}+\|(\sn_X^{\le 1}\sn_L, \sn_L \sn_X) \bA_{g,1}\|_{L^2_\Sigma}\\
&\les \l t\r^{-2}  (\log \l t\r)^{\f12\M+\frac{3}{2}}\big(\La_0+\Delta_0^\frac{5}{4}\big).
\end{align*}
It follows by using Proposition \ref{7.15.5.22}  and (\ref{9.12.2.22}) that
\begin{align*}
&\sum_{a=0}^1\|\sn_X^{1-a}(\tir \ud\bA\sn \sn_X^a \bA_{g,1})\|_{L^2_\Sigma}\\
&\les \|\sn_X\ud \bA\|_{L_u^2 L_\omega^4}\|\tir^2\bb^\f12 \sn\bA_{g,1}\|_{L_u^\infty L_\omega^4}+\|\bb\ud \bA\|_{L_u^\infty L_\omega^4}\|\bb^{-1}\tir (\sn \sn_X \bA_{g,1}, \sn_X \sn\bA_{g,1})\|_{L^2_\Sigma}\\
&\les \l t\r^{-2+2\delta}\Delta_0^2.
\end{align*}
Similarly,
\begin{align*}
\|\sn_X^{1-a}(\tir\ud \bA\bA_{g,1}\sn_X^a \bA_{g,1})\|_{L^2_\Sigma}\les \l t\r^{-3+3\delta}\Delta_0^2. 
\end{align*}
(\ref{3.18.9.24}) follows by summarizing the above estimates. 
  
 We next consider (\ref{3.18.8.24}). Due to (\ref{3.18.4.24}),
\begin{align}
I_a&\les\|\sn_X^{1-a}\Big(\pioh_{\Lb A}\sn_A \sn_X^a\bA_{g,1}, \ud\bA(\Omega) \sn_\bN \sn_X^a \bA_{g,1}\Big), \nn\\
&O(\tir)\sn_X^{\le 1-a}\Big(((\bA+\fB+\tir^{-1})\ud\bA+\sn\Lb \varrho+\Lb \sn\varrho)\sn_X^a \bA_{g,1}\Big) \|_{L^2_\Sigma}.\label{6.15.4.24}
\end{align}
Note that, due to (\ref{LbBA2'}) and the proved part in (\ref{3.16.13.24}), 
\begin{align}\label{3.18.7.24}
\|\sn_X \sn[\Lb\Phi], \sn_X \sn_\Lb \bA_{g,1}\|_{L_u^2 L_\omega^4}\les\log \l t\r^{\f12\M^2+2}\l t\r^{-2}(\Delta_0^\frac{5}{4}+\La_0). 
\end{align}
We first consider the most crucial term in terms of decay behavior in the commutators. 
Using (\ref{1.27.5.24}), (\ref{3.9.7.24}), (\ref{2.20.4.24}), (\ref{9.14.3.22}) and the above estimate, 
\begin{align*}
&\|\sn_X(\ud\bA(\Omega) \sn_\bN\bA_{g,1})\|_{L^2_\Sigma}\\
 &\les \|\bb^\f12\sn_X \ud\bA(\Omega)\|_{L_u^2 L_\omega^4}\|\tir\sn_\bN\bA_{g,1}\|_{L_\omega^4}+\|\bb^\f12 \ud\bA(\Omega)\|_{L_\omega^4}\|\tir\sn_X \sn_\bN \bA_{g,1}\|_{L_u^2 L_\omega^4}\\
 &\les \log \l t\r^{\f12(\M^2+\M)}\l t\r^{-1}(\La_0+\Delta_0^\frac{5}{4})^2. 
\end{align*}
Using (\ref{3.18.7.24}), Lemma \ref{5.13.11.21} (5), the last estimate in Proposition \ref{3.14.4.24} and (\ref{10.11.2.23}), we obtain
\begin{align*}
\l t\r\sum_{a=0}^1\|\sn_X^{1-a}\big((\sn \Lb \varrho+\sn_\Lb \sn \varrho)\sn_X^a \bA_{g,1}\big)\|_{L^2_\Sigma}&\les \l t\r\{\|\sn\Lb\varrho+\sn_\Lb\sn \varrho\|_{L_\omega^4}\|\bb^\f12\tir\sn_X \bA_{g,1}\|_{L_u^2 L_\omega^4}\\
&+\|\sn_X(\sn\Lb \varrho+\sn_\Lb \sn \varrho)\|_{L_u^2 L_\omega^4}\|\bb^\f12\tir\bA_{g,1}\|_{L_\omega^4}\}\\
&\les\log \l t\r^{\f12\M^2+\frac{\M}{4}+3}\l t\r^{-2}(\Delta_0^\frac{5}{4}+\La_0)^2.
\end{align*} 
Moreover, due to (\ref{5.21.1.21}) and Proposition \ref{7.15.5.22}
\begin{align*}
\sum_{a=0}^1\|\sn_X^{1-a}(\pioh_{\Lb A}\sn \sn_X^a \bA_{g,1})\|_{L^2_\Sigma}&\les \l t\r^{-\frac{3}{4}+\delta}\Delta_0^\f12(\|\sn_X^{1-a}\sn \sn_X^a \bA_{g,1}\|_{L^2_\Sigma}+\|\bb^\f12 \tir \sn \bA_{g,1}\|_{L_u^2 L_\omega^4})\\
&\les  \l t\r^{-2-\frac{3}{4}+2\delta} \Delta_0^\frac{3}{2}.
\end{align*}
For simplicity, we neglect the remaining lower order terms and conclude (\ref{3.18.8.24}) for $a=0$.
Therefore, we proved the third estimate in (\ref{3.16.13.24}) by combing the $a=0$ case in (\ref{3.18.8.24}) with (\ref{3.18.9.24}). 
To  bound $I_1$, it suffices to consider the second term on the right-hand side of the first line of (\ref{6.15.4.24}), since the remaining terms have been treated. Recasting the term below, also using (\ref{9.12.2.22}), (\ref{LbBA2'}) and the third estimate in (\ref{3.16.13.24}), we deduce
\begin{align*}
 \|\ud\bA(\Omega) \sn_\bN \sn_X\bA_{g,1}\|_{L^2_\Sigma}&\les\|\bb^\f12\ud \bA(\Omega)\|_{L_u^\infty L_\omega^4}\|\tir\sn_\bN \sn_X \bA_{g,1}\|_{L^2_u L_\omega^4}\\
 &\les \Delta_0\l t\r \log \l t\r \|\sn_\Omega \sn_\bN \sn_X \bA_{g,1}\|^\f12_{L_u^2 L_\omega^2}\|\sn_\bN\sn_X\bA_{g,1}\|^\f12_{L_u^2 L_\omega^2}\\
 &\les \l t\r^{-1}\log \l t\r^{\f12\M^2+1}(\La_0+\Delta_0^\frac{5}{4})\Delta_0.
\end{align*}
Therefore, we can conclude the estimate of $I_1$ in (\ref{3.18.8.24}). Combining the estimates of $I_1$, $II_1$, and the third estimate in (\ref{3.16.13.24}) gives the last estimate in (\ref{3.16.13.24}).

(\ref{3.18.2.24}) is a consequence of (\ref{3.16.7.24}) and the last estimate of (\ref{3.16.13.24}) by using (\ref{9.20.1.23}). Similar to the proof of the last estimate in (\ref{9.19.5.23}), (\ref{3.18.12.24}) can be obtained by using (\ref{3.16.3.24}) and (\ref{3.18.2.24}).

Next we consider (\ref{9.11.3.22}). Using (\ref{3.16.14.24}), (\ref{3.16.15.24}), Proposition \ref{9.8.6.22} and Proposition \ref{12.21.1.21}, we derive
\begin{align*}
\|\tir\Omega \tr\chi\|_{L_\omega^4}\les \l t\r^{-\frac{3}{4}}\log \l t\r^{\M^2-\f12\M}(\La_0+\Delta_0^\frac{5}{4}).
\end{align*}
This gives the first estimate in (\ref{9.11.3.22}). 

In view of (\ref{1.21.2.22}), we derive by using (\ref{9.19.5.23}), $\chi,\theta, \bp\Phi^\dagger=O(\l t\r^{-1})$ and Proposition \ref{1steng} that 
\begin{equation*}
\|\bR_{ABCL}\|_{L_\omega^4}\les \l t\r^{-3}\log \l t\r^{\f12\M}(\Delta_0^\frac{3}{2}+\La_0).
\end{equation*}
Combining the above two estimates, we infer by using (\ref{8.25.2.23}) and (\ref{10.11.2.23}),
\begin{equation*}
\|\tir\sn\chih, \chih\|_{L_\omega^4}\les \l t\r^{-\frac{7}{4}}\log \l t\r^{\M^2-\f12\M}(\La_0+\Delta_0^\frac{5}{4}).
\end{equation*}
To improve the $L^4_\omega$ estimate of $\chih$, we apply (\ref{9.20.1.23}) to $F=\chih$ to obtain 
\begin{align*}
\|\tir \chih\|^2_{L_\omega^4}\les\|\tir \chih\|^2_{L_\omega^4(S_{t, u_*})}&+\left(\|\tir^{-1} \chih\|_{L^2_{\Sigma_t}\cap[u, u_*]}+\|\sn_\bN \chih\|_{L^2_{\Sigma_t}\cap[u, u_*]}\right)\nn\\
&\times \|\sn^{\le 1}_\Omega \chih, \Delta_0 \log \l t\r\bb^{-1}\chih\|_{L^2_{\Sigma_t}\cap[u, u_*]}.
\end{align*} 
Substituting (\ref{3.30.2.24}), the lower order estimate in (\ref{8.8.6.22'}), the estimate of $\sn_\Omega^{\le 1}\chih$ in Proposition \ref{3.14.4.24} to the above Sobolev inequality, we conclude 
\begin{align*}
\|\tir \chih\|^2_{L_\omega^4}\les \l t\r^{-2} \log\l t\r^{\M^2+\frac{11}{2}}(\Delta_0^\frac{5}{4}+\La_0)^2.
\end{align*}
The $\|\bA_b\|_{L^4_\omega}$ estimate can be derived by refining (\ref{2.1.1.22}) with the help of (\ref{9.19.5.23}) and the estimate of $\Omega^2 \varrho$ in Proposition \ref{3.14.4.24}.
Hence we have completed the proof of (\ref{9.11.3.22}). (\ref{3.17.1.24}) follows as a consequence of (\ref{9.11.3.22}) due to Sobolev embedding. 

Next we prove (\ref{3.18.3.24}). For the first estimate, due to (\ref{3.16.3.24})  and $\fB=\Lb \varrho+[L\Phi]$, it suffices to prove
\begin{equation}\label{3.18.10.24}
\|X\Omega \bT\varrho\|_{L^4_\omega}\les (\log \l t\r)^{\frac{3\M^2}{4}-1}\l t\r^{-1}(\La_0+\Delta_0^\frac{5}{4}).
\end{equation}
If $X=\Omega$, this follows by applying Proposition \ref{10.30.4.21} (2), (\ref{2.14.1.24}) and (\ref{9.20.1.23}).
If $X=S$, we first show
\begin{align}\label{3.18.11.24}
\|\Omega S\Lb \varrho, S\Omega\Lb \varrho\|_{L^4_\omega}\les  \l t\r^{-1}\log \l t\r^{\f12\M+1}(\Delta_0^\frac{5}{4}+\La_0).
\end{align}
We use (\ref{9.11.3.22}), (\ref{9.14.3.22}) and the last estimate in (\ref{3.14.4.24}) to refine (\ref{3.10.6.24}) for any $f$ to be
\begin{align*}
X_2 S\Lb f&=X_2\big(\tir (\sD f-\Box_\bg f)\big)+O\Big(\log \l t\r^{\f12\M+1}(\Delta_0^\frac{5}{4}
+\La_0)\Big)_{L^4_\omega}(\Lb f, Lf)\\
&+O(1)X_2(\Lb f, Lf)+O(\l t\r^{-1+\delta}\Delta_0)_{L_\omega^4}\sn f.
\end{align*}
With $f=\varrho, X_2=\Omega$ in the above, using (\ref{9.14.3.22}) and (\ref{3.6.2.21}), we obtain
\begin{align*}
\|\Omega S\Lb\varrho\|_{L_\omega^4}&\les \l t\r^{-1}\log \l t\r^{\f12\M+1}(\Delta_0^\frac{5}{4}+\La_0)+\|\tir\Omega(\sD\varrho, \Box_\bg \varrho)\|_{L_\omega^4}. 
\end{align*}
In view of (\ref{10.22.5.22}) and (\ref{L4conn'}), we have
\begin{equation*}
\|\tir \Omega \Box_\bg \varrho\|_{L_\omega^4}\les \l t\r^{-2+2\delta}(\Delta_0^\frac{5}{4}+\La_0).
\end{equation*}
 Finally, combining the above estimates, also applying (\ref{3.16.15.24}) to control $\Omega \sD\varrho$, we conclude the first estimate in (\ref{3.18.11.24}). 
 
 Note, due to (\ref{5.13.10.21}) and (\ref{5.21.1.21}) 
\begin{equation*}
[\Omega, S]\Lb \varrho=O(\l t\r^{-\frac{3}{4}+\delta}\Delta_0^\f12)\c \Omega \Lb \varrho.
\end{equation*}
 By substituting (\ref{9.14.3.22}) to the last term in the above, it is direct to see that the commutator is a negligible error. Thus we conclude the second estimate in (\ref{3.18.11.24}) and (\ref{3.18.10.24}). The first estimate in (\ref{3.18.3.24}) is proved. 
 
  From (\ref{3.18.10.24}), (\ref{3.16.3.24}), (\ref{9.14.3.22}) and Sobolev embedding on spheres, we can obtain (\ref{3.18.13.24}).
 
 Next substituting the first estimate of (\ref{3.18.3.24}) to (\ref{10.29.1.22}), and using (\ref{3.9.7.24}) and Lemma \ref{3.17.2.22}, we derive 
\begin{align*}
&\|\tir(\tir \sn)\sn\log \bb\|_{L_\omega^4}\les\log \l t\r^{\frac{3\M^2}{4}}(\La_0+\Delta_0^\frac{5}{4})\\
&\|\tir(\tir\sn)^2( \sn\log \bb)\|_{L_u^2 L_\omega^2}\les \log \l t\r^{\M^2-\f12\M+\f12}(\La_0+\Delta_0^\frac{5}{4}).
\end{align*}
Combining Proposition \ref{3.14.4.24}, (\ref{2.21.2.24}) with the above estimates gives the estimates of $\ze$ in (\ref{3.19.4.24}). Due to (\ref{1.22.4.22}), the pointwise estimates of $\bA_{g,1}$ in (\ref{7.25.2.22}) and of $\ze$ in (\ref{3.19.4.24}), integrating along null geodesics gives the pointwise bound of $\b$. Thus we completed the proof of (\ref{3.19.4.24}). 
 
   To see the second estimate of (\ref{3.18.3.24}), we will consider the cases that $\bA_{g,1}=\sn \varrho, \ep, \eh$ separately.  
We first derive by (\ref{5.13.10.21}), (\ref{5.21.1.21}) and Proposition \ref{7.15.5.22} that 
  \begin{align*}
X[\Lb, \Omega]\varrho&=X(\pioh_{A\Lb}\sn \varrho+\Omega\log \bb \bN\varrho)\\
&=O(\l t\r^{-\frac{11}{4}+2\delta}\Delta_0^2)_{L_\omega^4}+X\Omega\log \bb \c \fB+\Omega\log \bb X\bN \varrho.
  \end{align*}
  If $X=S$, from (\ref{2.21.2.24}), (\ref{6.17.1.24}), (\ref{9.19.5.23}) and (\ref{8.23.1.23}), we have 
  \begin{equation*}
  S[\Lb, \Omega]\varrho=O(\l t\r^{-1}\log\l t\r^{\f12\M+2}(\La_0+\Delta_0^\frac{5}{4}))_{L_\omega^4}.
  \end{equation*}
  If $X=\Omega$, from (\ref{2.20.2.24}), (\ref{3.19.4.24}) and (\ref{8.23.1.23}), we derive
  \begin{equation*}
  \Omega[\Lb, \Omega]\varrho=O(\log \l t\r^{\frac{3\M^2}{4}}\l t\r^{-1}(\La_0+\Delta_0^\frac{5}{4}))_{L_\omega^4}.
  \end{equation*}
  Thus we conclude with the help of Lemma \ref{3.17.2.22} that 
  \begin{equation*}
  \tir\sn_X \sn_\Lb \sn\varrho=O(\log \l t\r^{\frac{3\M^2}{4}}\l t\r^{-1}(\La_0+\Delta_0^\frac{5}{4}))_{L_\omega^4}.
  \end{equation*}
  For the case that $\bA_{g,1}=\ep$, using (\ref{10.1.8.23}) with the help of the first estimate in (\ref{3.18.3.24}), (\ref{2.21.2.24}), (\ref{6.17.1.24}), (\ref{8.23.1.23}), (\ref{1.27.5.24}) and (\ref{3.18.13.24})
  \begin{align*}
  \sn_X\sn_\bN \ep&=\sn_X\big(\sn \fB+\ud\bA(\fB+\bAn)+(\bA+\tir^{-1}) \ep\big)\\
  &= O(\log \l t\r^{\frac{3\M^2}{4}}\l t\r^{-2}(\La_0+\Delta_0^\frac{5}{4}))_{L_\omega^4}
  \end{align*}
  where we also employed Proposition \ref{7.15.5.22} to control quadratic errors, and the last estimate in Proposition \ref{3.14.4.24} to control lower order terms. 
  
  Using (\ref{2.1.2.24}), (\ref{8.23.1.23}), (\ref{3.18.2.24})-(\ref{3.18.13.24}), Proposition \ref{7.15.5.22} and the last estimate in Proposition \ref{3.14.4.24}, we derive
  \begin{align*}
  \|\sn_X \sn_\bN \eh\|_{L^4_\omega}&\les \|\sn_X(\snc \hot \ep, \ud\bA\c \ep, \fB \bA_g), \tir^{-1}\sn_X^{\le 1}\eh\|_{L_\omega^4}+\l t\r^{-\frac{15}{4}+2\delta}\Delta_0^\frac{3}{2}\\
  &\les \l t\r^{-\frac{11}{4}}\log \l t\r^{\M^2}(\La_0+\Delta_0^\frac{5}{4}).
\end{align*}
Summarizing the above three cases, the second estimate in (\ref{3.18.3.24}) is proved. 

Next we prove (\ref{3.19.2.24}).
Using  Proposition \ref{10.30.4.21}, (\ref{3.31.2.22}) and (\ref{6.2.3.24})) we have
\begin{equation*}
\|\tir^2 (\tir\sn)^2\sF\|_{L^2_u L_\omega^2}\les (\log \l t\r)^{\M^2-\f12\M}(\La_0+\Delta_0^\frac{5}{4}).
\end{equation*}
Due to (\ref{12.22.1.24}), (\ref{3.16.7.24}),  Proposition \ref{10.30.4.21} (2), Proposition \ref{7.15.5.22} and (\ref{2.14.1.24}), we infer
\begin{equation*}
\|\tir \bb^{-1}\Omega^3\tr\chi\|_{L^2_u L_\omega^2}\les (\log \l t\r)^{\M^2-\f12\M}(\La_0+\Delta_0^\frac{5}{4}).
\end{equation*}
Thus we proved the second estimate in (\ref{3.19.2.24}). 

Due to (\ref{3.18.3.24}), (\ref{3.18.13.24}) and (\ref{1.27.5.24}), 
\begin{equation}\label{6.17.2.24}
\|\tir^2\Omega(\bb^{-1}\sn\Xi_4)\|_{L_\omega^4}\les \log \l t\r^{\frac{3\M^2}{4}}(\La_0+\Delta_0^\frac{5}{4}).
\end{equation}
Applying (\ref{L4conn'}), (\ref{9.14.3.22}), (\ref{6.17.2.24}), (\ref{3.18.13.24}), in view of (\ref{9.6.3.22}), we derive
\begin{align*}
\|\tir^2 \sF\|_{L_\omega^4}&\les \log \l t\r^{\f12\M+1}(\Delta_0^\frac{5}{4}
+\La_0)\\
\|\tir^3(\tir\sn)\sF\|_{L_\omega^4} &\les\La_0+\int_0^t \tir^2 \{\|(\tir \sn) (\bb^{-1}\sn \Xi_4)\|_{L_\omega^4}\\
&+\l t\r^{-1}\Delta_0^\f12(\|\bb^{-2}(\tir \sn) \Xi_4\|_{L_\omega^\infty}+\|\tir^2\sn\sF\|_{L_\omega^4})\}+\l t\r^{\frac{1}{4}+2\delta}\Delta_0^2\\ &+\|\tir^2\bb^{-1}(\tir\sn)^{\ell}(\sn[L\Phi])\|_{L_t^1 L_\omega^4}+\int_0^t \l t\r^{1+\delta}\Delta_0\|\sF\|_{L_\omega^4}\\
&\les\log \l t\r^{\frac{3\M^2}{4}}\l t\r(\La_0+\Delta_0^\frac{5}{4})\end{align*}
where to obtain the last line, we ran an auxiliary bootstrap argument, by using the smallness of $\Delta_0$.
 Hence, in view of (\ref{6.17.2.24}), we have obtained
\begin{align*}
\|\tir^2(\tir\sn)(\bb^{-1}\sn\tr\chi)\|_{L_\omega^4}\les\log \l t\r^{\frac{3\M^2}{4}}(\La_0+\Delta_0^\frac{5}{4}).
\end{align*}
Using (\ref{L4conn}) and (\ref{3.6.2.21}), we can obtain the $L^4$ estimate in (\ref{3.19.2.24}) by using the above estimate.

Using the pointwise estimate of $\ze$ in (\ref{3.19.4.24}), Proposition \ref{1steng}, the estimate of $\Omega\tr\chi$ in Proposition \ref{3.14.4.24} and Proposition \ref{9.8.6.22}, we improve the following estimates from (\ref{5.25.1.21}) to
\begin{align*}
\bJ[\Omega]_L&=O(\l t\r^{-2}\log \l t\r^{\M^2}(\La_0+\Delta_0^\frac{5}{4}) )_{L_u^2 L_\omega^2}\\
\bJ[\Omega]_A&=O(\fB)+O(\l t\r^{-1+\delta}\Delta_0)+O(\Delta_0\l t\r^{-2+2\delta})_{L_\omega^4}
\end{align*}
 and improve (\ref{2.14.3.24}) to
\begin{align*}
\sP[\Omega, \Phi]&+\f12 {}\rp{a}\pih_{LA}\bd^2_{\Lb A} \Phi+\f12\Lb \Phi \bJ[\Omega]_L\\
&= O(\log \l t\r^{\frac{3\M^2}{4}}(\La_0+\Delta_0^\frac{5}{4}))(\sn^2 \Phi+\bd^2_{LL}\Phi)\\
& +O(\l t\r^{-\frac{3}{4}+\delta}\Delta_0^\f12)\bd^2_{LA} \Phi+O(\l t\r^{-1+\delta}\Delta_0)_{L_\omega^4}( L \Phi+\l t\r^{-1+\delta}\sn\Phi)\nn\\
&+O(\fB+\l t\r^{-1+\delta}\Delta_0)\sn \Phi\\
&=O(\l t\r^{-2+\delta}\log \l t\r^{\f12\M}(\La_0+\Delta_0^\frac{5}{4}))_{L^2_\Sigma};
\end{align*}
and
\begin{align*}
\Lb \Phi \bJ[\Omega]_L&=O(\l t\r^{-2}\log \l t\r^{\M^2}(\La_0+\Delta_0^\frac{5}{4}))_{L^2_\Sigma};\\
{}\rp{a}\pih_{LA}\bd^2_{\Lb A} \Phi&=O(\log \l t\r^{\f12\M+\frac{3}{2}}\l t\r^{-2+\delta}(\Delta_0^\frac{5}{4}+\La_0))_{L^2_\Sigma} 
\end{align*}
where the last line follows from using (\ref{7.26.2.22}) and (\ref{12.25.3.23}).

Moreover, we can improve (\ref{6.6.1.23}) by using Proposition \ref{1steng} and Proposition \ref{9.8.6.22} to be
\begin{equation*}
X\Box_\bg \Phi=O(\l t\r^{-2}\log \l t\r^{\f12\M+1}(\Delta_0^\frac{5}{4}+\La_0))_{L^2_\Sigma}.
 \end{equation*}
 Combining the above estimates, we conclude (\ref{wave_ass'}).

\end{proof}

(\ref{3.12.1.21}) is improved by (\ref{8.21.4.21}), (\ref{8.25.2.21}), Proposition \ref{3.14.4.24}, and the first three estimates in Proposition \ref{10.30.4.21} (1); (\ref{L2BA2}) and (\ref{LbBA2}) are improved to (\ref{3.16.7.24}), (\ref{3.16.13.24}) and (\ref{LbBA2'}) with the lowest order estimate included in (\ref{8.21.4.21}); (\ref{L4BA1}) is improved by  (\ref{10.11.2.23}), (\ref{3.15.3.24}), (\ref{3.16.3.24}), (\ref{3.18.2.24}), (\ref{3.18.3.24}), (\ref{L4conn'}), the last estimate in Proposition \ref{3.14.4.24}, (\ref{9.14.3.22}) and the $|\Phi|$ bound in (\ref{7.25.2.22}); (\ref{L4conn}) is improved to (\ref{9.11.3.22}) ; (\ref{L2conndrv}) is improved to (\ref{L2conn'}), the third and the fourth estimates in Proposition \ref{3.14.4.24}, (\ref{8.8.6.22'})  and (\ref{L2conndrv'}).
(\ref{ConnH}) is improved by (\ref{3.19.2.24}); (\ref{zeh}) is improved by (\ref{3.19.4.24}); (\ref{wave_ass}) is improved to (\ref{wave_ass'}).

 The improvement can be seen in that the growth rate $\l t\r^\delta$ is improved by $(\log \l t\r)^{\M^2+1}$, with $15\le \M\approx \A_0$ to be confirmed in Section \ref{10.24.1.23}; and $\Delta_0$ in the assumptions is proved to be $C(\Delta_0^\frac{5}{2}+\La_0^2)^\f12$. In particular, we can choose $\Delta_0=C_6 \La_0$, with $C_6=2C$ and $0<\La_0<\frac{1}{8C^5}$, then $C^2(\Delta_0^\frac{5}{2}+\La_0^2)<\frac{3\Delta_0^2}{4}<\Delta_0^2$. The smallness of $\La_0$ is achieved by choosing $\ve_0$ sufficiently small, which is illustrated in Section \ref{9.25.2.22}. From now on, when employing the estimates from Proposition \ref{1steng}, Proposition \ref{8.29.8.21}, Proposition \ref{9.8.6.22}, Proposition \ref{10.30.4.21}, Proposition \ref{3.14.4.24}, Proposition \ref{imp_decay}, etc, we will use  $\Delta_0\approx \La_0$ to simplify the form of the bounds. 

Note that  (\ref{6.20.2.21+}) has been derived in the proof of Lemma \ref{5.13.11.21} as an improvement of (\ref{6.20.2.21}), (\ref{3.7.1.21-}) is improved by (\ref{3.17.1.24}) to (\ref{3.7.1.21}); (\ref{1.25.1.22}) has been improved by Lemma \ref{3.19.5.24} (3) and  Proposition \ref{10.16.1.22}. With the improvement of $\Delta_0$ to $C(\Delta_0^\frac{5}{4}+\La_0)$, we can obtain the estimates in (\ref{9.14.1.24}). It only remains to show (\ref{6.5.1.21_+}), which will be completed in Section \ref{10.24.1.23} after controlling the total energy for $\forall t<T_*$ in Section \ref{10.24.2.23}. 
\section{The continuation of the solution to $T_*$}\label{10.24.2.23}
In this section, we bound $\bp^{\le 3}\Phi\in L^2(\Sigma_t)$ uniformly for any $t<T_*<\infty$. By the standard local well-posedness result, the solution continues uniquely beyond $T_*$. 

Due to  (\ref{6.24.1.21}) and $\tir\bb \fB=O(1)$  in (\ref{3.11.4.21}), we have obtained the lower order estimates
\begin{equation*}
\|\bp\Phi\|_{L^2_\Sigma}\les 1,\quad |\Phi|\les \l t\r^{-1+\delta}.
\end{equation*}
 Next we control the higher order total energies. 
\begin{proposition}[Control of Total energies]\label{10.13.5.23}
Assuming (\ref{1.12.1.22})-(\ref{exist}) \begin{footnote}{Recall from Proposition \ref{12.21.1.21} that $\La_0\approx\ve$.}\end{footnote} and (\ref{6.5.1.21}), there hold the following estimates for any $0<t<T_*$ that
\begin{align}
&\|\bb^3 \bp^2\Phi\|_{L_\omega^4}\les\l t\r^{-1},  \|\bp^2\Phi\|_{L^2_\Sigma}\les 1\label{10.12.2.23}\\
& E[\bp^2\Phi](t)\les \log \l t\r^{\M+3}, \quad E[\bp^3\Phi](t)\les \log \l t\r^{2\M}\label{9.1.1.22}
\end{align}
as long as $\M>8$.
\end{proposition}
\begin{proof}
 By using (\ref{3.6.1.21}), (\ref{6.22.1.21}), Proposition \ref{9.8.6.22}, we first prove (\ref{10.12.2.23}).
Note by using (\ref{9.14.3.22}), (\ref{9.19.5.23}), (\ref{6.22.1.21}) and (\ref{3.6.1.21}), we derive 
$$
\|\bb^3\p \bT \varrho\|_{L_\omega^4}\les\|\bb^3(\Lb \Lb\varrho, L^2\varrho, L\Lb\varrho, \Lb L\varrho)\|_{L_\omega^4}+\|\bb^3\sn \bT\varrho\|_{L_\omega^4}\les \l t\r^{-1}.
$$
Next we decompose a scalar function by using the radial frame on $\Sigma$ 
\begin{align}\label{8.31.4.22}
|\nab^2 f|&\les |\sn^2 f|+ |\sn_\bN \sn f, \sn\sn_\bN f|+|\bN\bN f|+|\theta\nab f|+|\ud\bA\sn f|.
\end{align}
Applying the above formula to $f=\varrho$ and using (\ref{3.6.1.21}), (\ref{6.22.1.21}), Proposition \ref{9.8.6.22} and the symbolic identity $\bN \bN f=\Lb \Lb f+\Lb L f+L \Lb f+L^2 f$ gives  $$\|\bb^3 \nab^2 \varrho\|_{L_\omega^4}\les \l t\r^{-1}.$$
Note due to (\ref{4.23.1.19}), schematically, $\p\bT \Phi, \bT \p \Phi=\p^2 \varrho+\p \bT \varrho+ \p\Phi\bp\Phi$ where we dropped the factors of $c^m$ paired to for each terms. Hence using the above estimates of $\nab^2 \varrho$, $\p \bT\varrho$ and $\bb\fB=O(\l t\r^{-1})$, we infer
\begin{equation*}
\bb^3 (\p\bT \Phi, \bT \p \Phi)=O(\l t\r^{-1})_{L_\omega^4}.
\end{equation*}
Thus the case of $\Phi=\varrho$ in (\ref{10.12.2.23}), and $\Lb \Lb \Phi, \bT \bT\Phi$ and fully tangential derivatives of $\Phi$ in (\ref{10.12.2.23}) have been proved. It only remains to check $\nab^2 v$.

Similar to the treatment of $\nab^2\varrho$, using (\ref{8.31.4.22}), (\ref{3.19.4.24}), (\ref{8.23.1.23}) and  (\ref{6.22.1.21})  we derive
\begin{align*}
\|\bb^3\nab^2 v\|_{L_\omega^4}&\les \|\bb^3(\sn^2 v, L^2 v)\|_{L_\omega^4}+\|\bb^3(\sn_\bN \sn v, \sn\sn_\bN v, \Lb^2 v, \Lb L v, L\Lb v)\|_{L_\omega^4}+\l t\r^{-1+\frac{\delta}{2}}\|\bb^3\nab v\|_{L_\omega^4}\\
&\les \l t\r^{-1}
\end{align*}
Hence the first estimate in (\ref{10.12.2.23}) is proved and the second is a consequence. 

To prove (\ref{9.1.1.22}) we rely on the following formula
\begin{equation}
\Box_\bg \Phi=(\bT \Phi)^2+\nab \Phi\c\nab \Phi\label{9.19.7.23}
 \end{equation}
 which is recast from (\ref{4.10.1.19}) and (\ref{4.10.2.19}) symbolically. 

To prove the energy estimate with $n=2$ in (\ref{9.1.1.22}), using (\ref{9.19.7.23}), it is direct to obtain that
\begin{align}
\Box_\bg \bp^2\Phi&=\bg\c \Big(\bp\Phi \bp^3 \Phi+\bp^2 \Phi\c(\bp^2\Phi+\bp\Phi\c \bp\Phi)\Big)\label{7.4.1.24}\\
\Box_\bg \bp^3 \Phi&=\bg\c \Big(\bp^4\Phi\bp\Phi+\bp^3 \Phi((\bp\Phi)^2+\bp^2\Phi)+\bp^2\Phi(\bp^2\Phi\bp\Phi+(\bp\Phi)^3)\Big)\label{7.4.2.24}.
\end{align}
Applying (\ref{6.21.2.21}) to $\psi=\bp^2\Phi$ with the help of (\ref{7.4.1.24}), $O(\bb\bp\Phi)=\l t\r^{-1}$ and  (\ref{10.12.2.23}), we obtain by using (\ref{11.11.2.23}), Sobolev embedding and Gronwall's inequality 
\begin{align*}
E[\bp^2\Phi](t)&\les E[\bp^2\Phi](0)+\int_{\D_{u_0}^{t}}|\bT \bp^2 \Phi||\Box_\bg \bp^2\Phi|\\
 &\les E[\bp^2\Phi](0)+\int_{\D_{u_0}^{t}}|\bT \bp^2\Phi| (|\bp\Phi \bp^3 \Phi|+|\bp^2 \Phi\c(\bp^2\Phi+\bp\Phi\c \bp\Phi)|)\\
 &\les \log \l t\r^{\M+3}.
\end{align*}
Similarly, for the case $n=3$, using (\ref{10.12.2.23}), (\ref{7.4.2.24}) and the first estimate in (\ref{9.1.1.22}), we derive by using (\ref{11.11.2.23}), Sobolev embedding and Gronwall's inequality that
\begin{align*}
E[\bp^3 \Phi](t)&\les E[\bp^3\Phi](0)+\int_{\D_{u_0}^{t}}|\bT \bp^3 \Phi||\Box_\bg \bp^3\Phi|\\
&\les E[\bp^3\Phi](0)+\int_{\D_{u_0}^{t}}|\bT \bp^3\Phi||\bp^4\Phi\bp\Phi+\bp^3 \Phi((\bp\Phi)^2+\bp^2\Phi)+\bp^2\Phi(\bp^2\Phi\bp\Phi+(\bp\Phi)^3)|\\
&\les \log \l t\r^{\M+8}.
\end{align*}
Hence, the proof of (\ref{9.1.1.22}) is complete. 
\end{proof}

\section{The Completion of the bootstrap argument and the formation of rarefaction}\label{10.24.1.23}
\subsection{Control of $\Lb \varrho$}
In this subsection, we improve the bootstrap assumption (\ref{6.5.1.21}) to (\ref{6.5.1.21_+}), which will be completed in Proposition \ref{11.9.1.23} and Proposition \ref{10.9.3.22}.

Since all other estimates in Assumption \ref{5.13.11.21+} have been shown to hold beyond $T_*$, we recall (\ref{6.5.1.21}) in Assumption \ref{5.13.11.21+} that   $T_*=\sup\{t':\mbox{ the following inequality holds  on } (0, t')\}$ with a constant $\M_0>1$ depending on $\A_0$ and $c_*$, 
\begin{equation*}
-\varrho< -\tir \Lb \varrho\le \M_0 \big(1+\f12\wp\log(\frac{\l t\r}{2})\big)^{-1}.
\end{equation*}
  By fixing $\M_0=2\ti C(\mathfrak{C}+1)$, with the universal constants $\ti C, \mathfrak{C}$ specified in Proposition \ref{10.9.3.22}, (\ref{6.5.1.21}) will be improved by (\ref{6.5.1.21_+}), shown in Proposition \ref{11.9.1.23} (1) and Proposition \ref{10.9.3.22}.

We first prove
\begin{lemma}\label{11.10.3.23}
Under the assumptions (\ref{1.12.1.22})-(\ref{exist}) and the assumption (\ref{6.5.1.21}), there holds for $y=\tir \Lb\varrho$ that 
\begin{equation*}
Ly=\f12 \wp\tir^{-1}y^2+G, \quad 0<t<T_*
\end{equation*}
with
\begin{align}
G
&=\tir\Big((\frac{\wp-2}{2}L \varrho-L v_\bN+\bA_b)\Lb \varrho+\sD\varrho-\Box_\bg \varrho+2\zb^A \sn_A \varrho-\hb L \varrho\Big),\label{10.9.2.22}
\end{align}
and  
\begin{equation}\label{9.21.4.23}
|G|\le C_3\l t\r^{-\frac{7}{4}+\delta}\ve^\f12
\end{equation}
where $C_3>1$ is a universal constant. 
\end{lemma}
\begin{proof}
Recall from (\ref{6.30.2.19}) that
\begin{equation*}
L\Lb \varrho-k_{\bN\bN} \Lb \varrho=\sD\varrho-h\Lb \varrho-\hb L\varrho+2\zb^A \sn_A \varrho-\Box_\bg \varrho.
\end{equation*}
Also by using (\ref{3.22.1.21}) and decomposing $2\bT=L+\Lb$ , we calculate
\begin{align*}
L(\tir \Lb \varrho)&=\tir L \Lb \varrho+\Lb \varrho\\
&=\Lb \varrho+\tir((k_{\bN\bN}-h)\Lb\varrho-\hb L\varrho+2\zb^A \sn_A \varrho+\sD\varrho-\Box_\bg \varrho)\\
&=\Lb \varrho(1- h\tir)+\tir\big((\wp \bT \varrho-L \varrho-L v_\bN)\Lb \varrho-\hb L \varrho+2\zb \sn\varrho+\sD\varrho-\Box_\bg \varrho\big)\\
\displaybreak[0]
&=\f12\wp\tir (\Lb \varrho)^2+\Lb \varrho(1-h\tir)+\tir\cdot \Big(-\hb L \varrho+2\zb\sn \varrho\\
&+\sD\varrho-\Box_\bg \varrho+(\frac{\wp-2}{2}L \varrho-Lv_\bN)\Lb \varrho\Big).
\end{align*}
We then conclude (\ref{10.9.2.22}).
 The bound of $G$ follows by using the decay estimates in (\ref{7.25.2.22}), (\ref{3.18.2.24}), (\ref{3.17.1.24}), Proposition \ref{geonul_5.23_23} and $|\bb\tir\Lb \varrho|\les 1$.  
\end{proof}

  The improvement of the left-hand side inequality of (\ref{6.5.1.21}) is achieved by showing  the following result.
\begin{proposition}\label{11.9.1.23}
Under the assumptions (\ref{1.12.1.22})-(\ref{exist}) and  (\ref{6.5.1.21}), we have for $y=\tir \Lb \varrho$ that 
\begin{enumerate}
\item there is a positive-valued function $\bbf\approx 1$  and the universal constants $C, C_1, C_2>1$ such that 
\begin{equation}\label{9.29.3.23}
-(\bbf^{-1} (y-\varrho)+\a_1)> C_2^{-1}\big(-f^{-1}+\f12 C_1 \wp \log (\f12\l t\r)\big)^{-1},\quad 0<t<T_*
\end{equation}
where $\a_1=2C\ve \l t\r^{-\frac{3}{4}+\delta}$, and $f(0,u,\omega)=(y-\varrho+\a_1)(0)<0$. This gives
\begin{align*}
\a_1\les -(y-\varrho), \quad  0\le t<T_*.
\end{align*}
\item there holds 
\begin{align}\label{11.9.7.23}
 \varrho\le-C_8\l t\r^{-1}(\int_u^{u_*}\inf_\omega C(t,u, \omega)+\a_1(t))+C_7\ve\l t\r^{-1+\f12\delta}\log \l t\r^{\frac{\M}{4}} \quad 0<t<T_*
\end{align}
with $C(t,u,\omega)$ the right-hand side of (\ref{9.29.3.23}) and $C_7, C_8>0$  positive universal constants.
\end{enumerate}
\end{proposition}
\begin{proof}
With $y=\tir \Lb \varrho$, using $h+\hb=-2\tr k$, (\ref{k1}) and $\tr\eta=[L\Phi]$,
 we first recast (\ref{10.9.2.22}) below
\begin{align*}
Ly-L\varrho-\f12\wp \tir^{-1}y^2&=y([L\Phi]+\bA_b)+\tir([L\Phi]^2+ L\varrho\bA_b+\zb\sn \varrho+\sD\varrho-\Box_\bg \varrho).
\end{align*}
Using (\ref{null}) and the fact that $[\Lb\Phi]=\Lb\varrho+[L\Phi]$, symbolically we derive 
\begin{equation}\label{9.29.1.23}
\Box_\bg \varrho=\Lb \varrho\c[L\Phi]+[\sn\Phi]^2+|\eh|^2+[L\Phi]^2.
\end{equation}
Combining the above two formulas yields 
\begin{align}\label{11.8.1.23}
Ly-L\varrho-\f12 \wp y^2\tir^{-1}=(\bA_b+[L\Phi])y+G_2,
\end{align}
where we keep the terms on the left-hand side precise, and 
\begin{align*}
G_2&=\tir([L\Phi](\bA_b+[L\Phi])+\bA_{g,1}^2+\sD\varrho)\approx O(\ve) \l t\r^{-\frac{7}{4}+\delta},
\end{align*} 
where  the last estimate is derived in the same way as for (\ref{9.21.4.23}). 

We then rewrite 
\begin{align*}
L(y-\varrho)-\f12 \tir^{-1}\wp(y-\varrho)^2=(\tir^{-1}\wp\varrho+\bA_b+[L\Phi])(y-\varrho)+\wt{G_2}
\end{align*}
where
\begin{align*}
\wt{G_2}=\varrho(\bA_b+[L\Phi])+G_2+\frac{1}{2}\tir^{-1}\wp{\varrho}^2=O(\ve\l t\r^{-\frac{7}{4}+\delta})
\end{align*}
with the last estimate obtained by using (\ref{7.25.2.22}), (\ref{3.17.1.24}) and the bound of $G_2$. 

Consequently, with $\bbf:=\exp\Big(\int_0^t (\tir^{-1}\wp\varrho+\bA_b+[L\Phi]) dt'\Big)=\exp(O(\ve^\f12))$, we obtain 
\begin{equation*}
L(\bbf^{-1}(y-\varrho))=\f12\bbf\wp\big(\bbf^{-1}(y-\varrho)\big)^2\tir^{-1}+\bbf^{-1}\wt{G_2}. 
\end{equation*}
Following the scheme in Section \ref{5.9.1.23}, by setting $\ti z=\bbf^{-1}(y-\varrho)+\a_1$, we write
\begin{align}\label{9.29.2.23}
L(\bbf^{-1}(y-\varrho)+\a_1)=\f12 \bbf \wp\ti z^2 \tir^{-1}-\bbf \a_1\wp \ti z \tir^{-1}+\f12 \bbf \wp \a_1^2 \tir^{-1}+L\a_1+\bbf^{-1}\wt{G_2}.
\end{align} 
Since $\La_0\approx \ve$, we choose $\a_1=2C\ve \l t\r^{-\frac{3}{4}+\delta}$ such that 
\begin{align}\label{11.9.6.23}
(y-\varrho+\a_1)(0)=f(0)<0, \quad \f12 \bbf \wp \a_1^2 \tir^{-1}+L\a_1+\bbf^{-1}\wt{G_2}<0.
\end{align}
This is achievable if the universal constant $C>1$ is sufficiently large and $\ve$, no greater than a constant fraction of $q_0$, is sufficiently small.  
  Then we derive from (\ref{9.29.2.23}) that
\begin{equation}\label{11.10.4.23}
\ti z'-\f12\bbf \wp \ti z^2 \tir^{-1}+\bbf \a_1\wp \ti z\tir^{-1}<0.
\end{equation}
Suppose $t_1=\sup\{0\le  t<T_*: \ti z<0 \mbox{ on } (0, t)\}$. By continuity, $t_1>0$. Next we will show $t_1=T_*$. 

With $0<t<t_1$, multiplying (\ref{11.10.4.23}) by $-\ti z^{-2}$, followed with integrating along a null geodesic $\Upsilon_{u,\omega}$, gives
\begin{align*}
\ti z^{-1}\exp(\int_0^t -\bbf \a_1 \wp \tir^{-1})&>\ti z^{-1}(0)-\f12\int_0^t \bbf \wp {\tir'}^{-1}\exp(\int_0^{t'} -\bbf \a_1 \wp \tir^{-1})\\
&>\ti z^{-1}(0)-\f12 C_1\wp \log(\dfrac{\l t\r}{2}).
\end{align*}
Hence
\begin{equation*}
-\ti z\exp(\int_0^t \bbf \a_1\wp \tir^{-1})>\big(-\ti z^{-1}(0)+\f12 C_1\wp \log (\f12\l t\r)\big)^{-1}.
\end{equation*}
Noting that $\exp(\int_0^t \bbf \a_1\wp \tir^{-1})=\exp(O(\ve))<C_2\approx 1$, we derive
\begin{equation*}
-\ti z>C_2^{-1}(-f(0,u,\omega)^{-1}+\f12 C_1\wp \log (\f12\l t\r))^{-1},\quad 0<t<t_1.
\end{equation*}
Denote the right-hand side by $C(t,u,\omega)$, which is a positive-valued, decreasing function about $t$ with $u, \omega$ fixed. By compactness of $[u_0, u_*]$ and continuity of $f(0,u,\omega)$, $-f(0,u,\omega)>0$ has strictly positive lower and upper bounds on $[u_0, u_*]\times {\mathbb S^2}$. This shows $t_1=T_*$ and $-\ti z$ is uniformly bounded from below by a strictly positive number for $t<T_*, u_0\le u\le u_*$.
This gives (\ref{9.29.3.23}) and $\a_1\les -y+\varrho$ for $0\le t<T_*$. Thus the proof of (1) is complete.

Next we prove (\ref{11.9.7.23}).
Recall that $\bb \bN=\p_u-\b_A \frac{\p}{\p\omega^A}$. We directly compute
\begin{align}
-2\bb \bN (\stackrel{\circ}v_t^\frac{1}{4}\varrho)&=-2\bb\stackrel{\circ}v_t^\frac{1}{4}(\bN \varrho+\frac{\varrho}{2\tir}+\bA_b\varrho).\label{11.9.3.23}
\end{align} 
Using (\ref{L2BA2'}) and Sobolev embedding, (\ref{7.25.2.22}), (\ref{9.19.5.23}), (\ref{3.28.3.21}) and (\ref{3.17.1.24}), we bound
\begin{equation}\label{11.9.4.23}
\|L\varrho\|_{L_u^2 L_\omega^\infty}\les \l t\r^{-2}\log \l t\r^{\f12\M}\ve, |\b\sn\varrho|\les \l t\r^{-1+\frac{\delta}{2}}(\log \l t\r)^{\frac{\M}{4}}\ve^2,\, |\bb\bA_b \varrho|_{L_u^2 L_\omega^\infty}\les \l t\r^{-\frac{11}{4}+\delta}\ve.
\end{equation}
Integrating (\ref{11.9.3.23}) from $u$ to $u_*$, using $\vc_t\approx \l t\r^2$, (\ref{11.9.3.23}) and the last two estimates in (\ref{11.9.4.23}), we obtain
\begin{align*}
2\stackrel{\circ}v_t^\frac{1}{4}\varrho&=\int_u^{u_*}-2\p_u\left( \stackrel{\circ}v_t^\frac{1}{4}\varrho\right)=\int_u^{u_*}-2\bb \stackrel{\circ}v_t^\frac{1}{4}(\bN \varrho+\frac{1}{2\tir}\varrho)+O(\ve)\l t\r^{-\f12+\f12\delta}\log \l t\r^{\frac{\M}{4}}.
\end{align*}
Using (\ref{9.29.3.23}) we derive in view of $L-\Lb=2\bN$ that 
 \begin{align*}
-2 \tir(\bN \varrho+\f12\tir^{-1} \varrho)<-\bbf \a_1-\tir L\varrho-\bbf C(t,u,\omega)
 \end{align*}
 where $C(t,u,\omega)$ is the right-hand side of (\ref{9.29.3.23}). 
  
  Combining the above two estimates, also using (\ref{11.13.3.23}), leads to
 \begin{align*}
\stackrel{\circ}v_t^\frac{1}{4}\varrho\les\l t\r^{-\f12}\int_u^{u_*}\{-\bbf \a_1-\bbf C(t,u,\omega)\}+\log \l t \r\l t\r^{\f12}\int_u^{u_*} | L\varrho|+O(\ve)\l t\r^{-\f12+\f12\delta}\log \l t\r^{\frac{\M}{4}}.
 \end{align*}
 Hence using the first estimate in (\ref{11.9.4.23}), substituting the formula of $\a_1$ and using the fact that $\bbf \approx 1$, we have
 \begin{align*}
 \stackrel{\circ}v_t^\frac{1}{4}\varrho\les-\l t\r^{-\f12} (\int_u^{u_*}\inf_\omega C(t,u, \omega)+\a_1(t))+O(\ve)\l t\r^{-\f12+\f12\delta}\log \l t\r^{\frac{\M}{4}}.
 \end{align*}
Hence we obtained (\ref{11.9.7.23}) due to $\stackrel{\circ}v_t^\frac{1}{4}\approx \tir^\f12$. The proof of Proposition \ref{11.9.1.23} is complete. 
\end{proof}

The right-hand side inequality in (\ref{6.5.1.21}) will be improved by the following result
\begin{proposition}\label{10.9.3.22}
Under the assumptions (\ref{1.12.1.22})-(\ref{exist}), there are constants $C_4, C_5, C_6, \mathfrak{C}\ge1$ such that for $t<T_*$ and $u_0\le u\le u_*$  
\begin{align}
&-\a<-\tir \Lb \varrho\le C_5 \left(\frac{1}{\mathfrak{C}+3C\ve^\f12}+\f12 C_4^{-1}\wp \log(\f12 \l t\r)\right)^{-1},\label{11.10.1.23}\\
&-\l t\r^{-1}\log \l t\r^{\frac{\M}{2}}\ve^\f12\les -\tir \Lb \varrho\les \left(\frac{1}{\mathfrak{C}+3C\ve^\f12}+\f12 C_4^{-1}\wp \log(\f12 \l t\r)\right)^{-1}\label{11.10.2.23}
\end{align}
where $\a=3 C_6 \l t\r^{-\frac{3}{4}+\delta} \ve^\f12$ and the constant $C_6\ge \max(C_3, \frac{C_7}{3})$ with $C_3, C_7$ the constants appeared in (\ref{9.21.4.23}) and (\ref{11.9.7.23})
\end{proposition}
\begin{remark}\label{7.3.1.24}
From the right-hand side inequality in (\ref{11.10.1.23}), we derive
\begin{align*}
-\tir \Lb \varrho&\le C_5 C_4\big(\frac{C_4}{\mathfrak{C}+1}+\f12 \wp\log(\f12 \l t\r)\big)^{-1}\le C_5 C_4(\mathfrak{C}+1)(1+\f12(\mathfrak{C}+1)\wp \log(\f12\l t\r))^{-1}\\
 &\le \ti C (\mathfrak{C}+1)\left(1+\f12\wp \log(\f12\l t\r)\right)^{-1}
\end{align*}
which is the right-hand side inequality (\ref{6.5.1.21_+}) with $\ti C=C_4 C_5$. The left-hand side inequality of (\ref{6.5.1.21_+}) has been proved in Proposition \ref{11.9.1.23} (1). The proof of (\ref{6.5.1.21_+}) is therefore complete. 
\end{remark}
\begin{proof}
 Recall from (\ref{9.21.4.23}) that  $|G|\le C_3\ve^\f12 \l t\r^{-\frac{7}{4}+\delta},$ for $G$ in (\ref{10.9.2.22}).
Let
$$y=\tir \Lb \varrho, \quad\a=3 C_6 \l t\r^{-\frac{3}{4}+\delta} \ve^\f12.$$ We recast the asymptotic equation in Lemma \ref{11.10.3.23} 
\begin{align*}
L(-y+\a)&=-\f12 \wp \tir^{-1}(-y)^2+L\a-G\\
&=-\f12\wp\tir^{-1}\big((-y+\a)^2-\a^2+2y \a\big)+L\a-G\\
&=-\f12 \wp \tir^{-1}\big((-y+\a)^2+2(y-\a)\a\big)+L \a-G-\f12\wp\tir^{-1}\a^2.
\end{align*}
It is straightforward to see that $L\a-G-\f12 \wp \tir^{-1}\a^2<0$. Hence
\begin{equation*}
L(-y+\a)<-\f12 \wp \tir^{-1}\big((-y+\a)^2+2(y-\a)\a\big).
\end{equation*}
Let $z=-y+\a$. Then we obtain from the above inequality that
\begin{equation*}
Lz<\f12 \wp \tir^{-1}z(-z+2\a).
\end{equation*}
Recasting (\ref{11.9.7.23}) for $0<t<T_*$
\begin{align*}
\varrho=-y+\varrho+\a-z\le-C_8\l t\r^{-1}(\int_u^{u_*}\inf_\omega C(t,u, \omega)+\a_1(t))+C_7\ve\l t\r^{-1+\f12\delta}\log \l t\r^{\frac{\M}{4}},
\end{align*}
followed with rearranging terms in view of  the result that $-y+\varrho\ges\a_1>0$ in Proposition \ref{11.9.1.23}, we obtain
\begin{align*}
z\ge C_9^{-1} \a_1-C_7\ve\l t\r^{-1+\f12\delta}\log \l t\r^{\frac{\M}{4}}+\a\ges\a_1>0
\end{align*}
due to the choice of $\a$, where $C_9>1$ is a constant. 

Therefore
\begin{align*}
L(z^{-1})&=-z^{-2}Lz>-z^{-2}\wp \tir^{-1}z(-\f12 z+\a)\\
&=-\wp \tir^{-1}(-\f12+z^{-1}\a).
\end{align*}
We then derive that
\begin{equation*}
L(z^{-1}\exp\int_0^t \wp \a \tir^{-1})>\f12 \wp \tir^{-1}\exp\int_0^t \wp \a \tir^{-1}, \quad t<T_*.
\end{equation*}
Integrating from $S_{0,u}$ along $\H_u$ in $t$ leads to 
\begin{align*}
z^{-1}\exp \int_0^t \wp \a \tir^{-1} >z^{-1}(0)+\f12\wp \int_0^t {\tir'}^{-1}\exp\int_0^{t'} \wp \a \tir^{-1}
\end{align*}
which implies
\begin{align*}
\frac{\exp \int_0^t \wp \a \tir^{-1}}{-\tir \Lb \varrho+\a}>\frac{1}{(-\tir\Lb \varrho)(0)+3C_6\ve^\f12}+\f12 C_4^{-1}\wp \int_0^t \l t'\r^{-1},
\end{align*}
where we used $\tir\le  C_4\l t\r$ and $C_4\ge 1$ is a universal constant, depending on $c_*$. Due to continuity of $\ckc^{-1}\Lb \varrho$ and compactness of the initial slice, also using (\ref{9.30.15.23}), we have $-\mathfrak{C}\le \ckc^{-1}u\Lb \varrho(0)\les\ve$ on $u_0\le u\le u_*$, with $\mathfrak{C}>1$ a certain constant. Hence
\begin{equation*}
\frac{\exp \int_0^t \wp \a \tir^{-1}}{-\tir\Lb \varrho+\a}>\frac{1}{\mathfrak{C}+3C_6\ve^\f12}+\f12 C_4^{-1}\wp \int_0^t \l t'\r^{-1}.
\end{equation*}
Hence we conclude
\begin{align*}
0<-\tir \Lb\varrho+\a&< \left(\frac{1}{\mathfrak{C}+3C_6\ve^\f12}+\f12 C_4^{-1}\wp \log(\f12 \l t\r)\right)^{-1}\exp\int_0^t \wp \a \tir^{-1}\\
&<C_5 \left(\frac{1}{\mathfrak{C}+3C_6\ve^\f12}+\f12 C_4^{-1}\wp \log(\f12 \l t\r)\right)^{-1}.
\end{align*}
This gives (\ref{11.10.1.23}) for $t<T_*$. Consequently
\begin{align*}
-\l t\r^{-\frac{3}{4}+\delta}\ve^\f12\les -\tir \Lb \varrho\les (\frac{1}{\mathfrak{C}+3C_6\ve^\f12}+\f12 C_4^{-1}\wp \log(\f12 \l t\r))^{-1}.
\end{align*}
We can further use Proposition \ref{11.9.1.23} (1) and (\ref{7.25.2.22}) to improve the left-hand side inequality. Thus (\ref{11.10.2.23}) is proved.  
\end{proof}


\subsection{The rarefaction at null infinity}\label{brate}
Assuming (\ref{1.12.1.22})-(\ref{rarif}), we will show  (2) in Theorem \ref{mainthm1} holds, i.e. the rarefaction forms at the null infinity. 

Using (\ref{11.28.4.23}), Proposition \ref{9.8.6.22} and Proposition \ref{imp_decay}, we have
\begin{equation*}
|L(\bb \tir \Lb \varrho)|\les \ve^\f12 \l t\r^{-\frac{7}{4}+\delta}.
\end{equation*}
Integrating the above inequality along the null geodesic $\Upsilon_{u,\omega}$ in $t$ implies
\begin{equation}\label{9.10.4.22}
(\ckc^{-1} c^{-1}r\Lb \varrho)(0,u,\omega)-C\ve^\f12<\bb \tir \Lb \varrho(t,u,\omega)< (\ckc^{-1} c^{-1}r\Lb \varrho)(0,u,\omega)+C\ve^\f12.
\end{equation}
Suppose there exist $u_1\in[u_0, u_*]$ and $\omega_1\in{\mathbb S^2}$ such that 
\begin{equation*}
(\ckc^{-1} c^{-1} r\Lb \varrho)(0,u_1,\omega_1)<-2C\ve^\f12.
\end{equation*}
This implies
$$
\bb \tir \Lb \varrho(t,u_1,\omega_1)<- C\ve^\f12,\quad \forall t>0,
$$
where $C>0$ a constant. Due to $\bb>0$,  the term $-\tir \Lb \varrho(t,u_1,\omega_1)>0$.  In view of (\ref{11.10.2.23}), we derive 
 \begin{equation*}
\bb(t, u_1, \omega_1)\ges \ve^\f12\big(\frac{1}{\mathfrak{C}+3C_6\ve^\f12}+\f12 C_4^{-1} \wp\log (\f12\l t\r)\big).
\end{equation*}
With $\ve_0$ sufficiently small, we obtain the growth rate stated in Theorem \ref{mainthm1} (2).  
Consequently, along the null geodesic $\Upsilon(t,u_1,\omega_1)$, $\bb\rightarrow \infty$ as $t\rightarrow \infty$, which confirms the formation of rarefaction of the flow at the null infinity.

\section{The Control of the initial slice}\label{9.25.2.22}
The purpose of this section is to prove Proposition \ref{12.21.1.21}. The main tool we rely on is the following Sobolev inequality. 
\begin{lemma}
 At $t=0$, there holds the Sobolev inequality on $\Sigma_0$ to $f$ which vanishes at $u=u_*$,
\begin{align}
\int_{S_{0, u}}|f|^4 r^2&\les \int_{\Sigma_0\cap [u,u_*]}|\hN^{\le 1} f|^2\int_{\Sigma_0\cap [u,u_*]}|\Omega^{\le 1}  f|^2.\label{9.22.2.22}
\end{align}
\end{lemma}
\begin{proof}
(\ref{9.22.2.22}) follows from the following Sobolev inequality.
\begin{equation*}
\int_{S_{0,u}}|f|^4 r^2\les \int_{S_{0,u_*}}|r f|^4 r^{-2} + \int_{\Sigma_0\cap [u,u_*]}|\hN(r f)|^2 r^{-2} \cdot \sum_{l\le 1}\int_{\Sigma_0\cap [u,u_*]}|\Omega^l  f|^2.
\end{equation*}
\end{proof}
\subsection{Preliminary estimates}
(\ref{9.25.1.22}) follows from using (\ref{1.12.1.22}) and Sobolev embedding. (\ref{12.20.4.21}) are facts of the standard radial foliation at the initial slice. In this part, we give the proof of (\ref{5.14.1.23})-(\ref{12.6.2.23}) in Proposition \ref{12.21.1.21} except the last estimate in (\ref{9.30.1.23}), and other preliminary estimates which can be obtained directly from (\ref{1.12.1.22})-(\ref{9.22.1.22}). 
 
 Recall at $t=0$, $\hat \bN=\p_u$, $\bN=c\p_u$, and $\bp\Phi=O(1)$ from (\ref{9.25.1.22}), which will be used as facts without mentioning explicitly.

Applying (\ref{9.22.2.22}) to $f=Y^l\Phi, l=1,2,3, Y=\Omega, S$, also using (\ref{9.22.1.22}) and (\ref{1.12.1.22}), yields
\begin{equation*} 
\int_{S_{0,u}}r^2|Y^l\Phi|^4\les \int_{\Sigma_0}|\hN^{\le 1} Y^l \Phi|^2\int_{\Sigma_0}|\Omega^{\le 1}Y^l\Phi|^2\les\ve^2.
\end{equation*} 
Then using Sobolev embedding on spheres and (\ref{9.22.1.22}), we derive
\begin{equation}\label{9.30.10.23}
\|Y^{1+\le 2} \Phi(0)\|_{L^p_u L^4_\omega}\le \ve^{\f12+\frac{1}{p}}, \, \|Y^{1+\le 1}\Phi\|_{L_u^p L_\omega^\infty}\les \ve^{\f12+\frac{1}{p}}, p=2,\infty
\end{equation}
as stated in (\ref{5.14.1.23}). The estimates for $\Phi$ in (\ref{5.14.1.23}) are to be proved shortly.

(\ref{9.29.4.23}) is a direct consequence of (\ref{9.22.1.22}). 
Moreover, due to (\ref{1.12.1.22}), we can easily see that $\Phi(0), c(0)=O(1)$, which will be used frequently without mentioning.

Next we prove (\ref{9.18.1.23}) and (\ref{9.30.1.23}) except the last estimate in (\ref{9.30.1.23}).
Using (\ref{7.04.9.19}) and (\ref{9.29.6.23}), we have
\begin{equation*}
u^{-1}\overline{v_\hN}=\overline{[L\Phi]}.
\end{equation*} 
This implies 
\begin{equation}\label{9.30.9.23}
\overline{v_\hN}=O(\ve^\f12),\quad \|\overline{v_{\hN}}\|_{L^2_u}\les \ve.
\end{equation}
Note $\p_u\bar f=\overline{\hN f}$ at $t=0$. Using (\ref{9.29.5.23}), $\snc_\hN \hN=-\snc\log(\bb c)=0$ due to $\bb c=1$ at $t=0$, we write
\begin{align}\label{9.30.11.23}
\hN \overline{v_{\hN}}+u^{-1}\overline{v_\hN}=-\overline{\bT \varrho}.
\end{align}
Writing $\bT\varrho=-\bN \varrho+L\varrho$, in view of the estimate of $\overline{v_\hN}$ in (\ref{9.30.9.23}) and the second set of estimates in (\ref{9.30.10.23}),  we integrate (\ref{9.30.11.23}) at $t=0$ to derive
\begin{align}\label{12.10.2.23}
\overline{v_\hN}(u)+\int_u^{u*} \overline{\bN\varrho}=O(\ve)
\end{align}
due to $v(u_*)=0$. Note at $t=0$, since $\bb=c^{-1}$, $\p_u c=\f12(\ga-1)\bN \varrho$. We thus infer 
\begin{equation*}
\frac{2}{\ga-1}(\ckc-c_*)-\overline{v_\hN}=O(\ve), \ga>1.
\end{equation*}
If $\ga=1$, then $c\equiv c_*$. $\bN=c\p_u=c_*\p_u$. Hence integrating (\ref{12.10.2.23}) implies
\begin{equation*}
\overline{\varrho-c_*^{-1}v_\hN}=O(\ve).
\end{equation*}
We have shown that at $t=0$,  
\begin{equation}\label{9.30.13.23}
\ckc-c_*+O(1) \overline{v_\hN}=O(\ve), \mbox{ if } \ga>1; \quad \bar\varrho+O(1) \overline{v_\hN}=O(\ve), c=c_* \mbox{ if } \ga=1.
\end{equation}
Using the estimate of $\overline{v_\hN}$ in (\ref{9.30.9.23}), in both cases, as long as $\ve$ is sufficiently small, we infer
\begin{equation}\label{9.30.8.23}
\ckc>\frac{3}{4}c_*.
\end{equation} 
It follows due to \Poincare  inequality and the Sobolev embedding on spheres that for $f=c, \varrho$
\begin{equation}\label{9.30.14.23}
|\Osc{f}|\les \|(r\snc)^{\le 1}\Osc(f)\|_{L^4_\omega}\les\|r\snc f\|_{L^4_\omega}\les\ve^\f12.
\end{equation}
Using (\ref{9.30.13.23})-(\ref{9.30.14.23}) and the rough estimate in (\ref{9.30.9.23}), we further obtain 
\begin{equation}\label{12.19.1.21}
c(0)>\f12 c_*, \quad |c-c_*|\les \ve^\f12,
\end{equation}
as stated in (\ref{9.30.1.23}).

 
Next, we prove
\begin{align}\label{9.30.12.23}
v_A=O(\ve), v_\bN=O(\ve^\f12), \|\sn(v_\bN), \Osc(v_\bN), v_\bN\|_{L_u^2 L_\omega^\infty}\les \ve.
\end{align}
We will adapt the proof of (\ref{6.24.1.21}) to the initial slice. Firstly, we obtain 
\begin{equation*}
|v|\les\int_u^{u_*}|\hN v|\les 1
\end{equation*}
by using (\ref{12.19.1.21}) and (\ref{9.25.1.22}). 
Next, using (\ref{12.11.1.23}) and $\bb=c^{-1}$, integrating in $u$ at $t=0$ gives
\begin{align*}
|v_A(u)|_\ga^2\les\int_u^{u_*} c^{-1}(|\sn\log c \c v_\bN|+|[\sn v]|)|v_A| du'.
\end{align*}
 Using the rough bound $|v|\les 1$, (\ref{12.19.1.21}) and (\ref{9.30.10.23}), we obtain $v_A=O(\ve)$.

  Using this estimate and the last formula in (\ref{12.11.1.23}), we thus obtain
\begin{equation*}
\p_A(v_\bN)=[\p_A v]+v_B \thetac_{AB}=[\sn_A v]+O(\ve) 
\end{equation*}
Hence using (\ref{9.30.10.23}) again, we have $\sn_A(v_\bN)=O(\ve^\f12)$. This estimate together with (\ref{9.30.9.23}) gives the second estimate and the $\|\sn(v_\bN)\|_{L_u^2 L_\omega^\infty}$ bound in (\ref{9.30.12.23}). The last two estimates in (\ref{9.30.12.23}) follow as its consequences in view of (\ref{9.30.9.23}).

In (\ref{9.30.13.23}), we have controlled $\bar\varrho$ when $\ga=1$. Next, in the case of $\ga>1$, we derive the bound on $\bar\varrho$. Similar to the derivation of (\ref{9.30.11.23}), we have
\begin{align*}
\overline{\hN(c^{-1}v_{\hN})+c^{-1}(\tr\thetac+\hN\log c) v_\hN+c^{-1}\snc\log c \c v^A}=-\overline{c^{-1}\bT\varrho}.  
\end{align*}
 Integrating the above identity by using (\ref{9.25.1.22}), (\ref{9.30.10.23}), (\ref{9.30.12.23}) at $t=0$ in $u$ gives
\begin{equation*}
\overline{\varrho-c^{-1} v_{\hN}}=O(\ve).
\end{equation*}
Hence, $\overline{\varrho}=O(\ve^\f12)$ and $\|\overline{\varrho}\|_{L^2_u}\les \ve$ due to (\ref{9.30.9.23}). Therefore as long as $\ga\ge 1$, combining this estimate with the last estimate in (\ref{9.30.10.23}) and (\ref{9.30.14.23}), we have
\begin{equation}\label{7.3.2.24}
\varrho=O(\ve^\f12), \|\varrho\|_{L_u^2 L_\omega^\infty}\les\ve. 
\end{equation}
Combining the above estimate with (\ref{9.30.12.23}) and (\ref{9.30.10.23}), the second set of estimates in (\ref{5.14.1.23}) is proved.

Recall from (\ref{12.20.4.21}) and (\ref{1.6.1.21}) on $\Sigma_0$ that
\begin{equation}\label{1.6.2.21}
\tr\chi-\frac{2}{\tir}=\frac{2}{u}(c-\ckk c)+[L\Phi].
\end{equation}
 Symbolically, 
 \begin{equation}\label{9.30.16.23} 
\hat\theta(0)=0; \bA_b(0)=u^{-1}\Osc(c)+[L\Phi]; \bA_{g,1}(0)\in\sn\Phi^\dagger.
 \end{equation} 
 
 Hence using (\ref{9.30.10.23}) and (\ref{9.30.14.23}), we obtain 
\begin{equation}\label{10.5.2.23}
|\bA, \ze|\les \ve^\f12; \quad\|\bA, \ze\|_{L^2_u L_\omega^\infty}\les \ve,
\end{equation}
with the estimates of $\bA_{g,1}, \ze$ to be further improved. In summary, we have proved (\ref{9.18.1.23}) and (\ref{9.30.1.23}) except the last estimate of the latter.

Moreover, using (\ref{9.30.8.23}), $L v_t=\tr\chi v_t$, the above estimates of $\bA_b$ and (\ref{9.22.1.22}),
\begin{equation*}
W_2[X^m\Phi]^\f12(0)\les \sum_{Y=\Omega, S}\|Y^{\le 1} X^m\Phi\|_{L^2(\Sigma_0)}\les \ve, m=1,2,3.
\end{equation*}
 (\ref{5.14.1.23}) is proved.

\subsection{Completion of the proof of Proposition \ref{12.21.1.21}} We have obtained (\ref{9.25.1.22})-(\ref{9.30.1.23}) except that the last estimate in (\ref{9.30.1.23}), which will be derived from (\ref{12.22.1.21}) as a consequence of $\bA_{g,1}=O(\ve)_{L_\omega^4}$.

Next we derive a preliminary result. (\ref{10.1.6.23}) below is (\ref{10.22.3.23}); (\ref{10.1.9.23}) is (\ref{10.1.9.23'}).
\begin{proposition}\label{12.13.1.23}
Under the assumptions (\ref{1.12.1.22})-(\ref{9.22.1.22}), with $X\in\{\Omega, S\}$,
\begin{enumerate}
\item
there hold at $t=0$ that
\begin{align}
 &|\pioh(0)|\les\ve^\f12, \|\sn_\Omega^{\le 3}\pioh(0)\|_{L^2_\Sigma}\les \ve,\label{4.12.1.23}\\
&\sn_X X(0)=O(1); \sn_X^{\le 1}(X\bN)=O(1), \sn_X^{\le 1} (\Lb \hN)=O(\ve^\f12) \label{10.5.3.23}\\
&\sn_X^{\le 1}(\sn_Y\Omega_B)=O(1), Y=\sn, L, \Lb\label{10.5.5.23}
\end{align}
and (\ref{9.8.2.22}) holds at $t=0$;

\item there hold  at $t=0$ that
\begin{align}
&\sn_X{}(\rp{a}\pih_{L A}, {}\rp{a}\pih_{\Lb A})(0)=O(\ve^\f12), \, O(\ve)_{L^2_\Sigma}; \sn_X^{\le 2} \pioh_{A L}=O(\ve)_{L^2_\Sigma}\label{10.1.3.23}
\end{align} 
and with $l\le 2$ 
\begin{equation}\label{10.1.6.23}
\begin{split}
&\|\Sc(X^{1+\le l}\Phi), \Ac(X^{1+\le l}\Phi)\|_{L_\omega^4}+|\Sc(X^{\le l}\Phi), \Ac(X^{\le l}\Phi)|\les \ve^\f12,\\
 &\|\Sc(X^{\le l}\Omega\Phi),\Ac(X^{\le l}\Omega\Phi), \tir\Sc(X^{\le l} L\Phi), \tir\Ac(X^{\le l}L \Phi)\|_{L^2_\Sigma}\les \ve, l\le 3.
 \end{split} 
\end{equation}
\begin{equation}\label{3.21.6.24}
\|\Omega^{1+\le 1}\bT\varrho, \Omega^{1+\le 2}\varrho\|_{L_\omega^4}\les \ve
\end{equation}
\item 
\begin{align}
&\|\sn_X^l (\bA_b, \bA_{g,2})(0)\|_{L^2_\Sigma}+\|\sn_X^{\le l+1}(\bA_{g,1}, [L\Phi]), (\tir\sn)^{\le l+1}(\ze, \tr\chi)\|_{L^2(\Sigma_0)}\le \ve, l\le 2.\label{9.30.6.23}\\
&|\sn_S^{\le 1}\sn\log \bb(0)|\les \ve, \label{3.22.1.24}
\end{align}
\begin{align}
&\|\sta{\Omega^m, \Lb}\varrho, \Omega^m\bN\varrho(0)\|_{L^2_\Sigma}+\|\Omega^m \fB(0)\|_{L^2_\Sigma}+|\Omega \bT\varrho(0), \Omega^{1+\le 1}\varrho(0)|\les \ve, m=1,2,3\label{10.1.9.23}
\end{align}
\end{enumerate}
\end{proposition}
\begin{remark}
We will rely on the formulas in Lemma \ref{3.22.6.21}. With the formulas therein, the estimate of 
\begin{equation}\label{10.5.4.23}
\snc_\Omega^l \Omega=O(u)\approx O(1),\, l\le 3
\end{equation}
 can be directly derived by the Euclidean geometry  at $t=0$. To include the mixed tangential derivatives of $\Omega$, we need to control ${}\rp{a}\pih$.
 
 (\ref{3.21.6.24}) and the pointwise estimates in (\ref{10.1.9.23}) are direct consequence of (\ref{9.22.1.22}) by virtue of Sobolev embedding. 
\end{remark} 

\begin{proof} 
Recall from Proposition \ref{3.22.6.21} at $t=0$  and  Proposition \ref{8.18.3.21}, the following components  of the deformation tensor of ${}\rp{a}\Omega$
\begin{equation}\label{9.22.7.22}
\begin{split}
&{}\rp{a}\pih_{LL}, {}\rp{a}\pih_{A\bN}=0,\, {}\rp{a}\pi_{L\Lb}=-\f12{}\rp{a}\pih_{\Lb\Lb}=2\Omega\log c=-\fm{{}\rp{a}\Omega}\\
&{}\rp{a}\pih_{\bT A}=v^j \tensor{\ud\ep}{^a_j^l} {e_A}_l c^{-2}-{}\rp{a}\Omega(v^l) {e_A}_l c^{-2}\\
&{}\rp{a}\pih_{AB}=-{}\rp{a}\Omega\log c \ga_{AB}.
\end{split}
\end{equation}
Hence, using (\ref{9.18.1.23}), we obtain the first estimate of (\ref{4.12.1.23}). Then similar to the proof of (\ref{1.25.2.22}), the lower order estimate in (\ref{10.5.3.23}) follows from using the first estimate of (\ref{4.12.1.23}) and (\ref{10.5.4.23}). Hence (\ref{9.8.2.22}) also holds at $t=0$.

Using (\ref{11.30.2.23}) and (\ref{10.5.2.23}), we obtain $X\bN=O(1)$ and $\Lb \bN=\ud \bA=O(\ve^\f12)$ as desired in (\ref{10.5.3.23}). 
Then using (\ref{5.23.1.23}), $l=1$ cases in the $L^4_\omega$ estimate and the second line, and the pointwise estimate in (\ref{10.1.6.23}) are obtained by using (\ref{5.14.1.23}) and the lower order estimates in (\ref{10.5.3.23}).

Next, we give more estimates by using structure equations. 
Note $\chih(0)=-\eh(0)$. With $Y=S, \tir\Lb$, recasting (\ref{s2}) and (\ref{3chi}) at $t=0$ in view of (\ref{1.30.2.22}) and (\ref{5.21.3.23}), we write 
\begin{align}\label{10.1.4.23}
\begin{split}
\tir^{-1}\sn_Y \chih&=(\tr\chi, \tr\chib, \fB)\eh+(\sn, \sn_L)\bA_{g,1}+\bA_{g,1}(\bA+\fB+\tir^{-1})\\
\sn_X\sn_L\chih&=\sn_X((\tr\chi+\fB)\chih)+\sn_X\widehat{\bR_{4A4B}}
\end{split}
\end{align}
where we used the fact $\ud\bA(0)=\bA_{g,1}(0)$.

Using (\ref{6.3.1.23}), with $l=0,1$, we write
\begin{align}\label{10.2.3.23}
\sn_X^l L(\tr\chi-\frac{2}{\tir})=\sn_X^l(\tir^{-1}\bA_b+\bA_b^2)+\sn_X^l(\wt{L\Xi_4})+\sn_X^l(|\chih|^2+\tr\chi[L\Phi]+\N(\Phi, \bp\Phi)).
\end{align}
Hence, with $l=0$ in (\ref{10.2.3.23}), using the $l=1$ estimates of (\ref{10.1.6.23}) and the pointwise estimates therein, also using (\ref{10.5.2.23}) to treat the lower order terms, we obtain
\begin{align}\label{10.7.2.23}
\sn_L \bA_b(0)=O(\ve)_{L^2_\Sigma}, O(\ve^\f12). 
\end{align}
We also have $\sn_\Omega \bA_b(0)=O(\ve)_{L^2_\Sigma}, O(\ve^\f12)$ by using the $l=1$ results in (\ref{10.1.6.23}) and its pointwise estimates, and (\ref{1.6.2.21}).

Using (\ref{10.1.4.23}), the $l=1$ results and the pointwise results in (\ref{10.1.6.23}), (\ref{10.5.2.23}) and $\fB=O(1)$, we can obtain 
\begin{equation*}
\sn_Y \chih=O(\ve)_{L^2_\Sigma}, O(\ve^\f12), \quad Y=S, \Omega, \tir\Lb.
\end{equation*}
Using (\ref{8.26.1.23}) and the pointwise results in (\ref{10.1.6.23}), we have
\begin{equation}\label{10.7.4.23}
 L[\Lb\Phi], \Lb [L\Phi]=O(\fB).
\end{equation}
Hence we have obtained at $t=0$ with $X=\Omega, S$
\begin{align}
&\fB, \tr\chi, \tr\chib, \tr k=O(1), \label{9.23.9.22}\\
&S\fB, \Lb[L\Phi]=O(1)\fB, \Omega \bT\varrho=O(\ve)\label{12.11.3.23}\\
&\sn_X\tr\chi=\vs(X)O(1)+O(\ve^\f12), \sn_S \tr\chib=O(1)+O(\ve^\f12),\label{4.13.2.23}\\
&\sn_X(\bA_b, \chih), \sn_\Lb \chih=O(\ve)_{L^2_\Sigma}, O(\ve^\f12)\label{10.2.1.23}\\
&\mho, \bA, \sn^{\le 1}\ze=O(\ve^\f12), O(\ve)_{L^2_u L_\omega^\infty}\label{4.13.3.23}
\end{align}
where the estimate of $\mho(0)$ follows by repeating the proof of (\ref{mhoini}). Here the estimate of $\sn_S \tr\chib$ follows by using the result of $\sn_S\tr\chi$, (\ref{12.11.3.23}) with the help of $\tr\chi+\tr\chib=\fB$.

Recall from (\ref{5.13.10.21}) that
\begin{equation}\label{6.22.5.24}
[\Lb, \Omega]\varrho=\ud \bA(\Omega) \bN\varrho+\pioh_{A\Lb}\sn\varrho.
\end{equation}
Using (\ref{4.13.3.23}), $\Omega \Lb \varrho=O(\ve^\f12), O(\ve)_{L_u^2 L_\omega^\infty}$ (due to the first pointwise estimates in (\ref{10.1.9.23}) and the pointwise estimates in (\ref{10.1.6.23})),  and the first estimate in (\ref{4.12.1.23}), we have
\begin{equation}\label{6.22.2.24}
\Lb\Omega \varrho=O(\ve^\f12), O(\ve)_{L_u^2 L_\omega^\infty}
\end{equation}

At $t=0$, using (\ref{5.23.1.23}), the first estimate in (\ref{10.5.3.23}), (\ref{4.13.2.23})-(\ref{4.13.3.23}) and the pointwise estimate in (\ref{10.1.6.23}), we have
\begin{equation*}
\sn_X(X\bN)=O(1)
 \end{equation*}
 as stated in the second estimate in (\ref{10.5.3.23}).
Thus the $l=2$ cases  for both the $L_\omega^4$ estimate and  the second line of (\ref{10.1.6.23})  can be proved  in the same way as in Lemma \ref{6.30.4.23}.

Noting that $\Omega^{\le 1}\la(0)=0$ and $L \la=\bA_{g,1}\c\Omega$, using (\ref{7.16.2.22}), (\ref{3.28.3.24}), (\ref{9.18.1.23}) and the estimates of $\bA_{g,1}$ in (\ref{10.5.2.23}), we derive
\begin{equation}\label{3.21.3.24}
\sn_X^n{}\rp{a}\pih_{L A}=\sn_X^n(\eta(\Omega)+c^{-1}\sn\la)+O(1)(X^{\le n} v)+O(\ve), n=1.
\end{equation}
By (\ref{cmu_2}) and $\sn\la(0)=0$, $\sn_S\sn\la(0)=\sn S\la(0)$. Using the pointwise estimate in (\ref{10.1.6.23}), we have $\sn_X {}\rp{a}\pih_{L A}=O(\ve^\f12)$. 

 Similarly, using the $L^2_\Sigma$ estimates with $l=1$ in (\ref{10.1.6.23}), we derive 
\begin{equation*}
\|\sn_X {}\rp{a}\pih_{L A}(0)\|_{L^2_\Sigma}\les \|X\Omega^{\le 1} \Phi, X^{\le 1} v, \sn_X\sn\la\|_{L^2_\Sigma}+\ve\les \ve.
\end{equation*} 
as stated in (\ref{10.1.3.23}).  The estimates of $\sn_X{}\rp{a}\pih_{\Lb A}$ can be obtained similarly.

Using the pointwise estimates of $\sn_X^{\le 1}(\pioh_{AL},\pioh_{A\Lb})$, the pointwise estimates in (\ref{10.1.6.23}) and (\ref{10.1.9.23}), we can obtain (\ref{10.5.5.23}) similar to  (\ref{4.22.4.22})  at $t=0$.

Using the proved estimates in (\ref{10.5.3.23}) and (\ref{10.1.6.23}), we have
\begin{equation}
\|\sn_X^{\le 2} (\zb, \eh, \sn\varrho, [L\Phi]), (\tir\sn)^{\le 2}\ze\|_{L^2_\Sigma}\les \ve. \label{10.4.7.23}
\end{equation}
 Thus, due to  (\ref{9.30.16.23}), (\ref{1.6.2.21}) and (\ref{10.4.7.23}), we obtain
\begin{equation}\label{10.4.3.23}
\|\sn_\Omega^{\le 2} \bA(0)\|_{L^2_\Sigma}\les \ve,
\end{equation}
which is as desired in (\ref{9.30.6.23}). 

Next we prove 
 \begin{equation}\label{6.22.1.24}
 \vs^+(X^2)\|\sn_X^2 \bA(0)\|_{L^2_\Sigma}\les \ve
 \end{equation}
  in (\ref{9.30.6.23}).
  In particular we only need to prove it for $\bA=\bA_b, \bA_{g,2}$. Firstly, due to (\ref{3.21.6.24}), (\ref{10.2.1.23}) and the $L^4_\omega$ estimates in (\ref{10.1.6.23}), we have
 \begin{equation}\label{12.12.4.23}
  \|\Omega (\fB,\bA, [L\Phi])\|_{L^4_\omega}\les \ve^\f12.
  \end{equation}
Note, at $t=0$, for an $S_{t,u}$ tangent tensor or scalar $F$, we can bound 
\begin{equation*}
\|\sn_X(\fB F)\|_{L^2_\Sigma}\les \|\sn_X F\|_{L^2_\Sigma}+\|X\fB\|_{L_\omega^4}\|r F\|_{L_u^2 L_\omega^4}.
\end{equation*}
By using (\ref{10.4.7.23})-(\ref{12.12.4.23}) and Sobolev embedding on spheres,  we have 
$$\|\sn_\Omega\Big(\fB \c(\bA,[L\Phi])\Big)\|_{L^2_\Sigma}\les \ve,$$
Similarly, if $X=S$, it follows by using (\ref{12.11.3.23}), (\ref{10.2.1.23}) and (\ref{4.13.3.23}) that
\begin{equation*}
\|\sn_S\Big(\fB (\bA, [L\Phi])\Big)\|_{L^2_\Sigma}\les  \ve.
\end{equation*}
In both cases, we have
\begin{equation}\label{12.11.4.23}
\|\sn_X\Big(\fB (\bA, [L\Phi])\Big)\|_{L^2_\Sigma}\les  \ve.
\end{equation} 
Using $\snc\log \bb=-\snc\log c$ at $t=0$, we derive
\begin{equation*}
\snc_\hN \snc \log \bb=-\snc_\hN \snc \log c.
\end{equation*}
Since we can obtain from (\ref{10.1.6.23}) by using Sobolev embedding that 
 $L\Omega \varrho=O(\ve)_{L_u^2 L_\omega^\infty}$, 
   due to (\ref{6.22.2.24}), (\ref{4.13.3.23}) and $\sn_\hN \Omega=O(1)$, we infer 
\begin{align*}
\|\snc_\hN \snc \varrho\|_{L_u^2 L_\omega^\infty}\les \ve.
\end{align*}
Thus  $\sn\log \bb=O(\ve)$ by using the above estimate by virtue of (\ref{3.28.2.21}) and  $\b(0)=0$.
 
Using this estimate, in view of (\ref{1.31.2.24}), using $\bb^{-1}(0)=c(0)$, 
and the pointwise estimates of (\ref{10.1.9.23}), we can obtain the bound of $\sn_S \sn\log\bb$ in (\ref{3.22.1.24}). Hence (\ref{3.22.1.24}) is proved.

 By using the proved parts in (\ref{10.1.6.23}), (\ref{12.11.4.23}), (\ref{9.23.9.22})-(\ref{4.13.3.23}) and (\ref{10.4.7.23}), we claim
\begin{align}
&\|\sn_X^{\le 1} \widehat{\bR_{4A4B}}, \sn_X^{\le 1}\widehat{\bR_{A34B}}, \sn_X^{\le 1}\N(\Phi, \bp\Phi),  \sn_X^{\le 1}\wt{L \Xi_4}\|_{L^2_\Sigma}\les \ve. \label{10.2.2.23}
\end{align}
Indeed, in addition to the aforementioned estimates, the first two estimates follow by using (\ref{1.30.2.22}), (\ref{5.21.3.23}) and (\ref{10.7.4.23}). The estimate of the null form can be obtained by repeating the analysis in Proposition \ref{6.24.10.23}. Using the null form estimate and (\ref{3.20.1.22}), we can obtain the estimate for $\sn_X^{\le 1}\wt{L \Xi_4}$.  

By using second line in (\ref{10.1.4.23}), (\ref{10.2.3.23}), (\ref{9.23.9.22})-(\ref{4.13.3.23}), (\ref{10.2.2.23}) and (\ref{12.11.4.23}), we derive
\begin{equation}\label{10.6.2.23}
\|\sn_X L\bA_b, \sn_X \sn_L \chih, \sn_L \sn_X \bA_b, \sn_L \sn_X \chih\|_{L^2_\Sigma}\les \ve
\end{equation}
where the last two estimates follow from using the first two, the estimates of (\ref{7.17.6.21}) and (\ref{3.21.1.23}) at $t=0$. Thus we conclude (\ref{6.22.1.24}). 

In view of (\ref{5.23.1.23}), using the proved estimates in (\ref{9.30.6.23}) we can obtain 
\begin{equation*}
\sn_X^{\le 2}(X\bN)=O(1)+O(\ve)_{L^2_\Sigma}.
\end{equation*}
Using the above estimate and the first two estimates in (\ref{10.5.3.23}), we can get the $l=3$ case in the second line of  (\ref{10.1.6.23}). Thus (\ref{10.1.6.23}) is proved. 

Moreover, similar to (\ref{3.21.3.24}), we can obtain by using (\ref{10.1.6.23}) and (\ref{10.5.5.23})
\begin{equation*}
\sn_X^{\le 2} \pioh_{A L}=O(\ve)_{L^2_\Sigma}.
\end{equation*}

Therefore (\ref{10.1.3.23}) is proved. 

Using (\ref{10.1.3.23}), similar to (\ref{9.1.1.21}), we can obtain
\begin{equation*}
\sn_X^2 \sn_Z {\Omega^A}=O(1)+O(\ve)_{L^2_\Sigma}, Z= L, e_1, e_2.
\end{equation*}

It follows due to (\ref{10.1.6.23}) and  that
\begin{align*}
&\|\sn_X^{\le 3} (\zb, \eh, \sn\varrho, [L\Phi])\|_{L^2_\Sigma}\les \ve
\end{align*}
  and
  \begin{align*}
  \|(\tir\sn)^{\le 3}(\ze, \tr\chi)\|_{L^2_\Sigma}\les \ve, 
  \end{align*}
  where the second set of estimates follows also by using $\bb^{-1}=c$ and (\ref{9.30.16.23}) at $t=0$. Thus (\ref{9.30.6.23}) is proved. 
   
   Using (\ref{9.22.7.22}) and  (\ref{9.22.1.22}), we have  $\sn_\Omega^{\le 3} \pioh=O(\ve)_{L^2(\Sigma_0)}$ as stated in (\ref{4.12.1.23}). Using (\ref{4.12.1.23}), (\ref{10.1.3.23}), (\ref{9.30.6.23}), (\ref{4.13.3.23}), (\ref{6.22.2.24}) and  the pointwise estimates in (\ref{10.1.9.23}), we can obtain the first two estimates in (\ref{10.1.9.23}) inductively. Hence (\ref{10.1.9.23}) is proved.
   \end{proof}


Next we prove (\ref{1.12.4.22}) in the following result. 
\begin{lemma}
Let $X\in\{\Omega, S\}$. At $t=0$, under the assumptions (\ref{1.12.1.22})-(\ref{9.22.1.22}), we have
\begin{equation}\label{9.22.6.22}
\|X^{\le 1}\Box_\bg \Phi, X^{\le 1}\N(\Phi, \bp\Phi), \Box_\bg \Omega \Phi\|_{L^2_\Sigma}\les \ve.
\end{equation}
\end{lemma}  
\begin{proof}
At $t=0$, from (\ref{5.18.3.21}) and $\pioh_{\bN A}(0)=0$ (due to $\la(0)=0$), applying (\ref{5.02.3.21_1}) to $f=\Phi$ gives
\begin{align}\label{9.23.2.22}
\begin{split}
[\Box_\bg {},\rp{a}\Omega]\Phi&=\f12 \fm{{}\rp{a}\Omega}\Box_\bg \Phi+{}\rp{a}\pih^{AB}\sn_A \sn_B \Phi+\bJ[{}\rp{a}\Omega]^\mu\p_\mu \Phi\\
  &+\frac{1}{4}{}\rp{a}\pih_{\Lb\Lb}\bd^2_{LL}\Phi-{}\rp{a}\pih_{\bT A}\bd^2_{\bT A}\Phi.
  \end{split}
\end{align}
Using Proposition \ref{11.12.2.22}, Lemma \ref{3.23.2.23}, and (\ref{5.02.2.21}), we derive
\begin{align*}
\bJ[{}\rp{a}\Omega]_A&=-\f12\sn_A\fm{{}\rp{a}\Omega}-\sn_\bT{}\rp{a}\pih_{\bT A}+(\tr k+\eh+k_{\bN\bN}){}\rp{a}\pih_{\bT A}-\sn_B\log c{}\rp{a}\pih_{BA}\\
\bJ[{}\rp{a}\Omega]_L&=\sn_A{}\rp{a}\pih_{A\bT}+\zeta\c {}\rp{a}\pih_{A\bT}+L{}\rp{a}\Omega\varrho+c^{-1}\bA_{g,1}([L\Phi]+\tr\chi)\\
\bJ[{}\rp{a}\Omega]_\Lb&=\sn_A{}\rp{a}\pih_{A\bT}-\f12(L+\tr\chi-k_{\bN\bN}){}\rp{a}\pih_{\Lb \Lb}+\zb\pih_{A\bT}+\Lb{}\rp{a}\Omega\varrho\\
&+c^{-1}\bA_{g,1}([L\Phi]+\tr\chi).
\end{align*}
Moreover, it follows by using (\ref{6.28.6.21}), (\ref{12.16.1.23}), (\ref{12.11.3.23}),  (\ref{10.2.1.23}), (\ref{4.13.3.23}) and the pointwise estimate in (\ref{10.1.6.23}) that
\begin{equation*}
\sn_\bT {}\rp{a}\pih_{\bT A}=\sn_L(\la\ud \bA)+\sn_\bT(\eta(\Omega)+v^*)=O(1)
\end{equation*}
where we used the fact that $\bN(\la \ud \bA)$ is trivial terms on the initial slice, $\sn_\bT \Omega=O(1)$, and $|L\la\c \ud\bA|\les |\bA_{g,1}|^2\les \ve$ at the initial slice. 
We then conclude by using Proposition \ref{12.13.1.23} that
\begin{equation}\label{9.23.4.22}
\|\bJ[{}\rp{a}\Omega]_L, \bJ[{}\rp{a}\Omega]_\Lb\|_{L^2_\Sigma}\les \ve, \quad \bJ[{}\rp{a}\Omega]_A=O(1).
\end{equation}
By using Lemma \ref{comp}, (\ref{9.25.1.22})  and (\ref{9.22.1.22}), we have
\begin{equation}\label{4.17.2.23}
\|\sn^{1+\le 1}\Phi, \bd^2_{LL}\Phi,\bd^2_{L A}\Phi, \sn \Phi, L\Phi\|_{L^2_\Sigma}\les \ve,\quad |\bd^2_{\bN A}\Phi(0)|\les 1.
\end{equation}
Hence we conclude
\begin{equation*}
\|[\Box, \Omega]\Phi\|_{L^2_\Sigma}\les \ve.
\end{equation*}
Similar to (\ref{8.23.1.23}), $X\Lb \Phi=O(1)$ for $X=\Omega, S$. Then using Proposition \ref{geonul_5.23_23}, (\ref{12.11.4.23}), (\ref{9.25.1.22}) and (\ref{10.1.6.23}), we have
\begin{equation*}
\|X\Box_\bg\Phi, X^{\le 1}\N(\Phi, \bp\Phi), \Box_\bg\Phi\|_{L^2_\Sigma}\les \ve.
\end{equation*}
Similar to (\ref{8.23.1.23}), 
we then obtain the last estimates in (\ref{9.22.6.22}) by combining the two sets of estimates in the above. Thus the proof of the lemma is complete. 
\end{proof}
We prove (\ref{12.22.1.21}), (\ref{1.12.3.22}) and (\ref{3.21.6.24'}) by establishing the following result. 
\begin{lemma}
 Under the assumptions (\ref{1.12.1.22})-(\ref{9.22.1.22}), with $X\in\{\Omega, S\}$, there hold at $t=0$ that 
\begin{align}
&\sum_{Y=\tir\sn, \tir\Lb, L}\|\sn_X^l\sn_Y\bA_{g,1}\|_{L_\omega^4}+\|(\sn_X^l\sn, \sn\sn_X^l)([\Lb\Phi], [L\Phi]) \|_{L_\omega^4}\les \ve,\, l=0,1\label{10.4.1.23}\\
&\sum_{Y=\sn, \Lb, L}|\sn_Y^{\le 1} \bA_{g,1}|+|\sn^{\le 1}\bA_{g,2}, \sn \bA_b|\les \ve \label{6.22.4.24}\\ 
&\|\sn_X^{\le 1} \ze\|_{L_\omega^4}\les \ve, \|\sn_\Omega^{1+\le 1}\bA_b\|_{L_\omega^4}+\|\sn_X^{\le 2}\ze\|_{L^2_\Sigma}\les \ve\label{3.22.3.24}\\
&\| \sn_X^l(\tir\sn)([\Lb\Phi], [L\Phi]), \sn_X^l\sn_\Lb \bA_{g,1}, \sn_\Lb \sn_X^l \bA_{g,1}\|_{L^2_\Sigma}\les \ve,\, l\le 2\label{10.4.10.23}\\
&\|(\tir\sn)^{\le 1}\sF\|_{L_\omega^4}+\| (\tir\sn)^{\le 2}\sF\|_{L^2_\Sigma}\les \ve.\label{6.23.13.24}
\end{align}
\end{lemma}
\begin{remark}
To prove this lemma, we will use Proposition \ref{8.12.1.23} and Proposition \ref{10.4.4.23} which were proved by using elliptic estimates. They hold the same at $t=0$ without relying on Assumptions \ref{5.13.11.21+}.
\end{remark}
\begin{remark}
(\ref{1.12.3.22}) follows from  (\ref{10.4.10.23}), the last estimates in (\ref{6.23.13.24}) and (\ref{3.22.3.24}), and Proposition \ref{12.13.1.23} (3). (\ref{12.22.1.21}) is given by (\ref{3.22.1.24}), (\ref{10.4.1.23})-(\ref{3.22.3.24}) and the $L^4_\omega$ estimates in (\ref{6.23.13.24}). (\ref{3.21.6.24'}) can be obtained by using (\ref{3.21.6.24}) and (\ref{10.4.1.23}).
\end{remark}
\begin{proof}
 (\ref{6.22.4.24}) follows as the consequence of (\ref{10.4.1.23}) by using Sobolev embedding on spheres and (\ref{9.30.16.23}). We will prove (\ref{10.4.1.23}) together with (\ref{10.4.10.23}),  then prove (\ref{3.22.3.24}) and (\ref{6.23.13.24}).

{\bf $\bullet$ Step 1}: We consider the cases $l=0$ in (\ref{10.4.1.23}) and $l\le 1$ in (\ref{10.4.10.23}).

In view of (\ref{8.26.1.23}), by using (\ref{9.30.6.23}), (\ref{4.13.3.23}), (\ref{12.11.4.23}) and (\ref{9.22.6.22})
 we derive
\begin{equation}\label{12.14.1.23}
 \|\tir\sn \Lb[L\Phi](0)\|_{L^2_\Sigma}\les \ve.
\end{equation}
Note that, at $t=0$, for scalar functions $f$ there holds 
\begin{align*}
&\|[\Omega, \Lb]f, [\tir\sn, \Lb]f\|_{L^2_\Sigma}\les \|\sn f\|_{L^2_\Sigma}\ve^\f12+\ve\|\bN f\|_{L_u^\infty L_\omega^2}.
\end{align*}
Indeed, the above commutator estimates can  be obtained by using (\ref{5.13.10.21}), (\ref{7.03.1.19}), (\ref{4.12.1.23}) and (\ref{4.13.3.23}). 
Using (\ref{12.14.1.23}), the above estimates, (\ref{10.1.6.23}) and (\ref{12.11.3.23}), we obtain
\begin{align*}
\|\sn_\Lb(\tir\sn)[L\Phi](0)\|_{L^2_\Sigma}&\les \|\tir\sn \Lb[L\Phi](0)\|_{L^2_\Sigma}+\ve^\f12\|\sn[L\Phi]\|_{L^2_\Sigma}+\ve\|\bN[L\Phi]\|_{L_u^\infty L_\omega^2}\\
&\les \ve.
\end{align*}
Using (\ref{10.1.6.23}), we obtain 
\begin{equation}\label{3.21.5.24}
\Omega [L\Phi], \Omega \fB=O(\ve)_{L_\omega^4},
\end{equation}
where in view of $[\Lb\Phi]=\Lb\varrho+[L\Phi]$, the second estimate follows by using the first one and the first estimate in (\ref{3.21.6.24}). This gives the second set of estimates in (\ref{10.4.1.23}) for $l=0$.

Using (\ref{10.5.5.23}), (\ref{9.30.6.23}) and (\ref{10.1.9.23}), we have
\begin{equation*}
\|\sn_\Omega^{\le 1}(\tir \sn)\fB\|_{L^2_\Sigma}\les \ve. 
\end{equation*}
By using the above estimate, (\ref{10.1.8.23}), (\ref{9.30.6.23}), (\ref{4.13.3.23}) and (\ref{12.11.4.23}), we derive
\begin{align}\label{10.4.8.23}
\|\sn_\Omega^{\le 1}\sn_\bN \ep, \sn_\bN \sn_\Omega\ep\|_{L^2_\Sigma}\les \ve,
\end{align}
where the second estimate is obtained by using the first one, (\ref{9.30.6.23}) and (\ref{7.4.1.21}) with $\Delta_0$ therein replaced by $\ve^\f12$. 

Using Lemma \ref{2.9.3.23}, we obtain as a consequence of the first estimate in (\ref{4.12.1.23}), (\ref{10.1.9.23}) with $m=1$ and (\ref{10.1.6.23}) that
\begin{equation*}
\|[\Lb\Omega\Phi]\|_{L^2_\Sigma}\les \ve.
\end{equation*}
Using (\ref{10.5.5.23}), similar to (\ref{7.11.7.21}), at $t=0$ we infer
\begin{align*}
\|[\sn_S\sn_\Lb \sn\Phi]\|_{L^2_\Sigma}&\les \|[\tir\sD \Omega \Phi]+\tir[\Box_\bg \Omega \Phi]+[\Lb \Omega \Phi]+[L \Omega \Phi]\|_{L^2_\Sigma}+\|[\sn\Omega \Phi]\|_{L^2_\Sigma}\les \ve
\end{align*}
where we used (\ref{9.22.6.22}) and (\ref{10.1.6.23}) to derive the last estimate. 

 We then obtain by using (\ref{10.5.3.23}), similar to Lemma \ref{6.30.4.23}, that
\begin{equation*}
\|\sn_S^{\le 1} \sn_\Lb \zb\|_{L^2_\Sigma}\les\ve.
\end{equation*}
Combining the above estimate with (\ref{10.4.8.23}) and (\ref{9.30.6.23}) implies
\begin{equation*}
\|\sn_X\sn_\Lb \zb\|_{L^2_\Sigma}\les\ve.
\end{equation*}
It remains to consider the case that $\bA_{g,1}=\eh$ in (\ref{10.4.1.23}). Using (\ref{9.30.6.23}), (\ref{3.22.1.24}), (\ref{12.11.4.23}) and (\ref{1.30.2.24}) (with $\Delta_0$ replaced by $\ve$), we have
\begin{align*}
\|\sn_X^{\le 1} \sn_\Lb\eh\|_{L^2_\Sigma}\les\ve.
\end{align*}
Hence the case $\bA_{g,1}=\eh$ in (\ref{10.4.10.23}) is proved. Thus the estimate of $\sn_X\sn_\Lb \bA_{g,1}$ in (\ref{10.4.10.23}) is proved. 
Using this estimate, (\ref{7.4.1.21}) and (\ref{9.30.6.23}), we then obtain 
\begin{equation*}
\|\sn_\Lb \sn_X^{\le 1}\bA_{g,1}\|_{L^2_\Sigma}\les \ve
\end{equation*}
as stated in  (\ref{10.4.10.23}).  Thus the $l\le 1$ case in (\ref{10.4.10.23}) is proved.

Using the above estimate and (\ref{9.30.6.23}) we obtain
\begin{align}\label{3.22.2.24}
\|\sn_X \bA_{g,1}\|_{L^4_\omega}\les \ve, X=\Omega, S; |\bA_{g,1}|\les \ve
\end{align}
where the last estimate follows as a consequence of the first one by using Sobolev embedding. 

Next we prove (\ref{3.22.3.24}). For the $L^2_\Sigma$ estimate of $\ze$,  noting that the lower order estimates have been proved in (\ref{9.30.6.23}) and (\ref{3.22.1.24}), and the case of $\vs^+(X^2)=0$ is proved in (\ref{9.30.6.23}),  it suffices to consider only the second order estimates with $\vs^+(X^2)=1$. 
  
Using (\ref{6.23.6.24}), (\ref{4.13.2.23})-(\ref{4.13.3.23}), (\ref{3.22.1.24}), (\ref{10.4.10.23}) with $l\le 1$, we bound
\begin{align*}
\|\sn_X\sn_L\sn\log\bb\|_{L^2_\Sigma}\les\|\sn_X \sn k_{\bN\bN}\|_{L^2_\Sigma}+\ve\les \ve. 
\end{align*}
 Note that in view of (\ref{1.21.2.22}), (\ref{3.22.2.24}), the $l=0$ case in (\ref{10.4.1.23}) and the $l\le 1$ case in (\ref{10.4.10.23}), 
we have
\begin{equation}\label{6.23.1.24}
\bR_{AC4B}, \sn_X^{\le 1}(\chi \zb)=O(\ve)_{L_\omega^4}, \sn_X \bR_{AC4B}=O(\ve)_{L^2_\Sigma}
\end{equation}
where we also used (\ref{4.13.2.23}) and (\ref{10.2.1.23}) to bound $\sn_X^{\le 1}\chi$. 

 Similar to (\ref{3.21.1.23}) for $F$, using (\ref{6.23.1.24}), (\ref{4.12.1.23}) and (\ref{4.13.3.23}), we derive
 \begin{align}\label{6.23.7.24}
 [\tir\sn, L]F, [\Omega, L]F=O(\ve)_{L^4_\omega}F+O(\ve^\f12)\sn F.
 \end{align}
Applying the above estimate to $F=\sn\log \bb$, we derive
\begin{equation*}
[\Omega, L]\sn\log \bb=O(\ve)_{L_\omega^4}
\end{equation*}
where we used (\ref{4.13.3.23}).
Hence 
$$\|\sn_S\sn_\Omega \sn\log \bb\|_{L^2_\Sigma}\les \ve.$$ Therefore the last estimate in (\ref{3.22.3.24}) is proved with the help of (\ref{9.30.6.23}). 

Moreover, at $t=0$, the estimate of
$
\sn_X\ze=O(\ve)_{L_\omega^4}
$
 follows by using (\ref{3.22.2.24}) and (\ref{3.22.1.24}). Therefore, the $\ze$ estimates in (\ref{3.22.3.24}) are proved. 
 
Using the above estimate, (\ref{3.21.5.24}) and $\pioh=O(\ve^\f12)$ from (\ref{4.12.1.23}), 
recalling from (\ref{10.1.9.23}) and (\ref{6.22.5.24}), we have
\begin{equation*}
\Omega\Lb\varrho, \Lb\Omega\varrho=O(\ve)_{L_\omega^4}.
\end{equation*}
We then deduce from the above estimate, (\ref{10.5.5.23}) and $\fB=\Lb\varrho+[L\Phi]$ that
\begin{equation*}
\|\sn\fB, \sn_\Lb \sn\varrho\|_{L^4_\omega}\les \ve.
\end{equation*}

In view of (\ref{10.1.8.23}), the above estimate, (\ref{3.22.2.24}), we then obtain 
\begin{align*}
\|\sn_\bN \ep, \sn_\Lb \ep\|_{L^4_\omega}\les \ve. 
\end{align*}
Finally using (\ref{12.16.1.23}) and  (\ref{3.22.2.24}) we derive
\begin{align*}
\|\sn_\Lb \eh\|_{L_\omega^4}\les \ve.
\end{align*}
Therefore we have completed the proof of the $l=0$ case in (\ref{10.4.1.23}). In view of (\ref{9.30.16.23}), we also obtained $\sn \bA_b=O(\ve)_{L_\omega^4}$ by using the $l=0$ case in (\ref{10.4.1.23}). We will rely on the estimate of $\sn\sn_X [L\Phi]$ in (\ref{10.4.1.23}) to bound $\|\sn_\Omega^2 \bA_b\|_{L_\omega^4}$ which is the only remaining estimate in (\ref{3.22.3.24}). 

{\bf $\bullet$ Step 2:} $l=1$ in (\ref{10.4.1.23}) and $l=2$ in (\ref{10.4.10.23}).
 We consider the first estimate in (\ref{10.4.10.23}). 

{\bf Case (a) $\vs^+(X^2)=0$.}  Using the fact that $[\Lb\Phi]=\Lb \varrho+[L\Phi]$ and Proposition \ref{12.13.1.23}, we obtain
\begin{equation*}
\|\sn_\Omega^2\sn[\Lb\Phi]\|_{L^2_\Sigma}\les  \ve. 
\end{equation*}
{\bf Case (b) $\vs^+(X^2)=1$.} In view of (\ref{10.1.6.23}) and (\ref{10.5.5.23}), $\sn_X^2 (\tir \sn)[L\Phi]=\O(\ve)_{L^2_\Sigma}$. It only remains to consider $\sn_X^2 (\tir \sn)\Lb \varrho$. If $X_1=S$, we write
\begin{align*}
\sn_{X_2}\sn_S(\tir \sn)\Lb \varrho&=\sn_{X_2}[\sn_S, \tir\sn]\Lb \varrho+\sn_{X_2}(\tir \sn)S\Lb \varrho.
\end{align*} 
Note that (\ref{10.2.1.23}), (\ref{4.13.3.23}) and (\ref{10.1.6.23}) imply 
 $$\sn_X^{\le 1}\bA=O(\ve^\f12),\, X=\Omega, S$$
  In view of (\ref{cmu_2}) and the above estimate, we derive 
 yields
\begin{align*}
\|\sn_{X_2}[\sn_S, \tir\sn]\Lb \varrho\|_{L^2_\Sigma}&\les \sum_{X=\Omega, S}\|\sn_X^{\le 1} \bA\|_{L^\infty_x}\|\sn_{X_2}^{\le 1}\sn\Lb \varrho\|_{L^2_\Sigma}\les \ve^\f12\sum_{X=\Omega, X_2}\|\sn_{X_2}^{\le 1}\sn \Lb \varrho\|_{L^2_\Sigma}\les \ve^\frac{3}{2}
\end{align*}
where we used  the $l\le 1$ case in (\ref{10.4.10.23}) to obtain the last estimate. 

Next we consider $\|\sn_{X_2}(\tir \sn)S\Lb \varrho\|_{L^2_\Sigma}.$ Applying (\ref{3.22.5.24}) to $f=\varrho, a=0$, 
\begin{align*}
X_2\Omega S\Lb \varrho&=X_2\Omega(\tir (h-k_{\bN\bN})\Lb \varrho)+X\Omega(\tir \sD \varrho+\tir \Box_\bg \varrho)+X\Omega(h S\varrho+\tir\sn \varrho\c \zb).
\end{align*}
Using $l\le 1$ in (\ref{10.4.10.23}), (\ref{10.1.6.23}) and $l=0$ in (\ref{10.4.1.23}), we have
\begin{equation}\label{6.23.9.24}
X^{\le 1}\Omega(\Box_\bg \varrho)=O(\ve)_{L^2_\Sigma} 
\end{equation}
It follows by using the above estimate, (\ref{12.11.3.23}),  (\ref{10.1.6.23}), (\ref{9.30.6.23}), (\ref{9.22.6.22}) and the $l\le 1$ case in (\ref{10.4.10.23}) that 
\begin{equation*}
X_2\Omega S\Lb \varrho=O(\ve)_{L^2_\Sigma}. 
\end{equation*}
Hence
\begin{equation}\label{6.23.5.24}
\sn_{X_2}\sn_S(\tir \sn)\Lb \varrho=O(\ve)_{L^2_\Sigma}.
\end{equation}
If $X_2=S, X_1=\Omega$, we write
\begin{align*}
\sn_S\sn_\Omega (\tir \sn)\Lb \varrho=[S, \Omega](\tir \sn)\Lb \varrho+\sn_\Omega\sn_S(\tir \sn)\Lb \varrho. 
\end{align*}
The second term on the right-hand side has been treated in (\ref{6.23.5.24}). Applying (\ref{6.23.7.24}) to $F=\tir \sn \Lb \varrho$,  we obtain
\begin{align*}
[\sn_S, \sn_\Omega](\tir \sn\Lb \varrho)=O(\ve^\f12)\sn_\Omega\sn \Lb \varrho+O(\ve)_{L_\omega^4} \tir \sn \Lb \varrho=O(\ve^\frac{3}{2})_{L^2_\Sigma}
\end{align*}
where to obtain the last estimate, we used the $l\le 1$ result in (\ref{10.4.10.23}) and (\ref{3.21.5.24}).

Summarizing the above estimates, we conclude 
\begin{equation*}
\vs^+(X^2)\sn_X^2(\tir\sn)[\Lb\Phi]=O(\ve)_{L^2_\Sigma},
\end{equation*}
which closes {\bf Case (b)}. Thus combining the estimate of $\sn_X^{1+\le 2}([L\Phi])$ in (\ref{9.30.6.23}) with the results of the two cases, the first estimate in (\ref{10.4.10.23}) is proved. 

Next we consider the  estimate for $\sn_X^2\sn_\Lb \bA_{g,1}$  in (\ref{10.4.10.23}).

Note that it follows by using (\ref{3.22.2.24}), (\ref{5.13.10.21}) and $\la(0)=0$ that
\begin{equation}\label{9.21.1.24}
 \sn_X^{1+\le 1} \la=O(\ve)_{L_\omega^4},\, X\la=O(\ve).
 \end{equation}
 Using the above estimate,  (\ref{3.22.1.24}), (\ref{4.13.3.23}) and the last estimate of (\ref{3.22.3.24}), we have $X^2(\ud \bA\c \la)(0)=O(\ve)_{L^2_\Sigma}$. Using this estimate, the last estimate in (\ref{10.1.3.23}) and (\ref{10.1.6.23}), we have
\begin{equation*}
\sn_X^2\pioh_{A\Lb}(0)=O(\ve)_{L^2_\Sigma}.
\end{equation*}
Moreover similar to (\ref{8.23.1.23}), using (\ref{10.4.10.23}) with $l\le 1$, we have
\begin{equation}\label{3.22.6.24}
X^2\fB(0)=O(1)+O(\ve)_{L^2_\Sigma}.
\end{equation}
 Using the above estimates, and applying (\ref{1.13.1.23}) with $l=2$ to bound
\begin{equation}\label{3.23.1.24}
\sn_X^2\sn_\Lb \Omega(0)=O(1)+O(\ve)_{L^2_\Sigma}.
\end{equation}
Due to (\ref{3.23.1.24}), (\ref{10.5.5.23}), $l=1$ case in (\ref{10.4.10.23}) and (\ref{3.22.2.24}), we derive
\begin{equation*}
\|\sn_X^2 \sn_\Lb \sn\varrho\|_{L^2_\Sigma}\les \|\sn_X^2 \Lb \Omega\varrho\|_{L^2_\Sigma}+\ve.
\end{equation*}
If $X_1=S$, applying (\ref{3.22.5.24}) to $f=\varrho, a=1$,
\begin{align*}
X S\Lb\Omega \varrho&=X(\tir (h-k_{\bN\bN})\Lb \Omega \varrho)+X(\tir \sD \Omega \varrho+\tir \Box_\bg \Omega \varrho)+X(h S\Omega \varrho+\tir\sn \Omega \varrho\c \zb).
\end{align*}
Using (\ref{9.23.9.22})-(\ref{4.13.3.23}), the $l=0$ case in (\ref{10.4.1.23}), the $l=1$ case in (\ref{10.4.10.23}) and (\ref{10.1.6.23}), we have
\begin{equation}\label{3.22.7.24}
X S\Lb\Omega \varrho=X(\tir \Box_\bg \Omega \varrho)+O(\ve)_{L^2_\Sigma}.
\end{equation}
 
Using (\ref{xdjo}), the $l=0$ case in (\ref{10.4.1.23}), the $l=1$ case in (\ref{10.4.10.23}), the $\ze$ estimates in (\ref{3.22.3.24}), (\ref{9.21.1.24}) and Proposition \ref{12.13.1.23} we derive
\begin{align*}
&\sn_X\bJ[\Omega]_L, \sn_X \bJ[\Omega]_\Lb=O(\ve)_{L^2_\Sigma}, \sn_X\bJ_B=O(1)+O(\ve)_{L^2_\Sigma}\\
&\sn_X(\pioh_{AB}, \pioh_{\Lb\Lb}, \pioh_{A\bN}, \fm{\Omega})=O(\ve)_{L^4_\omega},  \sn_X \pioh_{\bT A}=O(\ve)_{L^2_\Sigma}.
\end{align*}
Similar to Corollary \ref{9.2.5.23} and Proposition \ref{7.16.1.21}, we have at $t=0$
\begin{align*}
&\sn_X\sn^2\Phi, \sn_X \bd^2_{LL}\Phi, \sn_X\bd^2_{LA}\Phi=O(\ve)_{L_u^2 L_\omega^4}, \sn_X \bd^2_{\Lb A}\Phi=O(1)+O(\ve)_{L^2_u L_\omega^4}\\
&\sn^2\Phi, \bd^2_{LL}\Phi, \bd^2_{LA}\Phi, \sn_X \sn\Phi,\sn_X L\Phi =O(\ve^\f12), O(\ve)_{L^2_\Sigma}, \bd^2_{\Lb A}\Phi, X\Lb\Phi=O(1).
\end{align*}
 Using $\fm{\Omega}=\Omega\log \bb$, (\ref{10.1.3.23}) and (\ref{4.12.1.23}), substituting the above four lines of estimates and (\ref{9.23.4.22}) into (\ref{5.02.3.21_1}) and (\ref{5.18.3.21}), we deduce
\begin{align*}
X[\Box_\bg, \Omega]\Phi&=\sn_X(\pioh_{AB}\sn^2\Phi+\pioh_{\Lb\Lb}\bd^2_{LL}\Phi+\pioh_{L A}\bd^2_{\Lb A}\Phi+\pioh_{\Lb A}\bd^2_{LA} \Phi)\\
&+X(\fm{\Omega}\Box_\bg \Phi)+X(\bJ[\Omega]^\mu\p_\mu \Phi)\\
&=O(\ve)_{L_\omega^4}\Box_\bg \Phi+O(\ve)X\Box_\bg \Phi+O(\ve)_{L^2_\Sigma}.
\end{align*} 
It follows by using (\ref{9.22.6.22}) and Sobolev embedding that
\begin{equation}\label{3.22.4.24}
X[\Box_\bg, \Omega]\Phi=O(\ve)_{L^2_\Sigma}.
\end{equation}
Substituting the above estimate and (\ref{6.23.9.24}) to (\ref{3.22.7.24}) gives 
$$
X S\Lb\Omega \varrho=O(\ve)_{L^2_\Sigma}.
$$ 

Hence we conclude
\begin{equation*}
\sn_X \sn_S\sn_\Lb \sn\varrho=O(\ve)_{L^2_\Sigma}.
\end{equation*}
Using (\ref{6.23.7.24}) and the $l\le 1$ case in (\ref{10.4.10.23}), we have
\begin{equation*}
\sn_S\sn_\Omega \sn_\Lb\sn \varrho=O(\ve)_{L^2_\Sigma}
\end{equation*}
It follows by using (\ref{10.1.9.23}), (\ref{10.5.5.23}) and (\ref{3.23.1.24}) that
\begin{equation*}
\sn_\Omega^2\sn_\Lb \sn\varrho=O(\ve)_{L^2_\Sigma}.
\end{equation*}
Next we apply the first estimate in (\ref{10.4.10.23}), (\ref{10.1.8.23}), (\ref{3.22.6.24}), (\ref{9.30.6.23}), (\ref{3.22.1.24}) and the $\ze$ estimates in (\ref{3.22.3.24})  to derive
\begin{equation*}
\sn_X^2 \sn_\Lb \ep=O(\ve)_{L^2_\Sigma},
\end{equation*}
for which we also used the preliminary results in (\ref{9.23.9.22})-(\ref{4.13.3.23}). Similarly, using (\ref{1.30.2.24}), we have
\begin{equation*}
\sn_X^2\sn_\Lb \eh=O(\ve)_{L^2_\Sigma}.
\end{equation*}
This completes the proof of the second estimate in (\ref{10.4.10.23}). The last estimate follows  from using the second estimate and the commutation formula (\ref{7.04.2.21}) with details skipped here. Using Sobolev inequality, we can obtain 
\begin{equation}\label{3.23.3.24}
\sn_X^{\le 1} (\sn, \sn_L)\bA_{g,1}=O(\ve)_{L_\omega^4}.
\end{equation}
This gives the $l=1$ case for the first estimate in (\ref{10.4.1.23}) except for the case that $Y=\Lb$.

Moreover, using the $l=0$ case in (\ref{10.4.1.23}) and (\ref{4.13.3.23}), we can refine (\ref{9.22.6.22}) to obtain
\begin{equation}\label{3.23.2.24}
\Omega\Box_\bg\varrho=O(\ve)_{L_\omega^4}.
\end{equation}
We then derive by using (\ref{6.30.2.19}) that 
\begin{align*}
\Omega L\Lb \varrho&=\Omega\Big(\sD \varrho-(h-k_{\bN\bN})\Lb \varrho-\hb L \varrho+2\zb^A \sn_A \varrho\Big)=O(\ve)_{L_\omega^4}
\end{align*}
where we used (\ref{3.21.6.24}), the $l=0$ case  in (\ref{10.4.1.23}),  $\sn_\Omega\bA_b=O(\ve)_{L_\omega^4}$ in (\ref{3.22.3.24}). Using (\ref{7.17.6.21}) and $\sn\Lb \varrho=O(\ve)_{L_\omega^4}$ in (\ref{10.4.1.23}) implies $S\Omega \Lb \varrho=O(\ve)_{L_\omega^4}$.

We infer from the above  and (\ref{10.5.5.23}) that
\begin{equation}\label{6.23.11.24}
(\sn_X \sn, \sn X) \Lb \varrho=O(\ve)_{L_\omega^4},\, X=S.
\end{equation}
As the case of $X=\Omega$,  $\sn\Omega\Lb\varrho=O(\ve)_{L_\omega^4}$
  follows by using (\ref{3.21.6.24}) and assuming
\begin{equation}\label{9.22.2.24}
X\Omega[L\Phi](0)=O(\ve)_{L^4_\omega}, X=S, \Omega.
\end{equation}
To confirm (\ref{9.22.2.24}), we first follow the same procedure for the proof of (\ref{6.23.10.24}) to derive
\begin{equation*}
\Lb X\Omega[L\Phi](0)=O(\ve)_{L^2_\Sigma},
\end{equation*}
with the detailed checking omitted here. Then we can obtain (\ref{9.22.2.24}) by using the first estimate of (\ref{10.4.10.23}) and  (\ref{9.22.2.22}).

By  (\ref{7.17.6.21}) and $\sn[L\Phi]=O(\ve)_{L_\omega^4}$, we then have from (\ref{9.22.2.24})
\begin{equation*}
\Omega S[L\Phi]=O(\ve)_{L^4_\omega}.
\end{equation*}
This allows us to conclude 
\begin{equation}\label{9.22.1.24}
 (\sn_X\sn, \sn X)[L\Phi]=O(\ve)_{L^4_\omega},\, X=S, \Omega.
 \end{equation}

Due to (\ref{9.22.1.24}) and (\ref{6.23.11.24}), we have
\begin{equation*}
(\sn_X \sn, \sn X)[\Lb\Phi]=O(\ve)_{L_\omega^4}.
\end{equation*}
Hence, by combining the above estimate with (\ref{9.22.1.24}), the $l=1$ case in the second set of estimates in (\ref{10.4.1.23}) is proved. As a consequence, the estimate of $\|\sn_\Omega^2 \bA_b\|_{L_\omega^4}$ in (\ref{3.22.3.24}) is also proved. The proof of (\ref{3.22.3.24}) is complete. 

It follows by using this result together with (\ref{10.1.8.23}) that
\begin{equation*}
\sn_X\sn_\Lb \ep=O(\ve)_{L_\omega^4}.
\end{equation*}
Using (\ref{2.1.2.24}) and the result in (\ref{3.23.3.24}), we have 
\begin{align*}
\sn_X\sn_\Lb \eh=O(\ve)_{L_\omega^4}.
\end{align*}
Next we prove 
\begin{equation}\label{6.23.12.24}
\sn_X \sn_\Lb \sn\varrho=O(\ve)_{L_\omega^4}.
\end{equation}
 Due to (\ref{6.23.11.24}) 
\begin{equation*}
\sn_X \sn_\Lb \sn\varrho=\sn_X\sn\Lb \varrho+\sn_X[\Lb, \sn]\varrho=O(\ve)_{L_\omega^4}+\sn_X[\Lb, \sn]\varrho.
\end{equation*}
It remains to estimate  the second term, which can be done by using (\ref{7.04.2.21})
\begin{align*}
\sn_X[\Lb, \sn]\varrho=\sn_X(\theta\sn \varrho+\ud \bA \bN \varrho)=O(\ve)_{L_\omega^4}
\end{align*}
where to derive the last estimate we adopted (\ref{9.23.9.22})-(\ref{4.13.3.23}), $\Omega L\varrho, \Omega \Lb \varrho=O(\ve)_{L_\omega^4}$ in (\ref{10.4.1.23}) and the first estimate in (\ref{3.22.3.24}). Thus (\ref{6.23.12.24}) is proved. We conclude 
\begin{equation*}
\sn_X \sn_\Lb \bA_{g,1}=O(\ve)_{L_\omega^4}.
\end{equation*}
Therefore (\ref{10.4.1.23}) is proved. The pointwise estimate of $\sn\bA_b$ in (\ref{6.22.4.24}) follow from (\ref{3.22.3.24}) and (\ref{4.13.3.23}). The remaining estimates in (\ref{6.22.4.24}) follow as a consequence of the first estimate of (\ref{10.4.1.23}) and $\bA_{g,2}(0)=\eh(0)$. 

Finally,  the estimates of (\ref{6.23.13.24}) can be obtained by using the estimates of $\sn^{1+\le 2}\tr\chi$, $\Omega^{1+\le 2}\Lb \varrho$, estimates of $\ze$ in Proposition \ref{12.13.1.23} (3), (\ref{3.22.3.24}) and (\ref{10.4.1.23}).
\end{proof}
\begin{proof}[Proof of (\ref{9.30.15.23}) and (\ref{12.6.2.23})]

To see (\ref{9.30.15.23}), we first write
\begin{equation*}
2\ckc^{-1}u\bN \varrho+\varrho=u\ckc^{-1}(L-\Lb) \varrho+\varrho.
\end{equation*}
Using $\bN=c\p_u$ again, we derive
\begin{align*}
2u\p_u \varrho+\varrho=-2\ckc^{-1}u\Osc(c)\p_u\varrho+u \ckc^{-1}L\varrho-\tir\Lb\varrho+\varrho.
\end{align*}
Using (\ref{rarefied}), (\ref{9.30.1.23}) and (\ref{9.25.1.22}), we infer that, at $t=0$, 
\begin{align}\label{10.5.1.23}
2u \p_u \varrho+\varrho+O(\ve)-u \ckc^{-1}L\varrho>q_0.
\end{align}
Using (\ref{9.30.10.23}), we derive 
\begin{align*}
\|u^\f12 \ckc^{-1} L\varrho\|_{L^1[u, u_*]}&\les \|L\varrho\|_{L_u^2 L_\omega^\infty}(\int_u^{u_*} u')^\f12\les \ve (u_*^2-u^2)^\f12\approx \ve(u_*-u)^\f12.  
\end{align*}
Hence, integrating (\ref{10.5.1.23}) in $u$ and using the above estimate, we derive 
\begin{equation*}
(u^\f12 \varrho)(u_*)-u^\f12\varrho\ges (q_0-C\ve)(u_*^\f12-u^\f12)-C'\ve(u_*-u)^\f12
\end{equation*}
where $C, C'>0$ are two universal constants, which implies
$$
u^\f12 \varrho+C_9 q_0(u_*^\f12-u^\f12)\le C_9(C\ve(u_*^\f12-u^\f12)+C'\ve(u_*-u)^\f12), 
$$
with some constants $C_9>0$. 
Thus $-\ve^\f12 \les \varrho(0)\les \ve$ due to $\varrho(0)=O(\ve^\f12)$ in (\ref{7.3.2.24}). (\ref{9.30.15.23}) is proved. 

In view of  (\ref{12.20.4.21}), (\ref{9.18.1.23}) and (\ref{9.30.1.23}), using (\ref{9.29.5.23}), we derive at $t=0$ that
\begin{align*}
\int_{u}^{u_*} \overline{-\bT\varrho} du=-\overline{v_\hN}(u)+\int_u^{u_*}\frac{2}{u'}\overline{ v_\hN} du'=O(\ve^\f12)
\end{align*}
which implies 
\begin{equation}\label{6.24.2.24}
\int_u^{u_*}\overline{-\Lb\varrho} du'=O(\ve^\f12).
\end{equation}
Using (\ref{rarefied}) and integrating in $u$, we derive
\begin{align*}
\int_u^{u_*}\overline{-\Lb \varrho+\tir^{-1}\varrho}\ge q_0\int_u^{u_*}\tir^{-1} 
\end{align*}
which gives $q_0\les \ve^\f12$ due to (\ref{6.24.2.24}) and using $\varrho(0)=O(\ve)_{L_u^2 L_\omega^\infty}$.

Due to (\ref{rarefied}), (\ref{9.25.1.22}) and $\varrho=O(\ve^\f12)$
\begin{equation*}
0<q_0\le -\tir\Lb \varrho+\varrho\les 1. 
\end{equation*}
Hence due to (\ref{9.30.15.23})
\begin{equation*}
q_0-C\ve<-\tir\Lb \varrho\les 1
\end{equation*}
This implies
\begin{equation}\label{6.24.1.24}
0\le [-\tir \Lb \varrho]_+\les 1, -\ve\les [-\tir \Lb \varrho]_-\le 0.
\end{equation}
Using the above estimate, it is direct to compute
\begin{align*}
\int_{u_0}^{u_*}\overline{-\Lb \varrho}^2&\le \int_{u_0}^{u_*}((\overline{[-\Lb\varrho]_+})^2+(\overline{[-\Lb\varrho]_-})^2)\le \int_{u_0}^{u_*}\overline{[-\Lb\varrho]_+}^2 du'+O(\ve^2)\\
&\les \int_{u_0}^{u_*}\overline{[-\Lb\varrho]_+} du'+O(\ve^2).
\end{align*}
We derive from (\ref{6.24.2.24}) and (\ref{6.24.1.24}) that 
\begin{equation*}
 \int_{u_0}^{u_*}\overline{[-\Lb\varrho]_+} du\les \ve^\f12-\int_{u_0}^{u_*}\overline{[-\Lb\varrho]_-} du\les \ve^\f12.
\end{equation*}
Hence 
$$
\int_{u_0}^{u_*}\overline{-\Lb \varrho}^2=O(\ve^\f12).
$$
Using \Poincare inequality and the first estimate in (\ref{10.1.9.23}), we conclude
\begin{equation*}
\int_{u_0}^{u_*}\|\Lb \varrho\|^2_{L_\omega^\infty}=O(\ve^\f12).
\end{equation*}
Hence  (\ref{12.6.2.23}) is proved.
\end{proof}
 Thus the proof of Proposition \ref{12.21.1.21} is completed.
 
\section{Decompositions of Riemann curvature components}\label{geostru2}
We crucially rely on geometric structures of various components of Riemann curvature in our analysis. In this section, we derive them in details.
\begin{proposition}\label{1.25.5.22}
Let $X$, $U$ and $V$ be in the null tetrad, $Y$ be either $L$ or $\Lb$. There hold symbolically
\begin{equation}\label{1.18.1.22}
U_1\bg_{\mu\la}U_2^\la U_3^\mu=[\sn\Phi],\quad X\bg_{\mu \la} e_A^\la Y^\mu=\sn \Phi^\dagger,\quad X\bg_{\mu \la} U^\la V^\mu=X \Phi^\dagger
\end{equation}
where $U_1, U_2, U_3$ are in null tetrad, and only one of them is $\{e_A\}_{A=1}^2$;
and
\begin{equation}\label{1.27.4.22}
\tensor{\Ga}{^L_{\Lb L}},  \tensor{\Ga}{^\Lb_L_L}=[L\Phi], \quad  \tensor{\Ga}{^\Lb_\Lb_L}=[\Lb \Phi],\quad \tensor{\Ga}{^L_L_L}=[\Lb\Phi] +[L\Phi].
\end{equation}
If $X, Y\in\{L, \Lb\}$  and $U$ is in null tetrad, then
\begin{equation}\label{1.18.3.22}
\tensor{\Ga}{^\a_A_\ga} X_\a Y^\ga,\tensor{\Ga}{^C_X_Y}=[\sn\Phi];\quad \tensor{\Ga}{^C_B_U}+U\log c\delta^C_B=(1-|\vs(U)|)\sn\log c\c \Pi, 
\end{equation}
\begin{equation}\label{1.21.1.22}\left\{
\begin{array}{lll}
 \tensor{\Ga}{^\Lb_A_B}+\f12L\log c\delta_{AB}=\f12 c^{-2}\sn_A v_B,\quad \tensor{\Ga}{^L_A_B}+\f12 \Lb \log c\delta_{AB}=\f12 c^{-2}\sn_A v_B,\\ \widehat{\tensor{\Ga}{^Y_A_B}}=-\f12\hk_{AB}, \quad Y=L, \Lb.
\end{array}\right.
\end{equation}
\begin{align}
&\tensor{\Ga}{^A_C_\la}\tensor{\Ga}{^\la_Y_B}-\tensor{\Ga}{^A_Y_\la}\tensor{\Ga}{^\la_C_B}=\bp \Phi^\dagger\c [\sn\Phi], \quad Y=L, \Lb\label{1.22.3.22}\\
&\tensor{\Ga}{^A_\Lb_\a}\tensor{\Ga}{^\a_B_L}- \tensor{\Ga}{^\a_\Lb_B}\tensor{\Ga}{^A_\a_L}=[\sn \Phi]^2.\label{1.28.3.22}
\end{align}
\begin{equation}\label{1.22.1.22}\left\{
\begin{array}{lll}
\sn (\tensor{\Ga}{^A_L_B})=-\sn L\log c \delta_{AB},\quad \sn_L (\tensor{\Ga}{^A_C_B})=\sn_L \sn\log c \c \Pi\\
 \widehat{\sn_L \tensor{\Ga}{^Y_B_A}}=-\f12\sn_L \hk_{AB},\quad \sn_B\tensor{\Ga}{^Y_L_A}=\sn[\sn\Phi], \, Y=L, \Lb.
\end{array}\right.
\end{equation}
\begin{align}
&\sn_\Lb \tensor{\Ga}{^A_B_L}=-\Lb L \log c \delta^A_B, \quad\sn_L\tensor{\Ga}{^A_B_\Lb}=-L \Lb\log c \delta^A_B\label{1.28.1.22}\\
&\bR_{ABCL}=(\sn(L \log c)+\sn_L \sn \log c)\c \Pi+[\sn\Phi](\chi, \chib, \bp\Phi^\dagger),
\label{1.21.2.22}\\
&\bR_{ABC\Lb}=(\sn(\Lb \log c)+\sn_\Lb \sn \log c)\c \Pi+[\sn\Phi](\chi, \chib, \bp\Phi^\dagger)+\ze \bp\Phi^\dagger,\label{9.2.1.22}\\
&\bR_{AB43}=(\zeta+[\sn\Phi])[\sn\Phi]\label{1.27.2.22}\\
&\widehat{\bR_{A4B4}}=\sn[\sn\Phi]+\sn_L \hk_{AB}+[\sn\Phi]^2+\hk_{AB}(\fB+\tr\chi)+(\chih, \chibh)(\fB+\hk_{AB})\label{1.30.2.22}\\
&\widehat{\bR_{A3B4}}=\sn[\sn\Phi]+\sn_L \hk_{AB}+[\sn\Phi]^2+\hk_{AB}(\fB+\tr\chi)+(\chih, \chibh)(\fB+\hk_{AB}),\label{5.21.3.23}\\
& \ga^{AB}\bR_{A43B}=\sn \bA_{g,1}+\varpi, \quad\varpi=(\bA+\frac{1}{\tir}+\fB)\fB\label{4.17.1.24}
\end{align}
where  for $S_{t,u}$ tangent $2$-tensor $F$, $\widehat{F}$ in the above standards for $\widehat{F}_{AB}=F_{AB}-\f12\ga_{AB}\ga^{CD}F_{CD}$.
\end{proposition}
\begin{proof}
 (\ref{1.18.1.22}) follows from (\ref{10.6.1.22}) and (1) in Proposition \ref{6.7.1.23}. We further recall from (\ref{10.6.1.22})
\begin{equation}\label{1.18.2.22}
X \bg_{\mu\la}e_A^\la Y^\mu=-c^{-2}X(v_i) e_A^i, \quad X \bg_{\mu\la}e_A^\la e_B^\mu=-2X\log c\ga_{AB}.
\end{equation}


Next, we prove (\ref{1.27.4.22}) by calculating Christoffel symbols in view of (\ref{10.6.1.22})
\begin{align*}
\tensor{\Ga}{^\rho_\mu_\nu}L_\rho \Lb^\mu L^\nu&=\f12 \bg^{L\a}(-\p_\a \bg_{\mu\nu}+\p_\mu\bg_{\a\nu}+\p_\nu \bg_{\a\mu})\Lb^\mu L^\nu\\
&=\f12 \bg^{L\Lb}\p_\nu \bg_{\a\mu}\Lb^\mu L^\nu \Lb^\a=-\frac{1}{4}\p_L \bg_{\a\mu}\Lb^\a \Lb^\mu=[L\Phi]
\end{align*}
and
\begin{align*}
&\tensor{\Ga}{^\rho_\mu_\nu} \Lb_\rho L^\mu L^\nu=-\frac{1}{4}\p_L \bg_{\a\nu} L^\nu L^\a=[L\Phi]\\
&\tensor{\Ga}{^\rho_\mu_\nu} \Lb_\rho \Lb^\mu L^\nu=-\frac{1}{4}\p_\Lb \bg_{\a\nu}L^\nu L^\a=[\Lb\Phi]\\
&\tensor{\Ga}{^\rho_\mu_\nu} L_\rho L^\mu L^\nu=-\frac{1}{4}(-\p_\Lb \bg_{\mu\nu}L^\mu L^\nu+2\p_L \bg_{\a\mu}\Lb^\a L^\mu)=[\Lb\Phi]+[L\Phi].
\end{align*}
Hence, (\ref{1.27.4.22}) is proved.

 Note
$$
\tensor{\Ga}{^\a_A_\ga}=\f12 \bg^{\a\b}(\p_s \bg_{\b \ga}+\p_\ga\bg_{\b s}-\p_\b \bg_{s\ga}) e_A^s,
$$
\begin{equation}\label{5.6.1.23}
\tensor{\Ga}{^C_\mu_\nu} =\f12 \bg^{C\a}(-\p_\a \bg_{\mu\nu}+\p_\mu \bg_{\a\nu}+\p_\nu\bg_{\a\mu}).
\end{equation}
Since $X_\a, Y^\ga$ are both null frames, the first identity in (\ref{1.18.3.22}) follows by using the first identity in (\ref{1.18.1.22}).
 Contracted by $ U^\mu e_B^\nu$ in (\ref{5.6.1.23}),
\begin{align*}
\tensor{\Ga}{^C_B_U}=\f12 \bg^{C\a}(-\p_\a \bg_{\mu\nu}+\p_\nu\bg_{\a\mu}) U^\mu e_B^\nu+\f12 e_C^\a e_B^\nu  U \bg_{\a\nu}.
\end{align*}
If $U$ is a null vector, the terms of bracket vanishes due to the first identity (\ref{1.18.2.22}) and $\sn_C v_B$ is symmetric. Using the first identity in (\ref{1.18.2.22}) we obtain the case that $U=L, \Lb$ of the last formula in (\ref{1.18.3.22}). If $U=e_D$, we apply the second identity of (\ref{1.18.2.22}) to obtain
 $
\tensor{\Ga}{^C_B_D}=-\sn\log c \c \Pi.
 $
Hence the second identity of (\ref{1.18.3.22}) is proved. 

Next we prove (\ref{1.21.1.22}). Note
\begin{align*}
\tensor{\Ga}{^L_A_B}=\f12 \bg^{L\a}(-\p_\a \bg_{\mu\nu}+\p_\mu \bg_{\a\nu}+\p_\nu \bg_{\a\mu})e_A^\mu e_B^\nu.
\end{align*}
By the first identity in (\ref{1.18.2.22}), for the last two terms, we obtain
\begin{equation*}
\bg^{L\a}(\p_\mu \bg_{\a\nu}+\p_\nu \bg_{\a\mu})e_A^\mu e_B^\nu= c^{-2}\sn_A v_B.
\end{equation*}
For the first term, we derive by the second identity in (\ref{1.18.2.22}) that
\begin{equation*}
\bg^{L\a}\p_\a \bg_{\mu\nu}e_A^\mu e_B^\nu=-\f12 \Lb \bg_{ij} e_A^i e_B^j= \Lb \log c \delta_{AB}.
\end{equation*}
This gives the second identity of (\ref{1.21.1.22}). Similarly, by changing $L$ to $\Lb$ in the above calculation and noting
\begin{equation*}
\bg^{\Lb\a}\p_\a \bg_{\mu\nu}e_A^\mu e_B^\nu=-\f12 L \bg_{ij} e_A^i e_B^j=L\log c \delta_{AB},
\end{equation*}
 we derive the first identity in (\ref{1.21.1.22}). The last formula in (\ref{1.21.1.22}) follows immediately as the consequence of the two identities. 

 To see (\ref{1.22.3.22}), we write
\begin{align*}
\tensor{\Ga}{^A_C_\la} \tensor{\Ga}{^\la_Y_B}-\tensor{\Ga}{^A_Y_\la}\tensor{\Ga}{^\la_C_B}&=\tensor{\Ga}{^A_C_L} \tensor{\Ga}{^L_Y_B}-\tensor{\Ga}{^A_Y_L}\tensor{\Ga}{^L_C_B}\\
&+\tensor{\Ga}{^A_C_\Lb} \tensor{\Ga}{^\Lb_Y_B}-\tensor{\Ga}{^A_Y_\Lb}\tensor{\Ga}{^\Lb_C_B}\\
&+\tensor{\Ga}{^A_C_D} \tensor{\Ga}{^D_Y_B}-\tensor{\Ga}{^A_Y_D}\tensor{\Ga}{^D_C_B}.
\end{align*}
It follows by
using (\ref{1.18.3.22}) and (\ref{1.21.1.22}) that the first two lines are $[\sn\Phi]\bp\Phi^\dagger$, and the last line vanishes.
Thus (\ref{1.22.3.22}) is proved.

To prove (\ref{1.28.3.22}), we derive
\begin{align*}
\tensor{\Ga}{^A_\Lb_\a}\tensor{\Ga}{^\a_B_L}- \tensor{\Ga}{^\a_\Lb_B}\tensor{\Ga}{^A_\a_L}&=\tensor{\Ga}{^A_\Lb_C}\tensor{\Ga}{^C_B_L}- \tensor{\Ga}{^C_\Lb_B}\tensor{\Ga}{^A_C_L}\\
&+\tensor{\Ga}{^A_\Lb_L}\tensor{\Ga}{^L_B_L}- \tensor{\Ga}{^L_\Lb_B}\tensor{\Ga}{^A_L_L}\\
&+\tensor{\Ga}{^A_\Lb_\Lb}\tensor{\Ga}{^\Lb_B_L}- \tensor{\Ga}{^\Lb_\Lb_B}\tensor{\Ga}{^A_\Lb_L}.
\end{align*}
The last two lines are $[\sn\Phi]^2$ by using (\ref{1.18.3.22}).
For the first pair of difference, due to (\ref{1.18.3.22}),
 $\tensor{\Ga}{^A_\Lb_C}=-\Lb\log c\delta^A_C, \tensor{\Ga}{^C_B_L}=-L \log c\delta^C_B$, thus this line vanishes.
  We hence conclude (\ref{1.28.3.22}).

Next we prove (\ref{1.22.1.22}). The first pair follows from (\ref{1.18.3.22}). The second line is a direct consequence of (\ref{1.21.1.22}) and the first identity in (\ref{1.18.3.22}). (\ref{1.28.1.22}) is a consequence of the second identity in (\ref{1.18.3.22}).

To prove (\ref{1.21.2.22}), we first recall that
\begin{equation}\label{1.27.3.22}
\tensor{\bR}{^\rho_\sigma_\mu_\nu}=\p_\mu\tensor{\Ga}{^\rho_\nu _\sigma}-\p_\nu \tensor{\Ga}{^\rho_\mu_\sigma}+\tensor{\Ga}{^\rho_\mu_\la}\tensor{\Ga}{^\la_\nu_\sigma}-\tensor{\Ga}{^\rho_\nu_\la}\tensor{\Ga}{^\la_\mu_\sigma}
\end{equation}
and \begin{align}\label{5.7.4.23}
\p_\mu \tensor{\Ga}{^\rho_\nu_\sigma}&=\bd_\mu\tensor{\Ga}{^\rho_\nu_\sigma}-\tensor{\Ga}{^\rho_\mu_\a}\tensor{\Ga}{^\a_\nu_\sigma}+\tensor{\Ga}{^\a_\mu_\nu}\tensor{\Ga}{^\rho_\a_\sigma}+\tensor{\Ga}{^\a_\mu_\sigma}\tensor{\Ga}{^\rho_\nu_\a}.
\end{align}
It is straightforward to derive
\begin{align*}
{e_A}_\rho e_B^\sigma e_C^\mu L^\nu \tensor{\bR}{^\rho_\sigma_\mu_\nu}&= {e_A}_\rho e_B^\sigma e_C^\mu L^\nu\left(\p_\mu\tensor{\Ga}{^\rho_\nu _\sigma}-\p_\nu \tensor{\Ga}{^\rho_\mu_\sigma}\right)+\tensor{\Ga}{^A_C_\la}\tensor{\Ga}{^\la_L_B}-\tensor{\Ga}{^A_L_\la}\tensor{\Ga}{^\la_C_B}\\
e_{A\rho} e_C^\mu e_B^\sigma L^\nu\bd_\nu \tensor{\Ga}{^\rho_\mu_\sigma}&=e_{A\rho} e_C^\mu e_B^\sigma L^\nu\p_\nu \tensor{\Ga}{^\rho_\mu_\sigma}+\tensor{\Ga}{^A_L_\a}\tensor{\Ga}{^\a_C_B}-\tensor{\Ga}{^\a_C_L}\tensor{\Ga}{^A_\a_B}-\tensor{\Ga}{^\a_L_B}\tensor{\Ga}{^A_C_\a}\\
e_{A\rho} e_C^\mu e_B^\sigma L^\nu\bd_\mu \tensor{\Ga}{^\rho_\nu_\sigma}&=e_{A\rho} e_C^\mu e_B^\sigma L^\nu\p_\mu \tensor{\Ga}{^\rho_\nu_\sigma}+\tensor{\Ga}{^A_C_\a}\tensor{\Ga}{^\a_L_B}-\tensor{\Ga}{^\a_C_L}\tensor{\Ga}{^A_\a_B}-\tensor{\Ga}{^\a_C_B}\tensor{\Ga}{^A_L_\a}.
\end{align*}
Hence we obtain
\begin{align}\label{1.22.2.22}
\tensor{\bR}{^A_B_C_L}&={e_A}_\rho e_B^\sigma e_C^\mu L^\nu\left(\bd_\mu\tensor{\Ga}{^\rho_\nu _\sigma}-\bd_\nu \tensor{\Ga}{^\rho_\mu_\sigma}\right)-(\tensor{\Ga}{^A_C_\la}\tensor{\Ga}{^\la_L_B}-\tensor{\Ga}{^A_L_\la}\tensor{\Ga}{^\la_C_B}).
\end{align}
The quadratic terms have been treated in (\ref{1.22.3.22}).
It remains to consider the second order terms.
We first prove
\begin{equation}\label{5.7.1.23}
\begin{split}
&e_{A\rho} e_C^\mu e_B^\sigma L^\nu\bd_\mu \tensor{\Ga}{^\rho_\nu_\sigma}-\sn_C \tensor{\Ga}{^A_L_B}=(\chi+\chib+L\log c\c \Pi)\c [\sn \Phi]\\
&e_{A\rho} e_C^\mu e_B^\sigma L^\nu \bd_\nu \tensor{\Ga}{^\rho_\mu_\sigma}-\sn_L\tensor{\Ga}{^A_C_B}=\bp \Phi^\dagger\c[\sn \Phi].
 \end{split}
\end{equation}
For the first identity, we note
\begin{align*}
\sn_C \tensor{\Ga}{^A_L_B}=e_{A\rho} e_C^\mu e_B^\sigma L^\nu\bd_\mu \tensor{\Ga}{^\rho_\nu_\sigma}+{e^A}_\mu\bd_C\Pi^\mu_\nu\tensor{\Ga}{^\nu_ L_B}+\bd_C L^\nu \tensor{\Ga}{^A _\nu_B}+\bd_C \Pi_\mu^\nu \tensor{\Ga}{^A_L_\nu}e_B^\mu.
\end{align*}
Using Proposition \ref{6.7con} and (\ref{1.18.3.22}), we derive
\begin{equation*}
{e^A}_\mu\bd_C\Pi^\mu_\nu\tensor{\Ga}{^\nu_ L_B}, \bd_C \Pi_\mu^\nu \tensor{\Ga}{^A_L_\nu}e_B^\mu=(\chi, \chib)\c [\sn\Phi]\, \quad \bd_C L^\nu \tensor{\Ga}{^A _\nu_B}=[\sn\Phi]\c (\chi, L\log c\c \Pi).
\end{equation*}
This gives the first line of (\ref{5.7.1.23}).

In view of (\ref{6.29.5.19}), symbolically,
\begin{equation*}
\sn_L \tensor{\Ga}{^A_C_B}-e_{A\rho} e_C^\mu e_B^\sigma L^\nu\bd_\nu \tensor{\Ga}{^\rho_\mu_\sigma}=\zb\c \bp\Phi^\dagger
\end{equation*}
where we  used the fact that $\Ga$ in null tetrad takes the form $\bp\Phi^\dagger$.
This gives the second identity in (\ref{5.7.1.23}).

Combining (\ref{1.22.3.22}), the first line of (\ref{1.22.1.22}), (\ref{5.7.1.23}) and (\ref{1.22.2.22}) gives (\ref{1.21.2.22}).  To see (\ref{9.2.1.22}), $Y$ in the higher order terms $\sn Y\log c, Y\sn \log c$ is changed from $L$ to $\Lb$. The calculations for the lower order terms can be done similarly by using Proposition \ref{6.7con} and (\ref{1.22.3.22}). This allows us to obtain (\ref{9.2.1.22}).

To prove (\ref{1.27.2.22}),
in view of (\ref{1.27.3.22}), we write
\begin{equation}\label{7.17.6.22}
\tensor{\bR}{^A_{B43}}=e_{A\rho}e_B^\sigma L^\mu\Lb^\nu(\p_\mu\tensor{\Ga}{^\rho_\nu_\sigma}-\p_\nu \tensor{\Ga}{^\rho_\mu_\sigma})+\tensor{\Ga}{^A_L_\la}\tensor{\Ga}{^\la_\Lb_B}-\tensor{\Ga}{^A_\Lb_\la}\tensor{\Ga}{^\la_L_B}.
\end{equation}
Due to (\ref{1.28.3.22}), the quadratic terms of Christoffel symbols are good terms $[\sn\Phi]^2$.
It only remains to consider the two second order derivatives. 
We apply (\ref{5.7.4.23}) and use (\ref{1.28.3.22}) to see that
\begin{align}
&(\p_\mu\tensor{\Ga}{^\rho_\nu_\sigma}-\bd_\mu\tensor{\Ga}{^\rho_\nu_\sigma})e_{A\rho}e_B^\sigma L^\mu\Lb^\nu-(\p_\nu \tensor{\Ga}{^\rho_\mu_\sigma}-\bd_\nu\tensor{\Ga}{^\rho_\mu_\sigma})e_{A\rho}e_B^\sigma L^\mu\Lb^\nu\nn\\
&=-2\tensor{\Ga}{^A_\Lb_\a}\tensor{\Ga}{^\a_B_L}+2 \tensor{\Ga}{^\a_\Lb_B}\tensor{\Ga}{^A_\a_L}=[\sn\Phi]^2.\label{1.28.2.22}
\end{align}
Hence we conclude from (\ref{7.17.6.22}) that
\begin{equation}\label{5.21.1.23}
\bR_{AB43}=e_{A\rho}e_B^\sigma L^\mu\Lb^\nu(\bd_\mu\tensor{\Ga}{^\rho_\nu_\sigma}-\bd_\nu \tensor{\Ga}{^\rho_\mu_\sigma})+[\sn\Phi]^2.
\end{equation}
Next we show
\begin{equation}\label{7.17.5.22}
\begin{split}
\sn_\Lb \tensor{\Ga}{^A_L_B}&-\tensor{\Ga}{^A_L_B} k_{\bN\bN}-\bd_\Lb\tensor{\Ga}{^A_L_B}, \sn_L \tensor{\Ga}{^A_\Lb_B}-\tensor{\Ga}{^A_\Lb_B} k_{\bN\bN}-\bd_L\tensor{\Ga}{^A_\Lb_B}\\
&=[\sn \Phi](\zeta+\zb).
\end{split}
\end{equation}

Indeed,
\begin{align*}
\sn_\Lb \tensor{\Ga}{^A_L_B}&= e_A^\mu\bd_\Lb(\Pi_{\mu\rho} \tensor{\Ga}{^\rho_\nu_\sigma}L^\nu \Pi^\sigma_{\sigma'})e_B^{\sigma'}\\
&=e_A^\mu (\bd_\Lb \Pi_{\mu\rho}\tensor{\Ga}{^\rho_\nu_\sigma}L^\nu\Pi^\sigma_{\sigma'}+\Pi_{\mu\rho}\bd_\Lb L^\nu\Pi^\sigma_{\sigma'} \tensor{\Ga}{^\rho_\nu_\sigma}\\
&+\Pi_{\mu\rho}\tensor{\Ga}{^\rho_\nu_\sigma}L^\nu \bd_\Lb \Pi^\sigma_{\sigma'})e_B^{\sigma'}+e_A^\mu\Pi_{\mu\rho}L^\nu \Pi^\sigma_{\sigma'}\bd_\Lb \tensor{\Ga}{^\rho_\nu_\sigma}e_B^{\sigma'}.
\end{align*}
In view of (\ref{6.29.6.19}) and (\ref{6.29.7.19}), symbolically,
\begin{align*}
e_A^\mu\bd_\Lb \Pi_{\mu\rho}\tensor{\Ga}{^\rho_\nu_\sigma}L^\nu e_B^\sigma&=\f12 \l \bd_\Lb \Lb, e_A\r \tensor{\Ga}{^L_L_B}+\f12\l \bd_\Lb L, e_A\r\tensor{\Ga}{^\Lb_L_B}
\end{align*}
which gives
\begin{equation*}
e_A^\mu\bd_\Lb \Pi_{\mu\rho}\tensor{\Ga}{^\rho_\nu_\sigma}L^\nu e_B^\sigma=(\zeta+\zb)\c [\sn \Phi].
\end{equation*}
Similarly, we obtain by using (\ref{6.29.6.19}) and (\ref{1.18.3.22}) that
\begin{align*}
{e_A}_\rho \tensor{\Ga}{^\rho_\nu_B}\bd_\Lb L^\nu=\tensor{\Ga}{^A_\nu_B}(2\zeta_C e_C^\nu+k_{\bN\bN}L^\nu)=k_{\bN\bN}\tensor{\Ga}{^A_L_B}+\zeta\c \sn\log c.
\end{align*}
By using (\ref{6.29.6.19}), (\ref{6.29.7.19}) and (\ref{1.18.3.22}), we obtain
\begin{align*}
{e_A}_\rho \tensor{\Ga}{^\rho_\nu_\sigma}L^\nu \bd_\Lb \Pi^\sigma_{\sigma'} e_B^{\sigma'}&=\f12 \tensor{\Ga}{^A_L_\sigma}\bd_\Lb(\Lb_{\sigma'}L^\sigma+L_{\sigma'}\Lb^\sigma)e_B^{\sigma'}\\
&=\f12 (\tensor{\Ga}{^A_L_L}\l \bd_\Lb \Lb, e_B\r+\tensor{\Ga}{^A_L_\Lb}\l \bd_\Lb L, e_B\r)\\
&=(\zeta_B+\zb_B)[\sn\Phi].
\end{align*}
Combining the above estimates implies
\begin{equation*}
\sn_\Lb \tensor{\Ga}{^A_L_B}-k_{\bN\bN}\tensor{\Ga}{^A_L_B}-\bd_\Lb\tensor{\Ga}{^A_L_B}=[\sn \Phi](\zeta+\zb).
\end{equation*}
This shows the first estimate in (\ref{7.17.5.22}), the other follows similarly by swapping $L$ and $\Lb$ in the above calculations.

Using (\ref{1.28.1.22}) we have
\begin{align*}
\sn_\Lb \tensor{\Ga}{^A_B_L}-\sn_L\tensor{\Ga}{^A_B_\Lb}=-[\Lb, L] \log c \delta^A_B.
\end{align*}
Moreover, applying (\ref{3.19.2}) to $f=\log c$ yields
\begin{equation*}
[L, \Lb]\log c=2(\zb^A-\zeta^A)\sn_A \log c-2k_{\bN\bN} \bN \log c.
\end{equation*}
From the above two lines, we deduce
\begin{align*}
\sn_\Lb \tensor{\Ga}{^A_B_L}&-\sn_L\tensor{\Ga}{^A_B_\Lb}+2k_{\bN\bN} \bN \log c \delta_B^A\\
&= 2(\zb^A-\zeta^A)\sn_A \log c+(\zeta+\zb)[\sn \Phi].
\end{align*}
Substituting (\ref{1.18.3.22}) and the above identity to (\ref{7.17.5.22}) leads to
\begin{align*}
(\bd_\mu\tensor{\Ga}{^\rho_\nu_\sigma}-\bd_\nu\tensor{\Ga}{^\rho_\mu_\sigma})e_{A\rho}e_B^\sigma L^\mu\Lb^\nu=(\zeta+\zb)[\sn \Phi].
\end{align*}
(\ref{1.27.2.22}) follows immediately as a consequence of the above line and (\ref{5.21.1.23}).

Next we prove (\ref{1.30.2.22}). Contracting (\ref{1.27.3.22}) by  $\Lb_\rho e_A^\sigma L^\mu e_B^\nu$, we have
\begin{equation*}
 \Lb_\rho e_A^\sigma L^\mu e_B^\nu \tensor{\bR}{^\rho_\sigma_\mu_\nu}= \Lb_\rho e_A^\sigma L^\mu e_B^\nu (\p_\mu\tensor{\Ga}{^\rho_\nu _\sigma}-\p_\nu \tensor{\Ga}{^\rho_\mu_\sigma})+\tensor{\Ga}{^\Lb_L_\la}\tensor{\Ga}{^\la_B_A}-\tensor{\Ga}{^\Lb_B_\la}\tensor{\Ga}{^\la_L_A}.
\end{equation*}
Moreover, applying (\ref{5.7.4.23}) again, we write
\begin{align}\label{5.7.2.23}
\begin{split}
&(\p_\mu \tensor{\Ga}{^\rho_\nu_\sigma}-\bd_\mu \tensor{\Ga}{^\rho_\nu_\sigma})L^\mu \Lb_\rho e_A^\sigma e_B^\nu=-\tensor{\Ga}{^\Lb_L_\a}\tensor{\Ga}{^\a_B_A}+\tensor{\Ga}{^\a_L_B} \tensor{\Ga}{^\Lb_\a_A}+\tensor{\Ga}{^\a_L_A} \tensor{\Ga}{^\Lb_B_\a}\\
&(\p_\nu \tensor{\Ga}{^\rho_\mu_\sigma}-\bd_\nu \tensor{\Ga}{^\rho_\mu_\sigma})L^\mu \Lb_\rho e_A^\sigma e_B^\nu=-\tensor{\Ga}{^\Lb_B_\a}\tensor{\Ga}{^\a_L_A}+\tensor{\Ga}{^\a_B_L}\tensor{\Ga}{^\Lb_\a_A}+\tensor{\Ga}{^\a_B_A}\tensor{\Ga}{^\Lb_L_\a}.
\end{split}
\end{align}
Hence
\begin{align*}
2\ti E_{AB}&:=\Big((\p_\mu \tensor{\Ga}{^\rho_\nu_\sigma}-\bd_\mu \tensor{\Ga}{^\rho_\nu_\sigma})-(\p_\nu \tensor{\Ga}{^\rho_\mu_\sigma}-\bd_\nu \tensor{\Ga}{^\rho_\mu_\sigma})\Big)L^\mu \Lb_\rho e_A^\sigma e_B^\nu\\
&=2\tensor{\Ga}{^\a_L_A} \tensor{\Ga}{^\Lb_B_\a}-2\tensor{\Ga}{^\a_B_A}\tensor{\Ga}{^\Lb_L_\a}.
\end{align*}
We then combine the above three calculations to obtain
\begin{equation}\label{5.7.5.23}
 \Lb_\rho e_A^\sigma L^\mu e_B^\nu \tensor{\bR}{^\rho_\sigma_\mu_\nu}= \Lb_\rho e_A^\sigma L^\mu e_B^\nu (\bd_\mu\tensor{\Ga}{^\rho_\nu _\sigma}-\bd_\nu \tensor{\Ga}{^\rho_\mu_\sigma})+\ti E_{AB}.
\end{equation}
Next we decompose $\ti E_{AB}=-C_{AB}+D_{AB}$  where each term is an $S_{t,u}$-tangent $2$-tensor
 $$C_{AB}=\tensor{\Ga}{^\a_B_A}\tensor{\Ga}{^\Lb_L_\a}; \quad D_{AB}=\tensor{\Ga}{^\a_L_A}\tensor{\Ga}{^\Lb_B_\a}.$$
Note
\begin{equation*}
C_{AB}=\sum_{Y=L, \Lb}\tensor{\Ga}{^Y_B_A}\tensor{\Ga}{^\Lb_L_Y}+\tensor{\Ga}{^C_B_A}\tensor{\Ga}{^\Lb_L_C}.
\end{equation*}
Using (\ref{1.18.3.22}), the last term is $\sn \log c\c \Pi\c [\sn\Phi]$.
Also using (\ref{1.21.1.22}) and (\ref{1.27.4.22}), the traceless part of the first term is $\hk_{AB} ([L\Phi]+[\Lb\Phi])$. Thus
\begin{equation*}
\widehat C_{AB}=\hk_{AB}([L \Phi]+[\Lb\Phi])+\sn \log c\c \Pi\c [\sn\Phi].
\end{equation*}
Using (\ref{1.18.3.22}) and (\ref{1.21.1.22}), we derive
\begin{align*}
D_{AB}&=\tensor{\Ga}{^C_L_A}\c \tensor{\Ga}{^\Lb_B_C}+\tensor{\Ga}{^L_L_A}\c\tensor{\Ga}{^\Lb_B_L}+\tensor{\Ga}{^\Lb_L_A}\c\tensor{\Ga}{^\Lb_B_\Lb}\\
&=-\f12 (L \log c)^2\delta_{AB}+L\log c k_{AB}+[\sn \Phi]^2,
\end{align*}
and hence
\begin{equation*}
\widehat{D}_{AB}=D_{AB}-\f12\ga_{AB}\ga^{CD}D_{CD}=L\log c \hk_{AB}+[\sn\Phi]^2.
\end{equation*}
Hence
\begin{equation}\label{5.21.2.23}
\widehat{\ti E}_{AB}=([L\Phi]+[\Lb\Phi])\hk_{AB}+[\sn\Phi]^2.
\end{equation}
Next by Proposition \ref{6.7con}, we derive
\begin{equation*}
\bd_L \tensor{\Ga}{^\Lb_B_A}=\sn_L \tensor{\Ga}{^\Lb_B_A}-\bd_L \Lb_\rho \tensor{\Ga}{^\rho_B_A}-\zb_B \tensor{\Ga}{^\Lb_L_A}-\zb_A \tensor{\Ga}{^\Lb_B_L}.
\end{equation*}
Symbolically, the traceless part of the second term on the right is $\zb\c \sn\log c \c \Pi+k_{\bN\bN} \hk_{BA}$ by using Proposition \ref{6.7con}, (\ref{1.18.3.22}) and (\ref{1.21.1.22}). The last two terms are $\zb\c[\sn\Phi]$ by using (\ref{1.18.3.22}).
Hence, using (\ref{1.22.1.22}) we have
\begin{align*}
\widehat{\bd_L \tensor{\Ga}{^\Lb_B_A}}&=\widehat{\sn_L\tensor{\Ga}{^\Lb_B_A}}+\zb\c \sn\log c \c \Pi+k_{\bN\bN}\hk_{BA}\\
&=-\f12\sn_L \hk_{AB}+\zb\c \sn\log c \c \Pi+k_{\bN\bN} \hk_{BA}.
\end{align*}
Similarly, it follows by using  (\ref{7.21.3.19}) and (\ref{6.29.7.19}) that
\begin{align*}
\bd_B\tensor{\Ga}{^\Lb_L_A}&=\sn_B \tensor{\Ga}{^\Lb_L_A}-\bd_B \Lb_\a \tensor{\Ga}{^\a_L_A}-\tensor{\Ga}{^\Lb_\a_A}\bd_B L^\a-\tensor{\Ga}{^\Lb_L_\a}(\bd_B e_A^\a-\sn_B e_A^\a)\\
&=\sn_B\tensor{\Ga}{^\Lb_L_A}-\chib_{BC}\tensor{\Ga}{^C_L_A}-\chi_{BC}\tensor{\Ga}{^\Lb_C_A}-\f12\chi_{AB}\tensor{\Ga}{^\Lb_L_\Lb}-\f12 \chib_{AB}\tensor{\Ga}{^\Lb_L_L}
\end{align*}
and, also by using (\ref{1.27.4.22})-(\ref{1.21.1.22}) and the fact that $\tr\eta=[L\Phi]$, 
\begin{equation*}
\widehat{\bd_B\tensor{\Ga}{^\Lb_L_A}}=\widehat{\sn_B\tensor{\Ga}{^\Lb_L_A}}+(\chih_{AB}, \chibh_{AB})([L\Phi]+[\Lb\Phi]+c^{-2}\eh)+\tr\chi c^{-2}\eh.
\end{equation*}
Using (\ref{1.22.1.22}), the first term on its right-hand side is $\sn [\sn\Phi]$. In view of (\ref{5.7.5.23}), (\ref{5.21.2.23}) and using the fact that $\chih+\chibh=-2\hk_{AB}$, we hence obtained the identity in (\ref{1.30.2.22}).

Next we prove (\ref{5.21.3.23}) by contracting (\ref{1.27.3.22}) with $L^\mu L_\rho e_A^\sigma e_B^\nu$.
Applying (\ref{5.7.4.23}) we write
\begin{align*}
&(\p_\mu \tensor{\Ga}{^\rho_\nu_\sigma}-\bd_\mu \tensor{\Ga}{^\rho_\nu_\sigma})L^\mu L_\rho e_A^\sigma e_B^\nu=-\tensor{\Ga}{^L_L_\a}\tensor{\Ga}{^\a_B_A}+\tensor{\Ga}{^\a_L_B} \tensor{\Ga}{^L_\a_A}+\tensor{\Ga}{^\a_L_A} \tensor{\Ga}{^L_B_\a}\\
&(\p_\nu \tensor{\Ga}{^\rho_\mu_\sigma}-\bd_\nu \tensor{\Ga}{^\rho_\mu_\sigma})L^\mu L_\rho e_A^\sigma e_B^\nu=-\tensor{\Ga}{^L_B_\a}\tensor{\Ga}{^\a_L_A}+\tensor{\Ga}{^\a_B_L}\tensor{\Ga}{^L_\a_A}+\tensor{\Ga}{^\a_B_A}\tensor{\Ga}{^L_L_\a}.
\end{align*}
Hence with $$
\ud E_{AB}=\tensor{\Ga}{^\a_L_A} \tensor{\Ga}{^L_B_\a}-\tensor{\Ga}{^\a_B_A}\tensor{\Ga}{^L_L_\a},
$$
we write
\begin{align*}
\big((\p_\mu &\tensor{\Ga}{^\rho_\nu_\sigma}-\bd_\mu \tensor{\Ga}{^\rho_\nu_\sigma})-(\p_\nu \tensor{\Ga}{^\rho_\mu_\sigma}-\bd_\nu \tensor{\Ga}{^\rho_\mu_\sigma})\big)L^\mu L_\rho e_A^\sigma e_B^\nu=2\ud E_{AB}.
\end{align*}
Similar to (\ref{5.7.5.23}), we obtain by the definition of curvature that
\begin{equation}\label{12.9.3.23}
 L_\rho e_A^\sigma L^\mu e_B^\nu \tensor{\bR}{^\rho_\sigma_\mu_\nu}= L_\rho e_A^\sigma L^\mu e_B^\nu (\bd_\mu\tensor{\Ga}{^\rho_\nu _\sigma}-\bd_\nu \tensor{\Ga}{^\rho_\mu_\sigma})+ \ud E_{AB}.
\end{equation}
Next we write $\ud E_{AB}=G_{AB}-H_{AB}$ with
\begin{equation*}
G_{AB}=\tensor{\Ga}{^\a_L_A} \tensor{\Ga}{^L_B_\a}, \quad H_{AB}=\tensor{\Ga}{^\a_B_A}\tensor{\Ga}{^L_L_\a}.
\end{equation*}
Similar to the treatment of $C_{AB}$ and $D_{AB}$,
we derive by using (\ref{1.27.4.22})-(\ref{1.21.1.22}) that
\begin{equation}\label{12.9.2.23}
\widehat{\ud E}_{AB}=\widehat{G}_{AB}-\widehat{H}_{AB}=[\sn\Phi]^2+\hk_{BA}([L\Phi]+[\Lb \Phi]).
\end{equation}
Next we follow the same proof as the previous case to obtain
\begin{align*}
\bd_L \tensor{\Ga}{^L_B_A}=\sn_L \tensor{\Ga}{^L_B_A}-\bd_L L_\rho \tensor{\Ga}{^\rho_B_A}-\zb_B \tensor{\Ga}{^L_L_A}-\zb_A \tensor{\Ga}{^L_B_L}.
\end{align*}
Taking the traceless part, in view of  (\ref{1.18.3.22}), (\ref{1.21.1.22}) and (\ref{1.22.1.22}) we deduce
\begin{equation}\label{12.9.1.23}
\widehat{\bd_L \tensor{\Ga}{^L_B_A}}=\sn_L\hk_{AB}+k_{\bN\bN}\hk_{AB}+\zb[\sn\Phi]
\end{equation}
Using Proposition \ref{6.7con} again, we derive
\begin{align*}
\bd_B\tensor{\Ga}{^L_L_A}&=\sn_B\tensor{\Ga}{^L_L_A}-\chi_{BC}\tensor{\Ga}{^C_L_A}-\chib_{BC}\tensor{\Ga}{^L_C_A}-\f12\chi_{AB}\tensor{\Ga}{^L_L_\Lb}-\f12 \chib_{AB}\tensor{\Ga}{^L_L_L}.
\end{align*}
Then using (\ref{1.18.3.22}), (\ref{1.21.1.22}) and (\ref{1.22.1.22}) again, we derive symbolically that
\begin{align*}
\widehat{\bd_B\tensor{\Ga}{^L_L_A}}&=\sn[\sn\Phi]+(\chih, \chibh)([L\Phi]+[\Lb\Phi]+c^{-2}\eh)+\tr\chi c^{-2}\eh.
\end{align*}
Combining the above formula with (\ref{12.9.3.23})-(\ref{12.9.1.23}), we obtain (\ref{5.21.3.23}).

We can refer to \cite[Section 8]{rough_fluid} for the proof of (\ref{4.17.1.24}). For completeness, we give the details below. We first compute with the help of Bianchi identity that
\begin{equation*}
\delta^{AB}\bR_{B43A}=\delta^{AB}(\bR_{AB}-\delta^{CD}\bR_{ACBD}).
\end{equation*}
 For $\delta^{AB}\bR_{AB}$ we use
$
\delta^{AB}\bR_{AB}=\sn_A (\Xi_A)+\varpi=\sn^2\varrho+\varpi
$
which follows from (\ref{ricc6.7.1})  together with (\ref{4.10.2.19}), with the formula of $\Xi_i$ given in (\ref{6.14.1.19}). 

Recall the direct decomposition of Riemann curvature, 
\begin{equation*}
\bR_{\a\b\ga \d}=\bd_\a \cir{\pi}_{\b\d\ga}+\bd_\b \cir{\pi}_{\a\ga\d}-\bd_\a \cir{\pi}_{\b\ga\d}-\bd_\b \cir{\pi}_{\d\a\ga}+E_{\a\b\ga\d}
\end{equation*}
with $E=\bg\c \bp \bg \c \bp \bg$ and $\cir{\pi}_{\a\b\ga}=\bp_\ga \bg_{\a\b}$. It is straightforward to check by using (\ref{10.6.1.22}) that $\cir{\pi}_{CD\bN}=\bN\log c \ga_{CD}$ and $\cir{\pi}_{C\bN D}=0$. 
Using this fact, contracting the above identity by the orthonormal frame on $\T S_{t,u}$ gives, in view of (\ref{10.6.1.22}) and Proposition \ref{6.7con}, we have 
  \begin{equation*}
  \bR_{ABCD}=\sn^2\log c+\varpi.
  \end{equation*}
Summarizing the above calculations gives (\ref{4.17.1.24}).


\end{proof}

\section{Some useful geometric calculations}\label{geocal}
In this section, we give some geometric calculations, which are based on the estimates established before Section \ref{mul_1} and the geometric formulas given in Section\ref{geostru2}. The material of this section is frequently used in treating the term $v^i \tensor{\ud\ep}{^a_i_l}\Pic^l_m$ and its derivatives appearing in the components of $\pioh_{\bT A}$, $\bJ[\Omega]$.   
\begin{lemma}
For $\Sigma$-tangent tensor $F$, with $Y=L, \Lb$ or  $\{e_A\}$,  there hold symbolically that 
\begin{equation}\label{10.17.1.23}
\begin{split}
\sn_Y F_A&=Y (F_{l'}) e_A^{l'}+ Y\log c F_A+c\Big(\bA_{g,1}\\
&+\max(-\vs(Y),0)\ze +(1-|\vs(Y)|)\thetac(Y)\Big) F_\bN.
\end{split}
\end{equation}
\end{lemma}
Indeed, for $\Sigma$-tangent tensor $F$, we directly compute
\begin{align*}
&\bd_Y F_A=Y F_{l'}e_A^{l'}-\tensor{\Ga}{^m_{Yl'}}F_m e_A^{l'},\quad \sn_Y\Pi_l^{l'}{e_A}^l=-\bd_Y \bN_A \bN^{l'}
\end{align*}
where we used Proposition \ref{6.7con} to obtain the second identity.
Hence using (\ref{1.18.3.22}) we have 
\begin{align}
\sn_Y F_A&=Y (F_{l'}) e_A^{l'}-\tensor{\Ga}{^C_{Yl'}}F_C e_A^{l'}-\tensor{\Ga}{^\bN_{Yl'}}F_\bN e_A^{l'}+\bd_Y \bN F_{\bN}\nn\\
&=Y (F_{l'}) e_A^{l'}+ Y \log c F_A+(\bd_Y \bN_A-\tensor{\Ga}{^\bN_{YA}})\c F_\bN\nn\\
&=Y (F_{l'}) e_A^{l'}+ Y\log c F_A+Y \bN_A \c F_\bN.\label{1.30.1.22}
\end{align}
Using (\ref{11.30.2.23}) we obtain (\ref{10.17.1.23}).

As a direct consequence, if $Y=L, \Lb$, using Lemma \ref{5.13.11.21} (1)-(2), (5), and (\ref{L2conndrv}) we bound
\begin{equation}\label{11.5.3.23}
\begin{array}{lll}
\sn_L(\tensor{\ud\ep}{^a_A_B})=O(L\log c)+O(\bA_{g,1}), \sn_\Lb(\tensor{\ud\ep}{^a_A_B})=O(\fB)+O(\ud \bA); \sn_\Omega\tensor{\ud\ep}{^a_A_B}=O(1)\\
\sn_S\sn_{\tir\Lb}(\tensor{\ud\ep}{^a_A_B})=O(1)+O(\log \l t\r\Delta_0)_{L^2_u L_\omega^2}, \quad\sn_X^2(\tensor{\ud\ep}{^a_A_B})=O(\l t\r^\delta\Delta_0)_{L^2_\Sigma}.
\end{array}
\end{equation}

Next we introduce the notations that
\begin{equation}\label{7.19.3.22}
{}\rp{a}v^*_l=v^i \tensor{\ud\ep}{^a_i_l}\Pic^l_m, \quad {}\rp{a}v^\sharp=v^i \tensor{\ud\ep}{^a_i_l}\hat \bN^l.
\end{equation}
For convenience we may write the contraction with $\tensor{\ud\ep}{^a_i^l}$ by $\sta{a}\wedge$, and drop the superscript $a$ in ${}\rp{a} v^*$ and ${}\rp{a}v^\sharp$ for short whenever there occurs no confusion.  We note that ${}\rp{a}\pih_b=c^{-1}{}\rp{a}v^*$. 
\begin{lemma}\label{4.1.1.23}
Under the assumptions of (\ref{3.12.1.21})-(\ref{1.25.1.22}), for ${}\rp{a}v^*$ and ${}\rp{a}v^\sharp$ defined in (\ref{7.19.3.22}), there hold
\begin{equation}\label{3.6.1.22}\left\{
\begin{array}{lll}
\snc_C {}\rp{a}v^*&=\snc_C(v^i \tensor{\ud\ep}{^a_i^l}\Pic_{ml})e_A^m=\snc_C(v^i) \tensor{\ud\ep}{^a_i^l} \Pic_{ml}e_A^m-\thetac_C^A \tensor{\ud\ep}{^a_i^l}\hat\bN_l v^i\\
&=-\f12\tr\thetac {}\rp{a}v^\sharp\Pic-\hat\thetac {}\rp{a}v^\sharp +\snc v\sta{a}{\wedge}\Pic\\
\snc{}\rp{a}v^\sharp&=v\sta{a}\wedge \thetac+\sn v\sta{a}\wedge \hat \bN=\f12\tr\thetac {}\rp{a}v^*+v\sta{a}\wedge \hat\thetac +\sn v\sta{a}\wedge \hat \bN\\
\divc{}\rp{a}v^*&=-\tr\thetac {}\rp{a}v^\sharp\\
\displaybreak[0]
\snc^2{}\rp{a}v^*&=\snc^2 v\sta{a}\wedge \Pic-\snc v\sta{a}\wedge \hat \bN \thetac-\snc(v\sta{a}\wedge\hat \bN \thetac)\\
&=-v\sta{a}\wedge\thetac\c\thetac-v\sta{a}\wedge \hat \bN \snc \thetac+\snc^2 v\sta{a}\wedge \Pic-2\snc v\sta{a}\wedge \hat \bN \thetac.
\end{array}\right.
\end{equation}
\begin{align}
\sn{}\rp{a}v^*&=\snc v\sta{a}{\wedge}\Pi+\snc\log c {}\rp{a}v^*+\thetac {}\rp{a}v^\sharp\label{6.24.2.21}\\
\sn_Y{}\rp{a}v^*&=Y v\sta{a}\wedge\Pi+Y\log c{}\rp{a}v^*+(\bA_{g,1}+\max(-\vs(Y), 0) \ze){}\rp{a}v^\sharp.\label{6.28.6.21}
\end{align}
where $Y=L, \Lb$. 
\begin{equation}\label{7.25.1.22}\left\{
\begin{array}{lll}
\hat \bN{}\rp{a}v^\sharp&=\hat \bN v\stackrel{a}\wedge \hat\bN+v\stackrel{a}\wedge\snc\log (\bb c)=O(\bA_{g,1}, \l t\r^{-1+\delta}\Delta_0^\f12\ud\bA)\\
\sn_L{}\rp{a}v^\sharp&= Lv\sta{a}\wedge\hat \bN+{}\rp{a}v^*[\sn\Phi]+{}\rp{a}v^\sharp L(c^{-1})=O(\l t\r^{-2+\delta}\Delta_0).
\end{array}\right.
\end{equation}
With $Y=L, \Lb$, and $X=S, \Omega$,
\begin{align}
\sn_X\sn_Y {}\rp{a}v^*&=XY\log c{}\rp{a}v^*+Y\log c\c \sn_X {}\rp{a}v^*+X Yv\sta{a}\wedge\Pi+Y v\sta{a}\wedge \hN(\bA_{g,1}\nn\\
&+(1-|\vs(X)|\thetac(X)))+\sn_X(\bA_{g,1}+\max(-\vs(Y),0)\ze){}\rp{a}v^\sharp\nn\\
&+(\bA_{g,1}+\max(-\vs(Y),0)\ze)X{}\rp{a}v^\sharp\label{7.19.6.22}
\end{align}
where in the terms on the right-hand side, we neglected the factors of $c^m, m\in \mathbb Z$.

In view of the above formulas, we have the following rough estimates
\begin{align}\label{3.28.3.24}
|\sn_X(v^*), X(v^\sharp)|\les |X^{\le 1}v|, \quad \sn_X^2(v^\sharp, v^*)=O(X^{\le 2} v)+O(\l t\r^{-2+2\delta}\Delta_0^\frac{3}{2})_{L^2_u L_\omega^2}.
\end{align}
\end{lemma}

\begin{proof}
The  identities in (\ref{3.6.1.22}) follow from direct calculations and  $\p_i v_j =\p_j v_i$.

(\ref{6.24.2.21}) follows in view of (\ref{1.27.1.22}), $\Ga(\gac)-\Ga(\ga)=\p_\omega\log c$ symbolically, and the symbolic formula for  $S_{t,u}$-tangent tensor $F$, $\sn F=\snc F+(\Ga(\gac)-\Ga(\ga))\c F$.

 (\ref{6.28.6.21}) is obtained by applying (\ref{1.30.1.22}) to $F=v^i \tensor{\ud\ep}{^a_i^l}$.

Next we prove (\ref{7.25.1.22}).
The first line follows by direct calculation with the help of $\hN(\hN)=\snc\log(\bb c)=\ud\bA$.
Note
\begin{equation*}
L {}\rp{a}v^\sharp=L v\sta{a}{\wedge} \hN+v\sta{a}{\wedge} L \hN.
\end{equation*}
For $L \hat \bN$ in the second term, it suffices to consider $c^{-1}L \bN^i$, with the extra term $L\varrho \c v^\sharp$ arisen from the second term. Due to Lemma \ref{dg} and Proposition \ref{6.7con},
\begin{equation*}
\bd_L \bN^i=L \bN^i+\Ga_{L\bN}^i=L\bN^i=[\sn\Phi].
\end{equation*}
Thus
\begin{equation*}
L {}\rp{a}v^\sharp= L v\sta{a}\wedge \hat\bN+ {}\rp{a}v^*[\sn\Phi]+{}\rp{a}v^\sharp L(c^{-1})
\end{equation*}
with the desired estimate followed from using Lemma \ref{5.13.11.21}. Hence (\ref{7.25.1.22}) is proved. 

Finally, (\ref{7.19.6.22}) is obtained by differentiating (\ref{6.28.6.21}) with the help of (\ref{10.17.1.23}).  The lower order estimate in (\ref{3.28.3.24}) is a direct consequence of the formulas for the first order derivatives in (\ref{3.6.1.22})-(\ref{7.25.1.22}), by using (\ref{3.6.2.21}) and (\ref{6.24.1.21}) in Lemma \ref{5.13.11.21}. The higher order estimate in (\ref{3.28.3.24}) can be derived by direct differentiation followed with using Lemma \ref{5.13.11.21} and (\ref{L2BA2}).
\end{proof}
\section{Appendix: Geometric calculations for the rotation vector-fields}\label{append}
In this section, we give the proofs of Proposition \ref{2.19.4.22} and Proposition \ref{3.22.6.21} in Section \ref{rotation}.
\begin{proof}[Proof of Proposition \ref{2.19.4.22}]

(\ref{3.22.5.21}) follows directly by combining (\ref{3.18.1.21}) and (\ref{3.22.4.21}).

We first prove (\ref{3.18.1.21}).
Note that
\begin{align*}
\nab_\bN {\Omega\rp{a}}^i&=\bN {\Omega\rp{a}}^i+\bN^m\Ga_{mn}^i{\Omega\rp{a}}^n\\
&=\bN {\Omega\rp{a}}^i+\bN^m(\p_m \varphi \delta_n^i+\p_n \varphi \delta_m^i-\nab^i \varphi g_{mn}) {\Omega\rp{a}}^n
\end{align*}
where $\varphi=-\log c$ and $\nab$ is the Levi-Civita connection of the induced metric $g$ on $\Sigma$.  We then obtain
\begin{equation*}
\l \nab_\bN {\Omega\rp{a}}, \bN\r=\l c^{-1}\bN {\Omega\rp{a}}, c^{-1}\bN\r_e-\Omega\rp{a}(\log c).
\end{equation*}
For the first term on the right-hand side of the above identity, noting that ${}\rp{a}\O$ is a killing vector field in Minkowski space, we deduce
\begin{align*}
0&=\l c^{-1} \bN {}\rp{a}\O, c^{-1}\bN\r_e=\l c^{-1}\bN(\Omega\rp{a}+\la\rp{a} \hat \bN), c^{-1}\bN\r_e\\
&=\l c^{-1}\bN(\Omega\rp{a}), c^{-1}\bN\r_e+c^{-1}\bN\la\rp{a} \l c^{-1}\bN, c^{-1}\bN\r_e+c^{-1}\la\rp{a}\l \bN(c^{-1}\bN), c^{-1}\bN\r_e\\
&=\l c^{-1}\bN(\Omega\rp{a}), c^{-1}\bN\r_e+c^{-1}\bN\la\rp{a} \l c^{-1}\bN, c^{-1}\bN\r_e.
\end{align*}
Computing the left-hand side implies
\begin{equation*}
\l\nab_\bN(\Omega\rp{a}), \bN\r=-\l\Omega\rp{a}, \nab_\bN\bN\r=\Omega\rp{a}\log \bb.
\end{equation*}
Hence we conclude (\ref{3.18.1.21}).

Next we prove (\ref{3.18.2.21}). By using the definition of $\la$,
\begin{align*}
\la^a\l \hat\bN, \hat \bN\r_e&=\l {}\rp{a}\O-\Omega\rp{a}, \hat \bN\r_e=\l {}\rp{a}\O, \hat \bN\r_e=c^{-1}\ep_{alk} x^l \bN^k.
\end{align*}
This gives (\ref{3.18.2.21}) in view of the definition of $y^k$.

Next we show (\ref{2.10.2.22}). In view of the definition (\ref{3.18.2.21}), it is straightforward to derive
\begin{align*}
\hat e_A(\la\rp{a})&=\hat e_A(x^l \tensor{\ud\ep}{^a_l_k}\c c^{-1}y^k)\\
&=\hat e_A (c^{-1}y^k) x^l \tensor{\ud\ep}{^a_l_k}+\hat e_A(x^l)\tensor{\ud\ep}{^a_l_k} c^{-1}y^k\\
&=(\thetac_A^k-\hat e_A(c^{-1}\frac{x^k}{\tir}))x^l \tensor{\ud\ep}{^a_l_k}+\hat e_A^l\tensor{\ud\ep}{^a_l_k}\c c^{-1}y^k.
\end{align*}
We hence conclude (\ref{2.10.2.22}).

Next we consider $\bT\la\rp{a}$. Note
\begin{equation*}
\l\bT {}\rp{a}\O, c^{-1}\bN\r_e +\l c^{-1}\bN {}\rp{a}\O, \bT\r_e=0.
\end{equation*}
By using (\ref{3.19.1.21}), we have
\begin{align*}
0=\l \bT{}\rp{a}\Omega, c^{-1}\bN\r_e+\l c^{-1}\bN {}\rp{a}\Omega, \bT\r_e+\l \bT(\la\rp{a} \hat \bN), c^{-1}\bN\r_e+\l \hat \bN(\la\rp{a}\hat \bN), \bT\r_e.
\end{align*}
This implies
\begin{align}
 \bT\la\rp{a}&=-\big(\l\bT{}\rp{a}\Omega, c^{-1}\bN\r_e+\l c^{-1}\bN {}\rp{a}\Omega, \bT\r_e\big)\nn\\
 &=\l{}\rp{a}\Omega, \bT(c^{-1}\bN)\r_e+\l {}\rp{a}\Omega, c^{-1}\bN \bT\r_e\nn\\
 &=\l{}\rp{a}\Omega, c^{-1}[\bT, \bN]\r_e+2\l{}\rp{a}\Omega, c^{-1}\bN\bT\r_e.\label{3.22.2.21}
\end{align}
Note due to Lemma \ref{dg}
\begin{equation*}
\bd_\bN \bT=\bN \bT+\bN^n \bT^\b\tensor{\Ga}{_{n\b}^i}=\bN \bT
\end{equation*}
Using (\ref{6.29.5.19})-(\ref{6.29.7.19}), we have $\l\bd_\bN \bT, \Omega\r=\zb({}\rp{a}\Omega)$. Thus, the second term on the right-hand side of (\ref{3.22.2.21}) equals $2c \zb({}\rp{a}\Omega)$.
Similarly,
\begin{align*}
\l \bd_\bT \bN, {}\rp{a}\Omega\r&=\frac{1}{4}\l \bd_L L+\bd_\Lb L-\bd_L \Lb-\bd_\Lb \Lb, {}\rp{a}\Omega\r\\
&=\frac{1}{4}\l 2\zeta_A e_A-2\zb_A e_A-(-2\zeta_A+2k_{\bN A}), {}\rp{a}\Omega\r=\zeta({}\rp{a}\Omega).
\end{align*}
Hence we can obtain from (\ref{3.22.2.21}) that
\begin{equation*}
\bT\la\rp{a}=c(\zeta({}\rp{a}\Omega)+\zb(\rp{a}\Omega))=c{}\rp{a}\Omega \log \bb.
\end{equation*}
This is (\ref{3.22.4.21}).
\end{proof}

 \begin{proof}[Proof of Proposition \ref{3.22.6.21}]
 (\ref{5.02.2.21}) follows as a direct consequence by using the formulas for ${}\rp{a}\pi$.

 The first identity follows directly from computation, by using (\ref{6.29.5.19}). Similarly, by using (\ref{6.29.5.19}) and (\ref{6.29.6.19})
 \begin{align*}
 {}\rp{a}\pi_{L\Lb}&=\l \bd_L {}\rp{a}\Omega, \Lb\r+\l \bd_\Lb {}\rp{a}\Omega, L\r\\
 &=-\l{}\rp{a}\Omega, \bd_L \Lb\r-\l {}\rp{a}\Omega, \bd_\Lb L\r=-2(\zeta+\zb)^A {}\rp{a}\Omega_A.
 \end{align*}
 By using (\ref{6.29.7.19}) and (\ref{3.19.1.21}), we have
 \begin{align*}
{}\rp{a}\pi_{\Lb\Lb}&=2\l\bd_\Lb {}\rp{a}\Omega, \Lb\r=-2\l {}\rp{a}\Omega, \bd_\Lb\Lb\r=4{}\rp{a}\Omega^A(\zeta_A-k_{\bN A}).
 \end{align*}
 \begin{align*}
 {}\rp{a}\pi_{\bT A}&=\l \bd_A {}\rp{a}\Omega, \bT\r+\l \bd_\bT {}\rp{a}\Omega, e_A\r\\
 &=\l \bd_\bT({}\rp{a}\O-\la\rp{a}\hat \bN), e_A\r-\l \bd_A \bT, {}\rp{a}\Omega\r\\
 &=\l\bd_ {{}\rp{a}\O}\bT+[\bT, {}\rp{a}\O], e_A\r-c^{-1}\la\rp{a} \zeta_A+k_{A{}\rp{a}\Omega}\\
 &=\l[\bT, {}\rp{a}\O], e_A\r-c^{-1}\la\rp{a} (\zeta_A+k_{A\bN}).
 \end{align*}
 Using ${}\rp{a}\O=\tensor{\ud\ep}{^a_i^j}x^i\p_j$, and $[\p_l, {}\rp{a}\O]=\tensor{\ud\ep}{^a_l^j}\p_j$,
\begin{align*}
[\bT, {}\rp{a}\O]&=[v^l\p_l, {}\rp{a}\O]=v^l[\p_l, {}\rp{a}\O]-{}\rp{a}\O(v^l)\p_l\\
&=v^j\tensor{\ud\ep}{^a_j^l}\p_l-{}\rp{a}\O(v^l)\p_l.
\end{align*}
Hence,
\begin{align*}
{}\rp{a}\pi_{\bT A}=v^j \tensor{\ud\ep}{^a_j^l} {e_A}_l -{}\rp{a}\O(v^l) {e_A}_l -c^{-1}\la\rp{a} (\zeta_A+k_{A  \bN}).
\end{align*}
 We first have by the formula of the conformal transform $g=e^{2\varphi}\delta_e$,
 \begin{equation*}
 {}\rp{X}\pi(g)=2X\varphi g+{}\rp{X}\pi(\delta_e).
 \end{equation*}
 Since $g=c^{-2}\delta_e$, we will apply the above formula with $\varphi=-\log c$. Using (\ref{5.17.1.21}), we deduce
\begin{align*}
{}\rp{a}\pi_{AB}&=\l\bd_A \Omega\rp{a}, e_B\r+\l \bd_B \Omega\rp{a}, e_A\r\\
&={}\rp{{}\rp{a}\O}\pi_{AB}-\{\l \bd_A (\la\rp{a}\hat \bN), e_B\r+\l \bd_B (\la\rp{a}\hat \bN), e_A\r\}\\
&=-2{}\rp{a}\O\log c \ga_{AB}-2c^{-1}\la\rp{a} \theta_{AB}\\
&=-2c^{-1}\la\rp{a} \hat\theta_{AB}-2({}\rp{a}\Omega \log c+c^{-1}\la\rp{a} L\log c) \ga_{AB}\\
&+(c^{-1} \la\rp{a} \p_C v_D \delta^{CD}-c^{-1}\la\rp{a} \tr\chi)\ga_{AB}\\
&=-2c^{-1}\la\rp{a} \hat \theta_{AB}+\f12 {}\rp{a}\ss\ga_{AB}
\end{align*}
where  ${}\rp{a}\ss:=\tr{}\rp{a}\sl{\pi}$. Due to ${}\rp{\O}\pi_{\bN A}=0$,  
\begin{align*}
{}\rp{a}\pi_{\bN A}&=\l \bd_\bN {}\rp{a}\Omega, e_A\r+\l \bd_A {}\rp{a}\Omega, \bN\r\\
&=\l \bd_\bN ({}\rp{a}\O-\la\rp{a} \hat \bN), e_A\r+\l \bd_A ({}\rp{a}\O-\la\rp{a} \hat\bN), \bN\r\\
&= \l \bd_\bN {}\rp{a}\O, e_A\r+\l\bd_A{}\rp{a}\O, \bN\r-c^{-1}\la\rp{a} \l\bd_\bN \bN, e_A\r -\sn_A(c^{-1}\la\rp{a})\\
&=c^{-1}\la\rp{a}\sn_A\log \bb -\sn_A(c^{-1}\la\rp{a}).
\end{align*}
(\ref{5.02.2.21})-(\ref{5.6.02.21}) are direct consequences of the formulas of ${}\rp{a}\pi_{AB}$ and (\ref{ricc_def}).

Next we prove (\ref{5.6.03.21}).
\begin{align*}
\sn_{k'}(\tensor{\ud\ep}{^i_l^m}x^l \Pi_{mn})\Pi_{n'}^n&=\Pi_{k'}^k\{\p_k(\tensor{\ud\ep}{^i_l^m}x^l \Pi_{mn})-\Ga_{kn}^{l'}\tensor{\ep}{^i_l^m}x^l \Pi_{ml'}\}\Pi_{n'}^n\\
\displaybreak[0]
&=\Pi_{k'}^k\{\tensor{\ud\ep}{^i_k^m}\Pi_{mn}+\tensor{\ud\ep}{^i_l^m}x^l \p_k \Pi_{mn}-\Ga_{kn}^{l'}\tensor{\ud\ep}{^i_l^m}x^l \Pi_{ml'}\}\Pi_{n'}^n\\
&=\Pi_{k'}^k\big(\tensor{\ud\ep}{^i_k^m}\Pi_{mn}+\tensor{\ud\ep}{^i_l^m}x^l (\nab_k \Pi_{mn}+\Ga_{km}^{l'}\Pi_{l'n})\big)\Pi_{n'}^n.
\end{align*}
For the second term on the right hand side, noting $\bN^m=c^2 \bN_m$, we calculate 
\begin{align*}
\tensor{\ud\ep}{^i_l^m}x^l \sn_{k'} \Pi_{mn}\Pi_{n'}^n&=-\tensor{\ud\ep}{^i_l^m}x^l(\sn_{k'} \bN_m \bN_n+\bN_m \sn_{k'} \bN_n)\Pi_{n'}^n\\
&=-\tensor{\ud\ep}{^i_l^m} x^l \bN_m\theta_{k'n}\Pi^n_{n'}=-c^{-1}\la\rp{i} \theta_{k'n} \Pi^n_{n'}.
\end{align*}
For the last term, we compute
\begin{align*}
\tensor{\ud\ep}{^i_l^m}x^l\Ga_{km}^{l'}\Pi_{l'n'}&=-\tensor{\ud\ep}{^i_l^m} x^l \p_m \log c \Pi_{kn'}-\tensor{\ud\ep}{^i_l^m}x^l \p_k \log c \Pi_{mn'}+\tensor{\ud\ep}{^i_l_k} x^l \nab^{l'}\log c \Pi_{l'n'}\\
&=-{}\rp{i}\O\log c\Pi_{kn'}-{}\rp{i}\Omega_{n'} \p_k \log c+\tensor{\ud\ep}{^i_l_k} x^l \sn_{n'} \log c.
\end{align*}
Thus
\begin{equation*}
\tensor{\ud\ep}{^i_l^m}x^l\Ga_{km}^{l'}\Pi_{l'n'}\Pi^{k}_{k'}=-{}\rp{i}\O\log c \Pi_{k'n'}-\sn_{k'} \log c {}\rp{i}\Omega_{n'}+\sn_{n'}\log c{}\rp{i}\Omega_{k'}.
\end{equation*}
Hence we can have
\begin{align*}
(\sn_k {}\rp{i}\Omega)_n&=\Pi_{kk'}\tensor{\ud\ep}{^{ik'm}}\Pi_{mn}-c^{-1}\la\rp{i}\theta_{kn}\nn\\
&-{}\rp{i}\O\log c \Pi_{kn}-\sn_k \log c {}\rp{i}\Omega_{n}+\sn_{n}\log c{}\rp{i}\Omega_k.
\end{align*}
Using ${}\rp{i}\O={}\rp{i}\Omega+c^{-1}\la \bN$, and $\theta_{AB}=\chi_{AB}+k_{AB}$, we can obtain (\ref{5.6.03.21}) by using (\ref{5.17.1.21}). (\ref{2.10.1.22}) and (\ref{6.8.1.22}) can be found in \cite{shock_demetrios}. (\ref{5.13.10.21}) can be obtained by direct calculations with the help of Proposition \ref{6.7con} and (\ref{lb}).

Differentiating (\ref{2.10.2.22}) implies
\begin{align*}
\snc_B \snc_A(\la\rp{a})&=\snc_B\left((\thetac_A^C-(c\tir)^{-1}\delta_A^C){}\rp{a}\Omega_C-\ud\ep_{akj}\frac{y^k}{c}\delta_A^j\right)\\
&=\snc_B (\thetac_A^C-(c\tir)^{-1}\delta_A^C){}\rp{a}\Omega_C+(\thetac_A^C-(c\tir)^{-1}\delta_A^C)\snc_B {}\rp{a}\Omega_C-\ud\ep_{akj}\snc_B(c^{-1}y^k)\hat e_A^j.
\end{align*}
Note
\begin{align*}
\ud\ep_{akj}\snc_B(c^{-1}y^k) \hat e_A^j&=\ud\ep_{akj}(\thetac_B^k-(c\tir)^{-1}\delta_B^k)\hat e_A^j-\ud\ep_{akj}\snc_B (c^{-1}) \frac{x^k}{\tir}\hat e_A^j.
\end{align*}
Also using (\ref{2.10.1.22}), combining the above calculations gives (\ref{2.10.3.22}).
\begin{align*}
\snc_B \snc_A\la\rp{a}&=\snc_B (\thetac_A^k) {}\rp{a}\Omega_k+(\thetac_A^k-(c\tir)^{-1}\delta_A^k)(\Pi_B^n\Pi_k^m \ep_{anm}-\la\rp{a} \thetac_{Bk})\\
&+(\thetac_B^k-(c\tir)^{-1}\delta_B^k)\Pi_A^n\Pi_k^m \ep_{anm}.
\end{align*}
This gives (\ref{2.10.3.22}). Taking trace and using $\stc\sdiv\thetac=\snc \tr\thetac$ (see \cite[Page 56 (3.1.2b)]{CK}), we obtain (\ref{12.17.1.23}).

Finally we prove (\ref{6.8.1.22}).
\begin{align*}
&\snc_A(\Pi_l^n \Pi_i^m \tensor{\ud\ep}{^a_n_m})e_A^l\\
&=-\snc_A\hat\bN_l \hat\bN^n \Pi_i^m \tensor{\ud\ep}{^a_n_m}e_A^l-(\snc_A \hat \bN_i \hat \bN^m+\hat\bN_i \snc_A\hat\bN^m)\tensor{\ud\ep}{^a_n_m} e_A^n\\
&=-\thetac_{Al}\hat \bN^n \Pi_i^m \tensor{\ud\ep}{^a_n_m}e_A^l-(\thetac_{Ai} \hat\bN^m +\thetac_A^m \hat\bN_i)\tensor{\ud\ep}{^a_n_m}e_A^n\\
&=-(r^{-1}{}\rp{a}\Omega_i+y'^n\Pi_i^m \tensor{\ud\ep}{^a_n_m}) \tr\thetac+e_A^l\Pi_l^n(\f12 \tr\thetac \Pic_{Ai}+\hat\thetac_{Ai})\hat \bN^m\tensor{\ud\ep}{^a_m_n} -e_A^n \hat\thetac_A^m \hat \bN_i \tensor{\ud\ep}{^a_n_m}.
\end{align*}
Since the last term vanishes, in view of (\ref{2.10.1.22}),  we obtain
\begin{align*}
\sDc\rp{a}\Omega_i=-\f12 \tr\thetac(r^{-1}{}\rp{a}\Omega_i+\tensor{\ud\ep}{^a_m_n}\Pi_i^n {y'}^m)+\hat \thetac_{Ai}\hat\bN^m \tensor{\ud\ep}{^a_m_n}e_A^n-\snc^B(\la\rp{a} \thetac_{Bi})
\end{align*}
as desired in (\ref{6.8.1.22}).
 \end{proof}

\subsection*{Acknowledgement}
The author would like to thank her colleague, Prof. Gui-Qiang Chen, for raising the enlightening question on the possibility of constructing a complementary counter scenario to the shock formation result in \cite{shock_demetrios} in a learning seminar in 2012 in OXPDE.
 
\end{document}